\documentclass[12pt]{amsbook}
\usepackage{amssymb}
\usepackage[all]{xy}

\usepackage[colorlinks=true,linkcolor=magenta,citecolor=magenta]{hyperref}
\makeindex

\newtheorem{theorem}{Theorem}[chapter]
\newtheorem{addendum}[theorem]{Addendum}
\newtheorem{advertisement}[theorem]{Advertisement}
\newtheorem{answer}[theorem]{Answer}
\newtheorem{bug}[theorem]{Bug}
\newtheorem{comment}[theorem]{Comment}
\newtheorem{conclusion}[theorem]{Conclusion}

\newtheorem{corollary}[theorem]{Corollary}
\newtheorem{definition}[theorem]{Definition}
\newtheorem{disclaimer}[theorem]{Disclaimer}
\newtheorem{exercise}[theorem]{Exercise}
\newtheorem{fact}[theorem]{Fact}
\newtheorem{goal}[theorem]{Goal}
\newtheorem{grandconclusion}[theorem]{Grand conclusion}

\newtheorem{principle}[theorem]{Principle}
\newtheorem{problem}[theorem]{Problem}
\newtheorem{proposition}[theorem]{Proposition}
\newtheorem{question}[theorem]{Question}
\newtheorem{questions}[theorem]{Questions}

\newtheorem{strategy}[theorem]{Strategy}
\newtheorem{thought}[theorem]{Thought}
\newtheorem{update}[theorem]{Update}
\newtheorem{warning}[theorem]{Warning}

\begin{document}

\title{Graphs and their symmetries}

\author{Teo Banica}
\address{Department of Mathematics, University of Cergy-Pontoise, F-95000 Cergy-Pontoise, France. {\tt teo.banica@gmail.com}}


\begin{abstract}
This is an introduction to graph theory, from a geometric and analytic viewpoint. A finite graph $X$ is described by its adjacency matrix $d\in M_N(0,1)$, which can be thought of as being a kind of discrete Laplacian, and we first discuss the basics of graph theory, by using $d$, and various linear algebra tools. Then we discuss the computation of the classical and quantum symmetry groups $G(X)\subset G^+(X)$, which must leave invariant the eigenspaces of $d$, with the quantum symmetry group $G^+(X)$ being in general bigger than the classical symmetry group $G(X)$.
\end{abstract}

\maketitle

\chapter*{Preface}

Graph theory is a wide topic, and there are many possible ways of getting into it. Basically any question which is of discrete nature, with a finite number of objects $N$ interacting between them, has an associated graph $X$, and the question can be reformulated, and hopefully solved, in terms of that graph $X$. A possible exception seems to come from discrete questions which are naturally encoded by a matrix, but isn't a $N\times N$ matrix $A$ some kind of graph too, having $1,\ldots,N$ as vertices, and with any pair of vertices $i,j$ producing the oriented edge $i\to j$, colored by corresponding the matrix entry $A_{ij}$.

\bigskip

Graphs appear as well in connection with continuous questions, via discretization methods. Discretizing is something commonplace in physics, and sometimes in pure mathematics too, with the continuous question, system or manifold $M$ appearing as the $N\to\infty$ limit of its discrete versions, which typically correspond to graphs $X_N$, which are often easier to solve, and with this philosophy having produced good results all around the spectrum, from common sense science and engineering questions up to fairly advanced topics, such as quantum gravity, dark matter and teleportation.

\bigskip

The purpose of this book is to talk about graphs, viewed from this perspective, as being discrete geometry objects. To be more precise, we will be rather mathematical, and the cornerstone of our philosophy will be the commonly accepted fact that the basic objects of mainstream, continuous mathematics are the Riemannian manifolds $M$, together with their Laplacians $\Delta$. In the discrete setting the analogues of these objects are the finite graphs $X$, together with their adjacency matrices $d$, and this will be our way to view the graphs, in this book, as being some kind of ``discrete Riemannian manifolds''.

\bigskip

This was for the philosophy. In practice, this will not prevent us from doing many elementary things, and our purpose will be to talk about graphs in a large sense, ranging from very basic, definition and elementary properties, to more advanced.

\bigskip

More in detail now, we will start in Part I with standard introductory material on graph theory, namely definition, main examples, and some computations and basic theorems. Our presentation here will be for the most of spectral flavor, by translating everything in terms of the adjacency matrix $d$, and proving the results via linear algebra. 

\bigskip

As a continuation of this, in Part II we will discuss various geometric aspects, first with a study of the Laplacian of graphs, and a discussion of related discretization methods, then with a look into trees and other planar graphs, and higher genus graphs as well, and finally with a look into knot invariants, constructed via plane projections.

\bigskip

Switching gears, in Part III we will decide that graph theory is most likely related to abstract algebra, and that what we are mostly interested in is the systematic study of the symmetry groups $G(X)$ of the finite graphs $X$. The material here will be mostly upper undergraduate level, sometimes erring on the graduate side.

\bigskip

Escalating difficulties, in Part IV we will include as well in our discussion the quantum symmetry groups $G^+(X)$ of the same finite graphs $X$, which are in general bigger than their classical counterparts $G(X)$, and with our motivation here coming from theoretical physics, and more specifically from quantum and statistical mechanics.

\bigskip

Part of this book is based on research work that I did some time ago on graphs and their symmetries, and I would like to thank my coworkers. Many thanks go as well to the young researchers in the area, who have recently pushed things to unexpected levels of deepness, that we will try to explain a bit here. Finally, many thanks go to my cats. When it comes to quick search over a graph, they are the kings of this universe.

\bigskip

\

{\em Cergy, October 2024}

\smallskip

{\em Teo Banica}

\baselineskip=15.95pt
\tableofcontents
\baselineskip=14pt

\part{Graph theory}

\ \vskip50mm

\begin{center}
{\em Fly robin fly

Fly robin fly

Fly robin fly

Up up to the sky}
\end{center}

\chapter{Finite graphs}

\section*{1a. Finite graphs}

We will be interested in this book in graph theory, which is the same as discrete geometry. Or perhaps vice versa. Personally I prefer the ``vice versa'' viewpoint, geometry is something basic, coming first, and graphs, which are more specialized, come after. And this, I hope, will agree with you too. I can only imagine that you have some knowledge of geometry, for instance that round thing that you played with as a kid is a sphere, and that's quality geometry, called Riemannian, you're already a bit expert in that. While in what regards graphs, you've probably heard of them, but have no idea what they are good for, and you are here, starting this book, for getting introduced to them.

\bigskip

So, discrete geometry. And the first question here is of course, why? Is there something wrong with the usual, continuous geometry?

\bigskip

Good question, and as a first answer, I would argue that a good quality digital music record is as good as an old vinyl one. Of course some people might claim the opposite, but these people can be proved, scientifically, to be wrong, for the simple reason that the human ear has a certain resolution, and once your digital technology passes that resolution, records are $100\%$ perfect to the human ear. And there is even more, because the same goes for vision, smell, taste and touch. In short, all our senses function like a computer, in a digital way, and we might be well living in a continuous world, but we will never be able to really sense this, and so fully benefit from that continuity.

\bigskip

Some further thinking on this leads to the scary perspective that our world might be actually discrete, at a resolution much finer than that of our human senses. You might argue that why not using some sharp scientific machinery instead, but that machinery has a certain resolution too. And so our scare is justified, and time to ask:

\index{geometry}
\index{quantization}

\begin{question}
Is the world continuous, or quantized?
\end{question}

Here we have used the word ``quantized'', which is a fancy scientific way of saying ``discrete'', as a matter of getting into serious physics. And here, hang on, opinions are split, but most physicists tend to favor the possibility that our world is indeed quantized. And at a resolution that is so fine, that is guaranteed to stay beyond the level of what can be directly observed with past, present and future scientific machinery.

\bigskip

But probably enough talking, time to get to work, remember that we are here for doing math and computations. Let us summarize this discussion by formulating:

\begin{conclusion}
Discrete geometry is worth a study, as a useful discretization of our usual continuous geometry. And why not as a replacement for it, in case it's wrong.
\end{conclusion}

In order now to get started, we first need to talk about continuous geometry. Generally speaking, that is about curves, surfaces and other shapes, called ``manifolds'' in $\mathbb R^N$. However, instead of getting into what a manifold exactly is, which can be a bit technical, let us just take this intuitively, $M$ is by definition a ``continuous shape'' in $\mathbb R^N$.

\bigskip

Now for discretizing such a manifold $M\subset\mathbb R^N$, the idea is very simple, namely placing a sort of net on $M$. To be more precise, assume that we found a finite set of points $X\subset M$, with some edges between them, denoted $i-j$, and corresponding to paths on $M$, all having the same lenght $\varepsilon>0$. Then, we have our discretization $X\subset M$.

\bigskip

As an example here, let us try to discretize a sphere $S\subset\mathbb R^3$. There are many ways of doing so, and a quite straightforward one is by using an inscribed cube, as follows:
$$
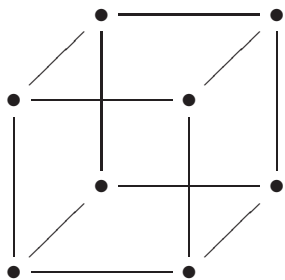
\xymatrix@R=20pt@C=20pt{
&\bullet\ar@{-}[rr]&&\bullet\\
\bullet\ar@{-}[rr]\ar@{-}[ur]&&\bullet\ar@{-}[ur]\\
&\bullet\ar@{-}[rr]\ar@{-}[uu]&&\bullet\ar@{-}[uu]\\
\bullet\ar@{-}[uu]\ar@{-}[ur]\ar@{-}[rr]&&\bullet\ar@{-}[uu]\ar@{-}[ur]
}$$

\smallskip

With this done, let us try now to forget about $M$ itself, and think at $X$, taken alone. This beast $X$ is a finite set of points, with edges between them. In addition, all edges have the same lenght $\varepsilon>0$, but now that $M$ is gone, we can assume if we want that we have $\varepsilon=1$, or simply forget about $\varepsilon$. Thus, we are led into the following definition:

\index{graph}
\index{finite graph}
\index{edges}

\begin{definition}
A graph $X$ is a finite set of points, with certain edges $i-j$ drawn between certain pairs of points $(i,j)$.
\end{definition}

All this might seem overly simplified, and I can hear some of you saying hey, but deep mathematics needs very complicated definitions, that no one really understands, and so on. Wrong. The above definition is in fact the good one, and the graphs themselves are very interesting objects. That will keep us busy, for the rest of this book.

\bigskip

Before getting into mathematics, let us explore a bit Definition 1.3, and see how the graphs look like. We would need here some sort of ``random graph'', in order to see what types of phenomena can appear, and here is an example, which is quite illustrating:
$$
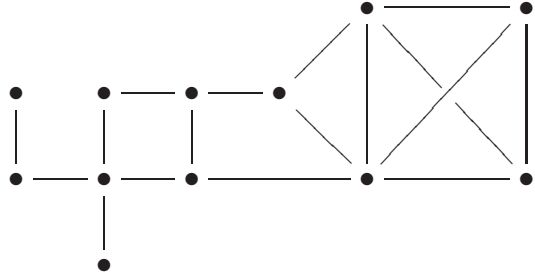
\xymatrix@R=20pt@C=20pt{
&&&&\bullet\ar@{-}[rr]\ar@{-}[dl]\ar@{-}[dd]\ar@{-}[dr]&&\bullet\ar@{-}[dd]\\
\bullet\ar@{-}[d]&\bullet\ar@{-}[r]&\bullet\ar@{-}[d]\ar@{-}[r]&\bullet\ar@{-}[dr]&&\ar@{-}[dr]\\
\bullet\ar@{-}[r]&\bullet\ar@{-}[r]\ar@{-}[u]\ar@{-}[d]&\bullet\ar@{-}[rr]&&\bullet\ar@{-}[rr]\ar@{-}[uurr]&&\bullet\\
&\bullet&
}$$

We can see that there are several things going on here, and here is a list of observations that can be formulated, just by looking at this graph:

\bigskip

(1) First of all, in relation with our discretization philosophy, $X\subset M$, there is a bit of a mess with dimensions here, because the whole middle of the graph seems to discretize a surface, or 2D manifold, but the left part is rather 1D, and the right part, with that crossing edges, suggests either a 2D manifold in $\mathbb R^3$ or higher, or a 3D manifold.

\bigskip

(2) So, this is one problem that we will have to face, when working with Definition 1.3, at that level of generality there is no indication about the ``dimension'' of what our graph is supposed to describe. But this is not an issue, because when this will be really needed, with some further axioms we can divide the graphs into classes and so on.

\bigskip

(3) As another observation now, our graph above is not that ``random'', because it is connected. If this annoys you, please consider that the graph contains as well three extra points, that I can even draw for you, here they are, $\bullet\bullet\bullet$\,, not connected to the others. But will this bring something interesting to our formalism, certainly not.

\bigskip

(4) In short, we will usually not bother with disconnected graphs $X$, a bit like geometers won't bother with disconnected manifolds $M$. This being said, at some point of this book, towards the end, when looking at graphs from a quantum perspective, we will reconsider this, because ``quantum'' allows all sorts of bizarre ``jumps''. More later.

\bigskip

(5) Finally, if you're an electric engineer you might be a bit deceived by all this, because what we're doing so far is obviously binary, and won't allow the installation of bulbs, resistors,  capacitors and so on, having continuous parameters. Criticism accepted, and we will expand our formalism, once our theory using Definition 1.3 will be ripe.

\bigskip

As a last comment, we are of course allowed to draw graphs as we find the most appropriate, and no rules here, free speech as we like it. However, as an important point, if you are a bit familiar with advanced mathematics, and know the difference between geometry and topology, you might have the impression that what we are doing is topology. But this is wrong, what we will be doing is definitely geometry. More on this later.

\bigskip

Let us record the main conclusions from this discussion, as follows:

\index{connected graph}

\begin{conclusion}
Our definition for the graphs looks good, and appears as a good compromise between:
\begin{enumerate}
\item More particular definitions, such as looking at connected graphs only, which is the main case of interest, geometrically speaking.

\item More general definitions, such as allowing machinery and symbols to be installed, either at the vertices, or on the edges, or both.
\end{enumerate}
\end{conclusion}

We will come back to this, when needed, and fine-tune our definition for the graphs, along the above lines, depending on the exact problems that we are interested in.

\bigskip

Moving ahead, now that we have our objects of study, the graphs, time to do some mathematics. For this purpose, we will mostly use linear algebra, and a bit of calculus too. At some point later, we will need as well some basic knowledge of group theory,  along with some basic functional analysis, and basic probability theory. All these can be learned from many places, and if looking for a compact package, talking a bit about everything that will be needed here, you can check my linear algebra book \cite{ba1}.

\bigskip

But let us start with some algebra and generalities. We surely know from Definition 1.3 what a graph is, and normally no need for more, we are good for starting some work. However, it is sometimes useful to have some alternative points of view on this. So, let us replace Definition 1.3 with the following definition, which can be quite useful:

\begin{definition}
A graph $X$ can be viewed as follows:
\begin{enumerate}
\item As a finite set of points $X$, with certain edges $i-j$ drawn between certain pairs of distinct points $(i,j)$.

\item As a finite vertex set $X$, given with an edge set $E\subset X\times X$, which must be symmetric, and avoiding the diagonal.

\item As a finite set $X$, given with a relation $i-j$ on it, which must be non-reflexive, meaning $i\not\!\!-\,i$, and symmetric.
\end{enumerate}
\end{definition}

To be more precise, here the formulation in (1) is the one in our old Definition 1.3, with the remark that we forgot to say there that $i\neq j$, and with this coming from the fact that, geometrically speaking, self-edges $i-i$ look like a pathology, to be avoided.

\bigskip

Regarding now (2), this is something clearly equivalent to (1). It is sometimes convenient to use the notation $X=(V,E)$, with the vertex set denoted $V$, but in what concerns us, we will keep using our policy of calling $X$ both the graph, and the vertex set.

\bigskip

Finally, regarding (3), this is something equivalent to (1) and (2), and with ``relation on $X$'' there not meaning anything in particular, and more specifically, just meaning ``subset of $X\times X$''. Observe that nothing is said about the transitivity of $-$.

\bigskip

All this is nice, and as a first mathematical question for us, let us clarify what happens in relation with the transitivity of $-$. To start with, as our previous examples of graphs indicate, the relation $-$ is not transitive, in general. In fact, $-$ can never be transitive, unless for graphs of type $\bullet\bullet\bullet$\,, without edges at all, because once we have an edge $i-j$, by symmetry and then transitivity we are led to the following wrong conclusion: 
$$i-j\ ,\ j-i\implies i-i$$

However, as a matter of recycling our question, we can ask if $-$, once completed with $i-i$ as to be symmetric, can be transitive. And the answer here is as follows:

\index{transitive relation}
\index{complete graph}
\index{disjoint union}

\begin{proposition}
The graphs $X$ having the property that $-$, once completed with $i-i$ as to be symmetric, is transitive, are exactly the graphs of type
$$
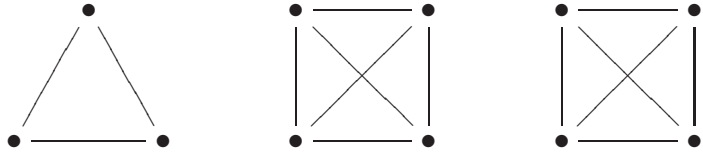
\xymatrix@R=15pt@C=15pt{
&\bullet\ar@{-}[ddl]\ar@{-}[ddr]&&&\bullet\ar@{-}[dd]\ar@{-}[ddrr]\ar@{-}[rr]&&\bullet\ar@{-}[dd]\ar@{-}[ddll]&&\bullet\ar@{-}[dd]\ar@{-}[ddrr]\ar@{-}[rr]&&\bullet\ar@{-}[dd]\ar@{-}[ddll]\\
\\
\bullet\ar@{-}[rr]&&\bullet&&\bullet\ar@{-}[rr]&&\bullet&&\bullet\ar@{-}[rr]&&\bullet
}$$
that is, are the disjoint unions of complete graphs.
\end{proposition} 

\begin{proof}
This is clear, because after thinking a bit, our question simply asks to draw the graph of an equivalence relation, and what you get is of course a picture as above, with various complete graphs corresponding to the various equivalence classes.
\end{proof}

In practice now, all three formulations of the notion of graph from Definition 1.5 can be useful, for certain pruposes, but we will keep using by default the formulation (1) there. Indeed, here is what happens for the cube, and judge yourself:

\index{cube}

\begin{proposition}
The cube graph is as follows:
\begin{enumerate}
\item Viewed as a set, with edges drawn between points:
$$
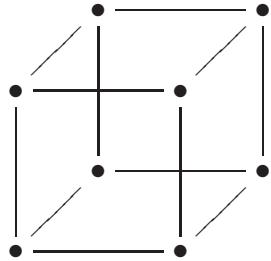
\xymatrix@R=18pt@C=18pt{
&\bullet\ar@{-}[rr]&&\bullet\\
\bullet\ar@{-}[rr]\ar@{-}[ur]&&\bullet\ar@{-}[ur]\\
&\bullet\ar@{-}[rr]\ar@{-}[uu]&&\bullet\ar@{-}[uu]\\
\bullet\ar@{-}[uu]\ar@{-}[ur]\ar@{-}[rr]&&\bullet\ar@{-}[uu]\ar@{-}[ur]
}$$

\item Viewed as a vertex set, plus edge set:
$$X=\{1,2,3,4,5,6,7,8\}$$
$$E=\left\{\begin{matrix}
12,14,15,21,23,26,32,34,37,41,43,48\\
51,56,58,62,65,67,73,76,78,84,85,87
\end{matrix}\right\}$$

\item Viewed as a set, with a relation on it: same as in $(2)$.
\end{enumerate} 
\end{proposition}

\begin{proof}
To start with, (1) is clear. For (2) however, we are a bit in trouble, because in order to figure out what the edge set is, we must first draw the cube, as in (1), so failure with respect to (1), and then label the vertices with numbers, say $1,\ldots,8$:
$$
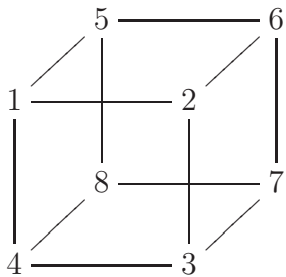
\xymatrix@R=16pt@C=20pt{
&5\ar@{-}[rr]&&6\\
1\ar@{-}[rr]\ar@{-}[ur]&&2\ar@{-}[ur]\\
&8\ar@{-}[rr]\ar@{-}[uu]&&7\ar@{-}[uu]\\
4\ar@{-}[uu]\ar@{-}[ur]\ar@{-}[rr]&&3\ar@{-}[uu]\ar@{-}[ur]
}$$

But with this done, we can then list the 12 edges, say in lexicographic order. Finally, (3) is obviously the same thing as (2), so failure too with respect to (1).
\end{proof}

Summarizing, and hope you agree with me, Definition 1.5 (1), which is more or less our original Definition 1.3, is the best, for general graphs. This being said, the cube example in Proposition 1.7 might look quite harsh, to the point that you might now be wondering if Definition 1.5 (2) and Definition 1.5 (3) have any uses at all.

\bigskip

Good point, and in order to answer, let us go back to the cube. With mea culpa to that cube, we can in fact do better when labeling the vertices, and we have:

\begin{proposition}
The cube can be best viewed, via edge and vertex set, as follows:
\begin{enumerate}
\item The vertices are the usual $3{\rm D}$ coordinates of the vertices of the unit cube, that is, $000$, $001$, $010$, $011$, $100$, $101$, $110$, $111$.

\item The edge set consists of the pairs $(abc,xyz)$ of such coordinates having the property that $abc\leftrightarrow xyz$ comes by modifying exactly one coordinate.
\end{enumerate}
\end{proposition}

\begin{proof}
This is clear from definitions, and no need to draw a cube or anything for this, in fact the geometric cube is there, in what we said in (1) and (2). Here is however an accompanying picture, in case you have troubles in seeing this right away:
$$
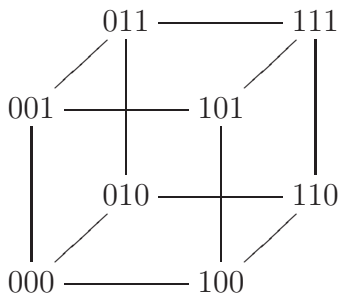
\xymatrix@R=18pt@C=11pt{
&011\ar@{-}[rr]&&111\\
001\ar@{-}[rr]\ar@{-}[ur]&&101\ar@{-}[ur]\\
&010\ar@{-}[rr]\ar@{-}[uu]&&110\ar@{-}[uu]\\
000\ar@{-}[uu]\ar@{-}[ur]\ar@{-}[rr]&&100\ar@{-}[uu]\ar@{-}[ur]
}$$

Thus, one way or another, we are led to the conclusion in the statement.
\end{proof}

As a comment here, we can further build along the above lines, with the ultimate conclusion being something which looks very good and conceptual, as follows:

\index{Cayley graph}

\begin{advertisement}
The cube is the Cayley graph of $\mathbb Z_2^3$. 
\end{advertisement}

Obviously, this looks like a quite deep statement, probably beating everything else that can be said, about the cube graph. We will talk about this later, when discussing in detail the Cayley graphs, but coming a bit in advance, some explanations on this:

\bigskip

(1) You surely know about the cyclic group $\mathbb Z_N$, having as elements the remainders $0,1,\ldots,N-1$ modulo $N$. But at $N=2$ this group is simply $\mathbb Z_2=\{0,1\}$, so in view of what we found in Proposition 1.8, and leaving now 3D geometry aside, for a trip into arithmetic, we can say that the vertices of the cube are the elements $g\in\mathbb Z_2^3$.

\bigskip

(2) Regarding now the edges $g-h$, we know from Proposition 1.8 that these appear precisely at the places where the passage $g\leftrightarrow h$ comes by modifying exactly one coordinate. But, and skipping some details here, with all this left to be discussed later in this book, this means precisely that the cube is the Cayley graph of $\mathbb Z_2^3$.

\bigskip

As an overall conclusion now, we have several equivalent definitions for the graphs, those in Definition 1.5, which can sometimes lead to interesting mathematics, and we have as well Conclusion 1.4, telling us how to fine-tune our graph formalism, when needed, depending on the precise questions that we are interested in. Very nice all this, foundations laid, and we are now good to go, for some tough mathematics on graphs.

\section*{1b. Adjacency matrix}

As a first true result about graphs, and we will call this theorem, not because of the difficulty of its proof, but because of its beauty and usefulness, we have:

\index{adjacency matrix}
\index{symmetric matrix}

\begin{theorem}
A graph $X$, with vertices labeled $1,\ldots,N$, is uniquely determined by its adjacency matrix, which is the matrix $d\in M_N(0,1)$ given by:
$$d_{ij}=\begin{cases}
1&{\rm if}\ i-j\\
0&{\rm if}\ {i\not\!\!-\,j}
\end{cases}$$
Moreover, the matrices $d\in M_N(0,1)$ which can appear in this way, from graphs, are precisely those which are symmetric, and have $0$ on the diagonal.
\end{theorem}

\begin{proof}
We have two things to be proved, the idea being as follows:

\medskip

(1) Given a graph $X$, we can construct a matrix $d\in M_N(0,1)$ as in the statement, and this matrix is obviously symmetric, and has $0$ on the diagonal.

\medskip

(2) Conversely, given a matrix $d\in M_N(0,1)$ which is symmetric, and has $0$ on the diagonal, we can associate to it the graph $X$ having as vertices the numbers $1,\ldots,N$, and with the edges between these vertices being defined as follows:
$$i-j\iff d_{ij}=1$$

It is then clear that the adjacency matrix of this graph $X$ is the matrix $d$ itself. Thus, we have established a correspondence $X\leftrightarrow d$, as in the statement.
\end{proof}

The above result is very useful for various purposes, a first obvious application being that we can now tell a computer what our favorite graph $X$ is, just by typing in the corresponding adjacency matrix $d\in M_N(0,1)$, which is something that the computer will surely understand. In fact, computers like 0-1 data, that's the language that they internally speak, and when that data comes in matrix form, they even get very happy and excited, and start doing all sorts of computations. But more on this later.

\bigskip

Speaking conversion between graphs $X$ and matrices $d$, this can work as well for the two of us. Here is my matrix $d$, and up to you to tell me what my graph $X$ was:
$$d=\begin{pmatrix}
0&1&1\\
1&0&1\\
1&1&0
\end{pmatrix}$$

Problem solved I hope, with the graph in question being a triangle. Here is now another matrix $d$, and in the hope that you will see here a square: 
$$d=\begin{pmatrix}
0&1&0&1\\
1&0&1&0\\
0&1&0&1\\
1&0&1&0
\end{pmatrix}$$

Here is now another matrix $d$, and in the hope that you will see the same square here, but this time with both its diagonals drawn, and you can call that tetrahedron too: 
$$d=\begin{pmatrix}
0&1&1&1\\
1&0&1&1\\
1&1&0&1\\
1&1&1&0
\end{pmatrix}$$

By the way, talking exercises, please allow me to record the solutions to the above exercises, not of course for you, who are young and enthusiastic and must train hard, but rather for me and my colleagues, who are often old and tired. Here they are:
$$
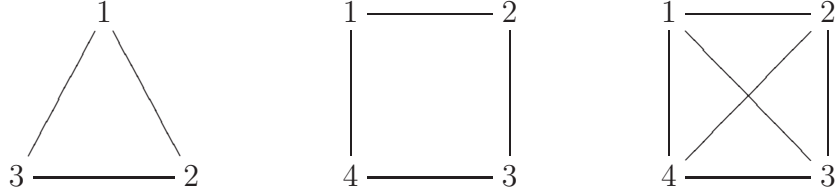
\xymatrix@R=20pt@C=20pt{
&1\ar@{-}[ddl]\ar@{-}[ddr]&&&1\ar@{-}[dd]\ar@{-}[rr]&&2\ar@{-}[dd]&&1\ar@{-}[dd]\ar@{-}[ddrr]\ar@{-}[rr]&&2\ar@{-}[dd]\ar@{-}[ddll]\\
\\
3\ar@{-}[rr]&&2&&4\ar@{-}[rr]&&3&&4\ar@{-}[rr]&&3
}$$

And to finish this discussion, at a more advanced level now, here is another matrix $d$, and in the hope that you will see something related to cycling here: 
$$d=\begin{pmatrix}
0&1&0&1&0&1\\
1&0&1&0&1&0\\
0&1&0&1&0&1\\
1&0&1&0&1&0\\
0&1&0&1&0&1\\
1&0&1&0&1&0
\end{pmatrix}$$

All this looks fun, and mathematically relevant too. Based on this, we will agree to sometimes speak directly in terms of $d$, which is quite practical. Especially for me, typing a matrix in Latex, the computer language used for typing math, and the present book, being easier than drawing a graph. Plus our computer friend will understand too.

\bigskip

Let us develop now some theory, for the general graphs. Once the number of vertices is fixed, $N\in\mathbb N$, what seems to most distinguish the graphs, for instance in connection with the easiness or pain of drawing them, is the total number of edges $|E|$. Equivalently, what distinguishes the graphs is the density of edges at each vertex, given by:
$$\rho=\frac{|E|}{N}$$

So, let us take a closer look at this quantity. It is convenient, for a finer study, to formulate our definition as follows:

\index{valence}
\index{regular graph}
\index{number of neighbors}

\begin{definition}
Given a graph $X$, with $X$ standing as usual for both the vertex set, and the graph itself, the valence of each vertex $i$ is the number of its neighbors:
$$v_i=\#\left\{j\in X\Big|i-j\right\}$$
We call $X$ regular when the valence function $v:X\to\mathbb N$ is constant. More specifically, in this case, we say that $X$ is $k$-regular, with $k$ being the common valence of vertices.
\end{definition}

At the level of examples, all the graphs pictured above, with our usual convention that 0-1 matrix means picture, are regular, with the valence being as follows:
$$\begin{matrix}
{\rm Graph}&{\rm Regularity}&{\rm Valence}\\
&-\\
{\rm triangle}&{\rm yes}&2\\
{\rm square}&{\rm yes}&2\\
{\rm tetrahedron}&{\rm yes}&3\\
{\rm hexagonal\ wheel}&{\rm yes}&3
\end{matrix}$$

Which leads us into the question, is there any interesting non-regular graph? Certainly yes, think for instance at all that beautiful pictures of snowflakes, such as:
$$
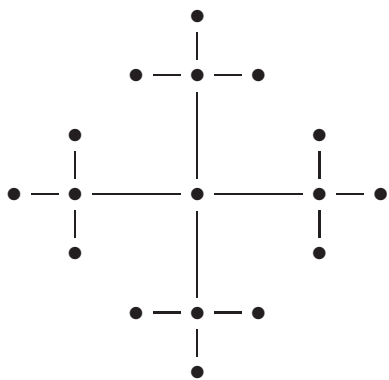
\xymatrix@R=10pt@C=10pt{
&&&\bullet\ar@{-}[d]\\
&&\bullet\ar@{-}[r]&\bullet\ar@{-}[dd]\ar@{-}[r]&\bullet\\
&\bullet\ar@{-}[d]&&&&\bullet\ar@{-}[d]\\
\bullet\ar@{-}[r]&\bullet\ar@{-}[rr]&&\bullet\ar@{-}[rr]\ar@{-}[dd]&&\bullet\ar@{-}[r]&\bullet\\
&\bullet\ar@{-}[u]&&&&\bullet\ar@{-}[u]\\
&&\bullet\ar@{-}[r]&\bullet\ar@{-}[d]\ar@{-}[r]&\bullet\\
&&&\bullet
}$$

Here the valence is 4 at the inner points, and 1 at the endpoints. Graphs as the above one are called in mathematics ``trees'', due to the fact that, thinking a bit, if you look at a usual tree from the above, say from a helicopter, what you see is a picture as above. But more on such graphs, which can be quite complicated, later.

\bigskip

For the moment, let us study the valence function, and the regular graphs. We can do some math, by using the adjacency matrix, as follows:

\index{all-one vector}
\index{stochastic matrix}
\index{valence function}
\index{row sums}

\begin{proposition}
Given a graph $X$, with adjacency matrix $d$:
\begin{enumerate}
\item The valence function $v$ is the row sum function for $d$.

\item Equivalently, $v=d\xi$, with $\xi$ being the all-$1$ vector.

\item $X$ is regular when $d$ is stochastic, meaning with constant row sums.

\item Equivalently, $X$ is $k$-regular, with $k\in\mathbb N$, when $d\xi=k\xi$.
\end{enumerate}
\end{proposition}

\begin{proof}
All this looks quite trivial, but let us make some effort, and write down a complete mathematical proof. Regarding (1), this follows from:
\begin{eqnarray*}
v_i
&=&\sum_{j-i}1\\
&=&\sum_{j-i}1+\sum_{j\not-i}0\\
&=&\sum_jd_{ij}
\end{eqnarray*}

Regarding (2), this follows from (1), because for any matrix $d$ we have:
$$d\xi=\begin{pmatrix}
d_{11}&\ldots&d_{1N}\\
\vdots&&\vdots\\
d_{N1}&\ldots&d_{NN}
\end{pmatrix}
\begin{pmatrix}
1\\
\vdots\\
1\end{pmatrix}
=\begin{pmatrix}
d_{11}+\ldots+d_{1N}\\
\vdots\\
d_{N1}+\ldots+d_{NN}\end{pmatrix}$$

Finally, (3) follows from (1), and (4) follows from (2).
\end{proof}

Observe that, in linear algebra terms, (4) above reformulates as follows:

\index{eigenvalue}

\begin{theorem}
A graph is regular precisely when the all-$1$ vector is an eigenvector of the adjacency matrix. In this case, the valence is the corresponding eigenvalue.
\end{theorem}

\begin{proof}
This is clear indeed from what we know from Proposition 1.12 (4). As a philosophical comment here, you might wonder how something previously labeled Proposition can suddenly become a Theorem. Welcome to mathematics, which is not an exact science, the general rule being that everything that looks fancy enough can be called Theorem. In fact, we have already met this phenomenon in the context of Proposition 1.8, which got converted into Advertisement 1.9, and with that advertisement being, you guessed right, an advertisement for a future Theorem, saying exactly the same thing.
\end{proof}

The above result is quite interesting, and suggests systematically looking at the eigenvalues and eigenvectors of $d$. Which is indeed a very good idea, but we will keep this for later, starting with chapter 2, once we will have a bit more general theory.

\section*{1c. Walks on graphs} 

Have you ever played chess, or simply observed a cat in action. There are all sorts of paths and combinations that can be followed, and the problem is that of quickly examining all these possibilities, with the series of conclusions being typically something on type damn, damn, damn, damn, damn, yes. With the yes followed by some action.

\bigskip

Mathematically speaking, this brings us into walks on graphs. And here, as usual, we get into the adjacency matrix $d$, thanks to the following key result:

\index{paths on graph}
\index{power of matrix}

\begin{theorem}
Given a graph $X$, with adjacency matrix $d\in M_N(0,1)$, we have:
$$(d^k)_{ij}=\#\Big\{i=i_0-i_1-\ldots-i_{k-1}-i_k=j\Big\}$$
That is, the $k$-th power of $d$ describes the length $k$ paths on $X$.
\end{theorem}

\begin{proof}
According to the usual rule of matrix multiplication, the formula for the powers of the adjacency matrix $d\in M_N(0,1)$ is as follows:
\begin{eqnarray*}
(d^k)_{i_0i_k}
&=&\sum_{i_1,\ldots,i_{k-1}}d_{i_0i_1}d_{i_1i_2}\ldots d_{i_{k-1}i_k}\\
&=&\sum_{i_1,\ldots,i_{k-1}}\delta_{i_0-i_1}\delta_{i_1-i_2}\ldots\delta_{i_{k-1}-i_k}\\
&=&\sum_{i_1,\ldots,i_{k-1}}\delta_{i_0-i_1-\ldots-i_{k-1}-i_k}\\
&=&\#\Big\{i_0-i_1-\ldots-i_{k-1}-i_k\Big\}
\end{eqnarray*}

Thus, we are led to the conclusion in the statement.
\end{proof}

Of particular interest are the paths which begin and end at the same point. These are called loops, and in the case of loops, Theorem 1.14 particularizes as follows:

\index{loops on graph}
\index{diagonalization}

\begin{theorem}
Given a graph $X$, with adjacency matrix $d\in M_N(0,1)$, we have:
$$(d^k)_{ii}=\#\Big\{k-{\rm loops\ based\ at\ }i\in I\Big\}$$
Also, the total number of $k$-loops on $X$, at various vertices, is the number
$$Tr(d^k)=\sum_i(d^k)_{ii}$$
which can be computed by diagonalizing $d$.
\end{theorem}

\begin{proof}
There are several things going on here, the idea being as follows:

\medskip

(1) The first assertion follows from Theorem 1.14, which at $i=j$ gives the following formula, which translates into the first formula in the statement:
$$(d^k)_{ii}=\#\Big\{i=i_0-i_1-\ldots-i_{k-1}-i_k=i\Big\}$$

(2) Regarding now the second assertion, this follows from the first one, simply by summing over all the vertices $i\in X$, which gives, as desired:
$$Tr(d^k)=\sum_i\#\Big\{k-{\rm loops\ based\ at\ }i\in I\Big\}$$

(3) Finally, the third assertion is something well-known from linear algebra, the idea being that once we diagonalize a matrix, $A=PDP^{-1}$, we have:
\begin{eqnarray*}
A=PDP^{-1}
&\implies&A^k=PD^kP^{-1}\\
&\implies&Tr(A^k)=Tr(D^k)
\end{eqnarray*}

Thus, back now to our graph case, if we denote by $\lambda_1,\ldots,\lambda_N\in\mathbb R$ the eigenvalues of $d$, then the formula of the trace of $d^k$, that we are interested in, is as follows:
$$Tr(d^k)=\lambda_1^k+\ldots+\lambda_N^k$$

Here we have used the fact that $d$, which is a real symmetric matrix, is indeed diagonalizable, and with real eigenvalues. You might perhaps know this from linear algebra, and if not do not worry, we will discuss all this in detail in chapter 2.
\end{proof}

All the above is quite interesting, and adds to our investigations from the previous section, suggesting as well to get into the diagonalization question for the adjacency matrix $d\in M_N(0,1)$. We will do this as soon as possible, and more precisely starting from chapter 2 below. But before that, we still have to discuss some examples for all this, plus a number of other topics of general interest, not needing diagonalization.

\bigskip

Let us start with the examples. Unfortunately, here everything is quite complicated, because even for very simple graphs $X$, go count the loops directly, via recurrences and so on, that is a lot of combinatorics, which is invariably non-trivial. Or go compute the powers $d^k$ without diagonalizing $d$, that is a lot of heavy work too.

\bigskip

So, as a first philosophical question, what are the simplest graphs $X$, that we can try to do some loop computations for? And here, we have 3 possible answers, as follows:

\index{circle graph}
\index{segment graph}

\begin{fact}
The following are graphs $X$, with a distinguished vertex $0\in X$:
\begin{enumerate}
\item The circle graph, having $N$ vertices, with $0$ being one of the vertices.

\item The segment graph, having $N$ vertices, with $0$ being the vertex at left.

\item The segment graph, having $2N+1$ vertices, with $0$ being in the middle.
\end{enumerate}
\end{fact}

So, let us start with these. For the circle, the computations are quite non-trivial, and you can try doing some, in order to understand what I am talking about. The problem comes from the fact that loops of length $k=0,2,4,6,\ldots$ are quite easy to count, but then, once we pass $k=N$, the loops can turn around the circle or not, and they can even turn several times, and so on, and all this makes the count too complicated. In addition, again due to loop turning, when $N$ is odd, we have as well loops of odd length.

\bigskip

As for the two segment graphs, here the computations look again complicated, and even more complicated than for the circle, because, again, once we pass $k=N$ many things can happen, and this makes the count too complicated. And here, again you can try doing some computations, in order to understand what I am talking about. 

\bigskip

So, shall we give up, in waiting for more advanced techniques, say coming from diagonalization? This would be a wise decision, but before that, let us pull an analysis trick, and formulate the following result, which is of course something informal, and modest:

\begin{theorem}
For the circle graph, having $N$ vertices, the number of length $k$ loops based at one of the vertices is approximately
$$L_k\simeq\frac{2^k}{N}$$
in the $k\to\infty$ limit, when $N$ is odd, and is approximately
$$L_k\simeq\begin{cases}
\frac{2^{k+1}}{N}&(k\ {\rm even})\\
0&(k\ {\rm odd})
\end{cases}$$
also with $k\to\infty$, when $N$ is even. However, in what regards the two segment graphs, we can expect here things to be more complicated.
\end{theorem}

\begin{proof}
This is something not exactly trivial, and with the way the statement is written, which is clearly informal, witnessing for that. The idea is as follows:

\medskip

(1) Consider the circle graph $X$, with vertices denoted $0,1,\ldots,N-1$. Since each vertex has valence 2, any length $k$ path based at 0 will consist of a binary choice at the beginning, then another binary choice afterwards, and so on up to a $k$-th binary choice at the end. Thus, there is a total of $2^k$ such paths, based at 0, and having length $k$. 

\medskip

(2) But now, based on the obvious ``uniformity'' of the circle, we can argue that, in the $k\to\infty$ limit, the endpoint of such a path will become random among the vertices $0,1,\ldots,N-1$. Thus, if we want this endpoint to be 0, as to have a loop, we have $1/N$ chances for this to happen, so the total number of loops is $L_k\simeq 2^k/N$, as stated.

\medskip

(3) With the remark, however, that the above argument works fine only when $N$ is odd. Indeed, when $N$ is even, the endpoint of a length $k$ path will be random among $0,2,\ldots,2N-2$ when $k$ is even, and random among $1,3,\ldots,2N-1$ when $k$ is odd. Thus for getting a loop we must assume that $k$ is even, and in this case the number of such loops is the total number of length $k$ paths, namely $2^k$, approximately divided by $N/2$, the number of points in $\{0,2,\ldots,2N-2\}$, which gives $L_k=2^k/(N/2)$, as stated.

\medskip

(4) All this was of course a bit borderline, I know, with respect to what rigorous mathematics is supposed to be, but honestly, I think that the argument is there, and good, in short I trust this proof. Needless to say, we will be back to all this later, with some better tools for attacking such problems, and with full rigor, at that time.

\medskip

(5) Moving ahead now to the segment graphs, it is pretty much clear that for both, we lack the ``uniformity'' needed in (2), and this due to the 2 endpoints of the segment. In fact, thinking well, these graphs are no longer 2-valent, again due to the 2 endpoints, each having valence 1, and so even (1) must be fixed. And so, we will stop here.
\end{proof}

All this is obviously not very good news, and so again, as question, shall we give up, in waiting for more advanced techniques, say coming from diagonalization? 

\bigskip

Well, instead of giving up, let us look face-to-face at the difficulties that we met. We are led this way, after analyzing the situation, to the following thought:

\begin{thought}
The difficulties that we met, with the circle and the two segments, come from the fact that our loops are not ``free to move'',
\begin{enumerate}
\item for the circle, because these can circle around the circle,

\item for the segments, obviously because of the endpoints,
\end{enumerate}
and so our difficulties will dissapear, and we will be able to do our exact loop count, once we find a graph $X$ where the loops are truly ``free to move''.
\end{thought}

Thinking some more, all this definitely buries the first interval graph, where the vertex $0$ is one of the endpoints. However, we can still try to recycle the circle, by unwrapping it, or extend our second interval graph up to $\infty$. But in both cases what we get is the graph $\mathbb Z$ formed by the integers. So, let us formulate the following definition:

\index{infinite graph}

\begin{definition}
An infinite graph is the same thing as a finite graph, but now with an infinity of vertices, $|X|=\infty$. As a basic example, we have $\mathbb Z$. We also have $\mathbb N$.
\end{definition}

Leaving aside now $\mathbb N$, which looks more complicated, let us try to count the length $k$ paths on $\mathbb Z$, based at $0$. At $k=1$ we have $2$ such paths, ending at $-1$ and $1$, and the count results can be pictured as follows, with everything being self-explanatory:
$$
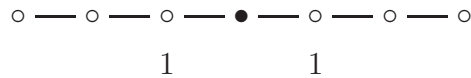
\xymatrix@R=5pt@C=15pt{
\circ\ar@{-}[r]&\circ\ar@{-}[r]&\circ\ar@{-}[r]&\bullet\ar@{-}[r]&\circ\ar@{-}[r]&\circ\ar@{-}[r]&\circ\\
&&1&&1
}$$

At $k=2$ now, we have 4 paths, one of which ends at $-2$, two of which end at 0, and one of which ends at 2. The results can be pictured as follows:
$$
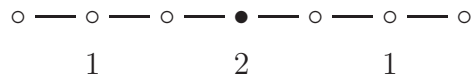
\xymatrix@R=5pt@C=15pt{
\circ\ar@{-}[r]&\circ\ar@{-}[r]&\circ\ar@{-}[r]&\bullet\ar@{-}[r]&\circ\ar@{-}[r]&\circ\ar@{-}[r]&\circ\\
&1&&2&&1
}$$

At $k=3$ now, we have 8 paths, the distribution of the endpoints being as follows:
$$
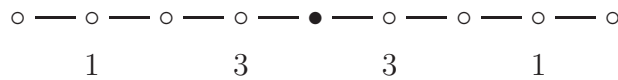
\xymatrix@R=5pt@C=15pt{
\circ\ar@{-}[r]&\circ\ar@{-}[r]&\circ\ar@{-}[r]&\circ\ar@{-}[r]&\bullet\ar@{-}[r]&\circ\ar@{-}[r]&\circ\ar@{-}[r]&\circ\ar@{-}[r]&\circ\\
&1&&3&&3&&1
}$$

As for $k=4$, here we have 16 paths, the distribution of the endpoints being as follows:
$$
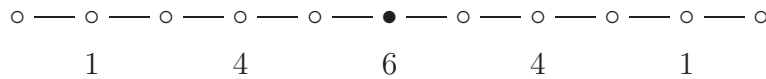
\xymatrix@R=5pt@C=15pt{
\circ\ar@{-}[r]&\circ\ar@{-}[r]&\circ\ar@{-}[r]&\circ\ar@{-}[r]&\circ\ar@{-}[r]&\bullet\ar@{-}[r]&\circ\ar@{-}[r]&\circ\ar@{-}[r]&\circ\ar@{-}[r]&\circ\ar@{-}[r]&\circ\\
&1&&4&&6&&4&&1
}$$

And good news, we can see in the above the Pascal triangle. Thus, eventually, we found the simplest graph ever, namely $\mathbb Z$, and we have the following result about it:

\index{binomial coefficient}
\index{Pascal triangle}
\index{central binomial coefficient}
\index{Stirling formula}
\index{paths on Z}

\begin{theorem}
The paths on $\mathbb Z$ are counted by the binomial coefficients. In particular, the $2k$-paths based at $0$ are counted by the central binomial coefficients,
$$\binom{2k}{k}\simeq\frac{4^k}{\sqrt{\pi k}}$$
with the estimate, in the $k\to\infty$ limit, coming from the Stirling formula.
\end{theorem}

\begin{proof}
This basically follows from the above discussion, as follows:

\medskip

(1) In what regards the count, we certainly have the Pascal triangle, as discovered above, and the rest is just a matter of finishing. There are many possible ways here, a straightforward one being that of arguing that the number $C_k^l$ of length $k$ loops $0\to l$  is subject, due to the binary choice at the end, to the following recurrence relation:
$$C_k^l=C_{k-1}^{l-1}+C_{k-1}^{l+1}$$

But this is exactly the recurrence for the Pascal triangle, so done with the count. 

\medskip

(2) In what regards the estimate, this follows indeed from Stirling, as follows:
\begin{eqnarray*}
\binom{2k}{k}
&=&\frac{(2k)!}{k!k!}\\
&\simeq&\left(\frac{2k}{e}\right)^{2k}\sqrt{4\pi k}\times\left(\frac{e}{k}\right)^{2k}\frac{1}{2\pi k}\\
&=&\frac{4^k}{\sqrt{\pi k}}
\end{eqnarray*}

Thus, we are led to the conclusions in the statement.
\end{proof}

Not bad all this. We will be back to other graphs, such as $\mathbb N$, which is still left, or the circle, or the two segments, and many more, later on. But before that, we still have to discuss a number of other topics of general interest, including coloring.

\section*{1d. Coloring questions}

There are actually two ways of coloring a graph, either by coloring the vertices, or by coloring the edges. Usually mathematicians are quite excited about coloring the vertices, and more on this in a moment. However, electric engineers for instance are more excited about coloring the edges, with machinery like bulbs, resistors and capacitors, or simply with arrows, describing the direction of the current through each edge. Although for machinery having several legs, such as transistors, we are back to vertex coloring.

\bigskip

In this book we will use the edge coloring convention, and this actually for rather mathematical reasons. To be more precise, recall from the very beginning of this book that we are philosophically interested in discretizing manifolds, $X\subset M$, and in fact this is how we run into graphs, by placing a net over such a manifold $M$. 

\bigskip

But now imagine that our net is allowed to stretch. In this case no more graph, what we have is simply a subset $X\subset M$, which technically is a finite metric space. Thus, beyond usual graphs, what we would mostly like to cover with our formalism are the finite metric spaces $X$. But this can be done by calling ``colored graph'' a graph with the vertices still uncolored, and with the edges colored. Indeed, any finite metric space $X$, consisting of $N$ points, can be viewed as the $N$-simplex, with each edge $i-j$ colored by its length $d_{ij}>0$. So, this will be our colored graph formalism in this book:

\index{colored graph}
\index{color set}
\index{finite metric space}
\index{metric space}

\begin{definition}
A colored graph is a usual graph $X$, which each edge $i-j$ colored by a symbol $d_{ij}\in C$, with $C$ being a set, called color set. We call the matrix
$$d\in M_N(C)$$
the adjacency matrix of such a colored graph $X$.
\end{definition}

In practice, as already mentioned in the above, we are mostly interested in the case $C=\{-1,0,1\}$, as to cover the oriented graphs, and also in the case $C=\mathbb R_+$, as to cover the finite metric spaces. However, the best is to use Definition 1.21 as stated, with an abstract color set $C$, that will quite often be a set of reals, $C\subset \mathbb R$.

\bigskip

As a first statement, regarding our enlarged graph formalism, we have:

\index{oriented graph}
\index{triangle inequality}

\begin{theorem}
The following are covered by our colored graph formalism:
\begin{enumerate}
\item The usual graphs. In fact, the usual graphs are precisely the colored versions of the simplex, with color set $C=\{0,1\}$.

\item The oriented graphs. In fact, the oriented graphs are precisely the colored versions of the simplex, with color set $C=\{-1,0,1\}$.

\item The finite metric spaces. These are the colored versions of the simplex, with color set $C=\mathbb R_+$, and with the coloring subject to the triangle inequality.
\end{enumerate}
\end{theorem}

\begin{proof}
All this is trivial, and self-explanatory, and the reason why we called this Theorem instead of Proposition only comes from its theoretical importance.
\end{proof}

As an interesting remark now, with colored graphs we are in fact not that far from the usual graphs, and this due to the following simple fact:

\index{formal matrix}
\index{color decomposition}
\index{color components}

\begin{theorem}
Any formal matrix $d\in M_N(C)$ has a color decomposition,
$$d=\sum_{c\in C}cd_c$$
with the color components $d_c\in M_N(0,1)$ being constructed as follows:
$$(d_c)_{ij}=\begin{cases}
1&{\rm if}\ d_{ij}=c\\
0&{\rm if}\ d_{ij}\neq c
\end{cases}$$
If $X$ is the colored graph having adjacency matrix $d$, these matrices $d_c$ are the adjacency matrices of the color components of $X$, which are usual graphs $X_c$.
\end{theorem}

\begin{proof}
As before with Theorem 1.22, this is something trivial, and self-explanatory, and called Theorem instead of Proposition just because its theoretical importance.
\end{proof}

Many things can be said about colored graphs, and especially about the oriented graphs and finite metric spaces from Theorem 1.22. We will be back to them, later.

\bigskip

All this being said, coloring the vertices of a graph $X$, and we will call such a beast a ``vertex-colored graph'', is something quite interesting too. Let us formulate:

\index{vertex-colored graph}

\begin{definition}
A vertex-colored graph is a usual graph $X$, with each vertex $i\in X$ colored by a symbol $c_i\in C$, with $C$ being a set, called color set. 
\end{definition}

As a first question regarding such vertex-colored graphs, pick whatever black and white geographical map, with countries and boundaries between them, then pick 4 colored pencils, and try coloring that map. Can you do that? What is your algorithm?

\bigskip

Well, the point here is that we have a deep theorem, as follows:

\index{4-color theorem}
\index{coloring maps}

\begin{theorem}
Any map can be colored with $4$ colors.
\end{theorem}

\begin{proof}
This is something non-trivial, which took mankind a lot of time to prove, and whose final proof is something very long, and with computers involved, at several key places. So, theorem coming without proof, and sorry for this. But, we will be back to this in chapter 7 below, when discussing planar graphs, with at least some explanations on what the word ``map'' in the statement exactly means, mathematically speaking.
\end{proof}

As a conclusion now, it looks like with our graph formalism, and our adjacency matrix technology, sometimes modified a bit, as above, we can deal with everything in discrete mathematics. However, before moving ahead with more mathematics, a warning:

\begin{warning}
Not everything in discrete mathematics having a picture, and an adjacency matrix, is a graph.
\end{warning}

This is something quite important, worth discussing in detail, with a good example. So, have you heard about projective geometry? In case you didn't yet, the general principle is that ``this is the wonderland where parallel lines cross''. Which might sound a bit crazy, and not very realistic, but take a picture of some railroad tracks, and look at that picture. Do that parallel railroad tracks cross, on the picture? Sure they do. So, we are certainly not into abstractions here, but rather into serious science.

\bigskip

Mathematically now, here are some axioms, to start with:

\index{projective space}

\begin{definition}
A projective space is a space consisting of points and lines, subject to the following conditions:
\begin{enumerate}
\item Each $2$ points determine a line.

\item Each $2$ lines cross, on a point.
\end{enumerate}
\end{definition}

As a basic example we have the usual projective space $P^2_\mathbb R$, which is best seen as being the space of lines in $\mathbb R^3$ passing through the origin. To be more precise, let us call each of these lines in $\mathbb R^3$ passing through the origin a ``point'' of $P^2_\mathbb R$, and let us also call each plane in $\mathbb R^3$ passing through the origin a ``line'' of $P^2_\mathbb R$. Now observe the following:

\bigskip

(1) Each $2$ points determine a line. Indeed, 2 points in our sense means 2 lines in $\mathbb R^3$ passing through the origin, and these 2 lines obviously determine a plane in $\mathbb R^3$ passing through the origin, namely the plane they belong to, which is a line in our sense.

\bigskip

(2) Each $2$ lines cross, on a point. Indeed, 2 lines in our sense means 2 planes in $\mathbb R^3$ passing through the origin, and these 2 planes obviously determine a line in $\mathbb R^3$ passing through the origin, namely their intersection, which is a point in our sense.

\bigskip

Thus, what we have is a projective space in the sense of Definition 1.27. More generally now, we can perform in fact this construction over an arbitrary field, as follows:

\index{finite field}

\begin{theorem}
Given a field $F$, we can talk about the projective space $P^2_F$, as being the space of lines in $F^3$ passing through the origin, having cardinality 
$$|P^2_F|=q^2+q+1$$
where $q=|F|$, in the case where our field $F$ is finite.
\end{theorem}

\begin{proof}
This is indeed clear from definitions, with the cardinality coming from:
$$|P^2_F|
=\frac{|F^3-\{0\}|}{|F-\{0\}|}
=\frac{q^3-1}{q-1}
=q^2+q+1$$

Thus, we are led to the conclusions in the statement.
\end{proof}

\index{Fano plane}

As an example, let us see what happens for the simplest finite field that we know, namely $F=\mathbb F_2$. Here our projective space, having $4+2+1=7$ points, and 7 lines, is a famous combinatorial object, called Fano plane, which looks as follows:
$$
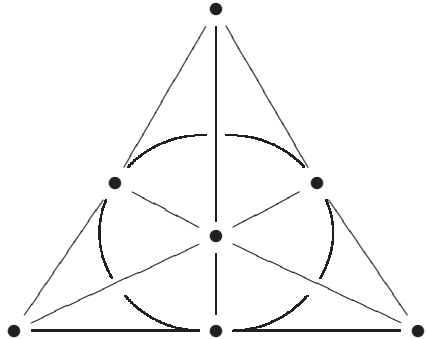
\xymatrix@R=8pt@C=9pt{
&&&&\bullet\ar@{-}[ddddd]\\
&&&&&&&\\
&&&&&&&\\
&&&&\ar@{-}@/^/[drr]\\
&&\bullet\ar@{-}[uuuurr]\ar@{-}@/^/[urr]\ar@{-}@/_/[dd]&&&&\bullet\ar@{-}[uuuull]&&\\
&&&&\bullet\ar@{-}[urr]\ar@{-}[ull]&&&&\\
&&\ar@{-}@/_/[drr]&&&&\ar@{-}@/^/[dll]\ar@{-}@/_/[uu]&&&\\
\bullet\ar@{-}[uuurr]\ar@{-}[rrrr]\ar@{-}[uurrrr]&&&&\bullet\ar@{-}[rrrr]\ar@{-}[uu]&&&&\bullet\ar@{-}[uuull]\ar@{-}[uullll]
}$$

Here the circle in the middle is by definition a line, and with this convention, the projective geometry axioms from Definition 1.27 are satisfied, in the sense that any two points determine a line, and any two lines determine a point. And isn't this magic. 

\bigskip

The question is now, in connection with Warning 1.26, as follows:

\begin{question}
What type of beast is the above Fano plane, with respect to graph theory, and its generalizations?
\end{question}

To be more precise, the Fano plane looks like a graph, but is not a graph, because, save for the circle in the middle and its precise conventions, it matters whether edges are aligned or not, and with this being the whole point with the Fano plane. 

\bigskip

However, not giving up, let us try to investigate this plane, inspired from what we know about graphs. A first interesting question is that of suitably labeling its vertices, and here we know from Theorem 1.28 and its proof that these are naturally indexed by the elements of the group $\mathbb Z_2^3$, with the neutral element $0=000$ excluded:
$$\mathbb Z_2^3-\{000\}=\{001,010,011,100,101,110,111\}$$ 

However, there is some sort of too much symmetry between these symbols, which will not help us much, so instead let us just label the vertices $1,2,\ldots,7$, as follows:
$$
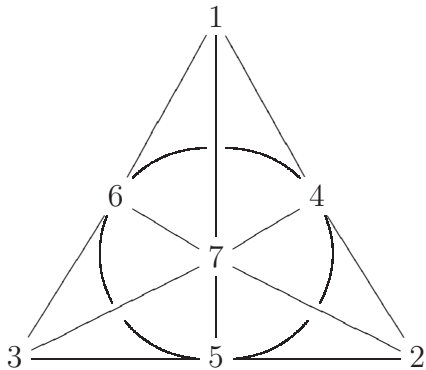
\xymatrix@R=8pt@C=9pt{
&&&&1\ar@{-}[ddddd]\\
&&&&&&&\\
&&&&&&&\\
&&&&\ar@{-}@/^/[drr]\\
&&6\ar@{-}[uuuurr]\ar@{-}@/^/[urr]\ar@{-}@/_/[dd]&&&&4\ar@{-}[uuuull]&&\\
&&&&7\ar@{-}[urr]\ar@{-}[ull]&&&&\\
&&\ar@{-}@/_/[drr]&&&&\ar@{-}@/^/[dll]\ar@{-}@/_/[uu]&&&\\
3\ar@{-}[uuurr]\ar@{-}[rrrr]\ar@{-}[uurrrr]&&&&5\ar@{-}[rrrr]\ar@{-}[uu]&&&&2\ar@{-}[uuull]\ar@{-}[uullll]
}$$

Next comes the question of suitably labeling the lines, and with here, we insist, these lines being not ``edges'', because the Fano plane is not a graph, as explained above. This is again an interesting question, of group theory flavor, and among the conclusions that we can come upon, by thinking at all this, we have the quite interesting fact that, when interchanging the 7 points and the 7 lines, the Fano plane stays the same.

\bigskip

Which looks like something quite deep, but which also teaches us that the labeling question for the lines is the same as the labeling question for the points, and so, in view of the above, that we have to give up. So, let us simply label the lines by letters $a,b,\ldots,g$, in a somewhat random fashion, a bit as we did for the vertices, as follows:
$$
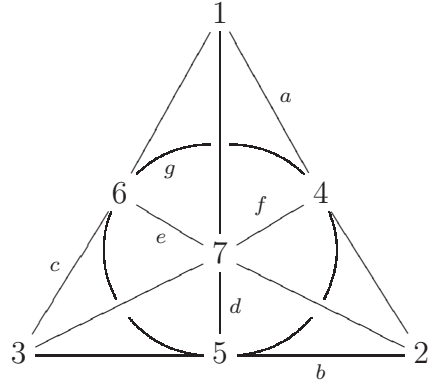
\xymatrix@R=8pt@C=9pt{
&&&&1\ar@{-}[ddddd]\\
&&&&&&&\\
&&&&&&&\\
&&&&\ar@{-}@/^/[drr]\\
&&6\ar@{-}[uuuurr]\ar@{-}@/^/[urr]_g\ar@{-}@/_/[dd]&&&&4\ar@{-}[uuuull]_a&&\\
&&&&7\ar@{-}[urr]^f\ar@{-}[ull]^e&&&&\\
&&\ar@{-}@/_/[drr]&&&&\ar@{-}@/^/[dll]\ar@{-}@/_/[uu]&&&\\
3\ar@{-}[uuurr]^c\ar@{-}[rrrr]\ar@{-}[uurrrr]&&&&5\ar@{-}[rrrr]_b\ar@{-}[uu]_d&&&&2\ar@{-}[uuull]\ar@{-}[uullll]
}$$

Now with our Fano plane fully labeled, we can answer Question 1.29, as follows:

\index{incidence matrix}

\begin{answer}
The Fano plane can be described by the $7\times7$ incidence matrix recording the matches between points and lines, which is
$$m=\begin{pmatrix}
1&0&1&1&0&0&0\\
1&1&0&0&1&0&0\\
0&1&1&0&0&1&0\\
1&0&0&0&0&1&1\\
0&1&0&1&0&0&1\\
0&0&1&0&1&0&1\\
0&0&0&1&1&1&0
\end{pmatrix}$$
with points on the vertical, and lines on the horizontal, labeled as above. However, unlike for graphs, this matrix is no longer symmetric, or having $0$ on the diagonal. 
\end{answer}

So, question answered, rather defavorably for graph theory, and with this being said, shall we give up? Well, never underestimate the graphs. Indeed, these can strike back, and we have the following alternative answer to Question 1.29:

\index{bipartite graph}
\index{tripartite graph}
\index{multipartite graph}
\index{block diagonal}
\index{rectangular matrix}
\index{transpose matrix}

\begin{answer}
The Fano plane can be described by the bipartite graph having as vertices the points and lines, and edges at matches, whose $14\times14$ adjacency matrix is
$$d=\begin{pmatrix}0&m\\ m^t&0\end{pmatrix}$$
with $m$ being the usual $7\times7$ incidence matrix between points and lines, constructed before, and with $m^t$ being the transpose of $m$.
\end{answer}

Not bad all this. By the way, many other interesting things can be said, about the bipartite graphs, that is, about the graphs whose vertices can be divided into 2 classes, with no edges within the same class. For instance, it is clear that these are precisely the graphs having the property that, with a suitable labeling of the vertices, the adjacency matrix looks as follows, with $m$ being a certain rectangular 0-1 matrix:
$$d=\begin{pmatrix}0&m\\ m^t&0\end{pmatrix}$$

Along the same lines, we can talk as well about tripartite graphs, with the adjacency matrices here being as follows, with $m,n,p$ being certain rectangular 0-1 matrices:
$$d=\begin{pmatrix}
0&m&n\\ 
m^t&0&p\\
n^t&p^t&0
\end{pmatrix}$$

More generally, we can talk about $N$-partite graphs for any $N\in\mathbb N$, among others with the trivial remark that with $|X|=N$, our graph $X$ is obviously $N$-partite.

\bigskip

We will be back to such graphs, and more specifically to the bipartite ones, which are the most useful and important, among multi-partite graphs, later in this book.

\section*{1e. Exercises}

Here are some exercises, in relation with what we did in this chapter, quite often rather difficult, but hey, if looking for an easy book, many other choices are available:

\begin{exercise}
Learn some motivating physics, in relation with quantization.
\end{exercise}

\begin{exercise}
What Cayley graph should mean, based on what we know so far?
\end{exercise}

\begin{exercise}
Try doing some general spectral theory for the trees.
\end{exercise}

\begin{exercise}
Clarify our asymptotic estimate for the circle graph.
\end{exercise}

\begin{exercise}
Work out loop numerics for the circle, and the segment graphs.
\end{exercise}

\begin{exercise}
Fully clarify the loop count on $\mathbb Z$, including learning Stirling.
\end{exercise}

\begin{exercise}
Learn more about the $4$-color theorem, and its proof.
\end{exercise}

\begin{exercise}
Learn more about projective spaces, and about the Paley biplane too.
\end{exercise}

As bonus exercise, learn some programming. Any graph expert knows some. So, find some free math software, download, and start playing with it.

\chapter{General theory}

\section*{2a. Linear algebra}

You probably know from linear algebra that the important question regarding any matrix $d\in M_N(\mathbb R)$ is its diagonalization. To be more precise, the first question is that of computing the eigenvalues, which are usually complex numbers, $\lambda\in\mathbb C$:
$$dv=\lambda v$$

Then comes the computation of the eigenvectors, which are usually complex vectors too, $v\in\mathbb C^N$, and then the diagonalization problem, amounting in finding a new basis of $\mathbb C^N$, formed by eigenvectors. This basis can exist or not, and in case it exists, we reach to a formula as follows, with $\Lambda=diag(\lambda)$, and $P$ being the basis change matrix:
$$d=P\Lambda P^{-1}$$

In the case of the graphs, or rather of the adjacency matrices $d\in M_N(0,1)$ of the graphs, these notions are quite important and intuitive, as shown by:

\index{eigenvectors}
\index{eigenvalues}
\index{harmonic function}
\index{sum over neighbors}
\index{average over neighbors}

\begin{theorem}
The eigenvectors of $d\in M_N(0,1)$, with eigenvalue $\lambda$, can be identified with the functions $f$ satisfying the following condition:
$$\lambda f(i)=\sum_{i-j}f(j)$$ 
That is, the value of $f$ at each vertex must be the rescaled average, over the neighbors. 
\end{theorem}

\begin{proof}
We have indeed the following computation, valid for any vector $f$:
\begin{eqnarray*}
(df)_i
&=&\sum_jd_{ij}f_j\\
&=&\sum_{i-j}d_{ij}f_j+\sum_{i\not-j}d_{ij}f_j\\
&=&\sum_{i-j}1\cdot f_j+\sum_{i\not-j}0\cdot f_j\\
&=&\sum_{i-j}f_j
\end{eqnarray*}

Thus, we are led to the conclusion in the statement.
\end{proof}

The above result is quite interesting, and as an illustration, when assuming that our graph is $k$-regular, for the particular value $\lambda=k$, the eigenvalue condition reads:
$$f(i)=\frac{1}{k}\sum_{i-j}f(j)$$ 

Thus, we can see here a relation with harmonic functions.  There are many things that can be said here, and we will be back to this later, when talking Laplace operators.

\bigskip

But let us pause now our study of graphs, and go back to linear algebra. Taking the scalars to be complex numbers, which is something very standard, we have the following general result, that you surely know about, from mathematics or physics classes:

\begin{theorem}
Given a matrix $d\in M_N(\mathbb C)$, consider its characteristic polynomial
$$P(x)=\det(d-x1_N)$$ 
then factorize this polynomial, by computing its complex roots, with multiplicities,
$$P(x)=(-1)^N(x-\lambda_1)^{m_1}\ldots(x-\lambda_k)^{m_k}$$
and finally compute the corresponding eigenspaces, for each eigenvalue found:
$$E_i=\left\{v\in\mathbb C^N\Big|dv=\lambda_iv\right\}$$
The dimensions of these eigenspaces satisfy then the following inequalities,
$$\dim(E_i)\leq m_i$$
and $d$ is diagonalizable precisely when we have equality for any $i$.
\end{theorem}

\begin{proof}
There are many things going on here, the idea being as follows:

\medskip

(1) Given $d\in M_N(\mathbb C)$, for any eigenvalue $\lambda\in\mathbb C$ we can define the corresponding eigenspace as being the vector space formed by the corresponding eigenvectors:
$$E_\lambda=\left\{v\in\mathbb C^N\Big|dv=\lambda v\right\}$$

These spaces $E_\lambda$ are easily seen to be in direct sum position, in the sense that given vectors $v_1\in E_{\lambda_1},\ldots,v_k\in E_{\lambda_k}$ coming from different eigenvalues $\lambda_1,\ldots,\lambda_k$, we have:
$$\sum_ic_iv_i=0\implies c_i=0$$

In particular we have the following estimate, with sum over all the eigenvalues, and our matrix is diagonalizable precisely when we have equality:
$$\sum_\lambda\dim(E_\lambda)\leq N$$

(2) Next, consider the characteristic polynomial of our matrix, given by:
$$P(x)=\det(d-x1_N)$$

Our claim is that the eigenvalues of $d$ are then the roots of $P$. Indeed, this follows from the following computation, using the elementary fact that a linear map is bijective precisely when the determinant of the associated matrix is nonzero:
\begin{eqnarray*}
\exists v,dv=\lambda v
&\iff&\exists v,(d-\lambda 1_N)v=0\\
&\iff&\det(d-\lambda 1_N)=0
\end{eqnarray*}

(3) Our next claim, which will lead to the result in the statement, is that we have an inequality as follows, where $m_\lambda$ is the multiplicity of $\lambda$, viewed as root of $P$:
$$\dim(E_\lambda)\leq m_\lambda$$

Indeed, for any eigenvalue $\lambda$, consider the dimension $n_\lambda=\dim(E_\lambda)$ of the corresponding eigenspace. By changing the basis of $\mathbb C^N$, as for $E_\lambda$ to be spanned by the first $n_\lambda$ basis elements, our matrix becomes as follows, with $e$ being a certain smaller matrix:
$$d\sim\begin{pmatrix}\lambda 1_{n_\lambda}&0\\0&e\end{pmatrix}$$

We conclude that the characteristic polynomial of $d$ is of the following form:
$$P_d
=P_{\lambda 1_{n_\lambda}}P_e
=(\lambda-x)^{n_\lambda}P_e$$

Thus the multiplicity of $\lambda$, as root of $P$, satisfies $m_\lambda\geq n_\lambda$, which proves our claim.

\medskip

(4) Getting now to what we wanted to prove, as a first observation, what is said in the theorem is well formulated, thanks to what we found above. Next, by summing the inequalities $\dim(E_\lambda)\leq m_\lambda$ found in (3), we obtain an inequality as follows:
$$\sum_\lambda\dim(E_\lambda)\leq\sum_\lambda m_\lambda\leq N$$

On the other hand, we know from (1) that our matrix is diagonalizable precisely when we have global equality. Thus, we are led to the conclusion in the statement.
\end{proof}

In practice, diagonalizing a matrix remains something quite complicated. Let us record as well a useful, algorithmic version of the above result, as follows:

\begin{theorem}
The square matrices $d\in M_N(\mathbb C)$ can be diagonalized as follows:
\begin{enumerate}
\item Compute the characteristic polynomial.

\item Factorize the characteristic polynomial.

\item Compute the eigenvectors, for each eigenvalue found.

\item If there are no $N$ eigenvectors, $d$ is not diagonalizable.

\item Otherwise, $d$ is diagonalizable, $d=P\Lambda P^{-1}$.
\end{enumerate}
\end{theorem}

\begin{proof}
This is an informal reformulation of Theorem 2.2, with (4) referring to the total number of linearly independent eigenvectors found in (3), and with $d=P\Lambda P^{-1}$ in (5) being the usual diagonalization formula, with $P,\Lambda$ being as before.
\end{proof}

Getting now towards our graph problematics, here is a key result:

\index{self-adjoint matrix}
\index{adjoint matrix}
\index{diagonalization}
\index{spectral theorem}

\begin{theorem}
Any matrix $d\in M_N(\mathbb C)$ which is self-adjoint, $d=d^*$, is diagonalizable, with the diagonalization being of the following type,
$$d=U\Lambda U^*$$
with $U\in U_N$, and with $\Lambda\in M_N(\mathbb R)$ diagonal. The converse holds too.
\end{theorem}

\begin{proof}
As a first remark, the converse trivially holds, because if we take a matrix of the form $d=U\Lambda U^*$, with $U$ unitary and $\Lambda$ diagonal and real, then we have:
\begin{eqnarray*}
d^*
&=&(U\Lambda U^*)^*\\
&=&U\Lambda^*U^*\\
&=&U\Lambda U^*\\
&=&d
\end{eqnarray*}

In the other sense now, assume that $d$ is self-adjoint, $d=d^*$.  Our first claim is that the eigenvalues are real. Indeed, assuming $dv=\lambda v$, we have:
\begin{eqnarray*}
\lambda<v,v>
&=&<\lambda v,v>\\
&=&<dv,v>\\
&=&<v,dv>\\
&=&<v,\lambda v>\\
&=&\bar{\lambda}<v,v>
\end{eqnarray*}

Thus we obtain $\lambda\in\mathbb R$, as claimed. Our next claim now is that the eigenspaces corresponding to different eigenvalues are pairwise orthogonal. Assume indeed that:
$$dv=\lambda v\quad,\quad 
dw=\mu w$$

We have then the following computation, using $\lambda,\mu\in\mathbb R$:
\begin{eqnarray*}
\lambda<v,w>
&=&<\lambda v,w>\\
&=&<dv,w>\\
&=&<v,dw>\\
&=&<v,\mu w>\\
&=&\mu<v,w>
\end{eqnarray*}

Thus $\lambda\neq\mu$ implies $v\perp w$, as claimed. In order now to finish, it remains to prove that the eigenspaces span $\mathbb C^N$. For this purpose, we will use a recurrence method. Let us pick an eigenvector, $dv=\lambda v$. Assuming $v\perp w$, we have:
\begin{eqnarray*}
<dw,v>
&=&<w,dv>\\
&=&<w,\lambda v>\\
&=&\lambda<w,v>\\
&=&0
\end{eqnarray*}

Thus, if $v$ is an eigenvector, then the vector space $v^\perp$ is invariant under $d$. In order to do the recurrence, it still remains to prove that the restriction of $d$ to the vector space $v^\perp$ is self-adjoint. But this comes from a general property of the self-adjoint matrices, that we will explain now. Our claim is that for any matrix $d\in M_N(\mathbb C)$, we have:
$$d=d^*\iff<dv,v>\in\mathbb R,\,\forall v$$

Indeed, this follows from the following computation:
\begin{eqnarray*}
<dv,v>-\overline{<dv,v>}
&=&<dv,v>-<v,dv>\\
&=&<dv,v>-<d^*v,v>\\
&=&<(d-d^*)v,v>\\
&=&\sum_{ij}(d-d^*)_{ij}\bar{v}_iv_j
\end{eqnarray*}

But this shows that the restriction of $d$ to any invariant subspace, and in particular to $v^\perp$, is self-adjoint. Thus, we can proceed by recurrence, and we obtain the result.
\end{proof}

In what concerns us, in relation with our graph problems, we will rather need the real version of the above result, which is also something well-known, as follows:

\index{spectral theorem}
\index{symmetric matrix}
\index{orthogonal matrix}

\begin{theorem}
Any matrix $d\in M_N(\mathbb R)$ which is symmetric, $d=d^t$, is diagonalizable, with the diagonalization being of the following type,
$$d=U\Lambda U^t$$
with $U\in O_N$, and with $\Lambda \in M_N(\mathbb R)$ diagonal. The converse holds too.
\end{theorem}

\begin{proof}
As before, the converse trivially holds, because if we take a matrix of the form $d=U\Lambda U^t$, with $U$ orthogonal, and $\Lambda$ diagonal and real, then we have:
\begin{eqnarray*}
d^t
&=&(U\Lambda U^t)^t\\
&=&U\Lambda^tU^t\\
&=&U\Lambda U^t\\
&=&d
\end{eqnarray*}

In the other sense now, this follows from Theorem 2.4, and its proof. Indeed, we know from there that the eigenvalues are real, and in what concerns the passage matrix, the arguments there carry over to the real case, and show that this matrix is real too.
\end{proof}

With the above results in hand, time now to get back to graphs. We have here the following particular case of Theorem 2.5, with the important drawback however that in what concerns the ``converse holds too'' part, that is unfortunately gone:

\index{trace of matrix}

\begin{theorem}
The adjacency matrix $d\in M_N(0,1)$ of any graph is diagonalizable, with the diagonalization being of the following type,
$$d=U\Lambda U^t$$
with $U\in O_N$, and with $\Lambda\in M_N(\mathbb R)$ diagonal. Moreover, we have $Tr(\Lambda)=0$.
\end{theorem}

\begin{proof}
Here the first assertion follows from Theorem 2.5, because $d$ is by definition real and symmetric. As for the last assertion, this deserves some explanations:

\medskip

(1) Generally speaking, in analogy with the last assertions in Theorem 2.4 and Theorem 2.5, which are something extremely useful, we would like to know under which assumptions on a rotation $U\in O_N$, and on a diagonal matrix $\Lambda\in M_N(\mathbb R)$, the real symmetric matrix $d=U\Lambda U^t$ has 0-1 entries, and 0 on the diagonal.

\medskip

(2) Unfortunately, both these questions are obviously difficult, there is no simple answer to them, and things are like that. So, gone the possibility of a converse. However, as a small consolation, we can make the remark that, with $d=U\Lambda U^t$, we have:
$$Tr(d)=Tr(U\Lambda U^t)=Tr(\Lambda)$$

Thus we have at least $Tr(\Lambda)=0$, as a necessary condition on $(U,\Lambda)$, as stated.
\end{proof}

In view of the above difficulties with the bijectivity, it is perhaps wise to formulate as well the graph particular case of Theorem 2.4. The statement here is as follows:

\begin{theorem}
The adjacency matrix $d\in M_N(0,1)$ of any graph is diagonalizable, with the diagonalization being of the following type,
$$d=U\Lambda U^*$$
with $U\in U_N$, and with $\Lambda\in M_N(\mathbb R)$ diagonal. Moreover, we have $Tr(\Lambda)=0$.
\end{theorem}

\begin{proof}
This follows from Theorem 2.4, via the various remarks from the proof of Theorem 2.5 and Theorem 2.6. But the simplest is to say that the statement itself is just a copy of Theorem 2.6, with $U\in O_N$ replaced by the more general $U\in U_N$.
\end{proof}

All the above is useful, and we will use these results on a regular basis, in what follows. However, before getting into more concrete things, let us formulate:

\begin{problem}
Find a geometric proof of Theorem 2.6, or of Theorem 2.7, based on the interpretation of eigenvalues and eigenvectors from Theorem 2.1.
\end{problem}

This question looks quite reasonable, at a first glance, after all what we have in Theorem 2.1 is all nice and gentle material, so do we really need all the above complicated linear algebra machinery in order to deal with all this. However, at a second look, meaning after studying some examples, the problem suddenly looks very complicated. So, homework for you, in case I forget to assign this, to come back to this problem, later.

\section*{2b. The simplex}

As an illustration for the above, let us diagonalize the adjacency matrix of the simplest graph that we know, namely the $N$-simplex. Let us start with:

\index{simplex}
\index{flat matrix}
\index{all-one matrix}

\begin{proposition}
The adjacency matrix of the $N$-simplex, having $0$ on the diagonal and $1$ elsewhere, is in matrix form:
$$d=\begin{pmatrix}
0&1&\ldots&1&1\\
1&0&\ldots&1&1\\
\vdots&\vdots&&\vdots&\vdots\\
1&1&\ldots&0&1\\
1&1&\ldots&1&0
\end{pmatrix}$$
We have the following formula for it, with $\mathbb I$ standing for the all-$1$ matrix, and with $1$ standing for the identity matrix, both of size $N$:
$$d=\mathbb I-1$$
Equivalently, $d=NP-1$, with $P$ being the projection on the all-$1$ vector $\xi\in\mathbb R^N$.
\end{proposition}

\begin{proof}
Here the first assertion is clear from definitions, and the second assertion is clear too. As for the last assertion, observe first that with $P=\mathbb I/N$ we have:
$$P^2=\left(\frac{\mathbb I}{N}\right)^2=\frac{\mathbb I^2}{N^2}=\frac{N\mathbb I}{N^2}=\frac{\mathbb I}{N}=P$$

Thus $P$ is a projection, and since we obviously have $P=P^t$, this matrix is an orthogonal projection. In order to find now the image, observe that for any vector $v\in\mathbb R^N$ we have the following formula, with $a\in\mathbb R$ being the average of the entries of $v$:
$$Pv=\begin{pmatrix}
a\\
\vdots\\
a\end{pmatrix}$$

We conclude that the image of $P$ is the vector space $\mathbb R\xi$, with $\xi\in\mathbb R^N$ being the all-1 vector, and so that $P$ is the orthogonal projection on $\xi$, as claimed.
\end{proof}

The above is very nice, in particular with $d=NP-1$ basically diagonalizing $d$ for us. However, thinking a bit, when it comes to explicitly diagonalize $d$, or, equivalently, $P$ or $\mathbb I$, things are quite tricky, and we run into the following strange problem:

\begin{problem}
In order to diagonalize $d$, we need solutions for
$$x_1+\ldots+x_N=0$$
and there are no standard such solutions, over the reals.
\end{problem}

So, this is the problem that we face, which might look a bit futile at a first glance, and in order for you to take me seriously here, let us work out some particular cases. At $N=2$ things are quickly solved, with the diagonalization being as follows, and with the passage matrix being easy to construct, I agree with you on this: 
$$\begin{pmatrix}0&1\\1&0\end{pmatrix}
\sim \begin{pmatrix}1&0\\0&-1\end{pmatrix}$$

However, things become suddenly complicated at $N=3$, where I challenge you to find the passage matrix for the following diagonalization:
$$\begin{pmatrix}0&1&1\\1&0&1\\1&1&0\end{pmatrix}
\sim \begin{pmatrix}2&0&0\\0&-1&0\\0&0&-1\end{pmatrix}$$

In the general case, $N\in\mathbb N$, the problem does not get any simpler, again with the challenge for you to find the passage matrix for the following diagonalization:
$$\begin{pmatrix}
0&1&\ldots&1&1\\
1&0&\ldots&1&1\\
\vdots&\vdots&&\vdots&\vdots\\
1&1&\ldots&0&1\\
1&1&\ldots&1&0
\end{pmatrix}\sim
\begin{pmatrix}
N-1&&&&0\\
&-1\\
&&\ddots&\\
&&&-1\\
0&&&&-1
\end{pmatrix}$$

In short, you got my point, Problem 2.10 is something real. Fortunately the complex numbers come to the rescue, via the following standard and beautiful result:

\index{roots of unity}
\index{barycenter}

\begin{theorem}
The roots of unity, $\{w^k\}$ with $w=e^{2\pi i/N}$, have the property
$$\frac{1}{N}\sum_{k=0}^{N-1}(w^k)^s=\delta_{N|s}$$
for any exponent $s\in\mathbb N$, where on the right we have a Kronecker symbol.
\end{theorem}

\begin{proof}
There are several possible proofs for this, as follows:

\medskip

(1) Nice proof. The numbers to be summed, when written more conveniently as $(w^s)^k$ with $k=0,\ldots,N-1$, form a certain regular polygon in the plane $P_s$. Thus, if we denote by $C_s$ the barycenter of this polygon, we have the following formula:
$$\frac{1}{N}\sum_{k=0}^{N-1}w^{ks}=C_s$$

Now observe that in the case $N\slash\hskip-1.6mm|\,s$ our polygon $P_s$ is non-degenerate, circling along the unit circle, and having center $C_s=0$. As for the case $N|s$, here the polygon is degenerate, lying at 1, and having center $C_s=1$. Thus, we have the following formula, as claimed:
$$C_s=\delta_{N|s}$$

(2) Ugly proof. Normally we won't bother with ugly proofs in this book, but these being mathematics too, at least in theory, here is the computation, for $N\slash\hskip-1.6mm|\,s$:
\begin{eqnarray*}
\sum_{k=0}^{N-1}(w^k)^s
&=&\sum_{k=0}^{N-1}(w^s)^k\\
&=&\frac{1-w^{sN}}{1-w^s}\\
&=&\frac{1-1^s}{1-w^s}\\
&=&0
\end{eqnarray*}

Thus, we are led again to the formula in the statement.
\end{proof}

Summarizing, we have the solution to our problem. In order now to finalize, let us start with the following definition, inspired by what happens in Theorem 2.11:

\index{Fourier matrix}

\begin{definition}
The Fourier matrix $F_N$ is the following matrix, with $w=e^{2\pi i/N}$:
$$F_N=
\begin{pmatrix}
1&1&1&\ldots&1\\
1&w&w^2&\ldots&w^{N-1}\\
1&w^2&w^4&\ldots&w^{2(N-1)}\\
\vdots&\vdots&\vdots&&\vdots\\
1&w^{N-1}&w^{2(N-1)}&\ldots&w^{(N-1)^2}
\end{pmatrix}$$
That is, $F_N=(w^{ij})_{ij}$, with indices $i,j\in\{0,1,\ldots,N-1\}$, taken modulo $N$.
\end{definition}

Here the name comes from the fact that $F_N$ is the matrix of the discrete Fourier transform, that over the cyclic group $\mathbb Z_N$, and more on this later, when talking Fourier analysis. As a first example now, at $N=2$ the root of unity is $w=-1$, and with indices as above, namely $i,j\in\{0,1\}$, taken modulo 2, our Fourier matrix is as follows:
$$F_2=\begin{pmatrix}1&1\\1&-1\end{pmatrix}$$

At $N=3$ now, the root of unity is $w=e^{2\pi i/3}$, and the Fourier matrix is:
$$F_3=\begin{pmatrix}1&1&1\\ 1&w&w^2\\ 1&w^2&w\end{pmatrix}$$

At $N=4$ now, the root of unit is $w=i$, and the Fourier matrix is:
$$F_4=\begin{pmatrix}
1&1&1&1\\
1&i&-1&-i\\
1&-1&1&-1\\
1&-i&-1&i
\end{pmatrix}$$

Finally, at $N=5$ the root of unity is $w=e^{2\pi i/5}$, and the Fourier matrix is:
$$F_5=\begin{pmatrix}
1&1&1&1&1\\
1&w&w^2&w^3&w^4\\
1&w^2&w^4&w&w^3\\
1&w^3&w&w^4&w^2\\
1&w^4&w^3&w^2&w
\end{pmatrix}$$

You get the point. Getting back now to the diagonalization problem for the flat matrix $\mathbb I$, this can be solved by using the Fourier matrix $F_N$, in the following way:

\index{flat matrix}
\index{Fourier matrix}

\begin{proposition}
The flat matrix $\mathbb I$ diagonalizes as follows,
$$\begin{pmatrix}
1&\ldots&\ldots&1\\
\vdots&&&\vdots\\
\vdots&&&\vdots\\
1&\ldots&\ldots&1\end{pmatrix}
=\frac{1}{N}\,F_N
\begin{pmatrix}
N&&&&0\\
&0\\
&&\ddots\\
&&&0\\
0&&&&0\end{pmatrix}F_N^*$$
with $F_N=(w^{ij})_{ij}$ being the Fourier matrix.
\end{proposition}

\begin{proof}
According to the last assertion in Proposition 2.9, we are left with finding the 0-eigenvectors of $\mathbb I$, which amounts in solving the following equation:
$$x_0+\ldots+x_{N-1}=0$$

For this purpose, we can use the root of unity $w=e^{2\pi i/N}$, and more specifically, the following standard formula, coming from Theorem 2.11:
$$\sum_{i=0}^{N-1}w^{ij}=N\delta_{j0}$$

This formula shows that for $j=1,\ldots,N-1$, the vector $v_j=(w^{ij})_i$ is a 0-eigenvector. Moreover, these vectors are pairwise orthogonal, because we have:
$$<v_j,v_k>
=\sum_iw^{ij-ik}
=N\delta_{jk}$$

Thus, we have our basis $\{v_1,\ldots,v_{N-1}\}$ of 0-eigenvectors, and since the $N$-eigenvector is $\xi=v_0$, the passage matrix $P$ that we are looking is given by:
$$P=\begin{bmatrix}v_0&v_1&\ldots&v_{N-1}\end{bmatrix}$$

But this is precisely the Fourier matrix, $P=F_N$. In order to finish now, observe that the above computation of $<v_i,v_j>$ shows that $F_N/\sqrt{N}$ is unitary, and so:
$$F_N^{-1}=\frac{1}{N}\,F_N^*$$

Thus, we are led to the diagonalization formula in the statement.
\end{proof}

By substracting now $-1$ from everything, we can formulate a final result, as follows:

\index{flat matrix}
\index{Fourier matrix}

\begin{theorem}
The adjacency matrix of the simplex diagonalizes as follows,
$$\begin{pmatrix}
0&1&\ldots&1&1\\
1&0&\ldots&1&1\\
\vdots&\vdots&&\vdots&\vdots\\
1&1&\ldots&0&1\\
1&1&\ldots&1&0
\end{pmatrix}
=\frac{1}{N}\,F_N
\begin{pmatrix}
N-1&&&&0\\
&-1\\
&&\ddots&\\
&&&-1\\
0&&&&-1
\end{pmatrix}F_N^*$$
with $F_N=(w^{ij})_{ij}$ being the Fourier matrix.
\end{theorem}

\begin{proof}
This follows as said above, from what we have in Proposition 2.13, by substracting $-1$ from everything. Alternatively, if you prefer a more direct proof, this follows from the various computations in the proof of Proposition 2.13.
\end{proof}

The above result was something quite tricky, and we will come back regularly to such things, in what follows. For the moment, let us formulate an interesting conclusion:

\begin{conclusion}
Theorem 2.7, telling us that $d$ diagonalizes over $\mathbb C$, is better than the stronger Theorem 2.6, telling us that $d$ diagonalizes over $\mathbb R$.
\end{conclusion}

And isn't this surprising. But after some thinking, after all no surprise, because the graphs, as we defined them in the beginning of chapter 1, are scalarless objects. So, when needing a field for studying them, we should just go with the mighty $F=\mathbb C$. 

\bigskip

By the way, regarding complex numbers, time to recommend some reading. Mathematically the book of Rudin \cite{ru2} is a must-read, and a pleasure to read. However, if you want to be really in love with complex numbers, and with this being an enormous asset, no matter what mathematics you want to do, nothing beats some physics reading.

\bigskip

The standard place for learning physics is the course of Feynman \cite{fe1}, \cite{fe2}, \cite{fe3}. If you already know a bit of physics, you can go as well with the lovely books of Griffiths \cite{gr1}, \cite{gr2}, \cite{gr3}. And if you know a bit more, good books are those of Weinberg \cite{we1}, \cite{we2}, \cite{we3}. In the hope that this helps, and I will not tell you of course what these books contain. Expect however lots of complex numbers, all beautiful, and used majestically.

\section*{2c. Further graphs} 

Motivated by the above, let us try now to diagonalize the adjacency matrices of other graphs, more complicated than the simplex. But, what graphs to start with?

\bigskip

We already faced this kind of question in chapter 1, when discussing random walks and other basic questions, and we will proceed in a similar way here, by using our intuition. Passed the simplices, the ``simplest'' graphs are most likely the unions of simplices:
$$
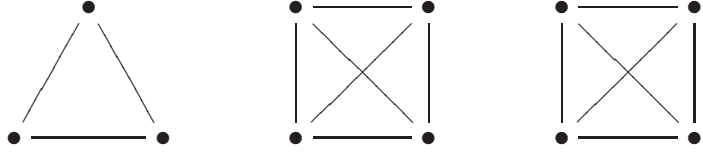
\xymatrix@R=15pt@C=15pt{
&\bullet\ar@{-}[ddl]\ar@{-}[ddr]&&&\bullet\ar@{-}[dd]\ar@{-}[ddrr]\ar@{-}[rr]&&\bullet\ar@{-}[dd]\ar@{-}[ddll]&&\bullet\ar@{-}[dd]\ar@{-}[ddrr]\ar@{-}[rr]&&\bullet\ar@{-}[dd]\ar@{-}[ddll]\\
\\
\bullet\ar@{-}[rr]&&\bullet&&\bullet\ar@{-}[rr]&&\bullet&&\bullet\ar@{-}[rr]&&\bullet
}$$

But the diagonalization question for such graphs is quickly settled, by using our results for the simplex from the previous section, and the following general fact:

\begin{theorem}
Given a disjoint union of connected graphs $X=X_1\cup\ldots\cup X_k$, the corresponding adjacency matrix is block-diagonal,
$$d=\begin{pmatrix}
d_1\\
&\ddots\\
&&d_k
\end{pmatrix}$$
and diagonalizes in the following way,
$$d=\begin{pmatrix}
U_1\\
&\ddots\\
&&U_k
\end{pmatrix}
\begin{pmatrix}
\Lambda_1\\
&\ddots\\
&&\Lambda_k
\end{pmatrix}
\begin{pmatrix}
U_1\\
&\ddots\\
&&U_k
\end{pmatrix}^*$$
with $d_i=U_i\Lambda_iU_i^*$ being diagonalizations of the adjacency matrices of the components $X_i$.
\end{theorem}

\begin{proof}
This is indeed something trivial and self-explanatory, coming from definitions, the idea being that when labeling the vertices of $X$ by starting with those of $X_1$, then with those of $X_2$, and so on up to those of $X_k$, the adjacency matrix $d$ takes the block-diagonal form in the statement. Thus, we are led to the above conclusions.
\end{proof}

Summarizing, we know how to diagonalize the adjacency matrix of $X=X_1\cup\ldots\cup X_k$, provided that we know how to diagonalize the adjacency matrix of each $X_i$. As an illustration here, for the graph pictured before Theorem 2.16, the diagonalization is as follows, with $F_N$ being as usual the Fourier matrix, and with $J_N=diag(N-1,-1,\ldots,-1)$ being the diagonal $N\times N$ matrix found in Theorem 2.14: 
$$d=\frac{1}{48}\begin{pmatrix}
F_3\\
&F_4\\
&&F_4
\end{pmatrix}
\begin{pmatrix}
J_3\\
&J_4\\
&&J_4
\end{pmatrix}
\begin{pmatrix}
F_3\\
&F_4\\
&&F_4
\end{pmatrix}^*$$

Back to the general case, as in Theorem 2.16, of particular interest is the case where all the components $X_i$ appear as copies of the same graph $Y$. And with the picture here, when $Y$ is a simplex, being something which looks quite nice, as follows:
$$
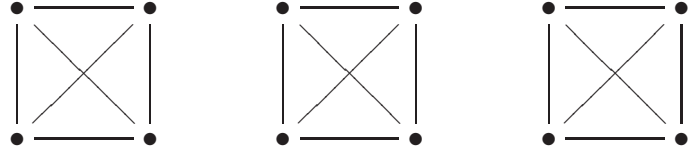
\xymatrix@R=15pt@C=15pt{
\bullet\ar@{-}[dd]\ar@{-}[ddrr]\ar@{-}[rr]&&\bullet\ar@{-}[dd]\ar@{-}[ddll]&&\bullet\ar@{-}[dd]\ar@{-}[ddrr]\ar@{-}[rr]&&\bullet\ar@{-}[dd]\ar@{-}[ddll]&&\bullet\ar@{-}[dd]\ar@{-}[ddrr]\ar@{-}[rr]&&\bullet\ar@{-}[dd]\ar@{-}[ddll]\\
\\
\bullet\ar@{-}[rr]&&\bullet&&\bullet\ar@{-}[rr]&&\bullet&&\bullet\ar@{-}[rr]&&\bullet
}$$

In this case we can further build on Theorem 2.16, and say a bit more about the diagonalization of the adjacency matrix, the final result being as follows:

\begin{theorem}
Given a disjoint union of type $X=kY$, with $Y$ being connected, the corresponding adjacency matrix can be written in the following way,
$$d=d_Y\otimes 1_k$$
and diagonalizes in the following way,
$$d=(U\otimes 1_k)(\Lambda\otimes 1_k)(U\otimes 1_k)^*$$
with $d_Y=U\Lambda U^*$ being the diagonalization of the adjacency matrix of $Y$.
\end{theorem}

\begin{proof}
This follows from Theorem 2.16, with the only issue being at the level of the labeling of vertices, and of the tensor product notations. In order to explain this, let us take as example the graph pictured before the statement, and label it as follows:
$$
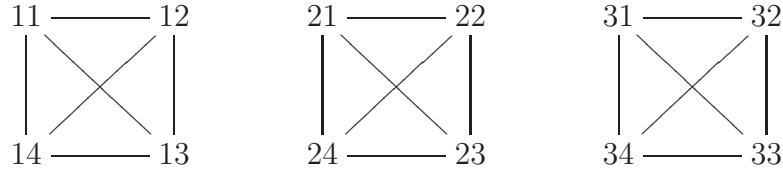
\xymatrix@R=15pt@C=15pt{
11\ar@{-}[dd]\ar@{-}[ddrr]\ar@{-}[rr]&&12\ar@{-}[dd]\ar@{-}[ddll]&&21\ar@{-}[dd]\ar@{-}[ddrr]\ar@{-}[rr]&&22\ar@{-}[dd]\ar@{-}[ddll]&&31\ar@{-}[dd]\ar@{-}[ddrr]\ar@{-}[rr]&&32\ar@{-}[dd]\ar@{-}[ddll]\\
\\
14\ar@{-}[rr]&&13&&24\ar@{-}[rr]&&23&&34\ar@{-}[rr]&&33
}$$

Observe that this is in agreement with what we did in the proof of Theorem 2.16, with the vertices of the first square coming first, then with the vertices of the second square coming second, and with the vertices of the third square coming third, according to:
$$11<12<13<14<21<22<23<24<31<32<33<34$$

Thus, we can use the conclusions of Theorem 2.16, as such. But since our adjacency matrix has now double indices, as matrix indices, we can switch to tensor product notations, according to the following standard rule, from tensor product calculus:
$$(M\otimes N)_{ia,jb}=M_{ij}N_{ab}$$

Thus, we are led to the conclusions in the statement.
\end{proof}

As before with Theorem 2.16, nothing better at this point than an illustration for all this. For the graph pictured before the statement, and in the proof too, we obtain:
$$d=\frac{1}{4}(F_4\otimes 1_k)(J_4\otimes 1_k)(F_4\otimes 1_k)^*$$

Moving forward, now that we understood the disjoint union operation, time to study the other basic operation for graphs, which is the complementation operation. However, things here are more tricky, with the general result on the subject being as follows:

\index{complement}

\begin{proposition}
The adjacency matrix of a graph $X$ and of its complement $X^c$ are related by the following formula, with $\mathbb I$ being the all-one matrix,
$$d_X+d_{X^c}=\mathbb I-1$$
so if we have a diagonalization formula of type $d_X=U\Lambda U^*$, with $U$ diagonalizing as well the all-one matrix $\mathbb I$, we obtain in this way a diagonalization for $d_X$.
\end{proposition}

\begin{proof}
This is something trivial, and which of course does not look very good, with our point coming from the following implication, which is immediate:
$$d_X=U\Lambda U^*\ ,\ \mathbb I=U\Sigma U^*\ \implies\ d_X=U(\Sigma-\Lambda-1)U^*$$

Thus, we are led to the conclusion in the statement.
\end{proof}

In practice now, the next step would be to look at the complements of various particular graphs, such as the disjoint unions $X=X_1\cup\ldots\cup X_k$ from Theorem 2.16. But this does not bring any simplification, as you can easily verify yourself. Moreover, when further particularizing, and looking at the complements of the graphs $X=kY$ from Theorem 2.17, there is no simplification either. However, when particularizing even more, and taking $Y=K_n$, we are led to something nice, which is good to know, namely:

\index{disjoint union}

\begin{theorem}
Given a disjoint union of simplices $X=kK_n$, the adjacency matrix of the complement $X^c$ can be written as
$$d=\mathbb I_n\otimes(\mathbb I_k-1_k)$$
and so diagonalizes in the following way,
$$d=\frac{1}{nk}(F_n\otimes F_k)\left(L_n\otimes(L_k-1_k)\right)(F_n\otimes F_k)^*$$
with $F_n$ being the Fourier matrix, and $L_n=diag(n,0,\ldots,0)$, as a $n\times n$ matrix.
\end{theorem}

\begin{proof}
By using the various labeling and tensor product conventions from Theorem 2.17 and its proof, the adjacency matrix of $(kK_n)^c$ is given by:
\begin{eqnarray*}
d
&=&\begin{pmatrix}
0&\mathbb I_n&\ldots&\mathbb I_n&\mathbb I_n\\
\mathbb I_n&0&\ldots&\mathbb I_n&\mathbb I_n\\
&&\ddots\\
\mathbb I_n&\mathbb I_n&\ldots&0&\mathbb I_n\\
\mathbb I_n&\mathbb I_n&\ldots&\mathbb I_n&0
\end{pmatrix}\\
&=&\mathbb I_n\otimes\begin{pmatrix}
0&1&\ldots&1&1\\
1&0&\ldots&1&1\\
&&\ddots\\
1&1&\ldots&0&1\\
1&1&\ldots&1&0
\end{pmatrix}\\
&=&\mathbb I_n\otimes(\mathbb I_k-1_k)
\end{eqnarray*}

Thus, we are led to the various conclusions in the statement.
\end{proof}

As before with other results of the same type, let us work out an illustration. Consider the following graph, that we already met before, in the context of Theorem 2.17:
$$
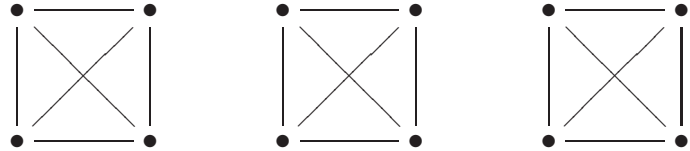
\xymatrix@R=15pt@C=15pt{
\bullet\ar@{-}[dd]\ar@{-}[ddrr]\ar@{-}[rr]&&\bullet\ar@{-}[dd]\ar@{-}[ddll]&&\bullet\ar@{-}[dd]\ar@{-}[ddrr]\ar@{-}[rr]&&\bullet\ar@{-}[dd]\ar@{-}[ddll]&&\bullet\ar@{-}[dd]\ar@{-}[ddrr]\ar@{-}[rr]&&\bullet\ar@{-}[dd]\ar@{-}[ddll]\\
\\
\bullet\ar@{-}[rr]&&\bullet&&\bullet\ar@{-}[rr]&&\bullet&&\bullet\ar@{-}[rr]&&\bullet
}$$

The adjacency matrix of the complementary graph diagonalizes then as follows:
$$d=\frac{1}{12}(F_4\otimes F_3)\left[
\begin{pmatrix}
4\\
&0\\
&&0\\
&&&0
\end{pmatrix}\otimes
\begin{pmatrix}
3\\
&-1\\
&&-1
\end{pmatrix}\right]
(F_4\otimes F_3)^*$$

And with this being actually not the end of the story, for this particular graph, because it is possible to prove that we have $F_{12}=F_4\otimes F_3$, with this coming from $F_{mn}=F_m\otimes F_n$, for $m,n$ prime to each other. But more on such phenomena later in this book.

\bigskip

Moving forward, now that we have some good understanding of the various operations for the graphs, time to get into the real thing, namely study of the graphs which are nice, simple and ``indecomposable''. Although we don't have yet a definition, for what indecomposable should mean, as a main example here, we certainly have the circle:
$$
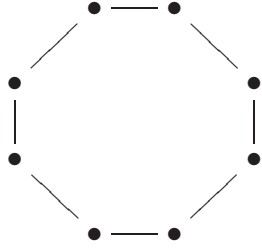
\xymatrix@R=16pt@C=17pt{
&\bullet\ar@{-}[r]\ar@{-}[dl]&\bullet\ar@{-}[dr]\\
\bullet\ar@{-}[d]&&&\bullet\ar@{-}[d]\\
\bullet\ar@{-}[dr]&&&\bullet\ar@{-}[dl]\\
&\bullet\ar@{-}[r]&\bullet}$$

So, in what follows we will study this circle graph, and its complement too. Let us first compute the eingenvalues and eigenvectors. This is easily done, as follows:

\index{circle graph}

\begin{proposition}
For the adjacency matrix $d$ of the circle graph, we have
$$d\begin{pmatrix}
1\\
w^k\\
w^{2k}\\
\vdots\\
w^{kN-k}
\end{pmatrix}
=(w^k+w^{-k})
\begin{pmatrix}
1\\
w^k\\
w^{2k}\\
\vdots\\
w^{kN-k}
\end{pmatrix}$$
with $w=e^{2\pi i/N}$, for any $k$, and this diagonalizes $d$.
\end{proposition}

\begin{proof}
In what regards the first assertion, this follows from:
\begin{eqnarray*}
d\begin{pmatrix}
1\\
w^k\\
w^{2k}\\
\vdots\\
w^{kN-k}
\end{pmatrix}
&=&
\begin{pmatrix}
0&1&0&\ldots&0&0&1\\
1&0&1&\ldots&0&0&0\\
&&\ddots&\ddots&\ddots\\
0&0&0&\ldots&1&0&1\\
1&0&0&\ldots&0&1&0
\end{pmatrix}
\begin{pmatrix}
1\\
w^k\\
w^{2k}\\
\vdots\\
w^{kN-k}
\end{pmatrix}\\
&=&\begin{pmatrix}
w^k+w^{kN-k}\\
1+w^{2k}\\
w^k+w^{3k}\\
\vdots\\
1+w^{kN-2k}
\end{pmatrix}\\
&=&(w^k+w^{-k})
\begin{pmatrix}
1\\
w^k\\
w^{2k}\\
\vdots\\
w^{kN-k}
\end{pmatrix}
\end{eqnarray*}

As for the second assertion, this comes from this, because with $k=0,1,\ldots,N-1$ the eigenvectors that we found, which are the rows of $F_N$, are linearly independent.
\end{proof}

We can now formulate a final result about the circle graph, as follows:

\begin{theorem}
The adjacency matrix of the circle graph diagonalizes as
$$d=\frac{1}{N}F_N
\begin{pmatrix}
2\\
&w+w^{-1}\\
&&w^2+w^{-2}\\
&&&\ddots\\
&&&&w^{N-1}+w^{1-N}
\end{pmatrix}F_N^*$$
where $F_N=(w^{ij})$ with $w=e^{2\pi i/N}$ is as usual the Fourier matrix.
\end{theorem}

\begin{proof}
In terms of the vector $\xi=(w^p)$, and as usual with indices $p=0,1,\ldots,N-1$, the eigenvector formula that we found in Proposition 2.20 reads:
$$d\xi^k=(w^k+w^{-k})\xi^k$$

By putting all these formulae together, with $k=0,1,\ldots,N-1$, we obtain:
\begin{eqnarray*}
&&d\begin{pmatrix}\xi^0&\xi^1&\xi^2&\ldots&\xi^{N-1}\end{pmatrix}\\
&=&\begin{pmatrix}2\xi^0&(w+w^{-1})\xi&(w^2+w^{-2})\xi^2&\ldots&(w^{N-1}+w^{1-N})\xi^{N-1}\end{pmatrix}\\
&=&\begin{pmatrix}\xi^0&\xi^1&\xi^2&\ldots&\xi^{N-1}\end{pmatrix}
\begin{pmatrix}
2\\
&w+w^{-1}\\
&&w^2+w^{-2}\\
&&&\ddots\\
&&&&w^{N-1}+w^{1-N}
\end{pmatrix}
\end{eqnarray*}

But, since we have $F_N=(\xi^k)$, as a square matrix, this formula reads:
$$dF_N=F_N\begin{pmatrix}
2\\
&w+w^{-1}\\
&&w^2+w^{-2}\\
&&&\ddots\\
&&&&w^{N-1}+w^{1-N}
\end{pmatrix}$$

Now by multipliying to the right by $F_N^{-1}=F_N^*/N$, we obtain the result.
\end{proof}

The above result is quite interesting, and as a first obvious conclusion, coming from this and from our previous results regarding the simplex, we have:

\begin{conclusion}
The adjacency matrices of both the simplex and the circle are diagonalized by the Fourier matrix $F_N$.
\end{conclusion}

We will be back to this phenomenon later in this book, with a systematic discussion of the graphs diagonalized by the Fourier matrix $F_N$, and of related topics.

\bigskip

In the meantime, as a corollary of this, remember Proposition 2.18, regarding the complementary graphs $X^c$, with the quite technical assumptions that we found there, seemingly bound for the trash can? Well, as good news, Conclusion 2.22 is exactly what we need, guaranteeing that Proposition 2.18 applies indeed, to the circle graph. So, as our next theorem, based on the above study for the circle, we have:

\index{complement of circle}

\begin{theorem}
The adjacency matrix of the complement of the circle graph diagonalizes as
$$d=\frac{1}{N}F_N
\begin{pmatrix}
N-3\\
&-1-w-w^{-1}\\
&&-1-w^2-w^{-2}\\
&&&\ddots\\
&&&&-1-w^{N-1}-w^{1-N}
\end{pmatrix}F_N^*$$
where $F_N=(w^{ij})$ with $w=e^{2\pi i/N}$ is as usual the Fourier matrix.
\end{theorem}

\begin{proof}
Consider the circle graph $X$, that we know well from the above:
$$
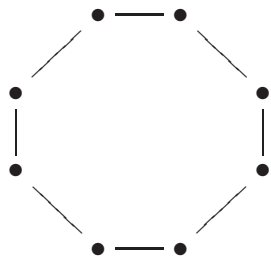
\xymatrix@R=17pt@C=18pt{
&\bullet\ar@{-}[r]\ar@{-}[dl]&\bullet\ar@{-}[dr]\\
\bullet\ar@{-}[d]&&&\bullet\ar@{-}[d]\\
\bullet\ar@{-}[dr]&&&\bullet\ar@{-}[dl]\\
&\bullet\ar@{-}[r]&\bullet}$$

Consider as well its complement $X^c$, which is a quite crowded and scary graph, of valence $N-3$, with the picture at $N=8$ being as follows:
$$
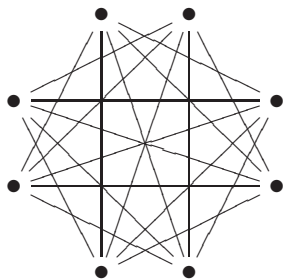
\xymatrix@R=20pt@C=20pt{
&\bullet\ar@{-}[ddrr]\ar@{-}[drr]\ar@{-}[ddd]\ar@{-}[dddr]&\bullet\ar@{-}[ddd]\ar@{-}[dll]\ar@{-}[ddll]\ar@{-}[dddl]\\
\bullet\ar@{-}[ddr]&&&\bullet\ar@{-}[lll]\ar@{-}[dlll]\ar@{-}[ddll]\\
\bullet\ar@{-}[uur]\ar@{-}[drr]&&&\bullet\ar@{-}[lll]\ar@{-}[uul]\ar@{-}[ulll]\\
&\bullet\ar@{-}[urr]&\bullet\ar@{-}[uur]\ar@{-}[uull]}$$

However, at the linear algebra level things are quite simple, because by combining Theorem 2.14 and Theorem 2.21, as indicated in Proposition 2.18, we obtain:
\begin{eqnarray*}
d_{X^c}
&=&\mathbb I_N-1_N-d_X\\
&=&\frac{1}{N}F_N\begin{pmatrix}
N-1&&&&0\\
&-1\\
&&\ddots&\\
&&&-1\\
0&&&&-1
\end{pmatrix}F_N^*\\
&&-\frac{1}{N}F_N
\begin{pmatrix}
2\\
&w+w^{-1}\\
&&w^2+w^{-2}\\
&&&\ddots\\
&&&&w^{N-1}+w^{1-N}
\end{pmatrix}F_N^*
\end{eqnarray*}

Thus, we are led to the formula in the statement.
\end{proof}

Very nice all this, so done with the circle and its complement, and on our to-do list, job for later to further explore Conclusion 2.22, and related Fourier analysis topics.

\section*{2d. The segment} 

As a continuation of the above study, let us discuss now the segment graph. This is a graph that we already met in chapter 1, in relation with random walk questions, and with our conclusions there being quite muddy, basically amounting in saying that this graph is something quite tough to fight with, knowing more mathematics than we do.

\bigskip

We will investigate here this graph, from the diagonalization point of view. As a first observation, since diagonalization solves the random walk question, as explained in chapter 1, we cannot expect much from our diagonalization study, done as we do, with bare hands. So, modesty, let us just start working on this, and we'll see what we get.

\bigskip

Some numerics first. These do not look very good, the result being as follows:

\begin{proposition}
The eigenvalues of the segment graph are as follows:
\begin{enumerate}
\item At $N=2$ we have $-1,1$.

\item At $N=3$ we have $0,\pm\sqrt{2}$.

\item At $N=4$ we have $\pm\frac{\sqrt{5}\pm1}{2}$.

\item At $N=5$ we have $0,\pm1,\pm\sqrt{3}$.

\item At $N=6$ we have the solutions of $x^6-5x^4+6x^2-1=0$.
\end{enumerate}
\end{proposition}

\begin{proof}
These are some straightforward linear algebra computations, with some tricks being needed only at $N=4$, the details being as follows:

\medskip

(1) At $N=2$ the adjacency matrix and its eigenvalues are as follows:
$$d=\begin{pmatrix}
0&1\\
1&0
\end{pmatrix}\quad,\quad
\begin{vmatrix}
x&-1\\
-1&x
\end{vmatrix}
=x^2-1=(x-1)(x+1)$$

(2) At $N=3$ the adjacency matrix and its eigenvalues are as follows:
$$d=\begin{pmatrix}
0&1&0\\
1&0&1\\
0&1&0
\end{pmatrix}\quad,\quad
\begin{vmatrix}
x&-1&0\\
-1&x&-1\\
0&-1&x
\end{vmatrix}
=x^3-2x=x(x^2-2)$$

(3) At $N=4$ the adjacency matrix and its characteristic polynomial are as follows:
$$d=\begin{pmatrix}
0&1&0&0\\
1&0&1&0\\
0&1&0&1\\
0&0&1&0
\end{pmatrix}\quad,\quad
\begin{vmatrix}
x&-1&0&0\\
-1&x&-1&0\\
0&-1&x&-1\\
0&0&-1&x
\end{vmatrix}
=x^4-3x^2+1$$

Now by solving the degree 2 equation, and fine-tuning the answer, we obtain:
$$x=\pm\sqrt{\frac{3\pm\sqrt{5}}{2}}=\pm\frac{\sqrt{5}\pm1}{2}$$

(4) At $N=5$ the adjacency matrix and its eigenvalues are as follows:
$$d=\begin{pmatrix}
0&1&0&0&0\\
1&0&1&0&0\\
0&1&0&1&0\\
0&0&1&0&1\\
0&0&0&1&0
\end{pmatrix}\quad,\quad
\begin{vmatrix}
x&-1&0&0&0\\
-1&x&-1&0&0\\
0&-1&x&-1&0\\
0&0&-1&x&-1\\
0&0&0&-1&x
\end{vmatrix}
=x(x^2-1)(x^2-3)$$

(5) At $N=6$ the adjacency matrix and its eigenvalues are as follows:
$$d=\begin{pmatrix}
0&1&0&0&0&0\\
1&0&1&0&0&0\\
0&1&0&1&0&0\\
0&0&1&0&1&0\\
0&0&0&1&0&1\\
0&0&0&0&1&0
\end{pmatrix}\ ,\ 
\begin{vmatrix}
x&-1&0&0&0&0\\
-1&x&-1&0&0&0\\
0&-1&x&-1&0&0\\
0&0&-1&x&-1&0\\
0&0&0&-1&x&-1\\
0&0&0&0&-1&x
\end{vmatrix}
=x^6-5x^4+6x^2-1$$

Thus, we are led to the formulae in the statement.
\end{proof}

All the above does not look very good. However, as a matter of having the dirty job fully done, with mathematical pride, let us look as well into eigenvectors. At $N=2$ things are quickly settled, with the diagonalization of the adjacency matrix being as follows:
$$\begin{pmatrix}
0&1\\
1&0
\end{pmatrix}
=\frac{1}{2}\begin{pmatrix}
1&1\\
1&-1
\end{pmatrix}
\begin{pmatrix}
1&0\\
0&-1
\end{pmatrix}
\begin{pmatrix}
1&1\\
1&-1
\end{pmatrix}$$

At $N=3$ the diagonalization formula becomes more complicated, as follows:
$$\begin{pmatrix}
0&1&0\\
1&0&1\\
0&1&0
\end{pmatrix}
=\frac{1}{4}\begin{pmatrix}
1&1&1\\
0&\sqrt{2}&-\sqrt{2}\\
-1&1&1
\end{pmatrix}
\begin{pmatrix}
0&0&0\\
0&\sqrt{2}&0\\
0&0&-\sqrt{2}
\end{pmatrix}
\begin{pmatrix}
2&0&-2\\
1&\sqrt{2}&1\\
1&-\sqrt{2}&1
\end{pmatrix}$$

At $N=4$ now, in view of the eigenvalue formula that we found, $x^4-3x^2+1=0$, we must proceed with care. The equation for the eigenvectors $dv=xv$ is as follows:
$$\begin{pmatrix}
0&1&0&0\\
1&0&1&0\\
0&1&0&1\\
0&0&1&0
\end{pmatrix}
\begin{pmatrix}
a\\
b\\
c\\
d
\end{pmatrix}
=x\begin{pmatrix}
a\\
b\\
c\\
d
\end{pmatrix}$$

In other words, the equations to be satisfied are as follows:
$$b=xa$$
$$a+c=xb$$
$$b+d=xc$$
$$c=xd$$

With the choice $a=1$, the solutions of these equations are as follows:
$$a=1$$
$$b=x$$
$$c=x^2-1$$
$$d=(x^2-1)/x$$

In order to compute now the $c$ components of the eigenvectors, we can use the formula $x=(\pm1\pm\sqrt{5})/2$ from Proposition 2.24. Indeed, this formula gives:
\begin{eqnarray*}
\left(\frac{\pm1\pm\sqrt{5}}{2}\right)^2-1
&=&\frac{6\pm 2\sqrt{5}}{4}-1\\
&=&\frac{3\pm\sqrt{5}}{2}-1\\
&=&\frac{1\pm\sqrt{5}}{2}
\end{eqnarray*}

Thus, almost done, and we deduce that the passage matrix is as follows:
$$P=\begin{pmatrix}
1&1&1&1\\
\frac{1+\sqrt{5}}{2}&\frac{1-\sqrt{5}}{2}&\frac{-1+\sqrt{5}}{2}&\frac{-1-\sqrt{5}}{2}\\
\frac{1+\sqrt{5}}{2}&\frac{1-\sqrt{5}}{2}&\frac{1-\sqrt{5}}{2}&\frac{1+\sqrt{5}}{2}\\
1&1&-1&-1
\end{pmatrix}$$

To be more precise, according to our equations above, the first row must consist of $a=1$ entries. Then the second row must consist of $b=x$ entries, with $x=(\pm1\pm\sqrt{5})/2$. Then the third row must consist of $c=x^2-1$ entries, but these are easily computed, as explained above. Finally, the fourth row must consist of $d=(x^2-1)/x$ entries, which means that the fourth row appears by dividing the third row by the second row, which is easily done too. In case you wonder why $d=\pm1$, here is another proof of this:
\begin{eqnarray*}
d=\pm1
&\iff&x^2-1=\pm x\\
&\iff&(x^2-1)^2=x^2\\
&\iff&x^4-2x^2+1=x^2\\
&\iff&x^4-3x^2+1=0
\end{eqnarray*}

Very nice all this, and leaving the computation of $P^{-1}$ for you, here are a few more observations, in relation with what we found in Proposition 2.24:

\bigskip

(1) In the case $N=5$ the eigenvectors can be computed too, and the diagonalization finished, via the standard method, namely system of equations, and with the numerics involving powers of the eigenvalues that we found. Exercise for you.

\bigskip

(2) The same stays true at $N=6$, again with the eigenvector numerics involving powers of the eigenvalues, and with these eigenvalues being explicitly computable, via the Cardano formula for degree 3 equations. Have fun with this too, of course.

\bigskip

(3) However, all this does not look very good, and at $N=7$ and higher we will certainly run into difficult questions, and save for some interesting remarks, normally depending on the parity of $N$, we will not be able to fully solve the diagonalization problem.

\bigskip

So, what to do? Work some more, of course, the hard way. Proposition 2.24 and its proof look quite trivial, but if you get into the full details of the computations, and let me assign here this, as a key exercise, you will certainly notice that, when computing that determinants, you can sort of use a recurrence method. And, this leads to:

\index{circle graph}
\index{Chebycheff polynomials}
\index{orthogonal polynomials}

\begin{theorem}
The characteristic polynomials $P_N$ of the segment graphs satisfy
$$P_0=1\quad,\quad P_1=x\quad,\quad P_{N+1}=xP_N-P_{N-1}$$
and are the well-known Chebycheff polynomials, enjoying lots of interesting properties.
\end{theorem} 

\begin{proof}
Obviously, many things going on here, ranging from precise to definitional, or even informal, the idea with all this being as follows:

\medskip

(1) By computing determinants, as indicated above, we are led to the recurrence formula in the statement. Here is the proof at $N=4$, the general case being similar:
\begin{eqnarray*}
P_5
&=&\begin{vmatrix}
x&-1&0&0&0\\
-1&x&-1&0&0\\
0&-1&x&-1&0\\
0&0&-1&x&-1\\
0&0&0&-1&x
\end{vmatrix}\\
&=&x\begin{vmatrix}
x&-1&0&0\\
-1&x&-1&0\\
0&-1&x&-1\\
0&0&-1&x
\end{vmatrix}
+\begin{vmatrix}
-1&-1&0&0\\
0&x&-1&0\\
0&-1&x&-1\\
0&0&-1&x
\end{vmatrix}\\
&=&x\begin{vmatrix}
x&-1&0&0\\
-1&x&-1&0\\
0&-1&x&-1\\
0&0&-1&x
\end{vmatrix}
-\begin{vmatrix}
x&-1&0\\
-1&x&-1\\
0&-1&x
\end{vmatrix}\\
&=&xP_4-P_3
\end{eqnarray*}

(2) Regarding now the intial values, according to Proposition 2.24 these should normally be $P_2=x^2-1$ and $P_3=x^3-2x$, but we can formally add the values $P_0=1$ and $P_1=x$, as for the final statement to look better, with this being justified by:
$$P_0=1$$
$$P_1=x$$
$$P_2=xP_1-P_0=x^2-1$$
$$P_3=xP_2-P_1=x^3-2x$$

(3) Thus, we have the recurrence formula in the statement, and with the initial values there, and an internet search, or some advanced calculus know-how, tells us that these are the well-known Chebycheff polynomials, enjoying lots of interesting properties.
\end{proof}

Very nice all this, and as a continuation, barring as usual the wish to fully solve the diagonalization problem, which remains something difficult, we can for instance go back now to the random walk computations from chapter 1, for the segment graph, with the various choices there for the distinguished vertex, and improve our results there.

\bigskip

However, while this is certainly something doable and interesting, we will not do it here, and rather keep the above as a mere remark, or as an exercise for you. This is because we will be back to all this in the next chapter, with some even better tools.

\section*{2e. Exercises}

We had a tough chapter here, mixing linear algebra with calculus, and with some probability too. Here are some exercises, in relation with the above:

\begin{exercise}
Learn the spectral theorem for general normal matrices.
\end{exercise}

\begin{exercise}
Furher build on the spectral characterization of graphs.
\end{exercise}

\begin{exercise}
Find a geometric proof of the spectral theorem, for graphs.
\end{exercise}

\begin{exercise}
View the Fourier matrix as matrix of a Fourier transform.
\end{exercise}

\begin{exercise}
Further build on our results regarding complementary graphs.
\end{exercise}

\begin{exercise}
Diagonalize the adjacency matrix of the segment, at $N=5,6,7$.
\end{exercise}

\begin{exercise}
Learn more about Chebycheff polynomials, and their properties.
\end{exercise}

\begin{exercise}
Study random walks on segments, using Chebycheff polynomials.
\end{exercise}

As bonus exercise, learn some probability. In what follows we will often get into probability, with the basics recalled, but nothing replaces some systematic learning.

\chapter{Random walks}

\section*{3a. Walks, measures}

We have learned so far the basics of graph theory and probability, and time now to see if this knowledge can be of any help, in relation with concrete questions. The question that we would like to discuss, which is something very basic, is as follows:

\begin{question}
Given a graph $X$, with a distinguished vertex $*$:
\begin{enumerate}
\item What is the number $L_k$ of length $k$ loops on $X$, based at $*$? 

\item Equivalently, what is the measure $\mu$ having $L_k$ as moments?
\end{enumerate}
\end{question}

To be more precise, we are mainly interested in the first question, counting loops on graphs, with this being notoriously related to many applied mathematics questions, of discrete type. As for the second question, this is a technical, useful probabilistic reformulation of the first question, that we will usually prefer, in what follows. 

\bigskip

Actually, in relation with this, the fact that a measure $\mu$ as above exists indeed is not exactly obvious. But comes from the following result, which is something rather elementary, and which can be very helpful for explicit computations:

\index{adjacency matrix}
\index{random walk}
\index{spectral measure}
\index{Dirac mass}

\begin{theorem}
Given a graph $X$, with adjacency matrix $d\in M_N(0,1)$, we have:
$$L_k=(d^k)_{**}$$
When writing $d=UDU^t$ with $U\in O_N$ and $D=diag(\lambda_1,\ldots,\lambda_N)$ with $\lambda_i\in\mathbb R$, we have
$$L_k=\sum_iU_{*i}^2\lambda_i^k$$
and the real probability measure $\mu$ having these numbers as moments is given by
$$\mu=\sum_iU_{*i}^2\delta_{\lambda_i}$$
with the delta symbols standing as usual for Dirac masses.
\end{theorem}

\begin{proof}
There are several things going on here, the idea being as follows:

\medskip

(1) According to the usual rule of matrix multiplication, the formula for the powers of the adjacency matrix $d\in M_N(0,1)$ is as follows:
\begin{eqnarray*}
(d^k)_{i_0i_k}
&=&\sum_{i_1,\ldots,i_{k-1}}d_{i_0i_1}d_{i_1i_2}\ldots d_{i_{k-1}i_k}\\
&=&\sum_{i_1,\ldots,i_{k-1}}\delta_{i_0-i_1}\delta_{i_1-i_2}\ldots\delta_{i_{k-1}-i_k}\\
&=&\sum_{i_1,\ldots,i_{k-1}}\delta_{i_0-i_1-\ldots-i_{k-1}-i_k}\\
&=&\#\Big\{i_0-i_1-\ldots-i_{k-1}-i_k\Big\}
\end{eqnarray*}

In particular, with $i_0=i_k=*$, we obtain the following formula, as claimed:
$$(d^k)_{**}=\#\Big\{\!*-\,i_1-\ldots-i_{k-1}-*\Big\}=L_k$$

(2) Now since the adjacency matrix $d\in M_N(0,1)$ is symmetric, by basic linear algebra, that we know from chapter 2, this matrix is diagonalizable, with the diagonalization being as follows, with $U\in O_N$, and $D=diag(\lambda_1,\ldots,\lambda_N)$ with $\lambda_i\in\mathbb R$:
$$d=UDU^t$$

By using this formula, we obtain the second formula in the statement:
\begin{eqnarray*}
L_k
&=&(d^k)_{**}\\
&=&(UD^kU^t)_{**}\\
&=&\sum_iU_{*i}\lambda_i^k(U^t)_{i*}\\
&=&\sum_iU_{*i}^2\lambda_i^k
\end{eqnarray*}

(3) Finally, the last assertion is clear from this, because the moments of the measure in the statement, $\mu=\sum_iU_{*i}^2\delta_{\lambda_i}$, are the following numbers:
\begin{eqnarray*}
M_k
&=&\int_\mathbb Rx^kd\mu(x)\\
&=&\sum_iU_{*i}^2\lambda_i^k\\
&=&L_k
\end{eqnarray*}

Observe also that $\mu$ is indeed of mass 1, because all rows of $U\in O_N$ must be of norm 1, and so $\sum_iU_{*i}^2=1$. Thus, we are led to the conclusions in the statement.
\end{proof}

At the level of examples now, we have already seen some computations in chapter 1, so let us first review that material. We first have the following result, from there:

\begin{theorem}
For the circle graph, having $N$ vertices, the number of length $k$ loops based at one of the vertices is approximately
$$L_k\simeq\frac{2^k}{N}$$
in the $k\to\infty$ limit, when $N$ is odd, and is approximately
$$L_k\simeq\begin{cases}
\frac{2^{k+1}}{N}&(k\ {\rm even})\\
0&(k\ {\rm odd})
\end{cases}$$
also with $k\to\infty$, when $N$ is even.
\end{theorem}

\begin{proof}
This is something that we know from chapter 1, the idea being as follows:

\medskip

(1) Consider the circle graph $X$, with vertices denoted $0,1,\ldots,N-1$. Since each vertex has valence 2, any length $k$ path based at 0 will consist of a binary choice at the beginning, then another binary choice afterwards, and so on up to a $k$-th binary choice at the end. Thus, there is a total of $2^k$ such paths, based at 0, and having length $k$. 

\medskip

(2) But now, based on the obvious ``uniformity'' of the circle, we can argue that, in the $k\to\infty$ limit, the endpoint of such a path will become random among the vertices $0,1,\ldots,N-1$. Thus, if we want this endpoint to be 0, as to have a loop, we have $1/N$ chances for this to happen, so the total number of loops is $L_k\simeq 2^k/N$, as stated.

\medskip

(3) With the remark, however, that the above argument works fine only when $N$ is odd. Indeed, when $N$ is even, the endpoint of a length $k$ path will be random among $0,2,\ldots,2N-2$ when $k$ is even, and random among $1,3,\ldots,2N-1$ when $k$ is odd. Thus for getting a loop we must assume that $k$ is even, and in this case the number of such loops is the total number of length $k$ paths, namely $2^k$, approximately divided by $N/2$, the number of points in $\{0,2,\ldots,2N-2\}$, which gives $L_k=2^k/(N/2)$, as stated.

\medskip

(4) So, this was what we already knew from chapter 1, but in view of our present more advanced knowledge of the segment graph, coming from the computations from chapter 2, we can now say more about all this. And, more on this later in this chapter.
\end{proof}

Another graph that we studied in chapter 1 was $\mathbb Z$, the result being as follows:

\index{central binomial coefficients}

\begin{theorem}
The paths on $\mathbb Z$ are counted by the binomial coefficients. In particular, the $2k$-paths based at $0$ are counted by the central binomial coefficients,
$$L_{2k}=\binom{2k}{k}$$
and $\mu$ is the centered measure having these numbers as even moments.
\end{theorem}

\begin{proof}
This is something elementary, the idea being as follows:

\medskip

(1) Let us try to count the length $k$ paths on $\mathbb Z$, based at $0$. At $k=1,2,3,4$ the count results can be pictured as follows, in a self-explanatory way:
$$
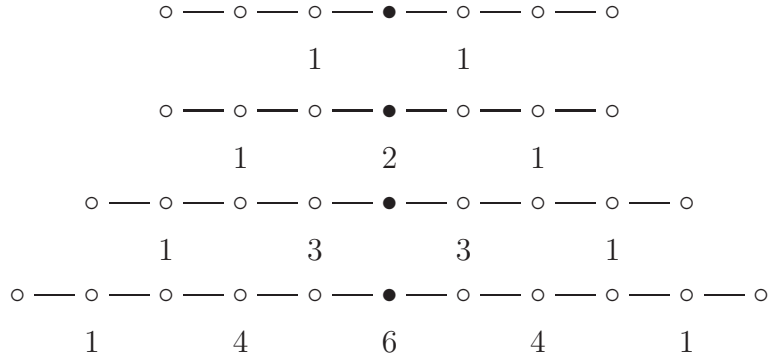
\xymatrix@R=5pt@C=15pt{
\circ\ar@{-}[r]&\circ\ar@{-}[r]&\circ\ar@{-}[r]&\bullet\ar@{-}[r]&\circ\ar@{-}[r]&\circ\ar@{-}[r]&\circ\\
&&1&&1
}$$
$$\xymatrix@R=5pt@C=15pt{
\circ\ar@{-}[r]&\circ\ar@{-}[r]&\circ\ar@{-}[r]&\bullet\ar@{-}[r]&\circ\ar@{-}[r]&\circ\ar@{-}[r]&\circ\\
&1&&2&&1
}$$
$$\xymatrix@R=5pt@C=15pt{
\circ\ar@{-}[r]&\circ\ar@{-}[r]&\circ\ar@{-}[r]&\circ\ar@{-}[r]&\bullet\ar@{-}[r]&\circ\ar@{-}[r]&\circ\ar@{-}[r]&\circ\ar@{-}[r]&\circ\\
&1&&3&&3&&1
}$$
$$\xymatrix@R=5pt@C=15pt{
\circ\ar@{-}[r]&\circ\ar@{-}[r]&\circ\ar@{-}[r]&\circ\ar@{-}[r]&\circ\ar@{-}[r]&\bullet\ar@{-}[r]&\circ\ar@{-}[r]&\circ\ar@{-}[r]&\circ\ar@{-}[r]&\circ\ar@{-}[r]&\circ\\
&1&&4&&6&&4&&1
}$$

(2) We recognize here the Pascal triangle, and the rest is just a matter of finishing. For instance we can argue that the number $C_k^l$ of length $k$ loops $0\to l$  is subject, due to the binary choice at the end, to the following recurrence relation:
$$C_k^l=C_{k-1}^{l-1}+C_{k-1}^{l+1}$$

But this is exactly the recurrence relation for the Pascal triangle, as desired. As for the last assertion, this is more of an empty statement, with $\mu$ still to be computed.
\end{proof}

\section*{3b. Catalan numbers}

As a third example, let us try to count the loops of $\mathbb N$, based at 0. This is something less obvious, and at the experimental level, the result is as follows:

\index{Catalan numbers}

\begin{proposition}
The Catalan numbers $C_k$, counting the loops on $\mathbb N$ based at $0$,
$$C_k=\#\Big\{0-i_1-\ldots-i_{2k-1}-0\Big\}$$
are numerically $1,2,5,14,42,132,429,1430,4862,16796,58786,\ldots$
\end{proposition}

\begin{proof}
To start with, we have indeed $C_1=1$, the only loop here being $0-1-0$. Then we have $C_2=2$, due to two possible loops, namely:
$$0-1-0-1-0$$
$$0-1-2-1-0$$

Then we have $C_3=5$, the possible loops here being as follows:
$$0-1-0-1-0-1-0$$
$$0-1-0-1-2-1-0$$
$$0-1-2-1-0-1-0$$
$$0-1-2-1-2-1-0$$
$$0-1-2-3-2-1-0$$

In general, the same method works, with $C_4=14$ being left to you, as an exercise, and with $C_5$ and higher to me, and I will be back with the solution, in due time.
\end{proof}

Obviously, computing the numbers $C_k$ is no easy task, and finding the formula of $C_k$, out of the data that we have, does not look as an easy task either. So, we will do what combinatorists do, let me teach you. The first step is to relax, then to look around, not with the aim of computing your numbers $C_k$, but rather with the aim of finding other objects counted by the same numbers $C_k$. With a bit of luck, among these objects some will be easier to count than the others, and this will eventually compute $C_k$.

\bigskip

This was for the strategy. In practice now, we first have the following result:

\index{Dyck paths}
\index{noncrossing pairings}
\index{noncrossing partitions}

\begin{theorem}
The Catalan numbers $C_k$ count:
\begin{enumerate}
\item The length $2k$ loops on $\mathbb N$, based at $0$.

\item The noncrossing pairings of $1,\ldots,2k$.

\item The noncrossing partitions of $1,\ldots,k$.

\item The length $2k$ Dyck paths in the plane.
\end{enumerate}
\end{theorem}

\begin{proof}
All this is standard combinatorics, the idea being as follows:

\medskip

(1) To start with, in what regards the various objects involved, the length $2k$ loops on $\mathbb N$ are the length $2k$ loops on $\mathbb N$ that we know, and the same goes for the noncrossing pairings of $1,\ldots,2k$, and for the noncrossing partitions of $1,\ldots,k$, the idea here being that you must be able to draw the pairing or partition in a noncrossing way. 

\medskip

(2) Regarding now the length $2k$ Dyck paths in the plane, these are by definition the paths from $(0,0)$ to $(k,k)$, marching North-East over the integer lattice $\mathbb Z^2\subset\mathbb R^2$, by staying inside the square $[0,k]\times[0,k]$, and staying as well under the diagonal of this square. As an example, here are the 5 possible Dyck paths at $n=3$:
$$
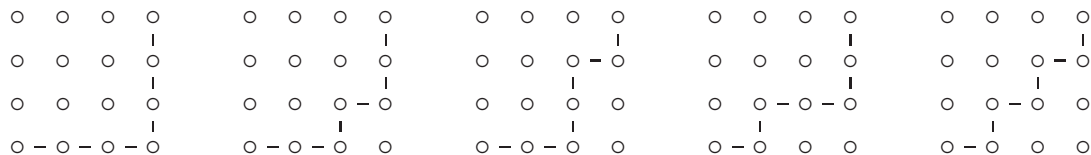
\xymatrix@R=4pt@C=4pt
{\circ&\circ&\circ&\circ\\
\circ&\circ&\circ&\circ\ar@{-}[u]\\
\circ&\circ&\circ&\circ\ar@{-}[u]\\
\circ\ar@{-}[r]&\circ\ar@{-}[r]&\circ\ar@{-}[r]&\circ\ar@{-}[u]}
\qquad
\xymatrix@R=4pt@C=4pt
{\circ&\circ&\circ&\circ\\
\circ&\circ&\circ&\circ\ar@{-}[u]\\
\circ&\circ&\circ\ar@{-}[r]&\circ\ar@{-}[u]\\
\circ\ar@{-}[r]&\circ\ar@{-}[r]&\circ\ar@{-}[u]&\circ}
\qquad
\xymatrix@R=4pt@C=4pt
{\circ&\circ&\circ&\circ\\
\circ&\circ&\circ\ar@{-}[r]&\circ\ar@{-}[u]\\
\circ&\circ&\circ\ar@{-}[u]&\circ\\
\circ\ar@{-}[r]&\circ\ar@{-}[r]&\circ\ar@{-}[u]&\circ}
\qquad
\xymatrix@R=4pt@C=4pt
{\circ&\circ&\circ&\circ\\
\circ&\circ&\circ&\circ\ar@{-}[u]\\
\circ&\circ\ar@{-}[r]&\circ\ar@{-}[r]&\circ\ar@{-}[u]\\
\circ\ar@{-}[r]&\circ\ar@{-}[u]&\circ&\circ}
\qquad
\xymatrix@R=4pt@C=4pt
{\circ&\circ&\circ&\circ\\
\circ&\circ&\circ\ar@{-}[r]&\circ\ar@{-}[u]\\
\circ&\circ\ar@{-}[r]&\circ\ar@{-}[u]&\circ\\
\circ\ar@{-}[r]&\circ\ar@{-}[u]&\circ&\circ}
$$

(3) Thus, we have definitions for all the objects involved, and in each case, if you start counting them, as we did in Proposition 3.5 with the loops on $\mathbb N$, you always end up with the same sequence of numbers, namely those found in Proposition 3.5:
$$1,2,5,14,42,132,429,1430,4862,16796,58786,\ldots$$

(4) In order to prove now that (1-4) produce indeed the same numbers, many things can be said. The idea is that, leaving aside mathematical brevity, and more specifically abstract reasonings of type $a=b,b=c\implies a=c$, what we have to do, in order to fully understand what is going on, is to etablish $\binom{4}{2}=6$ equalities, via bijective proofs.

\medskip

(5) But this can be done, indeed. As an example here, the noncrossing pairings of $1,\ldots,2k$ from (2) are in bijection with the noncrossing partitions of $1,\ldots,k$ from (3), via  fattening the pairings and shrinking the partitions. We will leave the details here as an instructive exercise, and exercise as well, to add (1) and (4) to the picture.

\medskip

(6) However, matter of having our theorem formally proved, I mean by me professor and not by you student, here is a less elegant argument, which is however very quick, and does the job. The point is that, in each of the cases (1-4) under consideration, the numbers $C_k$ that we get are easily seen to be subject to the following recurrence:
$$C_{k+1}=\sum_{a+b=k}C_aC_b$$ 

The initial data being the same, namely $C_1=1$ and $C_2=2$, in each of the cases (1-4) under consideration, we get indeed the same numbers.
\end{proof}

Now we can pass to the second step, namely selecting in the above list the objects that we find the most convenient to count, and count them. This leads to:

\index{Catalan numbers}

\begin{theorem}
The Catalan numbers are given by the formula
$$C_k=\frac{1}{k+1}\binom{2k}{k}$$
with this being best seen by counting the length $2k$ Dyck paths in the plane.
\end{theorem}

\begin{proof}
This is something quite tricky, the idea being as follows:

\medskip

(1) Let us count indeed the Dyck paths in the plane. For this purpose, we use a trick. Indeed, if we ignore the assumption that our path must stay under the diagonal of the square, we have $\binom{2k}{k}$ such paths. And among these, we have the ``good'' ones, those that we want to count, and then the ``bad'' ones, those that we want to ignore.

\medskip

(2) So, let us count the bad paths, those crossing the diagonal of the square, and reaching the higher diagonal next to it, the one joining $(0,1)$ and $(k,k+1)$. In order to count these, the trick is to ``flip'' their bad part over that higher diagonal, as follows:
$$
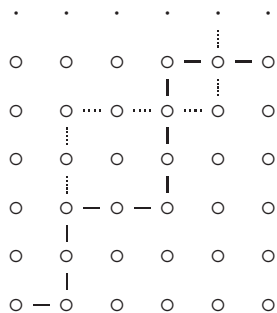
\xymatrix@R=6pt@C=6pt
{\cdot&\cdot&\cdot&\cdot&\cdot&\cdot\\
\circ&\circ&\circ&\circ\ar@{-}[r]&\circ\ar@{-}[r]\ar@{.}[u]&\circ\\
\circ&\circ\ar@{.}[r]&\circ\ar@{.}[r]&\circ\ar@{.}[r]\ar@{-}[u]&\circ\ar@{.}[u]&\circ\\
\circ&\circ\ar@{.}[u]&\circ&\circ\ar@{-}[u]&\circ&\circ\\
\circ&\circ\ar@{-}[r]\ar@{.}[u]&\circ\ar@{-}[r]&\circ\ar@{-}[u]&\circ&\circ\\
\circ&\circ\ar@{-}[u]&\circ&\circ&\circ&\circ\\
\circ\ar@{-}[r]&\circ\ar@{-}[u]&\circ&\circ&\circ&\circ}$$

(3) Now observe that, as it is obvious on the above picture, due to the flipping, the flipped bad path will no longer end in $(k,k)$, but rather in $(k-1,k+1)$. Moreover, more is true, in the sense that, by thinking a bit, we see that the flipped bad paths are precisely those ending in $(k-1,k+1)$. Thus, we can count these flipped bad paths, and so the bad paths, and so the good paths too, and so good news, we are done.

\medskip

(4) To finish now, by putting everything together, we have:
\begin{eqnarray*}
C_k
&=&\binom{2k}{k}-\binom{2k}{k-1}\\
&=&\binom{2k}{k}-\frac{k}{k+1}\binom{2k}{k}\\
&=&\frac{1}{k+1}\binom{2k}{k}
\end{eqnarray*}

Thus, we are led to the formula in the statement.
\end{proof}

We have as well another approach to all this, computation of the Catalan numbers, this time based on rock-solid standard calculus, as follows:

\index{Catalan numbers}

\begin{theorem}
The Catalan numbers have the following properties:
\begin{enumerate}
\item They satisfy $C_{k+1}=\sum_{a+b=k}C_aC_b$.

\item The series $f(z)=\sum_{k\geq0}C_kz^k$ satisfies $zf^2-f+1=0$.

\item This series is given by $f(z)=\frac{1-\sqrt{1-4z}}{2z}$.

\item We have the formula $C_k=\frac{1}{k+1}\binom{2k}{k}$.
\end{enumerate}
\end{theorem}

\begin{proof}
This is best viewed by using noncrossing pairings, as follows: 

\medskip

(1) Let us count the noncrossing pairings of $\{1,\ldots,2k+2\}$. Such a pairing appears by pairing 1 to an odd number, $2a+1$, and then inserting a noncrossing pairing of $\{2,\ldots,2a\}$, and a noncrossing pairing of $\{2a+2,\ldots,2k+2\}$. Thus we have, as claimed:
$$C_{k+1}=\sum_{a+b=k}C_aC_b$$ 

(2) Consider now the generating series of the Catalan numbers, $f(z)=\sum_{k\geq0}C_kz^k$. In terms of this generating series, the above recurrence gives, as desired:
\begin{eqnarray*}
zf^2
&=&\sum_{a,b\geq0}C_aC_bz^{a+b+1}\\
&=&\sum_{k\geq1}\sum_{a+b=k-1}C_aC_bz^k\\
&=&\sum_{k\geq1}C_kz^k\\
&=&f-1
\end{eqnarray*}

(3) By solving the equation $zf^2-f+1=0$ found above, and choosing the solution which is bounded at $z=0$, we obtain the following formula, as claimed:
$$f(z)=\frac{1-\sqrt{1-4z}}{2z}$$ 

(4) In order to compute this function, we use the generalized binomial formula, which is as follows, with $p\in\mathbb R$ being an arbitrary exponent, and with $|t|<1$:
$$(1+t)^p=\sum_{k=0}^\infty\binom{p}{k}t^k$$

To be more precise, this formula, which generalizes the usual binomial formula, holds indeed due to the Taylor formula, with the binomial coefficients being given by:
$$\binom{p}{k}=\frac{p(p-1)\ldots(p-k+1)}{k!}$$

(5) For the exponent $p=1/2$, the generalized binomial coefficients are:
\begin{eqnarray*}
\binom{1/2}{k}
&=&\frac{1/2(-1/2)(-3/2)\ldots(3/2-k)}{k!}\\
&=&(-1)^{k-1}\frac{1\cdot 3\cdot 5\ldots(2k-3)}{2^kk!}\\
&=&(-1)^{k-1}\frac{(2k-2)!}{2^{k-1}(k-1)!2^kk!}\\
&=&\frac{(-1)^{k-1}}{2^{2k-1}}\cdot\frac{1}{k}\binom{2k-2}{k-1}\\
&=&-2\left(\frac{-1}{4}\right)^k\cdot\frac{1}{k}\binom{2k-2}{k-1}
\end{eqnarray*}

(6) Thus the generalized binomial formula at exponent $p=1/2$ reads:
$$\sqrt{1+t}=1-2\sum_{k=1}^\infty\frac{1}{k}\binom{2k-2}{k-1}\left(\frac{-t}{4}\right)^k$$

With $t=-4z$ we obtain from this the following formula:
$$\sqrt{1-4z}=1-2\sum_{k=1}^\infty\frac{1}{k}\binom{2k-2}{k-1}z^k$$

(7) Now back to our series $f$, we obtain the following formula for it:
\begin{eqnarray*}
f(z)
&=&\frac{1-\sqrt{1-4z}}{2z}\\
&=&\sum_{k=1}^\infty\frac{1}{k}\binom{2k-2}{k-1}z^{k-1}\\
&=&\sum_{k=0}^\infty\frac{1}{k+1}\binom{2k}{k}z^k
\end{eqnarray*}

(8) Thus the Catalan numbers are given by the formula the statement, namely:
$$C_k=\frac{1}{k+1}\binom{2k}{k}$$

So done, and note in passing that I kept my promise, from the proof of Proposition 3.5. Indeed, with the above final formula, the numerics are easily worked out.
\end{proof}

Many other things can be said about the Catalan numbers, as a continuation of the above, and about the central binomial coefficients too. We will be back to this. 

\bigskip

In relation now with Question 3.1, we are led to the following questions:

\begin{question}
What are the following centered measures?
\begin{enumerate}
\item The measure having the central binomial coefficients as even moments.

\item The measure having the Catalan numbers as even moments.
\end{enumerate}
\end{question}

We will solve in what follows this question, which is of key importance.

\section*{3c. Stieltjes inversion}

As explained above, the problem is now, how to recover a probability measure out of its moments. And the answer here, which is something non-trivial, is as follows:

\index{Stieltjes inversion}
\index{Cauchy transform}
\index{moment problem}

\begin{theorem}
The density of a real probability measure $\mu$ can be recaptured from the sequence of moments $\{M_k\}_{k\geq0}$ via the Stieltjes inversion formula
$$d\mu (x)=\lim_{t\searrow 0}-\frac{1}{\pi}\,Im\left(G(x+it)\right)\cdot dx$$
where the function on the right, given in terms of moments by
$$G(\xi)=\xi^{-1}+M_1\xi^{-2}+M_2\xi^{-3}+\ldots$$
is the Cauchy transform of the measure $\mu$.
\end{theorem}

\begin{proof}
The Cauchy transform of our measure $\mu$ is given by:
\begin{eqnarray*}
G(\xi)
&=&\xi^{-1}\sum_{k=0}^\infty M_k\xi^{-k}\\\
&=&\int_\mathbb R\frac{\xi^{-1}}{1-\xi^{-1}y}\,d\mu(y)\\
&=&\int_\mathbb R\frac{1}{\xi-y}\,d\mu(y)
\end{eqnarray*}

Now with $\xi=x+it$, we obtain the following formula:
\begin{eqnarray*}
Im(G(x+it))
&=&\int_\mathbb RIm\left(\frac{1}{x-y+it}\right)d\mu(y)\\
&=&\int_\mathbb R\frac{1}{2i}\left(\frac{1}{x-y+it}-\frac{1}{x-y-it}\right)d\mu(y)\\
&=&-\int_\mathbb R\frac{t}{(x-y)^2+t^2}\,d\mu(y)
\end{eqnarray*}

By integrating over $[a,b]$ we obtain, with the change of variables $x=y+tz$:
\begin{eqnarray*}
\int_a^bIm(G(x+it))dx
&=&-\int_\mathbb R\int_a^b\frac{t}{(x-y)^2+t^2}\,dx\,d\mu(y)\\
&=&-\int_\mathbb R\int_{(a-y)/t}^{(b-y)/t}\frac{t}{(tz)^2+t^2}\,t\,dz\,d\mu(y)\\
&=&-\int_\mathbb R\int_{(a-y)/t}^{(b-y)/t}\frac{1}{1+z^2}\,dz\,d\mu(y)\\
&=&-\int_\mathbb R\left(\arctan\frac{b-y}{t}-\arctan\frac{a-y}{t}\right)d\mu(y)
\end{eqnarray*}

Now observe that with $t\searrow0$ we have:
$$\lim_{t\searrow0}\left(\arctan\frac{b-y}{t}-\arctan\frac{a-y}{t}\right)
=\begin{cases}
\frac{\pi}{2}-\frac{\pi}{2}=0& (y<a)\\
\frac{\pi}{2}-0=\frac{\pi}{2}& (y=a)\\
\frac{\pi}{2}-(-\frac{\pi}{2})=\pi& (a<y<b)\\
0-(-\frac{\pi}{2})=\frac{\pi}{2}& (y=b)\\
-\frac{\pi}{2}-(-\frac{\pi}{2})=0& (y>b)
\end{cases}$$

We therefore obtain the following formula:
$$\lim_{t\searrow0}\int_a^bIm(G(x+it))dx=-\pi\left(\mu(a,b)+\frac{\mu(a)+\mu(b)}{2}\right)$$

Thus, we are led to the conclusion in the statement.
\end{proof}

Before getting further, let us mention that the above result does not fully solve the moment problem, because we still have the question of understanding when a sequence of numbers $M_1,M_2,M_3,\ldots$ can be the moments of a measure $\mu$.  We have here:

\index{Hankel determinant}
\index{moment problem}

\begin{theorem}
A sequence of numbers $M_0,M_1,M_2,M_3,\ldots\in\mathbb R$, with $M_0=1$, is the series of moments of a real probability measure $\mu$ precisely when:
$$\begin{vmatrix}M_0\end{vmatrix}\geq0\quad,\quad 
\begin{vmatrix}
M_0&M_1\\
M_1&M_2
\end{vmatrix}\geq0\quad,\quad 
\begin{vmatrix}
M_0&M_1&M_2\\
M_1&M_2&M_3\\
M_2&M_3&M_4\\
\end{vmatrix}\geq0\quad,\quad 
\ldots$$
That is, the associated Hankel determinants must be all positive.
\end{theorem}

\begin{proof}
This is something a bit more advanced, the idea being as follows:

\medskip

(1) As a first observation, the positivity conditions in the statement tell us that the following associated linear forms must be positive:
$$\sum_{i,j=1}^nc_i\bar{c}_jM_{i+j}\geq0$$

(2) But this is something very classical, in one sense the result being elementary, coming from the following computation, which shows that we have positivity indeed:
\begin{eqnarray*}
\int_\mathbb R\left|\sum_{i=1}^nc_ix^i\right|^2d\mu(x)
&=&\int_\mathbb R\sum_{i,j=1}^nc_i\bar{c}_jx^{i+j}d\mu(x)\\
&=&\sum_{i,j=1}^nc_i\bar{c}_jM_{i+j}
\end{eqnarray*}

(3) As for the other sense, here the result comes once again from the above formula, this time via some standard functional analysis.
\end{proof}

As a basic application of the Stieltjes formula, let us solve the moment problem for the Catalan numbers $C_k$, and for the central binomial coefficients $D_k$. We first have:

\index{semicircle law}
\index{Wigner law}

\begin{theorem}
The real measure having as even moments the Catalan numbers, $C_k=\frac{1}{k+1}\binom{2k}{k}$, and having all odd moments $0$ is the measure
$$\gamma_1=\frac{1}{2\pi}\sqrt{4-x^2}dx$$
called Wigner semicircle law on $[-2,2]$.
\end{theorem}

\begin{proof}
In order to apply the inversion formula, our starting point will be the formula from Theorem 3.8 for the generating series of the Catalan numbers, namely:
$$\sum_{k=0}^\infty C_kz^k=\frac{1-\sqrt{1-4z}}{2z}$$

By using this formula with $z=\xi^{-2}$, we obtain the following formula:
\begin{eqnarray*}
G(\xi)
&=&\xi^{-1}\sum_{k=0}^\infty C_k\xi^{-2k}\\
&=&\xi^{-1}\cdot\frac{1-\sqrt{1-4\xi^{-2}}}{2\xi^{-2}}\\
&=&\frac{\xi}{2}\left(1-\sqrt{1-4\xi^{-2}}\right)\\
&=&\frac{\xi}{2}-\frac{1}{2}\sqrt{\xi^2-4}
\end{eqnarray*}

Now let us apply Theorem 3.10. The study here goes as follows:

\medskip

(1) According to the general philosophy of the Stieltjes formula, the first term, namely $\xi/2$, which is ``trivial'', will not contribute to the density. 

\medskip

(2) As for the second term, which is something non-trivial, this will contribute to the density, the rule here being that the square root $\sqrt{\xi^2-4}$ will be replaced by the ``dual'' square root $\sqrt{4-x^2}\,dx$, and that we have to multiply everything by $-1/\pi$. 

\medskip

(3) As a conclusion, by Stieltjes inversion we obtain the following density:
$$d\mu(x)
=-\frac{1}{\pi}\cdot-\frac{1}{2}\sqrt{4-x^2}\,dx
=\frac{1}{2\pi}\sqrt{4-x^2}dx$$

Thus, we have obtained the mesure in the statement, and we are done.
\end{proof}

We have the following version of the above result:

\index{Marchenko-Pastur law}

\begin{theorem}
The real measure having as sequence of moments the Catalan numbers, $C_k=\frac{1}{k+1}\binom{2k}{k}$, is the measure
$$\pi_1=\frac{1}{2\pi}\sqrt{4x^{-1}-1}\,dx$$
called Marchenko-Pastur law on $[0,4]$.
\end{theorem}

\begin{proof}
As before, we use the standard formula for the generating series of the Catalan numbers. With $z=\xi^{-1}$ in that formula, we obtain the following formula:
\begin{eqnarray*}
G(\xi)
&=&\xi^{-1}\sum_{k=0}^\infty C_k\xi^{-k}\\
&=&\xi^{-1}\cdot\frac{1-\sqrt{1-4\xi^{-1}}}{2\xi^{-1}}\\
&=&\frac{1}{2}\left(1-\sqrt{1-4\xi^{-1}}\right)\\
&=&\frac{1}{2}-\frac{1}{2}\sqrt{1-4\xi^{-1}}
\end{eqnarray*}

With this in hand, let us apply now the Stieltjes inversion formula, from Theorem 3.10. We obtain, a bit as before in Theorem 3.12, the following density:
$$d\mu(x)
=-\frac{1}{\pi}\cdot-\frac{1}{2}\sqrt{4x^{-1}-1}\,dx
=\frac{1}{2\pi}\sqrt{4x^{-1}-1}\,dx$$

Thus, we are led to the conclusion in the statement.
\end{proof}

Regarding now the central binomial coefficients, we have here:

\index{arcsine law}

\begin{theorem}
The real probability measure having as moments the central binomial coefficients, $D_k=\binom{2k}{k}$, is the measure
$$\alpha_1=\frac{1}{\pi\sqrt{x(4-x)}}\,dx$$
called arcsine law on $[0,4]$.
\end{theorem}

\begin{proof}
We have the following computation, using some standard formulae:
\begin{eqnarray*}
G(\xi)
&=&\xi^{-1}\sum_{k=0}^\infty D_k\xi^{-k}\\
&=&\frac{1}{\xi}\sum_{k=0}^\infty D_k\left(-\frac{t}{4}\right)^k\\
&=&\frac{1}{\xi}\cdot\frac{1}{\sqrt{1-4/\xi}}\\
&=&\frac{1}{\sqrt{\xi(\xi-4)}} 
\end{eqnarray*}

But this gives the density in the statement, via Theorem 3.10. 
\end{proof}

Finally, we have the following version of the above result:

\index{modified arcsine law}
\index{middle binomial coefficients}

\begin{theorem}
The real probability measure having as moments the middle binomial coefficients, $E_k=\binom{k}{[k/2]}$, is the following law on $[-2,2]$,
$$\sigma_1=\frac{1}{2\pi}\sqrt{\frac{2+x}{2-x}}\,dx$$
called modified arcsine law on $[-2,2]$.
\end{theorem}

\begin{proof}
In terms of the central binomial coefficients $D_k$, we have:
$$E_{2k}=\binom{2k}{k}=\frac{(2k)!}{k!k!}=D_k$$
$$E_{2k-1}=\binom{2k-1}{k}=\frac{(2k-1)!}{k!(k-1)!}=\frac{D_k}{2}$$

Standard calculus based on the Taylor formula for $(1+t)^{-1/2}$ gives:
$$\frac{1}{2x}\left(\sqrt{\frac{1+2x}{1-2x}}-1\right)=\sum_{k=0}^\infty E_kx^k$$

With $x=\xi^{-1}$ we obtain the following formula for the Cauchy transform:
\begin{eqnarray*}
G(\xi)
&=&\xi^{-1}\sum_{k=0}^\infty E_k\xi^{-k}\\
&=&\frac{1}{\xi}\left(\sqrt{\frac{1+2/\xi}{1-2/\xi}}-1\right)\\
&=&\frac{1}{\xi}\left(\sqrt{\frac{\xi+2}{\xi-2}}-1\right)
\end{eqnarray*}

By Stieltjes inversion we obtain the density in the statement.
\end{proof}

All this is very nice, and we are obviously building here, as this book goes by, some solid knowledge in classical probability. We will be back to all this later.

\section*{3d. Circular measures}

With the above done, we can come back now to walks on finite graphs, that we know from the above to be related to the eigenvalues of the adjacency matrix $d\in M_N(0,1)$. But here, we are led to the following philosophical question, to start with:

\begin{question}
What are the most important finite graphs, that we should do our computations for?
\end{question}

Not an easy question, you have to agree with me, with the answer to this obviously depending on your previous experience with mathematics, or physics, or chemistry, or computer science, or other branch of science that you are interested in, and also, on the specific problems that you are the most in love with, in that part of science.

\bigskip

So, we have to be subjective here. And with me writing this book, and doing some sort of complicated quantum physics, as daytime job, I will choose the ADE graphs. It is beyond our scope here to explain where these ADE graphs exactly come from, and what they are good for, but as a piece of advertisement for them, we have:

\index{ADE graphs}

\begin{advertisement}
The ADE graphs classify the following:
\begin{enumerate}
\item Basic Lie groups and algebras.

\item Subgroups of $SU_2$ and of $SO_3$.

\item Singularities of algebraic manifolds.

\item Basic invariants of knots and links.

\item Subfactors and planar algebras of small index.

\item Subgroups of the quantum permutation group $S_4^+$.

\item Basic quantum field theories, and other physics beasts.
\end{enumerate}
\end{advertisement}

Which sounds exciting, doesn't it. So, have a look at this, and with the comment that some heavy learning work is needed, in order to understand how all this works. And with the extra comment that, in view of (7), tough physics, no one really understands how all this works. A nice introduction to all this is the paper of Jones \cite{jo6}.

\bigskip

Getting to work now, we first need to know what the ADE graphs are. The A graphs, which are the simplest, are as follows, with the distinguished vertex being denoted $\bullet$, and with $A_n$ having $n\geq2$ vertices, and $\tilde{A}_{2n}$ having $2n\geq2$ vertices:
$$A_n=\bullet-\circ-\circ\cdots\circ-\circ-\circ\hskip18mm 
A_{\infty}=\bullet-\circ-\circ-\circ\cdots\hskip7mm$$
\vskip-3mm
$$\ \ \ \ \ \ \ \tilde{A}_{2n}=
\begin{matrix}
\circ&\!\!\!\!-\circ-\circ\cdots\circ-\circ-&\!\!\!\!\circ\\
|&&\!\!\!\!|\\
\bullet&\!\!\!\!-\circ-\circ-\circ-\circ-&\!\!\!\!\circ\\
\\
\\
\end{matrix}\hskip20mm 
\tilde{A}_\infty=
\begin{matrix}
\circ&\!\!\!\!-\circ-\circ-\circ\cdots\\
|&\\
\bullet&\!\!\!\!-\circ-\circ-\circ\cdots\\
\\
\\
\end{matrix}
\hskip15mm$$
\vskip-7mm

These A graphs do not actually look that scary, because we already met all of them in the above, and as a comment on them, summarizing the situation, we have:

\begin{comment}
With the $A$ graphs we are not really lost into quantum physics, because all these graphs are quite familiar to us, as follows:
\begin{enumerate}
\item $A_n$ is the segment.

\item $A_\infty$ is the $\mathbb N$ graph.

\item $\tilde{A}_{2n}$ is the circle.

\item $\tilde{A}_\infty$ is the $\mathbb Z$ graph.
\end{enumerate}
\end{comment}

You might probably say, why not stopping here, and doing our unfinished business for the segment and the circle, with whatever new ideas that we might have. Good point, but in answer, these ideas will apply as well, with minimal changes, to the D graphs, which are as follows, with $D_n$ having $n\geq3$ vertices, and $\tilde{D}_n$ having $n+1\geq5$ vertices:
$$D_n=\bullet-\circ-\circ\dots\circ-
\begin{matrix}\ \circ\\
\ |\\
\ \circ \\
\ \\
\  \end{matrix}-\circ\hskip71mm$$
\vskip-7mm
$$\hskip7mm\tilde{D}_n=\bullet-
\begin{matrix}\circ\\
|\\
\circ\\
\ \\
\ \end{matrix}-\circ\dots\circ-
\begin{matrix}\ \circ\\
\ |\\
\ \circ \\
\ \\
\  \end{matrix}-\circ\hskip18mm$$
\vskip-7mm
$$\hskip50mm D_\infty=\bullet-
\begin{matrix}\circ\\
|\\
\circ\\
\ \\
\ \end{matrix}-\circ-\circ\cdots$$
\vskip-7mm

As mentioned above, it is beyond our scope here to explain what the ADE graphs really stand for, but as an informal comment on these latter D graphs, we have:

\begin{comment}
The D graphs are not that scary either, and they can be thought of as being certain technical versions of the A graphs.
\end{comment}

So, this is the situation, you have to trust me here, and for more on all this, check for instance the paper of Jones \cite{jo6}. In what concerns us, we will just take the above D graphs as they come, and do our loop count work for them, without questions asked.

\bigskip

As another comment, the labeling conventions for the AD graphs, while very standard, can be a bit confusing. The first graph in each series is by definition as follows:
$$A_2=\bullet-\circ\hskip13mm 
\tilde{A}_2=\begin{matrix}
\circ\\
||\\
\bullet\\
&\\
&\\
\end{matrix}\hskip13mm 
D_3=\begin{matrix}\ \circ\\
\ |\\
\ \bullet \\
\ \\
\  \end{matrix}-\circ \hskip13mm
\tilde{D}_4=\bullet-\!\!\!\!\!\begin{matrix}
\circ\hskip5mm \circ\\
\backslash\ \,\slash\\
\circ\\
&\\
&\\
\end{matrix}\!\!\!\!\!\!\!\!\!\!-\circ$$
\vskip-7mm

Finally, there are also a number of exceptional ADE graphs. First we have:
$$E_6=\bullet-\circ-
\begin{matrix}\circ\\
|\\
\circ\\
\ \\
\ \end{matrix}-
\circ-\circ\hskip71mm$$
\vskip-13mm
$$E_7=\bullet-\circ-\circ-
\begin{matrix}\circ\\
|\\
\circ\\
\ \
\\
\ \end{matrix}-
\circ-\circ\hskip18mm$$
\vskip-15mm
$$\hskip30mm E_8=\bullet-\circ-\circ-\circ-
\begin{matrix}\circ\\
|\\
\circ\\
\ \\
\ \end{matrix}-
\circ-\circ$$
\vskip-5mm

Then, we have extended versions of the above exceptional graphs, as follows:
$$\tilde{E}_6=\bullet-\circ-\begin{matrix}
\circ\\
|
\\
\circ\\
|&\\
\circ&\!\!\!\!-\ \circ\\
\ \\
\   \\
\ \\
\ \end{matrix}-\circ\hskip71mm$$
\vskip-22mm
$$\tilde{E}_7=\bullet-\circ-\circ-
\begin{matrix}\circ\\
|\\
\circ\\
\ \\
\ \end{matrix}-
\circ-\circ-\circ\hskip18mm$$
\vskip-15mm
$$\hskip30mm \tilde{E}_8=\bullet-\circ-\circ-\circ-\circ-
\begin{matrix}\circ\\
|\\
\circ\\
\ \\
\ \end{matrix}-
\circ-\circ$$
\vskip-5mm

And good news, that is all. Hard job for me to come now with a comment on these latter E graphs, along the lines of Comments 3.18 and 3.19, and here is what I have: 

\begin{comment}
The E graphs naturally complement the AD series, by capturing the combinatorics of certain ``exceptional'' phenomena in mathematics and physics.
\end{comment}

So long for difficult definitions and related informal talk, and as already mentioned in the above, for more on all this, have a look at the paper of Jones \cite{jo6}. Getting now to work, we have some new graphs, and here is the problem that we would like to solve:

\begin{problem}
How to count loops on the ADE graphs?
\end{problem}

In answer, as mentioned in Comment 3.18, we are already familiar with two of the ADE graphs, namely $A_\infty$ and $\tilde{A}_\infty$, which are respectively the graphs that we previously called $\mathbb N$ and $\mathbb Z$. So, based on our work for these graphs, where the combinatorics naturally led us into generating series, let us formulate the following definition:

\index{Poincar\'e series}

\begin{definition}
The Poincar\'e series of a rooted bipartite graph $X$ is
$$f(z)=\sum_{k=0}^\infty L_{2k}z^k$$
where $L_{2k}$ is the number of $2k$-loops based at the root.
\end{definition}

To be more precise, observe that all the above ADE graphs are indeed bipartite. Now the point is that, for a bipartite graph, the loops based at any point must have even length. Thus, in order to study the loops on the ADE graphs, based at the root, we just have to count the above numbers $L_{2k}$. And then, considering the generating series $f(z)$ of these numbers, and calling this Poincar\'e series, is something very standard.

\bigskip

Before getting into computations, let us introduce as well:

\index{positive spectral measure}

\begin{definition}
The positive spectral measure $\mu$ of a rooted bipartite graph $X$ is the real probability measure having the numbers $L_{2k}$ as moments:
$$\int_\mathbb Rx^kd\mu(x)=L_{2k}$$
Equivalently, we must have the Stieltjes transform formula
$$f(z)=\int_\mathbb R\frac{1}{1-xz}\,d\mu(x)$$
where $f$ is the Poincar\'e series of $X$.
\end{definition}

Here the existence of $\mu$, and the fact that this is indeed a positive measure, meaning a measure supported on $[0,\infty)$, comes from the following simple fact:

\begin{theorem}
The positive spectral measure of a rooted bipartite graph $X$ is given by the following formula, with $d$ being the adjacency matrix of the graph,
$$\mu=law(d^2)$$
and with the probabilistic computation being with respect to the expectation 
$$A\to<A>$$
with $<A>$ being the $(*,*)$-entry of a matrix $A$, where $*$ is the root.
\end{theorem}

\begin{proof}
With the above conventions, we have the following computation:
\begin{eqnarray*}
f(z)
&=&\sum_{k=0}^\infty L_{2k}z^k\\
&=&\sum_{k=0}^\infty\left<d^{2k}\right>z^k\\
&=&\left<\frac{1}{1-d^2z}\right>
\end{eqnarray*}

But this shows that we have $\mu=law(d^2)$, as desired.
\end{proof}

The above result shows that computing $\mu$ might be actually a simpler problem than computing $f$, and in practice, this is indeed the case. So, in what follows we will rather forget about loops and Definition 3.22, and use Definition 3.23 instead, with our computations to follow being based on the concrete interpretation from Theorem 3.24.

\bigskip

However, even with this probabilistic trick in our bag, things are not exactly trivial. So, following now \cite{bdb}, let us introduce as well the following notion:

\index{circular measure}

\begin{definition}
The circular measure $\varepsilon$ of a rooted bipartite graph $X$ is given by
$$d\varepsilon(q)=d\mu((q+q^{-1})^2)$$
where $\mu$ is the associated positive spectral measure.
\end{definition}

To be more precise, we know from Theorem 3.24 that the positive measure $\mu$ is the spectral measure of a certain positive matrix, $d^2\geq0$, and it follows from this, and from basic spectral theory, that this measure is supported by the positive reals:
$$supp(\mu)\subset\mathbb R_+$$

But then, with this observation in hand, we can define indeed the circular measure $\varepsilon$ as above, as being the pullback of $\mu$ via the following map:
$$\mathbb R\cup\mathbb T\to\mathbb R_+\quad,\quad 
q\to (q+q^{-1})^2$$

As a basic example for this, to start with, assume that $\mu$ is a discrete measure, supported by $n$ positive numbers $x_1<\ldots<x_n$, with corresponding densities $p_1,\ldots,p_n$:
$$\mu=\sum_{i=1}^n p_i\delta_{x_i}$$

For each $i\in\{1,\ldots,n\}$ the equation $(q+q^{-1})^2=x_i$ has then four solutions, that we can denote $q_i,q_i^{-1},-q_i,-q_i^{-1}$. And with this notation, we have:
$$\varepsilon=\frac{1}{4}\sum_{i=1}^np_i\left(\delta_{q_i}+\delta_{q_i^{-1}}+\delta_{-q_i}+\delta_{-q_i^{-1}}\right)$$

In general, the basic properties of $\varepsilon$ can be summarized as follows:

\begin{theorem}
The circular measure has the following properties:
\begin{enumerate}
\item $\varepsilon$ has equal density at $q,q^{-1},-q,-q^{-1}$.

\item The odd moments of $\varepsilon$ are $0$.

\item The even moments of $\varepsilon$ are half-integers.

\item When $X$ has norm $\leq 2$, $\varepsilon$ is supported by the unit circle.

\item When $X$ is finite, $\varepsilon$ is discrete.

\item If $K$ is a solution of $d=K+K^{-1}$, then $\varepsilon=law(K)$. 
\end{enumerate}
\end{theorem}

\begin{proof}
These results can be deduced from definitions, the idea being that (1-5) are trivial, and that (6) follows from the formula of $\mu$ from Theorem 3.24.
\end{proof}

Getting now to computations, remember our struggle from the above, with the circle graph? We can now solve this question, majestically, as follows:

\index{circle graph}

\begin{theorem}
The circular measure of the basic index $4$ graph, namely 
$$\begin{matrix}
&\circ&\!\!\!\!-\circ-\circ\cdots\circ-\circ-&\!\!\!\!\circ\cr
\tilde{A}_{2n}=&|&&\!\!\!\!|\cr
&\bullet&\!\!\!\!-\circ-\circ-\circ-\circ-&\!\!\!\!\circ\cr\cr\cr\end{matrix}$$
\vskip-7mm

\noindent is the uniform measure on the $2n$-roots of unity.
\end{theorem}

\begin{proof}
Let us identify the vertices of $X=\tilde{A}_{2n}$ with the group $\{w^k\}$ formed by the $2n$-th roots of unity in the complex plane, where $w=e^{\pi i/n}$. The adjacency matrix of $X$ acts then on the functions $f\in C(X)$ in the following way:
$$df(w^s)=f(w^{s-1})+f(w^{s+1})$$

But this shows that we have $d=K+K^{-1}$, where $K$ is given by:
$$Kf(w^s)=f(w^{s+1})$$

Thus we can use Theorem 3.24 and Theorem 3.26 (6), and we get:
$$\varepsilon=law(K)$$

But this is the uniform measure on the $2n$-roots of unity, as claimed.
\end{proof}

All this is very nice, so, before going ahead with more computations, let us have an excursion into subfactor theory, and explain what is behind this trick. Following Jones \cite{jo7}, we can introduce the theta series of a graph $X$, as a version of the Poincar\'e series, via the change of variables $z^{-1/2}=q^{1/2}+q^{-1/2}$, as follows:

\index{theta series}
\index{Jones series}

\begin{definition}
The theta series of a rooted bipartite graph $X$ is
$$\Theta(q)=q+\frac{1-q}{1+q}f\left(\frac{q}{(1+q)^2}\right)$$
where $f$ is the Poincar\'e series.
\end{definition}

The theta series can be written as $\Theta(q)=\sum a_rq^r$, and it follows from the above formula, via some simple manipulations, that its coefficients are integers:
$$a_r\in\mathbb Z$$

In fact, we have the following explicit formula from Jones' paper \cite{jo7}, relating the coefficients of $\Theta(q)=\sum a_rq^r$ to those of the Poincar\'e series $f(z)=\sum c_kz^k$:
$$a_r=\sum_{k=0}^r(-1)^{r-k}\frac{2r}{r+k}\begin{pmatrix}r+k\cr r-k\end{pmatrix}c_k$$

As an important comment now, in the case where $X$ is the principal graph of a subfactor $A_0\subset A_1$ of index $N>4$, it is known from \cite{jo7} that the numbers $a_r$ are certain multiplicities associated to the planar algebra inclusion $TL_N\subset P$, as explained there. In particular, the coefficients of the theta series are in this case positive integers:
$$a_r\in\mathbb N$$

In relation now with the circular measure, the result here, which is quite similar to the Stieltjes transform formula from Definition 3.23, is as follows:

\begin{theorem}
We have the Stieltjes transform type formula
$$2\int\frac{1}{1-qu^2}\,d\varepsilon(u)=1+T(q)(1-q)$$
where the $T$ series of a rooted bipartite graph $X$ is by definition given by
$$T(q)=\frac{\Theta(q)-q}{1-q}$$
with $\Theta$ being the associated theta series.
\end{theorem}

\begin{proof}
This follows by applying the change of variables $q\to (q+q^{-1})^2$ to the fact that $f$ is the Stieltjes transform of $\mu$. Indeed, we obtain in this way:
\begin{eqnarray*}
2\int\frac{1}{1-qu^2}\,d\varepsilon(u)
&=&1+\frac{1-q}{1+q}f\left(\frac{q}{(1+q)^2}\right)\\
&=&1+\Theta(q)-q\\
&=&1+T(q)(1-q)
\end{eqnarray*}

Thus, we are led to the conclusion in the statement.
\end{proof}

Summarizing, we have a whole menagery of subfactor, planar algebra and bipartite graph invariants, which come in several flavors, namely series and measures, and which can be linear or circular, and which all appear as versions of the Poincar\'e series.

\bigskip

In order to discuss all this more systematically, let us introduce as well:

\index{cyclotomic series}

\begin{definition}
The series of the form
$$\xi(n_1,\ldots,n_s:m_1,\ldots,m_t)=\frac{(1-q^{n_1})\ldots(1-q^{n_s})}{(1-q^{m_1})\ldots(1-q^{m_t})}$$
with $n_i,m_i\in\mathbb N$ are called cyclotomic.
\end{definition}

It is technically convenient to allow as well $1+q^n$ factors, to be designated by $n^+$ symbols in the above writing. For instance we have, by definition:
$$\xi(2^+:3)=\xi(4:2,3)$$

Also, it is convenient in what follows to use the following notations:
$$\xi'=\frac{\xi}{1-q}\quad,\quad \xi''=\frac{\xi}{1-q^2}$$

The Poincar\'e series of the ADE graphs are given by quite complicated formulae. However, the corresponding $T$ series are all cyclotomic, as follows:

\index{T series}

\begin{theorem}
The $T$ series of the ADE graphs are as follows:
\begin{enumerate}
\item For $A_{n-1}$ we have $T=\xi(n-1:n)$.

\item For $D_{n+1}$ we have $T=\xi(n-1^+:n^+)$.

\item For $\tilde{A}_{2n}$ we have $T=\xi'(n^+:n)$.

\item For $\tilde{D}_{n+2}$ we have $T=\xi''(n+1^+:n)$.

\item For $E_6$ we have $T=\xi(8:3,6^+)$.

\item For $E_7$ we have $T=\xi(12:4,9^+)$.

\item For $E_8$ we have $T=\xi(5^+,9^+:15^+)$.

\item For $\tilde{E}_6$ we have $T=\xi(6^+:3,4)$.

\item For $\tilde{E}_7$ we have $T=\xi(9^+:4,6)$.

\item For $\tilde{E}_8$ we have $T=\xi(15^+:6,10)$.
\end{enumerate}
\end{theorem}

\begin{proof}
These formulae were obtained in \cite{bdb}, by counting loops, and then by making the following change of variables, and factorizing the resulting series:
$$z^{-1/2}=q^{1/2}+q^{-1/2}$$

An alternative proof for these formulae can be obtained by using planar algebra methods, along the lines of the paper of Jones \cite{jo7}. For details here, see \cite{bdb}.
\end{proof}

In order to formulate our final results, we will need more theory. First, we have:

\index{cyclotomic measure}

\begin{definition}
A cyclotomic measure is a probability measure $\varepsilon$ on the unit circle, having the following properties:
\begin{enumerate}
\item  $\varepsilon$ is supported by the $2n$-roots of unity, for some $n\in\mathbb N$.

\item $\varepsilon$ has equal density at $q,q^{-1},-q,-q^{-1}$.
\end{enumerate}
\end{definition}

As a first observation, it follows from Theorem 3.26 and from Theorem 3.31 that the circular measures of the finite ADE graphs are supported by certain roots of unity, hence are cyclotomic. We will be back to this in a moment, with details, and computations.

\bigskip

At the general level now, let us introduce as well the following notion:

\begin{definition}
The $T$ series of a cyclotomic measure $\varepsilon$ is given by
$$1+T(q)(1-q)=2\int\frac{1}{1-qu^2}\,d\varepsilon(u)$$
with $\varepsilon$ being as usual the circular spectral measure.
\end{definition}

Observe that this formula is nothing but the one in Theorem 3.29, written now in the other sense. In other words, if the cyclotomic measure $\varepsilon$ happens to be the circular measure of a rooted bipartite graph, then the $T$ series as defined above coincides with the $T$ series as defined before. This is useful for explicit computations.

\bigskip

Good news, with the above technology in hand, and with a computation already done, in Theorem 3.27, we are now ready to discuss the circular measures of all ADE graphs. The idea will be that these measures are all cyclotomic, of level $\leq 3$, and can be expressed in terms of the basic polynomial densities of degree $\leq 6$, namely:
$$\alpha=Re(1-q^2)$$
$$\beta=Re(1-q^4)$$
$$\gamma=Re(1-q^6)$$

To be more precise, we have the following final result on the subject, with $\alpha,\beta,\gamma$ being as above, with $d_n$ being the uniform measure on the $2n$-th roots of unity, and with $d_n'=2d_{2n}-d_n$ being the uniform measure on the odd $4n$-roots of unity:

\index{ADE graph}
\index{circular measure}

\begin{theorem}
The circular measures of the ADE graphs are given by:
\begin{enumerate}
\item $A_{n-1}\to\alpha_n$.

\item $\tilde{A}_{2n}\to d_n$.

\item $D_{n+1}\to\alpha_n'$.

\item $\tilde{D}_{n+2}\to (d_n+d_1')/2$.

\item $E_6\to\alpha_{12}+(d_{12}-d_6-d_4+d_3)/2$.

\item $E_7\to\beta_9'+(d_1'-d_3')/2$.

\item $E_8\to\alpha_{15}'+\gamma_{15}'-(d_5'+d_3')/2$.

\item $\tilde{E}_{n+3}\to (d_n+d_3+d_2-d_1)/2$.
\end{enumerate}
\end{theorem}

\begin{proof}
This is something which can be proved in three steps, as follows:

\medskip

(1) For the simplest graph, namely the circle $\tilde{A}_{2n}$, we already have the result, from Theorem 3.27, with the proof there being something elementary.

\medskip

(2) For the other non-exceptional graphs, that is, of type A and D, the same method works, namely direct loop counting, with some matrix tricks. See \cite{bdb}.

\medskip

(3) In general, this follows from the $T$ series formulae in Theorem 3.31, via some manipulations based on the general conversion formulae given above. See \cite{bdb}.
\end{proof}

We refer to \cite{bdb} and the subsequent literature for more on all this. Also, let us point out that all this leads to a more conceptual understanding of what we did before, for the graphs $\mathbb N$ and $\mathbb Z$. Indeed, even for these very basic graphs, using the unit circle and circular measures as above leads to a better understanding of the combinatorics.

\section*{3e. Exercises}

We had yet another tough chapter here, again mixing linear algebra with calculus, and with probability. Here are some exercises, in relation with the above:

\begin{exercise}
Work out all details for the bijections leading to Catalan numbers.
\end{exercise}

\begin{exercise}
Find some further interpretations of the Catalan numbers.
\end{exercise}

\begin{exercise}
Clarify all the details, in relation with our Dyck path counting.
\end{exercise}

\begin{exercise}
Learn more about Stieltjes inversion, and the moment problem.
\end{exercise}

\begin{exercise}
Learn more, from physicists, about the Wigner semicircle law.
\end{exercise}

\begin{exercise}
Learn as well, from probabilists, about Marchenko-Pastur.
\end{exercise}

\begin{exercise}
Do the missing computations for the A and D graphs.
\end{exercise}

\begin{exercise}
Compute the spectral measures of exceptional graphs.
\end{exercise}

As bonus exercise, learn some subfactor theory. More generally, learn, or at least get to know about, all mathematical theories featuring an ADE classification result.

\chapter{Transitive graphs}

\section*{4a. Circulant graphs}

You might have noticed already, some graphs look good, and some other look bad. This is of course a matter of taste, and there are several possible notions for what ``good'' should mean, and with this being a potential source of serious mathematical matter. As an example, we have already met a particularly beautiful graph before, namely:
$$
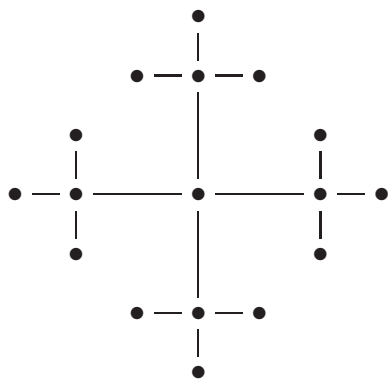
\xymatrix@R=10pt@C=10pt{
&&&\bullet\ar@{-}[d]\\
&&\bullet\ar@{-}[r]&\bullet\ar@{-}[dd]\ar@{-}[r]&\bullet\\
&\bullet\ar@{-}[d]&&&&\bullet\ar@{-}[d]\\
\bullet\ar@{-}[r]&\bullet\ar@{-}[rr]&&\bullet\ar@{-}[rr]\ar@{-}[dd]&&\bullet\ar@{-}[r]&\bullet\\
&\bullet\ar@{-}[u]&&&&\bullet\ar@{-}[u]\\
&&\bullet\ar@{-}[r]&\bullet\ar@{-}[d]\ar@{-}[r]&\bullet\\
&&&\bullet
}$$

Such a graph is called a tree, because when looking at a tree from the above, what you see is something like this. In general, trees can be axiomatized as follows:

\index{tree}
\index{no loops}

\begin{definition}
A tree is a graph having no cycles. That is, there is no loop
$$
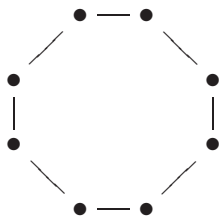
\xymatrix@R=12pt@C=12pt{
&\bullet\ar@{-}[r]\ar@{-}[dl]&\bullet\ar@{-}[dr]\\
\bullet\ar@{-}[d]&&&\bullet\ar@{-}[d]\\
\bullet\ar@{-}[dr]&&&\bullet\ar@{-}[dl]\\
&\bullet\ar@{-}[r]&\bullet}$$
having length $\geq3$, and distinct vertices, inside the graph.
\end{definition}

And aren't trees beautiful, hope you agree with me. But it is not about trees that we want to talk about here, these are in fact quite complicated mathematical objects, and we will keep them for later. As an alternative to them, we have the circulant graphs, which are equally beautiful, but in a somewhat opposite sense. Here is one, which is actually the most important graph in the history of mankind, along with the fire graph:
$$
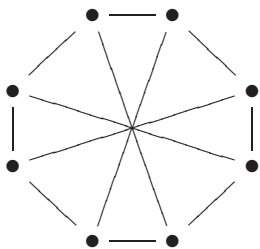
\xymatrix@R=16pt@C=17pt{
&\bullet\ar@{-}[r]\ar@{-}[dddr]\ar@{-}[dl]&\bullet\ar@{-}[dr]\ar@{-}[dddl]\\
\bullet\ar@{-}[d]\ar@{-}[drrr]&&&\bullet\ar@{-}[d]\\
\bullet\ar@{-}[dr]\ar@{-}[urrr]&&&\bullet\ar@{-}[dl]\\
&\bullet\ar@{-}[r]&\bullet}$$

Here is another circulant graph, again with 8 vertices, again with the picture suggesting the name ``circulant'', and of course, again beautiful as well:
$$
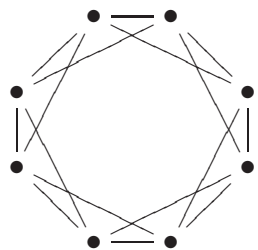
\xymatrix@R=16pt@C=16pt{
&\bullet\ar@{-}[r]\ar@{-}[drr]\ar@{-}[dl]&\bullet\ar@{-}[dr]\ar@{-}[dll]\\
\bullet\ar@{-}[d]\ar@{-}[ddr]&&&\bullet\ar@{-}[d]\\
\bullet\ar@{-}[dr]\ar@{-}[uur]\ar@{-}[drr]&&&\bullet\ar@{-}[dl]\ar@{-}[uul]\\
&\bullet\ar@{-}[r]\ar@{-}[urr]&\bullet\ar@{-}[uur]}$$

You get the point with these graphs, we are trying here to get in a sense which is opposite to that of Definition 4.1, with the cycles being not only welcome, but somehow mandatory. In general, the circulant graphs can be axiomatized as follows:

\index{circulant graph}
\index{wheel graph}

\begin{definition}
A graph is called circulant if, when drawn with its vertices on a circle, equally spaced, it is invariant under rotations.
\end{definition}

To be more precise, this is the definition of the circulant graphs, when the vertices are labeled in advance $1,\ldots,N$. In general, when the vertices are not labeled in advance, the convention is that the graph is called circulant when it is possible to label the vertices  $1,\ldots,N$, as for the graph to become circulant in the above sense.

\bigskip

In order to understand this, let us pick a graph which is obviously circulant, such as the wheel graph above, and mess up the labeling of the vertices, see what we get. For this purpose, let us put some random labels $1,2,\dots,8$ on our wheel graph, say as follows:
$$
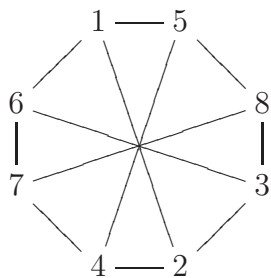
\xymatrix@R=16pt@C=18pt{
&1\ar@{-}[r]\ar@{-}[dddr]\ar@{-}[dl]&5\ar@{-}[dr]\ar@{-}[dddl]\\
6\ar@{-}[d]\ar@{-}[drrr]&&&8\ar@{-}[d]\\
7\ar@{-}[dr]\ar@{-}[urrr]&&&3\ar@{-}[dl]\\
&4\ar@{-}[r]&2}$$

Now let us redraw this graph, with the vertices $1,2,\ldots,8$ ordered on a circle, equally spaced, as Definition 4.2 requires. We get something not very beautiful, as follows:
$$
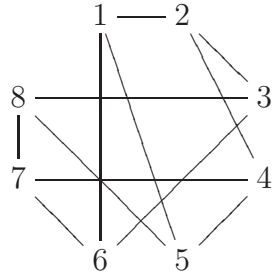
\xymatrix@R=16pt@C=18pt{
&1\ar@{-}[r]\ar@{-}[dddr]\ar@{-}[ddd]&2\ar@{-}[dr]\ar@{-}[ddr]\\
8\ar@{-}[d]&&&3\ar@{-}[lll]\ar@{-}[ddll]\\
7\ar@{-}[dr]&&&4\ar@{-}[dl]\ar@{-}[lll]\\
&6&5\ar@{-}[uull]}$$

So, here is the point. This graph, regarded as a graph with vertices labeled $1,2,\ldots,8$ is obviously not circulant, in the sense of Definition 4.2. However, when removing the labels, this graph does become circulant, as per our conventions above. 

\bigskip

All this might seem a bit confusing, when first seen, and you are probably in this situation, so here is a precise statement in this sense, coming with a full proof:

\begin{proposition}
The following graph is circulant, despite its bad look,
$$
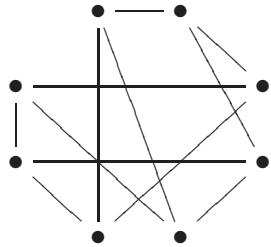
\xymatrix@R=16pt@C=18pt{
&\bullet\ar@{-}[r]\ar@{-}[dddr]\ar@{-}[ddd]&\bullet\ar@{-}[dr]\ar@{-}[ddr]\\
\bullet\ar@{-}[d]&&&\bullet\ar@{-}[lll]\ar@{-}[ddll]\\
\bullet\ar@{-}[dr]&&&\bullet\ar@{-}[dl]\ar@{-}[lll]\\
&\bullet&\bullet\ar@{-}[uull]}$$
in the sense that it can be put in circulant form, with a suitable labeling of the vertices.
\end{proposition}

\begin{proof}
As already mentioned, this normally follows from the above discussion, but let us prove this as well directly. The idea is as follows:

\medskip

(1) Let us label the vertices of our graph as follows, and I will explain in moment where this tricky labeling choice comes from:
$$
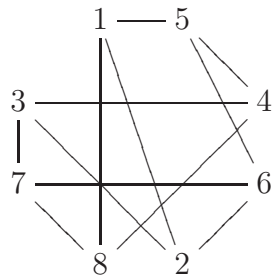
\xymatrix@R=16pt@C=18pt{
&1\ar@{-}[r]\ar@{-}[dddr]\ar@{-}[ddd]&5\ar@{-}[dr]\ar@{-}[ddr]\\
3\ar@{-}[d]&&&4\ar@{-}[lll]\ar@{-}[ddll]\\
7\ar@{-}[dr]&&&6\ar@{-}[dl]\ar@{-}[lll]\\
&8&2\ar@{-}[uull]}$$

Now let us redraw this graph, with the vertices $1,2,\ldots,8$ ordered on a circle, equally spaced, as Definition 4.2 requires. We get something very nice, as follows:
$$
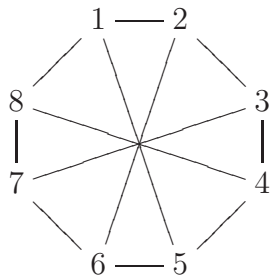
\xymatrix@R=16pt@C=18pt{
&1\ar@{-}[r]\ar@{-}[dddr]\ar@{-}[dl]&2\ar@{-}[dr]\ar@{-}[dddl]\\
8\ar@{-}[d]\ar@{-}[drrr]&&&3\ar@{-}[d]\\
7\ar@{-}[dr]\ar@{-}[urrr]&&&4\ar@{-}[dl]\\
&6\ar@{-}[r]&5}$$

Thus, our original graph was indeed circulant, as stated.

\medskip

(2) In order for everything to be fully clarified, we still must explain where the tricky labeling choice in (1) comes from. For this purpose, let us recall where the graph in the statement came from. Well, this graph was obtained by messing up the labeling of the vertices of the wheel graph, by using the following permutation: 
$$\sigma=\begin{pmatrix}
1&2&3&4&5&6&7&8\\
1&5&8&3&2&4&7&6
\end{pmatrix}$$

The point now is that, if we want to unmess our graph, we must use the inverse of the above permutation, obtained by reading things upside-down, which is given by:
$$\sigma^{-1}=\begin{pmatrix}
1&2&3&4&5&6&7&8\\
1&5&4&6&2&8&7&3
\end{pmatrix}$$

Thus, we must label our vertices $1,5,4,6,2,8,7,3$, precisely as done in (1).
\end{proof}

As a conclusion to all this, deciding whether a graph is circulant or not is not an easy business. In what follows we will rather focus on the graphs which come by definition in circulant form, and leave decision problems for arbitrary graphs for later.

\bigskip

As basic examples now of circulant graphs, we have the triangle, the square, the pentagon, the hexagon and so on. As usual, let us record the formula of the adjacency matrix for such graphs. In the general case of the $N$-gon, this is as follows:

\index{regular polygon}
\index{N-gon}

\begin{proposition}
The adjacency matrix of the $N$-gon is
$$d_{ij}=\delta_{|i-j|,1}$$
with the indices taken modulo $N$.
\end{proposition}

\begin{proof}
This is clear indeed from definitions, because with the indices taken modulo $N$, and arranged on a circle, being neighbor means $i=j\pm1$, and so $|i-j|=1$.
\end{proof}

Now observe that the adjacency matrix computed above best looks when written in usual matrix form, with its circulant nature being quite obvious. As an illustration for this, here is the adjacency matrix of the hexagon, which is obviously circulant:
$$d=\begin{pmatrix}
0&1&0&0&0&1\\
1&0&1&0&0&0\\
0&1&0&1&0&0\\
0&0&1&0&1&0\\
0&0&0&1&0&1\\
1&0&0&0&1&0
\end{pmatrix}$$

Passed the $N$-gons, we can have more examples of circulant graphs by adding spokes to the $N$-gons. Here is a hexagonal wheel, that we already met in chapter 1:
$$d=\begin{pmatrix}
0&1&0&1&0&1\\
1&0&1&0&1&0\\
0&1&0&1&0&1\\
1&0&1&0&1&0\\
0&1&0&1&0&1\\
1&0&1&0&1&0
\end{pmatrix}$$

Of course, not all circulant graphs appear in this way. As an example, here is another sort of ``hexagonal wheel'', without spokes, namely the Star of David:
$$d=\begin{pmatrix}
0&0&1&0&1&0\\
0&0&0&1&0&1\\
1&0&0&0&1&0\\
0&1&0&0&0&1\\
1&0&1&0&0&0\\
0&1&0&1&0&0
\end{pmatrix}$$

But probably enough examples, you got the point, and time now to do some theory. Inspired by the above examples of explicit adjacency matrices, we have:

\index{circulant matrix}

\begin{proposition}
A graph is circulant precisely when its adjacency matrix is circulant, in the sense that it is of the following form, with the indices taken modulo $N$:
$$d_{ij}=\gamma_{j-i}$$
In practice, and with indices $0,1,\ldots,N-1$, taken modulo $N$, this means that $d$ consists of a row vector $\gamma$, sliding downwards and to the right, in the obvious way.
\end{proposition}

\begin{proof}
This is something quite obvious. Indeed, at $N=4$ for instance, our assumption that $X$ is circulant means that the adjacency matrix must look as follows:
$$d=\begin{pmatrix}
x&y&z&t\\
t&x&y&z\\
z&t&x&y\\
y&z&t&x
\end{pmatrix}$$

Now let us call $\gamma$ the first row vector, $(x,y,z,t)$. With matrix indices $0,1,2,3$, taken modulo 4, as indicated in the statement, our matrix is then given by:
$$d=\begin{pmatrix}
\gamma_0&\gamma_1&\gamma_2&\gamma_3\\
\gamma_3&\gamma_0&\gamma_1&\gamma_2\\
\gamma_2&\gamma_3&\gamma_0&\gamma_1\\
\gamma_1&\gamma_2&\gamma_3&\gamma_0
\end{pmatrix}$$

In the general case, $N\in\mathbb N$, the situation is similar, and this leads to the result.
\end{proof}

Observe that the above result is not the end of the story, because we still have to discuss when a circulant matrix is symmetric, and has 0 on the diagonal. But this is something easy to do, and our final result is then as follows:

\begin{theorem}
The adjacency matrices of the circulant graphs are precisely the matrices of the following form, with indices $0,1,\ldots,N-1$, taken modulo $N$, 
$$d_{ij}=\gamma_{j-i}$$
with the vector $\gamma\in(0,1)^N$ being symmetric, $\gamma_i=\gamma_{-i}$, and with $\gamma_0=0$.
\end{theorem}

\begin{proof}
This follows from our observations in Proposition 4.5, because:

\medskip

(1) We know from there that we have $d_{ij}=\gamma_{j-i}$.

\medskip

(2) The symmetry of $d$ translates into the condition $\gamma_i=\gamma_{-i}$.

\medskip

(3) The fact that $d$ has 0 on the diagonal translates into the condition $\gamma_0=0$.
\end{proof}

All this is very nice, and with the circulant graphs being the same as the vectors $\gamma\in(0,1)^N$, subject to the simple assumptions above, this suggests that this might be the end of the story. Error. Recall that in chapter 2 we studied the $N$-simplex from a spectral point of view, and got into non-trivial things. But this $N$-simplex is, perhaps in rivalry with the $N$-gon, the simplest example of a circulant graph. 

\bigskip

So, let us first recall our findings from chapter 2. We have seen there that the adjacency matrix of the simplex, with a copy of the identity added for simplifying things, diagonalizes as follows, with $F_N=(w^{ij})_{ij}$ with $w=e^{2\pi i/N}$ being the Fourier matrix:
$$\begin{pmatrix}
1&\ldots&\ldots&1\\
\vdots&&&\vdots\\
\vdots&&&\vdots\\
1&\ldots&\ldots&1\end{pmatrix}=\frac{1}{N}\,F_N
\begin{pmatrix}
N\\
&0\\
&&\ddots\\
&&&0\end{pmatrix}F_N^*$$

More generally now, going towards the case of the general circulant graphs, we have the following result, which is something standard in discrete Fourier analysis:

\index{Fourier-diagonal}
\index{discrete Fourier transform}

\begin{theorem}
For a matrix $A\in M_N(\mathbb C)$, the following are equivalent,
\begin{enumerate}
\item $A$ is circulant, $A_{ij}=\xi_{j-i}$, for a certain vector $\xi\in\mathbb C^N$,

\item $A$ is Fourier-diagonal, $A=F_NQF_N^*$, for a certain diagonal matrix $Q$,
\end{enumerate}
and if so, $\xi=F_N^*q$, where $q\in\mathbb C^N$ is the vector formed by the diagonal entries of $Q$.
\end{theorem}

\begin{proof}
This follows from some basic computations with roots of unity, as follows:

\medskip

$(1)\implies(2)$ Assuming $A_{ij}=\xi_{j-i}$, the matrix $Q=F_N^*AF_N$ is indeed diagonal, as shown by the following computation:
\begin{eqnarray*}
Q_{ij}
&=&\sum_{kl}w^{-ik}A_{kl}w^{lj}\\
&=&\sum_{kl}w^{jl-ik}\xi_{l-k}\\
&=&\sum_{kr}w^{j(k+r)-ik}\xi_r\\
&=&\sum_rw^{jr}\xi_r\sum_kw^{(j-i)k}\\
&=&N\delta_{ij}\sum_rw^{jr}\xi_r
\end{eqnarray*}

$(2)\implies(1)$ Assuming $Q=diag(q_1,\ldots,q_N)$, the matrix $A=F_NQF_N^*$ is indeed circulant, as shown by the following computation:
$$A_{ij}
=\sum_kw^{ik}Q_{kk}w^{-jk}
=\sum_kw^{(i-j)k}q_k$$

To be more precise, in this formula the last term depends only on $j-i$, and so shows that we have $A_{ij}=\xi_{j-i}$, with $\xi$ being the following vector:
$$\xi_i
=\sum_kw^{-ik}q_k
=(F_N^*q)_i$$

Thus, we are led to the conclusions in the statement.
\end{proof}

The above result is something quite powerful, and useful, and suggests doing everything in Fourier, when dealing with the circulant matrices. And we can use here:

\index{discrete Fourier transform}

\begin{theorem}
The various basic sets of $N\times N$ circulant matrices are as follows, with the convention that associated to any $q\in\mathbb C^N$ is the matrix  $Q=diag(q_1,\ldots,q_N)$:
\begin{enumerate}
\item The set of all circulant matrices is:
$$M_N(\mathbb C)^{circ}=\left\{F_NQF_N^*\Big|q\in\mathbb C^N\right\}$$

\item The set of all circulant unitary matrices is:
$$U_N^{circ}=\left\{\frac{1}{N}F_NQF_N^*\Big|q\in\mathbb T^N\right\}$$

\item The set of all circulant orthogonal matrices is:
$$O_N^{circ}=\left\{\frac{1}{N}F_NQF_N^*\Big|q\in\mathbb T^N,\bar{q}_i=q_{-i},\forall i\right\}$$
\end{enumerate}
In addition, in this picture, the first row vector of $F_NQF_N^*$ is given by $\xi=F_N^*q$.
\end{theorem}

\begin{proof}
All this follows from Theorem 4.7, as follows:

\medskip

(1) This assertion, along with the last one, is Theorem 4.7 itself.

\medskip

(2) This is clear from (1), and from the fact that the rescaled matrix $F_N/\sqrt{N}$ is unitary, because the eigenvalues of a unitary matrix must be on the unit circle $\mathbb T$.

\medskip

(3) This follows from (2), because the matrix is real when $\xi_i=\bar{\xi}_i$, and in Fourier transform, $\xi=F_N^*q$, this corresponds to the condition $\bar{q}_i=q_{-i}$.
\end{proof}

There are many other things that can be said about the circulant matrices, and all this is quite interesting, in relation with the circulant graphs. We will be back to this.

\section*{4b. Transitive graphs}

We have seen so far that the circulant graphs are quite interesting objects. Our purpose now is to extend the theory that we have for them, to a more general class of graphs, the ``transitive'' ones. Let us start with a basic observation, as follows:

\begin{proposition}
A graph $X$, with vertices labeled $1,2,\ldots,N$, is circulant precisely when for any two vertices $i,j\in X$ there is a permutation $\sigma\in S_N$ such that:
\begin{enumerate}
\item $\sigma$ maps one vertex to another, $\sigma(i)=j$.

\item $\sigma$ is cyclic, $\sigma(k)=k+s$ modulo $N$, for some $s$.

\item $\sigma$ leaves invariant the edges, $k-l\iff\sigma(k)-\sigma(l)$.
\end{enumerate}
\end{proposition}

\begin{proof}
This is obvious from definitions, and with the remark that the number $s$ appearing in (2) is uniquely determined by (1), as being $s=j-i$, modulo $N$.
\end{proof}

The point now is that, with this picture of the circulant graphs in mind, it is quite clear that if we remove the assumption (2), that our permutation is cyclic, we will reach to a quite interesting class of graphs, generalizing them. So, let us formulate: 

\index{transitive graph}

\begin{definition}
A graph $X$, with vertices labeled $1,2,\ldots,N$, is called transitive when for any two vertices $i,j\in X$ there is a permutation $\sigma\in S_N$ such that:
\begin{enumerate}
\item $\sigma$ maps one vertex to another, $\sigma(i)=j$.

\item $\sigma$ leaves invariant the edges, $k-l\iff\sigma(k)-\sigma(l)$.
\end{enumerate}
\end{definition}

In short, what we did here is to copy the statement of Proposition 4.9, with the assumption (2) there removed, and call this a Definition. In view of this, obviously, any circulant graph is transitive. But, do we have other interesting examples?

\bigskip

As a first piece of answer to this question, which is very encouraging, we have:

\begin{theorem}
The cube graph, namely
$$
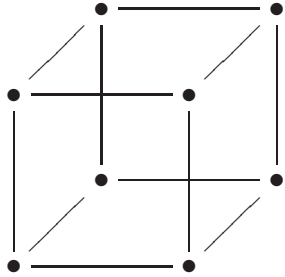
\xymatrix@R=20pt@C=20pt{
&\bullet\ar@{-}[rr]&&\bullet\\
\bullet\ar@{-}[rr]\ar@{-}[ur]&&\bullet\ar@{-}[ur]\\
&\bullet\ar@{-}[rr]\ar@{-}[uu]&&\bullet\ar@{-}[uu]\\
\bullet\ar@{-}[uu]\ar@{-}[ur]\ar@{-}[rr]&&\bullet\ar@{-}[uu]\ar@{-}[ur]
}$$
is transitive, but not circulant.
\end{theorem}

\begin{proof}
The fact that the cube is transitive is clear, because given any two vertices $i,j\in X$, we can certainly rotate the cube in 3D, as to have $i\to j$. As for the fact that the cube is not circulant, this is something more tricky, as follows:

\medskip

(1) As a first observation, when trying to draw the cube on a circle, in a somewhat nice and intuitive way, as to have it circulant, we reach to the following picture:
$$
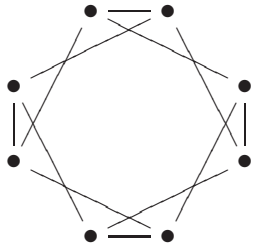
\xymatrix@R=16pt@C=16pt{
&\bullet\ar@{-}[r]\ar@{-}[drr]&\bullet\ar@{-}[dll]\\
\bullet\ar@{-}[d]\ar@{-}[ddr]&&&\bullet\ar@{-}[d]\\
\bullet\ar@{-}[uur]\ar@{-}[drr]&&&\bullet\ar@{-}[uul]\\
&\bullet\ar@{-}[r]\ar@{-}[urr]&\bullet\ar@{-}[uur]}$$

Thus, our cube is indeed not circulant, or at least not in an obvious way. 

\medskip

(2) However, this does not stand for a proof, and the problem of abstractly proving that the cube is not circulant remains. Normally this can be done by attempting to label the vertices in a circulant way. Indeed, up to some discussion here, that we will leave as an instructive exercise, we can always assume that $1,2$ are connected by an edge:
$$1-2$$

(3) But with this in hand, we can now start labeling the vertices of the cube, in a circulant way. Since $1-2$ implies via our circulant graph assumption $2-3$, $3-4$, and so on, in order to start our labeling, we must pick one vertex, and then follow a path on the cube, emanating from there. But, by some obvious symmetry reasons, this means that we can always assume that our first three vertices $1,2,3$ are as follows:
$$
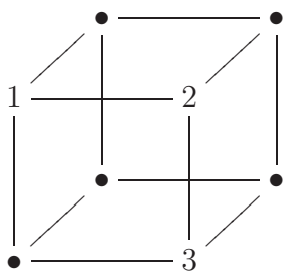
\xymatrix@R=17pt@C=20pt{
&\bullet\ar@{-}[rr]&&\bullet\\
1\ar@{-}[rr]\ar@{-}[ur]&&2\ar@{-}[ur]\\
&\bullet\ar@{-}[rr]\ar@{-}[uu]&&\bullet\ar@{-}[uu]\\
\bullet\ar@{-}[uu]\ar@{-}[ur]\ar@{-}[rr]&&3\ar@{-}[uu]\ar@{-}[ur]
}$$

(4) So, the question comes now, where the vertex 4 can be, as for all this to lead, in the end, to a circulant graph. And the point is that, among the two possible choices for the vertex 4, as new neighbors of 3, none works. Thus, our cube is indeed not circulant, and we will leave the remaining details here as an instructive exercise.
\end{proof}

Let us look now for some more examples. Thinking a bit, it is in fact not necessary to go up to the cube, which is a rather advanced object, in order to have an example of a transitive, non-circulant graph, and this because two triangles will do too:
$$
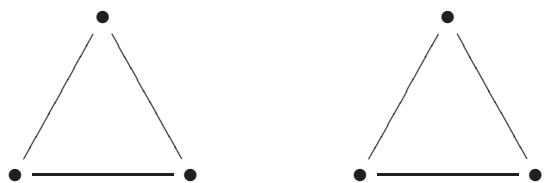
\xymatrix@R=20pt@C=20pt{
&\bullet\ar@{-}[ddl]\ar@{-}[ddr]&&&&\bullet\ar@{-}[ddl]\ar@{-}[ddr]\\
\\
\bullet\ar@{-}[rr]&&\bullet&\ &\bullet\ar@{-}[rr]&&\bullet
}$$

However, this latter graph is not connected, and so is not very good, as per our usual geometric philosophy. But, we can make it connected, by adding edges, as follows:
$$
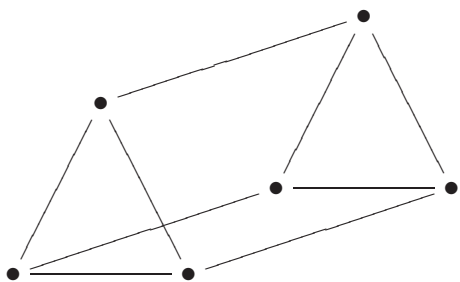
\xymatrix@R=20pt@C=20pt{
&&&&\bullet\ar@{-}[ddl]\ar@{-}[ddr]\\
&\bullet\ar@{-}[ddl]\ar@{-}[ddr]\ar@{-}[urrr]\\
&&&\bullet\ar@{-}[rr]&&\bullet\\
\bullet\ar@{-}[rr]\ar@{-}[urrr]&&\bullet\ar@{-}[urrr]
}$$

What we got here is a prism, and it is convenient now, for aesthetical and typographical reasons, to draw this prism on a circle, a bit like we did for the cube, at the beginning of the proof of Theorem 4.11. We are led in this way to the following statement:

\begin{theorem}
The prism, which is as follows when drawn on a circle,
$$
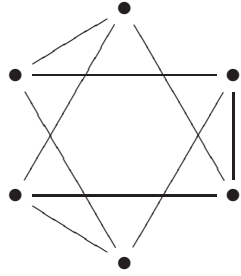
\xymatrix@R=13pt@C=28pt{
&\bullet\ar@{-}[dddl]\ar@{-}[dddr]\ar@{-}[dl]\\
\bullet\ar@{-}[rr]\ar@{-}[dddr]&&\bullet\ar@{-}[dddl]\ar@{-}[dd]\\
\\
\bullet\ar@{-}[rr]\ar@{-}[dr]&&\bullet\\
&\bullet
}$$
is transitive and non-circulant too, exactly as the cube was.
\end{theorem}

\begin{proof}
The fact that the prism is indeed transitive follows from the above discussion, but it is convenient to view this as well directly on the above picture. Indeed:

\medskip

-- The prism as drawn above has 3 obvious symmetry axes, allowing us to do many of the $i\to j$ operations required by the definition of transitivity.

\medskip

-- In addition, the prism is invariant as well by the $120^\circ$ and $240^\circ$ rotations, and when combining this with the above 3 symmetries, we have all that we need.

\medskip

Finally, the fact that the prism is indeed not circulant is quite clear, intuitively speaking, and this can be proved a bit as for the cube, as in the proof of Theorem 4.11.
\end{proof}

Summarizing, we have interesting examples, and our theory of transitive graphs seems worth developing. In order now to reach to something more conceptual, it is pretty much clear that we must get into group theory. So, let us formulate the following definition:

\index{group}
\index{permutation group}

\begin{definition}
A group of permutations is a subset $G\subset S_N$ which is stable under the composition of permutations, and under their inversion. We say that:
\begin{enumerate}
\item $G$ acts transitively on the set $\{1,\ldots,N\}$ if for any two points $i,j$ we can find $\sigma\in G$ mapping one point to another, $\sigma(i)=j$.

\item $G$ acts on a graph $X$ with vertices labeled $1,\ldots,N$ when each $\sigma\in G$ leaves invariant the edges, $k-l\iff\sigma(k)-\sigma(l)$.
\end{enumerate}
Also, we say that $G$ acts transitively on a graph $X$ with vertices labeled $1,\ldots,N$ when it acts on $X$ in the sense of (2), and the action is transitive in the sense of (1).
\end{definition}

All this might seem a bit heavy, but as we will soon discover, is worth the effort, because group theory is a powerful theory, and having it into our picture will be certainly a good thing, a bit similar to the update from the Stone Age to the Bronze Age. Or perhaps to the update from the Bronze Age to the Iron Age, because what we did so far in this book was sometimes non-trivial, and can be counted as Bronze Age weaponry.

\bigskip

As a first good surprise, once Definition 4.13 formulated and digested, our definition of the circulant and transitive graphs becomes something very simple, as follows:

\begin{theorem}
The following happen, for a graph $X$ having $N$ vertices:
\begin{enumerate}
\item $X$ is circulant when we have an action $\mathbb Z_N\curvearrowright X$.

\item $X$ is transitive when we have a transitive action $G\curvearrowright X$.
\end{enumerate}
\end{theorem}

\begin{proof}
This is something trivial and self-explanatory, and with the remark that in (1) we do not have to say something about transitivity, because the subgroup $\mathbb Z_N\subset S_N$ is transitive, in the sense of Definition 4.13. As usual, we have called this statement Theorem instead of Proposition simply due to its theoretical importance.
\end{proof}

As a second good surprise, our previous transitivity considerations regarding the cube and the prism take now a very simple form, in terms of groups, as follows:

\begin{proposition}
The following are transitive graphs:
\begin{enumerate}
\item The cube, due to an action $\mathbb Z_2^3\curvearrowright X$.

\item The prism, due to an action $\mathbb Z_2\times\mathbb Z_3\curvearrowright X$.
\end{enumerate}
\end{proposition}

\begin{proof}
As before with Theorem 4.14, this is trivial and self-explanatory, with the actions being the obvious ones, coming from our previous study of the cube and prism.
\end{proof}

Many things can be said about the transitive graphs, in general, but thinking well, what we would mostly like to have would be an extension of what we did in the previous section for the circulant graphs, including the Fourier transform material, which was something highly non-trivial and powerful, perhaps to a class of graphs smaller than that of the general transitive graphs. So, let us formulate the following definition:

\begin{definition}
A finite graph $X$ is called generalized circulant when it has a transitive action $G\curvearrowright X$, with $G$ being a finite abelian group.
\end{definition}

And this looks like a very good definition. Indeed, as examples we have the circulant graphs, but also the cube, and the prism, since products of abelian groups are obviously abelian. So, no interesting transitive graphs lost, when assuming that $G$ is abelian.

\bigskip

In order now to further build on this definition, and in particular to develop our generalized Fourier transform machinery, as hoped in the above, let us temporarily leave aside the graphs $X$, and focus on the finite abelian groups $G$. We first have:

\begin{proposition}
Given a finite abelian group $G$, the group morphisms
$$\chi:G\to\mathbb T$$
with $\mathbb T$ being the unit circle, called characters of $G$, form a finite abelian group $\widehat{G}$.
\end{proposition}

\begin{proof}
There are several things to be proved here, the idea being as follows:

\medskip

(1) Our first claim is that $\widehat{G}$ is a group, with the pointwise multiplication, namely:
$$(\chi\rho)(g)=\chi(g)\rho(g)$$

Indeed, if $\chi,\rho$ are characters, so is $\chi\rho$, and so the multiplication is well-defined on $\widehat{G}$. Regarding the unit, this is the trivial character, constructed as follows:
$$1:G\to\mathbb T\quad,\quad 
g\to1$$ 

Finally, we have inverses, with the inverse of $\chi:G\to\mathbb T$ being its conjugate:
$$\bar{\chi}:G\to\mathbb T\quad,\quad 
g\to\overline{\chi(g)}$$

(2) Our next claim is that $\widehat{G}$ is finite. Indeed, given a group element $g\in G$, we can talk about its order, which is smallest integer $k\in\mathbb N$ such that $g^k=1$. Now assuming that we have a character $\chi:G\to\mathbb T$, we have the following formula:
$$\chi(g)^k=1$$

Thus $\chi(g)$ must be one of the $k$-th roots of unity, and in particular there are finitely many choices for $\chi(g)$. Thus, there are finitely many choices for $\chi$, as desired.

\medskip

(3) Finally, the fact that $\widehat{G}$ is abelian follows from definitions, because the pointwise multiplication of functions, and in particular of characters, is commutative.
\end{proof}

The above construction is quite interesting, and we have:

\index{cyclic group}
\index{product of cyclic groups}
\index{self-dual group}

\begin{theorem}
The character group operation $G\to\widehat{G}$ for the finite abelian groups, called Pontrjagin duality, has the following properties:
\begin{enumerate}
\item The dual of a cyclic group is the group itself, $\widehat{\mathbb Z}_N=\mathbb Z_N$.

\item The dual of a product is the product of duals, $\widehat{G\times H}=\widehat{G}\times\widehat{H}$.

\item Any product of cyclic groups $G=\mathbb Z_{N_1}\times\ldots\times\mathbb Z_{N_k}$ is self-dual, $G=\widehat{G}$.
\end{enumerate}
\end{theorem}

\begin{proof}
We have several things to be proved, the idea being as follows:

\medskip

(1) A character $\chi:\mathbb Z_N\to\mathbb T$ is uniquely determined by its value $z=\chi(g)$ on the standard generator $g\in\mathbb Z_N$. But this value must satisfy:
$$z^N=1$$

Thus we must have $z\in\mathbb Z_N$, with the cyclic group $\mathbb Z_N$ being regarded this time as being the group of $N$-th roots of unity. Now conversely, any $N$-th root of unity $z\in\mathbb Z_N$ defines a character $\chi:\mathbb Z_N\to\mathbb T$, by setting, for any $r\in\mathbb N$:
$$\chi(g^r)=z^r$$

Thus we have an identification $\widehat{\mathbb Z}_N=\mathbb Z_N$, as claimed.

\medskip

(2) A character of a product of groups $\chi:G\times H\to\mathbb T$ must satisfy:
$$\chi(g,h)=\chi\left[(g,1)(1,h)\right]=\chi(g,1)\chi(1,h)$$

Thus $\chi$ must appear as the product of its restrictions $\chi_{|G},\chi_{|H}$, which must be both characters, and this gives the identification in the statement.

\medskip

(3) This follows from (1) and (2). Alternatively, any character $\chi:G\to\mathbb T$ is uniquely determined by its values $\chi(g_1),\ldots,\chi(g_k)$ on the standard generators of $\mathbb Z_{N_1},\ldots,\mathbb Z_{N_k}$, which must belong to $\mathbb Z_{N_1},\ldots,\mathbb Z_{N_k}\subset\mathbb T$, and this gives $\widehat{G}=G$, as claimed.
\end{proof}

At a more advanced level now, we have the following result:

\begin{theorem}
The finite abelian groups are the following groups,
$$G=\mathbb Z_{N_1}\times\ldots\times\mathbb Z_{N_k}$$
and these groups are all self-dual, $G=\widehat{G}$.
\end{theorem}

\begin{proof}
This is something quite tricky, the idea being as follows:

\medskip

(1) In order to prove our result, assume that $G$ is finite and abelian. For any prime number $p\in\mathbb N$, let us define $G_p\subset G$ to be the subset of elements having as order a power of $p$. Equivalently, this subset $G_p\subset G$ can be defined as follows:
$$G_p=\left\{g\in G\Big|\exists k\in\mathbb N,g^{p^k}=1\right\}$$

(2) It is then routine to check, based on definitions, that each $G_p$ is a subgroup. Our claim now is that we have a direct product decomposition as follows:
$$G=\prod_pG_p$$

(3) Indeed, by using the fact that our group $G$ is abelian, we have a morphism as follows, with the order of the factors when computing $\prod_pg_p$ being irrelevant:
$$\prod_pG_p\to G\quad,\quad (g_p)\to\prod_pg_p$$

Moreover, it is routine to check that this morphism is both injective and surjective, via some simple manipulations, so we have our group decomposition, as in (2).

\medskip

(4) Thus, we are left with proving that each component $G_p$ decomposes as a product of cyclic groups, having as orders powers of $p$, as follows:
$$G_p=\mathbb Z_{p^{r_1}}\times\ldots\times\mathbb Z_{p^{r_s}}$$

But this is something that can be checked by recurrence on $|G_p|$, via some routine computations, and we are led to the conclusion in the statement.

\medskip

(5) Finally, the fact that the finite abelian groups are self-dual, $G=\widehat{G}$, follows from the structure result that we just proved, and from Theorem 4.18 (3).
\end{proof}

In relation now with Fourier analysis, the result is as follows:

\begin{theorem}
Given a finite abelian group $G$, we have an isomorphism as follows, obtained by linearizing/delinearizing the characters,
$$C^*(G)\simeq C(\widehat{G})$$
where $C^*(G)$ is the algebra of functions $\varphi:G\to\mathbb C$, with convolution product, and $C(\widehat{G})$ is the algebra of functions $\varphi:\widehat{G}\to\mathbb C$, with usual product.
\end{theorem}

\begin{proof}
There are many things going on here, the idea being as follows:

\medskip

(1) Given a finite abelian group $G$, we can talk about the complex algebra $C(G)$ formed by the complex functions $\varphi:G\to\mathbb C$, with usual product, namely:
$$(\varphi\psi)(g)=\varphi(g)\psi(g)$$

Observe that we have $C(G)\simeq\mathbb C^N$ as an algebra, where $N=|G|$, with this being best seen via the basis of $C(G)$ formed by the Dirac masses at the points of $G$:
$$C(G)=\left\{\sum_{g\in G}\lambda_g\delta_g\Big|\lambda_g\in\mathbb C\right\}$$

(2) On the other hand, we can talk as well about the algebra $C^*(G)$ formed by the same functions $\varphi:G\to\mathbb C$, but this time with the convolution product, namely:
$$(\varphi*\psi)(g)=\sum_{h\in G}\varphi(gh^{-1})\psi(h)$$

Since we have $\delta_k*\delta_l=\delta_{kl}$ for any $k,l\in G$, as you can easily check by using the above formula, the Dirac masses $\delta_g\in C^*(G)$ behave like the group elements $g\in G$. Thus, we can view our algebra as follows, with multiplication given by $g\cdot h=gh$, and linearity:
$$C^*(G)=\left\{\sum_{g\in G}\lambda_gg\Big|\lambda_g\in\mathbb C\right\}$$

\medskip

(3) Now that we know what the statement is about, let us go for the proof. The first observation is that we have a morphism of algebras as follows:
$$C^*(G)\to C(\widehat{G})\quad,\quad g\to\left[\chi\to\chi(g)\right]$$

Now since on both sides we have vector spaces of dimension $N=|G|$, it is enough to check that this morphism is injective. But this is best done via Theorem 4.19, which shows that the characters $\chi\in\widehat{G}$ separate the points $g\in G$, as desired.
\end{proof}

In practice now, we can clearly feel that Theorem 4.20 is related to Fourier analysis, and more specifically to the Fourier transforms and series that we know from analysis, but also to the discrete Fourier transform from the beginning of this chapter. However, all this remains a bit difficult to clarify, and we have here the following statement:

\begin{fact}
The following happen, regarding the locally compact abelian groups:
\begin{enumerate}
\item What we did in the finite case, namely group characters, and construction and basic properties of the dual, can be extended to them.

\item As basic examples of this, besides what we have in the finite case, and notably $\widehat{\mathbb Z}_N=\mathbb Z_N$, we have $\widehat{\mathbb Z}=\mathbb T$, $\widehat{\mathbb T}=\mathbb Z$, and also $\widehat{\mathbb R}=\mathbb R$.

\item With some care for analytic aspects, $C^*(G)\simeq C(\widehat{G})$ remains true in this setting, and in the case $G=\mathbb R$, this isomorphism is the Fourier transform.
\end{enumerate}
\end{fact}

Obviously, all this is a bit heavy, but you get the point, there are 3 types of Fourier analysis in life, namely the ``standard'' one, that you might know from advanced calculus, corresponding to $G=\mathbb R$, then the ``Fourier series'' one, that you might know from advanced calculus too, corresponding to $G=\mathbb Z,\mathbb T$, and finally the ``discrete'' one that we started to learn in this book, over $G=\mathbb Z_N$ and other finite abelian groups.

\bigskip

In practice, all this is a bit complicated, and back now to the finite abelian groups, let us work out a softer version of all the above, which is what is really needed, in practice, when doing discrete Fourier analysis. We have here the following result:

\begin{theorem}
Given a finite abelian group $G$, with dual group $\widehat{G}=\{\chi:G\to\mathbb T\}$, consider the corresponding Fourier coupling, namely:
$$\mathcal F_G:G\times\widehat{G}\to\mathbb T\quad,\quad 
(i,\chi)\to\chi(i)$$
\begin{enumerate}
\item Via the standard isomorphism $G\simeq\widehat{G}$, this Fourier coupling can be regarded as a usual square matrix, $F_G\in M_G(\mathbb T)$.

\item This matrix $F_G\in M_G(\mathbb T)$ is complex Hadamard, in the sense that its entries are on the unit circle, and its rows are pairwise orthogonal.

\item In the case of the cyclic group $G=\mathbb Z_N$ we obtain in this way, via the standard identification $\mathbb Z_N=\{1,\ldots,N\}$, the Fourier matrix $F_N$.

\item In general, when using a decomposition $G=\mathbb Z_{N_1}\times\ldots\times\mathbb Z_{N_k}$, the corresponding Fourier matrix is given by $F_G=F_{N_1}\otimes\ldots\otimes F_{N_k}$.
\end{enumerate}
\end{theorem}

\begin{proof}
This follows indeed by using the above finite abelian group theory:

\medskip

(1) Via the identification $G\simeq\widehat{G}$, we have indeed a square matrix, given by:
$$(F_G)_{i\chi}=\chi(i)$$

(2) The scalar products between distinct rows are indeed zero, as shown by:
$$<R_i,R_j>
=\sum_\chi\chi(i)\overline{\chi(j)}
=\sum_\chi\chi(i-j)
=|G|\cdot\delta_{ij}$$

(3) This follows from the well-known and elementary fact that, via the identifications $\mathbb Z_N=\widehat{\mathbb Z_N}=\{1,\ldots,N\}$, the Fourier coupling here is as follows, with $w=e^{2\pi i/N}$:
$$(i,j)\to w^{ij}$$

(4) Observe first that $\widehat{H\times K}=\widehat{H}\times\widehat{K}$ gives, at the level of the Fourier couplings, $F_{H\times K}=F_H\otimes F_K$. Now by decomposing $G$ into cyclic groups, as in the statement, and using (3) for the cyclic components, we obtain the formula in the statement.
\end{proof}

So long for circulant graphs, finite group actions, transitive graphs, finite abelian groups, discrete Fourier analysis, and Hadamard matrices. We will be back to all this in Part III of the present book, where we will systematically investigate such things.

\section*{4c. Dihedral groups} 

As something more concrete now, which is a must-know, let us try to compute the dihedral group. This is a famous group, constructed as follows:

\index{dihedral group}
\index{regular polygon}

\begin{definition}
The dihedral group $D_N$ is the symmetry group of 
$$
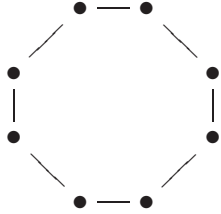
\xymatrix@R=12pt@C=12pt{
&\bullet\ar@{-}[r]\ar@{-}[dl]&\bullet\ar@{-}[dr]\\
\bullet\ar@{-}[d]&&&\bullet\ar@{-}[d]\\
\bullet\ar@{-}[dr]&&&\bullet\ar@{-}[dl]\\
&\bullet\ar@{-}[r]&\bullet}$$
that is, of the regular polygon having $N$ vertices.
\end{definition}

In order to understand how this works, here are the basic examples of regular $N$-gons, at small values of the parameter $N\in\mathbb N$, along with their symmetry groups:

\bigskip

\underline{$N=2$}. Here the $N$-gon is just a segment, and its symmetries are obviously the identity $id$, plus the symmetry $\tau$ with respect to the middle of the segment:
$$
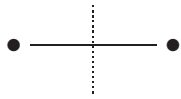
\xymatrix@R=10pt@C=20pt{
&\ar@{.}[dd]\\
\bullet\ar@{-}[rr]&&\bullet\\
&}$$

Thus we have $D_2=\{id,\tau\}$, which in group theory terms means $D_2=\mathbb Z_2$.

\bigskip

\underline{$N=3$}. Here the $N$-gon is an equilateral triangle, and we have 6 symmetries, the rotations of angles $0^\circ$, $120^\circ$, $240^\circ$, and the symmetries with respect to the altitudes: 
$$
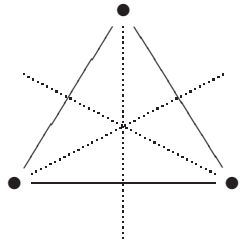
\xymatrix@R=13pt@C=28pt{
&\bullet\ar@{-}[dddr]\ar@{-}[dddl]\ar@{.}[dddd]\\
\ar@{.}[ddrr]&&\ar@{.}[ddll]\\
\\
\bullet\ar@{-}[rr]&&\bullet\\
&
}$$

Alternatively, we can say that the symmetries are all the $3!=6$ possible permutations of the vertices, and so that in group theory terms, we have $D_3=S_3$.

\bigskip

\underline{$N=4$}. Here the $N$-gon is a square, and as symmetries we have 4 rotations, of angles $0^\circ,90^\circ,180^\circ,270^\circ$, as well as 4 symmetries, with respect to the 4 symmetry axes, which are the 2 diagonals, and the 2 segments joining the midpoints of opposite sides:
$$
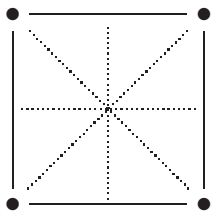
\xymatrix@R=26pt@C=26pt{
\bullet\ar@{-}[dd]\ar@{.}[ddrr]\ar@{-}[rr]&\ar@{.}[dd]&\bullet\ar@{-}[dd]\ar@{.}[ddll]\\
\ar@{.}[rr]&&\\
\bullet\ar@{-}[rr]&&\bullet
}$$

Thus, we obtain as symmetry group some sort of product between $\mathbb Z_4$ and $\mathbb Z_2$. Observe however that this product is not the usual one, our group being not abelian.

\bigskip

\underline{$N=5$}. Here the $N$-gon is a regular pentagon, and as symmetries we have 5 rotations, of angles $0^\circ,72^\circ,144^\circ,216^\circ,288^\circ$, as well as 5 symmetries, with respect to the 5 symmetry axes, which join the vertices to the midpoints of the opposite sides:
$$
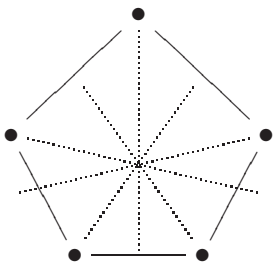
\xymatrix@R=13pt@C=11pt{
&&\bullet\ar@{-}[ddrr]\ar@{-}[ddll]\ar@{.}[dddd]\\
&&&&\\
\bullet\ar@{-}[ddr]\ar@{.}[drrrr]&&&&\bullet\ar@{-}[ddl]\ar@{.}[dllll]\\
&&&&\\
&\bullet\ar@{-}[rr]\ar@{.}[uuurr]&&\bullet\ar@{.}[uuull]&&
}$$

\underline{$N=6$}. Here the $N$-gon is a regular hexagon, and we have 6 rotations, of angles $0^\circ,60^\circ,120^\circ,180^\circ,240^\circ,300^\circ$, and 6 symmetries, with respect to the 6 symmetry axes, which are the 3 diagonals, and the 3 segments joining the midpoints of opposite sides:
$$
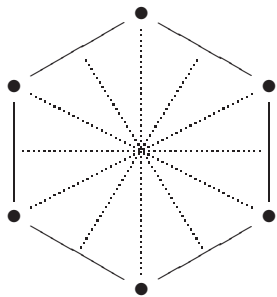
\xymatrix@R=4pt@C=14pt{
&&\bullet\ar@{-}[ddrr]\ar@{-}[ddll]\ar@{.}[dddddddd]\\
&\ar@{.}[ddddddrr]&&\ar@{.}[ddddddll]\\
\bullet\ar@{-}[dddd]\ar@{.}[ddddrrrr]&&&&\bullet\ar@{-}[dddd]\ar@{.}[ddddllll]\\
&&&&\\
\ar@{.}[rrrr]&&&&\\
&&&&\\
\bullet\ar@{-}[ddrr]&&&&\bullet\ar@{-}[ddll]\\
&&&&\\
&&\bullet
}$$

We can see from the above that the various dihedral groups $D_N$ have many common features, and that there are some differences as well. In general, we have:

\begin{proposition}
The dihedral group $D_N$ has $2N$ elements, as follows:
\begin{enumerate}
\item We have $N$ rotations $R_1,\ldots,R_N$, with $R_k$ being the rotation of angle $2k\pi/N$. When labelling the vertices $1,\ldots,N$, the formula is $R_k:i\to k+i$.

\item We have $N$ symmetries $S_1,\ldots,S_N$, with $S_k$ being the symmetry with respect to the $Ox$ axis rotated by $k\pi/N$. The symmetry formula is $S_k:i\to k-i$.
\end{enumerate}
\end{proposition}

\begin{proof}
Our group is indeed formed of $N$ rotations, of angles $2k\pi/N$ with $k=1,\ldots,N$, and then of the $N$ symmetries with respect to the $N$ possible symmetry axes, which are the $N$ medians of the $N$-gon when $N$ is odd, and are the $N/2$ diagonals plus the $N/2$ lines connecting the midpoints of opposite edges, when $N$ is even.
\end{proof}

With the above result in hand, we can talk about $D_N$ abstractly, as follows:

\begin{theorem}
The dihedral group $D_N$ is the group having $2N$ elements, $R_1,\ldots,R_N$ and $S_1,\ldots,S_N$, called rotations and symmetries, which multiply as follows,
$$R_kR_l=R_{k+l}\quad,\quad 
R_kS_l=S_{k+l}$$
$$S_kR_l=S_{k-l}\quad,\quad 
S_kS_l=R_{k-l}$$
with all the indices being taken modulo $N$.
\end{theorem}

\begin{proof}
With notations from Proposition 4.24, the various compositions between rotations and symmetries can be computed as follows:
$$R_kR_l\ :\ i\to l+i\to k+l+i\quad,\quad 
R_kS_l\ :\ i\to l-i\to k+l-i$$
$$S_kR_l\ :\ i\to l+i\to k-l-i\quad,\quad 
S_kS_l\ :\ i\to l-i\to k-l+i$$

But these are exactly the formulae for $R_{k+l},S_{k+l},S_{k-l},R_{k-l}$, as stated. Now since a group is uniquely determined by its multiplication rules, this gives the result.
\end{proof}

Observe that $D_N$ has the same cardinality as $E_N=\mathbb Z_N\times\mathbb Z_2$. We obviously don't have $D_N\simeq E_N$, because $D_N$ is not abelian, while $E_N$ is. So, our next goal will be that of proving that $D_N$ appears by ``twisting'' $E_N$. In order to do this, let us start with:

\begin{proposition}
The group $E_N=\mathbb Z_N\times\mathbb Z_2$ is the group having $2N$ elements, $r_1,\ldots,r_N$ and $s_1,\ldots,s_N$, which multiply according to the following rules,
$$r_kr_l=r_{k+l}\quad,\quad 
r_ks_l=s_{k+l}$$
$$s_kr_l=s_{k+l}\quad,\quad 
s_ks_l=r_{k+l}$$
with all the indices being taken modulo $N$.
\end{proposition}

\begin{proof}
With the notation $\mathbb Z_2=\{1,\tau\}$, the elements of the product group $E_N=\mathbb Z_N\times\mathbb Z_2$ can be labelled $r_1,\ldots,r_N$ and $s_1,\ldots,s_N$, as follows:
$$r_k=(k,1)\quad,\quad
s_k=(k,\tau)$$

These elements multiply then according to the formulae in the statement. Now since a group is uniquely determined by its multiplication rules, this gives the result.
\end{proof}

Let us compare now Theorem 4.25 and Proposition 4.26. In order to formally obtain $D_N$ from $E_N$, we must twist some of the multiplication rules of $E_N$, namely:
$$s_kr_l=s_{k+l}\to s_{k-l}\quad,\quad 
s_ks_l=r_{k+l}\to r_{k-l}$$

Informally, this amounts in following the rule ``$\tau$ switches the sign of what comes afterwards", and we are led in this way to the following definition:

\index{crossed product}

\begin{definition}
Given two groups $A,G$, with an action $A\curvearrowright G$, the crossed product
$$P=G\rtimes A$$
is the set $G\times A$, with multiplication $(g,a)(h,b)=(gh^a,ab)$.
\end{definition}

Now with this technology in hand, by getting back to the dihedral group $D_N$, we can improve Theorem 4.25, into a final result on the subject, as follows:

\index{dihedral group}
\index{crossed product decomposition}

\begin{theorem}
We have a crossed product decomposition as follows,
$$D_N=\mathbb Z_N\rtimes\mathbb Z_2$$
with $\mathbb Z_2=\{1,\tau\}$ acting on $\mathbb Z_N$ via switching signs, $k^\tau=-k$.
\end{theorem}

\begin{proof}
We have an action $\mathbb Z_2\curvearrowright\mathbb Z_N$ given by the formula in the statement, namely $k^\tau=-k$, so we can consider the corresponding crossed product group:
$$P_N=\mathbb Z_N\rtimes\mathbb Z_2$$

In order to understand the structure of $P_N$, we follow Proposition 4.26. The elements of $P_N$ can indeed be labelled $\rho_1,\ldots,\rho_N$ and $\sigma_1,\ldots,\sigma_N$, as follows:
$$\rho_k=(k,1)\quad,\quad 
\sigma_k=(k,\tau)$$

Now when computing the products of such elements, we basically obtain the formulae in Proposition 4.26, perturbed as in Definition 4.27. To be more precise, we have:
$$\rho_k\rho_l=\rho_{k+l}\quad,\quad 
\rho_k\sigma_l=\sigma_{k+l}$$
$$\sigma_k\rho_l=\sigma_{k+l}\quad,\quad 
\sigma_k\sigma_l=\rho_{k+l}$$

But these are exactly the multiplication formulae for $D_N$, from Theorem 4.25. Thus, we have an isomorphism $D_N\simeq P_N$ given by $R_k\to\rho_k$ and $S_k\to\sigma_k$, as desired.
\end{proof}

\section*{4d. Cayley graphs}

We have kept the best for the end. We have the following notion, that we already met in a vague form in chapter 1, and which is something quite far-reaching:

\index{Cayley graph}

\begin{definition}
Associated to any finite group $G=<S>$, with the generating set $S$ assumed to satisfy $1\notin S=S^{-1}$, is its Cayley graph, constructed as follows:
\begin{enumerate}
\item The vertices are the group elements $g\in G$.

\item Edges $g-h$ are drawn when $h=gs$, with $s\in S$.
\end{enumerate} 
\end{definition}

As a first observation, the Cayley graph is indeed a graph, because our assumption $S=S^{-1}$ on the generating set shows that we have $g-h\implies h-g$, as we should, and also because our assumption $1\notin S$ excludes the self-edges, $g\not\!\!-\,g$. 

\bigskip

We will see in what follows that most of the graphs that we met so far are Cayley graphs, and in addition, with $G$ being abelian in most cases. Before starting with the examples, however, let us point out that the Cayley graph depends, in a crucial way, on the generating set $S$. Indeed, if we choose for instance $S=G-\{1\}$, we obtain the complete graph with $N=|G|$ vertices, and this regardless of what our group $G$ is.

\bigskip

Thus, the Cayley graph as constructed above is not exactly associated to the group $G$, but rather to the group $G$ viewed as finitely generated group, $G=<S>$.

\bigskip

In order to construct now examples, let us start with the simplest finite groups that we know. The simplest such group is the cyclic group $\mathbb Z_N$, and we have:

\begin{proposition}
The $N$-gon graph, namely
$$
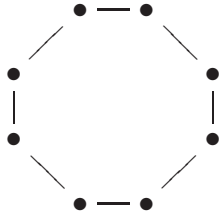
\xymatrix@R=12pt@C=12pt{
&\bullet\ar@{-}[r]\ar@{-}[dl]&\bullet\ar@{-}[dr]\\
\bullet\ar@{-}[d]&&&\bullet\ar@{-}[d]\\
\bullet\ar@{-}[dr]&&&\bullet\ar@{-}[dl]\\
&\bullet\ar@{-}[r]&\bullet}$$
is the Cayley graph of $\mathbb Z_N=<\pm 1>$, with $\mathbb Z_N$ being written here additively.
\end{proposition}

\begin{proof}
This is clear, because in additive notation, our condition for the edges $g-h$ reads $g=h\pm1$, so we are led to the $N$-gon graph in the statement.
\end{proof}

The next thing that we can do is to look at the products of cyclic groups. In the simplest case here, that of a product of the group $\mathbb Z_2$ with itself, we obtain:

\begin{theorem}
The square graph, namely
$$
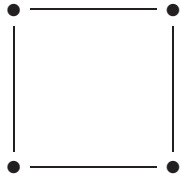
\xymatrix@R=20pt@C=20pt{
\bullet\ar@{-}[dd]\ar@{-}[rr]&&\bullet\ar@{-}[dd]\\
\\
\bullet\ar@{-}[rr]&&\bullet
}$$
is the Cayley graph of $\mathbb Z_2\times\mathbb Z_2$, with generating set $S=\{(1,0),(0,1)\}$.
\end{theorem}

\begin{proof}
This is something which requires a bit of thinking. We must first draw the 4 elements of $\mathbb Z_2\times\mathbb Z_2$, and this is best done by using their coordinates, as follows:
$$\xymatrix@R=17pt@C=17pt{
10&&11\\
\\
00&&10
}$$

Now let us construct the edges. In additive notation, our condition for the edges $(a,b)-(c,d)$ reads either $(a,b)=(c,d)+(1,0)$ or $(a,b)=(c,d)+(0,1)$, which amounts in saying that the passage $(a,b)\to(c,d)$ appears by modifying one coordinate, and keeping the other coordinate fixed. We conclude from this that the edges are as follows:
$$
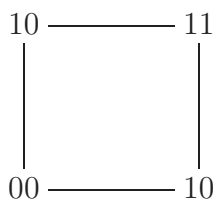
\xymatrix@R=20pt@C=20pt{
10\ar@{-}[dd]\ar@{-}[rr]&&11\ar@{-}[dd]\\
\\
00\ar@{-}[rr]&&10
}$$

Now by removing the vertex labels, we obtain the usual square, as claimed.
\end{proof}

What comes next? Obviously, more complicated products of cyclic groups. Here the situation ramifies a bit, and as our next basic example, we have:

\begin{proposition}
The prism graph, namely
$$
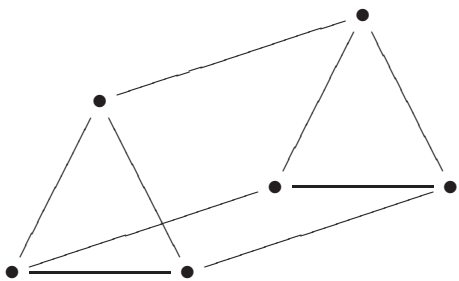
\xymatrix@R=20pt@C=20pt{
&&&&\bullet\ar@{-}[ddl]\ar@{-}[ddr]\\
&\bullet\ar@{-}[ddl]\ar@{-}[ddr]\ar@{-}[urrr]\\
&&&\bullet\ar@{-}[rr]&&\bullet\\
\bullet\ar@{-}[rr]\ar@{-}[urrr]&&\bullet\ar@{-}[urrr]
}$$
appears as Cayley graph of the group $\mathbb Z_2\times\mathbb Z_3$.
\end{proposition}

\begin{proof}
This follows a bit as before, for the square, and we will leave the details here, including of course finding the correct generating set $S$, as an instructive exercise.
\end{proof}

As a main example now, obtained by staying with $\mathbb Z_2$, we have:

\index{cube}

\begin{theorem}
The cube graph, namely
$$
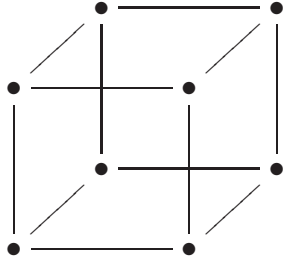
\xymatrix@R=18pt@C=20pt{
&\bullet\ar@{-}[rr]&&\bullet\\
\bullet\ar@{-}[rr]\ar@{-}[ur]&&\bullet\ar@{-}[ur]\\
&\bullet\ar@{-}[rr]\ar@{-}[uu]&&\bullet\ar@{-}[uu]\\
\bullet\ar@{-}[uu]\ar@{-}[ur]\ar@{-}[rr]&&\bullet\ar@{-}[uu]\ar@{-}[ur]
}$$
is the Cayley graph of $\mathbb Z_2^3$, with generating set $S=\{(1,0,0),(0,1,0),(0,0,1)\}$.
\end{theorem}

\begin{proof}
This is clear, as before for the square. In fact, we have already met this in chapter 1, with details, and this is our come-back to the subject, long overdue.
\end{proof}

Looking at what we have so far, it is pretty much clear that Theorem 4.31 and Theorem 4.33 are of the same nature. The joint extension of these theorems, along with the computation of the corresponding symmetry group, leads to the following statement:

\index{reflection group}
\index{hyperoctahedral group}

\begin{theorem}
The Cayley graph of the group $\mathbb Z_2^N$ is the hypercube 
$$\square_N\subset\mathbb R^N$$
and the symmetry group of $\square_N$, called hyperoctahedral group, is given by
$$H_N=\mathbb Z_2^N\rtimes S_N$$
with the permutations acting via $\sigma(e_1,\ldots,e_k)=(e_{\sigma(1)},\ldots,e_{\sigma(k)})$.
\end{theorem}

\begin{proof}
Regarding the Cayley graph assertion, this follows as for the square, or for the cube. As for the computation of $H_N$, this can be done as follows:

\medskip

(1) Consider the standard cube in $\mathbb R^N$, centered at 0, and having as vertices the points having coordinates $\pm1$. With this picture in hand, it is clear that the symmetries of the cube coincide with the symmetries of the $N$ coordinate axes of $\mathbb R^N$.

\medskip

(2) In order to count now these symmetries, observe first that we have as examples the $N!$ permutations of the $N$ coordinate axes of $\mathbb R^N$. But each of these permutations $\sigma\in S_N$ can be further ``decorated'' by a sign vector $e\in\{\pm1\}^N$, consisting of the possible $\pm1$ flips which can be applied to each coordinate axis, at the arrival. Thus, we have:
$$|H_N|
=|S_N|\cdot|\mathbb Z_2^N|
=N!\cdot2^N$$

(3) Now observe that at the level of the cardinalities, the above formula gives:
$$|H_N|=|\mathbb Z_2^N\times S_N|$$

To be more precise, given an element $g\in H_N$, the element $\sigma\in S_N$ is the corresponding permutation of the $N$ coordinate axes, regarded as unoriented lines in $\mathbb R^N$, and $e\in\mathbb Z_2^N$ is the vector collecting the possible flips of these coordinate axes, at the arrival. Now observe that the product formula for two such pairs $g=(e,\sigma)$ is as follows, with the permutations $\sigma\in S_N$ acting on the elements $f\in\mathbb Z_2^N$ as in the statement:
$$(e,\sigma)(f,\tau)=(ef^\sigma,\sigma\tau)$$

(4) Thus, we are precisely in the framework of Definition 4.27, and we conclude that we have a crossed product decomposition $H_N=\mathbb Z_2^N\rtimes S_N$, as in the statement.
\end{proof}

Many more things can be said here, and more on all this in Part III of this book.

\section*{4e. Exercises}

We had a rather easy and relaxing chapter here, and as exercises, we have:

\begin{exercise}
How to decide if an unlabeled graph is circulant or not?
\end{exercise}

\begin{exercise}
Learn more about permutation groups, their orbits and orbitals.
\end{exercise}

\begin{exercise}
Learn about Cayley embeddings $G\subset S_{|G|}$, and about $S_N\subset O_N$ too.
\end{exercise}

\begin{exercise}
Work out all details for the result on the finite abelian groups.
\end{exercise}

\begin{exercise}
Develop some theory for the generalized circulant graphs.
\end{exercise}

\begin{exercise}
Write down all possible decompositions of the group $D_4$.
\end{exercise}

\begin{exercise}
Compute the Cayley graphs of some non-abelian groups.
\end{exercise}

\begin{exercise}
If you teach, ask your students to compute $|H_3|$.
\end{exercise}

As bonus exercise, learn some group theory. We will need this, for Part III.

\part{Geometric aspects}

\ \vskip50mm

\begin{center}
{\em Lady, hear me tonight

Cause my feeling is just so right

As we dance by the moonlight

Can't you see you're my delight}
\end{center}

\chapter{The Laplacian}

\section*{5a. Laplace operator}

Welcome to geometry. Many things can be said here, but going straight to the point, there is nothing more geometric than waves. Indeed, waves produce physics, with each major branch of physics being guided by its own wave equation. And then physics produces, as we know it since Newton, calculus, also known as ``geometry and analysis''.

\bigskip

So, this will be our philosophy in what follows, waves, physics, calculus, geometry and analysis being more or less the same thing. And with our graphs being expected to help at some point, in relation with discretization questions, for waves, physics and so on. 

\bigskip

Getting started now, we first need to talk about waves. We have here:

\begin{fact}
Waves can be of many types, as follows:
\begin{enumerate}
\item Mechanical waves, such as the usual water waves, but also the sound waves, or the seismic waves. In all these cases, the wave propagates mechanically, via a certain medium, which can be solid, liquid or gaseous.

\item Electromagnetic waves, coming via a more complicated mechanism, namely an accelerating charge in the context of electromagnetism. These are the radio waves, microwaves, IR, visible light, UV, X-rays and $\gamma$-rays.

\item Other waves, including the heat waves, appearing in the context of heat diffusion through a certain material, and also the waves from quantum mechanics, describing the movements, which are wave-like, of the elementary particles.
\end{enumerate}
\end{fact}

Quite remarkably, the behavior of all the above mechanical and electromagnetic waves is basically described by the same wave equation, which looks as follows:
$$\ddot{\varphi}=v^2\Delta\varphi$$

As for the heat waves and quantum mechanical waves, these are described by some similar equations, namely $\dot{\varphi}=\alpha\Delta\varphi$, and versions of it. Which leads us into:

\begin{question}
What is the Laplace operator $\Delta$, appearing in the above, and making the world go round? 
\end{question}

Which might look a bit scary, I know, are we supposed now to stop mathematics, and take lots of physics classes, for a few years in a row, in order to answer our question. 

\bigskip

Fortunately for us, this world was made with some built-in beauty, among others with basic mathematics corresponding to basic physics, and vice versa. So, on purely philosophical grounds, expect the answer to Question 5.2 to be something very simple, mathematically speaking. And simple that mathematical answer indeed is, as follows:

\begin{answer}
The second derivative of a function $\varphi:\mathbb R^N\to\mathbb R$, making the formula
$$\varphi(x+h)\simeq\varphi(x)+\varphi'(x)h+\frac{<\varphi''(x)h,h>}{2}$$
work, is its Hessian matrix $\varphi''(x)\in M_N(\mathbb R)$, given by the following formula:
$$\varphi''(x)=\left(\frac{d^2\varphi}{dx_idx_j}\right)_{ij}$$
However, when needing a number, as second derivative, the trace of $\varphi''(x)$, denoted
$$\Delta\varphi=\sum_{i=1}^N\frac{d^2\varphi}{dx_i^2}$$
and called Laplacian of $\varphi$, is usually the correct quantity.
\end{answer}

Very nice all this, and as a plan now for the present chapter, we will first try to understand what Answer 5.3 says, from a purely mathematical perspective. Then, we will make a detour through lattices and graphs, further building on the graph theory developed in Part I. And finally, armed with all this knowledge, we will discuss waves.

\bigskip

Getting started for good now, as said in Answer 5.3, it is all about the Taylor formula at order 2. And here, in one variable, to start with, things are quite simple, and you certainly know them well. Indeed, given a function $\varphi:\mathbb R\to\mathbb R$, the first job is that of finding a quantity $\varphi'(x)\in\mathbb R$ making the following approximation formula work:
$$\varphi(x+h)\simeq\varphi(x)+\varphi'(x)h$$

But here, there are not so many choices, and the solution is that of defining the number $\varphi'(x)\in\mathbb R$ by the following formula, provided that the limit converges indeed:
$$\varphi'(x)=\lim_{h\to0}\frac{\varphi(x+h)-\varphi(x)}{h}$$

This number is called derivative of $\varphi$ at the point $x\in\mathbb R$, and as you surely know, geometrically, we have $\varphi'(x)=\tan\alpha$, with $\alpha$ being the slope of $\varphi$ at the point $x$. Now still in one variable, we can talk as well about second derivatives, as follows:

\begin{theorem}
The second derivative of a function $\varphi:\mathbb R\to\mathbb R$, making the formula
$$\varphi(x+h)\simeq\varphi(x)+\varphi'(x)h+\frac{\varphi''(x)h^2}{2}$$
work, is the derivative $\varphi''$ of the derivative $\varphi':\mathbb R\to\mathbb R$.
\end{theorem}

\begin{proof}
Assume indeed that $\varphi$ is twice differentiable at $x$, and let us try to construct an approximation of $\varphi$ around $x$ by a quadratic function, as follows:
$$\varphi(x+h)\simeq a+bh+ch^2$$

We must have $a=\varphi(x)$, and we also know that $b=\varphi'(x)$ is the correct choice for the coefficient of $h$. Thus, our approximation must be as follows:
$$\varphi(x+h)\simeq\varphi(x)+\varphi'(x)h+ch^2$$

In order to find the correct choice for $c\in\mathbb R$, observe that the function $\psi(h)=\varphi(x+h)$ matches with $P(h)=\varphi(x)+\varphi'(x)h+ch^2$ in what regards the value at $h=0$, and also in what regards the value of the derivative at $h=0$. Thus, the correct choice of $c\in\mathbb R$ should be the one making match the second derivatives at $h=0$, and this gives:
$$c=\frac{\varphi''(x)}{2}$$

We are therefore led to the formula in the statement. In order to prove now this formula, $\psi(h)\simeq P(h)$, we can use L'H\^opital's rule, and we obtain, as desired:
\begin{eqnarray*}
\frac{\psi(h)-P(h)}{h^2}
&\simeq&\frac{\psi'(h)-P'(h)}{2h}\\
&\simeq&\frac{\psi''(h)-P''(h)}{2}\\
&=&\frac{\varphi''(h)-\varphi''(h)}{2}\\
&=&0
\end{eqnarray*}

Thus, we are led to the conclusion in the statement.
\end{proof}

Many more things can be said about second derivatives, and let us record here:

\begin{proposition}
Intuitively speaking, the second derivative $\varphi''(x)\in\mathbb R$ computes how much different is $\varphi(x)$, compared to the average of $\varphi(y)$, with $y\simeq x$.
\end{proposition}

\begin{proof}
This is obviously something a bit heuristic, but which is good to know. Let us write the formula in Theorem 5.4, as such, and with $h\to-h$ too:
$$\varphi(x+h)\simeq\varphi(x)+\varphi'(x)h+\frac{\varphi''(x)}{2}\,h^2$$
$$\varphi(x-h)\simeq\varphi(x)-\varphi'(x)h+\frac{\varphi''(x)}{2}\,h^2$$

By making the average, we obtain the following formula:
$$\frac{\varphi(x+h)+\varphi(x-h)}{2}=\varphi(x)+\frac{\varphi''(x)}{2}\,h^2$$

Thus, thinking a bit, we are led to the conclusion in the statement. It is of course possible to say more here, but we will not really need all this, in what follows.
\end{proof}

Moving now to several variables, as a first job, given a function $\varphi:\mathbb R^N\to\mathbb R$, we would like to find a quantity $\varphi'(x)$ making the following approximation formula work:
$$\varphi(x+h)\simeq\varphi(x)+\varphi'(x)h$$

But here, again there are not so many choices, and the solution is that of defining $\varphi'(x)$ as being the row vector formed by the partial derivatives at $x$:
$$\varphi'(x)=\left(\frac{d\varphi}{dx_1}\ \ \ldots\ \ \frac{d\varphi}{dx_N}\right)$$

To be more precise, with this value for $\varphi'(x)$, our approximation formula $\varphi(x+h)\simeq\varphi(x)+\varphi'(x)h$ makes sense indeed, as an equality of real numbers, with $\varphi'(x)h\in\mathbb R$ being obtained as the matrix multiplication of the row vector $\varphi'(x)$, and the column vector $h$. As for the fact that our formula holds indeed, this follows by putting together the approximation properties of each of the partial derivatives $d\varphi/dx_i$, which gives:
$$\varphi(x+h)\simeq\varphi(x)+\sum_{i=1}^N\frac{d\varphi}{dx_i}\cdot h_i=\varphi(x)+\varphi'(x)h$$

Moving now to second derivatives, we have here, generalizing Theorem 5.4:

\begin{theorem}
The second derivative of a function $\varphi:\mathbb R^N\to\mathbb R$, making the formula
$$\varphi(x+h)\simeq\varphi(x)+\varphi'(x)h+\frac{<\varphi''(x)h,h>}{2}$$
work, is its Hessian matrix $\varphi''(x)\in M_N(\mathbb R)$, given by the following formula:
$$\varphi''(x)=\left(\frac{d^2\varphi}{dx_idx_j}\right)_{ij}$$
Moreover, this Hessian matrix is symmetric, $\varphi''(x)_{ij}=\varphi'(x)_{ji}$.
\end{theorem}

\begin{proof}
There are several things going on here, the idea being as follows:

\medskip

(1) As a first observation, at $N=1$ the Hessian matrix constructed above is simply the $1\times1$ matrix having as entry the second derivative $\varphi''(x)$, and the formula in the statement is something that we know well from Theorem 5.4, namely:
$$\varphi(x+h)\simeq\varphi(x)+\varphi'(x)h+\frac{\varphi''(x)h^2}{2}$$

(2) At $N=2$ now, we obviously need to differentiate $\varphi$ twice, and the point is that we come in this way upon the following formula, called Clairaut formula:
$$\frac{d^2\varphi}{dxdy}=\frac{d^2\varphi}{dydx}$$

But, is this formula correct or not? As an intuitive justification for it, let us consider a product of power functions, $\varphi(z)=x^py^q$. We have then:
$$\frac{d^2\varphi}{dxdy}=\frac{d}{dx}\left(\frac{dx^py^q}{dy}\right)=\frac{d}{dx}\left(qx^py^{q-1}\right)=pqx^{p-1}y^{q-1}$$
$$\frac{d^2\varphi}{dydx}=\frac{d}{dy}\left(\frac{dx^py^q}{dx}\right)=\frac{d}{dy}\left(px^{p-1}y^q\right)=pqx^{p-1}y^{q-1}$$

Thus, we can see that the Clairaut formula holds indeed, due to the fact that the functions in $x,y$ commute. Of course, this does not really prove our formula, in general. But exercise for you, to have this idea fully working, by using linear combinations and density arguments, or to look up the standard proof, using the mean value theorem.

\medskip

(3) Moving now to $N=3$ and higher, we can use here the Clairaut formula with respect to any pair of coordinates, and we obtain the Schwarz formula, namely:
$$\frac{d^2\varphi}{dx_idx_j}=\frac{d^2\varphi}{dx_jdx_i}$$

Thus, the second derivative, or Hessian matrix, is symmetric, as claimed.

\medskip

(4) Getting now to the main topic, namely approximation formula in the statement, let $y\in\mathbb R^N$, and consider the following function, with $r\in\mathbb R$:
$$f(r)=\varphi(x+ry)$$

We know from (1) that the Taylor formula for $f$, at the point $r=0$, reads:
$$f(r)\simeq f(0)+f'(0)r+\frac{f''(0)r^2}{2}$$

And our claim is that, with $h=ry$, this is precisely the formula in the statement.

\medskip

(5) So, let us see if our claim is correct. By using the chain rule, we have the following formula, with on the right, as usual, a row vector multiplied by a column vector:
$$f'(r)=\varphi'(x+ry)\cdot y$$

By using again the chain rule, we can compute the second derivative as well:
\begin{eqnarray*}
f''(r)
&=&(\varphi'(x+ry)\cdot y)'\\
&=&\left(\sum_i\frac{d\varphi}{dx_i}(x+ry)\cdot y_i\right)'\\
&=&\sum_i\sum_j\frac{d^2\varphi}{dx_idx_j}(x+ry)\cdot y_iy_j\\
&=&<\varphi''(x+ry)y,y>
\end{eqnarray*}

(6) Time now to conclude. We know that we have $f(r)=\varphi(x+ry)$, and according to our various computations above, we have the following formulae:
$$f(0)=\varphi(x)\quad,\quad 
f'(0)=\varphi'(x)\quad,\quad 
f''(0)=<\varphi''(x)y,y>$$

But with this data in hand, the usual Taylor formula for our one variable function $f$, at order 2, at the point $r=0$, takes the following form, with $h=ry$:
\begin{eqnarray*}
\varphi(x+ry)
&\simeq&\varphi(x)+\varphi'(x)ry+\frac{<\varphi''(x)y,y>r^2}{2}\\
&=&\varphi(x)+\varphi'(x)t+\frac{<\varphi''(x)h,h>}{2}
\end{eqnarray*}

Thus, we have obtained the formula in the statement.
\end{proof}

Getting back now to what we wanted to do, namely understand Answer 5.3, it remains to talk about the Laplace operator $\Delta$. Things are quite tricky here, basically requiring some physics that we still need to develop, but as something mathematical to start with, we have the following higher dimensional analogue of Proposition 5.5:

\begin{proposition}
Intuitively, the following quantity, called Laplacian of $\varphi$,
$$\Delta\varphi=\sum_{i=1}^N\frac{d^2\varphi}{dx_i^2}$$
computes how much different is $\varphi(x)$, compared to the average of $\varphi(y)$, with $y\simeq x$.
\end{proposition}

\begin{proof}
As before with Proposition 5.5, this is something a bit heuristic, but good to know. Let us write the formula in Theorem 5.6, as such, and with $h\to-h$ too:
$$\varphi(x+h)\simeq\varphi(x)+\varphi'(x)h+\frac{<\varphi''(x)h,h>}{2}$$
$$\varphi(x-h)\simeq\varphi(x)-\varphi'(x)h+\frac{<\varphi''(x)h,h>}{2}$$

By making the average, we obtain the following formula:
$$\frac{\varphi(x+h)+\varphi(x-h)}{2}=\varphi(x)+\frac{<\varphi''(x)h,h>}{2}$$

Thus, thinking a bit, we are led to the conclusion in the statement. It is of course possible to say more here, but we will not really need all the details, in what follows.
\end{proof}

With this understood, the problem is now, what can we say about the mathematics of $\Delta$? Inspired by linear algebra, let us formulate a basic question, as follows: 

\begin{question}
The Laplace operator being a linear operator,
$$\Delta(a\varphi+b\psi)=a\Delta\varphi+b\Delta\psi$$
what can we say about it, inspired by usual linear algebra?
\end{question}

In answer now, the space of functions $\varphi:\mathbb R^N\to\mathbb R$, on which $\Delta$ acts, being infinite dimensional, the usual tools from linear algebra do not apply as such, and we must be extremely careful. For instance, we cannot really expect to diagonalize $\Delta$.

\bigskip

Thinking some more, there is actually a real bug too with our problem, because at $N=1$ this problem becomes ``what can we say about the second derivatives $\varphi'':\mathbb R\to\mathbb R$ of the functions $\varphi:\mathbb R\to\mathbb R$, inspired by linear algebra'', with answer ``not much''.

\bigskip

And by thinking even more, still at $N=1$, there is a second bug too, because if $\varphi:\mathbb R\to\mathbb R$ is twice differentiable, nothing will guarantee that its second derivative $\varphi'':\mathbb R\to\mathbb R$ is twice differentiable too. Thus, we have some issues with the domain and range of $\Delta$, regarded as linear operator, and these problems will persist at higher $N$.

\bigskip

So, shall we trash Question 5.8? Not so quick, because, very remarkably, we have:

\begin{fact}
The functions $\varphi:\mathbb R^N\to\mathbb R$ which are $0$-eigenvectors of $\Delta$,
$$\Delta\varphi=0$$
called harmonic functions, have the following properties:
\begin{enumerate}
\item At $N=1$, nothing spectacular, these are just the linear functions.

\item At $N=2$, these are, locally, the real parts of holomorphic functions.

\item At $N\geq 3$, these still share many properties with the holomorphic functions. 
\end{enumerate}
\end{fact}

In order to understand this, or at least get introduced to it, let us first look at the case $N=2$. Here, any function $\varphi:\mathbb R^2\to\mathbb R$ can be regarded as function $\varphi:\mathbb C\to\mathbb R$, depending on $z=x+iy$. But, in view of this, it is natural to enlarge to attention to the functions $\varphi:\mathbb C\to\mathbb C$, and ask which of these functions are harmonic, $\Delta\varphi=0$. And here, we have the following remarkable result, making the link with complex analysis:

\begin{theorem}
Any holomorphic function $\varphi:\mathbb C\to\mathbb C$, when regarded as function
$$\varphi:\mathbb R^2\to\mathbb C$$
is harmonic. Moreover, the conjugates $\bar{\varphi}$ of holomorphic functions are harmonic too.
\end{theorem}

\begin{proof}
The first assertion comes from the following computation, with $z=x+iy$:
\begin{eqnarray*}
\Delta z^n
&=&\frac{d^2z^n}{dx^2}+\frac{d^2z^n}{dy^2}\\
&=&\frac{d(nz^{n-1})}{dx}+\frac{d(inz^{n-1})}{dy}\\
&=&n(n-1)z^{n-2}-n(n-1)z^{n-2}\\
&=&0
\end{eqnarray*}

As for the second assertion, this follows from $\Delta\bar{\varphi}=\overline{\Delta\varphi}$, which is clear from definitions, and which shows that if $\varphi$ is harmonic, than so is its conjugate $\bar{\varphi}$.
\end{proof}

Many more things can be said, along these lines, notably a proof of the assertion (2) in Fact 5.9, which is however a quite tough piece of mathematics, and then with a clarification of the assertion (3) too, from that same principle, which again requires some substantial mathematics. But, in what follows, we will not really need all this.

\section*{5b. Graphs, lattices}

Back now to graphs, we have the following definition, inspired by the above:

\index{Laplacian}

\begin{definition}
We call Laplacian of a graph $X$ the matrix
$$L=v-d$$
with $v$ being the diagonal valence matrix, and $d$ being the adjacency matrix.
\end{definition} 

This definition is inspired by differential geometry, or just by multivariable calculus, and more specifically by the well-known Laplace operator there, given by:
$$\Delta f=\sum_{i=1}^N\frac{d^2f}{dx_i^2}$$

More on this in a moment, but as a word regarding terminology, this is traditionally confusing in graph theory, and impossible to fix in a decent way, according to:

\index{negative Laplacian}
\index{positive Laplacian}

\begin{warning}
The graph Laplacian above is in fact the negative Laplacian,
$$L=-\Delta$$
with our preference for it, negative, coming from the fact that it is positive, $L\geq0$.
\end{warning}

Which sounds like a bad joke, but this is how things are, and more on this a moment. In practice now, the graph Laplacian is given by the following formula:
$$L_{ij}=\begin{cases}
v_i&{\rm if}\ i=j\\
-1&{\rm if}\ i-j\\
0&{\rm otherwise}
\end{cases}$$

Alternatively, we have the following formula, for the entries of the Laplacian:
$$L_{ij}=\delta_{ij}v_i-\delta_{i-j}$$

With these formulae in hand, we can formulate, as our first result on the subject:

\index{harmonic function}

\begin{proposition}
A function on a graph is harmonic, $Lf=0$, precisely when
$$f_i=\frac{1}{v_i}\sum_{i-j}f_j$$
that is, when the value at each vertex is the average over the neighbors.
\end{proposition}

\begin{proof}
We have indeed the following computation, for any function $f$:
\begin{eqnarray*}
(Lf)_i
&=&\sum_jL_{ij}f_j\\
&=&\sum_j(\delta_{ij}v_i-\delta_{i-j})f_j\\
&=&v_if_i-\sum_{i-j}f_j
\end{eqnarray*}

Thus, we are led to the conclusions in the statement.
\end{proof}

Summarizing, we have some good reasons for calling $L$ the Laplacian, because the solutions of $Lf=0$ satisfy what we would expect from a harmonic function, namely having the ``average over the neighborhood'' property. With the remark however that the harmonic functions on graphs are something trivial, due to the following fact:

\begin{proposition}
A function on a graph $X$ is harmonic in the above sense precisely when it is constant over the connected components of $X$.
\end{proposition}

\begin{proof}
This is clear from the equation that we found in Proposition 5.13, namely:
$$f_i=\frac{1}{v_i}\sum_{i-j}f_j$$

Indeed, based on this, we can say for instance that $f$ cannot have variations over a connected component, and so must be constant on these components, as stated.
\end{proof}

At a more advanced level now, let us try to understand the relation with the usual Laplacian from analysis $\Delta$, which is given by the following formula:
$$\Delta f=\sum_{i=1}^N\frac{d^2f}{dx_i^2}$$

In one dimension, $N=1$, the Laplacian is simply the second derivative, $\Delta f=f''$. Now let us recall that the first derivative of a one-variable function is given by:
$$f'(x)=\lim_{t\to 0}\frac{f(x+t)-f(x)}{t}$$

We deduce from this, or from the Taylor formula at order 2, to be fully correct, that the second derivative of a one-variable function is given by the following formula:
\begin{eqnarray*}
f''(x)
&=&\lim_{t\to 0}\frac{f'(x+t)-f'(x)}{t}\\
&=&\lim_{t\to 0}\frac{f(x+t)-2f(x)+f(x-t)}{t^2}
\end{eqnarray*}

Now since $\mathbb R$ can be thought of as appearing as the continuum limit, $t\to0$, of the graphs $t\mathbb Z\simeq\mathbb Z$, this suggests defining the Laplacian of $\mathbb Z$ by the following formula:
$$\Delta f(x)=\frac{f(x+1)-2f(x)+f(x-1)}{1^2}$$

But this is exactly what we have in Definition 5.11, up to a sign switch, the graph Laplacian of $\mathbb Z$, as constructed there, being given by the following formula:
$$Lf(x)=2f(x)-f(x+1)-f(x-1)$$

Summarizing, we have reached to the formula in Warning 5.12, namely:
$$L=-\Delta$$

In arbitrary dimensions now, everything generalizes well, and we have:

\index{lattice model}

\begin{theorem}
The Laplacian of graphs is compatible with the usual Laplacian,
$$\Delta f=\sum_{i=1}^N\frac{d^2f}{dx_i^2}$$
via the following formula, showing that our $L$ is in fact the negative Laplacian,
$$L=-\Delta$$
via regarding $\mathbb R^N$ as the continuum limit, $t\to0$, of the graphs $t\mathbb Z^N\simeq\mathbb Z^N$.
\end{theorem}

\begin{proof}
Indeed, at $N=2$, to start with, the formula that we need, coming from standard multivariable calculus, or just from the $N=1$ formula, is as follows:
\begin{eqnarray*}
\Delta f(x,y)
&=&\frac{d^2f}{dx^2}+\frac{d^2f}{dy^2}\\
&=&\lim_{t\to 0}\frac{f(x+t,y)-2f(x,y)+f(x-t,y)}{t^2}\\
&+&\lim_{t\to 0}\frac{f(x,y+t)-2f(x,y)+f(x,y-t)}{t^2}\\
&=&\lim_{t\to 0}\frac{f(x+t,y)+f(x-t,y)+f(x,y+t)+f(x,y-t)-4f(x,y)}{t^2}
\end{eqnarray*}

Now since $\mathbb R^2$ can be thought of as appearing as the continuum limit, $t\to0$, of the graphs $t\mathbb Z^2\simeq\mathbb Z^2$, this suggests defining the Laplacian of $\mathbb Z^2$ as follows:
$$\Delta f(x)
=\frac{f(x+1,y)+f(x-1,y)+f(x,y+1)+f(x,y-1)-4f(x,y)}{1^2}$$

But this is exactly what we have in Definition 5.11, up to a sign switch, the graph Laplacian of $\mathbb Z^2$, as constructed there, being given by the following formula:
$$Lf(x)=4f(x,y)-f(x+1,y)-f(x-1,y)-f(x,y+1)-f(x,y-1)$$

At higher $N\in\mathbb N$ the proof is similar, and we will leave this as an exercise.
\end{proof}

All this is quite interesting, and suggests doing all sorts of geometric and analytic things, with our graphs and their Laplacians. For instance, we can try to review the above with $\mathbb R^N$ replaced by more general manifolds, having a certain curvature. Also, we can now do PDE over our graphs, by using the negative Laplacian constructed above.

\bigskip

But more on this later. Now back to our general graph questions, and to Definition 5.11 as it is, the Laplacian of graphs as constructed there has the following properties:

\index{bistochastic}

\begin{theorem}
The graph Laplacian $L=v-d$ has the following properties:
\begin{enumerate}
\item It is symmetric, $L=L^t$.

\item It is positive definite, $L\geq 0$.

\item It is bistochastic, with row and column sums $0$.

\item It has $0$ as eigenvalue, with the other eigenvalues being positive.

\item The multiplicity of $0$ is the number of connected components.

\item In the connected case, the eigenvalues are $0<\lambda_1\leq\ldots\leq\lambda_{N-1}$.
\end{enumerate}
\end{theorem}

\begin{proof}
All this is straightforward, the idea being as follows:

\medskip

(1) This is clear from $L=v-d$, both $v,d$ being symmetric.

\medskip

(2) This follows from the following computation, for any function $f$ on the graph:
\begin{eqnarray*}
<Lf,f>
&=&\sum_{ij}L_{ij}f_if_j\\
&=&\sum_{ij}(\delta_{ij}v_i-\delta_{i-j})f_if_j\\
&=&\sum_iv_if_if_j-\sum_{i-j}f_if_j\\
&=&\sum_{i-j}f_i^2-\sum_{i-j}f_if_j\\
&=&\frac{1}{2}\sum_{i-j}(f_i-f_j)^2\\
&\geq&0
\end{eqnarray*}

(3) This is again clear from $L=v-d$, and from the definition of $v,d$.

\medskip

(4) Here the first assertion comes from (3), and the second one, from (2).

\medskip

(5) Given an arbitrary graph, we can label its vertices inceasingly, over the connected components, and this makes the adjacency matrix $d$, so the Laplacian $L$ as well, block diagonal. Thus, we are left with proving that for a connected graph, the multiplicity of 0 is precisely 1. But this follows from the formula from the proof of (2), namely:
$$<Lf,f>=\frac{1}{2}\sum_{i-j}(f_i-f_j)^2$$

Indeed, this formula shows in particular that we have the following equivalence:
$$Lf=0\iff f_i=f_j,\forall i-j$$

Now since our graph was assumed to be connected, as per the above beginning of proof, the condition on the right means that $f$ must be constant. Thus, the 0-eigenspace of the Laplacian follows to be 1-dimensional, spanned by the all-1 vector, as desired.

\medskip

(6) This follows indeed from (4) and (5), and with the remark that in fact we already proved this, in the proof of (5), with the formulae there being very useful in practice.
\end{proof}

\section*{5c. Discrete waves}

Getting now into discretization and physics, we have the following key result:

\index{wave equation}
\index{Laplace operator}
\index{lattice model}
\index{Hooke law}
\index{Newton law}

\begin{theorem}
The wave equation in $\mathbb R^N$ is
$$\ddot{\varphi}=v^2\Delta\varphi$$
where $v>0$ is the propagation speed.
\end{theorem}

\begin{proof}
The equation in the statement is of course what comes out of physics experiments. However, allowing us a bit of imagination, and trust in this imagination, we can mathematically ``prove'' this equation, by discretizing, as follows:

\medskip

(1) Let us first consider the 1D case. In order to understand the propagation of waves, we will model $\mathbb R$ as a network of balls, with springs between them, as follows:
$$\cdots\times\!\!\!\times\!\!\!\times\bullet\times\!\!\!\times\!\!\!\times\bullet\times\!\!\!\times\!\!\!\times\bullet\times\!\!\!\times\!\!\!\times\bullet\times\!\!\!\times\!\!\!\times\bullet\times\!\!\!\times\!\!\!\times\cdots$$

Now let us send an impulse, and see how the balls will be moving. For this purpose, we zoom on one ball. The situation here is as follows, $l$ being the spring length:
$$\cdots\cdots\bullet_{\varphi(x-l)}\times\!\!\!\times\!\!\!\times\bullet_{\varphi(x)}\times\!\!\!\times\!\!\!\times\bullet_{\varphi(x+l)}\cdots\cdots$$

We have two forces acting at $x$. First is the Newton motion force, mass times acceleration, which is as follows, with $m$ being the mass of each ball:
$$F_n=m\cdot\ddot{\varphi}(x)$$

And second is the Hooke force, displacement of the spring, times spring constant. Since we have two springs at $x$, this is as follows, $k$ being the spring constant:
\begin{eqnarray*}
F_h
&=&F_h^r-F_h^l\\
&=&k(\varphi(x+l)-\varphi(x))-k(\varphi(x)-\varphi(x-l))\\
&=&k(\varphi(x+l)-2\varphi(x)+\varphi(x-l))
\end{eqnarray*}

We conclude that the equation of motion, in our model, is as follows:
$$m\cdot\ddot{\varphi}(x)=k(\varphi(x+l)-2\varphi(x)+\varphi(x-l))$$

(2) Now let us take the limit of our model, as to reach to continuum. For this purpose we will assume that our system consists of $B>>0$ balls, having a total mass $M$, and spanning a total distance $L$. Thus, our previous infinitesimal parameters are as follows, with $K$ being the spring constant of the total system, which is of course lower than $k$:
$$m=\frac{M}{B}\quad,\quad k=KB\quad,\quad l=\frac{L}{B}$$

With these changes, our equation of motion found in (1) reads:
$$\ddot{\varphi}(x)=\frac{KB^2}{M}(\varphi(x+l)-2\varphi(x)+\varphi(x-l))$$

Now observe that this equation can be written, more conveniently, as follows:
$$\ddot{\varphi}(x)=\frac{KL^2}{M}\cdot\frac{\varphi(x+l)-2\varphi(x)+\varphi(x-l)}{l^2}$$

With $N\to\infty$, and therefore $l\to0$, we obtain in this way:
$$\ddot{\varphi}(x)=\frac{KL^2}{M}\cdot\frac{d^2\varphi}{dx^2}(x)$$

We are therefore led to the wave equation in the statement, which is $\ddot{\varphi}=v^2\varphi''$ in our present $N=1$ dimensional case, the propagation speed being $v=\sqrt{K/M}\cdot L$.

\medskip

(3) In $2$ dimensions now, the same argument carries on. Indeed, we can use here a lattice model as follows, with all the edges standing for small springs:
$$
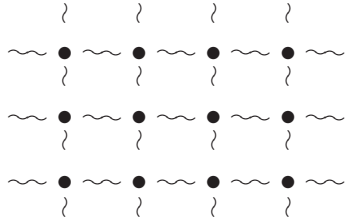
\xymatrix@R=12pt@C=15pt{
&\ar@{~}[d]&\ar@{~}[d]&\ar@{~}[d]&\ar@{~}[d]\\
\ar@{~}[r]&\bullet\ar@{~}[r]\ar@{~}[d]&\bullet\ar@{~}[r]\ar@{~}[d]&\bullet\ar@{~}[r]\ar@{~}[d]&\bullet\ar@{~}[r]\ar@{~}[d]&\\
\ar@{~}[r]&\bullet\ar@{~}[r]\ar@{~}[d]&\bullet\ar@{~}[r]\ar@{~}[d]&\bullet\ar@{~}[r]\ar@{~}[d]&\bullet\ar@{~}[r]\ar@{~}[d]&\\
\ar@{~}[r]&\bullet\ar@{~}[r]\ar@{~}[d]&\bullet\ar@{~}[r]\ar@{~}[d]&\bullet\ar@{~}[r]\ar@{~}[d]&\bullet\ar@{~}[r]\ar@{~}[d]&\\
&&&&&}$$

As before in one dimension, we send an impulse, and we zoom on one ball. The situation here is as follows, with $l$ being the spring length:
$$
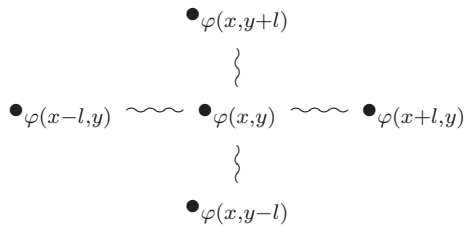
\xymatrix@R=20pt@C=20pt{
&\bullet_{\varphi(x,y+l)}\ar@{~}[d]&\\
\bullet_{\varphi(x-l,y)}\ar@{~}[r]&\bullet_{\varphi(x,y)}\ar@{~}[r]\ar@{~}[d]&\bullet_{\varphi(x+l,y)}\\
&\bullet_{\varphi(x,y-l)}}$$

We have two forces acting at $(x,y)$. First is the Newton motion force, mass times acceleration, which is as follows, with $m$ being the mass of each ball:
$$F_n=m\cdot\ddot{\varphi}(x,y)$$

And second is the Hooke force, displacement of the spring, times spring constant. Since we have four springs at $(x,y)$, this is as follows, $k$ being the spring constant:
\begin{eqnarray*}
F_h
&=&F_h^r-F_h^l+F_h^u-F_h^d\\
&=&k(\varphi(x+l,y)-\varphi(x,y))-k(\varphi(x,y)-\varphi(x-l,y))\\
&+&k(\varphi(x,y+l)-\varphi(x,y))-k(\varphi(x,y)-\varphi(x,y-l))\\
&=&k(\varphi(x+l,y)-2\varphi(x,y)+\varphi(x-l,y))\\
&+&k(\varphi(x,y+l)-2\varphi(x,y)+\varphi(x,y-l))
\end{eqnarray*}

We conclude that the equation of motion, in our model, is as follows:
\begin{eqnarray*}
m\cdot\ddot{\varphi}(x,y)
&=&k(\varphi(x+l,y)-2\varphi(x,y)+\varphi(x-l,y))\\
&+&k(\varphi(x,y+l)-2\varphi(x,y)+\varphi(x,y-l))
\end{eqnarray*}

(4) Now let us take the limit of our model, as to reach to continuum. For this purpose we will assume that our system consists of $B^2>>0$ balls, having a total mass $M$, and spanning a total area $L^2$. Thus, our previous infinitesimal parameters are as follows, with $K$ being the spring constant of the total system, taken to be equal to $k$:
$$m=\frac{M}{B^2}\quad,\quad k=K\quad,\quad l=\frac{L}{B}$$

With these changes, our equation of motion found in (3) reads:
\begin{eqnarray*}
\ddot{\varphi}(x,y)
&=&\frac{KB^2}{M}(\varphi(x+l,y)-2\varphi(x,y)+\varphi(x-l,y))\\
&+&\frac{KB^2}{M}(\varphi(x,y+l)-2\varphi(x,y)+\varphi(x,y-l))
\end{eqnarray*}

Now observe that this equation can be written, more conveniently, as follows:
\begin{eqnarray*}
\ddot{\varphi}(x,y)
&=&\frac{KL^2}{M}\times\frac{\varphi(x+l,y)-2\varphi(x,y)+\varphi(x-l,y)}{l^2}\\
&+&\frac{KL^2}{M}\times\frac{\varphi(x,y+l)-2\varphi(x,y)+\varphi(x,y-l)}{l^2}
\end{eqnarray*}

With $N\to\infty$, and therefore $l\to0$, we obtain in this way:
$$\ddot{\varphi}(x,y)=\frac{KL^2}{M}\left(\frac{d^2\varphi}{dx^2}+\frac{d^2\varphi}{dy^2}\right)(x,y)$$

Thus, we are led in this way to the following wave equation in two dimensions, with $v=\sqrt{K/M}\cdot L$ being the propagation speed of our wave:
$$\ddot{\varphi}(x,y)=v^2\left(\frac{d^2\varphi}{dx^2}+\frac{d^2\varphi}{dy^2}\right)(x,y)$$

But we recognize at right the Laplace operator, and we are done. As before in 1D, there is of course some discussion to be made here, arguing that our spring model in (3) is indeed the correct one. But do not worry, experiments confirm our findings.

\medskip

(5) In 3 dimensions now, which is the case of the main interest, corresponding to our real-life world, the same argument carries over, and the wave equation is as follows:
$$\ddot{\varphi}(x,y,z)=v^2\left(\frac{d^2\varphi}{dx^2}+\frac{d^2\varphi}{dy^2}+\frac{d^2\varphi}{dz^2}\right)(x,y,z)$$

(6) Finally, the same argument, namely a lattice model, carries on in arbitrary $N$ dimensions, and the wave equation here is as follows:
$$\ddot{\varphi}(x_1,\ldots,x_N)=v^2\sum_{i=1}^N\frac{d^2\varphi}{dx_i^2}(x_1,\ldots,x_N)$$

Thus, we are led to the conclusion in the statement.
\end{proof}

In order to reach to some further insight into our spring models above, we must get deeper into elasticity. We have here the following result:

\begin{theorem}
The wave equation can be understood as well directly, as a wave propagating through a linear elastic medium, via stress.
\end{theorem}

\begin{proof}
This is indeed something very standard, the idea being as follows:

\medskip

(1) In the 1D case, assume that we have a bar of length $L$, made of linear elastic material. The stiffness of the bar is then the following quantity, with $A$ being the cross-sectional area, and with $E$ being the Young modulus of the material:
$$K=\frac{EA}{L}$$

Now when sending a pulse, this propagates as follows, $M$ being the total mass:
$$\ddot{\varphi}=\frac{EAL}{M}\cdot\varphi''(x)$$

Bur since $V=AL$ is the volume, with $\rho=M/V$ being the density, we have:
$$\ddot{\varphi}=\frac{E}{\rho}\cdot\varphi''(x)$$

Thus, as a conclusion, the wave propagates with speed $v=\sqrt{E/\rho}$.

\medskip

(2) In two or more dimensions, the study, and final result, are similar.
\end{proof}

Getting now to mathematics, with a bit of work, we can fully solve the 1D wave equation. In order to explain this, we will need a standard calculus result, as follows:

\begin{proposition}
The derivative of a function of type
$$\varphi(x)=\int_{g(x)}^{h(x)}f(s)ds$$
is given by the formula $\varphi'(x)=f(h(x))h'(x)-f(g(x))g'(x)$.
\end{proposition}

\begin{proof}
Consider a primitive of the function that we integrate, $F'=f$. We have:
\begin{eqnarray*}
\varphi(x)
&=&\int_{g(x)}^{h(x)}f(s)ds\\
&=&\int_{g(x)}^{h(x)}F'(s)ds\\
&=&F(h(x))-F(g(x))
\end{eqnarray*}

By using now the chain rule for derivatives, we obtain from this:
\begin{eqnarray*}
\varphi'(x)
&=&F'(h(x))h'(x)-F'(g(x))g'(x)\\
&=&f(h(x))h'(x)-f(g(x))g'(x)
\end{eqnarray*}

Thus, we are led to the formula in the statement.
\end{proof}

Now back to the 1D waves, the general result here, due to d'Alembert, along with a little more, in relation with our lattice models above, is as follows:

\begin{theorem}
The solution of the 1D wave equation with initial value conditions $\varphi(x,0)=f(x)$ and $\dot{\varphi}(x,0)=g(x)$ is given by the d'Alembert formula, namely:
$$\varphi(x,t)=\frac{f(x-vt)+f(x+vt)}{2}+\frac{1}{2v}\int_{x-vt}^{x+vt}g(s)ds$$
In the context of our lattice model discretizations, what happens is more or less that the above d'Alembert integral gets computed via Riemann sums.
\end{theorem}

\begin{proof}
There are several things going on here, the idea being as follows:

\medskip

(1) Let us first check that the d'Alembert solution is indeed a solution of the wave equation $\ddot{\varphi}=v^2\varphi''$. The first time derivative is computed as follows:
$$\dot{\varphi}(x,t)=\frac{-vf'(x-vt)+vf'(x+vt)}{2}+\frac{1}{2v}(vg(x+vt)+vg(x-vt))$$

The second time derivative is computed as follows:
$$\ddot{\varphi}(x,t)=\frac{v^2f''(x-vt)+v^2f(x+vt)}{2}+\frac{vg'(x+vt)-vg'(x-vt)}{2}$$

Regarding now space derivatives, the first one is computed as follows:
$$\varphi'(x,t)=\frac{f'(x-vt)+f'(x+vt)}{2}+\frac{1}{2v}(g'(x+vt)-g'(x-vt))$$

As for the second space derivative, this is computed as follows:
$$\varphi''(x,t)=\frac{f''(x-vt)+f''(x+vt)}{2}+\frac{g''(x+vt)-g''(x-vt)}{2v}$$

Thus we have indeed $\ddot{\varphi}=v^2\varphi''$. As for the initial conditions, $\varphi(x,0)=f(x)$ is clear from our definition of $\varphi$, and $\dot{\varphi}(x,0)=g(x)$ is clear from our above formula of $\dot{\varphi}$.

\medskip

(2) Conversely now, we must show that our solution is unique, but instead of going here into abstract arguments, we will simply solve our equation, which among others will doublecheck out computations in (1). Let us make the following change of variables:
$$\xi=x-vt\quad,\quad\eta=x+vt$$

With this change of variables, which is quite tricky, mixing space and time variables, our wave equation $\ddot{\varphi}=v^2\varphi''$ reformulates in a very simple way, as follows:
$$\frac{d^2\varphi}{d\xi d\eta}=0$$

But this latter equation tells us that our new $\xi,\eta$ variables get separated, and we conclude from this that the solution must be of the following special form:
$$\varphi(x,t)=F(\xi)+G(\eta)=F(x-vt)+G(x+vt)$$

Now by taking into account the intial conditions $\varphi(x,0)=f(x)$ and $\dot{\varphi}(x,0)=g(x)$, and then integrating, we are led to the d'Alembert formula in the statement.

\medskip

(3) In regards now with our discretization questions, by using a 1D lattice model with balls and springs as before, what happens to all the above is more or less that the above d'Alembert integral gets computed via Riemann sums, in our model, as stated.
\end{proof}

In $N\geq2$ dimensions things are more complicated. We will be back to this.

\section*{5d. Discrete heat}

The general heat equation is quite similar to the one for the waves, as follows:

\index{heat equation}
\index{Laplace operator}

\begin{theorem}
Heat diffusion in $\mathbb R^N$ is described by the heat equation
$$\dot{\varphi}=\alpha\Delta\varphi$$
where $\alpha>0$ is the thermal diffusivity of the medium, and $\Delta$ is the Laplace operator.
\end{theorem}

\begin{proof}
The study here is quite similar to the study of waves, as follows:

\medskip

(1) To start with, as an intuitive explanation for the equation, since the second derivative $\varphi''$ in one dimension, or the quantity $\Delta\varphi$ in general, computes the average value of a function $\varphi$ around a point, minus the value of $\varphi$ at that point, the heat equation as formulated above tells us that the rate of change $\dot{\varphi}$ of the temperature of the material at any given point must be proportional, with proportionality factor $\alpha>0$, to the average difference of temperature between that given point and the surrounding material.

\medskip

(2) The heat equation as formulated above is of course something approximative, and several improvements can be made to it, first by incorporating a term accounting for heat radiation, and then doing several fine-tunings, depending on the material involved. But more on this later, for the moment let us focus on the heat equation above.

\medskip

(3) In relation with our modelling questions, we can recover this equation a bit as we did for the wave equation before, by using a basic lattice model. Indeed, let us first assume, for simplifying, that we are in the one-dimensional case, $N=1$. Here our model looks as follows, with distance $l>0$ between neighbors:
$$\xymatrix@R=10pt@C=20pt{
\ar@{-}[r]&\circ_{x-l}\ar@{-}[r]^l&\circ_x\ar@{-}[r]^l&\circ_{x+l}\ar@{-}[r]&
}$$

In order to model heat diffusion, we have to implement the intuitive mechanism explained above, namely ``the rate of change of the temperature of the material at any given point must be proportional, with proportionality factor $\alpha>0$, to the average difference of temperature between that given point and the surrounding material''.

\medskip

(4) In practice, this leads to a condition as follows, expressing the change of the temperature $\varphi$, over a small period of time $\delta>0$:
$$\varphi(x,t+\delta)=\varphi(x,t)+\frac{\alpha\delta}{l^2}\sum_{x\sim y}\left[\varphi(y,t)-\varphi(x,t)\right]$$

To be more precise, we have made several assumptions here, as follows:

\medskip

-- General heat diffusion assumption: the change of temperature at any given point $x$ is proportional to the average over neighbors, $y\sim x$, of the differences $\varphi(y,t)-\varphi(x,t)$ between the temperatures at $x$, and at these neighbors $y$.

\medskip

-- Infinitesimal time and length conditions: in our model, the change of temperature at a given point $x$ is proportional to small period of time involved, $\delta>0$, and is inverse proportional to the square of the distance between neighbors, $l^2$.

\medskip

(5) Regarding these latter assumptions, the one regarding the proportionality with the time elapsed $\delta>0$ is something quite natural, physically speaking, and mathematically speaking too, because we can rewrite our equation as follows, making it clear that we have here an equation regarding the rate of change of temperature at $x$:
$$\frac{\varphi(x,t+\delta)-\varphi(x,t)}{\delta}=\frac{\alpha}{l^2}\sum_{x\sim y}\left[\varphi(y,t)-\varphi(x,t)\right]$$

As for the second assumption that we made above, namely inverse proportionality with $l^2$, this can be justified on physical grounds too, but again, perhaps the best is to do the math, which will show right away where this proportionality comes from. 

\medskip

(6) So, let us do the math. In the context of our 1D model the neighbors of $x$ are the points $x\pm l$, and so the equation that we wrote above takes the following form:
$$\frac{\varphi(x,t+\delta)-\varphi(x,t)}{\delta}=\frac{\alpha}{l^2}\Big[(\varphi(x+l,t)-\varphi(x,t))+(\varphi(x-l,t)-\varphi(x,t))\Big]$$

Now observe that we can write this equation as follows:
$$\frac{\varphi(x,t+\delta)-\varphi(x,t)}{\delta}
=\alpha\cdot\frac{\varphi(x+l,t)-2\varphi(x,t)+\varphi(x-l,t)}{l^2}$$

(7) As it was the case with the wave equation before, we recognize on the right the usual approximation of the second derivative, coming from calculus. Thus, when taking the continuous limit of our model, $l\to 0$, we obtain the following equation:
$$\frac{\varphi(x,t+\delta)-\varphi(x,t)}{\delta}
=\alpha\cdot\varphi''(x,t)$$

Now with $t\to0$, we are led in this way to the heat equation, namely:
$$\dot{\varphi}(x,t)=\alpha\cdot\varphi''(x,t)$$

Summarizing, we are done with the 1D case, with our proof being quite similar to the one for the wave equation, from the previous section. 

\medskip

(8) In practice now, there are of course still a few details to be discussed, in relation with all this, for instance at the end, in relation with the precise order of the limiting operations $l\to0$ and $\delta\to0$ to be performed, but these remain minor aspects, because our equation makes it clear, right from the beginning, that time and space are separated, and so that there is no serious issue with all this. And so, fully done with 1D.

\medskip

(9) With this done, let us discuss now 2 dimensions. Here, as before for the waves, we can use a lattice model as follows, with all lengths being $l>0$, for simplifying:
$$
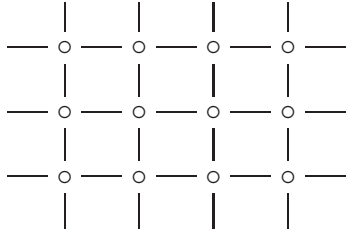
\xymatrix@R=12pt@C=15pt{
&\ar@{-}[d]&\ar@{-}[d]&\ar@{-}[d]&\ar@{-}[d]\\
\ar@{-}[r]&\circ\ar@{-}[r]\ar@{-}[d]&\circ\ar@{-}[r]\ar@{-}[d]&\circ\ar@{-}[r]\ar@{-}[d]&\circ\ar@{-}[r]\ar@{-}[d]&\\
\ar@{-}[r]&\circ\ar@{-}[r]\ar@{-}[d]&\circ\ar@{-}[r]\ar@{-}[d]&\circ\ar@{-}[r]\ar@{-}[d]&\circ\ar@{-}[r]\ar@{-}[d]&\\
\ar@{-}[r]&\circ\ar@{-}[r]\ar@{-}[d]&\circ\ar@{-}[r]\ar@{-}[d]&\circ\ar@{-}[r]\ar@{-}[d]&\circ\ar@{-}[r]\ar@{-}[d]&\\
&&&&&
}$$

(10) We have to implement now the physical heat diffusion mechanism, namely ``the rate of change of the temperature of the material at any given point must be proportional, with proportionality factor $\alpha>0$, to the average difference of temperature between that given point and the surrounding material''. In practice, this leads to a condition as follows, expressing the change of the temperature $\varphi$, over a small period of time $\delta>0$:
$$\varphi(x,y,t+\delta)=\varphi(x,y,t)+\frac{\alpha\delta}{l^2}\sum_{(x,y)\sim(u,v)}\left[\varphi(u,v,t)-\varphi(x,y,t)\right]$$

In fact, we can rewrite our equation as follows, making it clear that we have here an equation regarding the rate of change of temperature at $x$:
$$\frac{\varphi(x,y,t+\delta)-\varphi(x,y,t)}{\delta}=\frac{\alpha}{l^2}\sum_{(x,y)\sim(u,v)}\left[\varphi(u,v,t)-\varphi(x,y,t)\right]$$

(11) So, let us do the math. In the context of our 2D model the neighbors of $x$ are the points $(x\pm l,y\pm l)$, so the equation above takes the following form:
\begin{eqnarray*}
&&\frac{\varphi(x,y,t+\delta)-\varphi(x,y,t)}{\delta}\\
&=&\frac{\alpha}{l^2}\Big[(\varphi(x+l,y,t)-\varphi(x,y,t))+(\varphi(x-l,y,t)-\varphi(x,y,t))\Big]\\
&+&\frac{\alpha}{l^2}\Big[(\varphi(x,y+l,t)-\varphi(x,y,t))+(\varphi(x,y-l,t)-\varphi(x,y,t))\Big]
\end{eqnarray*}

Now observe that we can write this equation as follows:
\begin{eqnarray*}
\frac{\varphi(x,y,t+\delta)-\varphi(x,y,t)}{\delta}
&=&\alpha\cdot\frac{\varphi(x+l,y,t)-2\varphi(x,y,t)+\varphi(x-l,y,t)}{l^2}\\
&+&\alpha\cdot\frac{\varphi(x,y+l,t)-2\varphi(x,y,t)+\varphi(x,y-l,t)}{l^2}
\end{eqnarray*}

(12) As it was the case when modeling the wave equation before, we recognize on the right the usual approximation of the second derivative, coming from calculus. Thus, when taking the continuous limit of our model, $l\to 0$, we obtain the following equation:
$$\frac{\varphi(x,y,t+\delta)-\varphi(x,y,t)}{\delta}
=\alpha\left(\frac{d^2\varphi}{dx^2}+\frac{d^2\varphi}{dy^2}\right)(x,y,t)$$

Now with $t\to0$, we are led in this way to the heat equation, namely:
$$\dot{\varphi}(x,y,t)=\alpha\cdot\Delta\varphi(x,y,t)$$

Finally, in arbitrary $N$ dimensions the same argument carries over, namely a straightforward lattice model, and gives the heat equation, as formulated in the statement.
\end{proof}

Observe that we can use if we want different lenghts $l>0$ on the vertical and on the horizontal, because these will simplify anyway due to proportionality. Also, for some further mathematical fun, we can build our model on a cylinder, or a torus.

\bigskip

Also, as mentioned before, our heat equation above is something approximative, and several improvements can be made to it, first by incorporating a term accounting for heat radiation, and also by doing several fine-tunings, depending on the material involved. Some of these improvements can be implemented in the lattice model setting.

\bigskip

Regarding now the mathematics of the heat equation, many things can be said. As a first result here, often used by mathematicians, as to assume $\alpha=1$, we have:

\begin{proposition}
Up to a time rescaling, we can assume $\alpha=1$, as to deal with
$$\dot\varphi=\Delta\varphi$$
called normalized heat equation.
\end{proposition}

\begin{proof}
This is clear physically speaking, because according to our model, changing the parameter $\alpha>0$ will result in accelerating or slowing the heat diffusion, in time $t>0$. Mathematically, this follows via a change of variables, for the time variable $t$.
\end{proof}

Regarding now the resolution of the heat equation, we have here:

\index{heat kernel}

\begin{theorem}
The heat equation, normalized as $\dot\varphi=\Delta\varphi$, and with initial condition $\varphi(x,0)=f(x)$, has as solution the function
$$\varphi(x,t)=(K_t*f)(x)$$
where the function $K_t:\mathbb R^N\to\mathbb R$, called heat kernel, is given by
$$K_t(x)=(4\pi t)^{-N/2}e^{-||x||^2/4t}$$
with $||x||$ being the usual norm of vectors $x\in\mathbb R^N$.
\end{theorem}

\begin{proof}
According to the definition of the convolution operation $*$, we have to check that the following function satisfies $\dot\varphi=\Delta\varphi$, with initial condition $\varphi(x,0)=f(x)$:
$$\varphi(x,t)=(4\pi t)^{-N/2}\int_{\mathbb R^N}e^{-||x-y||^2/4t}f(y)dy$$

But both checks are elementary, coming from definitions.
\end{proof}

Regarding now our discretization questions, things here are quite tricky, and related to the Central Limit Theorem (CLT) from probability theory, which produces the normal laws, in dimension $N=1$, but also in general, in arbitrary $N\geq1$ dimensions.

\section*{5e. Exercises}

Exciting physics chapter that we had here, and as exercises, we have:

\begin{exercise}
Clarify our intuitive intepretations of $\varphi''(x)$ and $\Delta\varphi$.
\end{exercise}

\begin{exercise}
Read about harmonic functions, say from Rudin.
\end{exercise}

\begin{exercise}
Is there a better terminology for the Laplacian of graphs?
\end{exercise}

\begin{exercise}
Meditate, and read, about discretization of manifolds.
\end{exercise}

\begin{exercise}
Work out the discretization of the d'Alembert formula.
\end{exercise}

\begin{exercise}
Learn about the heat kernel, the CLT, and related mathematics.
\end{exercise}

\begin{exercise}
Have a look too into Poincar\'e, the Ricci flow, and Perelman.
\end{exercise}

\begin{exercise}
Learn about quantum waves, and the Schr\"odinger equation.
\end{exercise}

As bonus exercise, and no surprise here, learn more physics, as much as you can.

\chapter{Trees, counting}

\section*{6a. Trees, examples}

Trees, eventually. These are something beautiful, and difficult too, at the core of advanced graph theory. As starting point, we have the following definition:

\index{tree}

\begin{definition}
A tree is a connected graph with no cycles. That is, there is no loop
$$
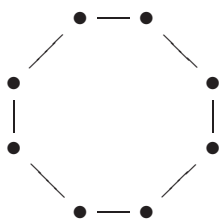
\xymatrix@R=12pt@C=12pt{
&\bullet\ar@{-}[r]\ar@{-}[dl]&\bullet\ar@{-}[dr]\\
\bullet\ar@{-}[d]&&&\bullet\ar@{-}[d]\\
\bullet\ar@{-}[dr]&&&\bullet\ar@{-}[dl]\\
&\bullet\ar@{-}[r]&\bullet}$$
having length $\geq3$, and distinct vertices, inside the graph.
\end{definition}

As a basic example here, which is quite illustrating for how trees look, in general, we have the following beautiful graph, that we already met in chapter 4:
$$
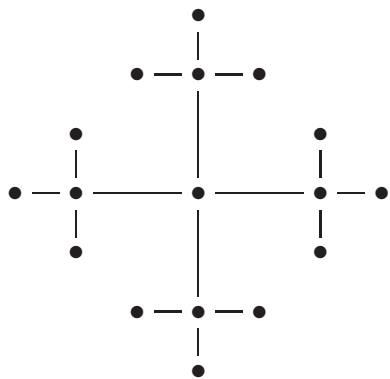
\xymatrix@R=10pt@C=10pt{
&&&\bullet\ar@{-}[d]\\
&&\bullet\ar@{-}[r]&\bullet\ar@{-}[dd]\ar@{-}[r]&\bullet\\
&\bullet\ar@{-}[d]&&&&\bullet\ar@{-}[d]\\
\bullet\ar@{-}[r]&\bullet\ar@{-}[rr]&&\bullet\ar@{-}[rr]\ar@{-}[dd]&&\bullet\ar@{-}[r]&\bullet\\
&\bullet\ar@{-}[u]&&&&\bullet\ar@{-}[u]\\
&&\bullet\ar@{-}[r]&\bullet\ar@{-}[d]\ar@{-}[r]&\bullet\\
&&&\bullet
}$$

As already mentioned in chapter 4, such a graph is called indeed a tree, because when looking at a tree from the above, what you see is something like this. 

\bigskip

Getting to work now, we can see right away, from Definition 6.1, where the difficulty in dealing with trees comes from. Indeed, with respect to what we did so far in this book, problems investigated, and solutions found for them, the situation is as follows:

\bigskip

(1) Most of our results so far were based on the correspondence between the graphs $X$ and their adjacency matrices $d$, which is something well understood, in general. However, for trees this correspondence is something quite tricky, because the loops of length $k$ are counted by the diagonal entries of $d^k$, and the condition that these entries must come only from trivial loops is obviously something complicated, that we cannot control well.

\bigskip

(2) Another thing that we did quite a lot, in the previous chapters, which was again of rather analytic type, in relation with the associated adjacency matrix $d$, was that of counting walks on our graphs. But here, for trees, things are again a bit complicated, as we have seen in chapter 1 for the simplest tree, namely $\mathbb Z$, and then in chapter 3 for the next simplest tree, namely $\mathbb N$. So, again, we are led into difficult questions here.

\bigskip

In view of this, we will take a different approach to the trees, by inventing and developing the mathematics which is best adapted to them. As a first observation, both $\mathbb N$, $\mathbb Z$ and the snowflake graph above are special kinds of trees, having a root, like normal, real-life trees do. It is possible to talk about rooted trees in general, and also about leaves, cherries, and so on, and do many other tree-inspired things. Let us record some:

\index{rooted tree}
\index{leaf}

\begin{definition}
In relation with our notion of tree:
\begin{enumerate}
\item We call rooted tree a tree with a distinguished vertex.

\item We call leaf of a tree a vertex having exactly $1$ neighbor.

\item With the convention that when the tree is rooted, the root is not a leaf.
\end{enumerate}
\end{definition}

In practice, it is often convenient to choose a root, and draw the tree oriented upwards, as the usual trees grow. Here is an example, with the root and leaves highlighted:
$$\xymatrix@R=18pt@C=12pt{
\ast&\ast&\ast&&&\ast\\
\ast&\circ\ar@{-}[ul]\ar@{-}[u]\ar@{-}[ur]&\ast&\ast&\ast&\circ\ar@{-}[u]&\ast\\
&\circ\ar@{-}[ul]\ar@{-}[u]&\ast&\circ\ar@{-}[ul]\ar@{-}[u]&&\circ\ar@{-}[ul]\ar@{-}[u]\ar@{-}[ur]\\
&&\circ\ar@{-}[ul]\ar@{-}[u]\ar@{-}[ur]&&\circ\ar@{-}[ur]\\
&&&\bullet\ar@{-}[ul]\ar@{-}[ur]
}$$

Along the same lines, simple facts that can be said about trees, we have as well:

\begin{theorem}
A tree has the following properties:
\begin{enumerate}
\item It is naturally a bipartite graph.

\item The number of edges is the number of vertices minus $1$.

\item Any two vertices are connected by a unique minimal path.
\end{enumerate}
\end{theorem}

\begin{proof}
All this is elementary, and is best seen by choosing a root, and then drawing the tree oriented upwards, as above, the idea being as follows:

\medskip

(1) We can declare the root, lying on the ground 0, to be of type ``even'', then all neighbors of the root, lying at height 1, to be of type ``odd'', then all neighbors of these neighbors, at height 2, to be of type ``even'', and so on. Thus, our tree is indeed bipartite. Here is an illustration, for how this works, for the tree pictured before the statement:
$$
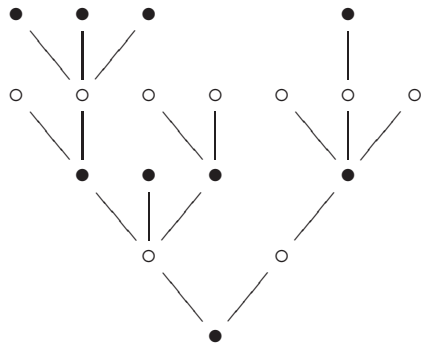
\xymatrix@R=18pt@C=12pt{
\bullet&\bullet&\bullet&&&\bullet\\
\circ&\circ\ar@{-}[ul]\ar@{-}[u]\ar@{-}[ur]&\circ&\circ&\circ&\circ\ar@{-}[u]&\circ\\
&\bullet\ar@{-}[ul]\ar@{-}[u]&\bullet&\bullet\ar@{-}[ul]\ar@{-}[u]&&\bullet\ar@{-}[ul]\ar@{-}[u]\ar@{-}[ur]\\
&&\circ\ar@{-}[ul]\ar@{-}[u]\ar@{-}[ur]&&\circ\ar@{-}[ur]\\
&&&\bullet\ar@{-}[ul]\ar@{-}[ur]
}$$

(2) This is again clear by choosing a root, and then drawing our tree oriented upwards. Indeed, in this picture any vertex, except for the root, comes from the edge below it, and this correspondence between non-root vertices and edges is bijective. Here is an illustration, for how this works, for the same example of tree as before, with each edge being represented by an arrow, pointing towards the vertex it is in bijection with:
$$
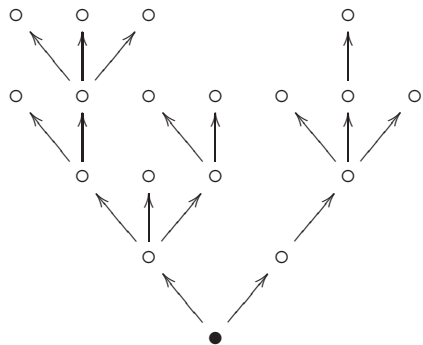
\xymatrix@R=18pt@C=12pt{
\circ&\circ&\circ&&&\circ\\
\circ&\circ\ar[ul]\ar[u]\ar[ur]&\circ&\circ&\circ&\circ\ar[u]&\circ\\
&\circ\ar[ul]\ar[u]&\circ&\circ\ar[ul]\ar[u]&&\circ\ar[ul]\ar[u]\ar[ur]\\
&&\circ\ar[ul]\ar[u]\ar[ur]&&\circ\ar[ur]\\
&&&\bullet\ar[ul]\ar[ur]
}$$

(3) Again, this is best seen by choosing a root, and then drawing our tree oriented upwards. Indeed, in this picture, traveling from one vertex to another, via the shortest path, is done uniquely, by going down as long as needed, and then up, in the obvious way. Here is an illustration for how this works, as usual for our favorite example of tree:
$$
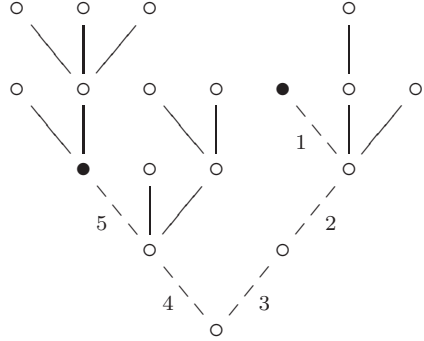
\xymatrix@R=18pt@C=12pt{
\circ&\circ&\circ&&&\circ\\
\circ&\circ\ar@{-}[ul]\ar@{-}[u]\ar@{-}[ur]&\circ&\circ&\bullet&\circ\ar@{-}[u]&\circ\\
&\bullet\ar@{-}[ul]\ar@{-}[u]&\circ&\circ\ar@{-}[ul]\ar@{-}[u]&&\circ\ar@{--}[ul]^1\ar@{-}[u]\ar@{-}[ur]\\
&&\circ\ar@{--}[ul]^5\ar@{-}[u]\ar@{-}[ur]&&\circ\ar@{--}_2[ur]\\
&&&\circ\ar@{--}[ul]^4\ar@{--}[ur]_3
}$$

Thus, we are led to the conclusions in the statement.
\end{proof}

Getting started for good now, let us look at the snowflake graph pictured before, with the convention that the graph goes on, in a 4-valent way, in each possible direction:
$$
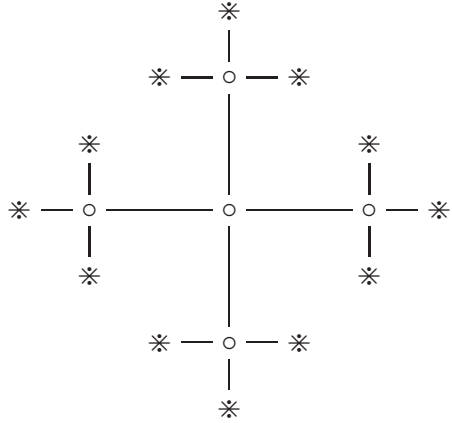
\xymatrix@R=10pt@C=10pt{
&&&\divideontimes\ar@{-}[d]\\
&&\divideontimes\ar@{-}[r]&\circ\ar@{-}[dd]\ar@{-}[r]&\divideontimes\\
&\divideontimes\ar@{-}[d]&&&&\divideontimes\ar@{-}[d]\\
\divideontimes\ar@{-}[r]&\circ\ar@{-}[rr]&&\circ\ar@{-}[rr]\ar@{-}[dd]&&\circ\ar@{-}[r]&\divideontimes\\
&\divideontimes\ar@{-}[u]&&&&\divideontimes\ar@{-}[u]\\
&&\divideontimes\ar@{-}[r]&\circ\ar@{-}[d]\ar@{-}[r]&\divideontimes\\
&&&\divideontimes
}$$

So, where does this beautiful graph come from, mathematically speaking? And a bit of thinking here leads to the following simple, and beautiful too, answer:

\index{free group}
\index{snowflake graph}

\begin{theorem}
Consider the free group $F_N$, given by the following formula, with $\emptyset$ standing for the lack of relations between generators:
$$F_N=\left<g_1,\ldots,g_N\,\Big|\,\emptyset\right>$$
The Cayley graph of $F_N$, with respect to the generating set $S=\{g_i,g_i^{-1}\}$, is then the infinite rooted tree having valence $2N$ everywhere.
\end{theorem}

\begin{proof}
This is something quite basic, the idea being as follows:

\medskip

(1) At $N=1$, our free group on one generator is by definition given by:
$$F_1=<g\,|\,\emptyset>$$

But, what does this mean? This means that $F_1$ is generated by a variable $g$ without relations of type $g^s=1$, between this variable and itself, preventing our group to be cyclic. We conclude that $F_1$ is the usual group formed by the integers:
$$F_1=\mathbb Z$$

Thus the associated Cayley graph is the usual $\mathbb Z$ graph, that we know well from the previous chapters, with edges drawn between the neighbors:
$$\xymatrix@R=10pt@C=10pt{
\cdots\ar@{-}[r]&\circ\ar@{-}[r]&\circ\ar@{-}[r]&\circ\ar@{-}[r]&\bullet\ar@{-}[r]&\circ\ar@{-}[r]&\circ\ar@{-}[r]&\circ\ar@{-}[r]&\cdots}$$

Here we have used, as before, a solid circle for the root. Now this being obviously the unique rooted tree having valence $2$ everywhere, we are done with the $N=1$ case.

\medskip

(2) At $N=2$ now, our free group on two generators is by definition given by:
$$F_2=<g,h\,|\,\emptyset>$$

So, question now, what does this mean? This means, by definition, that $F_2$ is generated by two variables $g,h$, without relations between them. In other words, our generators $g,h$ are as free as possible, in the sense that the following must be satisfied:
$$g^s\neq1\quad,\quad h^t\neq1$$
$$g^s\neq h^t$$
$$g^sh^t\neq h^tg^s$$
$$\vdots$$

Which is of course something a bit abstract, and dealing with all these conditions does not look like an easy task. However, after some thinking, this is in fact something very simple, because we can list the group elements, according to their length, as follows:
$$1$$
$$g,h,g^{-1},h^{-1}$$
$$g^2,gh,gh^{-1},hg,h^2,hg^{-1},g^{-1}h,g^{-2},g^{-1}h^{-1},h^{-1}g,h^{-1}g^{-1},h^{-2}$$
$$\vdots$$

Observe that we have one element of length 0, namely the unit 1, then 4 elements of length 1, then $16-4=12$ elements of length 2, then $64-16=48$ elements of length 3, and so on. Now when drawing the Cayley graph, the picture is as follows, with the convention that the graph goes on, in a 4-valent way, in each possible direction:
$$
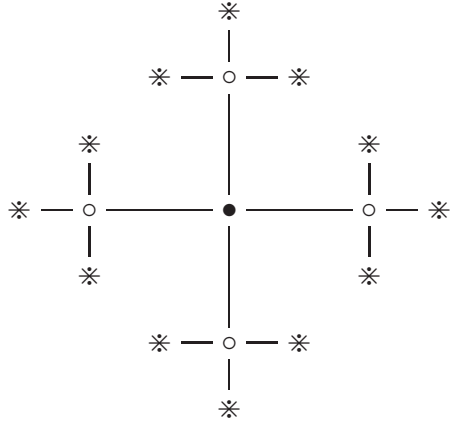
\xymatrix@R=10pt@C=10pt{
&&&\divideontimes\ar@{-}[d]\\
&&\divideontimes\ar@{-}[r]&\circ\ar@{-}[dd]\ar@{-}[r]&\divideontimes\\
&\divideontimes\ar@{-}[d]&&&&\divideontimes\ar@{-}[d]\\
\divideontimes\ar@{-}[r]&\circ\ar@{-}[rr]&&\bullet\ar@{-}[rr]\ar@{-}[dd]&&\circ\ar@{-}[r]&\divideontimes\\
&\divideontimes\ar@{-}[u]&&&&\divideontimes\ar@{-}[u]\\
&&\divideontimes\ar@{-}[r]&\circ\ar@{-}[d]\ar@{-}[r]&\divideontimes\\
&&&\divideontimes
}$$

But this graph being obviously the unique infinite rooted tree having valence $4$ everywhere, we are done with the $N=2$ case too.

\medskip

(3) At $N=3$ now, our free group on 3 generators is by definition given by:
$$F_3=<g,h,k\,|\,\emptyset>$$

As before at $N=2$, this is something a bit abstract, but we can list the group elements according to their length, and this provides us with some understanding, on what our group $F_3$ exactly is. And then, when drawing the Cayley graph, we are again in a situation which is similar to the one at $N=2$, but this time with our graph being $6$-valent everywhere, and more specifically, being the unique 6-valent rooted tree:
$$
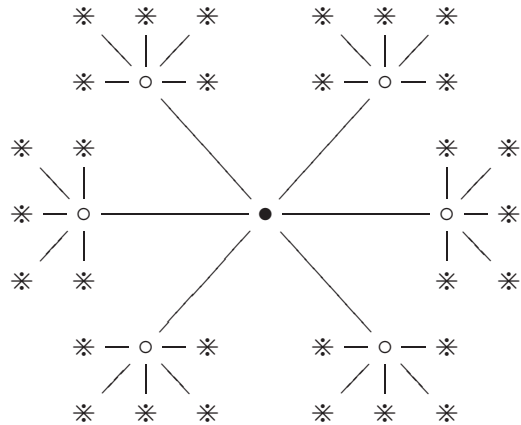
\xymatrix@R=10pt@C=7pt{
&\divideontimes\ar@{-}[dr]&\divideontimes\ar@{-}[d]&\divideontimes\ar@{-}[dl]&&\divideontimes\ar@{-}[dr]&\divideontimes\ar@{-}[d]&\divideontimes\ar@{-}[dl]\\
&\divideontimes\ar@{-}[r]&\circ\ar@{-}[ddrr]\ar@{-}[r]&\divideontimes&&\divideontimes\ar@{-}[r]&\circ\ar@{-}[ddll]\ar@{-}[r]&\divideontimes\\
\divideontimes\ar@{-}[dr]&\divideontimes\ar@{-}[d]&&&&&&\divideontimes\ar@{-}[d]&\divideontimes\ar@{-}[dl]\\
\divideontimes\ar@{-}[r]&\circ\ar@{-}[rrr]&&&\bullet\ar@{-}[rrr]\ar@{-}[ddrr]\ar@{-}[ddll]&&&\circ\ar@{-}[r]&\divideontimes\\
\divideontimes\ar@{-}[ur]&\divideontimes\ar@{-}[u]&&&&&&\divideontimes\ar@{-}[u]&\divideontimes\ar@{-}[ul]\\
&\divideontimes\ar@{-}[r]&\circ\ar@{-}[r]&\divideontimes&&\divideontimes\ar@{-}[r]&\circ\ar@{-}[r]&\divideontimes\\
&\divideontimes\ar@{-}[ur]&\divideontimes\ar@{-}[u]&\divideontimes\ar@{-}[ul]&&\divideontimes\ar@{-}[ur]&\divideontimes\ar@{-}[u]&\divideontimes\ar@{-}[ul]
}$$

(4) At an arbitrary $N\in\mathbb N$ now, the situation is quite similar to what we have seen in the above at $N=1,2,3$, and this leads to the conclusion in the statement.
\end{proof}

Many other things can be said about free groups, and other free products of groups, and their Cayley graphs. We will be back to these topics on several occasions, first at the end of the present chapter, with a discussion of the random walks on the above snowflake graphs, then in chapter 7 below, with a key occurrence of $F_N$ in topology, and then on a more systematic basis when doing groups, starting from chapter 9 below.

\section*{6b. Counting, Cayley}

Moving ahead, let us discuss now counting questions for trees. We already talked a bit about counting questions for various types of graphs before, but in the case of the trees, many interesting things can be said, that we will explore in what follows. 

\bigskip

We are mainly interested in counting the trees having $N$ vertices. However, this is something quite tricky, and it is better instead to try to count the tree with vertices labeled $1,\ldots,N$. And, regarding these latter trees, we have the following key fact:

\index{Pr\"ufer sequence}
\index{Pr\"ufer encoding}

\begin{definition}
To any tree with vertices labeled $1,\ldots,N$ we can associate its Pr\"ufer sequence, $(a_1,\ldots,a_{N-2})$ with $a_i\in\{1,\ldots,N\}$, constructed as follows:
\begin{enumerate}
\item At step $i$, remove the leaf with the smallest label,

\item Add the label of that leaf's neighbor to your list,

\item And stop when you have only $2$ vertices left.
\end{enumerate}
\end{definition}

This is something quite self-explanatory, and as an illustration for how the Pr\"ufer encoding works, consider for instance the following tree:
$$
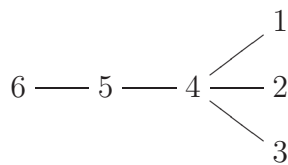
\xymatrix@R=10pt@C=20pt{
&&&1\\
6\ar@{-}[r]&5\ar@{-}[r]&4\ar@{-}[r]\ar@{-}[ur]\ar@{-}[dr]&2\\
&&&3}$$

So, let us see how the algorithm works for this tree. The details are as follows:

\bigskip

-- Following the algorithm, we first have to remove 1, and put 4 on our list. 

\bigskip

-- Then we have to remove 2, and then 3 too, again by putting 4 on our list, twice. 

\bigskip

-- Then, when left with $4,5,6$, we have to remove 4, and put 5 on our list.

\bigskip

-- And with this we are done, because we have only 2 vertices left, so we stop. 

\bigskip

As a conclusion to this, for the above tree the Pr\"ufer sequence is as follows:
$$p=(4,4,4,5)$$

So, this is how the Pr\"ufer encoding works, and feel free to work out some more examples, on your own. Quite remarkably, we have as well an inverse algorithm, as follows:

\index{Pr\"ufer decoding}

\begin{definition}
To any Pr\"ufer sequence, $(a_1,\ldots,a_{N-2})$ with $a_i\in\{1,\ldots,N\}$, we can associate a tree with vertices labeled $1,\ldots,N$, as follows:
\begin{enumerate}
\item Start with vertices $1,\ldots,N$, and declare the valence of $i\in\{1,\ldots,N\}$ to be the number of times $i$ appears in the Pr\"ufer sequence, plus $1$,

\item Read the Pr\"ufer sequence, and for any number on it $i\in\{1,\ldots,N\}$, connect $i$ to the smallest vertex $j$ having valence $1$, and lower by one the valences of $i,j$,

\item At the end, you will have two vertices left of valence $1$. Connect them.
\end{enumerate}
\end{definition}

Again, this is something quite self-explanatory, and as an illustration for how this works, consider the following Pr\"ufer sequence, that we already know from before:
$$p=(4,4,4,5)$$

So, let us see how the algorithm works for this sequence. The details are as follows:

\bigskip

-- Following the algorithm, we start with vertices, with valences, as follows:
$$1-\quad,\quad 2-\quad,\quad 3-\quad,\quad >4<\quad,\quad -5-\quad,\quad 6-$$

-- Then we read the Pr\"ufer sequence, as indicated. The first occurrence of 4 will lead to an edge $1-4$, then the second occurrence of 4 will lead to an edge $2-4$, and then the third occurrence of 4 will lead to an edge $3-4$. Thus, at this point of the algorithm, the sequence of vertices, with valences, and some edges added, will look as follows:
$${\ }^1\!\!>4<^2_3\quad,\quad -5-\quad,\quad 6-$$

-- We keep reading the Pr\"ufer sequence, by reading the last number there, the 5 at the end. This will produce an edge $4-5$, so our updated graph in the making becomes:
$${\ }^{\ 1}_{-5}\!>4<^2_3\quad,\quad 6-$$

-- Thus, bulk of the algorithm executed, and what we have now is what is indicated above, namely a graph, save for two valence 1 vertices, which are still to be connected. By connecting these two vertices, with an edge $5-6$, our graph becomes:
$${\ }^{\ \ 1}_{6-5}\!>4<^2_3$$

And good news, this is indeed a tree. And as further good news, this is in fact a tree that we know from before, having $p=(4,4,4,5)$ as Pr\"ufer sequence, namely:
$$
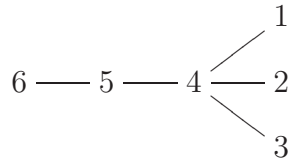
\xymatrix@R=10pt@C=20pt{
&&&1\\
6\ar@{-}[r]&5\ar@{-}[r]&4\ar@{-}[r]\ar@{-}[ur]\ar@{-}[dr]&2\\
&&&3}$$

Very nice all this, and, obviously, it looks like we are close to a big discovery. So, let us formulate now our main theorem, inspired from the above:

\begin{theorem}
The trees with vertices labeled $1,\ldots,N$ are bijectively encoded by their Pr\"ufer sequences, $(a_1,\ldots,a_{N-2})$ with $a_i\in\{1,\ldots,N\}$, constructed as follows:
\begin{enumerate}
\item At step $i$, remove the leaf with the smallest label,

\item Add the label of that leaf's neighbor to your list,

\item And stop when you have only $2$ vertices left.
\end{enumerate}
\end{theorem}

\begin{proof}
This follows from the above, and from some more thinking, as follows:

\medskip

(1) To start with, what we have in the statement is a copy of Definition 6.5, that is, of the Pr\"ufer encoding algorithm, along with claim that this algorithm is bijective.

\medskip

(2) In order to prove the bijectivity, we will use the Pr\"ufer decoding algorithm, from Definition 6.6. We have already seen in the above that these algorithms are inverse to each other, for a certain tree and sequence, and the claim is that this holds in general.

\medskip

(3) In order to verify this claim, nothing better than checking it first for some sort of ``random tree''. So, let us pick our favorite example of tree, that we used many times in the beginning of this chapter, and label its vertices in a somewhat random way:
$$
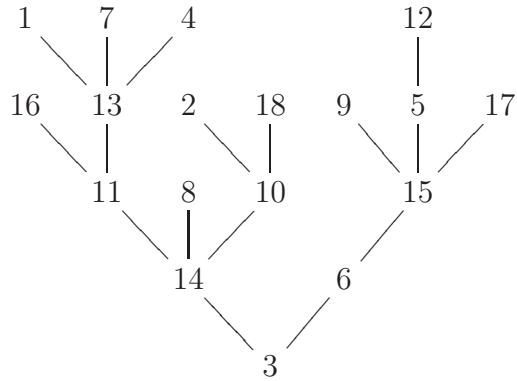
\xymatrix@R=18pt@C=12pt{
1&7&4&&&12\\
16&13\ar@{-}[ul]\ar@{-}[u]\ar@{-}[ur]&2&18&9&5\ar@{-}[u]&17\\
&11\ar@{-}[ul]\ar@{-}[u]&8&10\ar@{-}[ul]\ar@{-}[u]&&15\ar@{-}[ul]\ar@{-}[u]\ar@{-}[ur]\\
&&14\ar@{-}[ul]\ar@{-}[u]\ar@{-}[ur]&&6\ar@{-}[ur]\\
&&&3\ar@{-}[ul]\ar@{-}[ur]
}$$

(4) The associated Pr\"ufer sequence is then easy to compute, and is given by:
$$p=(13,10,13,13,14,15,5,15,11,11,14,15,6,3,14,10)$$

(5) Now let us compute the tree associated to this Pr\"ufer sequence. We follow here the algorithm in Definition 6.6, by starting with vertices, with valences, as follows:
$$1-\quad,\quad2-\quad,\quad-3-\quad,\quad4-\quad,\quad-5-\quad,\quad-6-$$
$$7-\quad,\quad8-\quad,\quad9-\quad,\quad-10<\quad,\quad-11<\quad,\quad12-$$
$$>13<\quad,\quad>14<\quad,\quad>15<\quad,\quad16-\quad,\quad17-\quad,\quad18-$$

We read now the Pr\"ufer sequence, and perform the operations in Definition 6.6. To start with, after reading the entries $13,10,13,13,14,15$, our picture becomes:
$$-3-\quad,\quad-5-\quad,\quad-6-\quad,\quad-10<^2\quad,\quad-11<\quad,\quad12-$$
$${\ }^7\!\!>13<^1_4\quad,\quad>14<^8\quad,\quad>15<^9\quad,\quad16-\quad,\quad17-\quad,\quad18-$$

Then, after reading the next entries $5,15,11$, our picture becomes:
$$-3-\quad,\quad-6-\quad,\quad-10<^2\quad,\quad{\ }^{\ \ 7}_{>11}\!\!>13<^1_4$$
$$>14<^8\quad,\quad>15<^9_{5-12}\quad,\quad16-\quad,\quad17-\quad,\quad18-$$

Then, after reading the next entry, which is an $11$, our picture becomes:
$$-3-\quad,\quad-6-\quad,\quad-10<^2\quad,\quad{\ }^{\ \ \ \ \ 7}_{{\ }_{16}>11}\!\!>13<^1_4$$
$$>14<^8\quad,\quad>15<^9_{5-12}\quad,\quad17-\quad,\quad18-$$

The next entry is a 14, and changing our graphics, our picture becomes:
$$\xymatrix@R=18pt@C=12pt{
1&7&4&&&12\\
16&13\ar@{-}[ul]\ar@{-}[u]\ar@{-}[ur]&&&9&5\ar@{-}[u]&&2&\\
&11\ar@{-}[ul]\ar@{-}[u]&8&&&15\ar@{-}[ul]\ar@{-}[u]\ar@{-}[ur]&&&10\ar@{-}[ul]\ar@{-}[u]\ar@{-}[dl]\\
&&14\ar@{-}[ul]\ar@{-}[u]\ar@{-}[ur]\ar@{-}[d]&&\ar@{-}[ur]&&&\\
&&&&&-3-&-6-&17-&18-
}$$

The next entries are $15,6,3$, and after reading them, our picture becomes:
$$\xymatrix@R=18pt@C=12pt{
1&7&4&&&12\\
16&13\ar@{-}[ul]\ar@{-}[u]\ar@{-}[ur]&&&9&5\ar@{-}[u]&17&2&\\
&11\ar@{-}[ul]\ar@{-}[u]&8&&&15\ar@{-}[ul]\ar@{-}[u]\ar@{-}[ur]&&&10\ar@{-}[ul]\ar@{-}[u]\ar@{-}[dl]\\
&&14\ar@{-}[ul]\ar@{-}[u]\ar@{-}[ur]\ar@{-}[d]&&6\ar@{-}[ur]&&&\\
&&&3\ar@{-}[ur]\ar@{-}[u]&&&18-
}$$

The last two entries are $14,10$, and after reading them, our picture becomes:
$$
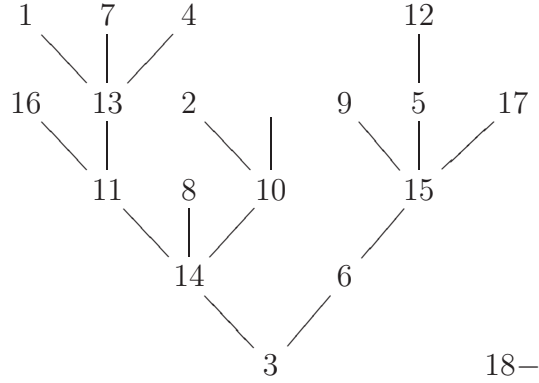
\xymatrix@R=18pt@C=12pt{
1&7&4&&&12\\
16&13\ar@{-}[ul]\ar@{-}[u]\ar@{-}[ur]&2&&9&5\ar@{-}[u]&17\\
&11\ar@{-}[ul]\ar@{-}[u]&8&10\ar@{-}[ul]\ar@{-}[u]&&15\ar@{-}[ul]\ar@{-}[u]\ar@{-}[ur]\\
&&14\ar@{-}[ul]\ar@{-}[u]\ar@{-}[ur]&&6\ar@{-}[ur]\\
&&&3\ar@{-}[ul]\ar@{-}[ur]&&&18-
}$$

It remains to connect the remaining 2 vertices, and we get what we want, namely: 
$$
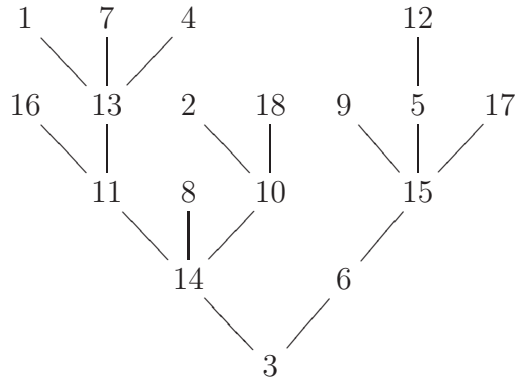
\xymatrix@R=18pt@C=12pt{
1&7&4&&&12\\
16&13\ar@{-}[ul]\ar@{-}[u]\ar@{-}[ur]&2&18&9&5\ar@{-}[u]&17\\
&11\ar@{-}[ul]\ar@{-}[u]&8&10\ar@{-}[ul]\ar@{-}[u]&&15\ar@{-}[ul]\ar@{-}[u]\ar@{-}[ur]\\
&&14\ar@{-}[ul]\ar@{-}[u]\ar@{-}[ur]&&6\ar@{-}[ur]\\
&&&3\ar@{-}[ul]\ar@{-}[ur]
}$$

(6) Summarizing, good work that we did here, and with our theorem proved for this beast, I would say that the confidence rate for our theorem jumps to a hefty $99\%$. As for the remaining $1\%$, I will leave this to you, as an instructive exercise. And with the comment that we will come back to this, later, with some alternative methods too.
\end{proof}

As a first consequence of Theorem 6.7, we have the following famous formula:

\index{Cayley formula}

\begin{theorem}
The number of trees with vertices labeled $1,\ldots,N$ is
$$T_N=N^{N-2}$$
called Cayley formula.
\end{theorem}

\begin{proof}
This is clear from the bijection in Theorem 6.7, because the number of possible Pr\"ufer sequences, $(a_1,\ldots,a_{N-2})$ with $a_i\in\{1,\ldots,N\}$, that is, sequences obtained by picking $N-2$ times numbers from $\{1,\ldots,N\}$, is obviously $T_N=N^{N-2}$.

\medskip

As a number of comments, however, which have to be made here, we have:

\medskip

(1) First of all, in order to avoid confusions, what we just counted are trees up to obvious isomorphism. For instance at $N=2$ we only have 1 tree, namely:
$$1-2$$

Indeed, the other possible tree, namely $2-1$, is clearly isomorphic to it.

\medskip

(2) At $N=3$ now, the trees are sequences of type $a-b-c$, with $\{a,b,c\}=\{1,2,3\}$. Now since such a tree is uniquely determined, up to graph isomorphism, by the middle vertex $b$, we have 3 exactly trees up to isomorphism, namely:
$$1-2-3\quad,\quad 
2-3-1\quad,\quad 
3-1-2$$

(3) At $N=4$, we first have $2\binom{4}{2}=12$ bivalent trees, which can be listed by choosing the ordered outer vertices, and putting the smallest label left at position $\#2$:
$$1-3-4-2\quad,\quad 
1-2-4-3\quad,\quad 
1-2-3-4$$
$$2-3-4-1\quad,\quad 
2-1-4-3\quad,\quad 
2-1-3-4$$
$$3-2-4-1\quad,\quad 
3-1-4-2\quad,\quad 
3-1-2-4$$
$$4-2-3-1\quad,\quad 
4-1-3-2\quad,\quad 
4-1-2-3$$

And then, for completing the count, we have $\binom{4}{1}=4$ trees with a trivalent vertex, which again, by using some sort of lexicographic order, look as follows:
$$2-1<^3_4\quad,\quad 1-2<^3_4\quad,\quad 1-3<^2_4\quad,\quad 1-4<^2_3$$

(4) And so on, and it is actually instructive here to try yourself listing the $5^3=125$ trees at $N=5$, in order to convince you that the Cayley formula is something quite subtle, hiding inside plenty of binomials and factorials.

\medskip

(5) Which leads us to the question, what is the simplest proof of the Cayley formula. Well, there are many possible proofs here, for all tastes, with the above one, using Pr\"ufer sequences, being probably the simplest one. We will present as well a second proof, a bit later, due to Kirchoff, which is based on linear algebra, and is non-trivial as well.
\end{proof}

Moving ahead, as a more specialized consequence of Theorem 6.7, we have:

\index{multinomial coefficient}

\begin{theorem}
The number of trees with vertices labeled $1,\ldots,N$, having respective valences $v_1,\ldots,v_N$, is the multinomial coefficient
$$T_{v_1,\ldots,v_N}=\binom{N-2}{v_1-1,\ldots,v_N-1}$$
and with this implying as well the Cayley formula, $T_N=N^{N-2}$.
\end{theorem}

\begin{proof}
This follows again from Theorem 6.7, the idea being as follows:

\medskip

(1) As a first observation, the formula in the statement makes sense indeed, with what we have there being indeed a multinomial coefficient, and this because by using the fact that the number of edges is $E=N-1$, that we know from Proposition 6.3, we have:
\begin{eqnarray*}
\sum_i(v_i-1)
&=&\left(\sum_iv_i\right)-N\\
&=&2E-N\\
&=&2(N-1)-N\\
&=&N-2
\end{eqnarray*}

(2) In what regards now the proof, this follows from Theorem 6.7, the point being that in the Pr\"ufer sequence, each number $i\in\{1,\ldots,N\}$ appears exactly $v_i-1$ times. 

\medskip

(3) As for the last assertion, this comes from this, and the multinomial formula:
\begin{eqnarray*}
T_N
&=&\sum_{v_1,\ldots,v_N}T_{v_1,\ldots,v_N}\\
&=&\sum_{v_1,\ldots,v_N}\binom{N-2}{v_1-1,\ldots,v_N-1}\\
&=&\sum_{v_1,\ldots,v_N}\binom{N-2}{v_1-1,\ldots,v_N-1}\times1^{v_1-1}\ldots1^{v_N-1}\\
&=&(1+\ldots+1)^{N-2}\\
&=&N^{N-2}
\end{eqnarray*}

Thus, we are led to the conclusions in the statement.
\end{proof}

\section*{6c. Kirchoff formula}

We discuss in this section another proof of the Cayley formula, this time due to Kirchoff, by using linear algebra techniques, of analytic flavor. As a bonus, we will see that the main result of Kirchoff goes well beyond what the Cayley formula says. 

\bigskip

The Kirchoff approach is based on the following simple observations:

\index{spanning tree}

\begin{proposition}
The following happen:
\begin{enumerate}
\item Any connected graph has a spanning tree, meaning a tree subgraph, making use of all vertices. 

\item For the complete graph $K_N$, with vertices labeled $1,\ldots,N$, the spanning trees are exactly the trees with vertices labeled $1,\ldots,N$.
\end{enumerate}
\end{proposition}

\begin{proof}
Both the assertions are trivial, the idea being as follows:

\medskip

(1) The fact that any connected graph has indeed a spanning tree is something which is very intuitive, clear on pictures, and we will leave the formal proof, which is not difficult, as an exercise. As an illustration for this, here is a picture of a quite random graph, which, after removal of some of the edges, the dotted ones, becomes indeed a tree:
$$
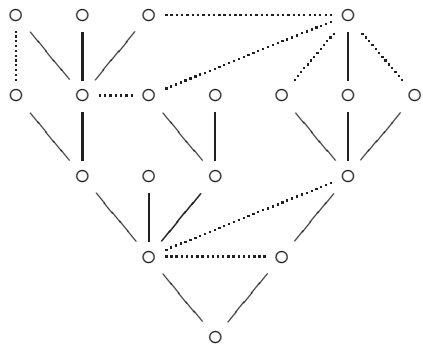
\xymatrix@R=18pt@C=12pt{
\circ\ar@{.}[d]&\circ&\circ\ar@{.}[rrr]&&&\circ\ar@{.}[dr]\ar@{.}[dlll]\ar@{.}[dl]\\
\circ&\circ\ar@{-}[ul]\ar@{-}[u]\ar@{-}[ur]\ar@{.}[r]&\circ&\circ&\circ&\circ\ar@{-}[u]&\circ\\
&\circ\ar@{-}[ul]\ar@{-}[u]&\circ&\circ\ar@{-}[ul]\ar@{-}[u]&&\circ\ar@{-}[ul]\ar@{-}[u]\ar@{-}[ur]\\
&&\circ\ar@{-}[ul]\ar@{-}[u]\ar@{-}[ur]\ar@{.}[rr]\ar@{-}[ur]\ar@{.}[urrr]&&\circ\ar@{-}[ur]\\
&&&\circ\ar@{-}[ul]\ar@{-}[ur]
}$$

(2) As for the second assertion, this is something which is clear too, and again we will leave the formal proof, which is not difficult at all, as an exercise.
\end{proof}

In view of the above, the following interesting question appears:

\begin{question}
Given a connected graph $X$, with vertices labeled $1,\ldots,N$, how to count its spanning trees? And, for the complete graph $K_N$, do we really get $N^{N-2}$ such spanning trees, in agreement with the Cayley formula, by using this method?
\end{question}

So, hope you get the point, we are trying here to do what was advertised in the above, namely new proof of the Cayley formula, plus more. Getting to work now, following Kirchoff, the idea will be that of connecting the spanning trees of a connected graph $X$ to the combinatorics of the Laplacian of $X$, given by the following formula:
$$L=v-d$$

We have already seen many things in chapter 5, regarding $L$, and in order to get started, we will just need the fact, that we know well, and which is trivial, that $L$ is bistochastic, with zero row and column sums. Indeed, this makes the link with the following basic linear algebra fact, that we can use afterwards, for our Laplacian $L$:

\index{bistochastic matrix}
\index{minors}

\begin{proposition}
For a matrix $L\in M_N(\mathbb R)$ which is bistochastic, with zero row and column sums, the signed minors
$$T_{ij}=(-1)^{i+j}\det(L^{ij})$$
do not depend on the choice of the indices $i,j$.
\end{proposition}

\begin{proof}
This is something very standard, the idea being as follows:

\medskip

(1) Before anything, let us do a quick check at $N=2$. Here the bistochastic matrices, with zero row and column sums, are as follows, with $a\in\mathbb R$:
$$L=\begin{pmatrix}a&-a\\ -a&a\end{pmatrix}$$

But this gives the result, with the number in question being $T_{ij}=a$.

\medskip

(2) Let us do as well a quick check at $N=3$. Here the bistochastic matrices, with zero row and column sums, are as follows, with $a,b,c,d\in\mathbb R$:
$$L=\begin{pmatrix}
a&b&-a-b\\
c&d&-c-d\\
-a-c&-b-d&a+b+c+d
\end{pmatrix}$$

But this gives again the result, with the number in question being $T_{ij}=ad-bc$.

\medskip

(3) In the general case now, where $N\in\mathbb N$ is arbitrary, the bistochastic matrices with zero row and column sums are as follows, with $A\in M_n(\mathbb R)$ with $n=N-1$ being an arbitary matrix, and with $R_1,\ldots,R_n$ and $C_1,\ldots,C_n$ being the row and column sums of this matrix, and $S=\sum R_i=\sum C_i$ being the total sum of this matrix:
$$L=\begin{pmatrix}
A_{11}&\ldots&A_{1n}&-R_1\\
\vdots&&\vdots&\vdots\\
A_{n1}&\ldots&A_{nn}&-R_n\\
-C_1&\ldots&-C_n&S
\end{pmatrix}$$

We want to prove that the signed minors of $L$ coincide, and by using the symmetries of the problem, it is enough to prove that the following equality holds:
$$L^{n+1,n}=-L^{n+1,n+1}$$

But, what we have on the right is $-\det A$, and what we have on the left is:
\begin{eqnarray*}
L^{n+1,n}
&=&\begin{vmatrix}
A_{11}&\ldots&A_{1,n-1}&-R_1\\
\vdots&&\vdots&\vdots\\
A_{n1}&\ldots&A_{n,n-1}&-R_n
\end{vmatrix}\\
&=&\begin{vmatrix}
A_{11}&\ldots&A_{1,n-1}&A_{11}+\ldots+A_{1,n-1}-R_1\\
\vdots&&\vdots&\vdots\\
A_{n1}&\ldots&A_{n,n-1}&A_{n1}+\ldots+A_{n,n-1}-R_n
\end{vmatrix}\\
&=&\begin{vmatrix}
A_{11}&\ldots&A_{1,n-1}&-A_{1n}\\
\vdots&&\vdots&\vdots\\
A_{n1}&\ldots&A_{n,n-1}&-A_{nn}
\end{vmatrix}\\
&=&-\det A
\end{eqnarray*}

Thus, we are led to the conclusion in the statement.
\end{proof}

We can now formulate, following Kirchoff, the following key result:

\index{Kirchoff formula}

\begin{theorem}
Given a connected graph $X$, with vertices labeled $1,\ldots,N$, the number of spanning trees inside $X$, meaning tree subgraphs using all vertices, is
$$T_X=(-1)^{i+j}\det(L^{ij})$$
with $L=v-d$ being the Laplacian, with this being independent on the chosen minor.
\end{theorem}

\begin{proof}
This is something non-trivial, the idea being as follows:

\medskip

(1) We know from Proposition 6.12 that the signed minors of $L$ coincide. In other words, we have a common formula as follows, with $T\in\mathbb Z$ being a certain number:
$$(-1)^{i+j}\det(L^{ij})=T$$

Our claim, which will prove the result, is that the number of spanning trees $T_X$ is precisely this common number $T$. That is, with $i=j=1$, our claim is that we have:
$$T_X=\det(L^{11})$$

(2) In order to prove our claim, which is non-trivial, we use a trick. We orient all the edges $e=(ij)$ of our graph as to have $i<j$, and we define the ordered incidence matrix of our graph, which is a rectangular matrix, with the vertices $i$ as row indices, and the oriented edges $e=(ij)$ as column indices, by the following formula:
$$E_{ie}=\begin{cases}
1&{\rm if}\ e=(ij)\\
-1&{\rm if}\ e=(ji)\\
0&{\rm otherwise}
\end{cases}$$

The point is that, in terms of this matrix, the Laplacian decomposes as follows:
$$L=EE^t$$

(3) Indeed, let us compute the matrix on the right. We have, by definition:
$$(EE^t)_{ij}=\sum_eE_{ie}E_{je}$$

Let us first compute the contributions of type $1\times1$, to the above sum. These come from the edges $e$ having the property $E_{ie}=E_{je}=1$. But $E_{ie}=1$ means $e=(ik)$ with $i<k$, and $E_{je}=1$ means $e=(jl)$ with $j<l$. Thus, our condition $E_{ie}=E_{je}=1$ means $i=j$, and $e=(ik)$ with $i<k$, so the contributions of type $1\times1$ are given by:
$$C_{1\times1}=\delta_{ij}\#\left\{k\Big|i<k,\ i-k\right\}$$

Similarly, the contributions of type $(-1)\times(-1)$ to our sum come from the equations $E_{ie}=E_{je}=-1$, which read $i=j$ and $e=(ki)$ with $k<i$, so are given by:
$$C_{(-1)\times(-1)}=\delta_{ij}\#\left\{k\Big|k<i,\ i-k\right\}$$

Now observe that by summing, the total $1$ contributions to our sum, be them of type $1\times1$ or $(-1)\times(-1)$, are given by the following formula, $v$ being the valence function:
\begin{eqnarray*}
C_1
&=&C_{1\times1}+C_{(-1)\times(-1)}\\
&=&\delta_{ij}\#\left\{k\Big|i<k,\ i-k\right\}+\delta_{ij}\#\left\{k\Big|k<i,\ i-k\right\}\\
&=&\delta_{ij}\#\left\{k\Big|i-k\right\}\\
&=&\delta_{ij}v_i
\end{eqnarray*}

(4) It remains to compute the total $-1$ contributions to our sum. But here, we first have the contributions of type $1\times(-1)$, coming from the equations $E_{ie}=1,E_{je}=-1$. Now since $E_{ie}=1$ means $e=(ik)$ with $i<k$, and $E_{je}=-1$ means $e=(lj)$ with $l<j$, our equations $E_{ie}=1,E_{je}=-1$ amount in saying that $e=(ij)$ with $i<j$. We conclude that the contributions of type $(-1)\times1$ to our sum are given by:
$$C_{1\times(-1)}=\delta_{i-j}\delta_{i<j}$$

Similarly, the contributions of type $(-1)\times1$ to our sum come from the equations $E_{ie}=-1,E_{je}=1$, which read $e=(ij)$ with $i<j$, so these are given by:
$$C_{(-1)\times1}=\delta_{i-j}\delta_{i>j}$$

Now by summing, the total $-1$ contributions to our sum, be them of type $1\times(-1)$ or $(-1)\times1$, are given by the following formula, $d$ being the adjacency matrix:
\begin{eqnarray*}
C_{-1}
&=&C_{1\times(-1)}+C_{(-1)\times1}\\
&=&\delta_{i-j}\delta_{i<j}+\delta_{i-j}\delta_{i>j}\\
&=&\delta_{i-j}\\
&=&d_{ij}
\end{eqnarray*}

(5) But with this, we can now finish the proof of our claim in (2), as follows:
\begin{eqnarray*}
(EE^t)_{ij}
&=&\sum_eE_{ie}E_{je}\\
&=&C_1-C_{-1}\\
&=&\delta_{ij}v_i-d_{ij}\\
&=&(v-d)_{ij}\\
&=&L_{ij}
\end{eqnarray*}

Thus, we have $EE^t=L$, and claim proved. Note in passing that our formula $EE^t=L$ gives another proof of the property $L\geq0$, that we know from chapter 5.

\medskip

(6) Getting now towards minors, if we denote by $F$ the submatrix of $E$ obtained by deleting the first row, the one coming from the vertex 1, we have, for any $i,j>1$:
\begin{eqnarray*}
(FF^t)_{ij}
&=&\sum_eF_{ie}F_{je}\\
&=&\sum_eE_{ie}E_{je}\\
&=&(EE^t)_{ij}\\
&=&L_{ij}\\
&=&(L^{11})_{ij}
\end{eqnarray*}

We conclude that we have the following equality of matrices:
$$L^{11}=FF^t$$

(7) The point now is that, in order to compute the determinant of this latter matrix, we can use the Cauchy-Binet formula from linear algebra. To be more precise, Cauchy-Binet says that given rectangular matrices $A,B$, of respective sizes $M\times N$ and $N\times M$, we have the following formula, with $A_S,B_S$ being both $M\times M$ matrices, obtained from $A,B$ by cutting, in the obvious way, with respect to the set of indices $S$:
$$\det(AB)=\sum_{|S|=M}\det(A_S)\det(B_S)$$

Observe that this formula tells us at $M>N$ that we have $\det(AB)=0$, as it should be, and at $M=N$ that we have $\det(AB)=\det A\det B$, again as it should be. At $M<N$, which is the interesting case, the Cauchy-Binet formula holds indeed, with the proof being a bit similar to that of the formula $\det(AB)=\det A\det B$ for the square matrices, which itself is not exactly a trivial business. For more on all this, details of the proof, examples, and other comments, we refer to any advanced linear algebra book.

\medskip

(8) Now back to our questions, in the context of our formula $L^{11}=FF^t$ from (6), we can apply Cauchy-Binet to the matrices $A=F$ and $B=F^t$, having respective sizes $(N-1)\times N$ and $N\times(N-1)$. We are led in this way to the following formula, with $S$ ranging over the subsets of the edge set having size $N-1$, and with $F_S$ being the corresponding square submatrix of $E$, having size $(N-1)\times(N-1)$, obtained by restricting the attention to the columns indexed by the subset $S$:
\begin{eqnarray*}
\det(L^{11})
&=&\det(FF^t)\\
&=&\sum_S\det(F_S)\det(F_S^t)\\
&=&\sum_S\det(F_S)^2
\end{eqnarray*}

(9) Now comes the combinatorics. The sets $S$ appearing in the above computation specify in fact $N-1$ edges of our graph, and so specify a certain subgraph $X_S$. But, in this picture, our claim is that we have the following formula:
$$\det(F_S)=\begin{cases}
\pm1&{\rm if}\ X_S\ {\rm is\ a\ spanning\ tree}\\
0&{\rm otherwise}
\end{cases}$$

Indeed, since the subgraph $X_S$ has $N$ vertices and $N-1$ edges, it must be either a spanning tree, or have a cycle, and the study here goes as follows:

\medskip

-- In the case where $X_S$ is a spanning tree, we pick a leaf of this tree, in theory I mean, by leaving it there, on the tree. The corresponding row of $F_S$ consists then of a $\pm1$ entry, at the neighbor of the leaf, and of $0$ entries elsewhere. Thus, by developing $\det(F_S)$ over that row, we are led to a recurrence, which gives $\det(F_S)=\pm1$, as claimed above.

\medskip

-- In the other case, where $X_S$ has a cycle, the sum of the columns of $F_S$ indexed by the vertices belonging to this cycle must be $0$. We conclude that in this case we have $\det(F_S)=0$, again as claimed above, and this finishes the proof of our claim.

\medskip

(10) By putting now everything together, we obtain the following formula:
$$\det(L^{11})=T_X$$

Thus, we are led to the conclusions in the statement.
\end{proof}

As a basic application of the Kirchoff formula, let us apply it to the complete graph $K_N$. We are led in this way to another proof of the Cayley formula, as follows:

\begin{theorem}
The number of spanning trees of the complete graph $K_N$ is
$$T_{K_N}=N^{N-2}$$
in agreement with the Cayley formula.
\end{theorem}

\begin{proof}
This is something which is clear from the Kirchoff formula, but let us prove this slowly, as an illustration for the various computations above:

\medskip

(1) At $N=2$ the Laplacian of the segment $K_2$ is given by the following fomula:
$$L=\begin{pmatrix}1&-1\\-1&1\end{pmatrix}$$

Thus the common cofactor is 1, which equals the number of spanning trees, $2^0=1$.

\medskip

(2) At $N=3$ the Laplacian of the triangle $K_3$ is given by the following fomula:
$$L=\begin{pmatrix}
2&-1&-1\\
-1&2&-1\\
-1&-1&2
\end{pmatrix}$$

Thus the common cofactor is 3, which equals the number of spanning trees, $3^1=3$.

\medskip

(3) At $N=4$ the Laplacian of the tetrahedron $K_4$ is given by the following fomula:
$$L=\begin{pmatrix}
3&-1&-1&-1\\
-1&3&-1&-1\\
-1&-1&3&-1\\
-1&-1&-1&3
\end{pmatrix}$$

Here the cofactor is $27-11=16$, which is the number of spanning trees, $4^2=16$.

\medskip

(4) In general, for the complete graph $K_N$, the Laplacian is as follows:
$$L=\begin{pmatrix}
N-1&-1&\ldots&-1&-1\\
-1&N-1&\ldots&-1&-1\\
\vdots&\vdots&&\vdots&\vdots\\
-1&-1&\ldots&N-1&-1\\
-1&-1&\ldots&-1&N-1
\end{pmatrix}$$

Thus, the common cofactor is $N^{N-2}$, in agreement with the Cayley formula.
\end{proof}

Very nice all this, so we have now a good understanding of the Cayley formula. However, the story is not over here, because in what regards the counting of trees having $N$ vertices, this time without labeled vertices, things are far more complicated, and there is no formula available, for the number of such trees. We refer here to the literature.

\section*{6d. More Kirchoff}

Our purpose now will be to further build on Theorem 6.13, with an analytic interpretation of the quantity $T_X=(-1)^{i+j}\det(L^{ij})$ found there. So, let us first go back to the general setting of Proposition 6.12. As a continuation of the material there, we have: 

\begin{proposition}
For a matrix $L\in M_N(\mathbb R)$ which is bistochastic, with zero row and column sums, we have, independently of the chosen indices $i,j$, 
$$(-1)^{i+j}\det(L^{ij})=\frac{\lambda_1\ldots\lambda_{N-1}}{N}$$
with $\{0,\lambda_1,\ldots,\lambda_{N-1}\}$ being the eigenvalues of $L$.
\end{proposition}

\begin{proof}
This is something quite standard, the idea being as follows:

\medskip

(1) To start with, since our matrix $L$ is bistochastic, with zero sums, 0 is indeed an eigenvalue of it, with the all-one vector as eigenvector, so our statement makes indeed sense. Observe also that, according to Proposition 6.12, the numbers $(-1)^{i+j}\det(L^{ij})$ in the statement do not depend on the choice of the chosen indices $i,j$.

\medskip

(2) In order to see what is going on, let us do a quick check at $N=2$. Here the bistochastic matrices, with zero row and column sums, are as follows, with $a\in\mathbb R$:
$$L=\begin{pmatrix}a&-a\\ -a&a\end{pmatrix}$$

The eigenvalues are then $0,2a$, and the formula in the statement holds indeed, as:
$$a=\frac{2a}{2}$$

(3) Let us do as well a quick check at $N=3$. Here the bistochastic matrices, with zero row and column sums, are as follows, with $a,b,c,d\in\mathbb R$:
$$L=\begin{pmatrix}
a&b&-a-b\\
c&d&-c-d\\
-a-c&-b-d&a+b+c+d
\end{pmatrix}$$

Thus, good news, time for some Sarrus, which gives, after 5 minutes or so:
\begin{eqnarray*}
\det(x-L)
&=&\begin{vmatrix}
x-a&-b&a+b\\
-c&x-d&c+d\\
a+c&b+d&x-a-b-c-d
\end{vmatrix}\\
&=&x^3-(2a+2d+b+c)x^2+3(ad-bc)x
\end{eqnarray*}

But $\det(x-L)=x(x-\lambda_1)(x-\lambda_2)$, so our formula holds again, as:
$$ad-bc=\frac{3(ad-bc)}{3}$$

(4) In the general case now, where $N\in\mathbb N$ is arbitrary, let us write:
$$\det(x-L)=x(x-\lambda_1)\ldots(x-\lambda_{N-1})$$

Based on this, we have the following computation, obtained by fully developing the determinant, and by using Proposition 6.12 at the end: 
\begin{eqnarray*}
\lambda_1\ldots\lambda_{N-1}
&=&(-1)^{N-1}\#\big(x\in \det(x-L)\big)\\
&=&(-1)^{N-1}\big(\det(-L^{11})+\ldots+\det(-L^{NN})\big)\\
&=&\det(L^{11})+\ldots+\det(L^{NN})\\
&=&N\det(L^{11})
\end{eqnarray*}

Thus, we are led to the formula in the statement.
\end{proof}

The above result is quite interesting, suggesting to reformulate the Kirchoff formula from Theorem 6.13 in a more analytic way, by using the eigenvalues of $L$, as follows:

\index{Kirchoff formula}

\begin{theorem}
Given a connected graph $X$, with vertices labeled $1,\ldots,N$, the number of spanning trees inside $X$ is
$$T_X=\frac{\lambda_1\ldots\lambda_{N-1}}{N}$$
with $\lambda_1,\ldots,\lambda_{N-1}$ being the nonzero eigenvalues of the Laplacian 
$L=v-d$.
\end{theorem}

\begin{proof}
This follows indeed by putting together Theorem 6.13 and Proposition 6.15, and with the remark that, as explained in chapter 5, for a connected graph the multiplicity of the 0 eigenvalue is 1, and so we have $\lambda_1,\ldots,\lambda_{N-1}>0$, as stated.
\end{proof}

As before with our previous version of the Kirchoff theorem, we can apply this to the complete graph, and we are led to yet another proof of the Cayley formula, as follows:

\begin{theorem}
The number of spanning trees of the complete graph $K_N$ is
$$T_{K_N}=N^{N-2}$$
in agreement with the Cayley formula.
\end{theorem}

\begin{proof}
We know from chapter 2 that the adjacency matrix of $K_N$ diagonalizes as follows, with $F_N=(w^{ij})$ with $w=e^{2\pi i/N}$ being as usual the Fourier matrix:
$$\begin{pmatrix}
0&1&\ldots&1&1\\
1&0&\ldots&1&1\\
\vdots&\vdots&&\vdots&\vdots\\
1&1&\ldots&0&1\\
1&1&\ldots&1&0
\end{pmatrix}
=\frac{1}{N}\,F_N
\begin{pmatrix}
N-1&&&&0\\
&-1\\
&&\ddots&\\
&&&-1\\
0&&&&-1
\end{pmatrix}F_N^*$$

We deduce from this that the corresponding Laplacian matrix, which is given by the formula $L=(N-1)1_N-d$, diagonalizes as follows:
$$\begin{pmatrix}
N-1&-1&\ldots&-1&-1\\
-1&N-1&\ldots&-1&-1\\
\vdots&\vdots&&\vdots&\vdots\\
-1&-1&\ldots&N-1&-1\\
-1&-1&\ldots&-1&N-1
\end{pmatrix}
=\frac{1}{N}\,F_N
\begin{pmatrix}
0&&&&0\\
&N\\
&&\ddots&\\
&&&N\\
0&&&&N
\end{pmatrix}F_N^*$$

Now by applying Theorem 6.16, the number of spanning trees follows to be:
$$T_{K_N}=\frac{N^{N-1}}{N}=N^{N-2}$$

Thus, we are led to the conclusions in the statement.
\end{proof}

In order to have now more examples, we have to browse through our diagonalization results from chapter 2, and afterwards, and we are led in this way into:

\begin{theorem}
For a connected circulant graph with vertices labeled $1,\ldots,N$, and having valence $k$, with the corresponding adjacency matrix diagonalized as
$$d=\frac{1}{N}F_N\begin{pmatrix}
k&&&&0\\
&\mu_1\\
&&\ddots&\\
0&&&&\mu_{N-1}
\end{pmatrix}F_N^*$$ 
with $\mu_1,\ldots,\mu_{N-1}\in\mathbb C$, the number of spanning trees inside $X$ is given by:
$$T_X=\frac{(k-\mu_1)\ldots(k-\mu_{N-1})}{N}$$
In the case of the complete graph $K_N$, we recover in this way the Cayley formula.
\end{theorem}

\begin{proof}
Here the formulation of the statement heavily relies on our knowledge of circulant graphs, from chapter 4. As for the result itself, this follows from Theorem 6.16, by using $L=(k1_n-d)$, and the diagonalization formula for $d$ in the statement. 
\end{proof}

As an illustration for the above result, we have:

\begin{corollary}
The number of spanning trees of the circle is $N$.
\end{corollary}

\begin{proof}
As a die-hard physicist, I never use Lemmas and Corollaries in my classes or books, way too scary and mathematical, but good time now for a Corollary. Indeed, the general policy is to end each chapter in beauty, with a Theorem, but the above statement being a total triviality, as you can instantly see by contemplating the circle graph, calling it Theorem would have been inappropriate. Anyway, and regarding now the proof, based on Theorem 6.18, which is after all something quite instructive, here that is:

\medskip

(1) We know from chapter 2 that the adjacency matrix of the circle $X$ diagonalizes as follows, with $F_N=(w^{ij})$ with $w=e^{2\pi i/N}$ being as usual the Fourier matrix:
$$d=\frac{1}{N}F_N
\begin{pmatrix}
2\\
&w+w^{-1}\\
&&w^2+w^{-2}\\
&&&\ddots\\
&&&&w^{N-1}+w^{1-N}
\end{pmatrix}F_N^*$$

(2) Now by using Theorem 6.18, and some root of unity know-how, we get:
\begin{eqnarray*}
T_X
&=&\frac{1}{N}\prod_{k=1}^{N-1}(2-w^k-w^{-k})\\
&=&\frac{1}{N}\prod_{k=1}^{N-1}w^{-k}(2w^k-w^{2k}-1)\\
&=&\frac{1}{N}\prod_{k=1}^{N-1}-w^{-k}(1-w^k)^2\\
&=&\frac{1}{N}\prod_{k=1}^{N-1}(-w^{-k})\left(\prod_{k=1}^{N-1}(1-w^k)\right)^2\\
&=&\frac{1}{N}\times1\times N^2\\
&=&N
\end{eqnarray*}

(3) To be more precise, at the end we have used $x^N-1=(x-1)(x-w)\ldots(x-w^{N-1})$, first with $x=0$ in order to get the first needed formula, and then with $x=1$, after dividing both sides by $x-1$, in order to get the second needed formula. Thus, Corollary proved.
\end{proof}

\section*{6e. Exercises}

We had an exciting chapter here, but as a downside, many exercises left for you:

\begin{exercise}
Find an explicit realization of $F_2$, as some sort of concrete group.
\end{exercise}

\begin{exercise}
Apart from free groups $F_N$, how can you get trees as Cayley graphs?
\end{exercise}

\begin{exercise}
Check all the details for the Pr\"ufer encoding theorem.
\end{exercise}

\begin{exercise}
Regarding Pr\"ufer sequences, write a code for them, too.
\end{exercise}

\begin{exercise}
Fully learn the Cauchy-Binet formula, proof and examples.
\end{exercise}

\begin{exercise}
Clarify the other details in the proof of the Kirchoff formula.
\end{exercise}

\begin{exercise}
Further build on our results for the circulant graphs.
\end{exercise}

\begin{exercise}
Try to investigate the random walks on snowflake graphs.
\end{exercise}

As a bonus exercise, learn more about the free groups $F_N$, which, despite their seemingly abstract and rather boring definition, are quite fascinating objects.

\chapter{Genus, planarity}

\section*{7a. Planar graphs} 

Remember the discussion from the opening of chapter 4, regarding good and bad graphs. We have seen there that a nice and fruitful notion of ``goodness'', coming from simple and beautiful things like the circulant graphs, is the notion of transitivity. 

\bigskip

On the other hand, we have just seen that trees are beautiful and good as well, from an opposite aesthetic to that of the circulant graphs. So, question now for the two of us, subjective as they come: thinking deeply, what makes the beauty of a tree?

\bigskip

You might agree with this or not, but here is an answer to this question:

\index{crossings}
\index{torus}

\begin{answer}
Graphs fall into three classes:
\begin{enumerate}
\item Trees and other graphs which can be drawn without crossings are good. 

\item If we can still do this, but on a torus, the graph is bad. 

\item And the rest is evil.
\end{enumerate}
\end{answer}

Here the fact that trees are indeed planar is obvious, and as an illustration, here is some sort of ``random'' tree, which is clearly planar, no question about it:
$$\xymatrix@R=18pt@C=12pt{
\bullet&\bullet&\bullet&&&\bullet\\
\bullet&\bullet\ar@{-}[ul]\ar@{-}[u]\ar@{-}[ur]&\bullet&\bullet&\bullet&\bullet\ar@{-}[u]&\bullet\\
&\bullet\ar@{-}[ul]\ar@{-}[u]&\bullet&\bullet\ar@{-}[ul]\ar@{-}[u]&&\bullet\ar@{-}[ul]\ar@{-}[u]\ar@{-}[ur]\\
&&\bullet\ar@{-}[ul]\ar@{-}[u]\ar@{-}[ur]&&\bullet\ar@{-}[ur]\\
&&&\bullet\ar@{-}[ul]\ar@{-}[ur]
}$$

Of course, there are many other interesting examples of planar graphs. Consider for instance our beloved cube graph, that we met on many occasions, in the above:
$$
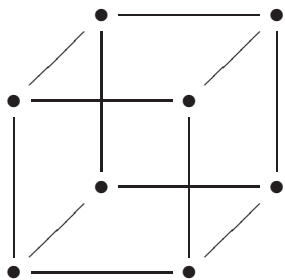
\xymatrix@R=20pt@C=20pt{
&\bullet\ar@{-}[rr]&&\bullet\\
\bullet\ar@{-}[rr]\ar@{-}[ur]&&\bullet\ar@{-}[ur]\\
&\bullet\ar@{-}[rr]\ar@{-}[uu]&&\bullet\ar@{-}[uu]\\
\bullet\ar@{-}[uu]\ar@{-}[ur]\ar@{-}[rr]&&\bullet\ar@{-}[uu]\ar@{-}[ur]
}$$

At the first glance, this graph does not look very planar. However, after thinking a bit, we can draw it as follows, making it clear that this graph is planar:
$$
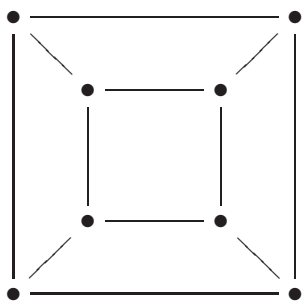
\xymatrix@R=15pt@C=15pt{
\bullet\ar@{-}[rrrr]&&&&\bullet\\
&\bullet\ar@{-}[rr]\ar@{-}[ul]&&\bullet\ar@{-}[ur]\\
\\
&\bullet\ar@{-}[uu]\ar@{-}[dl]\ar@{-}[rr]&&\bullet\ar@{-}[uu]\ar@{-}[dr]\\
\bullet\ar@{-}[rrrr]\ar@{-}[uuuu]&&&&\bullet\ar@{-}[uuuu]
}$$

Of course, not all graphs are planar. In order to find basic examples of non-planar graphs, we can look at simplices, and we are led to the following result:

\begin{proposition}
When looking at simplices, the segment $K_2$, the triangle $K_3$ and the tetrahedron $K_4$ are planar. However, the next simplex $K_5$, namely 
$$
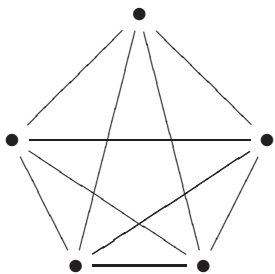
\xymatrix@R=14pt@C=11pt{
&&\bullet\ar@{-}[ddrr]\ar@{-}[ddll]\ar@{-}[ddddr]\ar@{-}[ddddl]\\
&&&&\\
\bullet\ar@{-}[ddr]\ar@{-}[rrrr]&&&&\bullet\ar@{-}[ddl]\ar@{-}[ddlll]\\
&&&&\\
&\bullet\ar@{-}[rr]\ar@{-}[uurrr]&&\bullet\ar@{-}[uulll]&&
}$$
is not planar. Nor are the higher simplices, $K_N$ with $N\geq6$, planar.
\end{proposition}

\begin{proof}
This is something quite elementary and intuitive, as follows:

\medskip

(1) The graphs $K_2,K_3,K_4$ are indeed planar, with this being clear for $K_2,K_3$, and with the planarity of $K_4$ being shown by the following picture for it:
$$
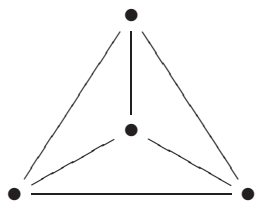
\xymatrix@R=12pt@C=12pt{
&&\bullet\ar@{-}[dd]\ar@{-}[dddrr]\ar@{-}[dddll]\\
\\
&&\bullet\ar@{-}[drr]\ar@{-}[dll]\\
\bullet\ar@{-}[rrrr]&&&&\bullet}$$

(2) Regarding now the non-planarity of $K_5$, let us try to manufacture an intuitive proof for this. In order to draw $K_5$ in a planar way, we first have to draw its subgraph $K_4$ in a planar way, and it is pretty much clear that this can only be done as a variation of the above picture, from (1), with curved edges this time, as follows:
$$
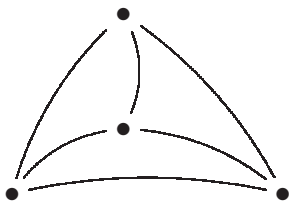
\xymatrix@R=12pt@C=11pt{
&&\bullet\ar@/^/@{-}[dd]\ar@/^/@{-}[dddrrr]\ar@/_/@{-}[dddll]\\
\\
&&\bullet\ar@/^/@{-}[drrr]\ar@/_/@{-}[dll]\\
\bullet\ar@/^/@{-}[rrrrr]&&&&&\bullet}$$

But with this in hand, it is clear that there is no room in the plane for our 5th vertex, as to avoid crossings. Indeed, we have 4 possible regions in the plane for this 5th vertex, and each of them is forbidden by the edge towards a certain vertex, as follows:
$$
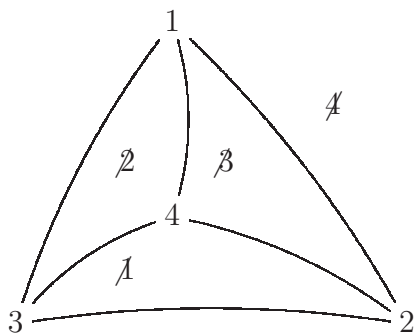
\xymatrix@R=4pt@C=4pt{
&&&&1\ar@/^/@{-}[dddd]\ar@/^/@{-}[ddddddrrrrrr]\ar@/_/@{-}[ddddddllll]\\
\\
&&&&&&&&\not\!4\\
&&&\not2&&\not3\\
&&&&4\ar@/^/@{-}[ddrrrrrr]\ar@/_/@{-}[ddllll]\\
&&&\not1\\
3\ar@/^/@{-}[rrrrrrrrrr]&&&&&&&&&&2}$$

(3) Finally, the fact that the graphs $K_N$ with $N\geq6$ are not planar either follows from the fact that their subgraphs $K_5$ are not planar, that we know from (2).
\end{proof}

In order to find some further examples of non-planar graphs, we can look as well at the bipartite simplices, and we are led to the following result:

\begin{proposition}
When looking at bipartite simplices, the square $K_{2,2}$ is planar, and so are all the graphs $K_{2,N}$. However, the next such graph, namely $K_{3,3}$,
$$
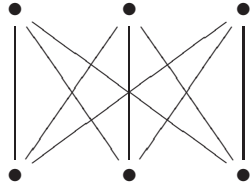
\xymatrix@R=50pt@C=30pt{
\bullet\ar@{-}[d]\ar@{-}[dr]\ar@{-}[drr]&\bullet\ar@{-}[dl]\ar@{-}[d]\ar@{-}[dr]&\bullet\ar@{-}[dll]\ar@{-}[dl]\ar@{-}[d]\\
\bullet&\bullet&\bullet}$$
called ``utility graph'' is not planar. Nor are planar the graphs $K_{M,N}$, for any $M,N\geq3$.
\end{proposition}

\begin{proof}
Again, this is something elementary and intuitive, as follows:

\medskip

(1) In what regards the first bipartite simplex, which is $K_{2,2}$, this is indeed the square, which is of course a planar graph, as shown by the following equality:
$$
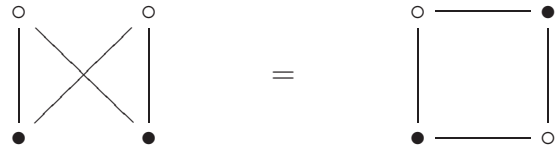
\xymatrix@R=12pt@C=36pt{
\circ\ar@{-}[dd]\ar@{-}[ddr]&\circ\ar@{-}[dd]\ar@{-}[ddl]&&\circ\ar@{-}[dd]\ar@{-}[r]&\bullet\ar@{-}[dd]\\
&&=\\
\bullet&\bullet&&\bullet\ar@{-}[r]&\circ}$$

(2) Regarding now the bipartite simplex $K_{2,N}$ with $N\geq2$ arbitrary, this graph looks at follows, with $N$ vertices in the lower row:
$$
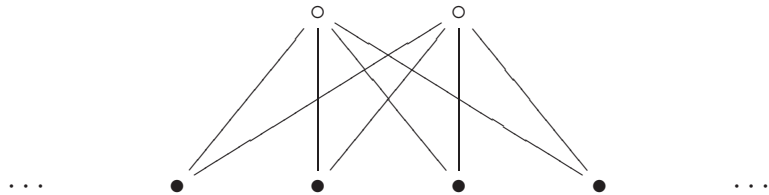
\xymatrix@R=53pt@C=40pt{
&&\circ\ar@{-}[d]\ar@{-}[dr]\ar@{-}[drr]\ar@{-}[dl]&\circ\ar@{-}[dll]\ar@{-}[dl]\ar@{-}[d]\ar@{-}[dr]\\
\cdots&\bullet&\bullet&\bullet&\bullet&\cdots}$$

But this graph is planar too, because we can draw it in the following way:
$$
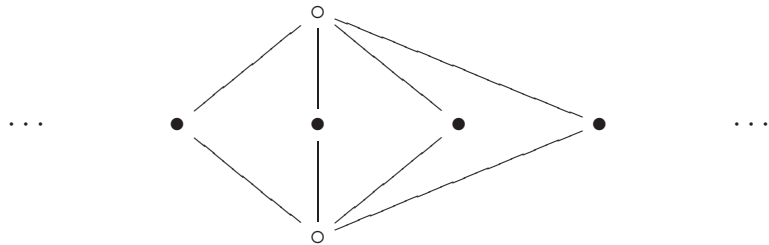
\xymatrix@R=30pt@C=40pt{
&&\circ\ar@{-}[d]\ar@{-}[dr]\ar@{-}[drr]\ar@{-}[dl]\\
\cdots&\bullet&\bullet&\bullet&\bullet&\cdots\\
&&\circ\ar@{-}[u]\ar@{-}[ur]\ar@{-}[urr]\ar@{-}[ul]}$$

(3) Regarding now $K_{3,3}$, as before with the simplex $K_5$, the result here is quite clear by thinking a bit, and drawing pictures. To be more precise, reasoning by contradiction, we first have to draw its subgraph $K_{2,3}$ in a planar way, and this is done as follows:
$$
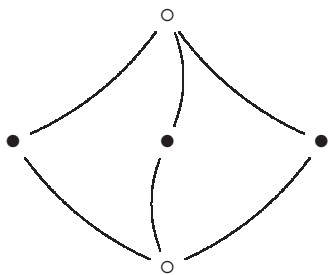
\xymatrix@R=35pt@C=45pt{
&\circ\ar@/^/@{-}[d]\ar@/_/@{-}[dr]\ar@/^/@{-}[dl]\\
\bullet&\bullet&\bullet\\
&\circ\ar@/^/@{-}[u]\ar@/_/@{-}[ur]\ar@/^/@{-}[ul]}$$

(4) But now, as before with $K_5$, it is clear that there is no room in the plane for our 6th vertex, as to avoid crossings. Indeed, we have 3 regions in the plane for this 6th vertex, and each of them is forbidden by the edge towards a certain vertex, as follows:
$$
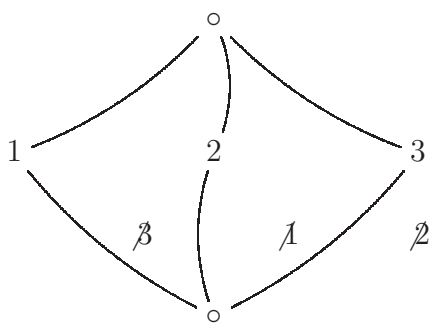
\xymatrix@R=15pt@C=13pt{
&&&\circ\ar@/^/@{-}[dd]\ar@/_/@{-}[ddrrr]\ar@/^/@{-}[ddlll]\\
\\
1&&&2&&&3\\
&&\not3&&\not1&&\not2\\
&&&\circ\ar@/^/@{-}[uu]\ar@/_/@{-}[uurrr]\ar@/^/@{-}[uulll]}$$

Thus, theorem proved for the utility graph $K_{3,3}$, via the same method as for $K_5$.

\medskip

(5) Still talking $K_{3,3}$, let us mention that this is called indeed ``utility graph'', as said above, due to a certain story with it. The story involves 3 companies, selling gas, water and electricity to 3 customers, and looking for a way to arrange their underground tubes and wires as not to cross. Thus, they are looking to implement their ``utillity graph'', which is $K_{3,3}$, in a planar way, and unfortunately, this is not possible.

\medskip

(6) And as further comments on $K_{3,3}$, quite remarkably, in recent years the stocks of the above-mentioned 3 companies have skyrocketed, apparently due to very good business done by their Saturn ring branches, which were able to considerably cut from their costs. But are we here for talking about economy, or about mathematics.

\medskip

(7) Finally, the bipartite simplex $K_{M,N}$ with $M,N\geq3$ is not planar either, because it contains $K_{3,3}$. Thus, we are led to the conclusions in the statement.
\end{proof}

So long for planar graphs, which are the ``good'' ones, in the sense of Answer 7.1. But, what about the graphs which are bad, or evil, in the sense of Answer 7.1? It is pretty much clear that their mathematics is quite complicated, and more on this later, but in the meantime, in regards with the graphs $K_N$ and $K_{M,N}$, let us formulate:

\begin{theorem}
The simplices $K_N$ and bipartite simplices $K_{M,N}$ with $M\leq N$ are as follows, with respect to our good, bad and evil classification scheme:
\begin{enumerate}
\item The good, planar graphs are $K_2,K_3,K_4$ and $K_{2,N}$.

\item The bad graphs, that can be drawn however on a torus, are $K_5,K_6,K_7$ and $K_{3,3},K_{3,4},K_{3,5},K_{3,6}$, and $K_{4,4}$.

\item As for the remaining graphs, these are all evil, not drawable on a torus.
\end{enumerate}
\end{theorem}

\begin{proof}
This does not look easy at all, so basically theorem coming without proof, and more on this later. However, in the meantime, some remarks about all this:

\medskip

(1) What we know so far, from Proposition 6.2 and Proposition 6.3, is that (1) holds. Not bad, and it remains to decide, for the remaining graphs, if these are bad or evil.

\medskip

(2) Regarding the fact that the graphs in (2) are indeed bad, that is, drawable on a torus, this is normally not difficult, just requiring some silence, thinking, and patience. Indeed, to start with, there is no really need to bother with 3D pictures for the torus, because this torus appears from a rectangle, by gluing the opposite sides, as follows:
$$
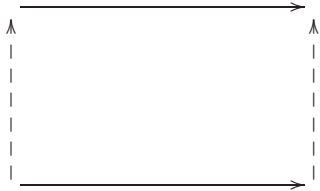
\xymatrix@R=60pt@C=50pt{
\ar[rr]&&\\
\ar@{-->}[u]\ar[rr]&&\ar@{-->}[u]}$$

(3) With this drawing method in hand, let us first look at $K_5$. After some thinking, this graph is indeed toral, because it can be drawn on a torus as follows, with the convention that the dotted edges connect, according to our gluing conventions above:
$$
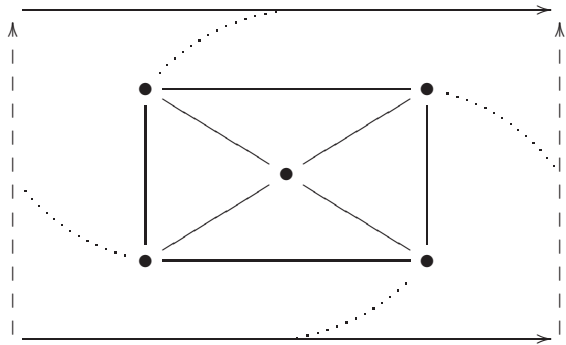
\xymatrix@R=20pt@C=40pt{
\ar[rrrr]&&&&\\
&\bullet\ar@{-}[rr]\ar@{-}[dd]\ar@/^/@{.}[ur]&&\bullet\ar@{-}[dd]\ar@/^/@{.}[dr]\\
&&\bullet\ar@{-}[ur]\ar@{-}[ul]\ar@{-}[dr]\ar@{-}[dl]&&\\
&\bullet\ar@{-}[rr]\ar@/^/@{.}[ul]&&\bullet\\
\ar@{-->}[uuuu]\ar[rrrr]&&\ar@/_/@{.}[ur]&&\ar@{-->}[uuuu]}$$

(4) A similar method works for the bipartite simplex $K_{3,3}$, with this graph being drawable on a torus as follows, with our various pictorial conventions above:
$$
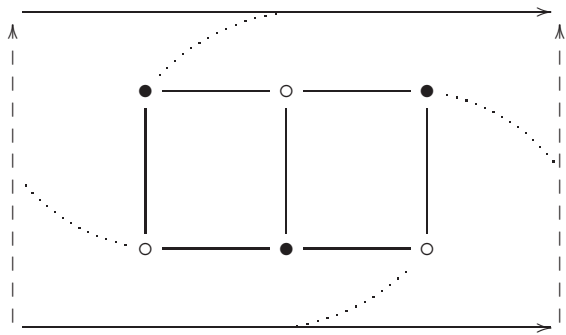
\xymatrix@R=20pt@C=40pt{
\ar[rrrr]&&&&\\
&\bullet\ar@{-}[r]\ar@{-}[dd]\ar@/^/@{.}[ur]&\circ\ar@{-}[r]\ar@{-}[dd]&\bullet\ar@{-}[dd]\ar@/^/@{.}[dr]\\
&&&&\\
&\circ\ar@{-}[r]\ar@/^/@{.}[ul]&\bullet\ar@{-}[r]&\circ\\
\ar@{-->}[uuuu]\ar[rrrr]&&\ar@/_/@{.}[ur]&&\ar@{-->}[uuuu]}$$

(5) Thus, theorem fully proved for the basic non-planar graphs, namely $K_5,K_{3,3}$. In what regards $K_6,K_7$, and then $K_{3,4},K_{3,5},K_{3,6}$, and then $K_{4,4}$, we will leave drawing them on a torus as an instructive exercise, that is, 6 graphs, say 1/2 hour each. Of course you might say that it is enough to do it for $K_7,K_{3,6},K_{4,4}$, but don't be greedy, if you try directly these 3 graphs, these might well cost you 2 hours each, so no gain.

\medskip

(6) As for the graphs left, here we just have to prove that $K_8,K_{3,7},K_{4,5}$ are evil. This is something a bit more difficult, requiring extreme silence, thinking, and patience, and as before, we will leave this as an instructive exercise, of rather difficult type.

\medskip

(7) We will be of course back to this, later in this chapter, but in meantime, you might wonder if there is something more conceptual, say some kind of ``formula'', behind all this. And in answer, yes there is a formula. The idea is that associated to any graph is a certain positive integer  $g\in\mathbb N$, called genus, which is as follows:
$$g=\begin{cases}
0&{\rm (planar)}\\
1&{\rm (toral)}\\
\geq2&{\rm (other)}
\end{cases}$$

(8) In general, the genus is quite hard to compute, and more on this later. However, coming a bit in advance, the genus of $K_N$ is indeed computable, given by the following formula, with $\lceil x\rceil$ being the ceiling function, that is, the smallest integer $\geq x$:
$$g(K_N)=\left\lceil\frac{(N-3)(N-4)}{12}\right\rceil$$

And for $N=4,5,6,7,8$, this gives the following values, which are what we need:
$$0,1,1,1,2$$

(9) As for the bipartite simplices, $K_{M,N}$ with $M\leq N$, here again the genus is explicitly computable, with appropriate tools, the formula being as follows:
$$g(K_{M,N})=\left\lceil\frac{(M-2)(N-2)}{4}\right\rceil$$

And at $M=3$ and $N=2,3,4,5,6,7$, and $M=4$ and $N=4,5$ we get, as needed:
$$0,1,1,1,1,2\qquad,\qquad 1,2$$

But more on such things, which are quite technical, later in this chapter.
\end{proof}

\section*{7b. Main theorems} 

Let us get back now to the planar graphs, and try to develop some systematic theory for them. As a first main result about these graphs, we have:

\begin{theorem}
The fact that a graph $X$ is non-planar can be checked as follows:
\begin{enumerate}
\item Kuratowski criterion: $X$ contains a subdivision of $K_5$ or $K_{3,3}$.

\item Wagner criterion: $X$ has a minor of type $K_5$ or $K_{3,3}$.
\end{enumerate}
\end{theorem}

\begin{proof}
This is obviously something quite powerful, when thinking at the potential applications, and non-trivial to prove as well, the idea being as follows:

\medskip

(1) Regarding the Kuratowski criterion, the convention is that ``subdivision'' means graph obtained by inserting vertices into edges, e.g. replacing $\bullet-\bullet$ with $\bullet-\bullet-\bullet$.

\medskip

(2) Regarding the Wagner criterion, the convention there is that ``minor'' means graph obtained by contracting certain edges into vertices.

\medskip

(3) Regarding now the proofs, the Kuratowski and Wagner criteria are more or less equivalent, and their proof is via standard, although long, recurrence methods.

\medskip

(4) In short, non-trivial, but rather routine results that we have here, and we will leave finding and studying their complete proofs as an instructive exercise.

\medskip

(5) Finally, let us mention that, often in practice, Wagner works a bit better than Kuratowski. More on this in a moment, when discussing examples.
\end{proof}

Regarding now the applications of the Kuratowski and Wagner criteria, things are quite tricky here, because most of the graphs that we met so far in this book are trees and other planar graphs, for which these criteria are not needed. We have as well the graphs $K_N$ and $K_{M,N}$, to which these criteria apply trivially. Thus, for illustrations, we have to go to more complicated graphs, and as a standard example here, we have:

\begin{proposition}
The Petersen graph $P$, namely
$$
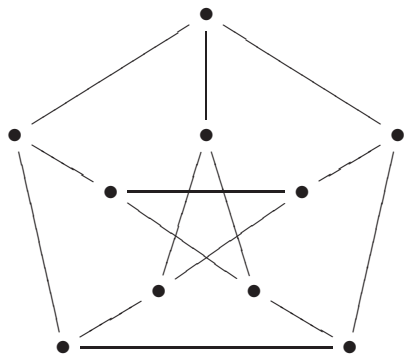
\xymatrix@R=1pt@C=5pt{
&&&&\bullet\ar@{-}[dddddrrrr]\ar@{-}[dddddllll]\\
\\
\\
\\
\\
\bullet\ar@{-}[ddddddddr]&&&&\bullet\ar@{-}[uuuuu]\ar@{-}[ddddddl]\ar@{-}[ddddddr]&&&&\bullet\ar@{-}[ddddddddl]\\
\\
&&
\bullet\ar@{-}[uull]\ar@{-}[ddddrrr]\ar@{-}[rrrr]&&&&\bullet\ar@{-}[uurr]\ar@{-}[ddddlll]\\
\\
\\
\\
&&&\bullet&&\bullet\\
\\
&\bullet\ar@{-}[rrrrrr]\ar@{-}[uurr]&&&&&&\bullet\ar@{-}[uull]}
$$
is not planar, the reasons for this being as follows:
\begin{enumerate}
\item Kuratowski: $P$ contains no subdivision of $K_5$, but contains a subdivision of $K_{3,3}$.

\item Wagner: $P$ has both $K_5$ and $K_{3,3}$ as minors.
\end{enumerate}
\end{proposition}

\begin{proof}
We have four things to be proved, all instructive, the idea being as follows:

\medskip

(1) To start with, $P$ contains no subdivision of $K_5$, because $P$ has valence 3, while $K_5$ has valence $4>3$. Thus, game over with the Kuratowski criterion using $K_5$.

\medskip

(2) On the other hand, regarding $K_{3,3}$, this has valence 3, exactly as $P$, so there is a chance for the Kuratowski criterion using $K_{3,3}$ to apply to $P$. And this is indeed the case, showing that $P$ is not planar, with the subdivision of $K_{3,3}$ being obtained as follows:
$$
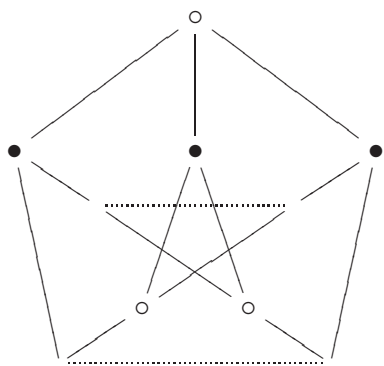
\xymatrix@R=2pt@C=7pt{
&&&&\circ\ar@{-}[dddddrrrr]\ar@{-}[dddddllll]\\
\\
\\
\\
\\
\bullet\ar@{-}[ddddddddr]&&&&\bullet\ar@{-}[uuuuu]\ar@{-}[ddddddl]\ar@{-}[ddddddr]&&&&\bullet\ar@{-}[ddddddddl]\\
\\
&&
\ar@{-}[uull]\ar@{-}[ddddrrr]\ar@{.}[rrrr]&&&&\ar@{-}[uurr]\ar@{-}[ddddlll]\\
\\
\\
\\
&&&\circ&&\circ\\
\\
&\ar@{.}[rrrrrr]\ar@{-}[uurr]&&&&&&\ar@{-}[uull]}
$$

To be more precise, ignoring the dotted edges, what we have here is indeed a subdivision of $K_{3,3}$, obtained from $K_{3,3}$ by inserting 4 vertices into 4 certain edges.

\medskip

(3) Regarding now Wagner, in contrast with Kuratowski, and better than it, this applies to $P$ by using $K_5$, with the $K_5$ minor of $P$ being obtained as follows:
$$
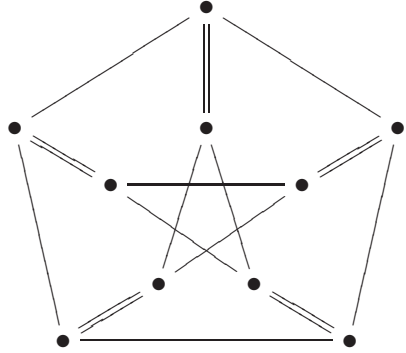
\xymatrix@R=1pt@C=5pt{
&&&&\bullet\ar@{-}[dddddrrrr]\ar@{-}[dddddllll]\\
\\
\\
\\
\\
\bullet\ar@{-}[ddddddddr]&&&&\bullet\ar@{=}[uuuuu]\ar@{-}[ddddddl]\ar@{-}[ddddddr]&&&&\bullet\ar@{-}[ddddddddl]\\
\\
&&
\bullet\ar@{=}[uull]\ar@{-}[ddddrrr]\ar@{-}[rrrr]&&&&\bullet\ar@{=}[uurr]\ar@{-}[ddddlll]\\
\\
\\
\\
&&&\bullet&&\bullet\\
\\
&\bullet\ar@{-}[rrrrrr]\ar@{=}[uurr]&&&&&&\bullet\ar@{=}[uull]}
$$

To be more precise, the convention here is that we identify the vertices joined by = edges, and this procedure obviously producing the graph $K_5$, we have $K_5$ as minor. 

\medskip

(4) Finally, still regarding Wagner, and adding to the power of this criterion, this applies to $P$ as well by using $K_{3,3}$, with the $K_{3,3}$ minor of $P$ being obtained as follows, again with the above conventions, namely identifying the vertices joined by = edges:
$$
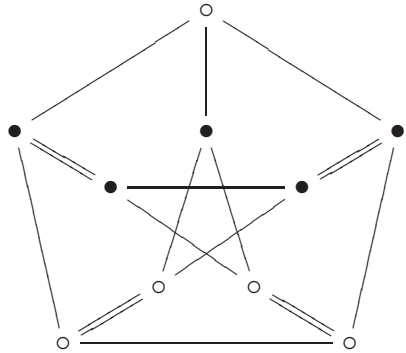
\xymatrix@R=1pt@C=5pt{
&&&&\circ\ar@{-}[dddddrrrr]\ar@{-}[dddddllll]\\
\\
\\
\\
\\
\bullet\ar@{-}[ddddddddr]&&&&\bullet\ar@{-}[uuuuu]\ar@{-}[ddddddl]\ar@{-}[ddddddr]&&&&\bullet\ar@{-}[ddddddddl]\\
\\
&&
\bullet\ar@{=}[uull]\ar@{-}[ddddrrr]\ar@{-}[rrrr]&&&&\bullet\ar@{=}[uurr]\ar@{-}[ddddlll]\\
\\
\\
\\
&&&\circ&&\circ\\
\\
&\circ\ar@{-}[rrrrrr]\ar@{=}[uurr]&&&&&&\circ\ar@{=}[uull]}
$$

Thus, we are led to the various conclusions in the statement.
\end{proof}

The above result is quite interesting, the Petersen graph being a key object in combinatorics, usually providing counterexamples to all sorts of graph theory statements that can be made. We will be back to this graph later, with a proof that it is toral.

\bigskip

As a second main result now about the planar graphs, we have:

\begin{theorem}
For a connected planar graph we have the Euler formula
$$v-e+f=2$$
with $v,e,f$ being the number of vertices, edges and faces. 
\end{theorem}

\begin{proof}
This is something very standard, the idea being as follows:

\medskip

(1) Regarding the precise statement, given a connected planar graph, drawn in a planar way, without crossings, we can certainly talk about the numbers $v$ and $e$, as for any graph, and also about $f$, as being the number of faces that our graph has, in our picture, with these including by definition the outer face too, the one going to $\infty$. With these conventions, the claim is that the Euler formula $v-e+f=2$ holds indeed.

\medskip

(2) As a first illustration for how this formula works, consider a triangle:
$$
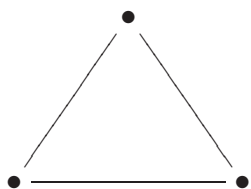
\xymatrix@R=50pt@C=30pt{
&\bullet\ar@{-}[dl]\ar@{-}[dr]\\
\bullet\ar@{-}[rr]&&\bullet}$$

Here we have $v=e=3$, and $f=2$, with this accounting for the interior and exterior, and we conclude that the Euler formula holds indeed in this case, as follows:
$$3-3+2=2$$

(3) More generally now, let us look at an arbitrary $N$-gon graph:
$$
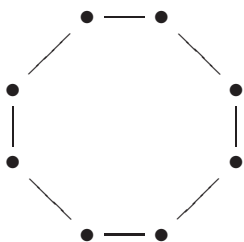
\xymatrix@R=15pt@C=15pt{
&\bullet\ar@{-}[r]\ar@{-}[dl]&\bullet\ar@{-}[dr]\\
\bullet\ar@{-}[d]&&&\bullet\ar@{-}[d]\\
\bullet\ar@{-}[dr]&&&\bullet\ar@{-}[dl]\\
&\bullet\ar@{-}[r]&\bullet}$$

Then, for this graph, the Euler formula holds indeed, as follows:
$$N-N+2=2$$

(4) With these examples discussed, let us look now for a proof. The idea will be to proceed by recurrence on the number of faces $f$. And here, as a first observation, the result holds at $f=1$, where our graph must be planar and without cycles, and so must be a tree. Indeed, with $N$ being the number of vertices, the Euler formula holds, as:
$$N-(N-1)+1=2$$

(5) At $f=2$ now, our graph must be an $N$-gon as above, but with some trees allowed to grow from the vertices, with an illustrating example here being as follows:
$$
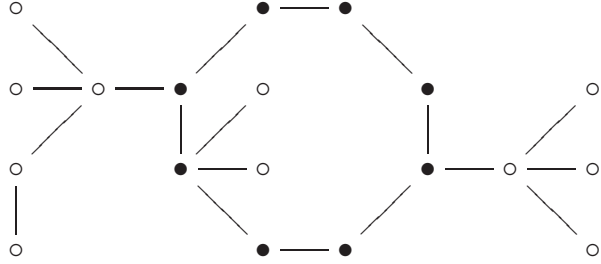
\xymatrix@R=18pt@C=18pt{
\circ\ar@{-}[dr]&&&\bullet\ar@{-}[r]\ar@{-}[dl]&\bullet\ar@{-}[dr]\\
\circ\ar@{-}[r]&\circ\ar@{-}[r]&\bullet\ar@{-}[d]&\circ&&\bullet\ar@{-}[d]&&\circ\\
\circ\ar@{-}[ur]\ar@{-}[d]&&\bullet\ar@{-}[dr]\ar@{-}[ur]\ar@{-}[r]&\circ&&\bullet\ar@{-}[dl]\ar@{-}[r]&\circ\ar@{-}[r]\ar@{-}[ur]\ar@{-}[dr]&\circ\\
\circ&&&\bullet\ar@{-}[r]&\bullet&&&\circ}$$

But here we can argue, again based on the fact that for a rooted tree, the non-root vertices are in obvious bijection with the edges, that removing all these trees won't change the problem. So, we are left with the problem for the $N$-gon, already solved in (3).

\medskip

(6) And so on, the idea being that we can first remove all the trees, by using the argument in (5), and then we are left with some sort of agglomeration of $N$-gons, for which we can check the Euler formula directly, a bit as in (3), or by recurrence.

\medskip

(7) To be more precise, let us try to do the recurrence on the number of faces $f$. For this purpose, consider one of the faces of our graph, which looks as follows, with $v_i$ denoting the number of vertices on each side, with the endpoints excluded:
$$
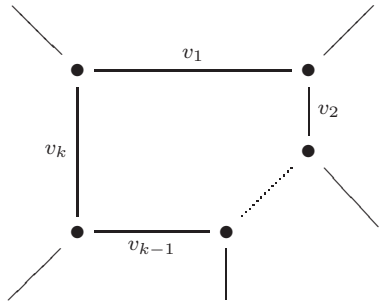
\xymatrix@R=18pt@C=18pt{
&&&&&\\
&\bullet\ar@{-}[dd]_{v_k}\ar@{-}[rrr]^{v_1}\ar@{-}[ul]&&&\bullet\ar@{-}[d]^{v_2}\ar@{-}[ur]\\
&&&&\bullet\ar@{.}[dl]\ar@{-}[dr]\\
&\bullet\ar@{-}[rr]_{v_{k-1}}\ar@{-}[dl]&&\bullet\ar@{-}[d]&&&\\
&&&&}$$

(8) Now let us collapse this chosen face to a single point, in the obvious way. In this process, the total number of vertices $v$, edges $e$, and faces $f$, evolves as follows:
$$v\to v-k+1-\sum v_i$$
$$e\to e-\sum(v_i+1)$$
$$f\to f-1$$

Thus, in this process, the Euler quantity $v-e+f$ evolves as follows:
\begin{eqnarray*}
v-e+f
&\to&v-k+1-\sum v_i-e+\sum(v_i+1)+f-1\\
&=&v-k+1-\sum v_i-e+\sum v_i+k+f-1\\
&=&v-e+f
\end{eqnarray*}

So, done with the recurrence, and the Euler formula is proved.
\end{proof}

As a famous application, or rather version, of the Euler formula, let us record:

\begin{proposition}
For a convex polyhedron we have the Euler formula
$$v-e+f=2$$
with $v,e,f$ being the number of vertices, edges and faces. 
\end{proposition}

\begin{proof}
This is more or less the same thing as Theorem 7.7, save for getting rid of the internal trees of the planar graph there, the idea being as follows:

\medskip

(1) In one sense, consider a convex polyhedron $P$. We can then enlarge one face, as much as needed, and then smash our polyhedron with a big hammer, as to get a planar graph $X$. As an illustration, here is how this method works, for a cube:
$$
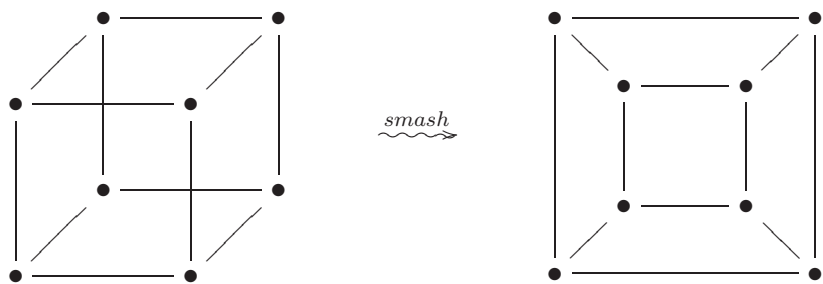
\xymatrix@R=20pt@C=20pt{
&\bullet\ar@{-}[rr]&&\bullet\\
\bullet\ar@{-}[rr]\ar@{-}[ur]&&\bullet\ar@{-}[ur]\\
&\bullet\ar@{-}[rr]\ar@{-}[uu]&&\bullet\ar@{-}[uu]\\
\bullet\ar@{-}[uu]\ar@{-}[ur]\ar@{-}[rr]&&\bullet\ar@{-}[uu]\ar@{-}[ur]
}
\qquad
\xymatrix@R=15pt@C=30pt{\\ \\ \ar@{~>}[r]^{smash}&\\ \\ }
\qquad
\xymatrix@R=13pt@C=13pt{
\bullet\ar@{-}[rrrr]&&&&\bullet\\
&\bullet\ar@{-}[rr]\ar@{-}[ul]&&\bullet\ar@{-}[ur]\\
\\
&\bullet\ar@{-}[uu]\ar@{-}[dl]\ar@{-}[rr]&&\bullet\ar@{-}[uu]\ar@{-}[dr]\\
\bullet\ar@{-}[rrrr]\ar@{-}[uuuu]&&&&\bullet\ar@{-}[uuuu]
}$$

But, in this process, each of the numbers $v,e,f$ stays the same, so we get the Euler formula for $P$, as a consequence of the Euler formula for $X$, from Theorem 7.7.

\medskip

(2) Conversely, consider a connected planar graph $X$. Then, save for getting rid of the internal trees, as explained in the proof of Theorem 7.7, we can assume that we are dealing with an agglomeration of $N$-gons, again as explained in the proof of Theorem 7.7. But now, we can inflate our graph as to obtain a convex polyhedron $P$, as follows:
$$
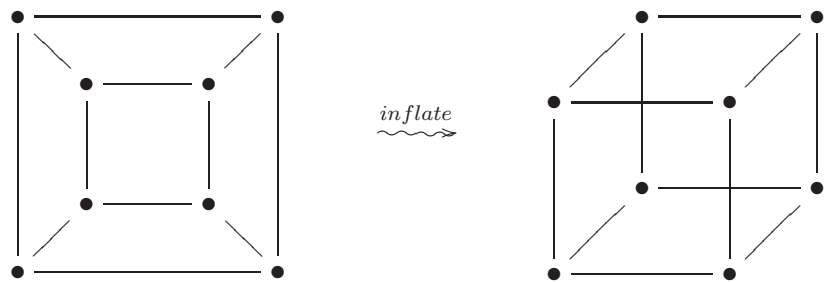
\xymatrix@R=13pt@C=13pt{
\bullet\ar@{-}[rrrr]&&&&\bullet\\
&\bullet\ar@{-}[rr]\ar@{-}[ul]&&\bullet\ar@{-}[ur]\\
\\
&\bullet\ar@{-}[uu]\ar@{-}[dl]\ar@{-}[rr]&&\bullet\ar@{-}[uu]\ar@{-}[dr]\\
\bullet\ar@{-}[rrrr]\ar@{-}[uuuu]&&&&\bullet\ar@{-}[uuuu]}
\qquad
\xymatrix@R=15pt@C=30pt{\\ \\ \ar@{~>}[r]^{inflate}&\\ \\ }
\qquad
\xymatrix@R=20pt@C=20pt{
&\bullet\ar@{-}[rr]&&\bullet\\
\bullet\ar@{-}[rr]\ar@{-}[ur]&&\bullet\ar@{-}[ur]\\
&\bullet\ar@{-}[rr]\ar@{-}[uu]&&\bullet\ar@{-}[uu]\\
\bullet\ar@{-}[uu]\ar@{-}[ur]\ar@{-}[rr]&&\bullet\ar@{-}[uu]\ar@{-}[ur]
}$$

Again, in this process, each of the numbers $v,e,f$ will stay the same, and so we get the Euler formula for $X$, as a consequence of the Euler formula for $P$.
\end{proof}

Summarizing, Euler formula understood, but as a matter of making sure that we didn't mess up anything with our mathematics, let us do some direct checks as well:

\begin{proposition}
The Euler formula $v-e+f=2$ holds indeed for the five possible regular polyhedra, as follows:
\begin{enumerate}
\item Tetrahedron: $4-6+4=2$.

\item Cube: $8-12+6=2$.

\item Octahedron: $6-12+8=2$.

\item Dodecahedron: $20-30+12=2$.

\item Isocahedron: $12-30+20=2$.
\end{enumerate}
\end{proposition} 

\begin{proof}
The figures in the statement are certainly the good ones for the tetrahedron and the cube. Regarding now the octahedron, again the figures are the good ones, by thinking in 3D, but as an interesting exercise for us, which is illustrating for the above, let us attempt to find a nice way of drawing the corresponding graph:

\medskip

(1) To start with, the ``smashing'' method from the proof of Proposition 7.8 provides us with a graph which is certainly planar, but which, even worse than before for the cube, sort of misses the whole point with the 3D octahedron, its symmetries, and so on:
$$
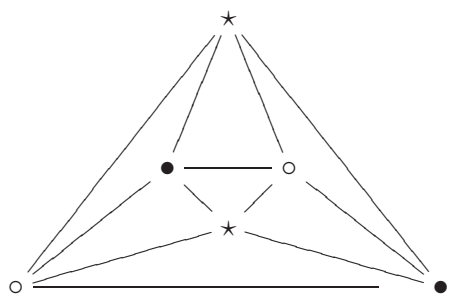
\xymatrix@R=10pt@C=10pt{
&&&&\star\ar@{-}[dddl]\ar@{-}[dddr]\ar@{-}[dddddllll]\ar@{-}[dddddrrrr]\\
\\
\\
&&&\bullet\ar@{-}[rr]&&\circ\\
&&&&\star\ar@{-}[ur]\ar@{-}[ul]\\
\circ\ar@{-}[rrrrrrr]\ar@{-}[urrrr]\ar@{-}[uurrr]&&&&&&&&\bullet\ar@{-}[ullll]\ar@{-}[uulll]}$$

(2) Much nicer, instead, is the following picture, which still basically misses the 3D beauty of the octahedron, but at least reveals some of its symmetries:
$$
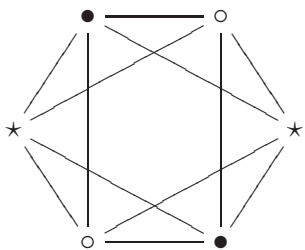
\xymatrix@R=30pt@C=15pt{
&\bullet\ar@{-}[dl]\ar@{-}[rr]\ar@{-}[dd]\ar@{-}[drrr]&&\circ\ar@{-}[dr]\ar@{-}[dd]\ar@{-}[dlll]\\
\star\ar@{-}[dr]&&&&\star\ar@{-}[dl]\\
&\circ\ar@{-}[rr]\ar@{-}[urrr]&&\bullet\ar@{-}[ulll]
}$$

In short, you get the point, quite subjective all this, and as a conclusion, drawing graphs in an appropriate way remains an art. As for the dodecahedron and isocahedron, exercise here for you, and if failing, take some drawing classes. Math is not everything.
\end{proof}

The Euler formula $v-e+f=2$, in both its above formulations, the graph one from Theorem 7.7, and the polyhedron one from Proposition 7.8, is something very interesting, at the origin of modern pure mathematics, and having countless other versions and generalizations. We will be back to it on several occasions, in what follows.

\bigskip

As a third main result now about the planar graphs, we have:

\begin{theorem}
Any planar graph has the following properties:
\begin{enumerate}
\item It is vertex $4$-colorable.

\item It is a $4$-partite graph.
\end{enumerate}
\end{theorem}

\begin{proof}
Heavy theorem that we have here, with our comments being as follows:

\medskip

(1) This is something that we talked about in chapter 1, by calling there ``maps'' the planar graphs, and with the comment that the proof is something extremely complicated. So, you might perhaps expect, with our graph learning doing quite well, that we might have some new things to say, about this. Unfortunately, no, this is definitely very difficult, beyond our reach, but do not hesitate to look it up, and learn more about it.

\medskip

(2) This being said, there are a few elementary things that we can do. First is the observation that 3 colors will not do, and that not all planar graphs are tripartite, and exercise here for you, to find some counterexamples. Also, it is possible to say of few other things, of theoretical nature, providing some evidence for the theorem, as stated.  

\medskip

(3) To be more precise, the fact that 3 colors will not do is clear for a flattened tetrahedron, and with the 4-coloring here $ABCD$ being, of course, as follows:
$$
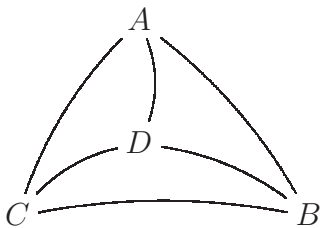
\xymatrix@R=12pt@C=11pt{
&&A\ar@/^/@{-}[dd]\ar@/^/@{-}[dddrrr]\ar@/_/@{-}[dddll]\\
\\
&&D\ar@/^/@{-}[drrr]\ar@/_/@{-}[dll]\\
C\ar@/^/@{-}[rrrrr]&&&&&B}$$

There are of course many things that can be said, along these lines.
\end{proof}

As a conclusion to all this, planar graphs and related topics, as you can see from the above, the theory of planar graphs can vary a lot, with:

\bigskip

(1) Theorem 7.5 being something quite tricky.

\bigskip

(2) Theorem 7.7 being something rather deep and trivial.

\bigskip

(3) Theorem 7.10 being something of extreme difficulty. 

\bigskip

Quite fascinating all this, hope you agree with me. In short, beware of planar graphs, you never know when coming across something trivial, or extremely complicated.

\section*{7c. Some topology} 

Switching topics now, let us get into the following question, which puts under the spotlight the graphs that we neglected so far, in this chapter:

\begin{question}
What are the graphs which are not planar, but can be however drawn on a torus? Also, what about graphs which can be drawn on higher surfaces, having $g\geq2$ holes, instead of the $g=0$ holes of the sphere, and the $g=1$ hole of the torus? 
\end{question}

All this looks quite interesting, but before going head-first into such questions, let us pause our study of graphs, and look instead at surfaces, and other manifolds. The same questions make sense here, does our manifold have holes, and how many holes, of which type, and so on. Let us start with something that we know well, namely:

\begin{definition}
A topological space $X$ is called path connected when any two points $x,y\in X$ can be connected by a path. That is, given any two points $x,y\in X$, we can find a continuous function $f:[0,1]\to X$ such that $f(0)=x$ and $f(1)=y$.
\end{definition}

The problem is now, given a connected space $X$, how to count its ``holes''. And this is quite subtle problem, because as examples of such spaces we have:

\medskip

(1) The sphere, the donut, the double-holed donut, the triple-holed donut, and so on. These spaces are quite simple, and intuition suggests to declare that the number of holes of the $N$-holed donut is, and you guessed right, $N$.

\medskip

(2) However, we have as well as example the empty sphere, I mean just the crust of the sphere, and while this obviously falls into the class of ``one-holed spaces'', this is not the same thing as a donut, its hole being of different nature.

\medskip

(3) As another example, consider again the sphere, but this time with two tunnels drilled into it, in the shape of a cross. Whether that missing cross should account for 1 hole, or for 2 holes, or for something in between, I will leave it up to you. 

\medskip

Summarizing, things are quite tricky, suggesting that the ``number of holes'' of a topological space $X$ is not an actual number, but rather something more complicated. Now with this in mind, let us formulate the following definition:

\index{homotopy group}
\index{hole}
\index{loop}
\index{topological space}
\index{loop on space}
\index{base point}
\index{null loop}

\begin{definition}
The homotopy group $\pi_1(X)$ of a connected space $X$ is the group of loops based at a given point $*\in X$, with the following conventions,
\begin{enumerate}
\item Two such loops are identified when one can pass continuously from one loop to the other, via a family of loops indexed by $t\in[0,1]$,

\item The composition of two such loops is the obvious one, namely is the loop obtained by following the first loop, then the second loop, 

\item The unit loop is the null loop at $*$, which stays there, and the inverse of a given loop is the loop itself, followed backwards,
\end{enumerate}
with the remark that the group $\pi_1(X)$ defined in this way does not depend on the choice of the given point $*\in X$, where the loops are based.
\end{definition}

This definition is obviously something non-trivial, based on some preliminary thinking on the subject, the technical details being as follows:

\bigskip

-- The fact that the set $\pi_1(X)$ defined as above is indeed a group is obvious, with all the group axioms being clear from definitions.

\bigskip

-- Obvious as well is the fact that, since $X$ is assumed to be connected, this group does not depend on the choice of the given point $*\in X$, where the loops are based.

\bigskip

As basic examples now, for spaces having ``no holes'', such as $\mathbb R$ itself, or $\mathbb R^N$, and so on, we have $\pi_1=\{1\}$. In fact, having no holes can only mean, by definition, $\pi_1=\{1\}$:

\begin{definition}
A space is called simply connected when:
$$\pi_1=\{1\}$$
That is, any loop inside our space must be contractible.
\end{definition}

So, this will be our starting definition, for the considerations in this section. As further illustrations for Definition 7.13, here are now a few basic computations:

\index{free group}

\begin{theorem}
We have the following computations of homotopy groups:
\begin{enumerate}
\item For the circle, we have $\pi_1=\mathbb Z$.

\item For the torus, we have $\pi_1=\mathbb Z\times\mathbb Z$.

\item For the disk minus $2$ points, we have $\pi_1=F_2$.

\item In fact, for the disk minus $N$ points, we have $\pi_1=F_N$.
\end{enumerate}
\end{theorem}

\begin{proof}
These results are all standard, as follows:

\medskip

(1) The first assertion is clear, because a loop on the circle must wind $n\in\mathbb Z$ times around the center, and this parameter $n\in\mathbb Z$ uniquely determines the loop, up to the identification in Definition 7.13. Thus, the homotopy group of the circle is the group of such parameters $n\in\mathbb Z$, which is of course the group $\mathbb Z$ itself.

\medskip

(2) In what regards now the second assertion, the torus being a product of two circles, we are led to the conclusion that its homotopy group must be some kind of product of $\mathbb Z$ with itself. But pictures show that the two standard generators of $\mathbb Z$, and so the two copies of $\mathbb Z$ themselves, commute, $gh=hg$, so we obtain the product of $\mathbb Z$ with itself, subject to commutation, which is the usual product $\mathbb Z\times\mathbb Z$:
$$\left<g,h\Big|gh=hg\right>=\mathbb Z\times\mathbb Z$$

It is actually instructive here to work out explicitly the proof of the commutation relation. We can use our usual drawing conventions for the torus, namely:
$$
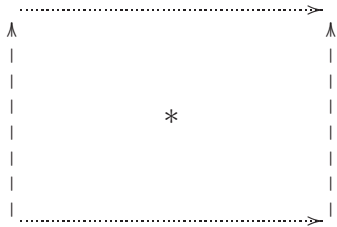
\xymatrix@R=30pt@C=50pt{
\ar@.[rr]&&\\
&\ast\\
\ar@{-->}[uu]\ar@.[rr]&&\ar@{-->}[uu]}$$

The standard generators $g,h$ of the homotopy group are then as follows:
$$
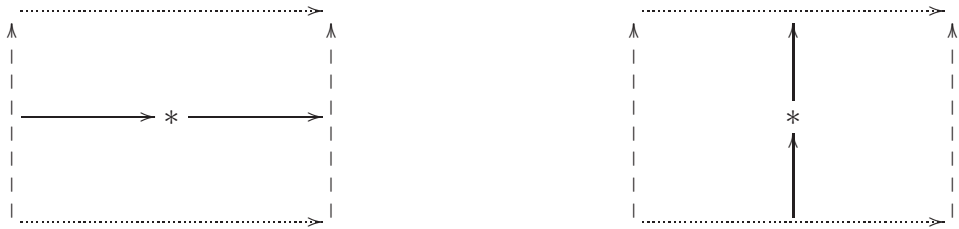
\xymatrix@R=30pt@C=50pt{
\ar@.[rr]&&&&\ar@.[rr]&&\\
\ar[r]&\ast\ar[r]&&&&\ast\ar[u]&\\
\ar@{-->}[uu]\ar@.[rr]&&\ar@{-->}[uu]&&\ar@{-->}[uu]\ar@.[rr]&\ar[u]&\ar@{-->}[uu]}$$

Regarding now the two compositions $gh,hg$, these are as follows:
$$
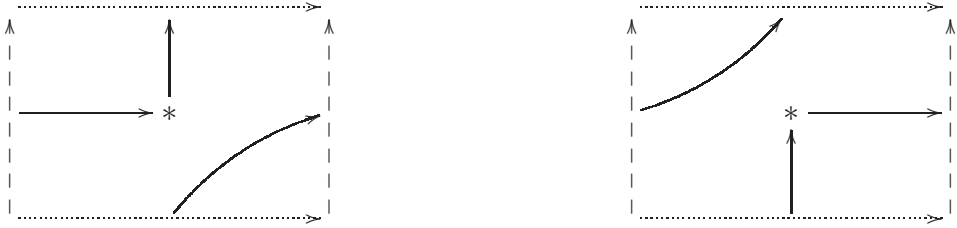
\xymatrix@R=30pt@C=50pt{
\ar@.[rr]&&&&\ar@.[rr]&&\\
\ar[r]&\ast\ar[u]&&&\ar@/_/[ur]&\ast\ar[r]&\\
\ar@{-->}[uu]\ar@.[rr]&\ar@/^/[ur]&\ar@{-->}[uu]&&\ar@{-->}[uu]\ar@.[rr]&\ar[u]&\ar@{-->}[uu]}$$

But these two pictures coincide, up to homotopy, with the following picture:
$$
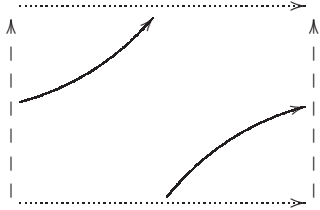
\xymatrix@R=30pt@C=50pt{
\ar@.[rr]&&\\
\ar@/_/[ur]&&\\
\ar@{-->}[uu]\ar@.[rr]&\ar@/^/[ur]&\ar@{-->}[uu]}$$

Thus we have indeed $gh=hg$, as desired, which gives the formula in (2).

\medskip

(3) Regarding now the disk minus $2$ points, the result here is quite clear, because the homotopy group is generated by the 2 loops around the 2 missing points, and these 2 loops are obviously free, algebrically speaking. Thus, we obtain a free product of the group $\mathbb Z$ with itself, which is by definition the free group on 2 generators $F_2$.

\medskip

(4) This is again clear, because the homotopy group is generated here by the $N$ loops around the $N$ missing points, which are free, algebrically speaking. Thus, we obtain a $N$-fold free product of $\mathbb Z$ with itself, which is the free group on $N$ generators $F_N$.
\end{proof}

We will be back to more such computations in chapter 8 below.

\bigskip

There are many other interesting things that can be said about homotopy groups. Also, another thing that can be done with the arbitrary topological spaces $X$, again in relation with studying their ``shape'', is that of looking at the fiber bundles over them, again up to continuous deformation. We are led in this way into a group, called $K_0(X)$. Moreover, both $\pi_1(X)$ and $K_0(X)$ have higher analogues $\pi_n(X)$ and $K_n(X)$ as well, and the general goal of algebraic topology is that of understanding all these groups, along with some other groups, of similar flavor, which can be constructed as well.

\bigskip

Moving ahead with some further topology, we have:

\index{planar graph}
\index{toral graph}
\index{genus of graph}
\index{Petersen graph}
\index{higher genus}
\index{Riemann surface}
\index{number of holes}
\index{algebraic topology}
\index{topology}

\begin{fact}
We can talk about the genus of a surface
$$g\in\mathbb N$$
as being its number of holes.
\end{fact}

To be more precise here, the genus of a sphere is by definition 0, then the genus of a torus is 1, then the genus of a doubly holed torus is 2, and so on.

\bigskip

At this point, I am sure that you are wondering what the genus exactly is, mathematically speaking. There are many answers here, ranging from elementary to advanced, depending on how much geometric you want to be, and the best answer, which is a bit complicated, involves complex analysis, and the notion of Riemann surface.

\section*{7d. Higher genus}

Getting back to graphs, with the above understood, we can now formulate the following key result, coming as a continuation and generalization of Theorem 7.7:

\begin{theorem}
For a connected graph of genus $g\in\mathbb N$ we have the Euler formula
$$v-e+f=2-2g$$
with $v,e,f$ being the number of vertices, edges and faces. 
\end{theorem}

\begin{proof}
This is something quite straightforward, as follows:

\medskip

(1) This comes as a continuation of Theorem 7.7, dealing with the case $g=0$, and assuming that you have worked out all the details of the proof there, you will certainly have no troubles now in understanding the present extension, to genus $g\in\mathbb N$.

\medskip

(2) So, instead of further insisting on the proof, which is something straightforward, let us see some illustrations, in genus $g\geq1$. Our simplest example of non-planar graph is the simplex $K_5$, and we have seen how to draw this graph on a torus, showing that the genus is $g=1$. But the picture that we used can now be recycled, in order to count the faces of $K_5$, when drawn on a torus, and there are 5 faces, as follows:
$$
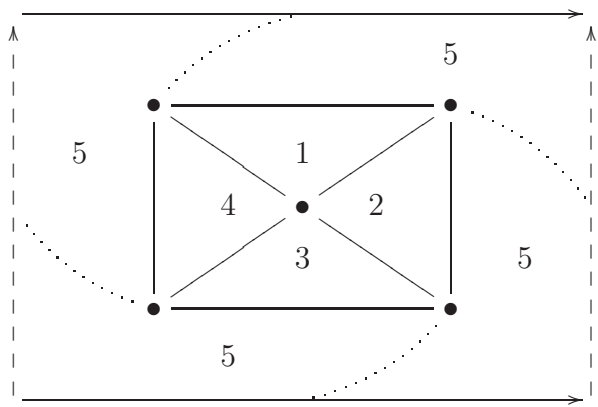
\xymatrix@R=5pt@C=15pt{
\ar[rrrrrrrr]&&&&&&&&\\
&&&&&&5\\
&&\bullet\ar@{-}[rrrr]\ar@{-}[dddd]\ar@/^/@{.}[uurr]&&&&\bullet\ar@{-}[dddd]\ar@/^/@{.}[ddrr]\\
&5&&&1\\
&&&4&\bullet\ar@{-}[uurr]\ar@{-}[uull]\ar@{-}[ddrr]\ar@{-}[ddll]&2&&&\\
&&&&3&&&5\\
&&\bullet\ar@{-}[rrrr]\ar@/^/@{.}[uull]&&&&\bullet\\
&&&5\\
\ar@{-->}[uuuuuuuu]\ar[rrrrrrrr]&&&&\ar@/_/@{.}[uurr]&&&&\ar@{-->}[uuuuuuuu]}$$

To be more precise, we have the four obvious faces $1,2,3,4$, and then everything else is in fact a single face 5, because by moving left and right, and up and down, according to our gluing conventions for the torus, we can travel from everywhere to everywhere, by avoiding that two edges, represented by dotted lines. Thus, we have 5 faces indeed, and as a conclusion to this, the Euler formula holds indeed for $K_5$, as:
$$5-10+5=0$$

(3) A similar discussion goes for the bipartite simplex $K_{3,3}$, with this graph being drawable on a torus, with 3 faces, as follows, with our various conventions above:
$$
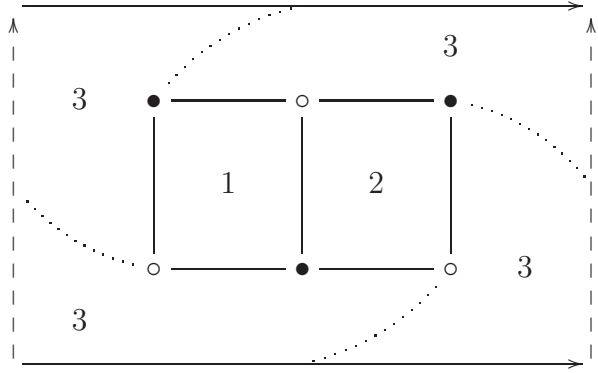
\xymatrix@R=5pt@C=15pt{
\ar[rrrrrrrr]&&&&&&&&\\
&&&&&&3\\
&3&\bullet\ar@{-}[rr]\ar@{-}[dddd]\ar@/^/@{.}[uurr]&&\circ\ar@{-}[rr]\ar@{-}[dddd]&&\bullet\ar@{-}[dddd]\ar@/^/@{.}[ddrr]\\
\\
&&&1&&2&&&\\
\\
&&\circ\ar@{-}[rr]\ar@/^/@{.}[uull]&&\bullet\ar@{-}[rr]&&\circ&3\\
&3\\
\ar@{-->}[uuuuuuuu]\ar[rrrrrrrr]&&&&\ar@/_/@{.}[uurr]&&&&\ar@{-->}[uuuuuuuu]}$$

Thus, again as good news, the Euler formula holds indeed for $K_{3,3}$, as:
$$6-9+3=0$$

(4) Regarding now the higher simplices $K_N$ and bipartite simplices $K_{M,N}$, things here are quite complicated, as already discussed in Theorem 7.4. To be more precise, in what regards the simplices, we know that $K_3,K_4$ have genus $g=0$, then $K_5,K_6,K_7$ have genus $g=1$, and then $K_8$ has genus $g=2$, and in all cases the Euler formula holds, as:
$$K_3\quad:\quad 3-3+2=2$$
$$K_4\quad:\quad 4-6+4=2$$
$$\ \, K_5\quad:\quad 5-10+5=0$$
$$\ \ K_6\quad:\quad 6-15+9=0$$
$$\ \ \ K_7\quad:\quad 7-21+14=0$$
$$\ \ \ \ \ K_8\quad:\quad 8-28+18=-2$$

In general, the verification of the Euler formula for the simplex $K_N$ is something quite complicated, which takes the following form, with $f$ being the number of faces, with this coming from the general genus formula mentioned in the proof of Theorem 7.4:
$$N-\binom{N}{2}+f=2-2\left\lceil\frac{(N-3)(N-4)}{12}\right\rceil$$

(5) As for the general bipartite simplices, $K_{M,N}$ with $M\leq N$, the situation here is quite similar. To start with, at $M=2$ these graphs are planar, as explained in the proof of Proposition 7.3, and the picture there shows that the Euler formula for them reads:
$$(N+2)-2N+N=2$$

At $M\geq3$ these graphs are no longer planar, and here are a few numerics for them, namely Euler formula, with input coming from Theorem 7.4:
$$K_{3,3}\quad:\quad 6-9+3=0$$
$$\ K_{3,4}\quad:\quad 7-12+5=0$$
$$\ \, K_{3,5}\quad:\quad 8-15+7=0$$
$$\ \, K_{3,6}\quad:\quad 9-18+9=0$$
$$\ \ \ \ \ \, K_{3,7}\quad:\quad 10-21+9=-2$$
$$\ \ K_{4,4}\quad:\quad 8-16+8=0$$
$$\ \ \ \ \, K_{4,5}\quad:\quad 9-20+9=-2$$

In general, the Euler formula holds indeed for $K_{M,N}$, in the following form, coming from the general genus formula mentioned in the proof of Theorem 7.4:
$$M+N-MN+f=2-2\left\lceil\frac{(M-2)(N-2)}{4}\right\rceil$$

(6) Summarizing, quite non-trivial all this, concrete applications of the Euler formula, and exercise of course for you, to learn more about such things.
\end{proof}

As another application of the notion of genus, and of the Euler formula, let us discuss now the Petersen graph. We have the following result, about it:

\begin{theorem}
The Petersen graph, namely
$$
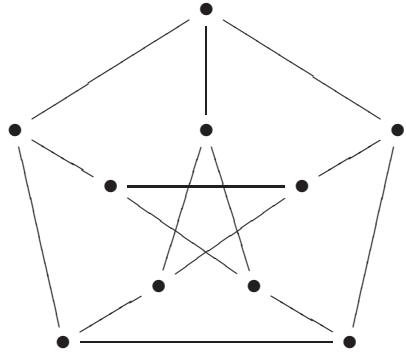
\xymatrix@R=1pt@C=5pt{
&&&&\bullet\ar@{-}[dddddrrrr]\ar@{-}[dddddllll]\\
\\
\\
\\
\\
\bullet\ar@{-}[ddddddddr]&&&&\bullet\ar@{-}[uuuuu]\ar@{-}[ddddddl]\ar@{-}[ddddddr]&&&&\bullet\ar@{-}[ddddddddl]\\
\\
&&
\bullet\ar@{-}[uull]\ar@{-}[ddddrrr]\ar@{-}[rrrr]&&&&\bullet\ar@{-}[uurr]\ar@{-}[ddddlll]\\
\\
\\
\\
&&&\bullet&&\bullet\\
\\
&\bullet\ar@{-}[rrrrrr]\ar@{-}[uurr]&&&&&&\bullet\ar@{-}[uull]}
$$
is toral, and the Euler formula for it reads $10-15+5=0$.
\end{theorem}

\begin{proof}
There are several things going on here, the idea being as follows:

\medskip

(1) The fact that this graph is indeed not planar can be best seen by using the Wagner criterion from Theorem 7.5, with both the graphs $K_5$ and $K_{3,3}$ being minors of it, and we have already talked about this, with full details, in Proposition 7.6. 

\medskip

(2) Regarding now the toral graph assertion, this requires some skill. By inverting the two pentagons, in the obvious way, the Petersen graph becomes as follows:
$$
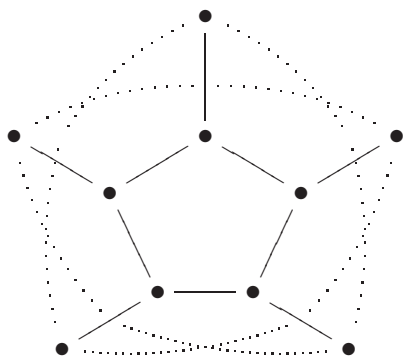
\xymatrix@R=1pt@C=5pt{
&&&&\bullet\ar@/^/@{.}[dddddrrr]\ar@/_/@{.}[dddddlll]\\
\\
\\
&&&&\ar@/^/@{.}[ddrrrr]\\
\\
\bullet\ar@/^/@{.}[uurrrr]\ar@/_/@{.}[dddddddrr]&\ar@/_/@{.}[dddddddd]&&&\bullet\ar@{-}[uuuuu]\ar@{-}[ddll]\ar@{-}[ddrr]&&&\ar@/^/@{.}[dddddddd]&\bullet\ar@/^/@{.}[dddddddll]\\
\\
&&\bullet\ar@{-}[uull]\ar@{-}[ddddr]&&&&\bullet\ar@{-}[uurr]\ar@{-}[ddddl]\\
\\
\\
\\
&&&\bullet\ar@{-}[rr]&&\bullet\\
&&\ar@/_/@{.}[drrrrr]&&&&\ar@/^/@{.}[dlllll]\\
&\bullet\ar@{-}[uurr]&&&&&&\bullet\ar@{-}[uull]}
$$

But now, we can keep the two pentagons and the solid edges, and send flying the various dotted edges, on suitable directions on a torus, as follows:
$$
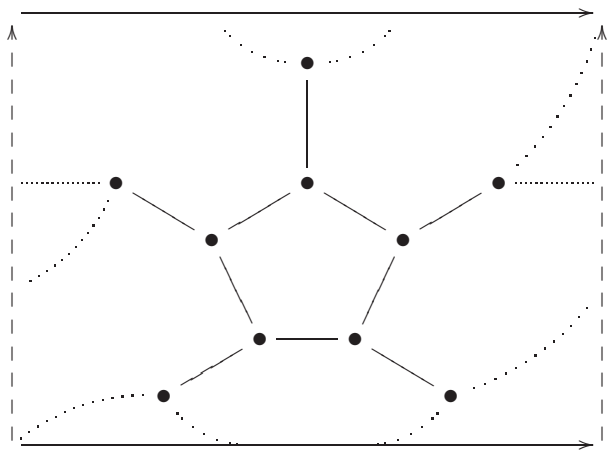
\xymatrix@R=1pt@C=5pt{
\ar[rrrrrrrrrrrrrr]&&&&&&&&&&&&&&\\
\\
&&&&&&&\bullet\ar@/_/@{.}[uurr]\ar@/^/@{.}[uull]\\
\\
\\
&&&&&&&\\
\\
\ar@{.}[rrr]&&&\bullet\ar@/^/@{.}[ddddlll]&&&&\bullet\ar@{-}[uuuuu]\ar@{-}[ddll]\ar@{-}[ddrr]&&&&\bullet\ar@/_/@{.}[uuuuuuurrr]\ar@{.}[rrr]&&&\\
\\
&&&&&\bullet\ar@{-}[uull]\ar@{-}[ddddr]&&&&\bullet\ar@{-}[uurr]\ar@{-}[ddddl]\\
\\
&&&&&&&&&&&&&&\ar@/^/@{.}[ddddllll]\\
\\
&&&&&&\bullet\ar@{-}[rr]&&\bullet\\
&&&&&&&&&\\
&&&&\bullet\ar@{-}[uurr]&&&&&&\bullet\ar@{-}[uull]\\
\\
\ar[rrrrrrrrrrrrrr]\ar@{-->}[uuuuuuuuuuuuuuuuu]\ar@/^/@{.}[uurrrr]&&&&&&\ar@/^/@{.}[uull]&&\ar@/_/@{.}[uurr]&&&&&&\ar@{-->}[uuuuuuuuuuuuuuuuu]
}$$

Observe the usage of the lower left vertex, which is identified with the upper right vertex, and in fact with the other two vertices of the rectangle as well, according to our gluing conventions for the torus. In any case, job done, and torality proved.

\medskip

(3) In order to finish, we still have to count the number of faces, in order to check the Euler formula. But there are 5 faces, as shown by the following picture:
$$
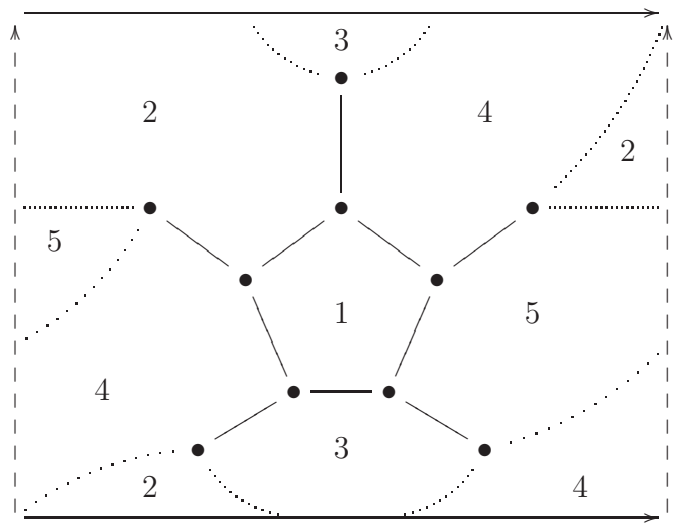
\xymatrix@R=0pt@C=5pt{
\ar[rrrrrrrrrrrrrr]&&&&&&&&&&&&&&\\
&&&&&&&3\\
&&&&&&&\bullet\ar@/_/@{.}[uurr]\ar@/^/@{.}[uull]\\
&&&2&&&&&&&4\\
&&&&&&&&&&&&&2\\
\\
\ar@{.}[rrr]&&&\bullet\ar@/^/@{.}[ddddlll]&&&&\bullet\ar@{-}[uuuu]\ar@{-}[ddll]\ar@{-}[ddrr]&&&&\bullet\ar@/_/@{.}[uuuuuurrr]\ar@{.}[rrr]&&&\\
&5\\
&&&&&\bullet\ar@{-}[uull]\ar@{-}[ddddr]&&&&\bullet\ar@{-}[uurr]\ar@{-}[ddddl]\\
&&&&&&&1&&&&5\\
&&&&&&&&&&&&&&\ar@/^/@{.}[ddddllll]\\
\\
&&4&&&&\bullet\ar@{-}[rr]&&\bullet\\
&&&&&&&&&\\
&&&&\bullet\ar@{-}[uurr]&&&3&&&\bullet\ar@{-}[uull]\\
&&&2&&&&&&&&&4\\
\ar[rrrrrrrrrrrrrr]\ar@{-->}[uuuuuuuuuuuuuuuu]\ar@/^/@{.}[uurrrr]&&&&&&\ar@/^/@{.}[uull]&&\ar@/_/@{.}[uurr]&&&&&&\ar@{-->}[uuuuuuuuuuuuuuuu]
}$$

Thus the Euler formula holds indeed in our embedding, as $10-15+5=0$.
\end{proof}

Many other things can be said about the Petersen graph, and other toral graphs, of similar type, or of general type. Exercise of course, read a bit here, if interested.

\section*{7e. Exercises}

We had yet another exciting chapter here, going into several directions, but as usual in such situations, as a downside to this, there are many exercises left for you:

\begin{exercise}
Prove that $K_6,K_7,K_{3,4}$ are toral.
\end{exercise}

\begin{exercise}
Prove that $K_{3,5},K_{3,6},K_{4,4}$ are toral.
\end{exercise}

\begin{exercise}
Prove that $K_8,K_{3,7},K_{4,5}$ are not toral.
\end{exercise}

\begin{exercise}
Read the full proofs of the Kuratowski and Wagner criteria.
\end{exercise}

\begin{exercise}
Learn more about regular polyhedra, and their various properties.
\end{exercise}

\begin{exercise}
Learn about genus, in all its flavors, including for Riemann surfaces.
\end{exercise}

\begin{exercise}
Work out the proof of the Euler formula, in general genus $g\in\mathbb N$.
\end{exercise}

\begin{exercise}
Learn about the genus of arbitrary simplices and bipartite simplices.
\end{exercise}

As a bonus exercise, learn some systematic algebraic topology, and more pure mathematics in general. All this material is old, very beautiful, and good to know.

\chapter{Knot invariants}

\section*{8a. Knots and links}

Leaving the graphs and related topological spaces aside, let us focus now on the simplest objects of topology, which are the knots. Knots are something very familiar, from the real life, and mathematically, it is most convenient to define them as follows:

\index{knot}
\index{link}
\index{random knot}
\index{tied knot}

\begin{definition}
A knot is a smooth closed curve in $\mathbb R^3$,
$$\gamma:\mathbb T\to\mathbb R^3$$
regarded modulo smooth transformations of $\mathbb R^3$.
\end{definition}

Observe that our knots are by definition oriented. The reverse knot $z\to\gamma(z^{-1})$ can be isomorphic or not to the original knot $z\to\gamma(z)$, and we will discuss this in a moment. At the level of examples, we first have the unknot, represented as follows:
$$
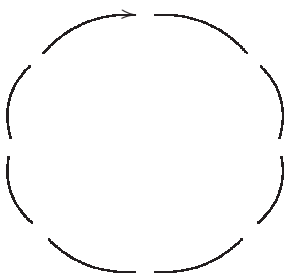
\xymatrix@R=9pt@C=3pt{
&&&&&\ar@{-}@/^/[drrrr]&\\
&\ar@{-}@/_/[ddl]\ar@/^/[urrrr]&&&&&&&&\ar@{-}@/^/[ddr]\\
&&&&&&&&&\\
\ar@{-}@/_/[ddr]&&&&&&&&&&&\\
&&&&&&&&&\\
&\ar@{-}@/_/[drrrr]&&&&&&&&\ar@{-}@/_/[uur]\\
&&&&&\ar@{-}@/_/[urrrr]}$$

For typographical reasons, it is most convenient to represent our knots by squarish diagrams, with these being far easier to type in Latex, the computer program used for writing math books, with the unknot for instance being represented as follows:
$$
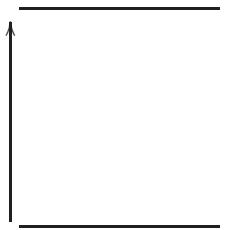
\xymatrix@R=75pt@C=75pt{
\ar@{-}[r]&\ar@{-}[d]\\
\ar[u]&\ar@{-}[l]}$$

The unknot is already a quite interesting mathematical object, suggesting a lot of exciting mathematical questions, for the most quite difficult, as follows:

\begin{questions}
In relation with the unknot:
\begin{enumerate}
\item Given a closed curve in $\mathbb R^3$, say given via its algebraic equations, can we decide if it is tied or not? 

\item Perhaps simpler, given the 2D picture of a knot, can we decide if the knot is tied or not?

\item Experience with cables and ropes shows that a random closed curve is usually tied. But, can we really prove this?
\end{enumerate}
\end{questions}

Obviously, difficult questions, and as you can see, knot theory is not an easy thing. But do not worry, we will manage to find our way through this jungle, and even come up with some mathematics for it. Going ahead now with examples, as the simplest possible true knot, meaning tied knot, we have the trefoil knot, which looks as follows:
$$
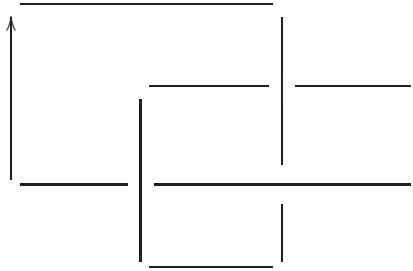
\xymatrix@R=18pt@C=40pt{
\ar@{-}[rr]&&\ar@{-}[dd]\\
&\ar@{-}[r]&|\ar@{-}[r]&\\
\ar@{-}[r]\ar[uu]&|\ar@{-}[rr]&-&\ar@{-}[u]\\
&\ar@{-}[r]\ar@{-}[uu]&\ar@{-}[u]
}$$

We also have the opposite trefoil knot, obtained by reversing the orientation, whose picture is identical to that of the trefoil knot, save for the orientation of the arrow:
$$
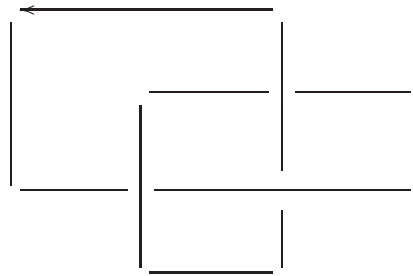
\xymatrix@R=18pt@C=40pt{
&&\ar[ll]\ar@{-}[dd]\\
&\ar@{-}[r]&|\ar@{-}[r]&\\
\ar@{-}[r]\ar@{-}[uu]&|\ar@{-}[rr]&-&\ar@{-}[u]\\
&\ar@{-}[r]\ar@{-}[uu]&\ar@{-}[u]
}$$

As before with the unknot, while the trefoil knot might look quite trivial, when it comes to formal mathematics regarding it, we are quickly led into delicate questions. Let us formulate a few intuitive observations about it, as follows:

\begin{fact}
In relation with the trefoil knot:
\begin{enumerate}
\item This knot is indeed tied, that is, not isomorphic to the unknot.

\item The trefoil knot and its opposite knot are not isomorphic.
\end{enumerate}
\end{fact}

To be more precise, here (1) is something which definitely holds, as we know it from real life, but if looking for a formal proof for this, based on Definition 8.1, we will certainly run into troubles. As for (2), here again we are looking for troubles, because when playing with two trefoil knots, made from rope, with opposite arrows marked on them, we certainly see that our two beasts are not identical, but go find a formal proof for that.

\bigskip

In short, as before with the unknot, modesty. For the moment, let us keep exploring the subject, by recording as Questions and Facts things that we see and feel, but cannot prove yet, mathematically, based on Definition 8.1 alone, due to a lack of tools.

\bigskip

Getting back now to Definition 8.1, as stated, it is convenient to allow, in relation with certain mathematical questions, links in our discussion:

\begin{definition}
A link is a collection of disjoint knots in $\mathbb R^3$, taken as usual oriented, and regarded as usual up to isotopy.
\end{definition}

As before with the knots, which can be truly knotted or not, there is a discussion here with respect to the links, which can be truly linked or not, and with orientation involved too. Drawing some pictures here, with some basic examples, is very instructive, the idea being that two or several basic unknots can be linked in many possible ways. For instance, as simplest non-trivial link, made of two unknots, which are indeed linked, we have: 
$$
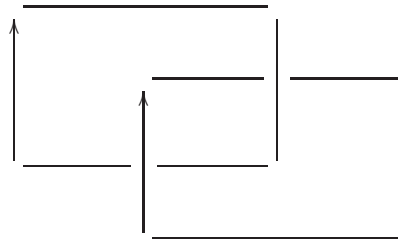
\xymatrix@R=14pt@C=40pt{
\ar@{-}[rr]&&\ar@{-}[dd]\\
&\ar@{-}[r]&|\ar@{-}[r]&\\
\ar@{-}[r]\ar[uu]&|\ar@{-}[r]&&\\
&\ar@{-}[rr]\ar[uu]&&\ar@{-}[uu]
}$$

By reversing the orientation of one unknot, we have as well the following link:
$$
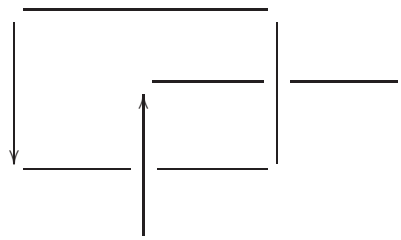
\xymatrix@R=14pt@C=40pt{
\ar@{-}[rr]\ar[dd]&&\ar@{-}[dd]\\
&\ar@{-}[r]&|\ar@{-}[r]&\\
\ar@{-}[r]&|\ar@{-}[r]&&\\
&\ar@{-}[rr]\ar[uu]&&\ar@{-}[uu]
}$$

This was for the story of two linked unknots, which is easily seen to stop here, with the above two links, but when trying to link $N$ unknots, with $N=3,4,5,\ldots\,$, many things can happen. Which leads us into the following philosophical question:

\begin{question}
Mathematically speaking, which are simpler to enumerate,
\begin{enumerate}
\item Usual knots, that is, links with one component,

\item Or links with several components, all being unknots,
\end{enumerate}
and this, in order to have some business started, for the links?
\end{question}

And with this being probably enough, as preliminary experimental work, time now to draw some conclusions. Obviously, what we have so far, namely Questions 8.2, Fact 8.3 and Question 8.5, is extremely interesting, at the core of everything that can be called ``geometry''. And by further thinking a bit, at how knots and links can be tied, in so many fascinating ways, we are led to the following philosophical conclusion:

\begin{conclusion}
Knots and links are to geometry and topology what prime numbers are to number theory.
\end{conclusion}

Very nice all this, we are now certainly motivated for studying the knots and links, and time for some mathematics. But the question is, how to get started? 

\bigskip

In view of the above, this is not an easy question. Fortunately, graph theory comes to the rescue, via to the following simple fact, which will be our starting point:

\index{projection of knot}
\index{tetravalent graph}
\index{crossings of knot}

\begin{fact}
The plane projection of a knot or link is something similar to an oriented graph with $4$-valent vertices, except for the fact that we have some extra data at the vertices, telling us, about the $2$ strings crossing there, which goes up and which goes down.
\end{fact}

Based on this, let us try to construct some knot invariants. A natural idea is that of defining the invariant on the 2D picture of the knot, that is, on a plane projection of the knot, and then proving that the invariant is indeed independent on the chosen plane. 

\bigskip

This method rests on the following technical result, which is well-known:

\index{Reidemeister moves}

\begin{theorem}
Two pictures correspond to plane projections of the same knot or link precisely when they differ by a sequence of Reidemeister moves, namely:
\begin{enumerate}
\item Moves of type I, given by $\propto\ \leftrightarrow\,|$.

\item Moves of type II, given by $)\hskip-1.7mm(\ \leftrightarrow\ )($.

\item Moves of type III, given by $\vartriangle\ \leftrightarrow\,\triangledown$.
\end{enumerate}
\end{theorem}

\begin{proof}
This is something very standard, as follows:

\medskip

(1) To start with, the Reidemeister moves of type I are by definition as follows:
$$
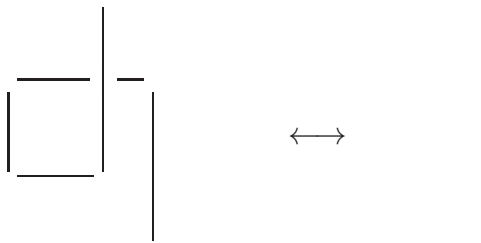
\xymatrix@R=6pt@C=10pt{
&&\ar@{-}[dddd]&&&&&&&\ar@{-}[dddddd]\\
\\
\ar@{-}[rr]\ar@{-}[dd]&&|\ar@{-}[r]&\ar@{-}[dddd]\\
&&&&&&\longleftrightarrow\\
\ar@{-}[rr]&&\\
\\
&&&&&&&&&
}$$

(2) Regarding the Reidemeister moves of type II, these are by definition as follows:
$$
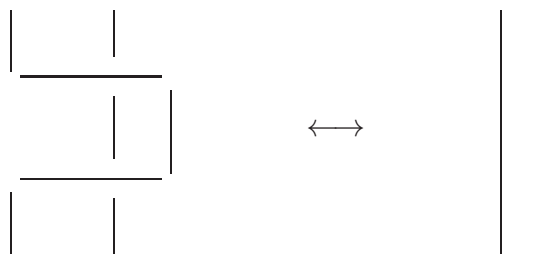
\xymatrix@R=6pt@C=10pt{
\ar@{-}[dd]&&\ar@{-}[dd]&&&&&&&\ar@{-}[dddddd]&\ar@{-}[dddddd]\\
\\
\ar@{-}[rrr]&&-\ar@{-}[dd]&\ar@{-}[dd]\\
&&&&&&\longleftrightarrow\\
\ar@{-}[dd]\ar@{-}[rrr]&&-\ar@{-}[dd]&\\
\\
&&&&&&&&&&
}$$

(3) As for the Reidemeister moves of type III, these are by definition as follows:
$$
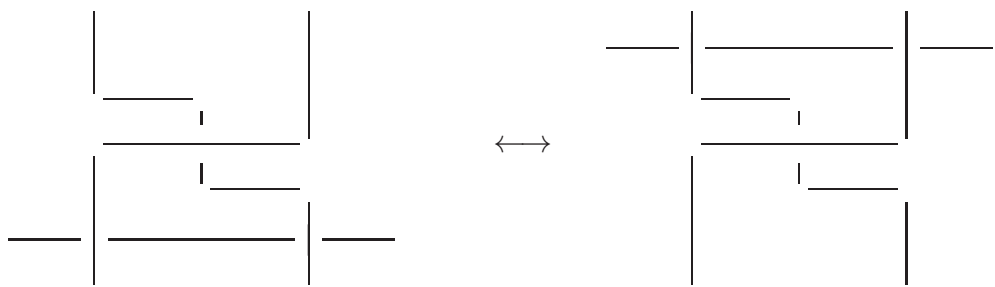
\xymatrix@R=6pt@C=10pt{
&&\ar@{-}[dd]&&&&\ar@{-}[ddd]&&&&&&&\ar@{-}[dd]&&&&\ar@{-}[ddd]\\
&&&&&&&&&&&\ar@{-}[rr]&&|\ar@{-}[rrrr]&&&&|\ar@{-}[rr]&&\\
&&\ar@{-}[rr]&&\ar@{-}[d]&&&&&&&&&\ar@{-}[rr]&&\ar@{-}[d]\\
&&\ar@{-}[rrrr]\ar@{-}[ddd]&&-\ar@{-}[d]&&&&&&\longleftrightarrow&&&\ar@{-}[rrrr]\ar@{-}[ddd]&&-\ar@{-}[d]&&\\
&&&&\ar@{-}[rr]&&\ar@{-}[dd]&&&&&&&&&\ar@{-}[rr]&&\ar@{-}[dd]\\
\ar@{-}[rr]&&|\ar@{-}[rrrr]&&&&|\ar@{-}[rr]&&\\
&&&&&&&&&&&&&&&&&
}$$

(4) This was for the precise statement of the theorem, and in what regards now the proof, this is somewhat clear from definitions, and in practice, this can be done by some sort of cut and paste procedure, or recurrence if you prefer, easy exercise for you.
\end{proof}

At a more advanced level now, we will need the following key observation, making the connection with group theory, and algebra in general, due to Alexander:

\index{braid group}

\begin{theorem}
Any knot or link can be thought of as being the closure of a braid,
$$
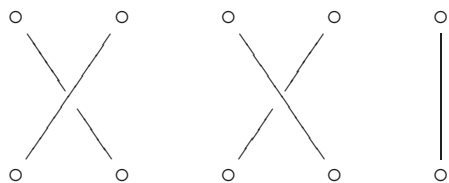
\xymatrix@C=10pt@R=20pt{
\circ\ar@{-}[dr]&&\circ\ar@{-}[ddll]&&\circ\ar@{-}[ddrr]&&\circ\ar@{-}[dl]&&\circ\ar@{-}[dd]\\
&\ar@{-}[dr]&&&&\ar@{-}[dl]&&&\\
\circ&&\circ&&\circ&&\circ&&\circ}
$$
with the braids forming a group $B_k$, called braid group.
\end{theorem}

\begin{proof}
Again, this is something quite self-explanatory, as follows:

\medskip

(1) Consider indeed the braids with $k$ strings, with the convention that things go from up to down. For instance the braid in the statement should be thought of as being:
$$
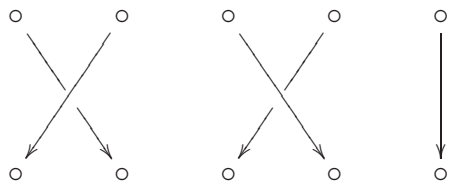
\xymatrix@C=10pt@R=20pt{
\circ\ar@{-}[dr]&&\circ\ar[ddll]&&\circ\ar[ddrr]&&\circ\ar@{-}[dl]&&\circ\ar[dd]\\
&\ar[dr]&&&&\ar[dl]&&&\\
\circ&&\circ&&\circ&&\circ&&\circ}
$$

But, with this convention, braids become some sort of permutations of $\{1,\ldots,k\}$, which are decorated at the level of crossings, with for instance the above braid corresponding to the following permutation of $\{1,2,3,4,5\}$, with due decorations: 
$$
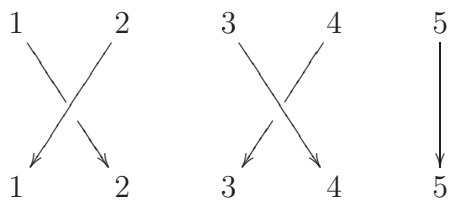
\xymatrix@C=10pt@R=20pt{
1\ar@{-}[dr]&&2\ar[ddll]&&3\ar[ddrr]&&4\ar@{-}[dl]&&5\ar[dd]\\
&\ar[dr]&&&&\ar[dl]&&&\\
1&&2&&3&&4&&5}
$$

In any case, we can see in this picture that $B_k$ is indeed a group, with composition law similar to that of the permutations in $S_k$, that is, going from up to down. 

\medskip

(2) Moreover, we can also see in this picture that we have a surjective group morphism $B_k\to S_k$, obtained by forgetting the decorations, at the level of crossings. For instance the braid pictured above is mapped in this way to the following permutation in $S_5$:
$$
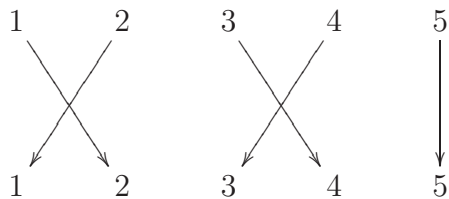
\xymatrix@C=10pt@R=20pt{
1\ar[ddrr]&&2\ar[ddll]&&3\ar[ddrr]&&4\ar[ddll]&&5\ar[dd]\\
&&&&&&&&\\
1&&2&&3&&4&&5}
$$

It is possible to do some more algebra here, in relation with the morphism $B_k\to S_k$, but we will not need this in what follows. We will keep in mind, from the above, the fact that ``braids are not exactly permutations, but they compose like permutations''.

\medskip

(3) Regarding now the closure operation in the statement, this consists by definition in adding semicircles at right, which makes our braid into a certain oriented link. As an illustration, the closure of the braid pictured above is the following link: 
$$
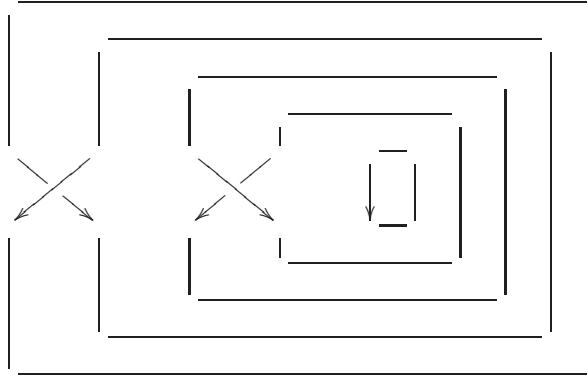
\xymatrix@C=10pt@R=7pt{
\ar@{-}[dddd]\ar@{-}[rrrrrrrrrrrrr]&&&&&&&&&&&&&\ar@{-}[dddddddddd]\\
&&\ar@{-}[ddd]\ar@{-}[rrrrrrrrrr]&&&&&&&&&&\ar@{-}[dddddddd]\\
&&&&\ar@{-}[dd]\ar@{-}[rrrrrrr]&&&&&&&\ar@{-}[dddddd]\\
&&&&&&\ar@{-}[d]\ar@{-}[rrrr]&&&&\ar@{-}[dddd]\\
\ar@{-}[dr]&&\ar[ddll]&&\ar[ddrr]&&\ar@{-}[dl]&&\ar[dd]\ar@{-}[r]&\ar@{-}[dd]\\
&\ar[dr]&&&&\ar[dl]&&&\\
&&&&&&&&\ar@{-}[r]&\\
&&&&&&\ar@{-}[u]\ar@{-}[rrrr]&&&&\\
&&&&\ar@{-}[uu]\ar@{-}[rrrrrrr]&&&&&&&\\
&&\ar@{-}[uuu]\ar@{-}[rrrrrrrrrr]&&&&&&&&&&\\
\ar@{-}[uuuu]\ar@{-}[rrrrrrrrrrrrr]&&&&&&&&&&&&&}
$$

(4) This was for the precise statement of the theorem, and in what regards now the proof, this can be done by some sort of cut and paste procedure, or recurrence if you prefer. As before with Theorem 8.8, we will leave this as an easy exercise for you.
\end{proof}

Many interesting things can be said about the braid group $B_k$, as for instance:

\begin{theorem}
The braid group $B_k$ has the following properties:
\begin{enumerate}
\item It is generated by variables $g_1,\ldots,g_{k-1}$, with the following relations:
$$g_ig_{i+1}g_i=g_{i+1}g_ig_{i+1}\quad,\quad 
g_ig_j=g_jg_i\ {\rm for}\ |i-j|\geq2$$

\item It is the homotopy group of $X=(\mathbb C^k-\Delta)/S_k$, with $\Delta\subset\mathbb C^k$ standing for the points $z$ satisfying $z_i=z_j$ for some $i\neq j$.
\end{enumerate}
\end{theorem}

\begin{proof}
These are things that we will not really need here, as follows:

\medskip

(1) In order to prove this assertion, due to Artin, consider the following braids:
$$\xymatrix@C=3pt@R=7pt{
&\circ\ar@{-}[dr]&&\circ\ar@{-}[ddll]&\circ\ar@{-}[dd]&\circ\ar@{-}[dd]&&\circ\ar@{-}[dd]&\circ\ar@{-}[dd]\\
g_1=&&\ar@{-}[dr]&&&&\ldots&\\
&\circ&&\circ&\circ&\circ&&\circ&\circ}
$$
$$\xymatrix@C=3pt@R=7pt{
&\circ\ar@{-}[dd]&\circ\ar@{-}[dr]&&\circ\ar@{-}[ddll]&\circ\ar@{-}[dd]&&\circ\ar@{-}[dd]&\circ\ar@{-}[dd]\\
g_2=&&&\ar@{-}[dr]&&&\ldots&\\
&\circ&\circ&&\circ&\circ&&\circ&\circ}
$$
$$\ \ \ \ \vdots$$
$$\xymatrix@C=3pt@R=7pt{
&\circ\ar@{-}[dd]&\circ\ar@{-}[dd]&&\circ\ar@{-}[dd]&\circ\ar@{-}[dd]&\circ\ar@{-}[dr]&&\circ\ar@{-}[ddll]\\
\ \ \ \ \ \ \ g_{k-1}=&&&\ldots&&&&\ar@{-}[dr]\\
&\circ&\circ&&\circ&\circ&\circ&&\circ&&&&}
$$

We have then $g_ig_j=g_jg_i$, for $|i-j|\geq2$. As for the relation $g_ig_{i+1}g_i=g_{i+1}g_ig_{i+1}$, by translation it is enough to check this at $i=1$. And here, we first have:
$$\xymatrix@C=3pt@R=7pt{
&\circ\ar@{-}[dr]&&\circ\ar@{-}[ddll]&&\circ\ar@{-}[dd]&\circ\ar@{-}[dd]&&\circ\ar@{-}[dd]&\circ\ar@{-}[dd]\\
&&\ar@{-}[dr]&&&&&\ldots&\\
&\circ\ar@{-}[dd]&&\circ\ar@{-}[dr]&&\circ\ar@{-}[ddll]&\circ\ar@{-}[dd]&&\circ\ar@{-}[dd]&\circ\ar@{-}[dd]\\
g_1g_2g_1=&&&&\ar@{-}[dr]&&&\ldots\\
&\circ\ar@{-}[dr]&&\circ\ar@{-}[ddll]&&\circ\ar@{-}[dd]&\circ\ar@{-}[dd]&&\circ\ar@{-}[dd]&\circ\ar@{-}[dd]\\
&&\ar@{-}[dr]&&&&&\ldots&\\
&\circ&&\circ&&\circ&\circ&&\circ&\circ}
$$

On the other hand, we have as well the following computation:
$$\xymatrix@C=3pt@R=7pt{
&\circ\ar@{-}[dd]&&\circ\ar@{-}[dr]&&\circ\ar@{-}[ddll]&\circ\ar@{-}[dd]&&\circ\ar@{-}[dd]&\circ\ar@{-}[dd]\\
&&&&\ar@{-}[dr]&&&\ldots\\
&\circ\ar@{-}[dr]&&\circ\ar@{-}[ddll]&&\circ\ar@{-}[dd]&\circ\ar@{-}[dd]&&\circ\ar@{-}[dd]&\circ\ar@{-}[dd]\\
g_2g_1g_2=&&\ar@{-}[dr]&&&&&\ldots&\\
&\circ\ar@{-}[dd]&&\circ\ar@{-}[dr]&&\circ\ar@{-}[ddll]&\circ\ar@{-}[dd]&&\circ\ar@{-}[dd]&\circ\ar@{-}[dd]\\
&&&&\ar@{-}[dr]&&&\ldots\\
&\circ&&\circ&&\circ&\circ&&\circ&\circ}$$

Now since the above two pictures are identical, up to isotopy, we have $g_1g_2g_1=g_2g_1g_2$, as desired. Thus, the braid group $B_k$ is indeed generated by elements $g_1,\ldots,g_{k-1}$ with the relations in the statement, and in what regards now the proof of universality, this can only be something quite routine, and we will leave this as an instructive exercise.

\medskip

(2) This is something quite self-explanatory, based on the general homotopy group material from chapter 7, and we will leave this as an easy exercise for you.

\medskip

(3) Finally, before leaving the subject, let us mention that the Artin relations in (1) are something very useful, in order to construct explicit matrix representations of $B_k$. For instance, it can be shown that the braid group $B_k$ is linear, and well, we will leave this as usual as an exercise for you, meaning either solve it, or look it up.
\end{proof}

\section*{8b. Temperley-Lieb}

Getting back now to knots and links, a quick comparison between our main results so far, namely Theorem 8.8 due to Reidemeister, and then Theorem 8.9 due to Alexander, suggests the following question, whose answer will certainly advance us:

\begin{question}
What is the analogue of the Reidemeister theorem, in the context of braids? That is, when do two braids produce, via closing, the same link?
\end{question} 

And this is, and we insist, a very good question, because assuming that we have an answer to it, no need afterwards to bother with plane projections, decorated graphs, Reidemeister moves, and amateurish topology in general, it will be all about groups and algebra. Which groups and algebra questions, you guessed right, we will eat them raw.

\bigskip

In answer now, we have the following theorem, due to Markov:

\begin{theorem}
Two elements of the full braid group, obtained as the increasing union of the various braid groups, with embeddings given by $\beta\to\beta\,|$,
$$B_\infty=\bigsqcup_{k=1}^\infty B_k$$
produce the same link, via closing, when one can pass from one to another via: 
\begin{enumerate}
\item Conjugation: $\beta\to\alpha\beta\alpha^{-1}$.

\item Markov move: $\beta\to g_k^{\pm1}\beta$.
\end{enumerate}
\end{theorem}

\begin{proof}
This is a version of the Reidemeister theorem, the idea being as follows:

\medskip

(1) To start with, it is clear that conjugating a braid, $\beta\to\alpha\beta\alpha^{-1}$, will produce the same link after closing, because we can pull the $\alpha,\alpha^{-1}$ to the right, in the obvious way, and there on the right, these $\alpha,\alpha^{-1}$ will annihilate, according to $\alpha\alpha^{-1}=1$.

\medskip

(2) Regarding now the Markov move from the statement, with $\beta\in B_k\subset B_{k+1}$ and with $g_1,\ldots,g_k\in B_{k+1}$ being the standard Artin generators, from Theorem 8.10 and its proof, this is the tricky move, which is worth a proof. Taking $k=3$ for an illustration, and representing $\beta\in B_3$ by a box, the link obtained by closing $g_4\beta$ is as follows, which is obviously the same link as the one obtained by closing $\beta$, and the same goes for $g_4^{-1}\beta$:
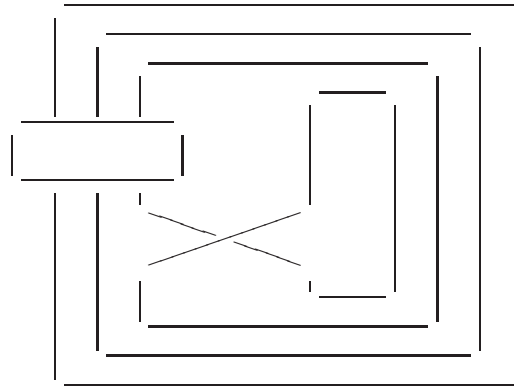
$$\xymatrix@C=9pt@R=4pt{
&\ar@{-}[dddd]\ar@{-}[rrrrrrrrrrr]&&&&&&&&&&&\ar@{-}[ddddddddddddd]\\
&&\ar@{-}[ddd]\ar@{-}[rrrrrrrrr]&&&&&&&&&\ar@{-}[ddddddddddd]\\
&&&\ar@{-}[dd]\ar@{-}[rrrrrrr]&&&&&&&\ar@{-}[ddddddddd]\\
&&&&&&&\ar@{-}[rr]\ar@{-}[dddd]&&\ar@{-}[ddddddd]\\
\ar@{-}[rrrr]\ar@{-}[dd]&&&&\ar@{-}[dd]\\
\\
\ar@{-}[rrrr]&\ar@{-}[ddddddd]&\ar@{-}[dddddd]&\ar@{-}[d]&&&\\
&&&\ar@{-}[drr]&&&&\ar@{-}[ddllll]\\
&&&&&\ar@{-}[drr]&\\
&&&\ar@{-}[dd]&&&&\ar@{-}[d]&\\
&&&&&&&\ar@{-}[rr]&&\\
&&&\ar@{-}[rrrrrrr]&&&&&&&\\
&&\ar@{-}[rrrrrrrrr]&&&&&&&&&\\
&\ar@{-}[rrrrrrrrrrr]&&&&&&&&&&&
}$$

(3) Thus, the links produced by braids are indeed invariant under the two moves in the statement. As for the proof of the converse, this comes from the Reidemeister theorem, applied in the context of the Alexander theorem, or perhaps simpler, by reasoning directly, a bit as in the proof of the Reidemeister theorem. We will leave this as an exercise.
\end{proof}

As explained before, the above kind of theorem is exactly what we need, in order to reformulate everything in terms of groups and algebra. To be more precise, looking now more in detail at what Theorem 8.12 exactly says, we are led to the following strategy:

\begin{strategy}
In order to construct numeric invariants for knots and links:
\begin{enumerate}
\item We must map $B_\infty$ somewhere, and then apply the trace.

\item And if the trace is preserved by Markov moves, it's a win.
\end{enumerate}
\end{strategy}

You get the point with all this, if we are do (1) then, by using the trace property $tr(ab)=tr(ba)$ of the trace, we will have $tr(\beta)=tr(\alpha\beta\alpha^{-1})$, in agreement with what Theorem 8.12 (1) requires. And if we do (2) too, whatever that condition exactly means, and more on this in a moment, we will have as well $tr(\beta)=tr(g_k^{\pm1}\beta)$, in agreement with what Theorem 8.12 (2) requires, so we will have our invariant for knots and links.

\bigskip

This sound very good, but before getting into details, let us be a bit megalomaniac, and add two more ambitious points to our war plan, as follows:

\begin{addendum}
Our victory will be total, with a highly reliable invariant, if:
\begin{enumerate}
\item The representation and trace are faithful as possible.

\item And they depend, if possible, on several parameters.
\end{enumerate}
\end{addendum}

Here (1) and (2) are obviously related, because the more parameters we have in (2), the more chances for our constructions in (1) to be faithful will be. In short, what we are wishing here for is an invariant which distinguishes well between various knots and links, and this can only come via a mixture of faithfulness, and parameters involved.

\bigskip

So long for the plan, and in practice now, getting back to what Strategy 8.13 says, we are faced right away with a problem, coming from the fact that $B_\infty$ is not that easy to represent. You might actually already know this, if you have struggled a bit with the exercise that I left for you, at the end of the previous section. So, we are led to:

\begin{question}
How to represent the braid group $B_\infty$?
\end{question}

So, this was the question that Reidemeister, Alexander, Markov, Artin and the others were fighting with, a long time ago, in the first half of the 20th century. Quite surprisingly, the answer to it came very late, in the 80s, from Jones \cite{jo2}, with inspiration from operator algebras, and more specifically, from his previous paper \cite{jo1} about subfactors.

\bigskip

Retrospectively looking at all this, what really matters in Jones' answer to Question 8.15 is the algebra constructed by Temperley and Lieb in \cite{tli}, in the context of questions from statistical mechanics. But then, by looking even more retrospectively at all this, we can even say that the answer to Question 8.15 comes from nothing at all, meaning basic category theory. So, this will be our approach in what follows, with our answer being:

\begin{answer}
Thinking well, $B_\infty$ is self-represented, without help from the outside.
\end{answer}

So, ready for some category theory? We first need objects, and our set of objects will be the good old $\mathbb N$. As for the arrows, somehow in relation with topology and braids, we will choose something very simple too, with our definition being as follows:

\index{meander}
\index{partition}
\index{noncrossing partition}
\index{Temperley-Lieb category}

\begin{definition}
The Temperley-Lieb category $TL_N$ has the positive integers $\mathbb N$ as objects, with the space of arrows $k\to l$ being the formal span
$$TL_N(k,l)=span(NC_2(k,l))$$
of noncrossing pairings between an upper row of $k$ points, and a lower row of $l$ points
$$
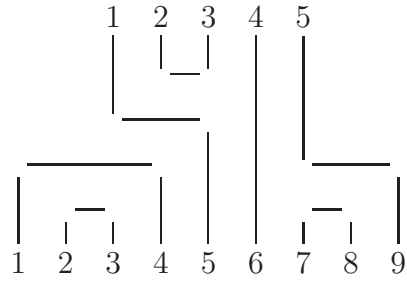
\xymatrix@R=10pt@C=5pt{
&&1\ar@{-}[dd]&2\ar@{-}[d]&3\ar@{-}[d]&4&5\ar@{-}[ddd]\\
&&&\ar@{-}[r]&&&&&\\
&&\ar@{-}[rr]&&&&&&\\
\ar@{-}[rrr]&&&&&&\ar@{-}[rr]&&\\
&\ar@{-}[r]&&&&&\ar@{-}[r]&&\\
1\ar@{-}[uu]&2\ar@{-}[u]&3\ar@{-}[u]&4\ar@{-}[uu]&5\ar@{-}[uuu]&6\ar@{-}[uuuuu]&7\ar@{-}[u]&8\ar@{-}[u]&9\ar@{-}[uu]
}$$
and with the composition of arrows appearing by composing the pairings, in the obvious way, with the rule $\bigcirc=N$, for the closed loops that might appear.
\end{definition}

This definition is something quite subtle, hiding several non-trivial things, and is worth a detailed discussion, our comments about it being as follows:

\bigskip

(1) First of all, our scalars in this chapter will be complex numbers, $\lambda\in\mathbb C$, and the ``formal span'' in the above must be understood in this sense, namely abstract complex vector space spanned by the elements of $NC_2(k,l)$. Of course it is possible to use an arbitrary field, at least at this stage of things, but remember that we are interested in quantum mechanics, and related mathematics, where the field of scalars is $\mathbb C$.

\bigskip

(2) Regarding the composition of arrows, this is by vertical concatenation, with our usual convention that things go ``from up to down''. And with this coming from care for our planet, and for entropy at the galactic level, I mean why pushing things from left to right, when we can have gravity work for us, pulling them from up to down:
$$\xymatrix@R=40pt@C=5pt{
up\ar[d]\\
down
}$$

(3) Less poetically, this ``from up to down'' convention is also useful for purely mathematical purposes, because the left-right direction will be reserved for the intervention of sums $\Sigma$ and scalars $\lambda\in\mathbb C$, while the up-down direction will be reserved for ``action''.

\bigskip

(4) Let us discuss now what happens with the closed circles, when concatenating. As an example, let us consider a full capping of noncrossing pairings, also called meander:
$$
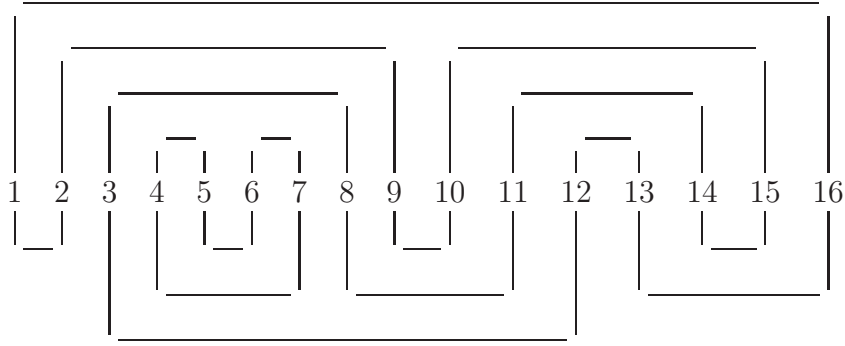
\xymatrix@R=10pt@C=5pt{
\ar@{-}[rrrrrrrrrrrrrrr]&&&&&&&&&&&&&&&\\
&\ar@{-}[rrrrrrr]&&&&&&&&\ar@{-}[rrrrr]&&&&&&\\
&&\ar@{-}[rrrrr]&&&&&&&&\ar@{-}[rrr]&&&&&\\
&&&\ar@{-}[r]&&\ar@{-}[r]&&&&&&\ar@{-}[r]&&&&\\
1\ar@{-}[uuuu]&2\ar@{-}[uuu]&3\ar@{-}[uu]&4\ar@{-}[u]&5\ar@{-}[u]&6\ar@{-}[u]&7\ar@{-}[u]&8\ar@{-}[uu]&9\ar@{-}[uuu]&10\ar@{-}[uuu]&11\ar@{-}[uu]&12\ar@{-}[u]&13\ar@{-}[u]&14\ar@{-}[uu]&15\ar@{-}[uuu]&16\ar@{-}[uuuu]\\
\ar@{-}[u]\ar@{-}[r]&\ar@{-}[u]&&&\ar@{-}[u]\ar@{-}[r]&\ar@{-}[u]&&&\ar@{-}[u]\ar@{-}[r]&\ar@{-}[u]&&&&\ar@{-}[u]\ar@{-}[r]&\ar@{-}[u]&\\
&&&\ar@{-}[uu]\ar@{-}[rrr]&&&\ar@{-}[uu]&\ar@{-}[uu]\ar@{-}[rrr]&&&\ar@{-}[uu]&&\ar@{-}[uu]\ar@{-}[rrr]&&&\ar@{-}[uu]\\
&&\ar@{-}[uuu]\ar@{-}[rrrrrrrrr]&&&&&&&&&\ar@{-}[uuu]&&&&
}$$

According to our conventions, this meander appears as the product $\pi\sigma\in NC_2(0,0)$ between the upper pairing $\sigma\in NC_2(0,16)$ and the lower pairing $\pi\in NC_2(16,0)$. But, what is the value of this product? We have two loops appearing, namely:
$$1-2-9-10-15-14-11-8-3-12-13-16$$
$$4-5-6-7$$

Thus, according to Definition 8.17, the value of this meander is $N^2$, with one $N$ for each  of the above loops, and with these two values of $N$ multiplying each other.

\bigskip

(5) The same discussion applies to an arbitrary composition $\pi\sigma\in NC_2(k,m)$ between an upper pairing $\sigma\in NC_2(k,l)$ and a lower pairing $\pi\in NC_2(l,m)$, with a certain number of loops appearing in this way, each contributing with a multiplicative factor $N$.

\bigskip

(6) Finally, in Definition 8.17 the value of the circle $N=\bigcirc$ can be pretty much anything, but due to some positivity reasons to become clear later, we will assume in what follows $N\in[1,\infty)$. Also, we will call this parameter $N$ the ``index'', with the precise reasons for calling this index to become clear later, too, as this book develops.

\bigskip

With all this discussed, what is next? More category theory I guess, and matter of having a theorem formulated too, instead of definitions only, let us formulate:

\index{tensor category}

\begin{theorem}
The Temperley-Lieb category $TL_N$ is a tensor $*$-category, with:
\begin{enumerate}
\item Composition of arrows: by vertical concatenation.

\item Tensoring of arrows: by horizontal concatenation.

\item Star operation: by turning the arrows upside-down.
\end{enumerate}
\end{theorem}

\begin{proof}
This is more of a definition, disguised as a theorem. To be more precise, we already know about (1), from Definition 8.17, and we can talk as well about (2) and (3), constructed as above, with (2) using of course multiplicativity with respect to the scalars, and with (3) using antimultiplicativity with respect to the scalars:
$$\left(\sum_i\lambda_i\pi_i\right)\otimes\left(\sum_j\mu_j\sigma_j\right)=\sum_{ij}\lambda_i\mu_j\pi_i\otimes\sigma_j$$
$$\left(\sum_i\lambda_i\pi_i\right)^*=\sum_i\bar{\lambda}_i\pi_i^*$$

And the point now is that our three operations are compatible with each other via all sorts of compatibility formulae, which are all clear from definitions, with the conclusion being that what we have a tensor $*$-category, as stated. We will leave the details here, basically amounting in figuring out what a tensor $*$-category exactly is, as an exercise.
\end{proof}

In order to further understand the category $TL_N$, let us focus on its diagonal part, formed by the End spaces of various objects. With the convention that these End spaces embed into each other by adding bars at right, this is a graded algebra, as follows:
$$\Delta TL_N=\bigcup_{k\geq0}TL_N(k,k)$$

Moreover, for further fine-tuning our study, let us actually focus on the individual components of this graded algebra. These components will play a key role in what follows, and they are worth a dedicated definition, and new notation and name, as follows:

\index{Temperley-Lieb algebra}

\begin{definition}
The Temperley-Lieb algebra $TL_N(k)$ is the formal span
$$TL_N(k)=span(NC_2(k,k))$$
with multiplication coming by concatenating, with the rule $\bigcirc=N$.
\end{definition}

In other words, $TL_N(k)$ appears as the formal span of the noncrossing pairings between an upper row of $k$ points, and a lower row of $k$ points, with multiplication coming by concatenating, with $\bigcirc=N$. As an example, here is a basis element of $TL_N(8)$:
$$
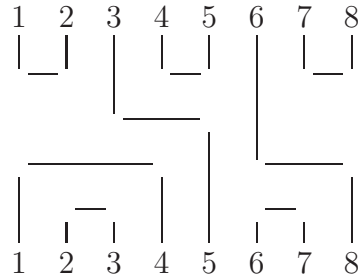
\xymatrix@R=10pt@C=5pt{
1\ar@{-}[d]&2\ar@{-}[d]&3\ar@{-}[dd]&4\ar@{-}[d]&5\ar@{-}[d]&6\ar@{-}[ddd]&7\ar@{-}[d]&8\ar@{-}[d]\\
\ar@{-}[r]&&&\ar@{-}[r]&&&\ar@{-}[r]&\\
&&\ar@{-}[rr]&&&&&\\
\ar@{-}[rrr]&&&&&\ar@{-}[rr]&&\\
&\ar@{-}[r]&&&&\ar@{-}[r]&&\\
1\ar@{-}[uu]&2\ar@{-}[u]&3\ar@{-}[u]&4\ar@{-}[uu]&5\ar@{-}[uuu]&6\ar@{-}[u]&7\ar@{-}[u]&8\ar@{-}[uu]
}$$

Getting back now to what we know about $TL_N$, from Theorem 8.18, the tensor product operation makes sense in the context of the diagonal algebra $\Delta TL_N$, but does not apply to its individual components $TL_N(k)$. However, the involution is useful, and we have:

\begin{proposition}
The Temperley-Lieb algebra $TL_N(k)$ is a $*$-algebra, with involution coming by turning the diagrams upside-down.
\end{proposition}

\begin{proof}
This is something trivial, which follows from Theorem 8.18, and can be verified as well directly, and we will leave this as an instructive exercise.
\end{proof}

Getting back now to knots and links, we first have to make the connection between braids and Temperley-Lieb diagrams. But this can be done as follows:

\begin{theorem}
The following happen:
\begin{enumerate}
\item We have a braid group representation $B_k\to TL_N(k)$, mapping standard generators to standard generators.

\item We have a trace $tr:TL_N(k)\to\mathbb C$, obtained by closing the diagrams, which is positive, and has a suitable Markov invariance property.
\end{enumerate}
\end{theorem}

\begin{proof}
Again, this is something quite intuitive, with the generators in (1) being by definition the standard ones, on both sides, and with the closing operation in (2) being similar to the one for braids, from Theorem 8.9. To be more precise:

\medskip

(1) The idea here is to map the Artin generators of the braid group to suitable modifications of the following Temperley-Lieb diagrams, called Jones projections:
$$\xymatrix@C=3pt@R=7pt{
&\circ\ar@/_/@{-}[rr]&&\circ&\circ\ar@{-}[dd]&\circ\ar@{-}[dd]&&\circ\ar@{-}[dd]&\circ\ar@{-}[dd]\\
e_1=&&&&&&\ldots&\\
&\circ\ar@/^/@{-}[rr]&&\circ&\circ&\circ&&\circ&\circ}
$$
$$\xymatrix@C=3pt@R=7pt{
&\circ\ar@{-}[dd]&\circ\ar@/_/@{-}[rr]&&\circ&\circ\ar@{-}[dd]&&\circ\ar@{-}[dd]&\circ\ar@{-}[dd]\\
e_2=&&&&&&\ldots&\\
&\circ&\circ\ar@/^/@{-}[rr]&&\circ&\circ&&\circ&\circ}
$$
$$\ \ \ \ \vdots$$
$$\xymatrix@C=3pt@R=7pt{
&\circ\ar@{-}[dd]&\circ\ar@{-}[dd]&&\circ\ar@{-}[dd]&\circ\ar@{-}[dd]&\circ\ar@/_/@{-}[rr]&&\circ\\
\ \ \ \ \ \ \ e_{k-1}=&&&\ldots&&&&\\
&\circ&\circ&&\circ&\circ&\circ\ar@/^/@{-}[rr]&&\circ&&&&}
$$

As a first observation, these diagrams satisfy $e_i^2=Ne_i$, with $N=\bigcirc$ being as usual the value of the circle, so it is rather the rescaled versions $f_i=e_i/N$ which are projections, but we will not bother with this, and use our terminology above. Next, our Jones projections certainly satisfy the Artin relations $e_ie_j=e_je_i$, for $|i-j|\geq2$. Our claim now is that is that we have as well the formula $e_ie_{i\pm1}e_i=e_i$. Indeed, by translation it is enough to check $e_ie_{i+1}e_i=e_i$ at $i=1$, and this follows from the following computation:
$$\xymatrix@C=3pt@R=7pt{
&\circ\ar@/_/@{-}[rr]&&\circ&&\circ\ar@{-}[dd]&\circ\ar@{-}[dd]&&\circ\ar@{-}[dd]&\circ\ar@{-}[dd]\\
&&&&&&&\ldots&\\
&\circ\ar@{-}[dd]\ar@/^/@{-}[rr]&&\circ\ar@/_/@{-}[rr]&&\circ&\circ\ar@{-}[dd]&&\circ\ar@{-}[dd]&\circ\ar@{-}[dd]\\
e_1e_2e_1=&&&&&&&\ldots&&&=e_1\\
&\circ\ar@/_/@{-}[rr]&&\circ\ar@/^/@{-}[rr]&&\circ\ar@{-}[dd]&\circ\ar@{-}[dd]&&\circ\ar@{-}[dd]&\circ\ar@{-}[dd]\\
&&&&&&&\ldots&\\
&\circ\ar@/^/@{-}[rr]&&\circ&&\circ&\circ&&\circ&\circ}
$$

As for the verification of the relation $e_2e_1e_2=e_2$, this is similar, as follows:
$$\xymatrix@C=3pt@R=7pt{
&\circ\ar@{-}[dd]&&\circ\ar@/_/@{-}[rr]&&\circ&\circ\ar@{-}[dd]&&\circ\ar@{-}[dd]&\circ\ar@{-}[dd]\\
&&&&&&&\ldots\\
&\circ\ar@/_/@{-}[rr]&&\circ\ar@/^/@{-}[rr]&&\circ\ar@{-}[dd]&\circ\ar@{-}[dd]&&\circ\ar@{-}[dd]&\circ\ar@{-}[dd]\\
e_2e_1e_2=&&&&&&&\ldots&&&=e_2\\
&\circ\ar@/^/@{-}[rr]&&\circ\ar@/_/@{-}[rr]&&&\circ\ar@{-}[dd]&&\circ\ar@{-}[dd]&\circ\ar@{-}[dd]\\
&&&&&&&\ldots\\
&\circ\ar@{-}[uu]&&\circ\ar@/^/@{-}[rr]&&\circ&\circ&&\circ&\circ}$$

Now with the relations $e_ie_{i\pm1}e_i=e_i$ in hand, let us try to reach to the Artin relations $g_ig_{i+1}g_i=g_{i+1}g_ig_{i+1}$. For this purpose, let us set $g_i=te_i-1$. We have then:
\begin{eqnarray*}
g_ig_{i+1}g_i
&=&(te_i-1)(te_{i+1}-1)(te_i-1)\\
&=&t^3e_i-t^2(Ne_i+e_ie_{i+1}+e_{i+1}e_i)+t(2e_i+e_{i+1})-1\\
&=&t(t^2-Nt+2)e_i+te_{i+1}-t^2(e_ie_{i+1}+e_{i+1}e_i)
\end{eqnarray*}

On the other hand, we have as well the following computation:
\begin{eqnarray*}
g_{i+1}g_ig_{i+1}
&=&(te_{i+1}-1)(te_i-1)(te_{i+1}-1)\\
&=&t^3e_{i+1}-t^2(Ne_{i+1}+e_ie_{i+1}+e_{i+1}e_i)+t(2e_{i+1}+e_i)-1\\
&=&t(t^2-Nt+2)e_{i+1}+te_i-t^2(e_ie_{i+1}+e_{i+1}e_i)
\end{eqnarray*}

Thus with $t^2-Nt+1=0$ we have a representation $B_k\to TL_N(k)$, as desired.

\medskip

(2) This is something more subtle, especially in what regards the positivity properties of the trace $tr:TL_N(k)\to\mathbb C$, which requires a bit more mathematics. So, no hurry with this, and we will discuss all this, and applications, in the remainder of this chapter.
\end{proof}

\section*{8c. Meanders, trace} 

Let us discuss now the positivity property of the Temperley-Lieb trace, constructed as indicated in Theorem 8.21. This is something quite subtle, which in the operator algebra context, that of the original paper of Jones \cite{jo2}, comes for free, or almost, and more on this in chapter 16 below. In the meantime, we will present a more pedestrian approach to the question, based on pure combinatorics, due to Di Francesco \cite{dif}.

\bigskip

The positivity will come from a systematic study of the partitions. Let us start with:

\index{lattice of partitions}
\index{order of partitions}
\index{supremum of partitions}

\begin{definition}
Let $P(k)$ be the set of partitions of $\{1,\ldots,k\}$, and $\pi,\sigma\in P(k)$.
\begin{enumerate}
\item We write $\pi\leq\sigma$ if each block of $\pi$ is contained in a block of $\sigma$.

\item We let $\pi\vee\sigma\in P(k)$ be the partition obtained by superposing $\pi,\sigma$.
\end{enumerate}
Also, we denote by $|.|$ the number of blocks of the partitions $\pi\in P(k)$.
\end{definition}

As an illustration here, at $k=2$ we have $P(2)=\{||,\sqcap\}$, and we have:
$$||\leq\sqcap$$

Also, at $k=3$ we have $P(3)=\{|||,\sqcap|,\sqcap\hskip-3.2mm{\ }_|\,,|\sqcap,\sqcap\hskip-0.7mm\sqcap\}$, and the order relation is as follows:
$$|||\ \leq\ \sqcap|\ ,\ \sqcap\hskip-3.2mm{\ }_|\ ,\ |\sqcap\ \leq\ \sqcap\hskip-0.7mm\sqcap$$

In order to study the Gram matrix $G_k(\pi,\sigma)=N^{|\pi\vee\sigma|}$, and more specifically to compute its determinant, we will use several standard facts about the partitions. We have:

\index{Gram matrix}
\index{M\"obius function}

\begin{definition}
The M\"obius function of any lattice, and so of $P$, is given by
$$\mu(\pi,\sigma)=\begin{cases}
1&{\rm if}\ \pi=\sigma\\
-\sum_{\pi\leq\tau<\sigma}\mu(\pi,\tau)&{\rm if}\ \pi<\sigma\\
0&{\rm if}\ \pi\not\leq\sigma
\end{cases}$$
with the construction being performed by recurrence.
\end{definition}

As an illustration here, for $P(2)=\{||,\sqcap\}$, we have by definition:
$$\mu(||,||)=\mu(\sqcap,\sqcap)=1$$

Also, $||<\sqcap$, with no intermediate partition in between, so we obtain:
$$\mu(||,\sqcap)=-\mu(||,||)=-1$$

Finally, we have $\sqcap\not\leq||$, and so we have as well the following formula:
$$\mu(\sqcap,||)=0$$

Thus, as a conclusion, we have computed the M\"obius matrix $M_2(\pi,\sigma)=\mu(\pi,\sigma)$ of the lattice $P(2)=\{||,\sqcap\}$, the formula being as follows:
$$M_2=\begin{pmatrix}1&-1\\ 0&1\end{pmatrix}$$

Back to the general case now, the main interest in the M\"obius function comes from the M\"obius inversion formula, which states that the following happens:
$$f(\sigma)=\sum_{\pi\leq\sigma}g(\pi)\quad
\implies\quad g(\sigma)=\sum_{\pi\leq\sigma}\mu(\pi,\sigma)f(\pi)$$

In linear algebra terms, the statement and proof of this formula are as follows:

\index{M\"obius inversion}

\begin{theorem}
The inverse of the adjacency matrix of $P(k)$, given by
$$A_k(\pi,\sigma)=\begin{cases}
1&{\rm if}\ \pi\leq\sigma\\
0&{\rm if}\ \pi\not\leq\sigma
\end{cases}$$
is the M\"obius matrix of $P$, given by $M_k(\pi,\sigma)=\mu(\pi,\sigma)$.
\end{theorem}

\begin{proof}
This is well-known, coming for instance from the fact that $A_k$ is upper triangular. Indeed, when inverting, we are led into the recurrence from Definition 8.23.
\end{proof}

Now back to our Gram matrix considerations, we have the following key result:

\begin{proposition}
The Gram matrix $G_{\pi\sigma}=N^{|\pi\vee\sigma|}$ decomposes as a product of upper/lower triangular matrices, $G_k=A_kL_k$, where
$$L_k(\pi,\sigma)=
\begin{cases}
N(N-1)\ldots(N-|\pi|+1)&{\rm if}\ \sigma\leq\pi\\
0&{\rm otherwise}
\end{cases}$$
and where $A_k$ is the adjacency matrix of $P(k)$.
\end{proposition}

\begin{proof}
We have indeed the following computation:
\begin{eqnarray*}
G_k(\pi,\sigma)
&=&N^{|\pi\vee\sigma|}\\
&=&\#\left\{i_1,\ldots,i_k\in\{1,\ldots,N\}\Big|\ker i\geq\pi\vee\sigma\right\}\\
&=&\sum_{\tau\geq\pi\vee\sigma}\#\left\{i_1,\ldots,i_k\in\{1,\ldots,N\}\Big|\ker i=\tau\right\}\\
&=&\sum_{\tau\geq\pi\vee\sigma}N(N-1)\ldots(N-|\tau|+1)
\end{eqnarray*}

According now to the definition of $A_k,L_k$, this formula reads:
\begin{eqnarray*}
G_k(\pi,\sigma)
&=&\sum_{\tau\geq\pi}L_k(\tau,\sigma)\\
&=&\sum_\tau A_k(\pi,\tau)L_k(\tau,\sigma)\\
&=&(A_kL_k)(\pi,\sigma)
\end{eqnarray*}

Thus, we are led to the formula in the statement.
\end{proof}

As an illustration for the above result, at $k=2$ we have $P(2)=\{||,\sqcap\}$, and the above decomposition $G_2=A_2L_2$ appears as follows:
$$\begin{pmatrix}N^2&N\\ N&N\end{pmatrix}
=\begin{pmatrix}1&1\\ 0&1\end{pmatrix}
\begin{pmatrix}N^2-N&0\\N&N\end{pmatrix}$$

We are led in this way to the following formula, due to Lindst\"om \cite{lin}:

\index{Gram matrix}
\index{Gram determinant}
\index{Lindst\"om formula}
\index{linear independence}

\begin{theorem}
The determinant of the Gram matrix $G_k$ is given by
$$\det(G_k)=\prod_{\pi\in P(k)}\frac{N!}{(N-|\pi|)!}$$
with the convention that in the case $N<k$ we obtain $0$.
\end{theorem}

\begin{proof}
If we order $P(k)$ as usual, with respect to the number of blocks, and then lexicographically, $A_k$ is upper triangular, and $L_k$ is lower triangular. Thus, we have:
\begin{eqnarray*}
\det(G_k)
&=&\det(A_k)\det(L_k)\\
&=&\det(L_k)\\
&=&\prod_\pi L_k(\pi,\pi)\\
&=&\prod_\pi N(N-1)\ldots(N-|\pi|+1)
\end{eqnarray*}

Thus, we are led to the formula in the statement.
\end{proof}

Getting now to what we wanted to do, namely computation of the Gram determinant for the lattice of noncrossing pairings, which will provide us with the desired positivity properties of the Temperley-Lieb trace, let us begin with some examples. We will need:

\index{fattening of partitions}
\index{shrinking of partitions}

\begin{proposition}
We have a bijection $NC(k)\simeq NC_2(2k)$, constructed by fattening and shrinking, as follows:
\begin{enumerate}
\item The application $NC(k)\to NC_2(2k)$ is the ``fattening'' one, obtained by doubling all the legs, and doubling all the strings too.

\item Its inverse $NC_2(2k)\to NC(k)$ is the ``shrinking'' application, obtained by collapsing pairs of consecutive neighbors.
\end{enumerate}
\end{proposition}

\begin{proof}
This is something self-explanatory, and in order to see how this works, let us discuss an example. Consider a noncrossing partition, say the following one:
$$
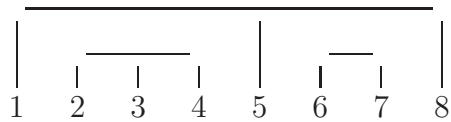
\xymatrix@R=10pt@C=10pt{
\ar@{-}[rrrrrrr]&&&&&&&\\
&\ar@{-}[rr]&&&&\ar@{-}[r]&&\\
1\ar@{-}[uu]&2\ar@{-}[u]&3\ar@{-}[u]&4\ar@{-}[u]&5\ar@{-}[uu]&6\ar@{-}[u]&7\ar@{-}[u]&8\ar@{-}[uu]
}$$

Now let us ``fatten'' this partition, by doubling everything, as follows:
$$
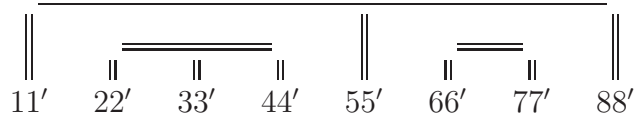
\xymatrix@R=10pt@C=10pt{
\ar@{=}[rrrrrrr]&&&&&&&\\
&\ar@{=}[rr]&&&&\ar@{=}[r]&&\\
11'\ar@{=}[uu]&22'\ar@{=}[u]&33'\ar@{=}[u]&44'\ar@{=}[u]&55'\ar@{=}[uu]&66'\ar@{=}[u]&77'\ar@{=}[u]&88'\ar@{=}[uu]
}$$

Now by relabeling the points $1,\ldots,16$, what we have is indeed a noncrossing pairing. As for the reverse operation, that is obviously obtained by ``shrinking'' our pairing, by collapsing pairs of consecutive neighbors, that is, by identifying $1=2$, then $3=4$, then $5=6$, and so on, up to $15=16$. Thus, we are led to the conclusion in the statement.
\end{proof}

At the level of the associated Gram matrices, the result is as follows:

\index{Gram matrix}

\begin{proposition}
The Gram matrices of $NC_2(2k)\simeq NC(k)$ are related by
$$G_{2k,n}(\pi,\sigma)=n^k(\Delta_{kn}^{-1}G_{k,n^2}\Delta_{kn}^{-1})(\pi',\sigma')$$
where $\pi\to\pi'$ is the shrinking operation, and $\Delta_{kn}$ is the diagonal of $G_{kn}$.
\end{proposition}

\begin{proof}
In the context of the bijection from Proposition 8.27, we have:
$$|\pi\vee\sigma|=k+2|\pi'\vee\sigma'|-|\pi'|-|\sigma'|$$

We therefore have the following formula, valid for any $n\in\mathbb N$:
$$n^{|\pi\vee\sigma|}=n^{k+2|\pi'\vee\sigma'|-|\pi'|-|\sigma'|}$$

Thus, we are led to the formula in the statement.
\end{proof}

Getting back now to our business, namely computation of the Gram determinant for the lattice of noncrossing pairings, we first have the following elementary result:

\begin{proposition}
The first Gram matrices and determinants for $NC_2$ are
$$\det\begin{pmatrix}N^2&N\\N&N^2\end{pmatrix}=N^2(N^2-1)$$
$$\det\begin{pmatrix}
N^3&N^2&N^2&N^2&N\\
N^2&N^3&N&N&N^2\\
N^2&N&N^3&N&N^2\\
N^2&N&N&N^3&N^2\\
N&N^2&N^2&N^2&N^3
\end{pmatrix}=N^5(N^2-1)^4(N^2-2)$$
with the matrices being written by using the lexicographic order on $NC_2(2k)$.
\end{proposition}

\begin{proof}
The formula at $k=2$, where $NC_2(4)=\{\sqcap\sqcap,\bigcap\hskip-4.9mm{\ }_\cap\,\}$, is clear. At $k=3$ however, things are tricky. We have $NC(3)=\{|||,\sqcap|,\sqcap\hskip-3.2mm{\ }_|\,,|\sqcap,\sqcap\hskip-0.7mm\sqcap\}$, and the corresponding Gram matrix and its determinant are, according to Theorem 8.26:
$$\det\begin{pmatrix}
N^3&N^2&N^2&N^2&N\\
N^2&N^2&N&N&N\\
N^2&N&N^2&N&N\\
N^2&N&N&N^2&N\\
N&N&N&N&N
\end{pmatrix}=N^5(N-1)^4(N-2)$$

By using Proposition 8.28, the Gram determinant of $NC_2(6)$ is given by:
\begin{eqnarray*}
\det(G_{6N})
&=&\frac{1}{N^2\sqrt{N}}\times N^{10}(N^2-1)^4(N^2-2)\times\frac{1}{N^2\sqrt{N}}\\
&=&N^5(N^2-1)^4(N^2-2)
\end{eqnarray*}

Thus, we have obtained the formula in the statement.
\end{proof}

In general, such tricks won't work, because $NC(k)$ is strictly smaller than $P(k)$ at $k\geq4$. However, following Di Francesco \cite{dif}, we have the following result:

\index{meander determinant}
\index{Gram determinant}

\begin{theorem}
The determinant of the Gram matrix for $NC_2$ is given by
$$\det(G_{kN})=\prod_{r=1}^{[k/2]}P_r(N)^{d_{k/2,r}}$$
where $P_r$ are the Chebycheff polynomials, given by
$$P_0=1\quad,\quad 
P_1=X\quad,\quad 
P_{r+1}=XP_r-P_{r-1}$$
and $d_{kr}=f_{kr}-f_{k,r+1}$, with $f_{kr}$ being the following numbers, depending on $k,r\in\mathbb Z$,
$$f_{kr}=\binom{2k}{k-r}-\binom{2k}{k-r-1}$$
with the convention $f_{kr}=0$ for $k\notin\mathbb Z$. 
\end{theorem}

\begin{proof}
This is something quite technical, obtained by using a decomposition as follows of the Gram matrix $G_{kN}$, with the matrix $T_{kN}$ being lower triangular:
$$G_{kN}=T_{kN}T_{kN}^t$$

Thus, a bit as in the proof of the Lindst\"om formula, we obtain the result, but the problem lies however in the construction of $T_{kN}$, which is non-trivial. See \cite{dif}.
\end{proof}

Let us record as well the following result, also from Di Francesco \cite{dif}:

\index{meander determinant}
\index{Gram determinant}

\begin{theorem}
The determinant of the Gram matrix for $NC$ is given by
$$\det(G_{kN})=(\sqrt{N})^{a_k}\prod_{r=1}^kP_r(\sqrt{N})^{d_{kr}}$$
where $P_r$ are the Chebycheff polynomials, given by
$$P_0=1\quad,\quad 
P_1=X\quad,\quad 
P_{r+1}=XP_r-P_{r-1}$$
and $d_{kr}=f_{kr}-f_{k,r+1}$, with $f_{kr}$ being the following numbers, depending on $k,r\in\mathbb Z$,
$$f_{kr}=\binom{2k}{k-r}-\binom{2k}{k-r-1}$$
with the convention $f_{kr}=0$ for $k\notin\mathbb Z$, and where $a_k=\sum_{\pi\in P(k)}(2|\pi|-k)$.
\end{theorem}

\begin{proof}
This follows indeed from Theorem 8.30, by using Proposition 8.28.
\end{proof}

We refer to the literature for more on the above, which is first class combinatorics.

\section*{8d. Jones polynomial} 

Now back to the knots and links, we have all the needed ingredients. Indeed, we can now put everything together, and we obtain, following Jones:

\index{Jones polynomial}
\index{Alexander polynomial}
\index{skein relations}
\index{Temperley-Lieb trace}
\index{Kauffmann polynomial}

\begin{theorem}
We can define the Jones polynomial of an oriented knot or link as being  the image of the corresponding braid producing it via the map
$$tr:B_k\to TL_N(k)\to\mathbb C$$
with the following change of variables: 
$$N=q^{1/2}+q^{-1/2}$$
We obtain a Laurent polynomial in $q^{1/2}$, which is an invariant, up to planar isotopy.
\end{theorem}

\begin{proof}
There is a long story here, the idea being as follows:

\medskip

(1) To start with, the result follows indeed by combining the above ingredients, the idea being that the various algebraic properties of $tr:TL_N(k)\to\mathbb C$ are exactly what is needed for the above composition, up to a normalization, to be invariant under the Reidemeister moves of type I, II, III, and so to produce indeed a knot invariant. 

\medskip

(2) More specifically, the result follows from Theorem 8.12, combined with what we have in Theorem 8.21, which is now fully proved, with the positivity part coming from Theorem 8.30, and with the change of variables $N=q^{1/2}+q^{-1/2}$ in the statement coming from the equation $t^2-Nt+1=0$ that we found in the proof of Theorem 8.21.

\medskip

(3) As an illustration for how this works, consider first the unknot:
$$
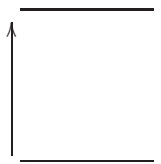
\xymatrix@R=50pt@C=50pt{
\ar@{-}[r]&\ar@{-}[d]\\
\ar[u]&\ar@{-}[l]}$$

For this knot, or rather unknot, the corresponding Jones polynomial is:
$$V=1$$

(4) Next, let us look at the link formed by two unlinked unknots:
$$
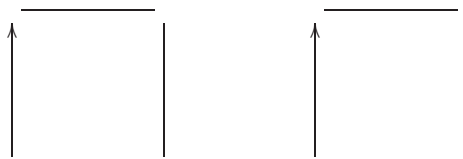
\xymatrix@R=50pt@C=50pt{
\ar@{-}[r]&\ar@{-}[d]&\ar@{-}[r]&\ar@{-}[d]\\
\ar[u]&\ar@{-}[l]&\ar[u]&\ar@{-}[l]
}$$

For this link, or rather unlink, the corresponding Jones polynomial is:
$$V=-q^{-1/2}-q^{1/2}$$

(5) Next, let us look at the link formed by two linked unknots, namely:
$$
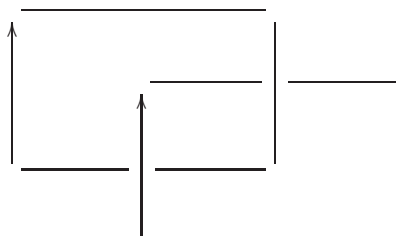
\xymatrix@R=14pt@C=40pt{
\ar@{-}[rr]&&\ar@{-}[dd]\\
&\ar@{-}[r]&|\ar@{-}[r]&\\
\ar@{-}[r]\ar[uu]&|\ar@{-}[r]&&\\
&\ar@{-}[rr]\ar[uu]&&\ar@{-}[uu]
}$$

For this link, the corresponding Jones polynomial is given by:
$$V=q^{1/2}+q^{5/2}$$

(6) Finally, let us look at the trefoil knot, which is as follows:
$$
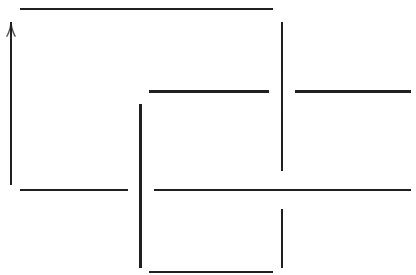
\xymatrix@R=18pt@C=40pt{
\ar@{-}[rr]&&\ar@{-}[dd]\\
&\ar@{-}[r]&|\ar@{-}[r]&\\
\ar@{-}[r]\ar[uu]&|\ar@{-}[rr]&-&\ar@{-}[u]\\
&\ar@{-}[r]\ar@{-}[uu]&\ar@{-}[u]
}$$

For this knot, the corresponding Jones polynomial is as follows:
$$V=q+q^3-q^4$$

Observe that, as previously for the unknot, this is a Laurent polynomial in $q$. This is part of a more general phenomenon, the point being that for knots, or more generally for links having an odd number of components, we get a Laurent polynomial in $q$.

\medskip

(7) In practice now, far more things can be said, about this. For instance the change of variables $N=q^{1/2}+q^{-1/2}$ in the statement, that we already used in chapter 3, in the random walk context for the small norm graphs, is something well-known in planar algebras, and with all this being related to operator algebras and subfactor theory. More on this in chapter 16 below, when discussing subfactors and planar algebras.

\medskip

(8) From a purely topological perspective, however, nothing beats the skein relation interpretation of the Jones polynomial $V_L(q)$, which is as follows, with $L_+,L_-,L_0$ being knots, or rather links, differing at exactly 1 crossing, in the 3 possible ways:
$$q^{-1}V_{L_+}-qV_{L_-}=(q^{1/2}+q^{-1/2})V_{L_0}$$

To be more precise, here are the conventions for $L_+,L_-,L_0$, that you need to know, in order to play with the above formula, and compute Jones polynomials at wish:
$$
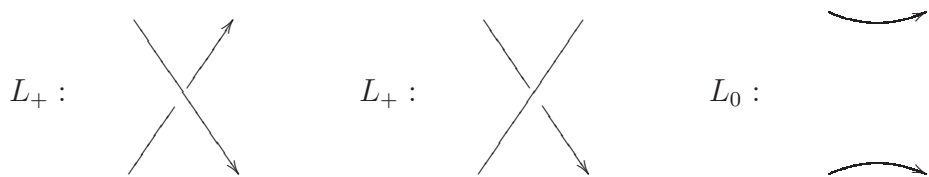
\xymatrix@R=20pt@C=15pt{
&\ar[ddrr]&&&&&\ar@{-}[dr]&&&&&\ar@/_/[rr]&&\\
L_+:&&\ar[ur]&&&L_+:&&\ar[dr]&&&L_0:\\
&\ar@{-}[ur]&&&&&\ar@{-}[uurr]&&&&&\ar@/^/[rr]&&}$$

As for the proof of the above formula, this comes from our definition of the Jones polynomial, because thinking well, ``unclosing'' links as to get braids, and then closing Temperley-Lieb diagrams as to get scalars, as required by the construction of $V_L(q)$, seemingly is some sort of identity operation, but the whole point comes from the fact that the Artin braids $g_1,\ldots,g_{k-1}$ and the Jones projections $e_1,\ldots,e_{k-1}$ differ precisely by a crossing being replaced by a non-crossing. Exercise for you, to figure out all this.

\medskip

(9) In short, up to you to learn all this, in detail, and its generalizations too, with link polynomials defined more generally via relations of the following type:
$$xP_{L_+}+yP_{L_-}+zP_{L_0}=0$$

Equivalently, we can define these more general invariants by using various versions of the Temperley-Lieb algebra. As usual, check here the papers of Jones \cite{jo1}, \cite{jo2}, \cite{jo3}.

\medskip

(10) With the comment here that, among all these invariants, Jones polynomial included, the first came, historially, the Alexander polynomial. However, from a modern point of view, the Alexander polynomial is something more complicated than the Jones polynomial, which remains the central invariant of knots and links.

\medskip 

(11) As another comment, with all this pure mathematics digested, physics strikes back, via a very interesting relation with statistical mechanics, happening in 2D as well, the idea being that ``interactions happen at crossings'', and it is these interactions that produce the knot invariant, as a kind of partition function. See Jones \cite{jo4}, \cite{jo5}.

\medskip

(12) Quite remarkably, the above invariants can be directly understood in 3D as well, in a purely geometric way, with elegance, and no need for 2D projection. But this is a more complicated story, involving ideas from quantum field theory. See Witten \cite{wit}.
\end{proof}

\section*{8e. Exercises} 

As usual, quite exciting chapter that we had here, but as a downside, many exercises left for you, in relation with all the above:

\begin{exercise}
Try to enumerate the links made of linked unknots.
\end{exercise}

\begin{exercise}
Clarify all details in relation with the Reidemeister moves.
\end{exercise}

\begin{exercise}
Clarify all the details in relation with the closure of braids.
\end{exercise}

\begin{exercise}
Learn more about the braid group, and its representations.
\end{exercise}

\begin{exercise}
Learn more about the Temperley-Lieb algebra, and its applications.
\end{exercise}

\begin{exercise}
Have some fun with meanders, and the Di Francesco formula.
\end{exercise}

\begin{exercise}
Clarify all details for the construction of the Jones polynomial.
\end{exercise}

\begin{exercise}
Compute the Jones polynomial, for knots and links of your choice.
\end{exercise}

As bonus exercise, learn some systematic topology, and more specifically some knot theory, from the papers of Jones \cite{jo2}, \cite{jo3}, \cite{jo4}, and from the paper of Witten \cite{wit}.

\part{Symmetry groups}

\ \vskip50mm

\begin{center}
{\em Here's my key

Philosophy

A freak like me

Just needs infinity}
\end{center}

\chapter{Symmetry groups}

\section*{9a. Finite groups}

We have seen that the study of the graphs $X$ is usually best done via the study of their adjacency matrices $d\in M_N(0,1)$, and with this latter study being typically a mix of linear algebra, combinatorics, calculus, probability, and some computer programming too. We will be back of course to this general principle, on numerous occasions.

\bigskip

However, as another interesting thing that we discovered, the group actions $G\curvearrowright X$ are something important too, making various group theoretical tools, coming from the study of the permutation groups $G\subset S_N$, useful for our graph questions. Remember for instance how many things we could say in the case where we have a transitive action $G\curvearrowright X$ of a cyclic or abelian group $G$, by using discrete Fourier analysis tools.

\bigskip

In this chapter, and in this whole Part III of this book, and in fact in the remainder of the whole book, we will take such things very seriously, and systematically develop the study of the group actions $G\curvearrowright X$, our aim being to solve the following question:

\begin{question}
Given a graph $X$, which groups $G$ act on it, $G\curvearrowright X$, and what can be said about $X$ itself, coming from the group theory of $G$?
\end{question} 

As already mentioned, we have a whole 100 pages, and later even more, for dealing with this question. So, we will go very slowly, at least in the beginning.

\bigskip

Let us start with some generalities. We already met groups $G$ and actions $G\curvearrowright X$ in Part I, and notably in chapter 4, which was dedicated to the transitive graphs. So, we are certainly not new to this subject. However, since there is no hurry with anything, let us first have a crash course in group theory. The starting point is of course:

\index{group}

\begin{definition}
A group is a set $G$ endowed with a multiplication operation 
$$(g,h)\to gh$$
which must satisfy the following conditions:
\begin{enumerate}
\item Associativity: we have $(gh)k=g(hk)$, for any $g,h,k\in G$.

\item Unit: there is an element $1\in G$ such that $g1=1g=g$, for any $g\in G$.

\item Inverses: for any $g\in G$ there is $g^{-1}\in G$ such that $gg^{-1}=g^{-1}g=1$.
\end{enumerate}
\end{definition}

As a first comment, the multiplication law is not necessarily commutative. In the case where it is, $gh=hg$ for any $g,h\in G$, we call $G$ abelian, and we usually denote its multiplication, unit and inverse operation in an additive way, as follows:
$$(g,h)\to g+h$$
$$0\in G$$
$$g\to-g$$

However, this is not a general rule, and rather the converse is true, in the sense that if a group is denoted as above, this means that the group must be abelian.

\bigskip

Let us first look at the abelian groups. Here as basic examples we have various groups of numbers, such as $\mathbb Z,\mathbb Q,\mathbb R,\mathbb C$, with the addition operation $+$. Observe that we have as well as $\mathbb Q^*,\mathbb R^*,\mathbb C^*$, and the unit circle $\mathbb T$ too, with the multiplication operation $\times$.

\bigskip

In order to reach to some theory, let us look into the finite group case, $|G|<\infty$. Here as basic examples we have the cyclic groups, constructed as follows:

\begin{proposition}
The following constructions produce the same group, denoted $\mathbb Z_N$, which is finite and abelian, and is called cyclic group of order $N$:
\begin{enumerate}
\item $\mathbb Z_N$ is the set of remainders modulo $N$, with operation $+$. 

\item $\mathbb Z_N\subset\mathbb T$ is the group of $N$-th roots of unity, with operation $\times$. 
\end{enumerate}
\end{proposition}

\begin{proof}
Here the equivalence between (1) and (2) is obvious. More complicated, however, is the question of finding the good philosophy and notation for this group. In what concerns us, we will be rather geometers, as usual, and we will often prefer the interpretation (2). As for the notation, we will use $\mathbb Z_N$, which is very natural.
\end{proof}

As a comment here, some algebraists prefer to reserve the notation $\mathbb Z_N$ for the ring of $p$-adic integers, when $N=p$ is prime, and use other notations for the cyclic group. Personally, although having a long love story with number theory, back in the days, I consider that the cyclic group is far more important, and gets the right to be denoted $\mathbb Z_N$. As for the $p$-adic integers, a reasonable notation for them, which does the job, is $\mathbf Z_N$.

\bigskip

As a basic thing to be known, about the abelian groups, still in the finite case, we can construct further examples of such groups by making products between various cyclic groups $\mathbb Z_N$. Quite remarkably, we obtain in this way all the finite abelian groups:

\begin{theorem}
The finite abelian groups are precisely the products of cyclic groups: 
$$G=\mathbb Z_{N_1}\times\ldots\times\mathbb Z_{N_k}$$
Moreover, there are technical extensions of this result, going beyond the finite case.
\end{theorem}

\begin{proof}
This is something that we already discussed in chapter 4, with the main ideas of the proof, and for further details on this, and for the technical extensions to the infinite groups as well, we recommend a solid algebra book, such as Lang \cite{lan}.
\end{proof}

Moving forward, let us look now into the general, non-abelian case. The first thought goes here to the $N\times N$ matrices with their multiplication, but these do not form a group, because we must assume $\det A\neq0$ in order for our matrix to be invertible.

\bigskip

So, let us call $GL_N(\mathbb C)$ the group formed by these latter matrices, with nonzero determinant, with $GL$ standing for ``general linear''. By further imposing the condition $\det A=1$ we obtain a subgroup $SL_N(\mathbb C)$, with $SL$ standing for ``special linear'', and then we can talk as well about the real versions of these groups, and also intersect everything with the group of unitary matrices $U_N$. We obtain in this way 8 groups, as follows:

\begin{theorem}
We have groups of invertible matrices as follows,
$$\xymatrix@R=22pt@C=5pt{
&GL_N(\mathbb R)\ar[rr]&&GL_N(\mathbb C)\\
O_N\ar[rr]\ar[ur]&&U_N\ar[ur]\\
&SL_N(\mathbb R)\ar[rr]\ar[uu]&&SL_N(\mathbb C)\ar[uu]\\
SO_N\ar[uu]\ar[ur]\ar[rr]&&SU_N\ar[uu]\ar[ur]
}$$
with $S$ standing for ``special'', meaning having determinant $1$.
\end{theorem}

\begin{proof}
This is clear indeed from the above discussion. As a comment, we can talk in fact about $GL_N(F)$ and $SL_N(F)$, once we have a ground field $F$, but in what regards the corresponding orthogonal and unitary groups, things here are more complicated.
\end{proof}

There are many other groups of matrices, besides the above ones, as for instance the symplectic groups $Sp_N\subset U_N$, appearing at $N\in2\mathbb N$. Generally speaking, the theory of Lie groups and algebras is in charge with the classification of such beasts.

\bigskip

Finally, and getting now to the point, let us discuss the general finite groups. As basic example here you have the symmetric group $S_N$, whose main properties are as follows:

\index{symmetric group}
\index{permutation group}
\index{signature}

\begin{theorem}
The permutations of $\{1,\ldots,N\}$ form a group, denoted $S_N$, and called symmetric group. This group has $N!$ elements. The signature map
$$\varepsilon:S_N\to\mathbb Z_2$$
can be regarded as being a group morphism, with values in $\mathbb Z_2=\{\pm1\}$.
\end{theorem}

\begin{proof}
The group property is indeed clear, and the count is clear as well, by recurrence on $N\in\mathbb N$. As for the last assertion, recall here the following formula for the signature of the permutations, which is something elementary to establish:
$$\varepsilon(\sigma\tau)=\varepsilon(\sigma)\varepsilon(\tau)$$

But this tells us precisely that $\varepsilon$ is a group morphism, as stated. 
\end{proof}

As graph theorists, we will be particularly interested in the subgroups $G\subset S_N$. We have already met a few such subgroups, as follows:

\begin{proposition}
We have finite non-abelian groups, as follows:
\begin{enumerate}
\item $S_N$, the group of permutations of $\{1,\ldots,N\}$.

\item $A_N\subset S_N$, the permutations having signature $1$.

\item $D_N\subset S_N$, the group of symmetries of the regular $N$-gon.
\end{enumerate}
\end{proposition}

\begin{proof}
The fact that we have indeed groups is clear from definitions, and the non-abelianity of these groups is clear as well, provided of course that in each case $N$ is chosen big enough, and with exercise for you to work out all this, with full details.
\end{proof}

For constructing further examples of finite non-abelian groups, the best is to ``look up'', by regarding $S_N$ as being the permutation group of the $N$ coordinate axes of $\mathbb R^N$. Indeed, this suggests looking at the symmetry groups of other geometric beasts inside $\mathbb R^N$, or even $\mathbb C^N$, and we end up with a whole menagery of groups, as follows: 

\begin{theorem}
We have groups of unitary matrices as follows,
$$
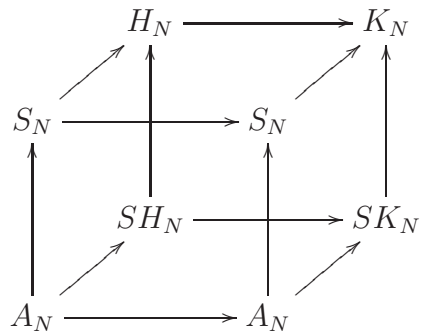
\xymatrix@R=20pt@C=16pt{
&H_N\ar[rr]&&K_N\\
S_N\ar[rr]\ar[ur]&&S_N\ar[ur]\\
&SH_N\ar[rr]\ar[uu]&&SK_N\ar[uu]\\
A_N\ar[uu]\ar[ur]\ar[rr]&&A_N\ar[uu]\ar[ur]
}$$
for the most finite, and non-abelian, called complex reflection groups.
\end{theorem}

\begin{proof}
This statement is of course something informal, and here are explanations on all this, including definitions for all the groups involved:

\medskip

(1) To start with, $S_N$ is the symmetric group $S_N$ that we know, but regarded now as permutation group of the $N$ coordinate axes of $\mathbb R^N$, and so as subgroup $S_N\subset O_N$.

\medskip

(2) Similarly, $A_N$ is the alternating group $A_N$ that we know, but coming now geometrically, as $A_N=S_N\cap SO_N$, with the intersection being computed inside $O_N$. 

\medskip

(3) Regarding $H_N\subset O_N$, this is a famous group, called hyperoctahedral group, appearing as the symmetry group of the hypercube $\square_N\subset\mathbb R^N$.

\medskip

(4) Regarding $K_N\subset U_N$, this is the complex analogue of $H_N$, consisting of the unitary matrices $U\in U_N$ having exactly one nonzero entry, on each row and each column.

\medskip

(5) We have as well on our diagram the groups $SH_N,SK_N$, with $S$ standing as usual for ``special'', that is, consisting of the matrices in  $H_N,K_N$ having determinant 1.

\medskip

(6) In what regards now the diagram itself, sure I can see that $S_N,A_N$ appear twice, but nothing can be done here, after thinking a bit, at how the diagram works.

\medskip

(7) Let us mention too that the groups $\mathbb Z_N,D_N$ have their place here, in $N$-dimensional geometry, but not exactly on our diagram, as being the symmetry groups of the oriented cycle, and unoriented cycle, with vertices at the simplex $X_N=\{e_i\}\subset\mathbb R^N$.

\medskip

(8) Finally, in what regards the finiteness, non-abelianity, and also the name ``complex reflection groups'', many things to be checked here, left to you as an exercise.
\end{proof}

Very nice all this. Let us summarize this group theory discussion as follows:

\begin{conclusion}
All groups, or almost, are best seen as being groups of matrices. And even as groups of unitary matrices, in most cases.
\end{conclusion}

Observe in particular that this justifies our choice in Definition 9.2, for the group operation to be denoted multiplicatively, $\times$. Indeed, in most cases, in view of the above general principle, that abstract multiplication is in fact a matrix multiplication.

\bigskip

Thinking a bit at the above principle, that is so useful in practice, taking us away from the abstraction of Definition 9.2, that we would like to have it as a theorem:

\begin{question}
Can we make our general ``most groups are in fact groups of matrices'' principle, into an abstract theorem?
\end{question}

In general, this is a quite subtle question, related to advanced group theory, such as Lie groups and algebras, and representation theory for them, and we will have a taste of such things in chapter 11 below, when developing the Peter-Weyl theory for the finite groups, that we will need at that time, in connection with our regular graph business. 

\bigskip

The point however is that, in the finite group case that we are interested in, the answer to Question 9.10 is indeed ``yes'', thanks to two well-known theorems, which are both elementary. First we have the Cayley embedding theorem, which is as follows:

\index{Cayley embedding}

\begin{theorem}
Given a finite group $G$, we have an embedding as follows,
$$G\subset S_N\quad,\quad g\to(h\to gh)$$
with $N=|G|$. Thus, any finite group is a permutation group. 
\end{theorem}

\begin{proof}
Given a group element $g\in G$, we can associate to it the following map:
$$\sigma_g:G\to G\quad,\quad 
h\to gh$$

Since $gh=gh'$ implies $h=h'$, this map is bijective, and so is a permutation of $G$, viewed as a set. Thus, with $N=|G|$, we can view this map as a usual permutation, $\sigma_G\in S_N$. Summarizing, we have constructed a map as follows:
$$G\to S_N\quad,\quad 
g\to\sigma_g$$

Our first claim is that this is a group morphism. Indeed, this follows from:
$$\sigma_g\sigma_h(k)
=\sigma_g(hk)
=ghk
=\sigma_{gh}(k)$$

It remains to prove that this group morphism is injective. But this follows from:
\begin{eqnarray*}
g\neq h
&\implies&\sigma_g(1)\neq\sigma_h(1)\\
&\implies&\sigma_g\neq\sigma_h
\end{eqnarray*}

Thus, we are led to the conclusion in the statement.
\end{proof}

As a comment, the above Cayley embedding theorem, while being certainly very beautiful at the theoretical level, has two weaknesses. First is the fact that the embedding $G\subset S_N$ that we constructed depends on a particular writing $G=\{g_1,\ldots,g_N\}$, which is needed in order to identify the permutations of $G$ with the elements of the symmetric group $S_N$. Which in practice, is of course something which is not very good.

\bigskip

As a second point of criticism, the Cayley theorem often provides us with a ``wrong embedding'' of our group. Indeed, as an illustration here, for the basic examples of groups that we know, the Cayley theorem provides us with embeddings as follows:
$$\mathbb Z_N\subset S_N\quad,\quad 
D_N\subset S_{2N}\quad,\quad 
S_N\subset S_{N!}\quad,\quad 
H_N\subset S_{2^NN!}$$

And, compare this with the ``good embeddings'' of these groups, which are:
$$\mathbb Z_N\subset S_N\quad,\quad 
D_N\subset S_N\quad,\quad 
S_N\subset S_N\quad,\quad 
H_N\subset S_{2N}$$

Thus, the first Cayley embedding is the good one, the second one is not the best possible one, but can be useful, and the third and fourth embeddings are obviously useless. So, as a conclusion, the Cayley theorem remains something quite theoretical.

\bigskip

Nevermind. We will fix this, once we will know more. Going ahead now, as previously mentioned, the Cayley theorem is just half of the story, the other half being:

\index{permutation matrix}

\begin{theorem}
We have a group embedding as follows, obtained by regarding $S_N$ as permutation group of the $N$ coordinate axes of $\mathbb R^N$,
$$S_N\subset O_N$$
which makes a permutation $\sigma\in S_N$ correspond to the matrix having $1$ on row $i$ and column $\sigma(i)$, for any $i$, and having $0$ entries elsewhere.
\end{theorem}

\begin{proof}
The first assertion is clear, because the permutations of the $N$ coordinate axes of $\mathbb R^N$ are isometries. Regarding now the explicit formula, we have by definition:
$$\sigma(e_j)=e_{\sigma(j)}$$

Thus, the permutation matrix corresponding to $\sigma$ is given by:
$$\sigma_{ij}=
\begin{cases}
1&{\rm if}\ \sigma(j)=i\\
0&{\rm otherwise}
\end{cases}$$

Thus, we are led to the formula in the statement.
\end{proof}

We can combine the above result with the Cayley theorem, and we obtain the following result, which is something very nice, of obvious theoretical importance:

\index{finite group}
\index{Cayley embedding}
\index{permutation group}

\begin{theorem}
Given a finite group $G$, we have an embedding as follows,
$$G\subset O_N\quad,\quad g\to(e_h\to e_{gh})$$
with $N=|G|$. Thus, any finite group is an orthogonal matrix group.
\end{theorem}

\begin{proof}
The Cayley theorem gives an embedding as follows:
$$G\subset S_N\quad,\quad g\to(h\to gh)$$

On the other hand, Theorem 9.12 provides us with an embedding as follows:
$$S_N\subset O_N\quad,\quad 
\sigma\to(e_i\to e_{\sigma(i)})$$

Thus, we are led to the conclusion in the statement.
\end{proof}

The same remarks as for the Cayley theorem apply. First, the embedding $G\subset O_N$ that we constructed depends on a particular writing $G=\{g_1,\ldots,g_N\}$. And also, for the basic examples of groups that we know, the embeddings that we obtain are as follows:
$$\mathbb Z_N\subset O_N\quad,\quad 
D_N\subset O_{2N}\quad,\quad 
S_N\subset O_{N!}\quad,\quad 
H_N\subset O_{2^NN!}$$

And, compare this with the ``good embeddings'' of these groups, which are:
$$\mathbb Z_N\subset O_N\quad,\quad 
D_N\subset O_N\quad,\quad 
S_N\subset O_N\quad,\quad 
H_N\subset O_N$$

As before with the abstract group embeddings $G\subset S_N$ coming from Cayley, the first abstract embedding is the good one, the second one is not the best possible one, but can be useful, and the third and fourth abstract embeddings are obviously useless. 

\bigskip

The problem is now, how to fix this, as to have a theorem providing us with good embeddings? And a bit of thinking at all the above leads to the following conclusion:

\begin{conclusion}
The weak point in all the above is Cayley.
\end{conclusion}

So, leaving aside now Cayley, or rather putting it in our pocket, matter of having it there, handy when meeting abstract questions, let us focus on Theorem 9.12, which is the only reasonable theorem that we have. In relation with the basic groups, we have:

\begin{theorem}
We have the following finite groups of matrices:
\begin{enumerate}
\item $\mathbb Z_N\subset O_N$, the cyclic permutation matrices.

\item $D_N\subset O_N$, the dihedral permutation matrices.

\item $S_N\subset O_N$, the permutation matrices.

\item $H_N\subset O_N$, the signed permutation matrices.
\end{enumerate}
\end{theorem}

\begin{proof}
This is something self-explanatory, the idea being that Theorem 9.12 provides us with embeddings as follows, given by the permutation matrices:
$$\mathbb Z_N\subset D_N\subset S_N\subset O_N$$

In addition, looking back at the definition of $H_N$, this group inserts into the embedding on the right, $S_N\subset H_N\subset O_N$. Thus, we are led to the conclusion that all our 4 groups appear as groups of suitable ``permutation type matrices''. To be more precise:

\medskip

(1) The cyclic permutation matrices are by definition the matrices as follows, with 0 entries elsewhere, and form a group, which is isomorphic to the cyclic group $\mathbb Z_N$:
$$U=\begin{pmatrix}
&&&1&&&\\
&&&&&\ddots&\\
&&&&&&1\\
1&&&&&\\
&\ddots&&&&&\\
&&1&&&
\end{pmatrix}$$

(2) The dihedral matrices are the above cyclic permutation matrices, plus some suitable symmetry permutation matrices, and form a group which is isomorphic to $D_N$.

\medskip

(3) The permutation matrices, which by Theorem 9.12 form a group which is isomorphic to $S_N$, are the $0-1$ matrices having exactly one 1 on each row and column.

\medskip

(4) Finally, regarding the signed permutation matrices, these are by definition the $(-1)-0-1$ matrices having exactly one nonzero entry on each row and column, and by Theorem 9.8 these matrices form a group, which is isomorphic to $H_N$.
\end{proof}

The above groups are all groups of orthogonal matrices. When looking into general unitary matrices, we led to the following interesting class of groups:

\index{complex reflection group}

\begin{definition}
The complex reflection group $H_N^s\subset U_N$, depending on parameters
$$N\in\mathbb N\quad,\quad s\in\mathbb N\cup\{\infty\}$$
is the group of permutation-type matrices with $s$-th roots of unity as entries,
$$H_N^s=M_N(\mathbb Z_s\cup\{0\})\cap U_N$$
with the convention $\mathbb Z_\infty=\mathbb T$, at $s=\infty$.
\end{definition}

Observe that at $s<\infty$, the above groups are finite. Also, at $s=1,2,\infty$ we obtain the following groups, that we already met in the above:
$$H_N^1=S_N\quad,\quad 
H_N^2=H_N\quad,\quad
H_N^\infty=K_N$$

We will be back later in this chapter with more details about $H_N^s$. However, at the philosophical level, we have extended our basic series of finite groups, as follows:
$$\mathbb Z_N\subset D_N\subset S_N\subset H_N\subset H_N^4\subset H_N^8\subset\ldots\ldots\subset U_N$$

But all this looks a bit complicated, so for being even more philosophers, let restrict the attention to the cases $s=1,2,\infty$, with $H_N^\infty=K_N$ having already been adopted, as a kind of ``almost'' finite group. With this convention, the conclusion is that we have extended our series of basic finite groups, regarded as unitary groups, as follows:
$$\mathbb Z_N\subset D_N\subset S_N\subset H_N\subset K_N\subset U_N$$

And good news, this will be our final say on the subject. Time now to formulate a grand conclusion to what we did so far in this chapter, as follows:

\begin{grandconclusion}
Group theory can be understood as follows:
\begin{enumerate}
\item Most interesting groups are groups of matrices, $G\subset GL_N(\mathbb C)$.

\item Quite often, these matrices can be chosen to be unitaries, $G\subset U_N$.

\item The Cayley theorem tells us that $G\subset S_N\subset O_N$, with $N=|G|$.

\item Most interesting finite groups appear as $G\subset U_N$, with $N<<|G|$.
\end{enumerate}
\end{grandconclusion}

Which is very nice, at least we know one thing, and with this type of talisman, we can now safely navigate through the abstractions of group theory. Of course, we will be back to this, once we will know more about groups, and regularly update our conclusions.

\section*{9b. Graph symmetries}

Moving now towards graph theory, that we are interested in, and which will in fact confirm and fine-tune the above conclusions, we have the following construction:

\begin{theorem}
Given a finite graph $X$, with vertices denoted $1,\ldots,N$, the symmetries of $X$, which are the permutations $\sigma\in S_N$ leaving invariant the edges, 
$$i-j\implies\sigma(i)-\sigma(j)$$
form a subgroup of the symmetric group, as follows, called symmetry group of $X$:
$$G(X)\subset S_N$$
As basic examples, for the empty graph, or for the simplex, we have $G(X)=S_N$.
\end{theorem}

\begin{proof}
Here the first assertion, regarding the group property of $G(X)$, is clear from definitions, because the symmetries of $X$ are stable under composition. The second assertion, regarding the empty graph and the simplex, is clear as well. So done, everything being trivial, and we have called this Theorem instead of Proposition because the construction $X\to G(X)$ will keep us busy, for the remainder of this book.
\end{proof}

Let us work out now some more examples. As a first result, dealing with the simplest graph ever, passed the empty graphs and the simplices, we have:

\begin{proposition}
The symmetry group of the regular $N$-gon
$$
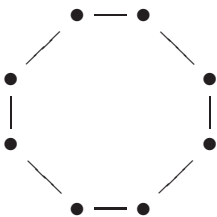
\xymatrix@R=12pt@C=12pt{
&\bullet\ar@{-}[r]\ar@{-}[dl]&\bullet\ar@{-}[dr]\\
\bullet\ar@{-}[d]&&&\bullet\ar@{-}[d]\\
\bullet\ar@{-}[dr]&&&\bullet\ar@{-}[dl]\\
&\bullet\ar@{-}[r]&\bullet}$$
is the dihedral group $D_N=\mathbb Z_N\rtimes\mathbb Z_2$.
\end{proposition}

\begin{proof}
This is something that we know well from chapter 4, and with the remark, which is something new, that the notation $D_N$ for the group that we get, which is the correct one, is justified by the general group theory discussion before, with $N$ standing for the natural ``dimensionality'' of this group. To be more precise, geometrically speaking, the regular $N$-gon is best viewed in $\mathbb R^N$, with vertices $1,\ldots,N$ at the standard basis:
$$1=(1,0,0,\ldots,0,0)$$ 
$$2=(0,1,0,\ldots,0,0)$$
$$\ \ \ \vdots$$
$$N=(0,0,0,\ldots,0,1)$$

But, with this interpretation in mind, we are led to an embedding as follows:
$$D_N\subset S_N\subset O_N$$

We conclude from this that $N$ is the correct dimensionality of our group, and so is the correct label to be attached to the dihedral symbol $D$. Of course, you might find this overly philosophical, or even a bit futile, but listen to this, there are two types of mathematicians in this world, those who use $D_N$ and those who use $D_{2N}$, and do not ask me why, but it is better to be in the first category, mathematicians using $D_N$.
\end{proof}

Moving ahead, the problem is now, is Proposition 9.19 good news, or bad news? I don't know about you, but personally I feel quite frustrated by the fact that the computation there leads to $D_N=\mathbb Z_N\rtimes\mathbb Z_2$, instead to $\mathbb Z_N$ itself. I mean, how can a theory be serious, if there is no room there, or even an Emperor's throne, for the cyclic group $\mathbb Z_N$.

\bigskip

So, let us fix this. It is obvious that the construction in Theorem 9.18 will work perfectly well for the oriented graphs, or for the colored graphs, so let us formulate:

\begin{definition}
Given a generalized graph $X$, with vertices denoted $1,\ldots,N$, the symmetries of $X$, which are the permutations $\sigma\in S_N$ leaving invariant the edges, 
$$i-j\implies\sigma(i)-\sigma(j)$$
with their orientations and colors, form a subgroup of the symmetric group
$$G(X)\subset S_N$$
called symmetry group of $X$.
\end{definition}

Here, as before with the construction in Theorem 9.18, the fact that we obtain indeed a group is clear from definitions. Now with this convention in hand, we have:

\begin{proposition}
The symmetry group of the oriented $N$-gon
$$
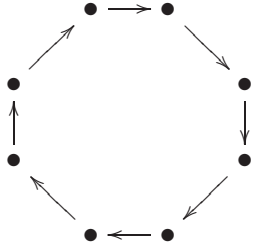
\xymatrix@R=16pt@C=16pt{
&\bullet\ar[r]&\bullet\ar[dr]\\
\bullet\ar[ur]&&&\bullet\ar[d]\\
\bullet\ar[u]&&&\bullet\ar[dl]\\
&\bullet\ar[ul]&\bullet\ar[l]}$$
is the cyclic group $\mathbb Z_N$.
\end{proposition}

\begin{proof}
This is clear from definitions, because once we choose a vertex $i$ and denote its image by $\sigma(i)=i+k$, the permutation $\sigma\in S_N$ leaving invariant the edges, with their orientation, must map $\sigma(i+1)=i+k+1$, $\sigma(i+2)=i+k+2$ and so on, and so must be an element of the cyclic group, in remainder modulo $N$ notation $\sigma=k\in\mathbb Z_N$.
\end{proof}

With this done, and the authority of $\mathbb Z_N$ restored, let us work out some general properties of the construction $X\to G(X)$. For simplicity we will restrict the attention to the usual graphs, as in Theorem 9.18, but pretty much everything will extend to the case of oriented or colored graphs. In fact, our policy in what follows will be that of saying nothing when things extend, and making a comment, when things do not extend.

\bigskip

As a first general result, coming as a useful complement to Theorem 9.18, we have:

\begin{theorem}
Having a group action on a graph $G\curvearrowright X$ is the same as saying that the action of $G$ leaves invariant the adjacency matrix $d$, in the sense that:
$$d_{ij}=d_{g(i)g(j)}\quad,\quad\forall g\in G$$
Equivalently, the action must preserve the spectral projections of $d$:
$$d=\sum_\lambda\lambda P_\lambda\implies (P_\lambda)_{ij}=(P_\lambda)_{g(i)g(j)}$$
Thus, the symmetry group $G(X)\subset S_N$ is the subgroup preserving the eigenspaces of $d$.
\end{theorem}

\begin{proof}
As before with Theorem 9.18, a lot of talking in the statement, with everything being trivial, coming from definitions, and with the statement itself being called Theorem instead of Proposition just due to its theoretical importance.
\end{proof}

Observe that Theorem 9.22 naturally leads us into colored graphs, because while the adjacency matrix is symmetric and binary, $d\in M_N(0,1)^{symm}$, the spectral projections $P_\lambda$ are also symmetric, but no longer binary, $P_\lambda\in M_N(\mathbb R)^{symm}$. Moreover, these spectral projections $P_\lambda$ can have 0 on the diagonal, pushing us into allowing self-edges in our colored graph formalism. We are led in this way to the following statement: 

\begin{theorem}
Having a group action on a colored graph $G\curvearrowright X$ is the same as saying that the action of $G$ leaves invariant the adjacency matrix $d$:
$$d_{ij}=d_{g(i)g(j)}\quad,\quad\forall g\in G$$
Equivalently, the action must preserve the spectral projections of $d$, as follows:
$$d=\sum_\lambda\lambda P_\lambda\implies (P_\lambda)_{ij}=(P_\lambda)_{g(i)g(j)}$$
Moreover, when allowing self-edges, each $P_\lambda$ will correspond to a colored graph $X_\lambda$.
\end{theorem}

\begin{proof}
This follows indeed from the above discussion, and with some extra discussion regarding the precise colors that we use, as follows:

\medskip

(1) When using real colors, the result follows from the linear algebra result regarding the diagonalization of real symmetric matrices, which tells us that the spectral projections of any such matrix $d\in M_N(\mathbb R)^{symm}$ are also real and symmetric, $P_\lambda\in M_N(\mathbb R)^{symm}$. 

\medskip

(2) When using complex colors, the result follows from the linear algebra result regarding the diagonalization of complex self-adjoint matrices, which tells us that the spectral projections of any such matrix $d\in M_N(\mathbb C)^{sa}$ are also self-adjoint, $P_\lambda\in M_N(\mathbb C)^{sa}$.
\end{proof}

The point with the perspective brought by the above results is that, when using permutation group tools for the study of the groups $G\subset S_N$ acting on our graph, $G\curvearrowright X$, what will eventually happen is that these tools, once sufficiently advanced, will become very close to the regular tools for the study of $d$, namely the same sort of mixture of linear algebra, calculus and probability, so in the end we will have a unified theory.

\bigskip

But probably too much talking, just trust me, we won't be doing groups and algebra here just because we are scared by analysis, and by the true graph problems. Quite the opposite. And we will see illustrations for this harmony and unity later on.

\bigskip

Leaving now the oriented or colored graphs aside, as per our general graph policy explained above, as a second general result about $X\to G(X)$, we have:

\index{complementation}

\begin{theorem}
The construction $X\to G(X)$ has the property
$$G(X)=G(X^c)$$
where $X\to X^c$ is the complementation operation.
\end{theorem}

\begin{proof}
This is clear from the construction of $G(X)$ from Theorem 9.18, and follows as well from the interpretation in Theorem 9.22, because the adjacency matrices of $X$, $X^c$ are related by the following formula, where $\mathbb I_N$ is the all-one matrix:
$$d_X+d_{X^c}=\mathbb I_N-1_N$$

Indeed, since on the right we have the adjacency matrix of the simplex, which commutes with everything, commutation with $d_X$ is equivalent to commutation with $d_{X^c}$, and this gives the result, via the interpretation of $G(X)$ coming from Theorem 9.22.
\end{proof}

In order to reach now to more advanced results, it is convenient to enlarge the attention to the colored graphs. Indeed, for the colored graphs, we can formulate:

\begin{theorem}
Having an action on a colored graph $G\curvearrowright X$ is the same as saying that the action leaves invariant the color components of $X$. Equivalently, with
$$d=\sum_{c\in C}cd_c$$
being the color decomposition of the adjacency matrix, with color components
$$(d_c)_{ij}=\begin{cases}
1&{\rm if}\ d_{ij}=c\\
0&{\rm otherwise}
\end{cases}$$
the action must leave invariant all these color components $d_c$. Thus, the symmetry group $G(X)\subset S_N$ is the subgroup which preserves all these matrices $d_c$.
\end{theorem}

\begin{proof}
As before with our other statements here, in the present section of this book, a lot of talking in the statement, with everything there being trivial.
\end{proof}

I have this feeling that you might get to sleep, on the occasion of the present section, which is overly theoretical, this is how things are, we have to have some theory started, right. But, in the case it is so, I have something interesting for you, in relation with the above. Indeed, by combining Theorem 9.23 with Theorem 9.25, both trivialities, we are led to the following enigmatic statement, which all of the sudden wakes us up:

\begin{theorem}
Given an adjacency matrix of a graph $X$, which can be taken in a colored graph sense, $d\in M_N(\mathbb C)$, or even binary as usual, 
$$d\in M_N(0,1)$$
a group action $G\curvearrowright X$ must preserve all ``spectral-color'' components of this matrix, obtained by succesively applying the spectral decomposition, and color decomposition.
\end{theorem}

\begin{proof}
This is clear indeed by combining Theorem 9.23 and Theorem 9.25, and with the remark that, indeed, even for a usual binary matrix $d\in M_N(0,1)$ this leads to something non-trivial, because the spectral components of this matrix are no longer binary, and so all of the sudden, we are into colors and everything.
\end{proof}

With the above result in hand, which is something quite unexpected, we are led into a quite interesting linear algebra question, which is surely new for you, namely:

\begin{question}
What are the spectral-color components of a matrix $d\in M_N(\mathbb C)$, or even of a usual binary matrix $d\in M_N(0,1)$?
\end{question}

This question is something non-trivial, and we will be back to this on several occasions, and notably at the end of this book, when talking planar algebras in the sense of Jones \cite{jo6}, which provide the good framework for the study of such questions. 

\bigskip

To be more precise, coming a bit in advance, we will see there that computing the spectral-color components of a matrix $d\in M_N(\mathbb C)$ is the same as computing the planar algebra generated by $d$, viewed as 2-box inside the tensor planar algebra.

\section*{9c. Reflection groups} 

We already know a number of things about the circulant graphs from Part I, and we also know that these usually generalize to the case where we have a transitive action of an abelian group on $X$. Both the group theory and the linear algebra aspects here can be quite subtle, involving the classification of finite abelian groups, and generalized Fourier matrices, and that discussion from Part I can be briefly summarized as follows:

\begin{fact}
The transitive abelian group actions on graphs, $G\curvearrowright X$ with
$$G=\mathbb Z_{N_1}\times\ldots\times\mathbb Z_{N_s}$$
can be investigated via discrete Fourier analysis, over the group $G$.
\end{fact} 

This is of course something very compact, and we refer to Part I for details. In order to advance now on all this, we have the following result, which is standard in discrete Fourier analysis, extending what we previously knew in the circulant case:

\index{Fourier-diagonal}
\index{discrete Fourier transform}

\begin{theorem}
For a matrix $A\in M_N(\mathbb C)$, the following are equivalent,
\begin{enumerate}
\item $A$ is $G$-invariant, $A_{ij}=\xi_{j-i}$, for a certain vector $\xi\in\mathbb C^N$,

\item $A$ is Fourier-diagonal, $A=F_GQF_G^*$, for a certain diagonal matrix $Q$,
\end{enumerate}
and if so, $\xi=F_G^*q$, where $q\in\mathbb C^N$ is the vector formed by the diagonal entries of $Q$.
\end{theorem}

\begin{proof}
This is something that we know from chapter 4 in the cyclic case, $G=\mathbb Z_N$, and the proof in general is similar, by using matrix indices as follows:
$$i,j\in G$$

To be more precise, in order to get started, with our generalization, let us decompose our finite abelian group $G$ as a product of cyclic groups, as follows:
$$G=\mathbb Z_{N_1}\times\ldots\times\mathbb Z_{N_s}$$

The corresponding Fourier matrix decomposes then as well, as follows:
$$F_G=F_{N_1}\otimes\ldots\otimes F_{N_s}$$

Now if we set $w_i=e^{2\pi i/N_i}$, this means that we have the following formula:
$$(F_G)_{ij}=w_1^{i_1j_1}\ldots w_s^{i_sj_s}$$

We can now prove the equivalence in the statement, as follows:

\medskip

$(1)\implies(2)$ Assuming $A_{ij}=\xi_{j-i}$, the matrix $Q=F_G^*AF_G$ is diagonal, as shown by the following computation, with all indices being group elements:
\begin{eqnarray*}
Q_{ij}
&=&\sum_{kl}\overline{(F_G)}_{ki}A_{kl}(F_G)_{lj}\\
&=&\sum_{kl}w_1^{-k_1i_1}\ldots w_s^{-k_si_s}\cdot\xi_{l-k}\cdot w_1^{l_1j_1}\ldots w_s^{l_sj_s}\\
&=&\sum_{kl}w_1^{l_1j_1-k_1i_1}\ldots w_s^{l_sj_s-k_si_s}\xi_{l-k}\\
&=&\sum_{kr}w_1^{(k_1+r_1)j_1-k_1i_1}\ldots w_s^{(k_s+r_s)j_s-k_si_s}\xi_r\\
&=&\sum_rw_1^{r_1j_1}\ldots w_s^{r_sj_s}\xi_r\sum_kw_1^{k_1(j_1-i_1)}\ldots w_s^{k_s(j_s-i_s)}\\
&=&\sum_rw_1^{r_1j_1}\ldots w_s^{r_sj_s}\xi_r\cdot N_1\delta_{i_1j_1}\ldots N_s\delta_{i_sj_s}\\
&=&N\delta_{ij}\sum_r(F_G)_{jr}\xi_r
\end{eqnarray*}

$(2)\implies(1)$ Assuming $Q=diag(q_1,\ldots,q_N)$, the matrix $A=F_GQF_G^*$ is $G$-invariant, as shown by the following computation, again with all indices being group elements:
\begin{eqnarray*}
A_{ij}
&=&\sum_{kl}(F_G)_{ik}Q_{kk}\overline{(F_G)}_{kj}\\
&=&\sum_kw_1^{i_1k_1}\ldots w_s^{i_sk_s}\cdot q_k\cdot w_1^{-j_1k_1}\ldots w_s^{-j_sk_s}\\
&=&\sum_kw_1^{(i_1-j_1)k_1}\ldots w_s^{(i_s-j_s)k_s}q_k
\end{eqnarray*}

To be more precise, in this formula the last term depends only on $j-i$, and so shows that we have $A_{ij}=\xi_{j-i}$, with $\xi$ being the following vector:
\begin{eqnarray*}
\xi_i
&=&\sum_kw_1^{-i_1k_1}\ldots w_s^{-i_sk_s}q_k\\
&=&\sum_k(F_G^*)_{ik}q_k\\
&=&(F_G^*q)_i
\end{eqnarray*}

Thus, we are led to the conclusions in the statement.
\end{proof}

As another generalization of what we did in chapter 4, in relation now with the dihedral groups, and then with the hyperoctahedral groups, we can investigate the complex reflection groups $H_N^s$ that we introduced in the above. Let us recall indeed that $H_N^s\subset U_N$ is the group of permutation-type matrices with $s$-th roots of unity as entries:
$$H_N^s=M_N(\mathbb Z_s\cup\{0\})\cap U_N$$

We know that at $s=1,2$, this group equals $S_N,H_N$, that we both know well. In analogy with what we know at $s=1,2$, we first have the following result:

\begin{proposition}
The number of elements of $H_N^s$ with $s\in\mathbb N$ is:
$$|H_N^s|=s^NN!$$
At $s=\infty$, the group $K_N=H_N^\infty$ that we obtain is infinite.
\end{proposition}

\begin{proof}
This is indeed clear from our definition of $H_N^s$, as a matrix group as above, because there are $N!$ choices for a permutation-type matrix, and then $s^N$ choices for the corresponding $s$-roots of unity, which must decorate the $N$ nonzero entries.
\end{proof}

Once again in analogy with what we know at $s=1,2$, we have as well:

\index{wreath product}
\index{hyperoctahedral group}

\begin{theorem}
We have a wreath product decomposition as follows,
$$H_N^s=\mathbb Z_s\wr S_N$$
which means by definition that we have a crossed product decomposition
$$H_N^s=\mathbb Z_s^N\rtimes S_N$$
with the permutations $\sigma\in S_N$ acting via $\sigma(e_1,\ldots,e_k)=(e_{\sigma(1)},\ldots,e_{\sigma(k)})$.
\end{theorem}

\begin{proof}
As explained in the proof of Proposition 9.30, the elements of $H_N^s$ can be identified with the pairs $g=(e,\sigma)$ consisting of a permutation $\sigma\in S_N$, and a decorating vector $e\in\mathbb Z_s^N$, so that at the level of the cardinalities, we have:
$$|H_N|=|\mathbb Z_s^N\times S_N|$$

Now observe that the product formula for two such pairs $g=(e,\sigma)$ is as follows, with the permutations $\sigma\in S_N$ acting on the elements $f\in\mathbb Z_s^N$ as in the statement:
$$(e,\sigma)(f,\tau)=(ef^\sigma,\sigma\tau)$$

Thus, we are in the framework of crossed products, and we obtain $H_N^s=\mathbb Z_s^N\rtimes S_N$. But this can be written, by definition, as $H_N^s=\mathbb Z_s\wr S_N$, and we are done.
\end{proof}

In relation with graph symmetries, the above groups appear as follows:

\begin{theorem}
The complex reflection group $H_N^s$ appears as symmetry group,
$$H_N^s=G(NC_s)$$
with $NC_s$ consisting of $N$ disjoint copies of the oriented cycle $C_s$.
\end{theorem}

\begin{proof}
This is something elementary, the idea being as follows:

\medskip

(1) Consider first the oriented cycle $C_s$, which looks as follows:
$$
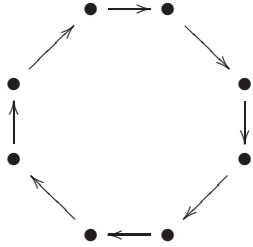
\xymatrix@R=16pt@C=16pt{
&\bullet\ar[r]&\bullet\ar[dr]\\
\bullet\ar[ur]&&&\bullet\ar[d]\\
\bullet\ar[u]&&&\bullet\ar[dl]\\
&\bullet\ar[ul]&\bullet\ar[l]}$$

It is then clear that the symmetry group of this graph is the cyclic group $\mathbb Z_s$.

\medskip

(2) In the general case now, where we have $N\in\mathbb N$ disjoint copies of the above cycle $C_s$, we must suitably combine the corresponding $N$ copies of the cyclic group $\mathbb Z_s$. But this leads to the wreath product group $H_N^s=\mathbb Z_s\wr S_N$, as stated.
\end{proof}

Normally, with the theory that we have so far, we can investigate the graphs having small number of vertices. But for more here, going beyond what we have, we need more product results, and we will develop the needed theory in the next chapter.

\section*{9d. Partial symmetries}

As a final topic for this theoretical chapter, let us discuss as well the partial symmetries of graphs. To be more precise, we will associate to any graph $X$ having $N$ vertices its semigroup of partial permutations $\widetilde{G}(X)\subset\widetilde{S}_N$, which is a quite interesting object.

\bigskip

In order to discuss all this, let us start with something well-known, namely:

\begin{definition}
$\widetilde{S}_N$ is the semigroup of partial permutations of $\{1\,\ldots,N\}$,
$$\widetilde{S}_N=\left\{\sigma:I\simeq J\Big|I,J\subset\{1,\ldots,N\}\right\}$$
with the usual composition operation for such partial permutations, namely
$$\sigma'\sigma:\sigma^{-1}(I'\cap J)\simeq\sigma'(I'\cap J)$$
being the composition of $\sigma':I'\simeq J'$ and $\sigma:I\simeq J$.
\end{definition}

Observe that $\widetilde{S}_N$ is not simplifiable, because the null permutation $\emptyset\in\widetilde{S}_N$, having the empty set as domain/range, satisfies the following formula, for any $\sigma\in\widetilde{S}_N$:
$$\emptyset\sigma=\sigma\emptyset=\emptyset$$

Observe also that our semigroup $\widetilde{S}_N$ has a kind of ``subinverse'' map, which is not a true inverse in the semigroup sense, sending a partial permutation $\sigma:I\to J$ to its usual inverse $\sigma^{-1}:J\to I$. Many other algebraic things can be said, along these lines.

\bigskip

As a first interesting result now about $\widetilde{S}_N$, which shows that we are dealing here with some non-trivial combinatorics, not really present in the $S_N$ context, we have:

\begin{theorem}
The number of partial permutations is given by
$$|\widetilde{S}_N|=\sum_{k=0}^Nk!\binom{N}{k}^2$$
that is, $1,2,7,34,209,\ldots\,$, and we have the formula
$$|\widetilde{S}_N|\simeq N!\sqrt{\frac{\exp(4\sqrt{N}-1)}{4\pi\sqrt{N}}}$$
in the $N\to\infty$ limit.
\end{theorem}

\begin{proof}
This is a mixture of trivial and non-trivial facts, as follows:

\medskip

(1) The first assertion is clear, because in order to construct a partial permutation $\sigma:I\to J$ we must choose an integer $k=|I|=|J|$, then we must pick two subsets $I,J\subset\{1,\ldots,N\}$ having cardinality $k$, and there are $\binom{N}{k}$ choices for each, and finally we must construct a bijection $\sigma:I\to J$, and there are $k!$ choices here. 

\medskip

(2) As for the estimate, which is non-trivial, this is something standard, and well-known, and exercise for you here to look that up, and read the proof.
\end{proof}

Another result, which is trivial, but is quite fundamental, is as follows:

\begin{proposition}
We have a semigroup embedding $u:\widetilde{S}_N\subset M_N(0,1)$, given by 
$$u_{ij}(\sigma)=
\begin{cases}
1&{\rm if}\ \sigma(j)=i\\
0&{\rm otherwise}
\end{cases}$$
whose image are the matrices having at most one nonzero entry, on each row and column.
\end{proposition}

\begin{proof}
This is something which is trivial from definitions, with the above embedding $u:\widetilde{S}_N\subset M_N(0,1)$ extending the standard embedding $u:S_N\subset M_N(0,1)$, given by the permutation matrices, that we have been heavily using, so far.
\end{proof}

Also at the algebraic level, we have the following simple and useful fact:

\begin{proposition}
Any partial permutation $\sigma:I\simeq J$ can be factorized as
$$\xymatrix@R=40pt@C=40pt
{I\ar[r]^{\sigma}\ar[d]_\gamma&J\\\{1,\ldots,k\}\ar[r]_\beta&\{1,\ldots,k\}\ar[u]_\alpha}$$
with $\alpha,\beta,\gamma\in S_k$ being certain non-unique permutations, where $k=|Dom(\sigma)|$.
\end{proposition}

\begin{proof}
We can choose indeed any two bijections $I\simeq\{1,\ldots,k\}$ and $\{1,\ldots,k\}\simeq J$, and then complete them up to permutations $\gamma,\alpha\in S_N$. The remaining permutation $\beta\in S_k$ is then uniquely determined by the formula $\sigma=\alpha\beta\gamma$.
\end{proof}

The above result is quite interesting, theoretically speaking, allowing us to reduce in principle questions about $\widetilde{S}_N$ to questions about the usual symmetric groups $S_k$, via some sort of homogeneous space procedure. We will be back to this, later in this book.

\bigskip

In relation now with the graphs, we have the following notion:

\begin{definition}
Associated to any graph $X$ is its partial symmetry semigroup
$$\widetilde{G}(X)\subset\widetilde{S}_N$$
consisting of the partial permutations $\sigma\in\widetilde{S}_N$ whose action preserves the edges.
\end{definition}

As a first observation, we have the following result, which provides an alternative to the above definition of the partial symmetry semigroup $\widetilde{G}(X)$:

\begin{proposition}
Given a graph $X$ with $N$ vertices, we have
$$\widetilde{G}(X)=\left\{\sigma\in\widetilde{S}_N\Big|d_{ij}=d_{\sigma(i)\sigma(j)},\ \forall i,j\in Dom(\sigma)\right\}$$
with $d\in M_N(0,1)$ being as usual the adjacency matrix of $X$.
\end{proposition}

\begin{proof}
The construction of $\widetilde{G}(X)$ from Definition 9.37 reformulates as follows, in terms of the usual adjacency relation $i-j$ for the vertices:
$$\widetilde{G}(X)=\left\{\sigma\in\widetilde{S}_N\Big|i-j,\exists\,\sigma(i),\exists\,\sigma(j)\implies \sigma(i)-\sigma(j)\right\}$$

But this leads to the formula in the statement, in terms of the adjacency matrix $d$.
\end{proof}

In order to discuss now some examples, let us make the following convention:

\begin{definition}
In the context of the partial permutations $\sigma:I\to J$, with $I,J\subset\{1,\ldots,N\}$, we decompose the domain set $I$ as a disjoint union
$$I=I_1\sqcup\ldots\sqcup I_p$$
with each $I_r$ being an interval consisting of consecutive numbers, and being maximal with this property, and with everything being taken cyclically.
\end{definition}

In other words, we represent the domain set $I\subset\{1,\ldots,N\}$ on a circle, with $1$ following $1,\ldots,N$, and then we decompose it into intervals, in the obvious way. With this convention made, in the case of the oriented cycle, we have the following result:

\begin{proposition}
For the oriented cycle $C_N^o$ we have
$$\widetilde{G}(C_N^o)=\widetilde{\mathbb Z}_N$$
with the semigroup on the right consisting of the partial permutations
$$\sigma:I_1\sqcup\ldots\sqcup I_p\to J$$
which are cyclic on any $I_r$, given there by $i\to i+k_r$, for a certain $k_r\in\{1,\ldots,N\}$.
\end{proposition}

\begin{proof}
According to the definition of $\widetilde{G}(X)$, we have the following formula:
$$\widetilde{G}(C_N^o)=\left\{\sigma\in\widetilde{S}_N\Big|d_{ij}=d_{\sigma(i)\sigma(j)},\ \forall i,j\in Dom(\sigma)\right\}$$

On the other hand, the adjacency matrix of $C_N^o$ is given by:
$$d_{ij}=\begin{cases}
1&{\rm if}\ j=i+1\\
0&{\rm otherwise}
\end{cases}$$

Thus, the condition defining $\widetilde{G}(C_N^o)$ is as follows:
$$j=i+1\iff\sigma(j)=\sigma(i)+1,\ \forall i,j\in Dom(\sigma)$$

But this leads to the conclusion in the statement.
\end{proof}

In the case of the unoriented cycle, the result is as follows:

\begin{proposition}
For the unoriented cycle $C_N$ we have
$$\widetilde{G}(C_N)=\widetilde{D}_N$$
with the semigroup on the right consisting of the partial permutations
$$\sigma:I_1\sqcup\ldots\sqcup I_p\to J$$
which are dihedral on any $I_r$, given there by $i\to \pm_ri+ k_r$, for a certain $k_r\in\{1,\ldots,N\}$, and a certain choice of the sign $\pm_r\in\{-1,1\}$.
\end{proposition}

\begin{proof}
The proof here is similar to the proof of Proposition 9.40. Indeed, the adjacency matrix of $C_N$ is given by:
$$d_{ij}=\begin{cases}
1&{\rm if}\ j=i\pm 1\\
0&{\rm otherwise}
\end{cases}$$

Thus, the condition defining $\widetilde{G}(C_N)$ is as follows:
$$j=i\pm 1\iff\sigma(j)=\sigma(i)\pm 1,\ \forall i,j\in Dom(\sigma)$$

But this leads to the conclusion in the statement.
\end{proof}

An interesting question is whether the semigroups $\widetilde{\mathbb Z}_N,\widetilde{D}_N$ are related by a formula similar to $D_N=\mathbb Z_N\rtimes\mathbb Z_2$. This is not exactly the case, at least with the obvious definition for the $\rtimes$ operation, because at the level of cardinalities we have:

\begin{theorem}
The cardinalities of $\widetilde{\mathbb Z}_N,\widetilde{D}_N$ are given by the formulae
$$|\widetilde{\mathbb Z}_N|=1+NK_1(N)+\sum_{p=2}^{[N/2]}N^pK_p(N)$$
$$|\widetilde{D}_N|=1+NK_1(N)+\sum_{p=2}^{[N/2]}(2N)^pK_p(N)$$
where $K_p(N)$ counts the sets having $p$ cyclic components, $I=I_1\sqcup\ldots\sqcup I_p$.
\end{theorem}

\begin{proof}
The first formula is clear from the description of $\widetilde{\mathbb Z}_N$ from Proposition 9.40, because for any domain set $I=I_1\sqcup\ldots\sqcup I_p$, we have $N$ choices for each scalar $k_r$, producing a cyclic partial permutation $i\to i+k_r$ on the interval $I_r$. Thus we have, as claimed:
$$|\widetilde{\mathbb Z}_N|=\sum_{p=0}^{[N/2]}N^pK_p(N)$$

In the case of $\widetilde{D}_N$ the situation is similar, with Proposition 9.41 telling us that the $N$ choices at the level of each interval $I_r$ must be now replaced by $2N$ choices, as to have a dihedral permutation $i\to \pm_ri+ k_r$ there. However, this is true only up to a subtlety, coming from the fact that at $p=1$ the choice of the $\pm1$ sign is irrelevant. Thus, we are led to the formula in the statement, with $2N$ factors everywhere, except at $p=1$.
\end{proof}

Summarizing, the partial symmetry group problematics leads to some interesting questions, even for simple graphs like $C_N^o$ and $C_N$. We will be back to this.

\section*{9e. Exercises}

Welcome to pure mathematics, we certainly had a pure mathematics chapter here, and of pure mathematics nature will be as well our exercises, as follows:

\begin{exercise}
Read more about Lie groups, and notably about $Sp_N$.
\end{exercise}

\begin{exercise}
Read more about finite groups, notably about Sylow theorems.
\end{exercise}

\begin{exercise}
And even more about finite groups, with complex reflection groups.
\end{exercise}

\begin{exercise}
Read about discrete groups too, and random walks on them.
\end{exercise}

\begin{exercise}
Meditate at the spectral-color decomposition, and its planar aspects.
\end{exercise}

\begin{exercise}
Meditate at our various graph formalisms, their pros and cons.
\end{exercise}

\begin{exercise}
Meditate at what else can be done, with discrete Fourier analysis.
\end{exercise}

\begin{exercise}
Find a proof for the asymptotic estimate for $|\widetilde{S}_N|$.
\end{exercise}

\begin{exercise}
Work out more examples of semigroups $\widetilde{G}(X)$.
\end{exercise}

As bonus exercise, learn some representation theory for finite groups. We will discuss this soon in this book, but the more you know a bit in advance, the better.

\chapter{Graph products}

\section*{10a. Small graphs}

We have seen in the previous chapter how to associate to any finite graph $X$ its symmetry group $G(X)$, and we have seen as well some basic properties of the correspondence $X\to G(X)$. Motivated by this, our goal in this chapter will be that of systematically computing the symmetry groups $G(X)$, for as many graphs $X$ that we can.

\bigskip

Easy task, you would say, because we already have $G(X)$ for all the $N$-gons, and since $\infty+n=\infty$, no need for new computations, in order to improve our results. Jokes left aside, we certainly need here some precise objectives and strategy, so let us formulate:

\begin{goal}
Take the integers one by one, $N=2,3,4,5,6,\ldots$ and look at all graphs of order $|X|=N$. For each such graph, find a decomposition of type
$$X=Y\times Z$$
with $\times$ being a certain graph product, and $|Y|,|Z|<N$, then work out a formula of type
$$G(X)=G(Y)\times G(Z)$$
with $\times$ being a certain group product, adapted to the above graph product, as to solve the problem by recurrence. When stuck, find some ad-hoc methods for dealing with $X$.
\end{goal}

This plan sounds quite reasonable, so let us see how this works. At very small values of $N=2,3,4,5,6,\ldots$ things are quite clear, and in fact no even need here for a plan, because the symmetry group is just obvious. However, since we will deal with higher $N$ afterwards, it is useful to resist the temptation of simply recording the value of $G(X)$, and stick to the above plan, or at least record the various graph products $\times$ that we meet on the way, and the corresponding adapted group products $\times$ too.

\bigskip

Getting started now, at $N=2,3$ everything is trivial, but let us record this:

\begin{proposition}
At $N=2,3$ we have, with no product operations involved:
\begin{enumerate}
\item Two points $\bullet\ \bullet$ and the segment $\bullet-\bullet$, with symmetry group $\mathbb Z_2$.

\item Three points $\bullet\ \bullet\ \bullet$ and the triangle $\triangle$, with symmetry group $S_3$.

\item The graphs $\bullet-\bullet\ \ \ \bullet$ and $\bullet-\bullet-\bullet$, with symmetry group $\mathbb Z_2$.
\end{enumerate}
\end{proposition}

\begin{proof}
All this is self-explanatory, but us record however a few observations:

\medskip

(1) All pairs of graphs in the above appear from each other via complementation, with this being obvious in all cases, save perhaps for (3), where the graphs are as follows:
$$
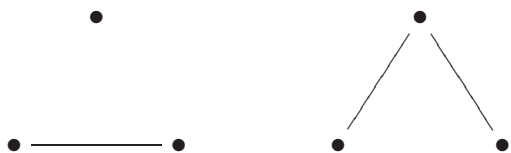
\xymatrix@R=36pt@C=18pt{
&\bullet&&\ &&\bullet\ar@{-}[dl]\ar@{-}[dr]\\
\bullet\ar@{-}[rr]&&\bullet&&\bullet&&\bullet
}$$

(2) Now since we have $G(X)=G(X^c)$, this suggests listing our graphs up to complementation. However, this is a quite bad idea, in practice, because believe me, you will end up sleeping bad at night, with thoughts of type damn, I forgot in my list this or that beautiful graph, only to realize later, after some computations in your head, that the graph in question was in fact complementary to a less beautiful graph, from your list.

\medskip

(3) Another comment, subjective too, concerns the labeling of the groups that we found. For instance $\mathbb Z_1=D_1=S_1$, and $\mathbb Z_2=D_2=S_2$, and $D_3=S_3$. Our policy will be that of regarding $\mathbb Z_N$ as the simplest group, because this is what this group is, and using it preferentially. Followed by $S_N$, because at least in questions regarding permutation groups, this is the second simplest group. As for $D_N$ and other more specialized groups, we will only use them when needed, and with $D_N$ being third on our list.

\medskip

(4) As a mathematical comment now, however trivial, the fact that we have no product operations involved conceptually comes from the fact that $2,3$ are both prime. 

\medskip

(5) However, in relation with this, observe that we have $S_3=D_3=\mathbb Z_3\rtimes\mathbb Z_2$, but this does not correspond to anything, at the level of the graphs $\bullet\ \bullet\ \bullet$ and $\triangle$. So, good to know, the world of graphs is somehow more rigid than that of the groups. 
\end{proof}

At $N=4$ now, which is both big enough, and composite, several new phenomena appear. First, we have the presence of ``ugly'' graphs, of the following type, that no one will be ever interested in, and which do not decompose as products either:
$$
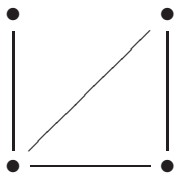
\xymatrix@R=45pt@C=45pt{
\bullet\ar@{-}[d]&\bullet\ar@{-}[d]\ar@{-}[dl]\\
\bullet\ar@{-}[r]&\bullet
}$$

So, leaving aside now this ugly graph, and its symmetry group $\mathbb Z_2$, and all sorts of other similar graphs, functioning on the same principle ``too many vertices, for not much symmetry for the buck'', we are left, quite obviously, with the transitive graphs, which do have lots of symmetry. So, let us update our Goal 10.1, as follows:

\begin{update}
In relation with our original program, we will restrict the attention to the transitive graphs. Moreover, we will list these transitive graphs following a 
$$\{0,N-1\}\quad,\quad\{1,N-2\}\quad,\quad\{2,N-3\}\quad,\quad\ldots$$
scheme, according to their valence, and coupled via complementation.
\end{update}

Here the scheme mentioned at the end is what comes out from Proposition 10.2 (1,2,3), and this is something quite self-explanatory, that will become clear as we work out more examples. Going back now to $N=4$, with the above update made, we have: 

\begin{proposition}
The transitive graphs at $N=4$ are at follows:
\begin{enumerate}
\item Valence $0,3$: the empty graph and the tetrahedron, with symmetry group $S_4$,
$$
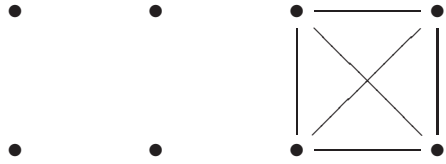
\xymatrix@R=40pt@C=40pt{
\bullet&\bullet&\bullet\ar@{-}[d]\ar@{-}[dr]&\bullet\ar@{-}[d]\ar@{-}[l]\ar@{-}[dl]\\
\bullet&\bullet&\bullet\ar@{-}[r]&\bullet
}$$

\item Valence $1,4$: the two segments and the square, namely
$$
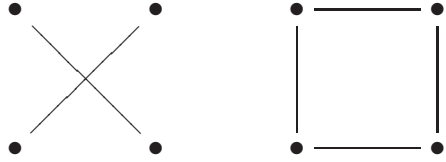
\xymatrix@R=40pt@C=40pt{
\bullet\ar@{-}[dr]&\bullet\ar@{-}[dl]&\bullet\ar@{-}[d]&\bullet\ar@{-}[d]\ar@{-}[l]\\
\bullet&\bullet&\bullet\ar@{-}[r]&\bullet
}$$
both products of a segment with itself, with symmetry group $D_4=\mathbb Z_4\rtimes\mathbb Z_2$.
\end{enumerate}
\end{proposition}

\begin{proof}
As before with Proposition 10.2, all this is trivial and self-explanatory, but some comments are in order, in regards with the last assertions, namely:

\medskip

(1) Regarding the ``products of a segment with itself'' assertion, this is of course something informal, and very intuitive, which remains of course to be clarified. Thus, to be added on our to-do list, two different product operations for graphs $\times$ to be constructed, as to make this work. And no worries for this, we will be back to it, very soon.

\medskip

(2) As another comment, even when assuming that we managed to solve (1), with some suitable product operations for graphs $\times$ constructed, it is quite unclear how we can get, via these operations, the formula $G(X)=\mathbb Z_4\rtimes\mathbb Z_2$. Thus, a potential source of worries, and again, this is just a comment, and we will be back to this, very soon.
\end{proof}

At $N=5$ now, a small prime, things a bit similar to $N=2,3$, and we have:

\begin{proposition}
The transitive graphs at $N=5$ are as follows:
\begin{enumerate}
\item Valence $0,4$: the empty graph and the simplex, with $G(X)=S_5$.

\item Valence $1,3$: none, since $N=5$ is odd.

\item Valence $2$: the pentagon, self-complementary, with $G(X)=D_5$.
\end{enumerate}
\end{proposition}

\begin{proof}
As before, everything is self-explanatory, with however the fact that the pentagon is self-dual being remarkable, definitely to be remembered, the picture being:
$$
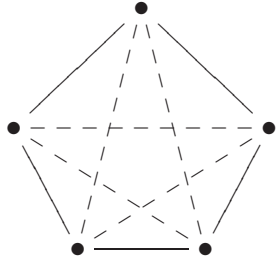
\xymatrix@R=13pt@C=11pt{
&&\bullet\ar@{-}[ddrr]\ar@{-}[ddll]\ar@{--}[ddddl]\ar@{--}[ddddr]\\
&&&&\\
\bullet\ar@{-}[ddr]\ar@{--}[ddrrr]&&&&\bullet\ar@{-}[ddl]\ar@{--}[ddlll]\ar@{--}[llll]\\
&&&&\\
&\bullet\ar@{-}[rr]&&\bullet&&
}$$

Thus, we are led to the conclusions in the statement.
\end{proof}

At $N=6$ now, composite number, things get interesting again, and we have:

\begin{theorem}
The transitive graphs at $N=6$ are as follows:
\begin{enumerate}
\item Valence $0,5$: the empty graph and the simplex, with $G(X)=S_6$.

\item Valence $1,4$: the $3$ segments and the star graph, with $G(X)=H_3$.

\item Valence $2,3$: the hexagon and the prism, with $G(X)=D_6$.

\item Valence $2,3$ too: the $2$ triangles and the wheel/utility graph, products.
\end{enumerate}
\end{theorem}

\begin{proof}
As before, everything here is self-explanatory, the idea being as follows:

\medskip

(1) Nothing special to be said about the empty graph and the simplex.

\medskip

(2) At valence 1, we obviously have as only solution the 3 segments. Regarding now the complementary graph, this looks as follows, like a star:
$$
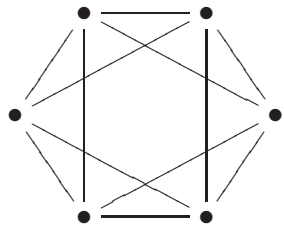
\xymatrix@R=26pt@C=13pt{
&\bullet\ar@{-}[dl]\ar@{-}[rr]\ar@{-}[dd]\ar@{-}[drrr]&&\bullet\ar@{-}[dr]\ar@{-}[dd]\ar@{-}[dlll]\\
\bullet\ar@{-}[dr]&&&&\bullet\ar@{-}[dl]\\
&\bullet\ar@{-}[rr]\ar@{-}[urrr]&&\bullet\ar@{-}[ulll]
}$$

As for the symmetry group assertion, this follows from what we know about the hyperoctahedral groups. Indeed, the hyperoctahedral group $H_3$ appears by definition as symmetry group of the centered hypercube $\square_3\subset\mathbb R^3$. But this group is also, obviously, the symmetry group of the space formed by the segments $[-1,1]$ along each coordinate axis. We conclude that the symmetry group of the 3 segments is indeed $H_3$.

\medskip

(3) At valence 2 now, which is the case left, along with the complementary case of valence 3, we have two solutions, namely the hexagon and the 2 triangles. Regarding the hexagon, have we have $G(X)=D_6$, and the only issue left is that of identifying the complementary graph. But this complementary graph is as follows, and when thinking a bit, by pulling one triangle, and then rotating it, you will see a prism here:
$$
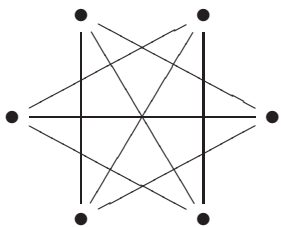
\xymatrix@R=26pt@C=13pt{
&\bullet\ar@{-}[ddrr]\ar@{-}[dd]\ar@{-}[drrr]&&\bullet\ar@{-}[ddll]\ar@{-}[dd]\ar@{-}[dlll]\\
\bullet&&&&\bullet\ar@{-}[llll]\\
&\bullet\ar@{-}[urrr]&&\bullet\ar@{-}[ulll]
}$$

\medskip

(4) Still at valence 2, and complementary valence 3, it remains to discuss the 2 triangles, and its complement. Nothing much to be said about the 2 triangles, but in what regards the complementary graph, there is something tricky here. Indeed, we can draw this complementary graph as follows, making it clear that we have a wheel:
$$
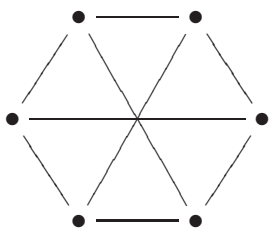
\xymatrix@R=26pt@C=12pt{
&\bullet\ar@{-}[dl]\ar@{-}[rr]\ar@{-}[ddrr]&&\bullet\ar@{-}[dr]\ar@{-}[ddll]\\
\bullet\ar@{-}[dr]\ar@{-}[rrrr]&&&&\bullet\ar@{-}[dl]\\
&\bullet\ar@{-}[rr]&&\bullet
}$$

But now, that we have this wheel, let us color the vertices as follows:
$$
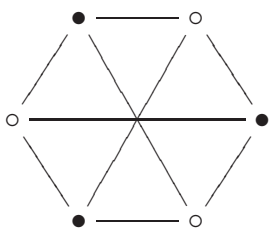
\xymatrix@R=26pt@C=12pt{
&\bullet\ar@{-}[dl]\ar@{-}[rr]\ar@{-}[ddrr]&&\circ\ar@{-}[dr]\ar@{-}[ddll]\\
\circ\ar@{-}[dr]\ar@{-}[rrrr]&&&&\bullet\ar@{-}[dl]\\
&\bullet\ar@{-}[rr]&&\circ
}$$

We can see here a bipartite graph, and by pulling apart the black and white vertices, we conclude that our graph is the utility graph, that we met in chapter 4:
$$
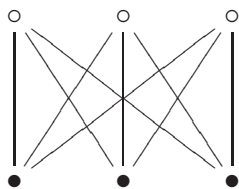
\xymatrix@R=50pt@C=28pt{
\circ\ar@{-}[d]\ar@{-}[dr]\ar@{-}[drr]&\circ\ar@{-}[dl]\ar@{-}[d]\ar@{-}[dr]&\circ\ar@{-}[dll]\ar@{-}[dl]\ar@{-}[d]\\
\bullet&\bullet&\bullet}$$

Finally, in what regards the computation of the symmetry group, the best here is to go back to the original 2 triangles, and draw them as follows:
$$
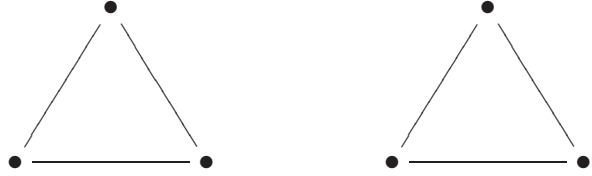
\xymatrix@R=46pt@C=23pt{
&\bullet\ar@{-}[dr]\ar@{-}[dl]&&\ &&\bullet\ar@{-}[dl]\ar@{-}[dr]\\
\bullet\ar@{-}[rr]&&\bullet&&\bullet\ar@{-}[rr]&&\bullet
}$$

But this is obviously a product graph, namely a product between a segment and a triangle, so the symmetry group must appear as some sort of product of $\mathbb Z_2$ and $S_3$. We will leave the details and computations here for a bit later, when discussing more in detail product operations, and the computation of the corresponding symmetry groups. 
\end{proof}

As an interesting conclusion coming from the above discussion, which is something useful and practice, that you always forget, let us record the following fact:

\begin{conclusion}
The utility graph is the wheel,
$$
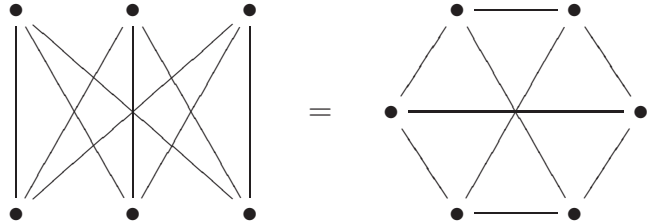
\xymatrix@R=26pt@C=12pt{
\bullet\ar@{-}[dd]\ar@{-}[ddrr]\ar@{-}[ddrrrr]&&\bullet\ar@{-}[ddll]\ar@{-}[dd]\ar@{-}[ddrr]&&\bullet\ar@{-}[ddllll]\ar@{-}[ddll]\ar@{-}[dd]&&&\bullet\ar@{-}[dl]\ar@{-}[rr]\ar@{-}[ddrr]&&\bullet\ar@{-}[dr]\ar@{-}[ddll]\\
&&&&&=&\bullet\ar@{-}[dr]\ar@{-}[rrrr]&&&&\bullet\ar@{-}[dl]\\
\bullet&&\bullet&&\bullet&&&\bullet\ar@{-}[rr]&&\bullet
}$$
with this being best seen via the $2$ triangles missing, in each case.
\end{conclusion}

Here the first assertion looks a bit like a joke, because what can be more useful to mankind, than a wheel. However, it is not from here that the name ``utility graph'' comes from. The story here involves 3 companies, selling gas, water and electricity to 3 customers, and looking for a way to arrange their underground tubes and wires as not to cross. Thus, they are looking to implement their ``utillity graph'', which is the above one, in a planar way, and as we know well from chapter 7, this is not possible, and even leads to some deep theorems about planar graphs, those of Kuratowski and Wagner.

\bigskip

Back to work now, and to our graph enumeration program, at $N=7$, which a small prime, things are a bit similar to those at $N=2,3,5$, as follows:

\begin{proposition}
The transitive graphs at $N=7$ are as follows:
\begin{enumerate}
\item Valence $0,6$: the empty graph and the simplex, with $G(X)=S_7$.

\item Valence $1,5$: none, since $N=7$ is odd.

\item Valence $2,4$: the heptagon and its complement, with $G(X)=D_7$.

\item Valence $3$: none, since $N=7$ is odd.
\end{enumerate}
\end{proposition}

\begin{proof}
As before, everything here is self-explanatory, coming from definitions, and from facts that we know well, with perhaps the only thing to be worked out being the picture of the complement of the heptagon, which looks as follows:
$$
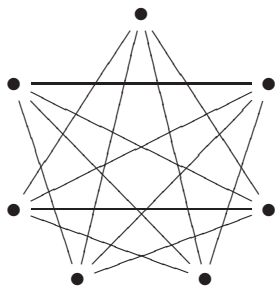
\xymatrix@R=14pt@C=2pt{
&&&&&\bullet\ar@{-}[ddddrr]\ar@{-}[ddddll]\ar@{-}[dddrrrr]\ar@{-}[dddllll]\\
&\bullet\ar@{-}[dddrrrrrr]\ar@{-}[ddrrrrrrrr]&&&&&&&&\bullet\ar@{-}[dddllllll]\ar@{-}[ddllllllll]\ar@{-}[dddll]\ar@{-}[llllllll]\\
\\
&\bullet\ar@{-}[drrrrrr]\ar@{-}[rrrrrrrr]&&&&&&&&\bullet\ar@{-}[dllllll]\\
&&&\bullet\ar@{-}[uuull]&&&&\bullet
}$$

Thus, we are led to the conclusions in the statement.
\end{proof}

\section*{10b. Medium graphs}

Let us discuss now the graphs with bigger number of vertices, $N\geq8$. In what regards the graphs with $8$ vertices, to start with, with $N=8$ being a multiply composite number, things get more complicated, with several new phenomena appearing, as follows:

\begin{theorem}
The transitive graphs at $N=8$ are as follows:
\begin{enumerate}
\item Valence $0,7$: the empty graph and the simplex, with $G(X)=S_8$.

\item Valence $1,6$: the $4$ segments and the thick tyre, products.

\item Valence $2,5$: the octagon and the big globe, with $G(X)=D_8$.

\item Valence $2,5$ too: the $2$ squares and the holy cross, products.

\item Valence $3,4$: the $2$ tetrahedra and the stop sign, products.

\item Valence $3,4$ too: the cube and the metal cross, with $G(X)=H_3$.

\item Valence $3,4$ too: the wheel and the big tent, with $G(X)=S_8$.
\end{enumerate}
\end{theorem}

\begin{proof}
As before, everything here is self-explanatory, the idea being as follows:

\medskip

(1) Nothing much to be said about the empty graph, and its complement.

\medskip

(2) At valence $1,6$ we have the $4$ segments, obviously a product, whose symmetry group we will compute later, and its complement, the thick tyre, which is as follows:
$$
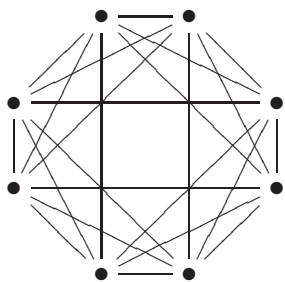
\xymatrix@R=20pt@C=20pt{
&\bullet\ar@{-}[r]\ar@{-}[drr]\ar@{-}[dl]\ar@{-}[ddrr]&\bullet\ar@{-}[dr]\ar@{-}[dll]\ar@{-}[ddll]\\
\bullet\ar@{-}[d]\ar@{-}[ddr]\ar@{-}[rrr]\ar@{-}[ddrr]&&&\bullet\ar@{-}[d]\ar@{-}[ddll]\\
\bullet\ar@{-}[dr]\ar@{-}[uur]\ar@{-}[drr]\ar@{-}[rrr]&&&\bullet\ar@{-}[dl]\ar@{-}[uul]\\
&\bullet\ar@{-}[r]\ar@{-}[urr]\ar@{-}[uuu]&\bullet\ar@{-}[uur]\ar@{-}[uuu]}$$

(3) At valence $2,5$ we have the octagon, whose symmetry group is $G(X)=D_8$, and its complement the big globe, which looks crowded too, as follows:
$$
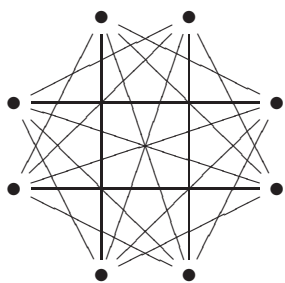
\xymatrix@R=20pt@C=20pt{
&\bullet\ar@{-}[ddrr]\ar@{-}[drr]\ar@{-}[ddd]\ar@{-}[dddr]&\bullet\ar@{-}[ddd]\ar@{-}[dll]\ar@{-}[ddll]\ar@{-}[dddl]\\
\bullet\ar@{-}[ddr]&&&\bullet\ar@{-}[lll]\ar@{-}[dlll]\ar@{-}[ddll]\\
\bullet\ar@{-}[uur]\ar@{-}[drr]&&&\bullet\ar@{-}[lll]\ar@{-}[uul]\ar@{-}[ulll]\\
&\bullet\ar@{-}[urr]&\bullet\ar@{-}[uur]\ar@{-}[uull]}$$

(4) At valence $2,5$ too we have the $2$ squares, obviously a product, whose symmetry group we will compute later, and its complement the holy cross, as follows:
$$
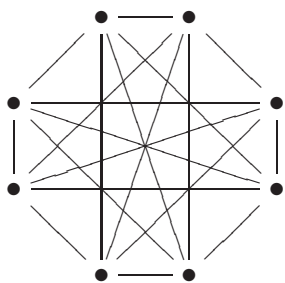
\xymatrix@R=20pt@C=20pt{
&\bullet\ar@{-}[r]\ar@{-}[dddr]\ar@{-}[dl]\ar@{-}[ddrr]&\bullet\ar@{-}[dr]\ar@{-}[ddll]\ar@{-}[dddl]\\
\bullet\ar@{-}[d]\ar@{-}[rrr]\ar@{-}[ddrr]&&&\bullet\ar@{-}[d]\ar@{-}[ddll]\ar@{-}[dlll]\\
\bullet\ar@{-}[dr]\ar@{-}[rrr]&&&\bullet\ar@{-}[dl]\ar@{-}[ulll]\\
&\bullet\ar@{-}[r]\ar@{-}[uuu]&\bullet\ar@{-}[uuu]}$$

(5) At valence $3,4$ now, we first have the $2$ tetrahedra, obviously a product, whose symmetry group we will compute later, and its complement, the stop sign:
$$
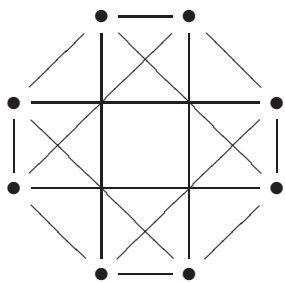
\xymatrix@R=20pt@C=20pt{
&\bullet\ar@{-}[r]\ar@{-}[dl]\ar@{-}[ddrr]&\bullet\ar@{-}[dr]\ar@{-}[ddll]\\
\bullet\ar@{-}[d]\ar@{-}[rrr]\ar@{-}[ddrr]&&&\bullet\ar@{-}[d]\ar@{-}[ddll]\\
\bullet\ar@{-}[dr]\ar@{-}[rrr]&&&\bullet\ar@{-}[dl]\\
&\bullet\ar@{-}[r]\ar@{-}[uuu]&\bullet\ar@{-}[uuu]}$$

(6) Still at valence $3,4$, we have as well the cube, whose symmetry group is $G(X)=H_3$, as we know well, and its complement the metal cross, which is as follows:
$$
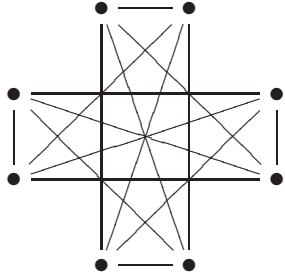
\xymatrix@R=20pt@C=20pt{
&\bullet\ar@{-}[r]\ar@{-}[ddrr]\ar@{-}[dddr]&\bullet\ar@{-}[ddll]\ar@{-}[dddl]\\
\bullet\ar@{-}[d]\ar@{-}[drrr]\ar@{-}[rrr]\ar@{-}[ddrr]&&&\bullet\ar@{-}[d]\ar@{-}[ddll]\\
\bullet\ar@{-}[rrr]\ar@{-}[urrr]&&&\bullet\\
&\bullet\ar@{-}[r]\ar@{-}[uuu]&\bullet\ar@{-}[uuu]}$$

(7) Finally, still at valence $3,4$, we have as well the wheel, whose symmetry group is $G(X)=S_8$, and its complement the big tent, which is as follows:
$$
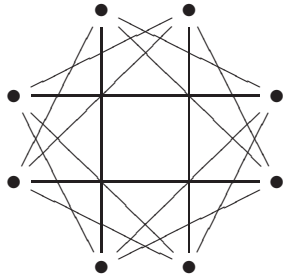
\xymatrix@R=20pt@C=20pt{
&\bullet\ar@{-}[ddrr]\ar@{-}[drr]\ar@{-}[ddd]&\bullet\ar@{-}[ddd]\ar@{-}[dll]\ar@{-}[ddll]\\
\bullet\ar@{-}[ddr]&&&\bullet\ar@{-}[lll]\ar@{-}[ddll]\\
\bullet\ar@{-}[uur]\ar@{-}[drr]&&&\bullet\ar@{-}[lll]\ar@{-}[uul]\\
&\bullet\ar@{-}[urr]&\bullet\ar@{-}[uur]\ar@{-}[uull]}$$

Thus, we are led to the conclusions in the statement.
\end{proof}

At $N=9$ now, things get again easier, but still non-trivial, as follows:

\begin{theorem}
The transitive graphs at $N=9$ are as follows:
\begin{enumerate}
\item Valence $0,8$: the empty graph and the simplex, with $G(X)=S_9$.

\item Valence $2,6$: the nonagon and its complement, with $G(X)=D_9$.

\item Valence $2,6$ too: the $3$ triangles and its complement, products.

\item Valence $4$: the torus graph and its complement, products.

\item Valence $4$ too: the wheel and its complement, with $G(X)=D_9$.
\end{enumerate}
\end{theorem}

\begin{proof}
As before, everything here is self-explanatory, the comments being:

\medskip

(1) Nothing much to be said about the empty graph, and its complement.

\medskip

(2) Nothing much to be said either about the nonagon, and its complement.

\medskip

(3) Regarding the 3 triangles, and its complement, we will compute here $G(X)$ later.

\medskip

(4) At valence $4$ we have an interesting graph which appears, namely the discrete torus, which is the simplest discretization of a torus, as follows, and whose symmetry group we will compute later, when discussing in detail the product operations:
$$
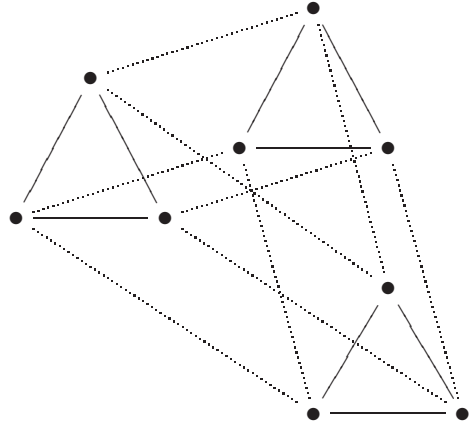
\xymatrix@R=14pt@C=15pt{
&&&&\bullet\ar@{-}[ddl]\ar@{-}[ddr]\ar@{.}[ddddr]\\
&\bullet\ar@{-}[ddl]\ar@{-}[ddr]\ar@{.}[urrr]\ar@{.}[dddrrrr]\\
&&&\bullet\ar@{-}[rr]\ar@{.}[ddddr]&&\bullet\ar@{.}[ddddr]\\
\bullet\ar@{-}[rr]\ar@{.}[urrr]\ar@{.}[dddrrrr]&&\bullet\ar@{.}[urrr]\ar@{.}[dddrrrr]\\
&&&&&\bullet\ar@{-}[ddl]\ar@{-}[ddr]\\
\\
&&&&\bullet\ar@{-}[rr]&&\bullet
}$$

(5) Also at valence 4, we have the wheel and its complement, with the wheel being drawn as follows, forced by the fact that 9 is odd, with two spokes at each vertex:
$$
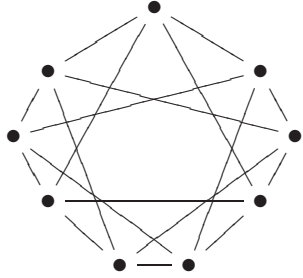
\xymatrix@R=12pt@C=0pt{
&&&&&\bullet\ar@{-}[drrrr]\ar@{-}[dllll]\ar@{-}[dddrrrr]\ar@{-}[dddllll]\\
&\bullet\ar@{-}[dl]\ar@{-}[drrrrrrrrr]&&&&&&&&\bullet\ar@{-}[dr]\ar@{-}[dddlll]\\
\bullet\ar@{-}[urrrrrrrrr]\ar@{-}[dr]&&&&&&&&&&\bullet\ar@{-}[dl]\ar@{-}[ddllllll]\\
&\bullet\ar@{-}[drrr]\ar@{-}[rrrrrrrr]&&&&&&&&\bullet\ar@{-}[dlll]\\
&&&&\bullet\ar@{-}[uuulll]&&\bullet\ar@{-}[ll]\ar@{-}[uullllll]
}$$

Thus, we are led to the conclusions in the statement.
\end{proof}

At $N=10$ now, things get really tough, and we have:

\begin{theorem}
The transitive graphs at $N=10$ are as follows:
\begin{enumerate}
\item Valence $0,9$: the empty graph and the simplex, with $G(X)=S_{10}$.

\item Valence $1,8$: the $5$ segments and its complement, products.

\item Valence $2,7$: the decagon and its complement, with $G(X)=D_{10}$.

\item Valence $2,7$ too: the $2$ pentagons and its complement, product.

\item Valence $3,6$: the wheel and its complement, with $G(X)=D_{10}$.

\item Valence $3,6$ too: the pentagon prism and its complement, products.

\item Valence $3,6$ too: the Petersen graph and its complement.

\item Valence $4,5$: the two $5$-simplices, and its complement, products.

\item Valence $4,5$ too: the $5$-simplex prism and the torch, products.

\item Valence $4,5$ too: the reinforced wheel and its complement, with $G(X)=D_{10}$.

\item Valence $4,5$ too: the alternative reinforced wheel and its complement.
\end{enumerate}
\end{theorem}

\begin{proof}
As before, everything here is self-explanatory, the idea being as follows:

\medskip

(1-6) Nothing much to be said about these graphs, and their complements.

\medskip

(7) At valence 3 we have as well the Petersen graph, which is as follows:
$$
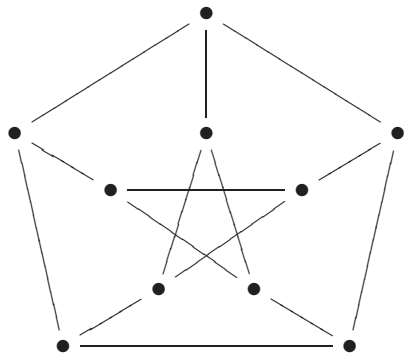
\xymatrix@R=1pt@C=5pt{
&&&&\bullet\ar@{-}[dddddrrrr]\ar@{-}[dddddllll]\\
\\
\\
\\
\\
\bullet\ar@{-}[ddddddddr]&&&&\bullet\ar@{-}[uuuuu]\ar@{-}[ddddddl]\ar@{-}[ddddddr]&&&&\bullet\ar@{-}[ddddddddl]\\
\\
&&
\bullet\ar@{-}[uull]\ar@{-}[ddddrrr]\ar@{-}[rrrr]&&&&\bullet\ar@{-}[uurr]\ar@{-}[ddddlll]\\
\\
\\
\\
&&&\bullet&&\bullet\\
\\
&\bullet\ar@{-}[rrrrrr]\ar@{-}[uurr]&&&&&&\bullet\ar@{-}[uull]}
$$

This graph seems to appear as some sort of product of the cycle $C_5$ with the segment $C_2$, but via a quite complicated procedure. We will be back to this, a bit later.

\medskip

(8) Nothing much to be said about the two $5$-simplices, and its complement, with the computation of the symmetry group here being left for later.

\medskip

(9) Here there is some discussion to be made, because, contrary to our policy so far, listing graphs with smaller valence than their complements first, here we did the opposite. To start with, the 5-simplex prism, having valence 5, is as follows:
$$
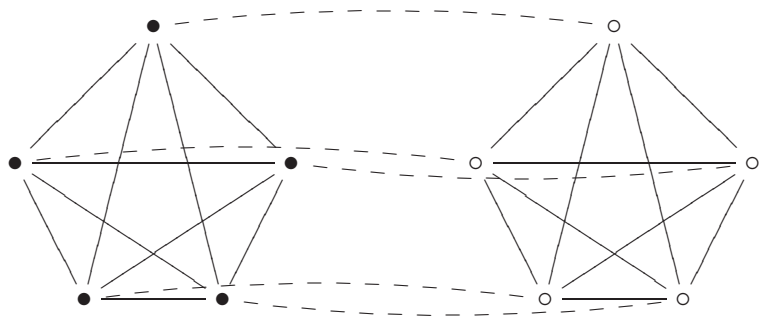
\xymatrix@R=16pt@C=13pt{
&&\bullet\ar@{-}[ddrr]\ar@{-}[ddll]\ar@{-}[ddddl]\ar@{-}[ddddr]\ar@{--}@/^/[rrrrrr]&&&&&&\circ\ar@{-}[ddrr]\ar@{-}[ddll]\ar@{-}[ddddl]\ar@{-}[ddddr]\\
&&&&\\
\bullet\ar@{-}[ddr]\ar@{-}[ddrrr]\ar@{--}@/^/[rrrrrr]&&&&\bullet\ar@{-}[ddl]\ar@{-}[ddlll]\ar@{-}[llll]\ar@{--}@/_/[rrrrrr]&\ \ \ \ \ \ &\circ\ar@{-}[ddr]\ar@{-}[ddrrr]&&&&\circ\ar@{-}[ddl]\ar@{-}[ddlll]\ar@{-}[llll]\\
&&&&\\
&\bullet\ar@{-}[rr]\ar@{--}@/^/[rrrrrr]&&\bullet\ar@{--}@/_/[rrrrrr]&&&&\circ\ar@{-}[rr]&&\circ&&
}$$

As for the complementary graph to this 5-simplex prism, this is the torch graph, having valence 4, which looks as follows:
$$
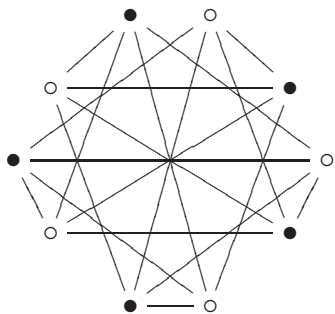
\xymatrix@R=15pt@C=1pt{
&&&&\bullet\ar@{-}[ddddrrr]\ar@{-}[dlll]\ar@{-}[ddrrrrrrr]\ar@{-}[dddlll]&&&\circ\ar@{-}[drrr]\ar@{-}[dddrrr]\ar@{-}[ddlllllll]\ar@{-}[ddddlll]\\
&\circ&&&&&&&&&\bullet\ar@{-}[dddlll]\ar@{-}[ddlllllllll]\ar@{-}[lllllllll]\\
\bullet\ar@{-}[dr]&&&&&&&&&&&\circ\ar@{-}[dl]\ar@{-}[ddlllllll]\ar@{-}[lllllllllll]\\
&\circ&&&&&&&&&\bullet\ar@{-}[uulllllllll]\ar@{-}[lllllllll]\\
&&&&\bullet\ar@{-}[uuulll]&&&\circ\ar@{-}[lll]\ar@{-}[uulllllll]
}$$

(10-11) There is some discussion to be made here, concerning the two possible reinforced wheels on $10$ vertices, having 2 spokes at each vertex. First we have a reinforced wheel as follows, whose symmetry group is $D_{10}$, as said in the statement:
$$
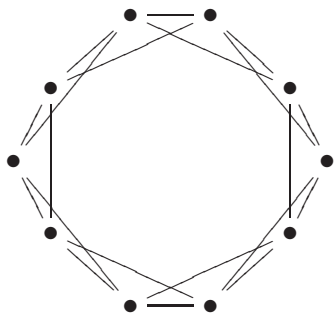
\xymatrix@R=15pt@C=1pt{
&&&&\bullet\ar@{-}[rrr]\ar@{-}[dlll]\ar@{-}[drrrrrr]\ar@{-}[ddllll]&&&\bullet\ar@{-}[drrr]\ar@{-}[ddrrrr]\ar@{-}[dllllll]\\
&\bullet\ar@{-}[dl]&&&&&&&&&\bullet\ar@{-}[dr]\ar@{-}[dd]\\
\bullet\ar@{-}[dr]&&&&&&&&&&&\bullet\ar@{-}[dl]\ar@{-}[ddllll]\\
&\bullet\ar@{-}[drrr]\ar@{-}[uu]&&&&&&&&&\bullet\ar@{-}[dlll]\ar@{-}[dllllll]\\
&&&&\bullet\ar@{-}[uullll]&&&\bullet\ar@{-}[lll]\ar@{-}[ullllll]
}$$

But then we have as well the following alternative reinforced wheel, whose symmetry group is a bit more complicated to compute, and we will discuss this later:
$$
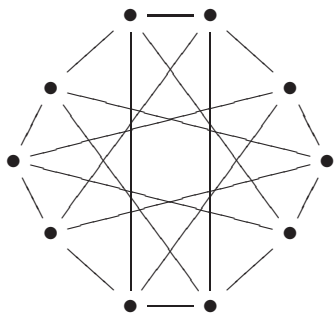
\xymatrix@R=15pt@C=1pt{
&&&&\bullet\ar@{-}[rrr]\ar@{-}[dlll]\ar@{-}[dddrrrrrr]\ar@{-}[dddd]&&&\bullet\ar@{-}[drrr]\ar@{-}[dddd]\ar@{-}[dddllllll]\\
&\bullet\ar@{-}[dl]&&&&&&&&&\bullet\ar@{-}[dr]\ar@{-}[dddllllll]\ar@{-}[dllllllllll]\\
\bullet\ar@{-}[dr]&&&&&&&&&&&\bullet\ar@{-}[dl]\ar@{-}[ullllllllll]\ar@{-}[dllllllllll]\\
&\bullet\ar@{-}[drrr]&&&&&&&&&\bullet\ar@{-}[dlll]\ar@{-}[ullllllllll]\\
&&&&\bullet&&&\bullet\ar@{-}[lll]\ar@{-}[uuullllll]
}$$

Thus, we are led to the conclusions in the statement.
\end{proof}

All the above is quite concerning, especially in what concerns the Petersen graph, which looks quite enigmatic. However, before leaving the subject, let us record as well the result at $N=11$. Here the order being a small prime, things get again easier:

\begin{proposition}
The transitive graphs at $N=11$ are as follows:
\begin{enumerate}
\item Valence $0,10$: the empty graph and the simplex, with $G(X)=S_{11}$.

\item Valence $2,8$: the hendecagon and its complement, with $G(X)=D_{11}$.

\item Valence $4,6$: the reinforced wheel and its complement, with $G(X)=D_{11}$.

\item Valence $4,6$: the other reinforced wheel and complement, with $G(X)=D_{11}$.
\end{enumerate}
\end{proposition}

\begin{proof}
Here again business as usual, with the only discussion concerning the two possible reinforced wheels, with 2 spokes at each vertex. A first such wheel is follows:
$$
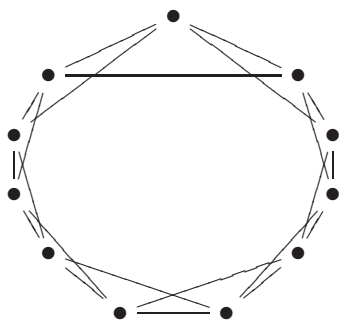
\xymatrix@R=10pt@C=0pt{
&&&&&&\bullet\ar@{-}[drrrrr]\ar@{-}[dlllll]\ar@{-}[ddllllll]\ar@{-}[ddrrrrrr]\\
&\bullet\ar@{-}[dl]\ar@{-}[ddl]&&&&&&&&&&\bullet\ar@{-}[dr]\ar@{-}[llllllllll]\ar@{-}[ddr]\\
\bullet\ar@{-}[d]\ar@{-}[ddr]&&&&&&&&&&&&\bullet\ar@{-}[d]\ar@{-}[ddl]\\
\bullet\ar@{-}[dr]\ar@{-}[ddrrrr]&&&&&&&&&&&&\bullet\ar@{-}[dl]\ar@{-}[ddllll]\\
&\bullet\ar@{-}[drrr]\ar@{-}[drrrrrrr]&&&&&&&&&&\bullet\ar@{-}[dlll]\ar@{-}[dlllllll]\\
&&&&\bullet&&&&\bullet\ar@{-}[llll]
}$$

But then, we have as well a second reinforced wheel, as follows:
$$
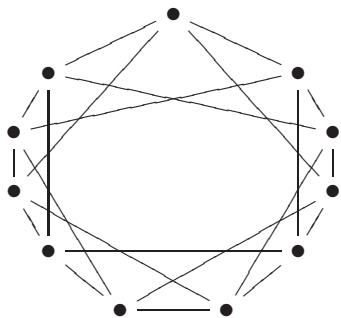
\xymatrix@R=10pt@C=0pt{
&&&&&&\bullet\ar@{-}[drrrrr]\ar@{-}[dlllll]\ar@{-}[dddllllll]\ar@{-}[dddrrrrrr]\\
&\bullet\ar@{-}[dl]\ar@{-}[ddd]\ar@{-}[drrrrrrrrrrr]&&&&&&&&&&\bullet\ar@{-}[dr]\ar@{-}[ddd]\ar@{-}[dlllllllllll]\\
\bullet\ar@{-}[d]\ar@{-}[dddrrrr]&&&&&&&&&&&&\bullet\ar@{-}[d]\ar@{-}[dddllll]\\
\bullet\ar@{-}[dr]\ar@{-}[ddrrrrrrrr]&&&&&&&&&&&&\bullet\ar@{-}[dl]\ar@{-}[ddllllllll]\\
&\bullet\ar@{-}[drrr]&&&&&&&&&&\bullet\ar@{-}[dlll]\ar@{-}[llllllllll]\\
&&&&\bullet&&&&\bullet\ar@{-}[llll]
}$$

Thus, we are led to the conclusions in the statement.
\end{proof}

And we will stop here our small $N$ study, for several reasons:

\bigskip

(1) First of all, because we did some good work in the above, and found many types of graph products, that we have now to investigate in detail.

\bigskip

(2) Second, because of the Petersen graph at $N=10$, which obviously brings us into uncharted territory, and needs some serious study, before going further. 

\bigskip

(3) And third, because at $N=12$, which is a dazzlingly composite number, things explode, with the number of transitive graphs here being a mighty $64$.

\section*{10c. Standard products}

For the transitive graphs, that we are mostly interested in, the point is that, according to the above, at very small values of the order, $N=2,\ldots,9$, these always decompose as products, via three main types of graph products, constructed as follows:

\index{direct product}
\index{Cartesian product}
\index{lexicographic product}

\begin{definition}
Let $X,Y$ be two finite graphs.
\begin{enumerate}
\item The direct product $X\times Y$ has vertex set $X\times Y$, and edges:
$$(i,\alpha)-(j,\beta)\Longleftrightarrow i-j,\, \alpha-\beta$$

\item The Cartesian product $X\,\square\,Y$ has vertex set $X\times Y$, and edges:
$$(i,\alpha)-(j,\beta)\Longleftrightarrow i=j,\, \alpha-\beta\mbox{ \rm{or} }i-j,\alpha=\beta$$

\item The lexicographic product $X\circ Y$ has vertex set $X\times Y$, and edges:
$$(i,\alpha)-(j,\beta)\Longleftrightarrow\alpha-\beta\mbox{ \rm{or} }\alpha=\beta,\,
 i-j$$
\end{enumerate}
\end{definition}

Several comments can be made here. First, the direct product $X\times Y$ is the usual one in a categorical sense, and we will leave clarifying this observation as an exercise. The Cartesian product $X\,\square\,Y$ is quite natural too from a geometric perspective, for instance because a product by a segment gives a prism. As for the lexicographic product $X\circ Y$, this is something interesting too, obtained by putting a copy of $X$ at each vertex of $Y$.

\bigskip

At the level of symmetry groups, several things can be said, and we first have:

\index{wreath product}

\begin{theorem}
We have group embeddings as follows, for any graphs $X,Y$,
$$G(X)\times G(Y)\subset G(X \times Y)$$
$$G(X)\times G(Y)\subset G(X\,\square\,Y)$$
$$G(X)\wr G(Y)\subset G(X\circ Y)$$
but these embeddings are not always isomorphisms.
\end{theorem}

\begin{proof}
The fact that we have indeed embeddings as above is clear from definitions. As for the counterexamples, in each case, these are easy to construct as well, provided by our study of small graphs, at $N=2,\ldots,11$, and we will leave this as an exercise.
\end{proof}

The problem now is that of deciding when the embeddings in Theorem 10.14 are isomorphisms. And this is something non-trivial, because there are both examples and counterexamples for these isomorphisms, coming from the various computations of symmetry groups that we did in the above, at $N=2,\ldots,11$. We will see, however, that the problem can be solved, via a technical study, of spectral theory flavor.

\bigskip

Now speaking technical algebra and spectral theory, it is good time to go back to the permutation groups $G\subset S_N$, as introduced and studied in chapter 9, and study them a bit more, from an algebraic perspective. We first have the following basic fact:

\begin{theorem}
Given a subgroup $G\subset S_N$, regarded as matrix group via 
$$G\subset S_N\subset O_N$$
the standard coordinates of the group elements, $u_{ij}(g)=g_{ij}$, are given by: 
$$u_{ij}=\chi\left(\sigma\in G\Big|\sigma(j)=i\right)$$
Moreover, these functions $u_{ij}:G\to\mathbb C$ generate the algebra $C(G)$.
\end{theorem}

\begin{proof}
Here the first assertion comes from the fact that the entries of the permutation matrices $\sigma\in S_N\subset O_N$, acting as $\sigma(e_i)=e_{\sigma(i)}$, are given by the following formula:
$$\sigma_{ij}=\begin{cases}
1&{\rm if}\ \sigma(j)=i\\
0&{\rm otherwise}
\end{cases}$$

As for the second assertion, this comes from the Stone-Weierstrass theorem, because the coordinate functions $u_{ij}:G\to\mathbb C$ obviously separate the group elements $\sigma\in G$.
\end{proof}

We are led in this way to the following definition:

\begin{definition}
The magic matrix associated to a permutation group $G\subset S_N$ is the $N\times N$ matrix of characteristic functions
$$u_{ij}=\chi\left(\sigma\in G\Big|\sigma(j)=i\right)$$
with the name ``magic'' coming from the fact that, on each row and each column, these characteristic functions sum up to $1$.
\end{definition}

The interest in this notion comes from the fact, that we know from Theorem 10.15, that the entries of the magic matrix generate the algebra of functions on our group:
$$C(G)=<u_{ij}>$$

We will talk more in detail later about such matrices, and their correspondence with the subgroups $G\subset S_N$, and what can be done with it, in the general framework of representation theory. However, for making our point, here is the general principle:

\begin{principle}
Everything that you can do with your group $G\subset S_N$ can be expressed in terms of the magic matrix $u=(u_{ij})$, quite often with good results.
\end{principle}

This principle comes from the above Stone-Weierstrass result, $C(G)=<u_{ij}>$. Indeed, when coupled with some basic spectral theory, and more specifically with the Gelfand theorem from operator algebras, this result tells us that our group $G$ appears as the spectrum of the algebra $<u_{ij}>$, therefore leading to the above principle. But more on this later in this book, when discussing spectral theory, and the Gelfand theorem.

\bigskip

For the moment, we will take Definition 10.16 as it is, something technical of group theory, of rather functional analysis flavor, that we can use in our proofs when needed. And we will take Principle 10.17 also as it is, namely a claim that this is useful indeed.

\bigskip

As an illustration for all this, in relation with the graphs, we have:

\begin{theorem}
Given a subgroup $G\subset S_N$, the transpose of its action map $X\times G\to X$ on the set $X=\{1,\ldots,N\}$, given by $(i,\sigma)\to\sigma(i)$, is given by:
$$\Phi(e_i)=\sum_je_j\otimes u_{ji}$$
Also, in the case where we have a graph with $N$ vertices, the action of $G$ on the vertex set $X$ leaves invariant the edges precisely when we have
$$du=ud$$
with $d$ being as usual the adjacency matrix of the graph.
\end{theorem}

\begin{proof}
There are several things going on here, the idea being as follows:

\medskip

(1) Given a subgroup $G\subset S_N$, if we set $X=\{1,\ldots,N\}$, we have indeed an action map as follows, and with the reasons of using $X\times G$ instead of the perhaps more familiar $G\times X$ being dictated by some quantum algebra, that we will do later in this book:
$$a:X\times G\to X\quad,\quad a(i,\sigma)=\sigma(i)$$

(2) Now by transposing this map, we obtain a morphism of algebras, as follows:
$$\Phi:C(X)\to C(X)\otimes C(G)\quad,\quad \Phi(f)(i,\sigma)=f(\sigma(i))$$

When evaluated on the Dirac masses, this map $\Phi$ is then given by:
$$\Phi(e_i)(j,\sigma)=e_i(\sigma(j))=\delta_{\sigma(j)i}$$

Thus, in tensor product notation, we have the following formula, as desired:
$$\Phi(e_i)(j,\sigma)=\left(\sum_je_j\otimes u_{ji}\right)(j,\sigma)$$

(3) Regarding now the second assertion, observe first that we have:
$$(du)_{ij}(\sigma)
=\sum_kd_{ik}u_{kj}(\sigma)
=\sum_kd_{ik}\delta_{\sigma(j)k}
=d_{i\sigma(j)}$$

On the other hand, we have as well the following formula:
$$(ud)_{ij}(\sigma)
=\sum_ku_{ik}(\sigma)d_{kj}
=\sum_k\delta_{\sigma(k)i}d_{kj}
=d_{\sigma^{-1}(i)j}$$

Thus $du=ud$ reformulates as $d_{ij}=d_{\sigma(i)\sigma(j)}$, which gives the result.
\end{proof}

So long for magic unitaries, and their basic properties, and we will be back to this, on several occasions, in what follows. In fact, the magic matrices will get increasingly important, as the present book develops, because not far away from now, when starting to talk about quantum permutation groups $G$, and their actions on the graphs $X$, these beasts will not really exist, as concrete objects $G$, but their associated magic matrices $u=(u_{ij})$ will exist, and we will base our whole study on them. More on this later.

\bigskip

Back to graphs now, we want to know when the embeddings in Theorem 10.14 are isomorphisms. In what regards the first two products, we have here the following result, coming with a proof from \cite{bbi}, obtained by using the magic matrix technology:

\index{connected graph}
\index{regular graph}

\begin{theorem}
Let $X$ and $Y$ be finite connected regular graphs. If their spectra $\{\lambda\}$ and $\{\mu\}$ do not contain $0$ and satisfy
$$\big\{\lambda_i/\lambda_j\big\}\cap\big\{\mu_k/\mu_l\big\}=\{1\}$$
then $G(X\times Y)=G(X)\times G(Y)$. Also, if their spectra satisfy
$$\big\{\lambda_i-\lambda_j\big\}\cap\big\{\mu_k-\mu_l\big\}=\{0\}$$
then $G(X\,\square\,Y)=G(X)\times G(Y)$.
\end{theorem}

\begin{proof}
This is something quite standard, the idea being as follows:

\medskip

(1) First, we know from Theorem 10.14 that we have embeddings as follows, valid for any two graphs $X,Y$, and coming from definitions:
$$G(X)\times G(Y)\subset G(X\times Y)$$
$$G(X)\times G(Y)\subset G(X\,\square\,Y)$$

(2) Now let $\lambda_1$ be the valence of $X$. Since $X$ is regular we have $\lambda_1\in Sp(X)$, with $1$ as eigenvector, and since $X$ is connected $\lambda_1$ has multiplicity 1. Thus if $P_1$ is the orthogonal projection onto $\mathbb C1$, the spectral decomposition of $d_X$ is of the following form:
$$d_X=\lambda_1P_1+\sum_{i\neq1}\lambda_iP_i$$

We have a similar formula for the adjacency matrix $d_Y$, namely:
$$d_Y=\mu_1Q_1+\sum_{j\neq1}\mu_jQ_j$$

(3) But this gives the following formulae for the graph products:
$$d_{X\times Y}=\sum_{ij}(\lambda_i\mu_j)P_{i}\otimes Q_{j}$$
$$d_{X\,\square\,Y}=\sum_{ij}(\lambda_i+\mu_i)P_i\otimes Q_j$$

Here the projections form partitions of unity, and the scalar are distinct, so these are spectral decompositions. The coactions will commute with any of the spectral projections, and so with both $P_1\otimes1$, $1\otimes Q_1$. In both cases the universal coaction $v$ is the tensor product of its restrictions to the images of $P_1\otimes1$, $1\otimes Q_1$, which gives the result.
\end{proof}

Regarding now the lexicographic product, things here are more tricky. Let us first recall that the lexicographic product of two graphs $X\circ Y$ is obtained by putting a copy of $X$ at each vertex of $Y$, the formula for the edges being as follows:
$$(i,\alpha)-(j,\beta)\Longleftrightarrow\alpha-\beta\mbox{ \rm{or} }\alpha=\beta,\,
 i-j$$

In what regards now the computation of the symmetry group, as before we must do here some spectral theory, and we are led in this way to the following result:

\begin{theorem}
Let $X,Y$ be regular graphs, with $X$ connected. If their spectra $\{\lambda_i\}$ and $\{\mu_j\}$ satisfy the condition
$$\big\{\lambda_1-\lambda_i\big|i\neq 1\big\}\cap\big\{-n\mu_j\big\}=\emptyset$$
where $n$ and $\lambda_1$ are the order and valence of $X$, then $G(X\circ Y)=G(X)\wr G(Y)$.   
\end{theorem}

\begin{proof}
This is something quite tricky, the idea being as follows:

\medskip

(1) First, we know from Theorem 10.14 that we have an embedding as follows, valid for any two graphs $X,Y$, and coming from definitions:
$$G(X)\wr G(Y)\subset G(X\circ Y)$$

(2) We denote by $P_i,Q_j$ the spectral projections corresponding to $\lambda_i,\mu_j$. Since $X$ is connected we have $P_1=\mathbb I/n$, and we obtain:
\begin{eqnarray*}
d_{X\circ Y}
&=&d_X\otimes 1+{\mathbb I}\otimes d_Y\\
&=&\left(\sum_i\lambda_iP_i\right)\otimes\left(\sum_jQ_j\right)+\left(nP_1\right)\otimes \left(\sum_i\mu_jQ_j\right)\\
&=&\sum_j(\lambda_1+n\mu_j)(P_1\otimes Q_j) + \sum_{i\not=1}\lambda_i(P_i\otimes 1)
\end{eqnarray*} 

In this formula the projections form a partition of unity and the scalars are distinct, so this is the spectral decomposition of $d_{X\circ Y}$. 

\medskip

(3) Now let $W$ be the universal magic matrix for $X\circ Y$. Then $W$ must commute with all spectral projections, and in particular:
$$[W,P_1\otimes Q_j]=0$$

Summing over $j$ gives $[W,P_1\otimes 1]=0$, so $1\otimes C(Y)$ is invariant under the coaction. So, consider the restriction of $W$, which gives a coaction of $G(X\circ Y)$ on $1\otimes C(Y)$, that we can denote as follows, with $y$ being a certain magic unitary:
$$W(1\otimes e_a)=\sum_b1\otimes e_b\otimes y_{ba}$$

(4) On the other hand, according to our definition of $W$, we can write:
$$W(e_i\otimes 1)=\sum_{jb}e_j\otimes e_b\otimes x_{ji}^b$$  

By multiplying by the previous relation, found in (3), we obtain:
\begin{eqnarray*}
W(e_i\otimes e_a)
&=&\sum_{jb}e_j\otimes e_b\otimes y_{ba}x_{ji}^b\\
&=&\sum_{jb}e_j \otimes e_b\otimes x_{ji}^by_{ba}
\end{eqnarray*}

But this shows that the coefficients of $W$ are of the following form:
$$W_{jb,ia}=y_{ba}x_{ji}^b=x_{ji}^b y_{ba}$$

(5) In order to advance, consider now the following matrix:
$$x^b=(x_{ij}^b)$$

Since the map $W$ above is a morphism of algebras, each row of $x^b$ is a partition  of unity. Also, by using the antipode map $S$, which is transpose to $g\to g^{-1}$, we have:
\begin{eqnarray*}
S\left(\sum_jx_{ji}^{b}\right)
&=&S\left(\sum_{ja}x_{ji}^{b}y_{ba}\right)\\
&=&S\left(\sum_{ja}W_{jb,ia}\right)\\
&=&\sum_{ja}W_{ia,jb}\\
&=&\sum_{ja}x_{ij}^ay_{ab}\\
&=&\sum_ay_{ab}\\
&=&1
\end{eqnarray*}

As a conclusion to this, the matrix $x^b$ constructed above is magic. 

\medskip

(6) We check now that both $x^a,y$ commute with $d_X,d_Y$. We have:
$$(d_{X\circ Y})_{ia,jb} = (d_X)_{ij}\delta_{ab} + (d_Y)_{ab}$$

Thus the two products between $W$ and $d_{X\circ Y}$ are given by:
$$(Wd_{X\circ Y})_{ia,kc}=\sum_j W_{ia,jc} (d_X)_{jk} + \sum_{jb}W_{ia,jb}(d_Y)_{bc}$$
$$(d_{X\circ Y}W)_{ia,kc}=\sum_j (d_X)_{ij} W_{ja,kc} + \sum_{jb}(d_Y)_{ab}W_{jb,kc}$$

(7) Now since the magic matrix $W$ commutes by definition with $d_{X\circ Y}$, the terms on the right in the above equations are equal, and by summing over $c$ we get:
$$\sum_j x_{ij}^a(d_X)_{jk} + \sum_{cb} y_{ab}(d_Y)_{bc}
= \sum_{j} (d_X)_{ij}x_{jk}^a + \sum_{cb} (d_Y)_{ab}y_{bc}$$

The second sums in both terms are equal to the valence of $Y$, so we get:
$$[x^a,d_X]=0$$

Now once again from the formula coming from $[W,d_{X\circ Y}]=0$, we get:
$$[y,d_Y] =0$$

(8) Summing up, the coefficients of $W$ are of the following form, where $x^b$ are magic unitaries commuting with $d_X$, and $y$ is a magic unitary commuting with $d_Y$: 
$$W_{jb,ia}=x_{ji}^by_{ba}$$

But this gives a morphism $C(G(X)\wr G(Y))\to G(X\circ Y)$ mapping $u_{ji}^{(b)}\to x_{ji}^b$ and $v_{ba}\to y_{ba}$, which is inverse to the morphism in (1), as desired.
\end{proof}

As before with the other graph products, there is some further spectral theory to be done here, plus working out examples. There are as well a number of finer results, for instance the theorem of Sabidussi, to be discussed. We will be back to this.

\section*{10d. Kneser graphs}

Good news, we have now enough tools in our bag for getting back to the transitive graphs of order $N=2,\ldots,11$, enumerated and partially studied in the beginning of this chapter, compute the missing symmetry groups there, and end up with a nice table. 

\bigskip

Let us begin with some notations, as follows:

\begin{definition}
We use the following notations for graphs:
\begin{enumerate}
\item $K_N$ is the complete graph on $N$ vertices.

\item $C_N$ is the cycle on $N$ vertices.

\item $P(X)$ the prism over a graph $X$.

\item $C_N^k$ is the cycle with chords of length $k$.
\end{enumerate}
\end{definition}

In what regards the transitive graphs of order $N=2,\ldots,9$, here we have all the symmetry groups already computed, except for 6 of them, which are as follows:

\begin{proposition}
The missing symmetry groups at $N\leq9$ are as follows:
\begin{enumerate}
\item For the two triangles $2K_3$ we obtain $S_3\wr\mathbb Z_2$.

\item For the four segments $4K_2$ we obtain $H_4$.

\item For the two squares $2C_4$ we obtain $H_2\wr\mathbb Z_2$.

\item For the two tetrahedra $2K_4$ we obtain $S_4\wr\mathbb Z_2$.

\item For the three triangles $3K_3$ we obtain $S_3\wr S_3$.

\item For the torus graph $K_3\times K_3$ we obtain $S_3\wr\mathbb Z_2$.
\end{enumerate}
\end{proposition}

\begin{proof}
This follows indeed by applying the various product results from the previous section, and we will leave the details here as an instructive exercise.
\end{proof}

Regarding the graphs of order $N=10,11$, things fine at $N=11$, but at $N=10$ we have 7 symmetry groups left. The first 6 of them can be computed as follows:

\begin{proposition}
The missing symmetry groups at $N=10$ are as follows:
\begin{enumerate}
\item For the five segments $5K_2$ we obtain $H_5$.

\item For the two pentagons $2C_5$ we obtain $D_5\wr\mathbb Z_2$.

\item For the pentagon prism $P(C_5)$ we obtain $D_{10}$.

\item For the two $5$-simplices $2K_5$ we obtain $S_5\wr\mathbb Z_2$.

\item For the $5$-simplex prism $P(K_5)$ we obtain $S_5\times\mathbb Z_2$.

\item For the second reinforced wheel $C_{10}^4$ we obtain $\mathbb Z_2\wr D_5$.
\end{enumerate}
\end{proposition}

\begin{proof}
This follows too by applying the various product results from the previous section, and we will leave the details here as an instructive exercise.
\end{proof}

We have kept the best for the end. Time now to face the only graph left, which is the famous Petersen graph, that we know well since chapter 7, which looks as follows:
$$
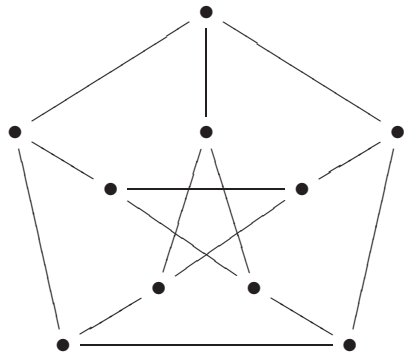
\xymatrix@R=1pt@C=5pt{
&&&&\bullet\ar@{-}[dddddrrrr]\ar@{-}[dddddllll]\\
\\
\\
\\
\\
\bullet\ar@{-}[ddddddddr]&&&&\bullet\ar@{-}[uuuuu]\ar@{-}[ddddddl]\ar@{-}[ddddddr]&&&&\bullet\ar@{-}[ddddddddl]\\
\\
&&
\bullet\ar@{-}[uull]\ar@{-}[ddddrrr]\ar@{-}[rrrr]&&&&\bullet\ar@{-}[uurr]\ar@{-}[ddddlll]\\
\\
\\
\\
&&&\bullet&&\bullet\\
\\
&\bullet\ar@{-}[rrrrrr]\ar@{-}[uurr]&&&&&&\bullet\ar@{-}[uull]}
$$

Intuition suggests that this graph should appear as a product of the cycle $C_5$ with the segment $C_2$. However, there are many bugs with this idea, which does not work.

\bigskip

In order to view the Petersen graph as part of a larger family, and have some general theory going, we have to proceed in a quite unexpected way, as follows:

\begin{proposition}
The Petersen graph is part of the Kneser graph family,
$$P_{10}=K(5,2)$$
with $K(n,s)$ having as vertices the $s$-element subsets of $\{1,\ldots,n\}$, and with the edges being drawn between distinct subsets.
\end{proposition}

\begin{proof}
Consider indeed the Kneser graph $K(n,s)$, as constructed above. At $n=5$, $s=2$ the vertices are the $\binom{5}{2}=10$ subsets of $\{1,\ldots,5\}$ having 2 elements, and since each such subset $\{p,q\}$ is disjoint from exactly 3 other such subsets $\{u,v\}$, our graph is trivalent, and when drawing the picture, we obtain indeed the Petersen graph.
\end{proof}

Many things can be said about the Kneser graphs, and among others, we have:

\begin{theorem}
The Kneser graphs $K(n,s)$ have the following properties:
\begin{enumerate}
\item $K(n,1)$ is the complete graph $K_n$.

\item $K(n,2)$ is the complement of the line graph $L(K_n)$.

\item $K(2m-1,m-1)$ is the so-called odd graph $O_m$.

\item The symmetry group of $K(n,s)$ is the symmetric group $S_n$.
\end{enumerate}
In particular, $P_{10}=K(5,2)=L(K_5)^c=O_3$, having symmetry group $S_5$.
\end{theorem}

\begin{proof}
All this is quite self-explanatory, the idea being as follows:

\medskip

(1) This is something trivial, coming from definitions.

\medskip

(2) This is clear too from definitions, with the line graph $L(X)$ of a given graph $X$ being by definition the incidence graph of the edges of $X$.

\medskip

(3) This stands for the definition of $O_m$, as being the graph $K(2m-1,m-1)$.

\medskip

(4) This is clear too, coming from the definition of the Kneser graphs.
\end{proof}

All this is very nice, and time to conclude. We have the following table, based on our various results above, containing all transitive graphs up to $N=11$ vertices, up to complementation, and their symmetry groups, written by using the conventions from Definition 10.21, and with the extra convention that $C_n^+$ with $n$ even denote the wheels:
\begin{center}\begin{tabular}[t]{|l|l|l|l|}\hline
Order&Graph&Symmetry group\\ 
\hline\hline 2&$K_2$&$\mathbb Z_2$\\
\hline\hline 3&$K_3$&$S_3$\\ 
\hline\hline 4&$2K_2$&$H_2$\\
\hline 4&$K_4$&$S_4$\\
\hline\hline 5&$C_5$&$D_5$\\ 
\hline5&$K_5$&$S_5$\\ 
\hline\hline 6&$C_6$&$D_6$\\ 
\hline 6&$2K_3$&$S_3\wr\mathbb Z_2$\\
\hline 6&$3K_2$&$H_3$\\ 
\hline 6&$K_6$&$S_6$\\
\hline\hline 7&$C_7$&$D_7$\\ 
\hline 7&$K_7$&$S_7$\\ 
\hline\hline 8&$C_8$, $C_8^+$&$D_8$\\
\hline 8&$P(C_4)$& $H_3$\\ 
\hline 8&$2K_4$&$S_4\wr \mathbb Z_2$\\
\hline 8&$2C_4$& $H_2\wr\mathbb Z_2$\\ 
\hline 8&$4K_2$&$H_4$\\ 
\hline 8&$K_8$&$S_8$\\
\hline\hline 9&$C_9$, $C_9^3$&$D_9$\\ 
\hline 9 & $K_3\times K_3$&$S_3\wr\mathbb Z_2$\\ 
\hline 9&$3K_3$&$S_3\wr S_3$\\ 
\hline 9&$K_9$&$S_9$\\ 
\hline \hline 10&$C_{10}$, $C_{10}^2$, $C_{10}^+$, $P(C_5)$&$D_{10}$\\ 
\hline 10 &$P(K_5)$&$S_5\times\mathbb  Z_2$\\ 
\hline10&$C_{10}^4$&$\mathbb Z_2\wr D_5$\\ 
\hline10&$2C_5$&$D_5\wr\mathbb  Z_2$\\ 
\hline10&$2K_{5}$&$S_5\wr\mathbb  Z_2$\\ 
\hline10&$5K_2$&$H_5$\\ 
\hline 10&$K_{10}$&$S_{10}$\\ 
\hline10&$P_{10}$&$S_5$\\
\hline\hline 11&$C_{11}$, $C_{11}^2$, $C_{11}^3$&$D_{11}$\\ 
\hline11&$K_{11}$&$S_{11}$\\ 
\hline\end{tabular}\end{center}

\smallskip

Afterwards, at $N=12$ and higher, the study becomes quite complicated, but we have results for several key series of graphs. We will be back to this, later in this book.

\section*{10e. Exercises}

We had a very exciting chapter here, with plenty of explicit graphs, nice pictures, and beautiful symmetry groups. As exercises on all this, we have:

\begin{exercise}
Draw all graphs with $N=9$ vertices, on a torus.
\end{exercise}

\begin{exercise}
Clarify which graphs at $N\leq11$ equal their complement.
\end{exercise}

\begin{exercise}
Enumerate all graphs on $N=12$ vertices, and good luck here.
\end{exercise}

\begin{exercise}
Try saying something nice, about the non-transitive graphs.
\end{exercise}

\begin{exercise}
Meditate on the magic matrices, and what can be done with them.
\end{exercise}

\begin{exercise}
Find more direct proofs for all our general product results.
\end{exercise}

\begin{exercise}
Work out all the concrete applications of our product results.
\end{exercise}

\begin{exercise}
Learn more about the Kneser graphs, and their various properties.
\end{exercise}

As bonus exercise, further meditate on the Petersen graph, and what can be done with it. This is indeed the typical counterexample to all sorts of questions, in graph theory.

\chapter{Spectral theory}

\section*{11a. Function algebras}

It is possible to do many other things, based on the techniques developed above, by further building on that material. But, my proposal now would be to slow down, temporarily forget about the graphs $X$, and focus instead on the finite groups $G$. These are our main tools, and before going to war, you have to sharpen your blades.

\bigskip

In order to discuss this, let us first have a second look at the magic unitaries, introduced in chapter 10. We have the following result, summarizing our knowledge from there:

\begin{theorem}
Given a subgroup $G\subset S_N\subset O_N$, the standard coordinates $u_{ij}\in C(G)$ generate $C(G)$, are given by the following formula, and form a magic matrix: 
$$u_{ij}=\chi\left(\sigma\in G\Big|\sigma(j)=i\right)$$
These coordinates appear as well as the coefficients of the transpose of the action map $a:X\times G\to X$ on the set $X=\{1,\ldots,N\}$, given by $(i,\sigma)\to\sigma(i)$, which is given by:
$$\Phi(e_i)=\sum_je_j\otimes u_{ji}$$
Also, in the case where we have a graph with $N$ vertices, the action of $G$ on the vertex set $X$ leaves invariant the edges precisely when $du=ud$.
\end{theorem}

\begin{proof}
This is something that we know from chapter 10, the idea being as follows:

\medskip

(1) Since the action of the group elements $\sigma\in G\subset S_N\subset O_N$ on the standard basis of $\mathbb R^N$ is given by $\sigma(e_i)=e_{\sigma(i)}$, we have $\sigma_{ij}=\delta_{\sigma(j)i}$, which gives the formula of $u_{ij}$.

\medskip

(2) The fact that our matrix $u=(u_{ij})$ is indeed magic, in the sense that its entries sum up to 1, on each row and each column, follows from our formula of $u_{ij}$. 

\medskip

(3) By Stone-Weierstrass we get $C(G)=<u_{ij}>$, since the coordinate functions $u_{ij}$ separate the points of $G$, in the sense that $\sigma\neq\pi$ needs $\sigma_{ij}\neq\pi_{ij}$, for some $i,j$.

\medskip

(4) Regarding the action $a:X\times G\to X$ and the coaction $\Phi:C(X)\to C(X)\otimes C(G)$, all this looks scary, but is in fact a triviality, as explained in chapter 10.

\medskip

(5) Finally, in what regards the last assertion, concerning $du=ud$, this again looks a bit abstract and scary, but is again a triviality, as explained in chapter 10.
\end{proof}

We have the following result, further building on the above:

\begin{theorem}
The symmetry group $G(X)$ of a graph $X$ having $N$ vertices is given by the following formula, at the level of the corresponding algebra of functions,
$$C(G(X))=C(S_N)\Big/\Big<du=ud\Big>$$
with $d\in M_N(0,1)$ being as usual the adjacency matrix of $X$.
\end{theorem}

\begin{proof}
This follows indeed from Theorem 11.1, and more specifically, is just an abstract reformulation of the last assertion there.
\end{proof}

In order to further build on all this, the idea will be that of getting rid of $S_N$, or rather of the corresponding algebra of functions $C(S_N)$, and formulating everything in terms of magic matrices. To be more precise, leaving aside $X$, we have the following question:

\begin{question}
With a suitable algebra formalism, do we have
$$C(S_N)=A\left((u_{ij})_{i,j=1,\ldots,N}\Big|u={\rm magic}\right)$$
with $A$ standing for ``universal algebra generated by"?
\end{question}

At the first glance, this question might seem overly theoretical and abstract, but the claim is that such things can be useful, in order to deal with graphs. Indeed, assuming a positive answer to this question, by Theorem 11.2 we would conclude that $C(G(X))$ is the universal algebra generated by the entries of a magic matrix commuting with $d$. 

\bigskip

Which is something very nice, and potentially useful, among others putting under a more conceptual light the various product computations from the previous chapter, done with magic matrices. In a word, we have seen in the previous chapter that the magic matrices can be very useful objects, so let us go now for it, and reformulate everything in terms of them, along the lines of Question 11.3, and of the comments afterwards.

\bigskip

Getting to work now, the algebra $C(S_N)$ from Question 11.3, that we want to axiomatize via magic matrices, is something very particular, and we have:

\begin{fact}
The function algebra $C(S_N)$ has the following properties:
\begin{enumerate}
\item It is a complex algebra.

\item It is finite dimensional.

\item It is commutative, $fg=gf$.

\item It has an involution, given by $f^*(x)=\overline{f(x)}$.

\item It has a norm, given by $||f||=\sup_x|f(x)|$.
\end{enumerate}
\end{fact}

So, which of these properties shall we choose, for our axiomatization? Definitely (1), and then common sense would suggest to try the combination (2+3). However, since we will be soon in need, in this book, of algebras which are infinite dimensional, or not commutative, or both, let us go instead with the combination (4+5), for most of our axiomatization work, and then with (2) or (3) added, if needed, towards the end. 

\bigskip

In short, trust me here, we need to do some algebra and this is what we will do, we will learn useful things in what follows, and here are some axioms, to start with:

\begin{definition}
A $C^*$-algebra is a complex algebra $A$, given with an involution $a\to a^*$ and a norm $a\to||a||$, such that:
\begin{enumerate}
\item The norm and involution are related by $||aa^*||=||a||^2$.

\item $A$ is complete, as metric space, with respect to the norm.
\end{enumerate}
\end{definition}

As a first basic class of examples, which are of interest for us, in relation with Question 11.3, we have the algebras of type $A=C(X)$, with $X$ being a finite, or more generally compact space, with the usual involution and norm of functions, namely:
$$f^*(x)=\overline{f(x)}\quad,\quad ||f||=\sup_x|f(x)|$$

Observe that such algebras are commutative, $fg=gf$, and also that both the conditions in Definition 11.5 are satisfied, with (1) being something trivial, and with (2) coming from the well-known fact that a uniform limit of continuous function is continuous.

\bigskip

Interestingly, and of guaranteed interest for many considerations to follow, in this book, as a second basic class of examples we have the matrix algebras $A=M_N(\mathbb C)$, with the usual involution and norm of the complex matrices, namely:
$$(M^*)_{ij}=\overline{M}_{ij}\quad,\quad ||M||=\sup_{||x||=1}||Mx||$$

Observe that such algebras are finite dimensional, and also that the two conditions in Definition 11.5 are satisfied, with (1) being a good linear algebra exercise for you, via double inequality, and with (2) being trivial, our algebra being finite dimensional.

\bigskip

Summarizing, good definition that we have, so let us develop now some theory, for the $C^*$-algebras. Inspired by the matrix algebra examples, we first have:

\begin{theorem}
Given an element $a\in A$ of a $C^*$-algebra, define its spectrum as:
$$\sigma(a)=\left\{\lambda\in\mathbb C\Big|a-\lambda\notin A^{-1}\right\}$$
The following spectral theory results hold, exactly as in the $A=M_N(\mathbb C)$ case:
\begin{enumerate}
\item We have $\sigma(ab)\cup\{0\}=\sigma(ba)\cup\{0\}$.

\item We have $\sigma(f(a))=f(\sigma(a))$, for any $f\in\mathbb C(X)$ having poles outside $\sigma(a)$.

\item The spectrum $\sigma(a)$ is compact, non-empty, and contained in $D_0(||a||)$.

\item The spectra of unitaries $(u^*=u^{-1})$ and self-adjoints $(a=a^*)$ are on $\mathbb T,\mathbb R$.

\item The spectral radius of normal elements $(aa^*=a^*a)$ is given by $\rho(a)=||a||$.
\end{enumerate}
In addition, assuming $a\in A\subset B$, the spectra of $a$ with respect to $A$ and to $B$ coincide.
\end{theorem}

\begin{proof}
Here the assertions (1-5), which are of course formulated a bit informally, are well-known for the matrix algebra $A=M_N(\mathbb C)$, and the proof in general is similar:

\medskip

(1) Assuming that $1-ab$ is invertible, with inverse $c$, we have $abc=cab=c-1$, and it follows that $1-ba$ is invertible too, with inverse $1+bca$. Thus $\sigma(ab),\sigma(ba)$ agree on $1\in\mathbb C$, and by linearity, it follows that $\sigma(ab),\sigma(ba)$ agree on any point $\lambda\in\mathbb C^*$.

\medskip

(2) The formula $\sigma(f(a))=f(\sigma(a))$ is clear for polynomials, $f\in\mathbb C[X]$, by factorizing $f-\lambda$, with $\lambda\in\mathbb C$. Then, the extension to the rational functions is straightforward, because $P(a)/Q(a)-\lambda$ is invertible precisely when $P(a)-\lambda Q(a)$ is.

\medskip

(3) By using $1/(1-b)=1+b+b^2+\ldots$ for $||b||<1$ we obtain that $a-\lambda$ is invertible for $|\lambda|>||a||$, and so $\sigma(a)\subset D_0(||a||)$. It is also clear that $\sigma(a)$ is closed, so what we have is a compact set. Finally, assuming $\sigma(a)=\emptyset$ the function $f(\lambda)=\varphi((a-\lambda)^{-1})$ is well-defined, for any $\varphi\in A^*$, and by Liouville we get $f=0$, contradiction.

\medskip

(4) Assuming $u^*=u^{-1}$ we have $||u||=1$, and so $\sigma(u)\subset D_0(1)$. But with $f(z)=z^{-1}$ we obtain via (2) that we have as well $\sigma(u)\subset f(D_0(1))$, and this gives $\sigma(u)\subset\mathbb T$. As for the result regarding the self-adjoints, this can be obtained from the result for the unitaries, by using (2) with functions of type $f(z)=(z+it)/(z-it)$, with $t\in\mathbb R$.

\medskip

(5) It is routine to check, by integrating quantities of type $z^n/(z-a)$ over circles centered at the origin, and estimating, that the spectral radius is given by $\rho(a)=\lim||a^n||^{1/n}$. But in the self-adjoint case, $a=a^*$, this gives $\rho(a)=||a||$, by using exponents of type $n=2^k$, and then the extension to the general normal case is straightforward.

\medskip 

(6) Regarding now the last assertion, the inclusion $\sigma_B(a)\subset\sigma_A(a)$ is clear. For the converse, assume $a-\lambda\in B^{-1}$, and set $b=(a-\lambda )^*(a-\lambda )$. We have then:
$$\sigma_A(b)-\sigma_B(b)=\left\{\mu\in\mathbb C-\sigma_B(b)\Big|(b-\mu)^{-1}\in B-A\right\}$$

Thus this difference in an open subset of $\mathbb C$. On the other hand $b$ being self-adjoint, its two spectra are both real, and so is their difference. Thus the two spectra of $b$ are equal, and in particular $b$ is invertible in $A$, and so $a-\lambda\in A^{-1}$, as desired.
\end{proof}

With these ingredients, we can now a prove a key result, as follows:

\begin{theorem}[Gelfand]
Any commutative $C^*$-algebra is of the form $A=C(X)$, with 
$$X=\Big\{\chi:A\to\mathbb C\,,\ {\rm normed\ algebra\ character}\Big\}$$
with topology making continuous the evaluation maps $ev_a:\chi\to\chi(a)$.
\end{theorem}

\begin{proof}
This is something quite tricky, the idea being as follows:

\medskip

(1) Given a commutative $C^*$-algebra $A$, let us define a space $X$ as in the statement. Then $X$ is compact, and $a\to ev_a$ is a morphism of algebras, as follows:
$$ev:A\to C(X)$$

(2) We first prove that $ev$ is involutive. We use the following formula, which is similar to the $z=Re(z)+iIm(z)$ decomposition formula for usual complex numbers:
$$a=\frac{a+a^*}{2}+i\cdot\frac{a-a^*}{2i}$$

Thus it is enough to prove $ev_{a^*}=ev_a^*$ for the self-adjoint elements $a$. But this is the same as proving that $a=a^*$ implies that $ev_a$ is a real function, which is in turn true, by Theorem 11.6 (4), because $ev_a(\chi)=\chi(a)$ is an element of $\sigma(a)$, contained in $\mathbb R$.

\medskip

(3) Since $A$ is commutative, each element is normal, so $ev$ is isometric:
$$||ev_a||
=\rho(a)
=||a||$$

It remains to prove that $ev$ is surjective. But this follows from the Stone-Weierstrass theorem, because $ev(A)$ is a closed subalgebra of $C(X)$, which separates the points.
\end{proof}

The above result is something truly remarkable, and we can now formulate:

\begin{definition}
Given an arbitrary $C^*$-algebra $A$, we can write it as
$$A=C(X)$$
with $X$ compact quantum space. When $A$ is commutative, $X$ is a usual compact space.
\end{definition}

Which is very nice, but obviously, a bit off-topic. More on quantum spaces a bit later in this book, and for the moment, just take this as it came, namely math professor supposed to lecture on something, and ending up in lecturing on something else.

\bigskip

More seriously now, the Gelfand theorem that we just learned is something very useful, and getting back to our original Question 11.3, we can answer it, as follows:

\begin{theorem}
The algebra of functions on $S_N$ has the following presentation,
$$C(S_N)=C^*_{comm}\left((u_{ij})_{i,j=1,\ldots,N}\Big|u={\rm magic}\right)$$
and the multiplication, unit and inversion map of $S_N$ appear from the maps
$$\Delta(u_{ij})=\sum_ku_{ik}\otimes u_{kj}\quad,\quad 
\varepsilon(u_{ij})=\delta_{ij}\quad,\quad 
S(u_{ij})=u_{ji}$$
defined at the algebraic level, of functions on $S_N$, by transposing.
\end{theorem}

\begin{proof}
The universal algebra $A$ in the statement being commutative, by the Gelfand theorem it must be of the form $A=C(X)$, with $X$ being a certain compact space. Now since we have coordinates $u_{ij}:X\to\mathbb R$, we have an embedding $X\subset M_N(\mathbb R)$. Also, since we know that these coordinates form a magic matrix, the elements $g\in X$ must be 0-1 matrices, having exactly one 1 entry on each row and each column, and so $X=S_N$. Thus we have proved the first assertion, and the second assertion is clear as well.
\end{proof}

In relation now with graphs, we have the following result:

\begin{theorem}
The symmetry group of a graph $X$ having $N$ vertices is given by the following formula, at the level of the corresponding algebra of functions,
$$C(G(X))=C^*_{comm}\left((u_{ij})_{i,j=1,\ldots,N}\Big|u={\rm magic},\ du=ud\right)$$
with $d\in M_N(0,1)$ being as usual the adjacency matrix of $X$, and the multiplication, unit and inversion map of $G(X)$ appear from the maps
$$\Delta(u_{ij})=\sum_ku_{ik}\otimes u_{kj}\quad,\quad 
\varepsilon(u_{ij})=\delta_{ij}\quad,\quad 
S(u_{ij})=u_{ji}$$
defined at the algebraic level, of functions on $G(X)$, by transposing.
\end{theorem}

\begin{proof}
This follows indeed from Theorem 11.9, by combining it with Theorem 11.2, which tells us that when dealing with a graph $X$, with adjacency matrix $d\in M_N(0,1)$, we simply must add the condition $du=ud$, on the corresponding magic matrix $u$.
\end{proof}

Let us discuss as well what happens in relation with the partial permutations, introduced at the end of chapter 9. We first have the following result:

\begin{theorem}
The algebra of functions on $\widetilde{S}_N$ has the following presentation, with submagic meaning formed of projections, pairwise orthogonal on rows and columns,
$$C(\widetilde{S}_N)=C^*_{comm}\left((u_{ij})_{i,j=1,\ldots,N}\Big|u={\rm submagic}\right)$$
and the multiplication and unit of $\widetilde{S}_N$ appear from the maps
$$\Delta(u_{ij})=\sum_ku_{ik}\otimes u_{kj}\quad,\quad 
\varepsilon(u_{ij})=\delta_{ij}$$
defined at the algebraic level, of functions on $\widetilde{S}_N$, by transposing.
\end{theorem}

\begin{proof}
This is very similar to the proof of Theorem 11.9, with the result again coming from the Gelfand theorem, applied to the universal algebra in the statement.
\end{proof}

In relation now with graphs, the result, which is quite tricky, is as follows:

\begin{theorem}
Given a graph $X$ having $N$ vertices, and adjacency matrix $d\in M_N(0,1)$, consider its partial automorphism semigroup, given by:
$$\widetilde{G}(X)=\left\{\sigma\in\widetilde{S}_N\Big|d_{ij}=d_{\sigma(i)\sigma(j)},\ \forall i,j\in Dom(\sigma)\right\}$$
We have then the following formula, with $R=diag(R_i)$, $C=diag(C_j)$, with $R_i,C_j$ being the row and column sums of the associated submagic matrix $u$:
$$C(\widetilde{G}(X))=C(\widetilde{S}_N)\Big/\Big<R(du-ud)C=0\Big>$$
Moreover, when using the relation $du=ud$ instead of the above one, we obtain a certain semigroup $\bar{G}(X)\subset\widetilde{G}(X)$, which can be strictly smaller.
\end{theorem}

\begin{proof}
This requires a bit of abstract thinking, the idea being as follows:

\medskip

(1) To start with, we will use the formula from Theorem 11.11, namely:
$$C(\widetilde{S}_N)=C^*_{comm}\left((u_{ij})_{i,j=1,\ldots,N}\Big|u={\rm submagic}\right)$$

(2) Getting now to graphs, the definition of $\widetilde{G}(X)$ in the statement reformulates as follows, in terms of the usual adjacency relation $i-j$ for the vertices:
$$\widetilde{G}(X)=\left\{\sigma\in\widetilde{S}_N\Big|i-j,\exists\,\sigma(i),\exists\,\sigma(j)\implies \sigma(i)-\sigma(j)\right\}$$

Indeed, this reformulation is something which is clear from definitions.

\medskip

(3) In view of this, we have the following product computation:
\begin{eqnarray*}
(du)_{ij}(\sigma)
&=&\sum_kd_{ik}u_{kj}(\sigma)\\
&=&\sum_{k\sim i}u_{kj}(\sigma)\\
&=&\begin{cases}
1&{\rm if}\ \sigma(j)-i\\
0&{\rm otherwise}
\end{cases}
\end{eqnarray*}

On the other hand, we have as well the following computation:
\begin{eqnarray*}
(ud)_{ij}(\sigma)
&=&\sum_ku_{ik}d_{kj}(\sigma)\\
&=&\sum_{k\sim j}u_{ik}(\sigma)\\
&=&\begin{cases}
1&{\rm if}\ \sigma^{-1}(i)-j\\
0&{\rm otherwise}
\end{cases}
\end{eqnarray*}

To be more precise, in the above two formulae the ``otherwise'' cases include by definition the cases where $\sigma(j)$, respectively $\sigma^{-1}(i)$, is undefined. 

\medskip

(4) On the other hand, we have as well the following formulae:
$$R_i(\sigma)=\sum_ju_{ij}(\sigma)
=\begin{cases}
1&{\rm if}\ \exists\,\sigma^{-1}(i)\\
0&{\rm otherwise}
\end{cases}$$
$$C_j(\sigma)=\sum_iu_{ij}(\sigma)
=\begin{cases}
1&{\rm if}\ \exists\,\sigma(j)\\
0&{\rm otherwise}
\end{cases}$$

(5) Now by multiplying the above formulae, we obtain the following formulae:
$$(R_i(du)_{ij}C_j)(\sigma)
=\begin{cases}
1&{\rm if}\ \sigma(j)-i\ {\rm and}\ \exists\,\sigma^{-1}(i)\ {\rm and}\ \exists\,\sigma(j)\\
0&{\rm otherwise}
\end{cases}$$
$$(R_i(ud)_{ij}C_j)(\sigma)
=\begin{cases}
1&{\rm if}\ \sigma^{-1}(i)-j\ {\rm and}\ \exists\,\sigma^{-1}(i)\ {\rm and}\ \exists\,\sigma(j)\\
0&{\rm otherwise}
\end{cases}$$

(6) We conclude that the relations in the statement, $R_i(du)_{ij}C_j=R_i(ud)_{ij}C_j$, when applied to a given element $\sigma\in\widetilde{S}_N$, correspond to the following condition:
$$\exists\,\sigma^{-1}(i),\ \exists\,\sigma(j)\implies[\sigma(j)-i\iff\sigma^{-1}(i)-j]$$

But with $i=\sigma(k)$, this latter condition reformulates as follows:
$$\exists\,\sigma(k),\ \exists\,\sigma(j)\implies[\sigma(j)-\sigma(k)\iff k-j]$$

Thus we must have $\sigma\in\widetilde{G}(X)$, and we obtain the presentation result for $\widetilde{G}(X)$. 

\medskip

(7) Regarding now the second assertion, the simplest counterexample here is simplex $X_N$, having $N$ vertices and edges everywhere. Indeed, the adjacency matrix of this simplex is $d=\mathbb I_N-1_N$, with $\mathbb I_N$ being the all-1 matrix, and so the commutation of this matrix with $u$ corresponds to the fact that $u$ must be bistochastic. Thus, $u$ must be in fact magic, and we obtain $\bar{G}(X_N)=S_N$, which is smaller than $\widetilde{G}(X_N)=\widetilde{S}_N$.
\end{proof}

Many interesting things can be said here, and we have of course many explicit examples, for the graphs having small number of vertices. We will be back to this.

\section*{11b. Representation theory}

Let us discuss now another advanced algebraic topic, namely the Peter-Weyl theory for the finite groups, and in particular for the permutation groups. The idea here will be that there are several non-trivial things that can be said about the group actions on graphs, $G\curvearrowright X$, by using the Peter-Weyl theory for finite groups, according to:

\begin{principle} 
Any finite group action on a finite graph $G\curvearrowright X$, with $|X|=N$, produces a unitary representation of $G$, obtained as
$$G\to S_N\subset U_N$$
that we can decompose and study by using the Peter-Weyl theory for $G$. And with this study leading to non-trivial results about the action $G\curvearrowright X$, and about $X$ itself. 
\end{principle}

Getting started now, with this program, we first have to forget about the finite graphs $X$, and develop the Peter-Weyl theory for the finite groups $G$. We first have:

\index{representation}
\index{character}
\index{finite group}
\index{character of representation}
\index{trace of representation}

\begin{definition}
A representation of a finite group $G$ is a morphism as follows:
$$u:G\to U_N$$
The character of such a representation is the function $\chi:G\to\mathbb C$ given by
$$g\to Tr(u_g)$$
where $Tr$ is the usual trace of the $N\times N$ matrices, $Tr(M)=\sum_iM_{ii}$.
\end{definition}

As a basic example here, for any finite group we always have available the trivial 1-dimensional representation, which is by definition as follows:
$$u:G\to U_1\quad,\quad 
g\to(1)$$

At the level of non-trivial examples now, most of the groups that we met so far, in chapter 5, naturally appear as subgroups $G\subset U_N$. In this case, the embedding $G\subset U_N$ is of course a representation, called fundamental representation:
$$u:G\subset U_N\quad,\quad 
g\to g$$

In this situation, there are many other representations of $G$, which are equally interesting. For instance, we can define the representation conjugate to $u$, as being:
$$\bar{u}:G\subset U_N\quad,\quad
g\to\bar{g}$$

In order to clarify all this, and see which representations are available, let us first discuss the various operations on the representations. The result here is as follows:

\index{sum of representations}
\index{product of representations}
\index{conjugate representation}
\index{spinned representation}

\begin{proposition}
The representations of a finite group $G$ are subject to:
\begin{enumerate}
\item Making sums. Given representations $u,v$, having dimensions $N,M$, their sum is the $N+M$-dimensional representation $u+v=diag(u,v)$.

\item Making products. Given representations $u,v$, having dimensions $N,M$, their tensor product is the $NM$-dimensional representation $(u\otimes v)_{ia,jb}=u_{ij}v_{ab}$.

\item Taking conjugates. Given a representation $u$, having dimension $N$, its complex conjugate is the $N$-dimensional representation $(\bar{u})_{ij}=\bar{u}_{ij}$.

\item Spinning by unitaries. Given a representation $u$, having dimension $N$, and a unitary $V\in U_N$, we can spin $u$ by this unitary, $u\to VuV^*$.
\end{enumerate}
\end{proposition}

\begin{proof}
The fact that the operations in the statement are indeed well-defined, among maps from $G$ to unitary groups, is indeed routine, and this gives the result.
\end{proof}

In relation now with characters, we have the following result:

\index{character}

\begin{proposition}
We have the following formulae, regarding characters
$$\chi_{u+v}=\chi_u+\chi_v\quad,\quad 
\chi_{u\otimes v}=\chi_u\chi_v\quad,\quad 
\chi_{\bar{u}}=\bar{\chi}_u\quad,\quad 
\chi_{VuV^*}=\chi_u$$
in relation with the basic operations for the representations.
\end{proposition}

\begin{proof}
All these assertions are elementary, by using the following well-known trace formulae, valid for any two square matrices $g,h$, and any unitary $V$:
$$Tr(diag(g,h))=Tr(g)+Tr(h)\quad,\quad 
Tr(g\otimes h)=Tr(g)Tr(h)$$
$$Tr(\bar{g})=\overline{Tr(g)}\quad,\quad
Tr(VgV^*)=Tr(g)$$

Thus, we are led to the conclusions in the statement.
\end{proof}

Assume now that we are given a finite group $G\subset U_N$. By using the above operations, we can construct a whole family of representations of $G$, as follows:

\index{Peter-Weyl representations}
\index{colored integer}

\begin{definition}
Given a finite group $G\subset U_N$, its Peter-Weyl representations are the tensor products between the fundamental representation and its conjugate:
$$u:G\subset U_N\quad,\quad 
\bar{u}:G\subset U_N$$ 
We denote these tensor products $u^{\otimes k}$, with $k=\circ\bullet\bullet\circ\ldots$ being a colored integer, with the colored tensor powers being defined according to the rules 
$$u^{\otimes\circ}=u\quad,\quad 
u^{\otimes\bullet}=\bar{u}\quad,\quad
u^{\otimes kl}=u^{\otimes k}\otimes u^{\otimes l}$$
and with the convention that $u^{\otimes\emptyset}$ is the trivial representation $1:G\to U_1$.
\end{definition}

Here are a few examples of such Peter-Weyl representations, namely those coming from the colored integers of length 2, to be often used in what follows:
$$u^{\otimes\circ\circ}=u\otimes u\quad,\quad 
u^{\otimes\circ\bullet}=u\otimes\bar{u}$$
$$u^{\otimes\bullet\circ}=\bar{u}\otimes u\quad,\quad 
u^{\otimes\bullet\bullet}=\bar{u}\otimes\bar{u}$$

Observe also that the characters of Peter-Weyl representations are given by the following formula, with the powers $\chi$ being given by
$\chi^\circ=\chi$, $\chi^\bullet=\bar{\chi}$ and multiplicativity:
$$\chi_{u^{\otimes k}}=(\chi_u)^k$$

In order now to advance, let us formulate the following key definition:

\index{Hom space}
\index{End space}
\index{Fix space}

\begin{definition}
Given a finite group $G$, and two of its representations,
$$u:G\to U_N\quad,\quad 
v:G\to U_M$$
we define the linear space of intertwiners between these representations as being 
$$Hom(u,v)=\left\{T\in M_{M\times N}(\mathbb C)\Big|Tu_g=v_gT,\forall g\in G\right\}$$
and we use the following conventions:
\begin{enumerate}
\item We use the notations $Fix(u)=Hom(1,u)$, and $End(u)=Hom(u,u)$.

\item We write $u\sim v$ when $Hom(u,v)$ contains an invertible element.

\item We say that $u$ is irreducible, and write $u\in Irr(G)$, when $End(u)=\mathbb C1$.
\end{enumerate}
\end{definition}

The terminology here is very standard, with Hom and End standing for ``homomorphisms'' and ``endomorphisms'', and with Fix standing for ``fixed points''. We have:

\index{tensor category}

\begin{theorem}
The following happen:
\begin{enumerate}
\item The intertwiners are stable under composition:
$$T\in Hom(u,v)\ ,\ 
S\in Hom(v,w)
\implies ST\in Hom(u,w)$$

\item The intertwiners are stable under taking tensor products:
$$S\in Hom(u,v)\ ,\ 
T\in Hom(w,t)\\
\implies S\otimes T\in Hom(u\otimes w,v\otimes t)$$

\item The intertwiners are stable under taking adjoints:
$$T\in Hom(u,v)
\implies T^*\in Hom(v,u)$$

\item Thus, the Hom spaces form a tensor $*$-category.
\end{enumerate}
\end{theorem}

\begin{proof}
All this is clear from definitions, the verifications being as follows:

\medskip

(1) This follows indeed from the following computation, valid for any $g\in G$:
$$STu_g=Sv_gT=w_gST$$

(2) Again, this is clear, because we have the following computation:
\begin{eqnarray*}
(S\otimes T)(u_g\otimes w_g)
&=&Su_g\otimes Tw_g\\
&=&v_gS\otimes t_gT\\
&=&(v_g\otimes t_g)(S\otimes T)
\end{eqnarray*}

(3) This follows from the following computation, valid for any $g\in G$:
\begin{eqnarray*}
Tu_g=v_gT
&\implies&u_g^*T^*=T^*v_g^*\\
&\implies&T^*v_g=u_gT^*
\end{eqnarray*}

(4) This is just an abstract conclusion of (1,2,3), with a tensor $*$-category being by definition an abstract beast satisfying these conditions (1,2,3). We will be back to tensor categories later on in this book, with more details on all this.
\end{proof}

As a main consequence of Theorem 11.19, we have:

\begin{theorem}
Given a representation $u:G\to U_N$, the linear space
$$End(u)\subset M_N(\mathbb C)$$
is a $*$-algebra, with respect to the usual involution of the matrices.
\end{theorem}

\begin{proof}
We know from Theorem 11.19 (1) that $End(u)$ is a subalgebra of $M_N(\mathbb C)$, and we know as well from Theorem 11.19 (3) that this subalgebra is stable under the involution $*$. Thus, what we have here is a $*$-subalgebra of $M_N(\mathbb C)$, as claimed.
\end{proof}

\bigskip

Our claim now is that Theorem 11.20 gives us everything that we need, in order to have some advanced representation theory started, for our finite groups $G$. Indeed, we can combine this result with the following standard fact, from matrix algebra:

\index{operator algebra}
\index{finite dimensional algebra}

\begin{theorem}
Let $A\subset M_N(\mathbb C)$ be a $*$-algebra.
\begin{enumerate}
\item We can write $1=p_1+\ldots+p_k$, with $p_i\in A$ being central minimal projections.

\item The linear spaces $A_i=p_iAp_i$ are non-unital $*$-subalgebras of $A$.

\item We have a non-unital $*$-algebra sum decomposition $A=A_1\oplus\ldots\oplus A_k$.

\item We have unital $*$-algebra isomorphisms $A_i\simeq M_{n_i}(\mathbb C)$, with $n_i=rank(p_i)$.

\item Thus, we have a $*$-algebra isomorphism $A\simeq M_{n_1}(\mathbb C)\oplus\ldots\oplus M_{n_k}(\mathbb C)$.
\end{enumerate}
\end{theorem}

\begin{proof}
This is something standard, whose proof is however quite long, as follows:

\medskip

(1) Consider an arbitrary $*$-algebra of the $N\times N$ matrices, $A\subset M_N(\mathbb C)$, as in the statement. Let us first look at the center of this algebra, $Z(A)=A\cap A'$. It is elementary to prove that this center, as an algebra, is of the following form:
$$Z(A)\simeq\mathbb C^k$$

Consider now the standard basis $e_1,\ldots,e_k\in\mathbb C^k$, and let  $p_1,\ldots,p_k\in Z(A)$ be the images of these vectors via the above identification. In other words, these elements $p_1,\ldots,p_k\in A$ are central minimal projections, summing up to 1:
$$p_1+\ldots+p_k=1$$

The idea is then that this partition of the unity will eventually lead to the block decomposition of $A$, as in the statement. We prove this in 4 steps, as follows:

\medskip

(2) We first construct the matrix blocks, our claim here being that each of the following linear subspaces of $A$ are non-unital $*$-subalgebras of $A$:
$$A_i=p_iAp_i$$

But this is clear, with the fact that each $A_i$ is closed under the various non-unital $*$-subalgebra operations coming from the projection equations $p_i^2=p_i=p_i^*$.

\medskip

(3) We prove now that the above algebras $A_i\subset A$ are in a direct sum position, in the sense that we have a non-unital $*$-algebra sum decomposition, as follows:
$$A=A_1\oplus\ldots\oplus A_k$$

As with any direct sum question, we have two things to be proved here. First, by using the formula $p_1+\ldots+p_k=1$ and the projection equations $p_i^2=p_i=p_i^*$, we conclude that we have the needed generation property, namely:
$$A_1+\ldots+ A_k=A$$

As for the fact that the sum is indeed direct, this follows as well from the formula $p_1+\ldots+p_k=1$, and from the projection equations $p_i^2=p_i=p_i^*$.

\medskip

(4) Our claim now, which will finish the proof, is that each of the $*$-subalgebras $A_i=p_iAp_i$ constructed above is a full matrix algebra. To be more precise here, with $n_i=rank(p_i)$, our claim is that we have isomorphisms, as follows:
$$A_i\simeq M_{n_i}(\mathbb C)$$

In order to prove this claim, recall that the projections $p_i\in A$ were chosen central and minimal. Thus, the center of each of the algebras $A_i$ reduces to the scalars:
$$Z(A_i)=\mathbb C$$

But this shows, either via a direct computation, or via the bicommutant theorem, that the each of the algebras $A_i$ is a full matrix algebra, as claimed.

\medskip

(5) We can now obtain the result, by putting together what we have. Indeed, by using the results from (3) and (4), we obtain an isomorphism as follows:
$$A\simeq M_{n_1}(\mathbb C)\oplus\ldots\oplus M_{n_k}(\mathbb C)$$

In addition to this, a careful look at the isomorphisms established in (4) shows that at the global level, that of the algebra $A$ itself, the above isomorphism simply comes by twisting the following standard multimatrix embedding, discussed in the beginning of the proof, step (2) above, by a certain unitary matrix $U\in U_N$:
$$M_{n_1}(\mathbb C)\oplus\ldots\oplus M_{n_k}(\mathbb C)\subset M_N(\mathbb C)$$

Now by putting everything together, we obtain the result.
\end{proof}

We can now formulate our first Peter-Weyl theorem, as follows:

\index{Peter-Weyl}

\begin{theorem}[PW1]
Let $u:G\to U_N$ be a representation, consider the algebra $A=End(u)$, and write its unit as above, with $p_i$ being central minimal projections:
$$1=p_1+\ldots+p_k$$
The representation $u$ decomposes then as a direct sum, as follows,
$$u=u_1+\ldots+u_k$$
with each $u_i$ being an irreducible representation, obtained by restricting $u$ to $Im(p_i)$.
\end{theorem}

\begin{proof}
This follows from Theorem 11.20 and Theorem 11.21, as follows:

\medskip

(1) As a first observation, by replacing $G$ with its image $u(G)\subset U_N$, we can assume if we want that our representation $u$ is faithful, $G\subset_uU_N$. However, this replacement will not be really needed, and we will keep using $u:G\to U_N$, as above.

\medskip

(2) In order to prove the result, we will need some preliminaries. We first associate to our representation $u:G\to U_N$ the corresponding action map on $\mathbb C^N$. If a linear subspace $V\subset\mathbb C^N$ is invariant, the restriction of the action map to $V$ is an action map too, which must come from a subrepresentation $v\subset u$. This is clear indeed from definitions, and with the remark that the unitaries, being isometries, restrict indeed into unitaries.

\medskip

(3) Consider now a projection $p\in End(u)$. From $pu=up$ we obtain that the linear space $V=Im(p)$ is invariant under $u$, and so this space must come from a subrepresentation $v\subset u$. It is routine to check that the operation $p\to v$ maps subprojections to subrepresentations, and minimal projections to irreducible representations.

\medskip

(4) With these preliminaries in hand, let us decompose the algebra $End(u)$ as in Theorem 11.21, by using the decomposition $1=p_1+\ldots+p_k$ into minimal projections. If we denote by $u_i\subset u$ the subrepresentation coming from the vector space $V_i=Im(p_i)$, then we obtain in this way a decomposition $u=u_1+\ldots+u_k$, as in the statement.
\end{proof}

Here is now our second Peter-Weyl theorem, complementing Theorem 11.22:

\index{Peter-Weyl}

\begin{theorem}[PW2]
Given a subgroup $G\subset_uU_N$, any irreducible representation 
$$v:G\to U_M$$
appears inside a tensor product of the fundamental representation $u$ and its adjoint $\bar{u}$.
\end{theorem}

\begin{proof}
We define the space of coefficients a representation $v:G\to U_M$ to be:
$$C_v=span\Big[g\to(v_g)_{ij}\Big]$$

The construction $v\to C_v$ is then functorial, in the sense that it maps subrepresentations into linear subspaces. Also, we have an inclusion of linear spaces as follows:
$$C_v\subset<g_{ij}>$$

On the other hand, by definition of the Peter-Weyl representations, we have:
$$<g_{ij}>=\sum_kC_{u^{\otimes k}}$$

Thus, we must have an inclusion as follows, for certain exponents $k_1,\ldots,k_p$:
$$C_v\subset C_{u^{\otimes k_1}\oplus\ldots\oplus u^{\otimes k_p}}$$

We conclude that we have an inclusion of representations, as follows:
$$v\subset u^{\otimes k_1}\oplus\ldots\oplus u^{\otimes k_p}$$

Together with Theorem 11.22, this leads to the conclusion in the statement.
\end{proof}

As a third Peter-Weyl theorem, which is something more advanced, we have: 

\index{Frobenius isomorphism}
\index{Peter-Weyl}

\begin{theorem}[PW3]
We have a direct sum decomposition of linear spaces
$$C(G)=\bigoplus_{v\in Irr(G)}M_{\dim(v)}(\mathbb C)$$
with the summands being pairwise orthogonal with respect to the scalar product
$$<a,b>=\int_Gab^*$$
where $\int_G$ is the averaging over $G$.
\end{theorem}

\begin{proof}
This is something more tricky, the idea being as follows:

\medskip

(1) By combining the previous two Peter-Weyl results, from Theorem 11.22 and Theorem 11.23, we deduce that we have a linear space decomposition as follows:
$$C(G)
=\sum_{v\in Irr(G)}C_v
=\sum_{v\in Irr(G)}M_{\dim(v)}(\mathbb C)$$

Thus, in order to conclude, it is enough to prove that for any two irreducible corepresentations $v,w\in Irr(A)$, the corresponding spaces of coefficients are orthogonal:
$$v\not\sim w\implies C_v\perp C_w$$ 

(2) We will need the basic fact, whose proof is elementary, that for any representation $v$ we have the following formula, where $P$ is the orthogonal projection on $Fix(v)$:
$$\left(id\otimes\int_G\right)v=P$$

(3) We will also need the basic fact, whose proof is elementary too, that for any two representations $v,w$ we have an isomorphism as follows, called Frobenius isomorphism:
$$Hom(v,w)\simeq Fix(v\otimes\bar{w})$$

(4) Now back to our orthogonality question from (1), let us set indeed:
$$P_{ia,jb}=\int_Gv_{ij}w_{ab}^*$$

Then $P$ is the orthogonal projection onto the following vector space:
$$Fix(v\otimes\bar{w})
\simeq Hom(v,w)
=\{0\}$$

Thus we have $P=0$, and this gives the result.
\end{proof}

Finally, we have the following result, completing the Peter-Weyl theory:

\index{Peter-Weyl}
\index{central function}
\index{algebra of characters}

\begin{theorem}[PW4]
The characters of irreducible representations belong to
$$C(G)_{central}=\left\{f\in C(G)\Big|f(gh)=f(hg),\forall g,h\in G\right\}$$
called algebra of central functions on $G$, and form an orthonormal basis of it.
\end{theorem}

\begin{proof}
We have several things to be proved, the idea being as follows:

\medskip

(1) Observe first that $C(G)_{central}$ is indeed an algebra, which contains all the characters. Conversely, consider a function $f\in C(G)$, written as follows:
$$f=\sum_{v\in Irr(G)}f_v$$

The condition $f\in C(G)_{central}$ states then that for any $v\in Irr(G)$, we must have:
$$f_v\in C(G)_{central}$$

But this means precisely that the coefficient $f_v$ must be a scalar multiple of $\chi_v$, and so the characters form a basis of $C(G)_{central}$, as stated. 

\medskip

(2) The fact that we have an orthogonal basis follows from Theorem 11.24. 

\medskip

(3) As for the fact that the characters have norm 1, this follows from:
$$\int_G\chi_v\chi_v^*
=\sum_{ij}\int_Gv_{ii}v_{jj}^*
=\sum_i\frac{1}{N}
=1$$

Here we have used the fact that the above integrals $\int_Gv_{ij}v_{kl}^*$ form the orthogonal projection onto the following vector space: 
$$Fix(v\otimes\bar{v})
\simeq End(v)
=\mathbb C1$$

Thus, the proof of our theorem is now complete.
\end{proof}

So long for Peter-Weyl theory. As a comment, our approach here, which was rather functional analytic, was motivated by what we will be doing later in this book, in relation with quantum groups. For a more standard presentation of the Peter-Weyl theory for finite groups, there are many good books available, such as the book of Serre \cite{ser}.

\section*{11c. Transitive actions}

With the above understood, let us get back to graphs. As already mentioned in the beginning of the previous section, the general principle is that representation theory and the Peter-Weyl theorems can help in understanding the group actions on graphs, $G\curvearrowright X$, by regarding these actions as unitary representations of $G$, as follows:
$$u:G\to S_N\subset U_N$$

To be more precise, from a finite group and representation theory perspective, the whole graph symmetry problematics can be reformulated as follows:

\begin{problem}
Given a finite group representation $u:G\to U_N$:
\begin{enumerate}
\item Is this representation magic, in the sense that it factorizes through $S_N$?

\item If so, what are the $0-1$ graph adjacency matrices $d\in End(u)$?

\item Among these latter matrices, which ones generate $End(u)$?
\end{enumerate}
\end{problem}

Generally speaking, these questions are quite difficult, ultimately leading into planar algebras in the sense of Jones \cite{jo6}, and despite having learned in this chapter some serious $C^*$-algebra theory, and Peter-Weyl theory, we are not ready yet for such things. All this will have to wait for the end of this book, when we will discuss planar algebras.

\bigskip

In the meantime we can, however, discuss some elementary applications of Peter-Weyl theory to the study of graphs. As starting point, we have the following result:

\begin{theorem}
Given a transitive action $G\curvearrowright X$, any group character
$$\chi:G\to\mathbb T$$
canonically produces an eigenfunction and eigenvalue of $X$, given by
$$f(g(0))=\chi(g)\quad,\quad\lambda=\sum_{g(0)-0}\chi(g)$$
with $0\in X$ being a chosen vertex of the graph.
\end{theorem}

\begin{proof}
This is something elementary, coming from definitions, as follows:

\medskip

(1) Let us fix a vertex $0\in X$, and consider the following set:
$$S=\left\{g\in G\Big|g(0)-0\right\}$$

Observe that this set $S\subset G$ satisfies $1\notin S$, and also $S=S^{-1}$, due to:
\begin{eqnarray*}
g(0)-0
&\implies&g^{-1}(g(0))-g^{-1}(0)\\
&\implies&0-g^{-1}(0)\\
&\implies&g^{-1}(0)-0
\end{eqnarray*}

(2) As a comment here, the condition $1\notin S=S^{-1}$ is something quite familiar, namely the Cayley graph assumption from chapter 4. More about this in chapter 12 below, where we will systematically discuss the Cayley graphs, as a continuation of that material.

\medskip

(3) Now given a character $\chi:G\to\mathbb T$ as in the statement, we can construct a function on the vertex set of our graph, $f:X\to\mathbb T$, according to the following formula:
$$f(g(0))=\chi(g)$$

Observe that this function is indeed well-defined, everywhere on the graph, thanks to our assumption that the group action $G\curvearrowright X$ is transitive.

\medskip

(4) Our claim now is that this function $f:X\to\mathbb T$ is an eigenfunction of the adjacency matrix of the graph. Indeed, we have the following computation, for any $g\in G$:
\begin{eqnarray*}
(df)(g(0))
&=&\sum_{g(0)-j}f(j)\\
&=&\sum_{g(0)-h(0)}f(h(0))\\
&=&\sum_{g(0)-h(0)}\chi(h)\\
&=&\sum_{s\in S}\chi(sg)\\
&=&\sum_{s\in S}\chi(s)\chi(g)\\
&=&\sum_{s\in S}\chi(s)f(g(0))
\end{eqnarray*}

Thus, we are led to the conclusion in the statement.
\end{proof}

The above result is quite interesting, and as a continuation of the story, we have:

\begin{theorem}
In the context of the eigenfunctions and eigenvalues constructed above, coming from transitive actions $G\curvearrowright X$ and characters $\chi:G\to\mathbb T$:
\begin{enumerate}
\item For the trivial character, $\chi=1$, we obtain the trivial eigenfunction, $f=1$.

\item When $G$ is assumed to be abelian, this gives the diagonalization of $d$.
\end{enumerate}
Moreover, it is possible to generalize this construction, by using arbitrary irreducible representations instead of characters, and this gives again the diagonalization of $d$.
\end{theorem}

\begin{proof}
There are several things going on here, the idea being as follows:

\medskip

(1) This is clear from the construction from the proof of Theorem 11.27, which for trivial character, $\chi=1$, gives $f=1$, and $\lambda=|S|$, with the notations there.

\medskip

(2) For an abelian group $G$, the characters form the dual group $\widehat{G}\simeq G$, and so we get a collection of $|\widehat{G}|=|G|=|X|$ eigenfunctions, as needed for diagonalizing $d$.

\medskip

(3) In what regards the last assertion, the construction there is straightforward, with the final conclusion coming from the following formula, coming itself from Peter-Weyl:
$$\sum_{r\in Irr(G)}(\dim r)^2=|G|$$

We will leave clarifying the details here as an instructive exercise, and we will come back to this in chapter 12, when discussing more systematically the Cayley graphs.
\end{proof}

\section*{11d. Asymptotic aspects}

We would like to end this chapter with something still advanced, in relation with representations, but more analytic and refreshing, featuring some probability.

\bigskip

It is about formal graphs $X_N$ and their formal symmetry groups $G_N\subset S_N$ that we want to talk about, in the $N\to\infty$ limit, a bit as the physicists do. But, how to do this? Not clear at all, because while it is easy to talk about series of graphs $X=(X_N)$, in an intuitive way, the corresponding symmetry groups $G_N\subset S_N$ are not necessarily compatible with each other, as to form an axiomatizable object $G=(G_N)$.

\bigskip

Long story short, we are into potentially difficult mathematics here, and as a more concrete question that we can attempt to solve, we have:

\begin{problem}
What families of groups $G_N\subset S_N$ have compatible representation theory invariants, in the $N\to\infty$ limit? And then, can we use this in order to talk about ``good'' families of graphs $X_N$, whose symmetry groups $G_N=G(X_N)$ are of this type?
\end{problem}

But probably too much talking, let us get to work. The simplest graphs are undoubtedly the empty graphs and the simplices, with symmetry group $S_N$, so as a first question, we would like to know if the symmetric groups $S_N$ themselves have compatible representation theory invariants, with some nice $N\to\infty$ asymptotics to work out.

\bigskip

And here, surprise, or miracle, the answer is indeed yes, with the result, which is something very classical, remarkable, and beautiful, being as follows:

\begin{theorem}
The probability for a random $\sigma\in S_N$ to have no fixed points is
$$P\simeq\frac{1}{e}$$
in the $N\to\infty$ limit, where $e=2.7182\ldots$ is the usual constant from analysis. More generally, the main character of $S_N$, which counts these permutations, given by
$$\chi=\sum_i\sigma_{ii}$$
via the standard embedding $S_N\subset O_N$, follows the Poisson law $p_1$, in the $N\to\infty$ limit. Even more generally, the truncated characters of $S_N$, given by
$$\chi=\sum_{i=1}^{[tN]}\sigma_{ii}$$
with $t>0$, follow the Poisson laws $p_t$, in the $N\to\infty$ limit. 
\end{theorem}

\begin{proof}
Obviously, many things going on here. The idea is as follows:

\medskip

(1) In order to prove the first assertion, which is the key one, and probably the most puzzling one too, we will use the inclusion-exclusion principle. Let us set:
$$S_N^k=\left\{\sigma\in S_N\Big|\sigma(k)=k\right\}$$

The set of permutations having no fixed points, called derangements, is then:
$$X_N=\left(\bigcup_kS_N^k\right)^c$$

Now the inclusion-exclusion principle tells us that we have:
\begin{eqnarray*}
|X_N|
&=&\left|\left(\bigcup_kS_N^k\right)^c\right|\\
&=&|S_N|-\sum_k|S_N^k|+\sum_{k<l}|S_N^k\cap S_N^l|-\ldots+(-1)^N\sum_{k_1<\ldots<k_N}|S_N^{k_1}\cup\ldots\cup S_N^{k_N}|\\
&=&N!-N(N-1)!+\binom{N}{2}(N-2)!-\ldots+(-1)^N\binom{N}{N}(N-N)!\\
&=&\sum_{r=0}^N(-1)^r\binom{N}{r}(N-r)!
\end{eqnarray*}

Thus, the probability that we are interested in, for a random permutation $\sigma\in S_N$ to have no fixed points, is given by the following formula:
$$P
=\frac{|X_N|}{N!}
=\sum_{r=0}^N\frac{(-1)^r}{r!}$$

Since on the right we have the expansion of $1/e$, this gives the result.

\medskip

(2) Let us construct now the main character of $S_N$, as in the statement. The permutation matrices being given by $\sigma_{ij}=\delta_{i\sigma(j)}$, we have the following formula:
$$\chi(\sigma)
=\sum_i\delta_{\sigma(i)i}
=\#\left\{i\in\{1,\ldots,N\}\Big|\sigma(i)=i\right\}$$

In order to establish now the asymptotic result in the statement, regarding these characters, we must prove the following formula, for any $r\in\mathbb N$, in the $N\to\infty$ limit:
$$P(\chi=r)\simeq\frac{1}{r!e}$$

We already know, from (1), that this formula holds at $r=0$. In the general case now, we have to count the permutations $\sigma\in S_N$ having exactly $r$ points. Now since having such a permutation amounts in choosing $r$ points among $1,\ldots,N$, and then permuting the $N-r$ points left, without fixed points allowed, we have:
\begin{eqnarray*}
\#\left\{\sigma\in S_N\Big|\chi(\sigma)=r\right\}
&=&\binom{N}{r}\#\left\{\sigma\in S_{N-r}\Big|\chi(\sigma)=0\right\}\\
&=&\frac{N!}{r!(N-r)!}\#\left\{\sigma\in S_{N-r}\Big|\chi(\sigma)=0\right\}\\
&=&N!\times\frac{1}{r!}\times\frac{\#\left\{\sigma\in S_{N-r}\Big|\chi(\sigma)=0\right\}}{(N-r)!}
\end{eqnarray*}

By dividing everything by $N!$, we obtain from this the following formula:
$$\frac{\#\left\{\sigma\in S_N\Big|\chi(\sigma)=r\right\}}{N!}=\frac{1}{r!}\times\frac{\#\left\{\sigma\in S_{N-r}\Big|\chi(\sigma)=0\right\}}{(N-r)!}$$

Now by using the computation at $r=0$, that we already have, from (1), it follows that with $N\to\infty$ we have the following estimate:
$$P(\chi=r)
\simeq\frac{1}{r!}\cdot P(\chi=0)
\simeq\frac{1}{r!}\cdot\frac{1}{e}$$

Thus, we obtain as limiting measure the Poisson law of parameter 1, as stated.

\medskip

(3) Finally, let us construct the truncated characters of $S_N$, as in the statement. As before in the case $t=1$, we have the following computation, coming from definitions:
$$\chi_t(\sigma)
=\sum_{i=1}^{[tN]}\delta_{\sigma(i)i}
=\#\left\{i\in\{1,\ldots,[tN]\}\Big|\sigma(i)=i\right\}$$

Also before in the proofs of (1) and (2), we obtain by inclusion-exclusion that:
\begin{eqnarray*}
P(\chi_t=0)
&=&\frac{1}{N!}\sum_{r=0}^{[tN]}(-1)^r\sum_{k_1<\ldots<k_r<[tN]}|S_N^{k_1}\cap\ldots\cap S_N^{k_r}|\\
&=&\frac{1}{N!}\sum_{r=0}^{[tN]}(-1)^r\binom{[tN]}{r}(N-r)!\\
&=&\sum_{r=0}^{[tN]}\frac{(-1)^r}{r!}\cdot\frac{[tN]!(N-r)!}{N!([tN]-r)!}
\end{eqnarray*}

Now with $N\to\infty$, we obtain from this the following estimate:
$$P(\chi_t=0)
\simeq\sum_{r=0}^{[tN]}\frac{(-1)^r}{r!}\cdot t^r
\simeq e^{-t}$$

More generally, by counting the permutations $\sigma\in S_N$ having exactly $r$ fixed points among $1,\ldots,[tN]$, as in the proof of (2), we obtain:
$$P(\chi_t=r)\simeq\frac{t^r}{r!e^t}$$

Thus, we obtain in the limit a Poisson law of parameter $t$, as stated.
\end{proof}

Escalating difficulties, let us discuss now the hyperoctahedral group $H_N$. Here the result is more technical, getting us into more advanced probability, as follows:

\begin{theorem}
For the hyperoctahedral group $H_N\subset O_N$, the law of the truncated character $\chi=g_{11}+\ldots +g_{ss}$ with $s=[tN]$ is, in the $N\to\infty$ limit, the measure
$$b_t=e^{-t}\sum_{r=-\infty}^\infty\delta_r\sum_{p=0}^\infty \frac{(t/2)^{|r|+2p}}{(|r|+p)!p!}$$
called Bessel law of parameter $t>0$.
\end{theorem}

\begin{proof}
We regard $H_N$ as being the symmetry group of the graph $I_N=\{I^1,\ldots ,I^N\}$ formed by $n$ segments. The diagonal coefficients are then given by:
$$u_{ii}(g)=\begin{cases}
\ 0\ \mbox{ if $g$ moves $I^i$}\\
\ 1\ \mbox{ if $g$ fixes $I^i$}\\
-1\mbox{ if $g$ returns $I^i$}
\end{cases}$$

Let us denote by $F_g,R_g$ the number of segments among $\{I^1,\ldots ,I^s\}$ which are fixed, respectively returned by an element $g\in H_N$. With this notation, we have:
$$u_{11}+\ldots+u_{ss}=F_g-R_g$$

Let us denote by $P_N$ the probabilities computed over $H_N$. The density of the law of the variable $u_{11}+\ldots+u_{ss}$ at a point $r\geq 0$ is then given by the following formula:
$$D(r)
=P_N(F_g-R_g=r)
=\sum_{p=0}^\infty P_N(F_g=r+p,R_g=p)$$

Assume first that we are in the case $t=1$. We have the following computation:
\begin{eqnarray*}
\lim_{N\to\infty}D(r)
&=&\lim_{N\to\infty}\sum_{p=0}^\infty(1/2)^{r+2p}\binom{r+2p}{r+p}P_N(F_g+R_g=r+2p)\\ 
&=&\sum_{p=0}^\infty(1/2)^{r+2p}\binom{r+2p}{r+p}\frac{1}{e(r+2p)!}\\
&=&\frac{1}{e}\sum_{p=0}^\infty\frac{(1/2)^{r+2p}}{(r+p)!p!}
\end{eqnarray*}

The general case $0<t\leq 1$ follows by performing some modifications in the above computation. Indeed, the asymptotic density can be computed as follows:
\begin{eqnarray*}
\lim_{N\to\infty}D(r)
&=&\lim_{N\to\infty}\sum_{p=0}^\infty(1/2)^{r+2p}\binom{r+2p}{r+p}P_N(F_g+R_g=r+2p)\\
&=&\sum_{p=0}^\infty(1/2)^{r+2p}\binom{r+2p}{r+p}\frac{t^{r+2p}}{e^t(r+2p)!}\\
&=&e^{-t}\sum_{p=0}^\infty\frac{(t/2)^{r+2p}}{(r+p)!p!}
\end{eqnarray*}

Together with $D(-r)=D(r)$, this gives the formula in the statement.
\end{proof}

In order to further discuss now all this, and extend the above results, we will need a number of standard probabilistic preliminaries. We have the following notion:

\begin{definition}
Associated to any compactly supported positive measure $\nu$ on $\mathbb C$, not necessarily of mass $1$, is the probability measure
$$p_\nu=\lim_{n\to\infty}\left(\left(1-\frac{t}{n}\right)\delta_0+\frac{1}{n}\nu\right)^{*n}$$
where $t=mass(\nu)$, called compound Poisson law.
\end{definition}

In what follows we will be mainly interested in the case where the measure $\nu$ is discrete, as is for instance the case for $\nu=t\delta_1$ with $t>0$, which produces the Poisson laws. The following standard result allows one to detect compound Poisson laws:

\begin{proposition}
For $\nu=\sum_{i=1}^st_i\delta_{z_i}$ with $t_i>0$ and $z_i\in\mathbb C$, we have
$$F_{p_\nu}(y)=\exp\left(\sum_{i=1}^st_i(e^{iyz_i}-1)\right)$$
where $F$ denotes the Fourier transform.
\end{proposition}

\begin{proof}
Let $\mu_n$ be the Bernoulli measure appearing in Definition 11.32, under the convolution sign. We have then the following Fourier transform computation:
\begin{eqnarray*}
F_{\mu_n}(y)=\left(1-\frac{t}{n}\right)+\frac{1}{n}\sum_{i=1}^st_ie^{iyz_i}
&\implies&F_{\mu_n^{*n}}(y)=\left(\left(1-\frac{t}{n}\right)+\frac{1}{n}\sum_{i=1}^st_ie^{iyz_i}\right)^n\\
&\implies&F_{p_\nu}(y)=\exp\left(\sum_{i=1}^st_i(e^{iyz_i}-1)\right)
\end{eqnarray*}

Thus, we have obtained the formula in the statement.
\end{proof}

Getting back now to the Poisson and Bessel laws, we have:

\begin{theorem}
The Poisson and Bessel laws are compound Poisson laws,
$$p_t=p_{t\delta_1}\quad,\quad b_t=p_{t\varepsilon}$$
where $\delta_1$ is the Dirac mass at $1$, and $\varepsilon$ is the centered Bernoulli law, $\varepsilon=(\delta_1+\delta_{-1})/2$.
\end{theorem}

\begin{proof}
We have two assertions here, the idea being as follows:

\medskip

(1) The first assertion, regarding the Poisson law $p_t$, is clear from Definition 11.22, which for $\nu=t\delta_1$ takes the following form:
$$p_\nu=\lim_{n\to\infty}\left(\left(1-\frac{t}{n}\right)\delta_0+\frac{t}{n}\,\delta_1\right)^{*n}$$

But this is a usual Poisson limit, producing the Poisson law $p_t$, as desired.

\medskip

(2) Regarding the second assertion, concerning the Bessel law $b_t$, by further building on the various formulae from Theorem 11.31 and its proof, we have:
$$F_{b_t}(y)=\exp\left(t\left(\frac{e^{iy}+e^{-iy}}{2}-1\right)\right)$$

On the other hand, the formula in Proposition 11.33 gives, for $\nu=t\varepsilon$, the same Fourier transform formula. Thus, with $\nu=t\varepsilon$ we have $p_\nu=b_t$, as claimed.
\end{proof}

More generally now, and in answer to Problem 11.29, for the groups $H_N^s$ from chapter 9 the asymptotic truncated characters follow the law $b_t^s=p_{t\varepsilon_s}$, with $\varepsilon_s$ being the uniform measure on the $s$-roots of unity, and if you are really picky about what happens with $N\to\infty$, these are the only solutions. For details on this, you can check my book \cite{ba1}.

\section*{11e. Exercises}

We had an exciting chapter here, mixing advanced algebra techniques with advanced analysis and probability. As exercises on all this, we have:

\begin{exercise}
Fine-tune your Peter-Weyl learning, with some character tables.
\end{exercise}

\begin{exercise}
Clarify what else can be done with magics, in the graph context.
\end{exercise}

\begin{exercise}
Learn more about planar algebras, and their relation with graphs.
\end{exercise}

\begin{exercise}
Learn about the representations of $S_N$, at fixed $N\in\mathbb N$.
\end{exercise}

\begin{exercise}
Have a look at the representations of $H_N$ too, at fixed $N\in\mathbb N$.
\end{exercise}

\begin{exercise}
Revise your probability knowledge, notably with the CLT and PLT.
\end{exercise}

\begin{exercise}
Learn more about compound Poisson laws, and their properties.
\end{exercise}

\begin{exercise}
Learn more about the Bessel laws, and their properties.
\end{exercise}

As bonus exercise, learn some advanced probability. Among others, because we will soon get into quantum, and there is no quantum without probability.

\chapter{Cayley graphs}

\section*{12a. Cayley graphs}

We discuss here the Cayley graphs of finite groups, and what can be done with them. We have already met these graphs in chapter 4, their definition being as follows:

\index{Cayley graph}

\begin{definition}
Associated to any finite group $G=<S>$, with the generating set $S$ assumed to satisfy $1\notin S=S^{-1}$, is its Cayley graph, constructed as follows:
\begin{enumerate}
\item The vertices are the group elements $g\in G$.

\item Edges $g-h$ are drawn when $h=gs$, with $s\in S$.
\end{enumerate} 
\end{definition}

As a first observation, the Cayley graph is indeed a graph, because our assumption $S=S^{-1}$ on the generating set shows that we have $g-h\implies h-g$, as we should, and also because our assumption $1\notin S$ excludes the self-edges, $g\not\!\!-\,g$.

\bigskip

Observe also that the Cayley graph depends, in a crucial way, on the generating set $S$ satisfying $1\notin S=S^{-1}$. Indeed, if we choose for instance $S=G-\{1\}$, we obtain the complete graph with $N=|G|$ vertices, and this regardless of what our group $G$ is. Thus, the Cayley graph as constructed above is not exactly associated to the group $G$, but rather to the group $G$ viewed as finitely generated group, $G=<S>$. 

\bigskip

In view of this, we will usually look for generating sets $S$ which are minimal, in order to perform our Cayley graph construction, and get non-trivial graphs. Here are now some examples, with minimal generating sets, that we know well from chapter 4:

\begin{proposition}
We have the following examples of Cayley graphs, each time with respect to the standard, minimal generating set satisfying $S=S^{-1}$:
\begin{enumerate}
\item For $\mathbb Z_N$ we obtain the cycle graph $C_N$.

\item For $\mathbb Z_2\times\mathbb Z_3$ we obtain the prism $P(C_3)$.

\item For $\mathbb Z_2^N$ we obtain the hypercube graph $\square_N$.
\end{enumerate}
\end{proposition}

\begin{proof}
This is something elementary, that we know well from chapter 4, with the generating set in question being $S=\{-1,1\}$ for the cyclic group $\mathbb Z_N$, written additively, unless we are in the case $N=2$, where this set simply becomes $S=\{1\}$, again written additively, and with the generating sets for products being obtained in the obvious minimal way, by taking the union of the generating sets of the components:

\medskip

(1) For the group $\mathbb Z_N=<-1,1>$, written additively, our condition for the edges $g-h$ reads $g=h\pm1$, so we are led to the cycle graph $C_N$, namely:
$$
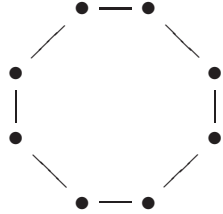
\xymatrix@R=12pt@C=12pt{
&\bullet\ar@{-}[r]\ar@{-}[dl]&\bullet\ar@{-}[dr]\\
\bullet\ar@{-}[d]&&&\bullet\ar@{-}[d]\\
\bullet\ar@{-}[dr]&&&\bullet\ar@{-}[dl]\\
&\bullet\ar@{-}[r]&\bullet}$$

(2) For the group $\mathbb Z_2\times\mathbb Z_3=<(1,0),(0,1),(0,2)>$, again written additively, our condition for the edges takes the following form:
$$(g,a)-(h,b)\Longleftrightarrow g=h,\, a=b\pm1\mbox{ \rm{or} }g=h+1,a=b$$

But this leads to the prism graph $P(C_3)$, which is as follows:
 $$
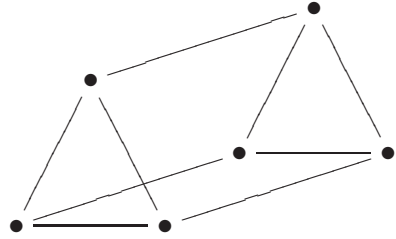
\xymatrix@R=15pt@C=15pt{
&&&&\bullet\ar@{-}[ddl]\ar@{-}[ddr]\\
&\bullet\ar@{-}[ddl]\ar@{-}[ddr]\ar@{-}[urrr]\\
&&&\bullet\ar@{-}[rr]&&\bullet\\
\bullet\ar@{-}[rr]\ar@{-}[urrr]&&\bullet\ar@{-}[urrr]
}$$

(3) Finally, for the group $\mathbb Z_2^N=<1_1,\ldots,1_N>$, with $1_i$ standing for the standard generators of the components, our condition for the edges takes the following form:
$$g-h\iff\exists!\,i,g_i\neq h_i$$

Now if we represent the elements of $\mathbb Z_2^N=(0,1)^N$ as points in $\mathbb R^N$, in the obvious way, we are led to the hypercube graph $\square_N$, which is as follows:
$$
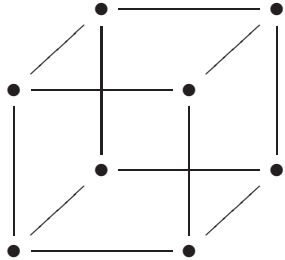
\xymatrix@R=18pt@C=20pt{
&\bullet\ar@{-}[rr]&&\bullet\\
\bullet\ar@{-}[rr]\ar@{-}[ur]&&\bullet\ar@{-}[ur]\\
&\bullet\ar@{-}[rr]\ar@{-}[uu]&&\bullet\ar@{-}[uu]\\
\bullet\ar@{-}[uu]\ar@{-}[ur]\ar@{-}[rr]&&\bullet\ar@{-}[uu]\ar@{-}[ur]
}$$

Thus, we are led to the conclusions in the statement. 
\end{proof}

The above result dates back to chapter 4, and time now to improve it, by using the product operations for graphs, that we learned in chapter 10. We have:

\begin{theorem}
Given two groups $G=<S>$ and $H=<T>$, we have
$$G\times H=<S\times1,1\times T>$$
and at the level of the corresponding Cayley graphs, we have the formula
$$X_{G\times H}=X_G\,\square\,X_H$$
involving the Cartesian product operation for graphs $\square$.
\end{theorem}

\begin{proof}
We have indeed a generating set, which satisfies the condition $1\notin S=S^{-1}$. Now observe that when constructing the Cayley graph, the edges are as follows:
$$(g,a)-(h,b)\Longleftrightarrow g=h,\, a-b\mbox{ \rm{or} }g-h,a=b$$

Thus, we obtain indeed a Cartesian product $X_{G\times H}=X_G\,\square\,X_H$, as claimed.
\end{proof}

Along the same lines, let us record as well the following result, which is however rather anecdotical, first because the subgroup $K\subset G\times H$ appearing below might not equal the group $G\times H$ itself, and also because, even when it does, the generating set $S\times T$ that we use is much bigger than the generating set $S\times1,1\times T$ from Theorem 12.3:

\begin{proposition}
Given two groups $G=<S>$ and $H=<T>$, if we set
$$K=<S\times T>\subset G\times H$$
then at the level of the corresponding Cayley graphs, we have the formula
$$X_K=X_G\times X_H$$
involving the direct product operation for graphs $\times$.
\end{proposition}

\begin{proof}
As before, we have a set satisfying the condition $1\notin S=S^{-1}$. Now observe that when constructing the Cayley graph, the edges are as follows:
$$(g,a)-(h,b)\Longleftrightarrow g-h,\, a-b$$

Thus, we obtain in this case a direct product $X_K=X_G\times X_H$, as claimed.
\end{proof}

As mentioned above, this latter result is something anecdotical, and we will not use it, in what follows. Thus, our convention for products will be the one in Theorem 12.3. We can now generalize the various computations from Proposition 12.2, as follows:

\begin{theorem}
For a finite abelian group, written as
$$G=\mathbb Z_{N_1}\times\ldots\times\mathbb Z_{N_k}$$
and with standard generating set, the corresponding Cayley graph is:
$$X_G=C_{N_1}\,\square\,\ldots\,\square\,C_{N_k}$$
That is, we obtain a Cartesian product of cycle graphs.
\end{theorem}

\begin{proof}
This follows indeed from Proposition 12.2 (1), and from Theorem 12.3.
\end{proof}

Getting now to symmetry groups, as a first observation, we have:

\begin{proposition}
Given a finite group, $G=<S>$ with $1\notin S=S^{-1}$, we have an action of this group on its associated Cayley graph,
$$G\curvearrowright X_G$$
and so $G\subset G(X_G)$. However, this inclusion is not an isomorphism, in general.
\end{proposition}

\begin{proof}
We have several assertions here, the idea being as follows:

\medskip

(1) Consider indeed the Cayley action of $G$ on itself, which is given by:
$$G\subset S_G\quad,\quad g\to[h\to gh]$$

(2) Thus $G$ acts on the vertices of its Cayley graph $X_G$, and our claim is that the edges are preserved by this action. Indeed, given an edge, $h-hs$ with $s\in S$, a group element $g\in G$ will transform it into $gh-ghs$, which is an edge too.

\medskip

(3) Thus, the first assertion holds indeed. As for the second assertion, this holds too, and in a very convincing way, for instance because for the groups in Proposition 12.2, the corresponding Cayley graphs and their symmetry groups are as follows:
$$\mathbb Z_N\to C_N\to D_N$$
$$\mathbb Z_2\times\mathbb Z_3\to P(C_3)\to D_6$$
$$\mathbb Z_2^N\to\square_N\to H_N$$

Thus, we are led to the conclusions in the statement.
\end{proof}

In relation with the above, a first question would be that of suitably fine-tuning the construction of the Cayley graph $X_G$, as to have as end result $G=G(X_G)$, which would be something nice. We will discuss this later, the idea being that this can be done indeed, by adding orientation, or colors, or both, to the construction of the Cayley graph $X_G$.

\bigskip

As a second question now, we can try to understand the operation $G\to G(X_G)$, with our present definition for the Cayley graph $X_G$. Many things can be said here, and to start with, we have the following result, generalizing our computation for $\mathbb Z_2^N$:

\begin{theorem}
The symmetry group of the Cayley graph of $G=\mathbb Z_s^N$, which is
$$X_G=C_s^{\,\square N}$$
according to our previous results, is the complex reflection group $H_N^s$:
$$G(X_G)=H_N^s$$
Thus, in this case, the operation $G\to G(X_G)$ is given by $\mathbb Z_s^N\to\mathbb Z_s^N\rtimes S_N$.
\end{theorem}

\begin{proof}
This is something that we discussed in chapter 4 at $s=2$, with the graph involved there being the hypercube $X_G=\square_N$. In general the proof is similar, by using the theory of the complex reflection groups $H_N^s$, that we developed in chapter 9.
\end{proof}

Many other things can be said, along these lines, for instance with a generalization of the computation $\mathbb Z_2\times\mathbb Z_3\to P(C_3)\to D_6$ too, and with a number of other results regarding the operation $G\to X_G\to G(X_G)$, for the finite abelian groups $G$. We will leave doing some study here as an instructive exercise.

\section*{12b. Orientation, colors}

We discuss now various versions of the Cayley graph construction $G\to X_G$, which are all technically useful, obtained by adding orientation, or colors, or both. Let us start with a straightforward oriented version of Definition 12.1, as follows:

\begin{definition}
Associated to any finite group $G=<S>$, with the generating set $S$ assumed to satisfy $1\notin S$, is its oriented Cayley graph, constructed as follows:
\begin{enumerate}
\item The vertices are the group elements $g\in G$.

\item Edges $g\to h$ are drawn when $h=gs$, with $s\in S$.
\end{enumerate} 
\end{definition}

Observe that the oriented Cayley graph is indeed a graph, because our assumption $1\notin S$ excludes the self-edges, $g\not\!\!-\,g$. Observe also that in the case $S=S^{-1}$ each edge $g\to h$ has as companion edge $h\to g$, so in this case, up to replacing the pairs of edges $g\leftrightarrow h$ by usual edges $g-h$, we obtain the previous, unoriented Cayley graph.

\bigskip

As in the case of unoriented Cayley graphs, our graph depends a lot on the chosen writing $G=<S>$, and the point is to look for generating sets $S$ which are minimal, in order to perform our Cayley graph construction, and get non-trivial graphs. 

\bigskip

In view of this, what we basically have to do, in order to get started with some theory, is to review the material from the previous section, by dropping the assumption $S=S^{-1}$ there, which in practice means to remove around $1/2$ of the generators.

\bigskip

Getting started now, as a first result, in relation with Proposition 12.2, we have:

\begin{proposition}
We have the following examples of oriented Cayley graphs, each time with respect to the standard, minimal generating set:
\begin{enumerate}
\item For $\mathbb Z_N$ we obtain the oriented cycle $C_N^o$.

\item For $\mathbb Z_2\times\mathbb Z_3$ we obtain the oriented prism $P(C_3^o)$.

\item For $\mathbb Z_2^N$ we obtain the hypercube graph $\square_N$.
\end{enumerate}
\end{proposition}

\begin{proof}
This is something elementary, with the convention that the generating set is $S=\{1\}$ for the cyclic group $\mathbb Z_N$, written additively, and with the generating sets for products being obtained by taking the union of the generating sets for components:

\medskip

(1) For the group $\mathbb Z_N=<1>$, written additively, our condition for the edges $g\to h$ reads $h=g+1$, so we are led to the oriented cycle graph $C_N^o$, namely:
$$
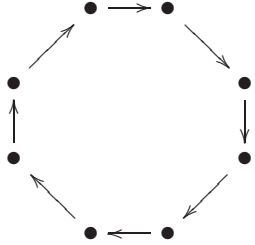
\xymatrix@R=16pt@C=16pt{
&\bullet\ar[r]&\bullet\ar[dr]\\
\bullet\ar[ur]&&&\bullet\ar[d]\\
\bullet\ar[u]&&&\bullet\ar[dl]\\
&\bullet\ar[ul]&\bullet\ar[l]}$$

As an observation here, that will be of importance later, in the particular case $N=2$ what we get is the graph $\bullet\longleftrightarrow\bullet$, which is the same as the usual segment, $\bullet-\!\!\!-\,\bullet$.

\medskip

(2) For the group $\mathbb Z_2\times\mathbb Z_3=<(1,0),(0,1)>$, again written additively, our condition for the edges takes the following form:
$$(g,a)\to(h,b)\Longleftrightarrow g=h,\, b=a+1\mbox{ \rm{or} }h=g+1,a=b$$

But this leads to the oriented prism graph $P(C_3^o)$, which is as follows, with the convention that the usual segments $\bullet-\!\!\!-\,\bullet$ stand for double edges $\bullet\longleftrightarrow\bullet$:
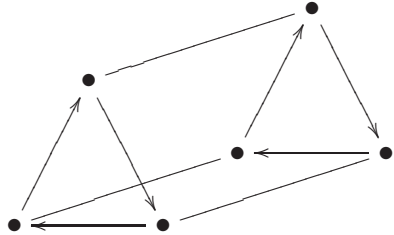
 $$\xymatrix@R=15pt@C=15pt{
&&&&\bullet\ar[ddr]\\
&\bullet\ar[ddr]\ar@{-}[urrr]\\
&&&\bullet\ar[uur]&&\bullet\ar[ll]\\
\bullet\ar@{-}[urrr]\ar[uur]&&\bullet\ar@{-}[urrr]\ar[ll]
}$$

(3) Finally, for the group $\mathbb Z_2^N=<1_1,\ldots,1_N>$, with $1_i$ standing for the standard generators of the components, written additively, there is no true orientation, due to the reason explained in (1), and we are led to the old hypercube graph $\square_N$, namely:
$$
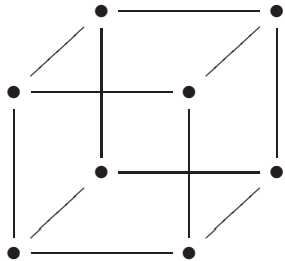
\xymatrix@R=18pt@C=20pt{
&\bullet\ar@{-}[rr]&&\bullet\\
\bullet\ar@{-}[rr]\ar@{-}[ur]&&\bullet\ar@{-}[ur]\\
&\bullet\ar@{-}[rr]\ar@{-}[uu]&&\bullet\ar@{-}[uu]\\
\bullet\ar@{-}[uu]\ar@{-}[ur]\ar@{-}[rr]&&\bullet\ar@{-}[uu]\ar@{-}[ur]
}$$

Thus, we are led to the conclusions in the statement.
\end{proof}

As before in the unoriented case, we can generalize these computations, as follows:

\begin{theorem}
Given two groups $G=<S>$ and $H=<T>$, we have
$$G\times H=<S\times1,1\times T>$$
and at the level of the corresponding oriented Cayley graphs, we have the formula
$$X_{G\times H}=X_G\,\square\,X_H$$
involving the Cartesian product operation for oriented graphs $\square$.
\end{theorem}

\begin{proof}
We have indeed a generating set, which satisfies the condition $1\notin S$. Now observe that when constructing the Cayley graph, the edges are as follows:
$$(g,a)\to(h,b)\Longleftrightarrow g=h,\, a\to b\mbox{ \rm{or} }g\to h,a=b$$

Thus, we obtain indeed a Cartesian product $X_{G\times H}=X_G\,\square\,X_H$, as claimed.
\end{proof}

At the level of symmetry groups now, of these Cayley graphs, things get more interesting, in the present oriented graph setting. Let us recall indeed that in the unoriented setting, the basic computations of type $G\to X_G\to G(X_G)$ were as follows:
$$\mathbb Z_N\to C_N\to D_N$$
$$\mathbb Z_2\times\mathbb Z_3\to P(C_3)\to D_6$$
$$\mathbb Z_2^N\to\square_N\to H_N$$

When adding orientation, according to Proposition 12.9, these computations become:
$$\mathbb Z_N\to C_N^o\to\mathbb Z_N$$
$$\mathbb Z_2\times\mathbb Z_3\to P(C_3^o)\to\mathbb Z_2\times\mathbb Z_3$$
$$\mathbb Z_2^N\to\square_N\to H_N$$

Thus, we are getting closer to a formula of type $G=G(X_G)$. In order to discuss this, let us start with a straightforward analogue of Proposition 12.6, as follows:

\begin{proposition}
Given a finite group, $G=<S>$ with $1\notin S$, we have an action of this group on its associated oriented Cayley graph,
$$G\curvearrowright X_G$$
and so $G\subset G(X_G)$. However, this inclusion is not an isomorphism, in general.
\end{proposition}

\begin{proof}
We have several assertions here, the idea being as follows:

\medskip

(1) Consider indeed the Cayley action of $G$ on itself, given by:
$$G\subset S_G\quad,\quad g\to[h\to gh]$$

(2) Thus $G$ acts on the vertices of its Cayley graph $X_G$, and our claim is that the edges are preserved by this action. Indeed, given an edge, $h\to hs$ with $s\in S$, a group element $g\in G$ will transform it into $gh\to ghs$, which is an edge too.

\medskip

(3) Thus, the first assertion holds indeed. As for the second assertion, this holds too, the counterexample coming from the computation $\mathbb Z_2^N\to\square_N\to H_N$, discussed above.
\end{proof}

Summarizing, as already said, we are certainly getting closer to a formula of type $G=G(X_G)$, but we are still not there, with some counterexamples still persisting. We will come back to this question in a moment, with a solution using colors.

\bigskip

Before leaving the subject, let us record the following result, due to Sabidussi:

\begin{theorem}
An oriented graph $X$ is the oriented Cayley graph of a given group $G$ precisely when it admits a simply transitive action of $G$.
\end{theorem}

\begin{proof}
This is something elementary, which is clear in one sense, from the proof of Proposition 12.11, and which in the other sense can be established as follows:

\medskip

-- Given a graph $X$ as in the statement, pick an arbitrary vertex, and label it 1.

\medskip

-- Then, label each vertex $v\in X$ by the unique $g\in G$ mapping $1\to v$.

\medskip

-- Finally, define $S\subset G$ as being the set of labels $i$, such that $1\to i$. 

\medskip

Indeed, with these operations performed, it follows from definitions that what we get is a generating set for our group, $G=<S>$, satisfying the condition $1\notin S$, and the Cayley graph of $G=<S>$ follows to be the original graph $X$ itself, as desired.
\end{proof}

Moving ahead, let us attempt now to further modify our Cayley graph formalism, as to have $G=G(X_G)$. As a first observation, this is certainly possible, due to:

\index{colored graph}

\begin{theorem}
Any finite group $G$ appears as $G=G(X_G)$, with $X_G$ being the complete oriented graph having $G$ as set of vertices, and with the edges being colored by
$$d_{hk}=h^{-1}k$$
according to the usual colored graph conventions, with color set $C=G$.
\end{theorem}

\begin{proof}
Consider indeed the Cayley action of $G$ on itself, which is given by:
$$G\subset S_G\quad,\quad g\to[h\to gh]$$

We have $d_{gh,gk}=d_{hk}$, which gives an action $G\curvearrowright X_G$, and so an inclusion $G\subset G(X_G)$. Conversely now, pick an arbitrary permutation $\sigma\in S_G$. We have then:
\begin{eqnarray*}
\sigma\in G(X_G)
&\implies&d_{\sigma(h)\sigma(k)}=d_{hk}\\
&\implies&\sigma(h)^{-1}\sigma(k)=h^{-1}k\\
&\implies&\sigma(1)^{-1}\sigma(k)=k\\
&\implies&\sigma(k)=\sigma(1)k\\
&\implies&\sigma\in G
\end{eqnarray*}

Thus, the inclusion $G\subset G(X_G)$ is an equality, as desired.
\end{proof}

Summarizing, we must add colors to our Cayley graph formalism, as follows:

\begin{definition}
Associated to any finite group $G=<S>$, with the generating set $S$ satisfying $1\notin S$, is its oriented colored Cayley graph, constructed as follows:
\begin{enumerate}
\item The vertices are the group elements $g\in G$.

\item Edges $g\to h$ are drawn when $h=gs$, with $s\in S$.

\item Such an edge $g\to h$ is colored by the element $s=g^{-1}h\in S$.
\end{enumerate} 
\end{definition}

As before, this oriented colored Cayley graph is indeed a graph, because our assumption $1\notin S$ excludes the self-edges, $g\not\!\!-\,g$. Observe also that in the case $S=S^{-1}$ each edge $g\to h$ has as companion edge $h\to g$, so in this case, up to replacing the pairs of edges $g\leftrightarrow h$ by usual edges $g=h$, we obtain an unoriented colored graph.

\bigskip

At the level of examples, the graph in Theorem 12.13 appears as in Definition 12.14, with generating set $S=G-\{1\}$. Also, in relation with Proposition 12.2 and Proposition 12.9, and with the associated symmetry groups too, we have the following result:

\begin{proposition}
We have the following oriented colored Cayley graphs, each time with respect to the standard, minimal generating set:
\begin{enumerate}
\item For $\mathbb Z_N$ we obtain $C_N^o$, having symmetry group $\mathbb Z_N$.

\item For $\mathbb Z_2\times\mathbb Z_3$ we obtain the bicolored $P(C_3^o)$, with symmetry group $\mathbb Z_2\times\mathbb Z_3$.

\item For $\mathbb Z_2^N$ we obtain the $N$-colored cube $\square_N$, having symmetry group $\mathbb Z_2^N$.
\end{enumerate}
\end{proposition}

\begin{proof}
This basically comes from our previous computations, as follows:

\medskip

(1) For the group $\mathbb Z_N=<1>$, written as usual additively, the number of colors is $|\{1\}|=1$, and so no colors, and the graph is the usual oriented cycle $C_N^o$, namely:
$$
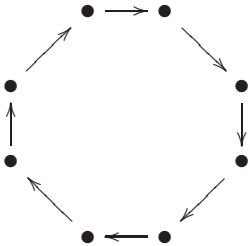
\xymatrix@R=16pt@C=16pt{
&\bullet\ar[r]&\bullet\ar[dr]\\
\bullet\ar[ur]&&&\bullet\ar[d]\\
\bullet\ar[u]&&&\bullet\ar[dl]\\
&\bullet\ar[ul]&\bullet\ar[l]}$$

But, as we know well from the above, the symmetry group of this graph is $\mathbb Z_N$.

\medskip

(2) For the group $\mathbb Z_2\times\mathbb Z_3=<(1,0),(0,1)>$, again written additively, we have 2 colors, and we obtain the bicolored version of the oriented prism $P(C_3^o)$, which is as follows:
 $$
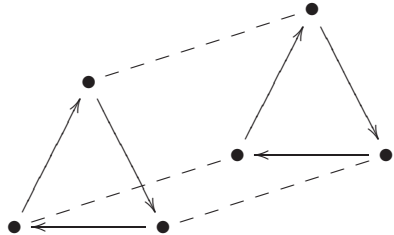
\xymatrix@R=15pt@C=15pt{
&&&&\bullet\ar[ddr]\\
&\bullet\ar[ddr]\ar@{--}[urrr]\\
&&&\bullet\ar[uur]&&\bullet\ar[ll]\\
\bullet\ar@{--}[urrr]\ar[uur]&&\bullet\ar@{--}[urrr]\ar[ll]
}$$

But the symmetry group of this bicolored oriented prism is $\mathbb Z_2\times\mathbb Z_3$, as claimed.

\medskip

(3) Finally, for the group $\mathbb Z_2^N=<1_1,\ldots,1_N>$, with $1_i$ standing for the standard generators of the components, written as usual additively, here there is no true orientation, due to $1=-1$, but we have however $N$ colors, those of the generators. We are led in this way to the $N$-colored version of the old hypercube graph $\square_N$, namely:
$$
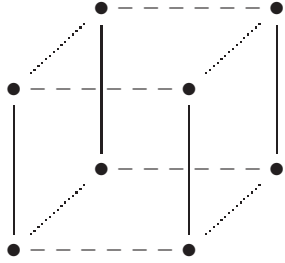
\xymatrix@R=18pt@C=20pt{
&\bullet\ar@{--}[rr]&&\bullet\\
\bullet\ar@{--}[rr]\ar@{.}[ur]&&\bullet\ar@{.}[ur]\\
&\bullet\ar@{--}[rr]\ar@{-}[uu]&&\bullet\ar@{-}[uu]\\
\bullet\ar@{-}[uu]\ar@{.}[ur]\ar@{--}[rr]&&\bullet\ar@{-}[uu]\ar@{.}[ur]
}$$

But the symmetry group of this latter graph is $\mathbb Z_2^N$, as claimed.
\end{proof}

All this is quite interesting, so let us generalize now the above results. In what regards the product operations, the result here, which is very standard, is as follows:

\begin{theorem}
Given two groups $G=<S>$ and $H=<T>$, we have
$$G\times H=<S\times1,1\times T>$$
and at the level of the corresponding oriented colored Cayley graphs, we have the formula
$$X_{G\times H}=X_G\,\square\,X_H$$
involving the Cartesian product operation for oriented colored graphs $\square$.
\end{theorem}

\begin{proof}
We have indeed a generating set, which satisfies the condition $1\notin S$. Now observe that when constructing the Cayley graph, the edges are as follows:
$$(g,a)\to(h,b)\Longleftrightarrow g=h,\, a\to b\mbox{ \rm{or} }g\to h,a=b$$

Thus, we obtain indeed a Cartesian product $X_{G\times H}=X_G\,\square\,X_H$, as claimed.
\end{proof}

As for the generalization of our symmetry group computations, this is as follows:

\begin{theorem}
Given a group $G=<S>$ with $1\notin S$, we have the formula 
$$G=G(X_G)$$
with $X_G$ being the corresponding oriented colored Cayley graph.
\end{theorem}

\begin{proof}
We use the same method as for Theorem 12.13, which corresponds to the case $S=G-\{1\}$. The adjacency matrix of the graph $X_G$ in the statement is given by:
$$d_{hk}
=\begin{cases}
h^{-1}k&{\rm if}\ h^{-1}k\in S\\
0&{\rm otherwise}
\end{cases}$$

Thus we have an action $G\curvearrowright X_G$, and so on inclusion $G\subset G(X_G)$. Conversely now, pick an arbitrary permutation $\sigma\in S_G$. We know that $\sigma$ must preserve all the color components of $d$, which are the following matrices, depending on a color $c\in S$:
$$d_{hk}^c
=\begin{cases}
1&{\rm if}\ h^{-1}k=c\\
0&{\rm otherwise}
\end{cases}$$

In other words, we have the following equivalences:
\begin{eqnarray*}
\sigma\in G(X)
&\iff&d^c_{\sigma(h)\sigma(k)}=d_{hk},\forall c\in S\\
&\iff&\sigma(h)^{-1}\sigma(k)=h^{-1}k,\forall h^{-1}k\in S
\end{eqnarray*}

Now observe that with $h=1$ we obtain from this that we have:
\begin{eqnarray*}
k\in S
&\implies&\sigma(1)^{-1}\sigma(k)=k\\
&\implies&\sigma(k)=\sigma(1)k
\end{eqnarray*}

Next, by taking $h\in S$, we obtain from the above formula that we have:
\begin{eqnarray*}
k\in hS
&\implies&\sigma(h)^{-1}\sigma(k)=h^{-1}k\\
&\implies&\sigma(k)=\sigma(h)h^{-1}k\\
&\implies&\sigma(k)=(\sigma(1)h)h^{-1}k\\
&\implies&\sigma(k)=\sigma(1)k
\end{eqnarray*}

Thus with $g=\sigma(1)$ we have the following formula, for any $k\in S$:
$$\sigma(k)=gk$$

But  the same method shows that this formula holds as well for any $k\in S^2$, then for any $k\in S^3$, any $k\in S^4$, and so on. Thus the above formula $\sigma(k)=gk$ holds for any $k\in G$, and so the inclusion $G\subset G(X_G)$ is an equality, as desired.
\end{proof}

The above results are not the end of the story, but rather the beginning of it. Indeed, at a more advanced level, we have the following classical result, due to Frucht:

\index{Frucht theorem}

\begin{theorem}[Frucht]
Any finite group $G$ appears as the symmetry group
$$G=G(X)$$
of a certain uncolored, unoriented graph $X$.
\end{theorem}

\begin{proof}
This is something quite tricky, the idea being as follows:

\medskip

(1) We can start with the graph in Definition 12.8, namely the associated oriented Cayley graph $X_G$, with the convention that the edges $g\to h$ exist when:
$$g^{-1}h\in S$$

(2) The point now is that we can make a suitable unoriented graph out of this graph, by replacing each edge with a copy of the following graph, with the height being in a chosen bijection with the corresponding element $g^{-1}h\in S$: 
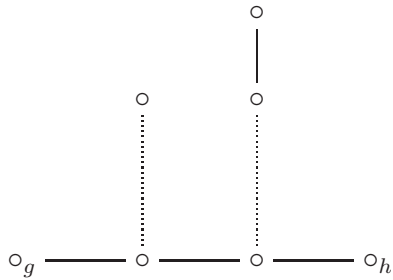
$$\xymatrix@R=20pt@C=30pt{
&&\circ\\
&\circ&\circ\ar@{-}[u]\\
\\
\circ_g\ar@{-}[r]&\circ\ar@{.}[uu]\ar@{-}[r]&\circ\ar@{.}[uu]\ar@{-}[r]&\circ_h}$$

(3) With these replacements made, and under suitable assumptions on the generating set $S$, namely the usual $1\notin S$ assumption, plus the fact that $S\cap S^{-1}$ must consist only of involutions, one can prove that $G$ appears indeed of the symmetry group of this graph $X$. We will leave checking the details here as an instructive exercise.
\end{proof}

\section*{12c. Brauer theorems}

The above results are certainly interesting from a graph theory perspective, but from a group theory perspective, they remain a bit anecdotical. We would like to present now a series of alternative results, going in the other sense, that is, featuring less graphs, or rather featuring some combinatorial objects which are more complicated and abstract than graphs, but which can be extremely useful for the study of finite groups.

\bigskip

The idea will be a bit the same as for the Frucht theorem, namely that of viewing an arbitrary finite group $G$ as symmetry group of a combinatorial object, $G=G(X)$. What will change, however, is the nature of $X$, the general principle being as follows:

\begin{principle}
Any finite, or even general compact group $G$ appears as the symmetry group of its corresponding Tannakian category $C_G$,
$$G=G(C_G)$$
and by suitably delinearizing $C_G$, say via a Brauer theorem of type $C_G=span(D_G)$, we can view $G$ as symmetry group of a certain combinatorial object $D_G$.
\end{principle}

Excited about this? Does not look easy, all this material, with both Tannaka and Brauer being quite scary names, in the context of group theory. But, believe me, all this is worth learning, and there is nothing better, when doing graphs or any kind of other algebraic discipline, to have in your bag some cutting-edge technology regarding the groups, such as the results of Tannaka and Brauer. So, we will go for this.

\bigskip

Getting started now, we will develop our theory as a continuation of the Peter-Weyl theory developed in chapter 11. That theory was developed for the finite groups, but with some minimal changes, and we will leave clarifying the details here to you, everything works in fact for a compact group $G$. So, our starting point will be:

\begin{theorem}
We have the following Peter-Weyl results, valid for any compact group $G\subset_u U_N$, with the orthogonality being with respect to the Haar integration:
\begin{enumerate}
\item Any representation decomposes as a sum of irreducible representations.

\item Each irreducible representation appears inside a certain tensor power $u^{\otimes k}$.

\item $C(G)=\overline{\bigoplus}_{v\in Irr(G)}M_{\dim(v)}(\mathbb C)$, the summands being pairwise orthogonal.

\item The characters of irreducible representations form an orthonormal system.
\end{enumerate}
\end{theorem}

\begin{proof}
As explained above, this is something that we know from chapter 11, in the finite group case, and in the general compact group case the proof is similar, with the only technical point being that of proving, somewhere between (1,2) and (3,4), the existence of the Haar measure. But this can be proved by using the arguments from chapter 11, with that chapter being written precisely with this idea in mind, namely that of having a quite straightforward extension to the compact group case, whenever needed.   
\end{proof}

Going now towards Tannakian duality, let us start with:

\begin{definition}
A tensor category over $H=\mathbb C^N$ is a collection $C=(C_{kl})$ of linear spaces $C_{kl}\subset\mathcal L(H^{\otimes k},H^{\otimes l})$ satisfying the following conditions:
\begin{enumerate}
\item $S,T\in C$ implies $S\otimes T\in C$.

\item If $S,T\in C$ are composable, then $ST\in C$.

\item $T\in C$ implies $T^*\in C$.

\item Each $C_{kk}$ contains the identity operator.

\item $C_{\emptyset k}$ with $k=\circ\bullet,\bullet\circ$ contain the operator $R:1\to\sum_ie_i\otimes e_i$.

\item $C_{kl,lk}$ with $k,l=\circ,\bullet$ contain the flip operator $\Sigma:a\otimes b\to b\otimes a$.
\end{enumerate}
\end{definition}

Here, as usual, the tensor powers $H^{\otimes k}$, which are Hilbert spaces depending on a colored integer $k=\circ\bullet\bullet\circ\ldots\,$, are defined by the following formulae, and multiplicativity:
$$H^{\otimes\emptyset}=\mathbb C\quad,\quad 
H^{\otimes\circ}=H\quad,\quad
H^{\otimes\bullet}=\bar{H}\simeq H$$

We have already met such categories, when dealing with the Tannakian categories of the closed subgroups $G\subset U_N$, and our knowledge can be summarized as follows:

\index{tensor category}
\index{Tannakian category}

\begin{proposition}
Given a closed subgroup $G\subset U_N$, its Tannakian category
$$C_{kl}=\left\{T\in\mathcal L(H^{\otimes k},H^{\otimes l})\Big|Tg^{\otimes k}=g^{\otimes l}T,\forall g\in G\right\}$$
is a tensor category over $H=\mathbb C^N$. Conversely, given a tensor category $C$ over $\mathbb C^N$,
$$G=\left\{g\in U_N\Big|Tg^{\otimes k}=g^{\otimes l}T,\forall k,l,\forall T\in C_{kl}\right\}$$
is a closed subgroup of $U_N$.
\end{proposition}

\begin{proof}
This is something that we basically know, the idea being as follows:

\medskip

(1) Regarding the first assertion, we have to check here the axioms (1-6) in Definition 12.21. The axioms (1-4) being all clear from definitions, let us establish (5). But this follows from the fact that each element $g\in G$ is a unitary, which can be reformulated as follows, with $R:1\to\sum_ie_i\otimes e_i$ being the map in Definition 12.21:
$$R\in Hom(1,g\otimes\bar{g})\quad,\quad 
R\in Hom(1,\bar{g}\otimes g)$$

Regarding now the condition in Definition 12.21 (6), this comes from the fact that the matrix coefficients $g\to g_{ij}$ and their conjugates $g\to\bar{g}_{ij}$ commute with each other.

\medskip

(2) Regarding the second assertion, we have to check that the subset $G\subset U_N$ constructed in the statement is a closed subgroup. But, assuming $g,h\in G$, we have $gh\in G$, due to the following computation, valid for any $k,l$ and any $T\in C_{kl}$:
\begin{eqnarray*}
T(gh)^{\otimes k}
&=&Tg^{\otimes k}h^{\otimes k}\\
&=&g^{\otimes l}Th^{\otimes k}\\
&=&g^{\otimes l}h^{\otimes l}T\\
&=&(gh)^{\otimes l}T
\end{eqnarray*}

Also, we have $1\in G$, trivially. And also, assuming $g\in G$, we have $g^{-1}\in G$, due to:
\begin{eqnarray*}
T(g^{-1})^{\otimes k}
&=&(g^{-1})^{\otimes l}[g^{\otimes l}T](g^{-1})^{\otimes k}\\
&=&(g^{-1})^{\otimes l}[Tg^{\otimes k}](g^{-1})^{\otimes k}\\
&=&(g^{-1})^{\otimes l}T
\end{eqnarray*}

Finally, the fact that our subgroup $G\subset U_N$ is closed is clear from definitions.
\end{proof}

Summarizing, we have so far precise axioms for the tensor categories $C=(C_{kl})$, given in Definition 12.21, as well as correspondences as follows:
$$G\to C_G\quad,\quad 
C\to G_C$$

We will prove in what follows that these correspondences are inverse to each other. In order to get started, we first have the following technical result:

\begin{proposition}
Consider the following conditions:
\begin{enumerate}
\item $C=C_{G_C}$, for any tensor category $C$.

\item $G=G_{C_G}$, for any closed subgroup $G\subset U_N$.
\end{enumerate}
We have then $(1)\implies(2)$. Also, $C\subset C_{G_C}$ is automatic.
\end{proposition}

\begin{proof}
Given $G\subset U_N$, we have $G\subset G_{C_G}$. On the other hand, by using (1) we have $C_G=C_{G_{C_G}}$. Thus, we have an inclusion of closed subgroups of $U_N$, which becomes an isomorphism at the level of the associated Tannakian categories, so $G=G_{C_G}$. Finally, the fact that we have an inclusion $C\subset C_{G_C}$ is clear from definitions.
\end{proof}

The point now is that it is possible to prove that we have $C_{G_C}\subset C$, by doing some abstract algebra, and we are led in this way to the following conclusion:

\index{Tannakian duality}

\begin{theorem}
The Tannakian duality constructions 
$$C\to G_C\quad,\quad 
G\to C_G$$
are inverse to each other.
\end{theorem}

\begin{proof}
This is something quite tricky, the idea being as follows:

\medskip

(1) According to Proposition 12.23, we must prove $C_{G_C}\subset C$. For this purpose, given a tensor category $C=(C_{kl})$ over a Hilbert space $H$, consider the following $*$-algebra:
$$E_C
=\bigoplus_{k,l}C_{kl}
\subset\bigoplus_{k,l}B(H^{\otimes k},H^{\otimes l})
\subset B\left(\bigoplus_kH^{\otimes k}\right)$$

Consider also, inside this $*$-algebra, the following $*$-subalgebra:
$$E_C^{(s)}
=\bigoplus_{|k|,|l|\leq s}C_{kl}
\subset\bigoplus_{|k|,|l|\leq s}B(H^{\otimes k},H^{\otimes l})
=B\left(\bigoplus_{|k|\leq s}H^{\otimes k}\right)$$

(2) It is then routine to check that we have equivalences as follows:
\begin{eqnarray*}
C_{G_C}\subset C
&\iff&E_{C_{G_C}}\subset E_C\\
&\iff&E_{C_{G_C}}^{(s)}\subset E_C^{(s)},\forall s\\
&\iff&E_{C_{G_C}}^{(s)'}\supset E_C^{(s)'},\forall s
\end{eqnarray*}

(3) Summarizing, we would like to prove that we have $E_C^{(s)'}\subset E_{C_{G_C}}^{(s)'}$. But this can be done by doing some abstract algebra, and we refer here to Malacarne \cite{mal}, or to the paper of Woronowicz \cite{wo2}. For more on all this, you have as well my book \cite{ba1}.
\end{proof}

With this piece of general theory in hand, let us go back to Principle 12.19, and develop the second idea there, namely delinearization and Brauer theorems. We have:

\begin{definition}
A category of crossing partitions is a collection $D=\bigsqcup_{k,l}D(k,l)$ of subsets $D(k,l)\subset P(k,l)$, having the following properties:
\begin{enumerate}
\item Stability under the horizontal concatenation, $(\pi,\sigma)\to[\pi\sigma]$.

\item Stability under vertical concatenation $(\pi,\sigma)\to[^\sigma_\pi]$, with matching middle symbols.

\item Stability under the upside-down turning $*$, with switching of colors, $\circ\leftrightarrow\bullet$.

\item Each set $P(k,k)$ contains the identity partition $||\ldots||$.

\item The sets $P(\emptyset,\circ\bullet)$ and $P(\emptyset,\bullet\circ)$ both contain the semicircle $\cap$.

\item The sets $P(k,\bar{k})$ with $|k|=2$ contain the crossing partition $\slash\hskip-2.0mm\backslash$.
\end{enumerate}
\end{definition} 

Observe the similarity with Definition 12.21, and more on this in a moment. In order now to construct a Tannakian category out of such a category, we will need:

\begin{proposition}
Each partition $\pi\in P(k,l)$ produces a linear map
$$T_\pi:(\mathbb C^N)^{\otimes k}\to(\mathbb C^N)^{\otimes l}$$
given by the following formula, with $e_1,\ldots,e_N$ being the standard basis of $\mathbb C^N$,
$$T_\pi(e_{i_1}\otimes\ldots\otimes e_{i_k})=\sum_{j_1\ldots j_l}\delta_\pi\begin{pmatrix}i_1&\ldots&i_k\\ j_1&\ldots&j_l\end{pmatrix}e_{j_1}\otimes\ldots\otimes e_{j_l}$$
and with the Kronecker type symbols $\delta_\pi\in\{0,1\}$ depending on whether the indices fit or not. The assignement $\pi\to T_\pi$ is categorical, in the sense that we have
$$T_\pi\otimes T_\sigma=T_{[\pi\sigma]}\quad,\quad 
T_\pi T_\sigma=N^{c(\pi,\sigma)}T_{[^\sigma_\pi]}\quad,\quad 
T_\pi^*=T_{\pi^*}$$
where $c(\pi,\sigma)$ are certain integers, coming from the erased components in the middle.
\end{proposition}

\begin{proof}
This is something elementary, the computations being as follows:

\medskip

(1) The concatenation axiom follows from the following computation:
\begin{eqnarray*}
&&(T_\pi\otimes T_\sigma)(e_{i_1}\otimes\ldots\otimes e_{i_p}\otimes e_{k_1}\otimes\ldots\otimes e_{k_r})\\
&=&\sum_{j_1\ldots j_q}\sum_{l_1\ldots l_s}\delta_\pi\begin{pmatrix}i_1&\ldots&i_p\\j_1&\ldots&j_q\end{pmatrix}\delta_\sigma\begin{pmatrix}k_1&\ldots&k_r\\l_1&\ldots&l_s\end{pmatrix}e_{j_1}\otimes\ldots\otimes e_{j_q}\otimes e_{l_1}\otimes\ldots\otimes e_{l_s}\\
&=&\sum_{j_1\ldots j_q}\sum_{l_1\ldots l_s}\delta_{[\pi\sigma]}\begin{pmatrix}i_1&\ldots&i_p&k_1&\ldots&k_r\\j_1&\ldots&j_q&l_1&\ldots&l_s\end{pmatrix}e_{j_1}\otimes\ldots\otimes e_{j_q}\otimes e_{l_1}\otimes\ldots\otimes e_{l_s}\\
&=&T_{[\pi\sigma]}(e_{i_1}\otimes\ldots\otimes e_{i_p}\otimes e_{k_1}\otimes\ldots\otimes e_{k_r})
\end{eqnarray*}

(2) The composition axiom follows from the following computation:
\begin{eqnarray*}
&&T_\pi T_\sigma(e_{i_1}\otimes\ldots\otimes e_{i_p})\\
&=&\sum_{j_1\ldots j_q}\delta_\sigma\begin{pmatrix}i_1&\ldots&i_p\\j_1&\ldots&j_q\end{pmatrix}
\sum_{k_1\ldots k_r}\delta_\pi\begin{pmatrix}j_1&\ldots&j_q\\k_1&\ldots&k_r\end{pmatrix}e_{k_1}\otimes\ldots\otimes e_{k_r}\\
&=&\sum_{k_1\ldots k_r}N^{c(\pi,\sigma)}\delta_{[^\sigma_\pi]}\begin{pmatrix}i_1&\ldots&i_p\\k_1&\ldots&k_r\end{pmatrix}e_{k_1}\otimes\ldots\otimes e_{k_r}\\
&=&N^{c(\pi,\sigma)}T_{[^\sigma_\pi]}(e_{i_1}\otimes\ldots\otimes e_{i_p})
\end{eqnarray*}

(3) Finally, the involution axiom follows from the following computation:
\begin{eqnarray*}
&&T_\pi^*(e_{j_1}\otimes\ldots\otimes e_{j_q})\\
&=&\sum_{i_1\ldots i_p}<T_\pi^*(e_{j_1}\otimes\ldots\otimes e_{j_q}),e_{i_1}\otimes\ldots\otimes e_{i_p}>e_{i_1}\otimes\ldots\otimes e_{i_p}\\
&=&\sum_{i_1\ldots i_p}\delta_\pi\begin{pmatrix}i_1&\ldots&i_p\\ j_1&\ldots& j_q\end{pmatrix}e_{i_1}\otimes\ldots\otimes e_{i_p}\\
&=&T_{\pi^*}(e_{j_1}\otimes\ldots\otimes e_{j_q})
\end{eqnarray*}

Summarizing, our correspondence is indeed categorical.
\end{proof}

We can now formulate a key theoretical result, as follows:

\begin{theorem}
Any category of crossing partitions $D\subset P$ produces a series of compact groups $G=(G_N)$, with $G_N\subset U_N$ for any $N\in\mathbb N$, via the formula
$$C_{kl}=span\left(T_\pi\Big|\pi\in D(k,l)\right)$$
for any $k,l$, and Tannakian duality. We call such groups easy.
\end{theorem}

\begin{proof}
Indeed, once we fix an integer $N\in\mathbb N$, the various axioms in Definition 12.25 show, via Proposition 12.26, that the following spaces form a Tannakian category:
$$span\left(T_\pi\Big|\pi\in D(k,l)\right)$$

Thus, Tannakian duality applies, and provides us with a closed subgroup $G_N\subset U_N$ such that the following equalities are satisfied, for any colored integers $k,l$:
$$C_{kl}=span\left(T_\pi\Big|\pi\in D(k,l)\right)$$

Thus, we are led to the conclusion in the statement.
\end{proof}

At the level of basic examples of easy groups, these are the real and complex rotation groups, coming from the following key theorem of Brauer: 

\begin{theorem}
We have the following results:
\begin{enumerate}
\item $U_N$ is easy, coming from the category of all matching pairings $\mathcal P_2$.

\item $O_N$ is easy too, coming from the category of all pairings $P_2$.
\end{enumerate}
\end{theorem}

\begin{proof}
This can be deduced from Tannakian duality, the idea being as follows:

\medskip

(1) The group $U_N$ being defined via the relations $u^*=u^{-1}$, $u^t=\bar{u}^{-1}$, the associated Tannakian category is $C=span(T_\pi|\pi\in D)$, with:
$$D
=<{\ }^{\,\cap}_{\circ\bullet}\,\,,{\ }^{\,\cap}_{\bullet\circ}>
=\mathcal P_2$$

(2) The group $O_N\subset U_N$ being defined by imposing the relations $u_{ij}=\bar{u}_{ij}$, the associated Tannakian category is $C=span(T_\pi|\pi\in D)$, with:
$$D
=<\mathcal P_2,|^{\hskip-1.32mm\circ}_{\hskip-1.32mm\bullet},|_{\hskip-1.32mm\circ}^{\hskip-1.32mm\bullet}>
=P_2$$
  
Thus, we are led to the conclusions in the statement.
\end{proof}

Moving now towards finite groups, we first have the following result:

\begin{theorem}
The symmetric group $S_N$, regarded as group of unitary matrices,
$$S_N\subset O_N\subset U_N$$
via the permutation matrices, is easy, coming from the category of all partitions $P$.
\end{theorem}

\begin{proof}
Consider indeed the group $S_N$, regarded as a group of unitary matrices, with each permutation $\sigma\in S_N$ corresponding to the associated permutation matrix:
$$\sigma(e_i)=e_{\sigma(i)}$$

Consider as well the easy group $G\subset O_N$ coming from the category of all partitions $P$. Since $P$ is generated by the one-block ``fork'' partition $Y\in P(2,1)$, we have:
$$C(G)=C(O_N)\Big/\Big<T_Y\in Hom(u^{\otimes 2},u)\Big>$$

The linear map associated to $Y$ is given by the following formula:
$$T_Y(e_i\otimes e_j)=\delta_{ij}e_i$$

In order to do the computations, we use the following formulae:
$$u=(u_{ij})_{ij}\quad,\quad 
u^{\otimes 2}=(u_{ij}u_{kl})_{ik,jl}\quad,\quad 
T_Y=(\delta_{ijk})_{i,jk}$$

We therefore obtain the following formula:
$$(T_Yu^{\otimes 2})_{i,jk}
=\sum_{lm}(T_Y)_{i,lm}(u^{\otimes 2})_{lm,jk}
=u_{ij}u_{ik}$$

On the other hand, we have as well the following formula:
$$(uT_Y)_{i,jk}
=\sum_lu_{il}(T_Y)_{l,jk}
=\delta_{jk}u_{ij}$$

Thus, the relation defining $G\subset O_N$ reformulates as follows:
$$T_Y\in Hom(u^{\otimes 2},u)\iff u_{ij}u_{ik}=\delta_{jk}u_{ij},\forall i,j,k$$

In other words, the elements $u_{ij}$ must be projections, which must be pairwise orthogonal on the rows of $u=(u_{ij})$. We conclude that $G\subset O_N$ is the subgroup of matrices $g\in O_N$ having the property $g_{ij}\in\{0,1\}$. Thus we have $G=S_N$, as desired. 
\end{proof}

The hyperoctahedral group $H_N$ is easy as well, the result here being as follows:

\begin{theorem}
The hyperoctahedral group $H_N$, regarded as group of matrices,
$$S_N\subset H_N\subset O_N$$
is easy, coming from the category of partitions with even blocks $P_{even}$.
\end{theorem}

\begin{proof}
This follows as usual from Tannakian duality. To be more precise, consider the following one-block partition, which, as the name indicates, looks like a $H$ letter:
$$H\in P(2,2)$$

The linear map associated to this partition is then given by:
$$T_H(e_i\otimes e_j)=\delta_{ij}e_i\otimes e_i$$

By using this formula, we have the following computation:
\begin{eqnarray*}
(T_H\otimes id)u^{\otimes 2}(e_a\otimes e_b)
&=&(T_H\otimes id)\left(\sum_{ijkl}e_{ij}\otimes e_{kl}\otimes u_{ij}u_{kl}\right)(e_a\otimes e_b)\\
&=&(T_H\otimes id)\left(\sum_{ik}e_i\otimes e_k\otimes u_{ia}u_{kb}\right)\\
&=&\sum_ie_i\otimes e_i\otimes u_{ia}u_{ib}
\end{eqnarray*}

On the other hand, we have as well the following computation:
\begin{eqnarray*}
u^{\otimes 2}(T_H\otimes id)(e_a\otimes e_b)
&=&\delta_{ab}\left(\sum_{ijkl}e_{ij}\otimes e_{kl}\otimes u_{ij}u_{kl}\right)(e_a\otimes e_a)\\
&=&\delta_{ab}\sum_{ij}e_i\otimes e_k\otimes u_{ia}u_{ka}
\end{eqnarray*}

We conclude from this that we have the following equivalence:
$$T_H\in End(u^{\otimes 2})\iff \delta_{ik}u_{ia}u_{ib}=\delta_{ab}u_{ia}u_{ka},\forall i,k,a,b$$

But the relations on the right tell us that the entries of $u=(u_{ij})$ must satisfy $\alpha\beta=0$ on each row and column of $u$, and so that the corresponding closed subgroup $G\subset O_N$ consists of the matrices $g\in O_N$ which are permutation-like, with $\pm1$ nonzero entries. Thus, the corresponding group is $G=H_N$, and as a conclusion to this, we have:
$$C(H_N)=C(O_N)\Big/\Big<T_H\in End(u^{\otimes 2})\Big>$$

But this means that the hyperoctahedral group $H_N$ is easy, coming from the category of partitions $D=<H>=P_{even}$. Thus, we are led to the conclusion in the statement.
\end{proof}

More generally now, we have in fact the following result, regarding the series of complex reflection groups $H_N^s$, which covers both the groups $S_N,H_N$:

\index{Brauer theorem}
\index{complex reflection group}

\begin{theorem}
The complex reflection group $H_N^s=\mathbb Z_s\wr S_N$ is easy, the corresponding category $P^s$ consisting of the partitions satisfying the condition
$$\#\circ=\#\bullet(s)$$
as a weighted sum, in each block. In particular, we have the following results:
\begin{enumerate}
\item $S_N$ is easy, coming from the category $P$.

\item $H_N=\mathbb Z_2\wr S_N$ is easy, coming from the category $P_{even}$.

\item $K_N=\mathbb T\wr S_N$ is easy, coming from the category $\mathcal P_{even}$.
\end{enumerate}
\end{theorem}

\begin{proof}
This is something that we already know at $s=1,2$, from Theorems 12.29 and 12.30. In general, the proof is similar, based on Tannakian duality. To be more precise, in what regards the main assertion, the idea here is that the one-block partition $\pi\in P(s)$, which generates the category of partitions $P^s$ in the statement, implements the relations producing the subgroup $H_N^s\subset S_N$. As for the last assertions, these are all elementary:

\medskip

(1) At $s=1$ we know that we have $H_N^1=S_N$. Regarding now the corresponding category, here the condition $\#\circ=\#\bullet(1)$ is automatic, and so $P^1=P$.

\medskip

(2) At $s=2$ we know that we have $H_N^2=H_N$. Regarding now the corresponding category, here the condition $\#\circ=\#\bullet(2)$ reformulates as follows:
$$\#\circ+\,\#\bullet=0(2)$$

Thus each block must have even size, and we obtain, as claimed, $P^2=P_{even}$.

\medskip

(3) At $s=\infty$ we know that we have $H_N^\infty=K_N$. Regarding now the corresponding category, here the condition $\#\circ=\#\bullet(\infty)$ reads:
$$\#\circ=\#\bullet$$

But this is the condition defining $\mathcal P_{even}$, and so $P^\infty=\mathcal P_{even}$, as claimed.
\end{proof}

Summarizing, we have many examples. In fact, our list of easy groups has currently become quite big, and here is a selection of the main results that we have so far: 

\begin{theorem}
We have a diagram of compact groups as follows,
$$\xymatrix@R=50pt@C=50pt{
K_N\ar[r]&U_N\\
H_N\ar[u]\ar[r]&O_N\ar[u]}$$
where $H_N=\mathbb Z_2\wr S_N$ and $K_N=\mathbb T\wr S_N$, and all these groups are easy.
\end{theorem}

\begin{proof}
This follows from the above results. To be more precise, we know that the above groups are all easy, the corresponding categories of partitions being as follows:
$$\xymatrix@R=16mm@C=18mm{
\mathcal P_{even}\ar[d]&\mathcal P_2\ar[l]\ar[d]\\
P_{even}&P_2\ar[l]}$$

Thus, we are led to the conclusion in the statement.
\end{proof}

Summarizing, we have reached to the conclusions formulated in Principle 12.19. All this remains of course a bit away from graph theory, but we will make a good use of what we learned here, later on in this book, especially when talking quantum groups.

\section*{12d. Infinite graphs}

We would like to end this chapter, and this Part III of the present book, with a discussion regarding the infinite graphs. Generally speaking, the subject is heavily analytic, and it is beyond our purposes here to really get into this, at least at this stage of things. However, we have already seen some interesting examples of infinite graphs in the above, and also some of our basic definitions in the above extend in a quite straightforward way to the infinite graph setting, and all this is certainly worth an informal discussion.

\bigskip

To start with, an infinite graph $X$ is the same thing as a finite graph, except for the fact that the set of vertices is infinite, $|X|=\infty$, assumed countable. As in the finite group case, many interesting examples appear as Cayley graphs, as follows:

\begin{definition}
Associated to any discrete group $G=<S>$, with the generating set $S$ assumed to satisfy $1\notin S=S^{-1}$, is its Cayley graph, constructed as follows:
\begin{enumerate}
\item The vertices are the group elements $g\in G$.

\item Edges $g-h$ are drawn when $h=gs$, with $s\in S$.
\end{enumerate} 
\end{definition}

Of course, we can talk as well about oriented Cayley graphs, and about colorings too, exactly as in the finite group case, but this discussion being quite informal anyway, we will focus on the most standard types of Cayley graphs, which are those above.

\bigskip

At the level of basic examples, we have two of them, as follows:

\bigskip

(1) Consider the group $\mathbb Z^N$, that is, the free abelian group on $N$ generators. We can represent the group elements as vectors in $\mathbb R^N$, in the obvious way, by using the embedding $\mathbb Z^N\subset\mathbb R^N$, and if we endow our group with its standard generating set, namely $S=\{\pm1_1,\ldots,\pm1_N\}$, written additively, the edges will appear precisely at the edges of the lattice $\mathbb Z^N\subset\mathbb R^N$. Thus, the Cayley graph that we get is precisely this lattice:
$$
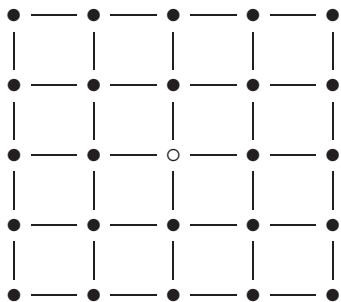
\xymatrix@R=14pt@C=17pt{
\bullet\ar@{-}[d]\ar@{-}[r]&\bullet\ar@{-}[d]\ar@{-}[r]&\bullet\ar@{-}[d]\ar@{-}[r]&\bullet\ar@{-}[d]\ar@{-}[r]&\bullet\ar@{-}[d]\\
\bullet\ar@{-}[d]\ar@{-}[r]&\bullet\ar@{-}[r]\ar@{-}[d]&\bullet\ar@{-}[r]\ar@{-}[d]&\bullet\ar@{-}[r]\ar@{-}[d]&\bullet\ar@{-}[d]\\
\bullet\ar@{-}[d]\ar@{-}[r]&\bullet\ar@{-}[r]\ar@{-}[d]&\circ\ar@{-}[r]\ar@{-}[d]&\bullet\ar@{-}[r]\ar@{-}[d]&\bullet\ar@{-}[d]\\
\bullet\ar@{-}[d]\ar@{-}[r]&\bullet\ar@{-}[r]\ar@{-}[d]&\bullet\ar@{-}[r]\ar@{-}[d]&\bullet\ar@{-}[r]\ar@{-}[d]&\bullet\ar@{-}[d]\\
\bullet\ar@{-}[r]&\bullet\ar@{-}[r]&\bullet\ar@{-}[r]&\bullet\ar@{-}[r]&\bullet}$$

(2) Consider now the free group $F_N=\mathbb Z^{*N}$ on $N$ generators, with its standard generating set, formed by the $N$ generators and their inverses, $S=\{g_1^{\pm1},\ldots,g_N^{\pm1}\}$. As explained in chapter 6, when discussing trees, at $N=2$ the corresponding Cayley graph is as follows, and in general the picture is similar, with valence 4 being replaced by valence $4N$:
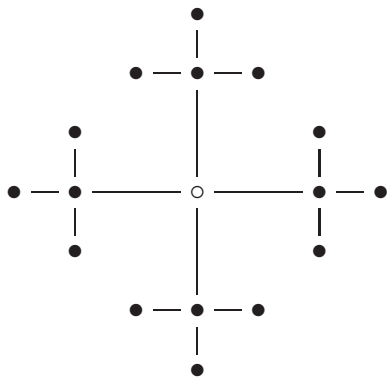
$$\xymatrix@R=10pt@C=10pt{
&&&\bullet\ar@{-}[d]\\
&&\bullet\ar@{-}[r]&\bullet\ar@{-}[dd]\ar@{-}[r]&\bullet\\
&\bullet\ar@{-}[d]&&&&\bullet\ar@{-}[d]\\
\bullet\ar@{-}[r]&\bullet\ar@{-}[rr]&&\circ\ar@{-}[rr]\ar@{-}[dd]&&\bullet\ar@{-}[r]&\bullet\\
&\bullet\ar@{-}[u]&&&&\bullet\ar@{-}[u]\\
&&\bullet\ar@{-}[r]&\bullet\ar@{-}[d]\ar@{-}[r]&\bullet\\
&&&\bullet
}$$

So, these will be our basic examples of infinite graphs. And with the comment that these are both excellent graphs, not only aesthetically, but mathematically too.

\bigskip

Regarding now the general theory for the infinite graphs $X$, we can certainly talk about the corresponding symmetry groups $G(X)\subset S_\infty$, and other algebraic aspects, but the main questions remain those of analytic nature, notably in relation with:

\begin{question}
Given an infinite rooted graph $X$:
\begin{enumerate}
\item What is the number $L_k$ of length $k$ loops based at the root? 

\item What about the probability measure $\mu$ having the numbers $L_k$ as moments?
\end{enumerate}
\end{question}

We refer to chapter 1 for the computation for $\mathbb Z$, that is, for the Cayley graph of $\mathbb Z=F_1$, and to chapter 3 for some further computations, for the graphs $\mathbb N$ and $D_\infty$. There are many interesting questions here, notably with the computation for the Cayley graphs of $\mathbb Z^N$ and $F_N$ at arbitrary $N\in\mathbb N$, leading to a lot of interesting mathematics, and for more on all this, we refer to any advanced probability book.

\bigskip

Finally, another set of interesting questions appears in relation with the Laplace operator introduced in chapter 5. Indeed, as explained there, for the Cayley graph of $\mathbb Z^N$, all this is related to the discretization of the basic equations of physics, via the finite element method. For more on all this, we refer to any good PDE book.

\bigskip

As a last topic of discussion, let us get back to the $N\to\infty$ considerations from the end of chapter 11. It is pretty much clear, from the discussion there, that, although this might seem related to the infinite graphs, the world of infinite graphs is in fact too narrow, for discussing such things. However, as good news, we have something new to say on that subject, by using the Brauer theorems from the previous section, namely:

\begin{theorem}
Given a family of easy groups $G=(G_N)$, coming from a category of crossing partitions $D\subset P$, we have the following formula, with $D_k=D(0,k)$:
$$\lim_{N\to\infty}\int_{G_N}\chi^k=|D_k|$$ 
More generally, we have the following formula, with $|.|$ being the number of blocks:
$$\lim_{N\to\infty}\int_{G_N}\chi_t^k=\sum_{\pi\in D_k}t^{|\pi|}$$ 
In the case of the groups $S=(S_N)$, $H=(H_N)$, and more generally $H^s=(H_N^s)$, we recover in this way, more conceptually, our previous probabilistic results.
\end{theorem} 

\begin{proof}
There are several things going on here, the idea being as follows:

\medskip

(1) To start with, we have the following formula, coming from Peter-Weyl, and then from easiness, with $u$ standing as usual for the fundamental representation:
\begin{eqnarray*}
\int_{G_N}\chi^k
&=&\int_{G_N}\chi_{u^{\otimes k}}\\
&=&\dim(Fix(u^{\otimes k}))\\
&=&\dim\big(span(T_\pi|\pi\in D_k)\big)
\end{eqnarray*}

(2) The point now is that, with $N\to\infty$, the above vectors $T_\pi$ become linearly independent. This is something not exactly trivial, the standard argument here being that it is enough to check this for the biggest possible category, $D=P$, and that for this category, the determinant of the Gram matrix of the vectors $T_\pi$ can be explicitly computed, thanks to a result of Lindst\"om, and follows to be nonzero, with $N\to\infty$.

\medskip

(3) Long story short, we have the first formula in the statement, at $t=1$. As for the general $t>0$ formula, this can be deduced as well, via more technical integration, called Weingarten formula, again by using Peter-Weyl and easiness.

\medskip

(4) Finally, in what regards the examples as the end, for $S=(S_N)$, where $D=P$, and at $t=1$, it is notorious that the measure having as moments the Bell numbers $B_k=|P_k|$ is the Poisson law $p_1$. Thus we have the general $H_N^s$ result at $s=1,t=1$, and the extension to the case of arbitrary parameters $s\in\mathbb N$ and $t>0$ is straightforward.
\end{proof}

And good news, that is all. In the hope that you liked this Part III of the present book, and of course with some apologies for going a bit off-topic on a number of occassions, as for instance with this Theorem 12.35 that we just proved, featuring no graphs. The point, however, with all this, is that all this learning will be very useful, for what comes next.
 
\section*{12e. Exercises}

We had a quite varied chapter here, but at the level of exercises, we would like to insist on the infinite graphs, which are certainly worth a more detailed look:

\begin{exercise}
Draw the Cayley graphs of other infinite groups $G$.
\end{exercise}

\begin{exercise}
Clarify the formalism of symmetry groups $G(X)\subset S_\infty$.
\end{exercise}

\begin{exercise}
Clarify too the partial symmetry groups $\widetilde{G}(X)\subset\widetilde{S}_\infty$.
\end{exercise}

\begin{exercise}
Do the loop computations for the Cayley graph of $\mathbb Z^N$.
\end{exercise}

\begin{exercise}
Do the loop computations for the Cayley graph of $F_N$ too.
\end{exercise}

\begin{exercise}
Do the $F_N$ computations appear as ``liberations'' of those for $\mathbb Z^N$?
\end{exercise}

\begin{exercise}
Study the wave equation on various infinite graphs.
\end{exercise}

\begin{exercise}
Study the heat equation too, on various infinite graphs.
\end{exercise}

As bonus exercise, learn some functional analysis. This would be useful in connection with the above, and also in connection with what we will be doing next, in this book.

\part{Quantum symmetry}

\ \vskip50mm

\begin{center}
{\em Think about the way

That we live today

Think about the way

How some people play}
\end{center}

\chapter{Quantum groups}

\section*{13a. Quantum spaces}

Welcome to quantum symmetry. Our purpose in what follows will be to look for hidden, quantum symmetries of graphs, according to the following principle:

\begin{principle}
The following happen, in the quantum world:
\begin{enumerate}
\item $S_N$ has a free analogue $S_N^+$, which is a compact quantum group.

\item This quantum group $S_N^+$ is infinite, and reminding $SO_3$, at $N\geq4$.

\item $S_N\to S_N^+$ can be however understood, using algebra and probability. 

\item $S_N^+$ is the quantum symmetry group $G^+(K_N)$ of the complete graph $K_N$.

\item In fact, any graph $X\subset K_N$ has a quantum symmetry group $G^+(X)\subset S_N^+$.

\item $G(X)\subset G^+(X)$ can be an isomorphism or not, depending on $X$.

\item $G(X)\to G^+(X)$ can be understood, via algebra and probability.
\end{enumerate}
\end{principle} 

Excited about this? We will learn this technology, in this chapter, and in the next one. To be more precise, in this chapter we will talk about Hilbert spaces, operator algebras, quantum spaces, quantum groups, and then about (1), with a look into (2,3) too. And then, in the next chapter, we will talk about (4,5), with a look into (6,7) too. 

\bigskip

Getting started now, we already know a bit about operator algebras and quantum spaces from chapter 11, but that material was explained in a hurry, time now to do this the right way. At the gates of the quantum world are the Hilbert spaces:

\index{Hilbert space}

\begin{definition}
A Hilbert space is a complex vector space $H$ with a scalar product $<x,y>$, which will be linear at left and antilinear at right,
$$<\lambda x,y>=\lambda<x,y>\quad,\quad <x,\lambda y>=\bar{\lambda}<x,y>$$
and which is complete with respect to corresponding norm
$$||x||=\sqrt{<x,x>}$$
in the sense that any sequence $\{x_n\}$ which is a Cauchy sequence, having the property $||x_n-x_m||\to0$ with $n,m\to\infty$, has a limit, $x_n\to x$.
\end{definition}

Here our convention for the scalar products, written $<x,y>$ and being linear at left, is one among others, often used by mathematicians, and also by certain professional quantum physicists, like myself. As further comments now on Definition 13.2, there is some mathematics encapsulated there, needing some discussion. First, we have:

\begin{theorem}
Given an index set $I$, which can be finite or not, the space
$$l^2(I)=\left\{(x_i)_{i\in I}\Big|\sum_i|x_i|^2<\infty\right\}$$
is a Hilbert space, with scalar product as follows:
$$<x,y>=\sum_ix_i\bar{y}_i$$
When $I$ is finite, $I=\{1,\ldots,N\}$, we obtain in this way the usual space $H=\mathbb C^N$.
\end{theorem}

\begin{proof}
All this is well-known and routine, the idea being as follows:

\medskip

(1) We know that $l^2(I)\subset\mathbb C^I$ is the space of vectors satisfying $||x||<\infty$. We want to prove that $l^2(I)$ is a vector space, that $<x,y>$ is a scalar product on it, that $l^2(I)$ is complete with respect to $||.||$, and finally that for $|I|<\infty$ we have $l^2(I)=\mathbb C^{|I|}$.

\medskip

(2) The last assertion, $l^2(I)=\mathbb C^{|I|}$ for $|I|<\infty$, is clear, because in this case the sums are finite, so the condition $||x||<\infty$ is automatic. So, we know at least one thing.

\medskip

(3) Next, we can use the Cauchy-Schwarz inequality, which is as follows, coming from the positivity of the degree 2 quantity $f(t)=||twx+y||^2$, with $t\in\mathbb R$ and $w\in\mathbb T$:
$$|<x,y>|\leq||x||\cdot||y||$$

(4) Indeed, with Cauchy-Schwarz in hand, everything is straightforward. We first obtain, by raising to the square and expanding, that for any $x,y\in l^2(I)$ we have:
$$||x+y||\leq||x||+||y||$$

(5) Thus $l^2(I)$ is indeed a vector space, and $<x,y>$ is surely a scalar product on it, because all the conditions for a scalar product are trivially satisfied.

\medskip

(6) Finally, the completness with respect to $||.||$ follows in the obvious way, the limit of a Cauchy sequence $\{x^n\}$ being the vector $y=(y_i)$ given by $y_i=\lim_{n\to\infty}x^n_i$.
\end{proof}

Going now a bit abstract, we have, more generally, the following result:

\begin{theorem}
Given a space $X$ with a positive measure $\mu$ on it, the space
$$L^2(X)=\left\{f:X\to\mathbb C\Big|\int_X|f(x)|^2\,d\mu(x)<\infty\right\}$$
with the convention $f=g$ when $||f-g||=0$, is a Hilbert space, with scalar product:
$$<f,g>=\int_Xf(x)\overline{g(x)}\,d\mu(x)$$
When $X=I$ is discrete, $\mu(\{x\})=1$ for any $x\in X$, we obtain the previous space $l^2(I)$.
\end{theorem}

\begin{proof}
This is something routine, remake of Theorem 13.3, as follows:

\medskip

(1) The proof of the first, and main assertion is something perfectly similar to the proof of Theorem 13.3, by replacing everywhere the sums by integrals.

\medskip

(2) As for the last assertion, when $\mu$ is the counting measure all our integrals here become usual sums, and so we recover in this way Theorem 13.3.
\end{proof}

As a third and last theorem about Hilbert spaces, that we will need, we have:

\index{Gram-Schmidt}
\index{separable space}

\begin{theorem}
Any Hilbert space $H$ has an orthonormal basis $\{e_i\}_{i\in I}$, which is by definition a set of vectors whose span is dense in $H$, and which satisfy
$$<e_i,e_j>=\delta_{ij}$$
with $\delta$ being a Kronecker symbol. The cardinality $|I|$ of the index set, which can be finite, countable, or worse, depends only on $H$, and is called dimension of $H$. We have
$$H\simeq l^2(I)$$
in the obvious way, mapping $\sum\lambda_ie_i\to(\lambda_i)$. The Hilbert spaces with $\dim H=|I|$ being countable, including $l^2(\mathbb N)$ and $L^2(\mathbb R)$, are all isomorphic, and are called separable.
\end{theorem}

\begin{proof}
We have many assertions here, the idea being as follows:

\medskip

(1) In finite dimensions an orthonormal basis $\{e_i\}_{i\in I}$ can be constructed by starting with any vector space basis $\{x_i\}_{i\in I}$, and using the Gram-Schmidt procedure. But the same method works in arbitrary dimensions, by using the Zorn lemma.

\medskip

(2) Regarding $L^2(\mathbb R)$, here we can argue that, since functions can be approximated by polynomials, we have a countable algebraic basis, namely $\{x^n\}_{n\in\mathbb N}$, called the Weierstrass basis, that we can orthogonalize afterwards by using Gram-Schmidt.
\end{proof}

Moving ahead, now that we know what our vector spaces are, we can talk about infinite matrices with respect to them. And the situation here is as follows:

\index{bounded operator}

\begin{theorem}
Given a Hilbert space $H$, consider the linear operators $T:H\to H$, and for each such operator define its norm by the following formula:
$$||T||=\sup_{||x||=1}||Tx||$$
The operators which are bounded, $||T||<\infty$, form then a complex algebra $B(H)$, which is complete with respect to $||.||$. When $H$ comes with a basis $\{e_i\}_{i\in I}$, we have
$$B(H)\subset\mathcal L(H)\subset M_I(\mathbb C)$$
where $\mathcal L(H)$ is the algebra of all linear operators $T:H\to H$, and $\mathcal L(H)\subset M_I(\mathbb C)$ is the correspondence $T\to M$ obtained via the usual linear algebra formulae, namely:
$$T(x)=Mx\quad,\quad M_{ij}=<Te_j,e_i>$$
In infinite dimensions, none of the above two inclusions is an equality.
\end{theorem}

\begin{proof}
This is something straightforward, the idea being as follows:

\medskip

(1) The fact that we have indeed an algebra, satisfying the product condition in the statement, follows from the following estimates, which are all elementary:
$$||S+T||\leq||S||+||T||\quad,\quad 
||\lambda T||=|\lambda|\cdot||T||\quad,\quad 
||ST||\leq||S||\cdot||T||$$

(2) Regarding now the completness assertion, if $\{T_n\}\subset B(H)$ is Cauchy then $\{T_nx\}$ is Cauchy for any $x\in H$, so we can define the limit $T=\lim_{n\to\infty}T_n$ by setting:
$$Tx=\lim_{n\to\infty}T_nx$$

Let us first check that the application $x\to Tx$ is linear. We have:
\begin{eqnarray*}
T(x+y)
&=&\lim_{n\to\infty}T_n(x+y)\\
&=&\lim_{n\to\infty}T_n(x)+T_n(y)\\
&=&\lim_{n\to\infty}T_n(x)+\lim_{n\to\infty}T_n(y)\\
&=&T(x)+T(y)
\end{eqnarray*}

Similarly, we have $T(\lambda x)=\lambda T(x)$, and we conclude that $T\in\mathcal L(H)$.

\medskip

(3) With this done, it remains to prove now that we have $T\in B(H)$, and that $T_n\to T$ in norm. For this purpose, observe that we have:
\begin{eqnarray*}
||T_n-T_m||\leq\varepsilon\ ,\ \forall n,m\geq N
&\implies&||T_nx-T_mx||\leq\varepsilon\ ,\ \forall||x||=1\ ,\ \forall n,m\geq N\\
&\implies&||T_nx-Tx||\leq\varepsilon\ ,\ \forall||x||=1\ ,\ \forall n\geq N\\
&\implies&||T_Nx-Tx||\leq\varepsilon\ ,\ \forall||x||=1\\
&\implies&||T_N-T||\leq\varepsilon
\end{eqnarray*}

But this gives both $T\in B(H)$, and $T_N\to T$ in norm, and we are done.

\medskip

(4) Regarding the embeddings, the correspondence $T\to M$ in the statement is indeed linear, and its kernel is $\{0\}$, so we have indeed an embedding as follows, as claimed:
$$\mathcal L(H)\subset M_I(\mathbb C)$$

In finite dimensions we have an isomorphism, because any $M\in M_N(\mathbb C)$ determines an operator $T:\mathbb C^N\to\mathbb C^N$, given by $<Te_j,e_i>=M_{ij}$. However, in infinite dimensions, we have matrices not producing operators, as for instance the all-one matrix. 

\medskip

(5) As for the examples of linear operators which are not bounded, these are more complicated, coming from logic, and we will not need them in what follows.
\end{proof}

Finally, as a second and last result regarding the operators, we will need:

\index{adjoint operator}

\begin{theorem}
Each operator $T\in B(H)$ has an adjoint $T^*\in B(H)$, given by: 
$$<Tx,y>=<x,T^*y>$$
The operation $T\to T^*$ is antilinear, antimultiplicative, involutive, and satisfies:
$$||T||=||T^*||\quad,\quad ||TT^*||=||T||^2$$
When $H$ comes with a basis $\{e_i\}_{i\in I}$, the operation $T\to T^*$ corresponds to
$$(M^*)_{ij}=\overline{M}_{ji}$$ 
at the level of the associated matrices $M\in M_I(\mathbb C)$.
\end{theorem}

\begin{proof}
This is standard too, and can be proved in 3 steps, as follows:

\medskip

(1) The existence of the adjoint operator $T^*$, given by the formula in the statement, comes from the fact that the function $\varphi(x)=<Tx,y>$ being a linear map $H\to\mathbb C$, we must have a formula as follows, for a certain vector $T^*y\in H$:
$$\varphi(x)=<x,T^*y>$$

Moreover, since this vector is unique, $T^*$ is unique too, and we have as well:
$$(S+T)^*=S^*+T^*\quad,\quad
(\lambda T)^*=\bar{\lambda}T^*\quad,\quad 
(ST)^*=T^*S^*\quad,\quad 
(T^*)^*=T$$

Observe also that we have indeed $T^*\in B(H)$, because:
\begin{eqnarray*}
||T||
&=&\sup_{||x||=1}\sup_{||y||=1}<Tx,y>\\
&=&\sup_{||y||=1}\sup_{||x||=1}<x,T^*y>\\
&=&||T^*||
\end{eqnarray*}

(2) Regarding now $||TT^*||=||T||^2$, which is a key formula, observe that we have:
$$||TT^*||
\leq||T||\cdot||T^*||
=||T||^2$$

On the other hand, we have as well the following estimate:
\begin{eqnarray*}
||T||^2
&=&\sup_{||x||=1}|<Tx,Tx>|\\
&=&\sup_{||x||=1}|<x,T^*Tx>|\\
&\leq&||T^*T||
\end{eqnarray*}

By replacing $T\to T^*$ we obtain from this $||T||^2\leq||TT^*||$, as desired.

\medskip

(3) Finally, when $H$ comes with a basis, the formula $<Tx,y>=<x,T^*y>$ applied with $x=e_i$, $y=e_j$ translates into the formula $(M^*)_{ij}=\overline{M}_{ji}$, as desired.
\end{proof}

Generally speaking, the theory of bounded operators can be developed in analogy with the theory of the usual matrices, and the main results can be summarized as follows:

\begin{fact}
The following happen, extending the spectral theorem for matrices:
\begin{enumerate}
\item Any self-adjoint operator, $T=T^*$, is diagonalizable.

\item More generally, any normal operator, $TT^*=T^*T$, is diagonalizable.

\item In fact, any family $\{T_i\}$ of commuting normal operators is diagonalizable.
\end{enumerate}
\end{fact}

You might wonder here, why calling this Fact instead of Theorem. In answer, this is something which is quite hard to prove, and in fact not only we will not prove this, but we will also find a way of short-circuiting all this. But more on this in a moment, for now, let us enjoy this. As a consequence of all this, we can formulate as well:

\begin{fact}
The following happen, regarding the closed $*$-algebras $A\subset B(H)$:
\begin{enumerate}
\item For $A=<T>$ with $T$ normal, we have $A=C(X)$, with $X=\sigma(T)$.

\item In fact, all commutative algebras are of the form $A=C(X)$, with $X$ compact.

\item In general, we can write $A=C(X)$, with $X$ being a compact quantum space.
\end{enumerate}
\end{fact}

To be more precise here, the first assertion is more or less part of the spectral theorems from Fact 13.8, with the spectrum of an operator $T\in B(H)$ being defined as follows:
$$\sigma(T)=\left\{\lambda\in\mathbb C\Big|T-\lambda\notin B(H)^{-1}\right\}$$

Regarding the second assertion, if we write $A=span(T_i)$, then the family $\{T_i\}$ consists of commuting normal operators, and this leads to the above conclusion, with $X$ being a certain compact space associated to the family $\{T_i\}$. As for the third assertion, which is something important to us, this is rather a philosophical conclusion, to all this.

\bigskip

Very good all this, so we have quantum spaces, and you would say, it remains to understand the proofs of all the above, and then we are all set, ready to go ahead with quantum groups, and the rest of our program. However, there is a bug with all this:

\begin{bug}
Besides the spectral theorem in infinite dimensions being something tough, the resulting notion of compact quantum spaces is not very satisfactory, because we cannot define operator algebras $A\subset B(H)$ with generators and relations, as we would love to.
\end{bug}

In short, nice try with the above, but time now to forget all this, and invent something better. And, fortunately, the solution to our problem exists, due to Gelfand, with the starting definition here, that we already met in chapter 11, being as follows:

\begin{definition}
A $C^*$-algebra is a complex algebra $A$, having a norm $||.||$ making it a Banach algebra, and an involution $*$, related to the norm by the formula 
$$||aa^*||=||a||^2$$
which must hold for any $a\in A$.
\end{definition}

As a basic example, the full operator algebra $B(H)$ is a $C^*$-algebra, and so is any norm closed $*$-subalgebra $A\subset B(H)$. We will see in a moment that a converse of this holds, in the sense that any $C^*$-algebra appears as an operator algebra, $A\subset B(H)$.

\bigskip

But, let us start with finite dimensions. We know that the matrix algebra $M_N(\mathbb C)$ is a $C^*$-algebra, with the usual matrix norm and involution of matrices, namely:
$$||M||=\sup_{||x||=1}||Mx||\quad,\quad 
(M^*)_{ij}=\bar{M}_{ji}$$

More generally, any $*$-subalgebra $A\subset M_N(\mathbb C)$ is automatically closed, and so is a $C^*$-algebra. In fact, in finite dimensions, the situation is as follows:

\index{finite dimensional algebra}
\index{multimatrix algebra}

\begin{theorem}
The finite dimensional $C^*$-algebras are exactly the algebras
$$A=M_{N_1}(\mathbb C)\oplus\ldots\oplus M_{N_k}(\mathbb C)$$
with norm $||(a_1,\ldots,a_k)||=\sup_i||a_i||$, and involution $(a_1,\ldots,a_k)^*=(a_1^*,\ldots,a_k^*)$.
\end{theorem}

\begin{proof}
This is something that we discussed in chapter 11, the idea being that this comes by splitting the unit of our algebra $A$ as a sum of central minimal projections, $1=p_1+\ldots+p_k$. Indeed, when doing so, each of the $*$-algebras $A_i=p_iAp_i$ follows to be a matrix algebra, $A_i\simeq M_{N_i}(\mathbb C)$, and this gives the decomposition in the statement.
\end{proof}

In order to develop more theory, we will need a technical result, as follows:

\begin{theorem}
Given an element $a\in A$ of a $C^*$-algebra, define its spectrum as:
$$\sigma(a)=\left\{\lambda\in\mathbb C\Big|a-\lambda\notin A^{-1}\right\}$$
The following spectral theory results hold, exactly as in the $A=B(H)$ case:
\begin{enumerate}
\item We have $\sigma(ab)\cup\{0\}=\sigma(ba)\cup\{0\}$.

\item We have $\sigma(f(a))=f(\sigma(a))$, for any $f\in\mathbb C(X)$ having poles outside $\sigma(a)$.

\item The spectrum $\sigma(a)$ is compact, non-empty, and contained in $D_0(||a||)$.

\item The spectra of unitaries $(u^*=u^{-1})$ and self-adjoints $(a=a^*)$ are on $\mathbb T,\mathbb R$.

\item The spectral radius of normal elements $(aa^*=a^*a)$ is given by $\rho(a)=||a||$.
\end{enumerate}
In addition, assuming $a\in A\subset B$, the spectra of $a$ with respect to $A$ and to $B$ coincide.
\end{theorem}

\begin{proof}
All the above assertions, which are of course formulated a bit informally, are well-known to hold for the full operator algebra $A=B(H)$, and the proof in general is similar. We refer here to chapter 11, where all this was already discussed.
\end{proof}

With these ingredients, we can now a prove a key result of Gelfand, as follows:

\index{Gelfand theorem}
\index{commutative algebra}
\index{compact space}

\begin{theorem}
Any commutative $C^*$-algebra $A$ is of the form
$$A=C(X)$$
with $X=Spec(A)$ being the space of Banach algebra characters $\chi:A\to\mathbb C$.
\end{theorem}

\begin{proof}
This is something that we know too from chapter 11, the idea being that with $X$ as in the statement, we have a morphism of algebras as follows:
$$ev:A\to C(X)\quad,\quad a\to ev_a=[\chi\to\chi(a)]$$

(1) Quite suprisingly, the fact that $ev$ is involutive is not trivial. But here we can argue that it is enough to prove that we have $ev_{a^*}=ev_a^*$ for the self-adjoint elements $a$, which in turn follows from Theorem 13.13 (4), which shows that we have:
$$ev_a(\chi)=\chi(a)\in\sigma(a)\subset\mathbb R$$

(2) Next, since $A$ is commutative, each element is normal, so $ev$ is isometric:
$$||ev_a||=\rho(a)=||a||$$

It remains to prove that $ev$ is surjective. But this follows from the Stone-Weierstrass theorem, because $ev(A)$ is a closed subalgebra of $C(X)$, which separates the points.
\end{proof}

In view of Theorem 13.14, we can formulate the following definition:

\index{quantum space}
\index{compact quantum space}
\index{noncommutative space}

\begin{definition}
Given an arbitrary $C^*$-algebra $A$, we can write 
$$A=C(X)$$
and call the abstract space $X$ a compact quantum space.
\end{definition}

In other words, we can define the category of compact quantum spaces $X$ as being the category of the $C^*$-algebras $A$, with the arrows reversed. A morphism $f:X\to Y$ corresponds by definition to a morphism $\Phi:C(Y)\to C(X)$, a product of spaces $X\times Y$ corresponds by definition to a product of algebras $C(X)\otimes C(Y)$, and so on.

\bigskip

All this is of course a bit speculative, and as a first true result, we have:

\index{finite quantum space}

\begin{theorem}
The finite quantum spaces are exactly the disjoint unions of type
$$X=M_{N_1}\sqcup\ldots\sqcup M_{N_k}$$
where $M_N$ is the finite quantum space given by $C(M_N)=M_N(\mathbb C)$.
\end{theorem}

\begin{proof}
For a compact quantum space $X$, coming from a $C^*$-algebra $A$ via the formula $A=C(X)$, being finite can only mean that the following number is finite:
$$|X|=\dim_\mathbb CA<\infty$$

Thus, by using Theorem 13.12, we are led to the conclusion that we must have:
$$C(X)=M_{N_1}(\mathbb C)\oplus\ldots\oplus M_{N_k}(\mathbb C)$$

But since direct sums of algebras $A$ correspond to disjoint unions of quantum spaces $X$, via the correspondence $A=C(X)$, this leads to the conclusion in the statement.
\end{proof}

Finally, at the general level, we have as well the following key result:

\begin{theorem}
Any $C^*$-algebra appears as an operator algebra:
$$A\subset B(H)$$
Moreover, when $A$ is separable, which is usually the case, $H$ can be taken separable.
\end{theorem}

\begin{proof}
This result, called GNS representation theorem after Gelfand-Naimark-Segal, comes as a continuation of the Gelfand theorem, the idea being as follows:

\medskip

(1) Let us first prove that the result holds in the commutative case, $A=C(X)$. Here, we can pick a positive measure on $X$, and construct our embedding as follows:
$$C(X)\subset B(L^2(X))\quad,\quad f\to[g\to fg]$$

(2) In general the proof is similar, the idea being that given a $C^*$-algebra $A$ we can construct a Hilbert space $H=L^2(A)$, and then an embedding as above:
$$A\subset B(L^2(A))\quad,\quad a\to[b\to ab]$$

(3) Finally, the last assertion is clear, because when $A$ is separable, meaning that it has a countable algebraic basis, so does the associated Hilbert space $H=L^2(A)$.
\end{proof}

Many other things can be said about the $C^*$-algebras, and we recommend here any operator algebra book. But for our purposes here, the above will do.

\section*{13b. Quantum groups}

We are ready now to introduce the quantum groups. The axioms here, due to Woronowicz \cite{wo1}, and slightly modified for our purposes, are as follows:

\index{Woronowicz algebra}
\index{comultiplication}
\index{counit}
\index{antipode}

\begin{definition}
A Woronowicz algebra is a $C^*$-algebra $A$, given with a unitary matrix $u\in M_N(A)$ whose coefficients generate $A$, such that the formulae
$$\Delta(u_{ij})=\sum_ku_{ik}\otimes u_{kj}\quad,\quad
\varepsilon(u_{ij})=\delta_{ij}\quad,\quad 
S(u_{ij})=u_{ji}^*$$
define morphisms of $C^*$-algebras $\Delta:A\to A\otimes A$, $\varepsilon:A\to\mathbb C$ and $S:A\to A^{opp}$, called comultiplication, counit and antipode. 
\end{definition}

Here the tensor product needed for $\Delta$ can be any $C^*$-algebra tensor product, and more on this later. In order to get rid of redundancies, coming from this and from amenability issues, we will divide everything by an equivalence relation, as follows:

\index{equivalence relation}
\index{amenability}

\begin{definition}
We agree to identify two Woronowicz algebras, $(A,u)=(B,v)$, when we have an isomorphism of $*$-algebras
$$<u_{ij}>\simeq<v_{ij}>$$
mapping standard coordinates to standard coordinates, $u_{ij}\to v_{ij}$.
\end{definition}

We say that $A$ is cocommutative when $\Sigma\Delta=\Delta$, where $\Sigma(a\otimes b)=b\otimes a$ is the flip. We have then the following key result, from \cite{wo1}, providing us with examples:

\index{cocommutative algebra}
\index{commutative algebra}
\index{compact Lie group}
\index{finitely generated group}

\begin{theorem}
The following are Woronowicz algebras, which are commutative, respectively cocommutative:
\begin{enumerate}
\item $C(G)$, with $G\subset U_N$ compact Lie group. Here the structural maps are:
$$\Delta(\varphi)=\big[(g,h)\to \varphi(gh)\big]\quad,\quad 
\varepsilon(\varphi)=\varphi(1)\quad,\quad
S(\varphi)=\big[g\to\varphi(g^{-1})\big]$$

\item $C^*(\Gamma)$, with $F_N\to\Gamma$ finitely generated group. Here the structural maps are:
$$\Delta(g)=g\otimes g\quad,\quad 
\varepsilon(g)=1\quad,\quad
S(g)=g^{-1}$$
\end{enumerate}
Moreover, we obtain in this way all the commutative/cocommutative algebras.
\end{theorem}

\begin{proof}
In both cases, we first have to exhibit a certain matrix $u$, and then prove that we have indeed a Woronowicz algebra. The constructions are as follows:

\medskip

(1) For the first assertion, we can use the matrix $u=(u_{ij})$ formed by the standard matrix coordinates of $G$, which is by definition given by:
$$g=\begin{pmatrix}
u_{11}(g)&\ldots&u_{1N}(g)\\
\vdots&&\vdots\\
u_{N1}(g)&\ldots&u_{NN}(g)
\end{pmatrix}$$

(2) For the second assertion, we can use the diagonal matrix formed by generators:
$$u=\begin{pmatrix}
g_1&&0\\
&\ddots&\\
0&&g_N
\end{pmatrix}$$

Finally, regarding the last assertion, in the commutative case this follows from the Gelfand theorem, and in the cocommutative case, we will be back to this.
\end{proof}

In order to get now to quantum groups, we will need as well:

\index{Pontrjagin duality}
\index{dual of group}
\index{Fourier transform}

\begin{proposition}
Assuming that $G\subset U_N$ is abelian, we have an identification of Woronowicz algebras $C(G)=C^*(\Gamma)$, with $\Gamma$ being the Pontrjagin dual of $G$:
$$\Gamma=\big\{\chi:G\to\mathbb T\big\}$$
Conversely, assuming that $F_N\to\Gamma$ is abelian, we have an identification of Woronowicz algebras $C^*(\Gamma)=C(G)$, with $G$ being the Pontrjagin dual of $\Gamma$:
$$G=\big\{\chi:\Gamma\to\mathbb T\big\}$$
Thus, the Woronowicz algebras which are both commutative and cocommutative are exactly those of type $A=C(G)=C^*(\Gamma)$, with $G,\Gamma$ being abelian, in Pontrjagin duality.
\end{proposition}

\begin{proof}
This follows from the Gelfand theorem applied to $C^*(\Gamma)$, and from the fact that the characters of a group algebra come from the characters of the group.
\end{proof}

In view of this result, and of the findings from Theorem 13.20 too, we have the following definition, complementing Definition 13.18 and Definition 13.19:

\index{compact quantum group}
\index{discrete quantum group}
\index{quantum group}
\index{dual group}
\index{dual quantum group}

\begin{definition}
Given a Woronowicz algebra, we write it as follows, and call $G$ a compact quantum Lie group, and $\Gamma$ a finitely generated discrete quantum group:
$$A=C(G)=C^*(\Gamma)$$
Also, we say that $G,\Gamma$ are dual to each other, and write $G=\widehat{\Gamma},\Gamma=\widehat{G}$.
\end{definition}

Let us discuss now some tools for studying the Woronowicz algebras, and the underlying quantum groups. First, we have the following result:

\index{square of antipode}
\index{Hopf algebra axioms}

\begin{proposition}
Let $(A,u)$ be a Woronowicz algebra.
\begin{enumerate} 
\item $\Delta,\varepsilon$ satisfy the usual axioms for a comultiplication and a counit, namely:
$$(\Delta\otimes id)\Delta=(id\otimes \Delta)\Delta$$
$$(\varepsilon\otimes id)\Delta=(id\otimes\varepsilon)\Delta=id$$

\item $S$ satisfies the antipode axiom, on the $*$-algebra generated by entries of $u$: 
$$m(S\otimes id)\Delta=m(id\otimes S)\Delta=\varepsilon(.)1$$

\item In addition, the square of the antipode is the identity, $S^2=id$.
\end{enumerate}
\end{proposition}

\begin{proof}
As a first observation, the result holds in the commutative case, $A=C(G)$ with $G\subset U_N$. Indeed, here we know from Theorem 13.20 that $\Delta,\varepsilon,S$ appear as functional analytic transposes of the multiplication, unit and inverse maps $m,u,i$:
$$\Delta=m^t\quad,\quad 
\varepsilon=u^t\quad,\quad 
S=i^t$$

Thus, in this case, the various conditions in the statement on $\Delta,\varepsilon,S$ simply come by transposition from the group axioms satisfied by $m,u,i$, namely:
$$m(m\times id)=m(id\times m)$$
$$m(u\times id)=m(id\times u)=id$$
$$m(i\times id)\delta=m(id\times i)\delta=1$$

Here $\delta(g)=(g,g)$. Observe also that the result holds as well in the cocommutative case, $A=C^*(\Gamma)$ with $F_N\to\Gamma$, trivially. In general now, the first axiom follows from:
$$(\Delta\otimes id)\Delta(u_{ij})=(id\otimes \Delta)\Delta(u_{ij})=\sum_{kl}u_{ik}\otimes u_{kl}\otimes u_{lj}$$

As for the other axioms, the verifications here are similar.
\end{proof}

In order to reach now to more advanced results, the idea will be that of doing representation theory. Following Woronowicz \cite{wo1}, let us start with the following definition:

\index{representation}
\index{corepresentation}

\begin{definition}
Given $(A,u)$, we call corepresentation of it any unitary matrix $v\in M_n(\mathcal A)$, with $\mathcal A=<u_{ij}>$, satisfying the same conditions as $u$, namely:
$$\Delta(v_{ij})=\sum_kv_{ik}\otimes v_{kj}\quad,\quad
\varepsilon(v_{ij})=\delta_{ij}\quad,\quad
S(v_{ij})=v_{ji}^*$$
We also say that $v$ is a representation of the underlying compact quantum group $G$.
\end{definition}

In the commutative case, $A=C(G)$ with $G\subset U_N$, we obtain in this way the finite dimensional unitary smooth representations $v:G\to U_n$, via the following formula:
$$v(g)=\begin{pmatrix}
v_{11}(g)&\ldots&v_{1n}(g)\\
\vdots&&\vdots\\
v_{n1}(g)&\ldots&v_{nn}(g)
\end{pmatrix}$$

With this convention, we have the following fundamental result, from \cite{wo1}:

\begin{theorem}
Any Woronowicz algebra has a unique Haar integration functional, 
$$\left(\int_G\otimes id\right)\Delta=\left(id\otimes\int_G\right)\Delta=\int_G(.)1$$
which can be constructed by starting with any faithful positive form $\varphi\in A^*$, and setting
$$\int_G=\lim_{n\to\infty}\frac{1}{n}\sum_{k=1}^n\varphi^{*k}$$
where $\phi*\psi=(\phi\otimes\psi)\Delta$. Moreover, for any corepresentation $v\in M_n(\mathbb C)\otimes A$ we have
$$\left(id\otimes\int_G\right)v=P$$
where $P$ is the orthogonal projection onto $Fix(v)=\{\xi\in\mathbb C^n|v\xi=\xi\}$.
\end{theorem}

\begin{proof}
Following \cite{wo1}, this can be done in 3 steps, as follows:

\medskip

(1) Given $\varphi\in A^*$, our claim is that the following limit converges, for any $a\in A$:
$$\int_\varphi a=\lim_{n\to\infty}\frac{1}{n}\sum_{k=1}^n\varphi^{*k}(a)$$

Indeed, by linearity we can assume that $a\in A$ is the coefficient of certain corepresentation, $a=(\tau\otimes id)v$. But in this case, an elementary computation gives the following formula, with $P_\varphi$ being the orthogonal projection onto the $1$-eigenspace of $(id\otimes\varphi)v$:
$$\left(id\otimes\int_\varphi\right)v=P_\varphi$$

(2) Since $v\xi=\xi$ implies $[(id\otimes\varphi)v]\xi=\xi$, we have $P_\varphi\geq P$, where $P$ is the orthogonal projection onto the fixed point space in the statement, namely:
$$Fix(v)=\left\{\xi\in\mathbb C^n\Big|v\xi=\xi\right\}$$

The point now is that when $\varphi\in A^*$ is faithful, by using a standard positivity trick, we can prove that we have $P_\varphi=P$, exactly as in the classical case.

\medskip

(3) With the above formula in hand, the left and right invariance of $\int_G=\int_\varphi$ is clear on coefficients, and so in general, and this gives all the assertions. See \cite{wo1}. 
\end{proof}

We can now develop, again following \cite{wo1}, the Peter-Weyl theory for the corepresentations of $A$. Consider the dense subalgebra $\mathcal A\subset A$ generated by the coefficients of the fundamental corepresentation $u$, and endow it with the following scalar product: 
$$<a,b>=\int_Gab^*$$

With this convention, we have the following result, also from \cite{wo1}:

\index{Peter-Weyl theorem}
\index{Frobenius isomorphism}
\index{finite dimensional algebra}
\index{character}
\index{orthonormal system}

\begin{theorem}
We have the following Peter-Weyl type results:
\begin{enumerate}
\item Any corepresentation decomposes as a sum of irreducible corepresentations.

\item Each irreducible corepresentation appears inside a certain $u^{\otimes k}$.

\item $\mathcal A=\bigoplus_{v\in Irr(A)}M_{\dim(v)}(\mathbb C)$, the summands being pairwise orthogonal.

\item The characters of irreducible corepresentations form an orthonormal system.
\end{enumerate}
\end{theorem}

\begin{proof}
This is something that we met in chapters 11-12, in the case where $G\subset U_N$ is a finite group, or more generally a compact group. In general, when $G$ is a compact quantum group, the proof is quite similar, using Theorem 13.25.
\end{proof}

Finally, no discussion about compact and discrete quantum groups would be complete without a word on amenability. The result here, again from \cite{wo1}, is as follows:

\begin{theorem}
Let $A_{full}$ be the enveloping $C^*$-algebra of $\mathcal A$, and $A_{red}$ be the quotient of $A$ by the null ideal of the Haar integration. The following are then equivalent:
\begin{enumerate}
\item The Haar functional of $A_{full}$ is faithful.

\item The projection map $A_{full}\to A_{red}$ is an isomorphism.

\item The counit map $\varepsilon:A_{full}\to\mathbb C$ factorizes through $A_{red}$.

\item We have $N\in\sigma(Re(\chi_u))$, the spectrum being taken inside $A_{red}$.
\end{enumerate}
If this is the case, we say that the underlying discrete quantum group $\Gamma$ is amenable.
\end{theorem}

\begin{proof}
This is well-known in the group dual case, $A=C^*(\Gamma)$, with $\Gamma$ being a usual discrete group. In general, the result follows by adapting the group dual case proof:

\medskip

$(1)\iff(2)$ This simply follows from the fact that the GNS construction for the algebra $A_{full}$ with respect to the Haar functional produces the algebra $A_{red}$.

\medskip

$(2)\iff(3)$ Here $\implies$ is trivial, and conversely, a counit $\varepsilon:A_{red}\to\mathbb C$ produces an isomorphism $\Phi:A_{red}\to A_{full}$, by slicing the map $\widetilde{\Delta}:A_{red}\to A_{red}\otimes A_{full}$.

\medskip

$(3)\iff(4)$ Here $\implies$ is clear, coming from $\varepsilon(N-Re(\chi (u)))=0$, and the converse can be proved by doing some functional analysis. See \cite{wo1}.
\end{proof}

This was for the basic theory of the quantum groups in the sense of Woronowicz, quickly explained. For more on all this, we have for instance my book \cite{ba3}.

\section*{13c. Quantum permutations}

Following Wang \cite{wan}, let us discuss now the construction and basic properties of the quantum permutation group $S_N^+$. Let us first look at $S_N$. We have here:

\index{magic matrix}
\index{magic unitary}
\index{symmetric group}
\index{Gelfand theorem}
\index{commutative algebra}

\begin{theorem}
The algebra of functions on $S_N$ has the following presentation,
$$C(S_N)=C^*_{comm}\left((u_{ij})_{i,j=1,\ldots,N}\Big|u={\rm magic}\right)$$
and the multiplication, unit and inversion map of $S_N$ appear from the maps
$$\Delta(u_{ij})=\sum_ku_{ik}\otimes u_{kj}\quad,\quad 
\varepsilon(u_{ij})=\delta_{ij}\quad,\quad 
S(u_{ij})=u_{ji}$$
defined at the algebraic level, of functions on $S_N$, by transposing.
\end{theorem}

\begin{proof}
This is something that we know from chapter 11, coming from the Gelfand theorem, applied to the universal algebra in the statement. Indeed, that algebra follows to be of the form $A=C(X)$, with $X$ being a certain compact space. Now since we have coordinates $u_{ij}:X\to\mathbb R$, we have an embedding $X\subset M_N(\mathbb R)$. Also, since we know that these coordinates form a magic matrix, the elements $g\in X$ must be 0-1 matrices, having exactly one 1 entry on each row and each column, and so $X=S_N$, as desired.
\end{proof}

Following now Wang \cite{wan}, we can liberate $S_N$, as follows:

\index{quantum permutation group}
\index{magic unitary}
\index{free symmetric group}

\begin{theorem}
The following universal $C^*$-algebra, with magic meaning as usual formed by projections $(p^2=p^*=p)$, summing up to $1$ on each row and each column,
$$C(S_N^+)=C^*\left((u_{ij})_{i,j=1,\ldots,N}\Big|u={\rm magic}\right)$$
is a Woronowicz algebra, with comultiplication, counit and antipode given by:
$$\Delta(u_{ij})=\sum_ku_{ik}\otimes u_{kj}\quad,\quad 
\varepsilon(u_{ij})=\delta_{ij}\quad,\quad 
S(u_{ij})=u_{ji}$$
Thus the space $S_N^+$ is a compact quantum group, called quantum permutation group.
\end{theorem}

\begin{proof}
As a first observation, the universal $C^*$-algebra in the statement is indeed well-defined, because the conditions $p^2=p^*=p$ satisfied by the coordinates give:
$$||u_{ij}||\leq1$$

In order to prove now that we have a Woronowicz algebra, we must construct maps $\Delta,\varepsilon,S$ given by the formulae in the statement. Consider the following matrices:
$$u^\Delta_{ij}=\sum_ku_{ik}\otimes u_{kj}\quad,\quad 
u^\varepsilon_{ij}=\delta_{ij}\quad,\quad 
u^S_{ij}=u_{ji}$$

Our claim is that, since $u$ is magic, so are these three matrices. Indeed, regarding $u^\Delta$, its entries are idempotents, as shown by the following computation:
$$(u_{ij}^\Delta)^2
=\sum_{kl}u_{ik}u_{il}\otimes u_{kj}u_{lj}
=\sum_{kl}\delta_{kl}u_{ik}\otimes\delta_{kl}u_{kj}
=u_{ij}^\Delta$$

These elements are self-adjoint as well, as shown by the following computation:
$$(u_{ij}^\Delta)^*
=\sum_ku_{ik}^*\otimes u_{kj}^*
=\sum_ku_{ik}\otimes u_{kj}
=u_{ij}^\Delta$$

The row and column sums for the matrix $u^\Delta$ can be computed as follows:
$$\sum_ju_{ij}^\Delta
=\sum_{jk}u_{ik}\otimes u_{kj}
=\sum_ku_{ik}\otimes 1
=1$$
$$\sum_iu_{ij}^\Delta
=\sum_{ik}u_{ik}\otimes u_{kj}
=\sum_k1\otimes u_{kj}
=1$$

Thus, $u^\Delta$ is magic. Regarding now $u^\varepsilon,u^S$, these matrices are magic too, and this for obvious reasons. Thus, all our three matrices $u^\Delta,u^\varepsilon,u^S$ are magic, so we can define $\Delta,\varepsilon,S$ by the formulae in the statement, by using the universality property of $C(S_N^+)$.
\end{proof}

Our first task now is to make sure that Theorem 13.29 produces indeed a new quantum group, which does not collapse to $S_N$. Following Wang \cite{wan}, we have:

\index{liberation}

\begin{theorem}
We have an embedding $S_N\subset S_N^+$, given at the algebra level by: 
$$u_{ij}\to\chi\left(\sigma\in S_N\Big|\sigma(j)=i\right)$$
This is an isomorphism at $N\leq3$, but not at $N\geq4$, where $S_N^+$ is not classical, nor finite.
\end{theorem} 

\begin{proof}
The fact that we have indeed an embedding as above follows from Theorem 13.28. Observe that in fact more is true, because Theorems 13.28 and 13.29 give:
$$C(S_N)=C(S_N^+)\Big/\Big<ab=ba\Big>$$

Thus, the inclusion $S_N\subset S_N^+$ is a ``liberation'', in the sense that $S_N$ is the classical version of $S_N^+$. We will often use this basic fact, in what follows. Regarding now the second assertion, we can prove this in four steps, as follows:

\medskip

\underline{Case $N=2$}. The fact that $S_2^+$ is indeed classical, and hence collapses to $S_2$, is trivial, because the $2\times2$ magic matrices are as follows, with $p$ being a projection:
$$U=\begin{pmatrix}p&1-p\\1-p&p\end{pmatrix}$$

Thus $C(S_2^+)$ is commutative, and equals its biggest commutative quotient, $C(S_2)$.

\medskip

\underline{Case $N=3$}. It is enough to check that $u_{11},u_{22}$ commute. But this follows from:
\begin{eqnarray*}
u_{11}u_{22}
&=&u_{11}u_{22}(u_{11}+u_{12}+u_{13})\\
&=&u_{11}u_{22}u_{11}+u_{11}u_{22}u_{13}\\
&=&u_{11}u_{22}u_{11}+u_{11}(1-u_{21}-u_{23})u_{13}\\
&=&u_{11}u_{22}u_{11}
\end{eqnarray*}

Indeed, by conjugating, $u_{22}u_{11}=u_{11}u_{22}u_{11}$, so $u_{11}u_{22}=u_{22}u_{11}$, as desired.

\medskip

\underline{Case $N=4$}. Consider the following matrix, with $p,q$ being projections:
$$U=\begin{pmatrix}
p&1-p&0&0\\
1-p&p&0&0\\
0&0&q&1-q\\
0&0&1-q&q
\end{pmatrix}$$ 

This matrix is magic, and we can choose $p,q\in B(H)$ as for the algebra $<p,q>$ to be noncommutative and infinite dimensional. We conclude that $C(S_4^+)$ is noncommutative and infinite dimensional as well, and so $S_4^+$ is non-classical and infinite, as claimed.

\medskip

\underline{Case $N\geq5$}. Here we can use the standard embedding $S_4^+\subset S_N^+$, obtained at the level of the corresponding magic matrices in the following way:
$$u\to\begin{pmatrix}u&0\\ 0&1_{N-4}\end{pmatrix}$$

Indeed, with this in hand, the fact that $S_4^+$ is a non-classical, infinite compact quantum group implies that $S_N^+$ with $N\geq5$ has these two properties as well.
\end{proof}

As a first observation, as a matter of doublechecking our findings, we are not wrong with our formalism, because as explained once again in \cite{wan}, we have as well:

\index{quantum permutation}
\index{counting measure}
\index{coaction}
\index{standard coaction}

\begin{theorem}
The quantum permutation group $S_N^+$ acts on the set $X=\{1,\ldots,N\}$, the corresponding coaction map $\Phi:C(X)\to C(X)\otimes C(S_N^+)$ being given by:
$$\Phi(e_i)=\sum_je_j\otimes u_{ji}$$
In fact, $S_N^+$ is the biggest compact quantum group acting on $X$, by leaving the counting measure invariant, in the sense that $(tr\otimes id)\Phi=tr(.)1$, where $tr(e_i)=\frac{1}{N},\forall i$.
\end{theorem}

\begin{proof}
Our claim is that given a compact matrix quantum group $G$, the following formula defines a morphism of algebras, which is a coaction map, leaving the trace invariant, precisely when the matrix $u=(u_{ij})$ is a magic corepresentation of $C(G)$: 
$$\Phi(e_i)=\sum_je_j\otimes u_{ji}$$

Indeed, let us first determine when $\Phi$ is multiplicative. We have:
$$\Phi(e_i)\Phi(e_k)
=\sum_{jl}e_je_l\otimes u_{ji}u_{lk}
=\sum_je_j\otimes u_{ji}u_{jk}$$
$$\Phi(e_ie_k)
=\delta_{ik}\Phi(e_i)
=\delta_{ik}\sum_je_j\otimes u_{ji}$$

We conclude that the multiplicativity of $\Phi$ is equivalent to the following conditions:
$$u_{ji}u_{jk}=\delta_{ik}u_{ji}\quad,\quad\forall i,j,k$$

Similarly, $\Phi$ is unital when $\sum_iu_{ji}=1$, $\forall j$. Finally, the fact that $\Phi$ is a $*$-morphism translates into $u_{ij}=u_{ij}^*$, $\forall i,j$. Summing up, in order for $\Phi(e_i)=\sum_je_j\otimes u_{ji}$ to be a morphism of $C^*$-algebras, the elements $u_{ij}$ must be projections, summing up to 1 on each row of $u$. Regarding now the preservation of the trace, observe that we have:
$$(tr\otimes id)\Phi(e_i)=\frac{1}{N}\sum_ju_{ji}$$

Thus the trace is preserved precisely when the elements $u_{ij}$ sum up to 1 on each of the columns of $u$. We conclude from this that $\Phi(e_i)=\sum_je_j\otimes u_{ji}$ is a morphism of $C^*$-algebras preserving the trace precisely when $u$ is magic, and this gives the result.
\end{proof}

\section*{13d. Liberation theory} 

In order to study $S_N^+$, and better understand the liberation operation $S_N\to S_N^+$, we can use representation theory. We have the following version of Tannakian duality:

\index{Tannakian duality}
\index{soft Tannakian duality}

\begin{theorem}
The following operations are inverse to each other:
\begin{enumerate}
\item The construction $A\to C$, which associates to any Woronowicz algebra $A$ the tensor category formed by the intertwiner spaces $C_{kl}=Hom(u^{\otimes k},u^{\otimes l})$.

\item The construction $C\to A$, which associates to a tensor category $C$ the Woronowicz algebra $A$ presented by the relations $T\in Hom(u^{\otimes k},u^{\otimes l})$, with $T\in C_{kl}$.
\end{enumerate}
\end{theorem}

\begin{proof}
This is something quite deep, going back to Woronowicz's paper \cite{wo2} in a slightly different form, with the idea being as follows:

\medskip

-- We have indeed a construction $A\to C$ as above, whose output is a tensor $C^*$-subcategory with duals of the tensor $C^*$-category of Hilbert spaces.

\medskip

-- We have as well a construction $C\to A$ as above, simply by dividing the free $*$-algebra on $N^2$ variables by the relations in the statement.

\medskip

Some elementary algebra shows then that $C=C_{A_C}$ implies $A=A_{C_A}$, and also that $C\subset C_{A_C}$ is automatic. Thus we are left with proving $C_{A_C}\subset C$, and this can be done by doing some algebra, and using von Neumann's bicommutant theorem. See \cite{ba3}.
\end{proof}

We will need as well, following the classical work of Weyl, Brauer and many others, the notion of ``easiness''. Let us start with the following definition:

\index{category of partitions}

\begin{definition}
Let $P(k,l)$ be the set of partitions between an upper row of $k$ points, and a lower row of $l$ points. A set $D=\bigsqcup_{k,l}D(k,l)$ with $D(k,l)\subset P(k,l)$ is called a category of partitions when it has the following properties:
\begin{enumerate}
\item Stability under the horizontal concatenation, $(\pi,\sigma)\to[\pi\sigma]$.

\item Stability under the vertical concatenation, $(\pi,\sigma)\to[^\sigma_\pi]$.

\item Stability under the upside-down turning, $\pi\to\pi^*$.

\item Each set $P(k,k)$ contains the identity partition $||\ldots||$.

\item The set $P(0,2)$ contains the semicircle partition $\cap$.
\end{enumerate}
\end{definition} 

Observe that this is precisely the definition that we used in chapter 12, with the condition there on the basic crossing $\slash\hskip-2.1mm\backslash$, which produces commutativity via Tannakian duality, removed. In relation with the quantum groups, we have the following notion:

\index{easiness}
\index{easy quantum group}

\begin{definition}
A compact quantum matrix group $G$ is called easy when
$$Hom(u^{\otimes k},u^{\otimes l})=span\left(T_\pi\Big|\pi\in D(k,l)\right)$$
for any colored integers $k,l$, for certain sets of partitions $D(k,l)\subset P(k,l)$, where
$$T_\pi(e_{i_1}\otimes\ldots\otimes e_{i_k})=\sum_{j_1\ldots j_l}\delta_\pi\begin{pmatrix}i_1&\ldots&i_k\\ j_1&\ldots&j_l\end{pmatrix}e_{j_1}\otimes\ldots\otimes e_{j_l}$$
with the Kronecker type symbols $\delta_\pi\in\{0,1\}$ depending on whether the indices fit or not. 
\end{definition}

Again, this is something coming as a continuation of the material from chapter 12. Many things can be said here, but getting now straight to the point, we have:

\index{Brauer theorem}

\begin{theorem}
We have the following results:
\begin{enumerate}
\item $S_N$ is easy, coming from the category of all partitions $P$.

\item $S_N^+$ is easy, coming from the category of all noncrossing partitions $NC$.
\end{enumerate}
\end{theorem}

\begin{proof}
This is something quite fundamental, with the proof, using the above Tannakian results and subsequent easiness theory, being as follows:

\medskip

(1) $S_N^+$. We know that this quantum group comes from the magic condition. In order to interpret this magic condition, consider the fork partition:
$$Y\in P(2,1)$$

By arguing as in chapter 8, we conclude that we have the following equivalence:
$$T_Y\in Hom(u^{\otimes 2},u)\iff u_{ij}u_{ik}=\delta_{jk}u_{ij},\forall i,j,k$$

The condition on the right being equivalent to the magic condition, we conclude that $S_N^+$ is indeed easy, the corresponding category of partitions being, as desired:
$$D
=<Y>
=NC$$

(2) $S_N$. Here there is no need for new computations, because we have:
$$S_N=S_N^+\cap O_N$$

At the categorical level means that $S_N$ is easy, coming from:
$$<NC,\slash\hskip-2.2mm\backslash>=P$$

Thus, we are led to the conclusions in the statement.
\end{proof}

Summarizing, we have now a good understanding of the liberation operation $S_N\to S_N^+$, the idea being that this comes, via Tannakian duality, from $P\to NC$. 

\bigskip

In order to go further in this direction, we will need the following result, with $*$ being the classical convolution, and $\boxplus$ being Voiculescu's free convolution operation \cite{vdn}:

\begin{theorem}
The following Poisson type limits converge, for any $t>0$,
$$p_t=\lim_{n\to\infty}\left(\left(1-\frac{1}{n}\right)\delta_0+\frac{1}{n}\delta_t\right)^{*n}$$
$$\pi_t=\lim_{n\to\infty}\left(\left(1-\frac{1}{n}\right)\delta_0+\frac{1}{n}\delta_t\right)^{\boxplus n}$$
the limiting measures being the Poisson law $p_t$, and the Marchenko-Pastur law $\pi_t$, 
$$p_t=\frac{1}{e^t}\sum_{k=0}^\infty\frac{t^k\delta_k}{k!}$$
$$\pi_t=\max(1-t,0)\delta_0+\frac{\sqrt{4t-(x-1-t)^2}}{2\pi x}\,dx$$
whose moments are given by the following formulae:
$$M_k(p_t)=\sum_{\pi\in P(k)}t^{|\pi|}\quad,\quad 
M_k(\pi_t)=\sum_{\pi\in NC(k)}t^{|\pi|}$$
The Marchenko-Pastur measure $\pi_t$ is also called free Poisson law.
\end{theorem}

\begin{proof}
This is something quite advanced, related to probability theory, free probability theory, and random matrices, the idea being as follows:

\medskip

(1) The first step is that of finding suitable functional transforms, which linearize the convolution operations in the statement. In the classical case this is the logarithm of the Fourier transform $\log F$, and in the free case this is Voiculescu's $R$-transform.

\medskip

(2) With these tools in hand, the above limiting theorems can be proved in a standard way, a bit as when proving the Central Limit Theorem. The computations give the moment formulae in the statement, and the density computations are standard as well.

\medskip

(3) Finally, in order for the discussion to be complete, what still remains to be explained is the precise nature of the ``liberation'' operation $p_t\to\pi_t$, as well as the random matrix occurrence of $\pi_t$. This is more technical, and we refer here to \cite{bpa}, \cite{mpa}, \cite{vdn}.
\end{proof}

Getting back now to quantum permutations, the results here are as follows:

\index{main character}
\index{truncated character}
\index{Poisson law}
\index{free Poisson law}
\index{Marchenko-Pastur law}
\index{Weingarten formula}

\begin{theorem}
The law of the main character, given by 
$$\chi=\sum_iu_{ii}$$
for $S_N/S_N^+$ becomes $p_1/\pi_1$ with $N\to\infty$. As for the truncated character 
$$\chi_t=\sum_{i=1}^{[tN]}u_{ii}$$
for $S_N/S_N^+$, with $t\in(0,1]$, this becomes $p_t/\pi_t$ with $N\to\infty$.
\end{theorem}

\begin{proof}
This is again something quite technical, the idea being as follows:

\medskip

(1) In the classical case this is well-known, and follows by using the inclusion-exclusion principle, and then letting $N\to\infty$, as explained in chapter 11. 

\medskip

(2) In the free case there is no such simple argument, and we must use what we know about $S_N^+$, namely its easiness property. We know from easiness that we have:
$$Fix(u^{\otimes k})=span(NC(k))$$

On the other hand, a direct computation shows that the partitions in $P(k)$, and in particular those in $NC(k)$, implemented as linear maps via the operation $\pi\to T_\pi$ from Definition 13.34, become linearly independent with $N\geq k$. Thus we have, as desired:
\begin{eqnarray*}
\int_{S_N^+}\chi^k
&=&\dim\left(Fix(u^{\otimes k})\right)\\
&=&\dim\left(span\left(T_\pi\Big|\pi\in NC(k)\right)\right)\\
&\simeq&|NC(k)|\\
&=&\sum_{\pi\in NC(k)}1^{|\pi|}
\end{eqnarray*}

(3) In the general case now, where our parameter is an arbitrary number $t\in(0,1]$, the above computation does not apply, but we can still get away with Peter-Weyl theory. Indeed, we know from Theorem 13.25 how to compute the Haar integration of $S_N^+$, out of the knowledge of the fixed point spaces $Fix(u^{\otimes k})$, and in practice, by using easiness, this leads to the following formula, called Weingarten integration formula:
$$\int_{S_N^+}u_{i_1j_1}\ldots u_{i_kj_k}=\sum_{\pi,\sigma\in NC(k)}\delta_\pi(i)\delta_\sigma(j)W_{kN}(\pi,\sigma)$$

Here the $\delta$ symbols are Kronecker type symbols, checking whether the indices fit or not with the partitions, and $W_{kN}=G_{kN}^{-1}$, with $G_{kN}(\pi,\sigma)=N^{|\pi\vee\sigma|}$, where $|.|$ is the number of blocks. Now by using this formula for computing the moments of $\chi_t$, we obtain:
\begin{eqnarray*}
\int_{S_N^+}\chi_t^k
&=&\sum_{i_1=1}^{[tN]}\ldots\sum_{i_k=1}^{[tN]}\int u_{i_1i_1}\ldots u_{i_ki_k}\\
&=&\sum_{\pi,\sigma\in NC(k)}W_{kN}(\pi,\sigma)\sum_{i_1=1}^{[tN]}\ldots\sum_{i_k=1}^{[tN]}\delta_\pi(i)\delta_\sigma(i)\\
&=&\sum_{\pi,\sigma\in NC(k)}W_{kN}(\pi,\sigma)G_{k[tN]}(\sigma,\pi)\\
&=&Tr(W_{kN}G_{k[tN]})
\end{eqnarray*}

(4) The point now is that with $N\to\infty$ the Gram matrix $G_{kN}$, and so the Weingarten matrix $W_{kN}$ too, becomes asymptotically diagonal. We therefore obtain:
$$\int_{S_N^+}\chi_t^k\simeq\sum_{\pi\in NC(k)}t^{|\pi|}$$

Thus, we are led to the conclusion in the statement. For details, see \cite{ba3}.
\end{proof}

\section*{13e. Exercises}

This was a pleasant chapter to write for me, because I've been doing such things for long, but probably quite hard to read, for you. As exercises, we have:

\begin{exercise}
Do Gram-Schmidt for spaces $L^2(X)$ of your choice.
\end{exercise}

\begin{exercise}
Clarify the details in the proof in the spectral radius formula.
\end{exercise}

\begin{exercise}
Learn more about the GNS embedding theorem, and its proof.
\end{exercise}

\begin{exercise}
Learn the details of Woronowicz's Peter-Weyl theory.
\end{exercise}

\begin{exercise}
Can we compute the Haar measure via Peter-Weyl, and how.
\end{exercise}

\begin{exercise}
Find a new, very simple proof for $S_3^+=S_3$.
\end{exercise}

\begin{exercise}
Prove that the quantum group $S_4^+$ is coamenable.
\end{exercise}

\begin{exercise}
Compute the representations of $S_N^+$, at $N\geq4$.
\end{exercise}

As bonus exercise, learn some quantum mechanics, from Feynman \cite{fe3}, or Griffiths \cite{gr2}, or Weinberg \cite{we2}. There ain't no quantum without quantum mechanics.

\chapter{Graph symmetries}

\section*{14a. Graph symmetries}

We can get back now to graphs. By using the quantum permutation group $S_N^+$ constructed in the previous chapter, we can perform the following construction:

\begin{theorem}
Given a finite graph $X$, with adjacency matrix $d\in M_N(0,1)$, the following construction produces a quantum permutation group, 
$$C(G^+(X))=C(S_N^+)\Big/\Big<du=ud\Big>$$
whose classical version $G(X)$ is the usual  automorphism group of $X$.
\end{theorem}

\begin{proof}
The fact that we have indeed a quantum group comes from the fact that $du=ud$ reformulates as $d\in End(u)$, which makes it clear that $\Delta,\varepsilon,S$ factorize. Regarding now the second assertion, we must establish here the following equality:
$$C(G(X))=C(S_N)\Big/\Big<du=ud\Big>$$

For this purpose, recall that we have $u_{ij}(\sigma)=\delta_{\sigma(j)i}$. We therefore obtain:
$$(du)_{ij}(\sigma)
=\sum_kd_{ik}u_{kj}(\sigma)
=\sum_kd_{ik}\delta_{\sigma(j)k}
=d_{i\sigma(j)}$$

On the other hand, we have as well the following formula:
$$(ud)_{ij}(\sigma)
=\sum_ku_{ik}(\sigma)d_{kj}
=\sum_k\delta_{\sigma(k)i}d_{kj}
=d_{\sigma^{-1}(i)j}$$

Thus $du=ud$ reformulates as $d_{ij}=d_{\sigma(i)\sigma(j)}$, which gives the result.
\end{proof}

Let us work out some examples. With the convention that $\hat{*}$ is the dual free product, obtained by diagonally concatenating the magic unitaries, we have:

\index{simplex}
\index{square graph}
\index{disconnected union}

\begin{proposition}
The construction $X\to G^+(X)$ has the following properties:
\begin{enumerate}
\item For the $N$-point graph, having no edges at all, we obtain $S_N^+$.

\item For the $N$-simplex, having edges everywhere, we obtain as well $S_N^+$.

\item We have $G^+(X)=G^+(X^c)$, where $X^c$ is the complementary graph.

\item For a disconnected union, we have $G^+(X)\,\hat{*}\,G^+(Y)\subset G^+(X\sqcup Y)$.

\item For the square, we obtain a non-classical, proper subgroup of $S_4^+$.
\end{enumerate}
\end{proposition}

\begin{proof}
All these results are elementary, the proofs being as follows:

\medskip

(1) This follows from definitions, because here we have $d=0$.

\medskip

(2) Here $d=\mathbb I-1$, where $\mathbb I$ is the all-one matrix, and the magic condition gives $u\mathbb I=\mathbb Iu=N\mathbb I$. We conclude that $du=ud$ is automatic, and so $G^+(X)=S_N^+$.

\medskip

(3) The adjacency matrices of $X,X^c$ being related by the following formula: $$d_X+d_{X^c}=\mathbb I-1$$

By using now the above formula $u\mathbb I=\mathbb Iu=N\mathbb I$, we conclude that $d_Xu=ud_X$ is equivalent to $d_{X^c}u=ud_{X^c}$. Thus, we obtain, as claimed, $G^+(X)=G^+(X^c)$.

\medskip

(4) The adjacency matrix of a disconnected union is given by:
$$d_{X\sqcup Y}=diag(d_X,d_Y)$$

Now let $w=diag(u,v)$ be the fundamental corepresentation of $G^+(X)\,\hat{*}\,G^+(Y)$. Then $d_Xu=ud_X$ and $d_Yv=vd_Y$, and we obtain, as desired, $d_{X\sqcup Y}w=wd_{X\sqcup Y}$.

\medskip

(5) We know from (3) that we have $G^+(\square)=G^+(|\ |)$. We know as well from (4) that we have $\mathbb Z_2\,\hat{*}\,\mathbb Z_2\subset G^+(|\ |)$. It follows that $G^+(\square)$ is non-classical. Finally, the inclusion $G^+(\square)\subset S_4^+$ is indeed proper, because $S_4\subset S_4^+$ does not act on the square.
\end{proof}

In order to further advance, and to explicitely compute various quantum automorphism groups, we can use the spectral decomposition of $d$, as follows:

\index{spectral decomposition}

\begin{theorem}
A closed subgroup $G\subset S_N^+$ acts on a graph $X$ precisely when
$$P_\lambda u=uP_\lambda\quad,\quad\forall\lambda\in\mathbb R$$
where $d=\sum_\lambda\lambda\cdot P_\lambda$ is the spectral decomposition of the adjacency matrix of $X$.
\end{theorem}

\begin{proof}
Since $d\in M_N(0,1)$ is a symmetric matrix, we can consider indeed its spectral decomposition, $d=\sum_\lambda\lambda\cdot P_\lambda$, and we have the following formula:
$$<d>=span\left\{P_\lambda\Big|\lambda\in\mathbb R\right\}$$

Thus $d\in End(u)$ when $P_\lambda\in End(u)$ for all $\lambda\in\mathbb R$, which gives the result.
\end{proof}

In order to exploit Theorem 14.3, we will often combine it with the following fact:

\begin{proposition}
Given a closed subgroup $G\subset S_N^+$, with associated coaction
$$\Phi:\mathbb C^N\to \mathbb C^N\otimes C(G)\quad,\quad e_i\to\sum_je_j\otimes u_{ji}$$
and a linear subspace $V\subset\mathbb C^N$, the following are equivalent:
\begin{enumerate}
\item The magic matrix $u=(u_{ij})$ commutes with $P_V$.

\item $V$ is invariant, in the sense that $\Phi(V)\subset V\otimes C(G)$.
\end{enumerate}
\end{proposition}

\begin{proof}
Let $P=P_V$. For any $i\in\{1,\ldots,N\}$ we have the following formula:
$$\Phi(P(e_i))
=\Phi\left(\sum_kP_{ki}e_k\right) 
=\sum_{jk}P_{ki}e_j\otimes u_{jk}
=\sum_je_j\otimes (uP)_{ji}$$

On the other hand the linear map $(P\otimes id)\Phi$ is given by a similar formula:
$$(P\otimes id)(\Phi(e_i))
=\sum_kP(e_k)\otimes u_{ki}
=\sum_{jk}P_{jk}e_j\otimes u_{ki}
=\sum_je_j\otimes (Pu)_{ji}$$

Thus $uP=Pu$ is equivalent to $\Phi P=(P\otimes id)\Phi$, and the conclusion follows.
\end{proof}

As an application of the above results, we have the following computation:

\index{square graph}
\index{cycle graph}
\index{circulant graph}

\begin{theorem}
The quantum automorphism group of the $N$-cycle is, at $N\neq4$:
$$G^+(X)=D_N$$
However, at $N=4$ we have $D_4\subset G^+(X)\subset S_4^+$, with proper inclusions.
\end{theorem}

\begin{proof}
We know from Proposition 14.2, and from $S_N=S_N^+$ at $N\leq3$, that the various assertions hold indeed at $N\leq4$. So, assume $N\geq5$. Given a $N$-th root of unity, $w^N=1$, the vector $\xi=(w^i)$ is an eigenvector of $d$, with eigenvalue as follows:
$$\lambda=w+w^{N-1}$$

Now by taking $w=e^{2\pi i/N}$, it follows that the are eigenvectors of $d$ are:
$$1,f,f^2,\ldots ,f^{N-1}$$

More precisely, the invariant subspaces of $d$ are as follows, with the last subspace having dimension 1 or 2 depending on the parity of $N$:
$$\mathbb C 1,\, \mathbb C f\oplus\mathbb C f^{N-1},\, \mathbb C f^2\oplus\mathbb C f^{N-2},\ldots$$

Assuming $G\subset G^+(X)$, consider the coaction $\Phi:\mathbb C^N\to \mathbb C^N\otimes C(G)$, and write:
$$\Phi(f)=f\otimes a+f^{N-1}\otimes b$$

By taking the square of this equality we obtain the following formula:
$$\Phi(f^2)=f^2\otimes a^2+f^{N-2}\otimes b^2+1\otimes(ab+ba)$$

It follows that $ab=-ba$, and that $\Phi(f^2)$ is given by the following formula:
$$\Phi(f^2)=f^2\otimes a^2+f^{N-2}\otimes b^2$$

By multiplying this with $\Phi(f)$ we obtain the following formula:
$$\Phi(f^3)=f^3\otimes a^3+f^{N-3}\otimes b^3+f^{N-1}\otimes ab^2+f\otimes ba^2$$

Now since $N\geq 5$ implies that $1,N-1$ are different from $3,N-3$, we must have $ab^2=ba^2=0$. By using this and $ab=-ba$, we obtain by recurrence on $k$ that:
$$\Phi(f^k)=f^k\otimes a^k+f^{N-k}\otimes b^k$$

In particular at $k=N-1$ we obtain the following formula:
$$\Phi(f^{N-1})=f^{N-1}\otimes a^{N-1}+f\otimes b^{N-1}$$

On the other hand we have $f^*=f^{N-1}$, so by applying $*$ to $\Phi(f)$ we get:
$$\Phi(f^{N-1})=f^{N-1}\otimes a^*+f\otimes b^*$$

Thus $a^*=a^{N-1}$ and $b^*=b^{N-1}$. Together with $ab^2=0$ this gives:
$$(ab)(ab)^*
=abb^*a^*
=ab^Na^{N-1}
=(ab^2)b^{N-2}a^{N-1}
=0$$

From positivity we get from this $ab=0$, and together with $ab=-ba$, this shows that $a,b$ commute. On the other hand $C(G)$ is generated by the coefficients of $\Phi$, which are powers of $a,b$, and so $C(G)$ must be commutative, and we obtain the result.
\end{proof}

The above result is quite suprising, but we will be back to this, with a more conceptual explanation for the fact that the square $\square$ has quantum symmetry. Back to theory now, we have the following useful result, complementary to Theorem 14.3:

\index{color decomposition}
\index{color components}

\begin{theorem}
Given a matrix $p\in M_N(\mathbb C)$, consider its ``color'' decomposition
$$p=\sum_{c\in\mathbb C}c\cdot p_c$$
with the color components $p_c\in M_N(0,1)$ with $c\in\mathbb C$ being constructed as follows:
$$(p_c)_{ij}=\begin{cases}
1&{\rm if}\ p_{ij}=c\\
0&{\rm otherwise}
\end{cases}$$
Then a magic matrix $u=(u_{ij})$ commutes with $p$ iff it commutes with all matrices $p_c$.
\end{theorem}

\begin{proof}
Consider the multiplication and counit maps of the algebra $\mathbb C^N$:
$$M:e_i\otimes e_j\to e_ie_j\quad,\quad 
C:e_i\to e_i\otimes e_i$$

Since $M,C$ intertwine $u,u^{\otimes 2}$, their iterations $M^{(k)},C^{(k)}$ intertwine $u,u^{\otimes k}$, and so:
$$M^{(k)}p^{\otimes k}C^{(k)}
=\sum_{c\in\mathbb C}c^kp_c
\in End(u)$$

Now since this formula holds for any $k\in\mathbb N$, we obtain the result.
\end{proof}

The above results can be combined, and we are led to the following statement:

\index{spectral-color decomposition}
\index{color-spectral decomposition}
\index{planar algebra}

\begin{theorem}
A closed subgroup $G\subset S_N^+$ acts on a graph $X$ precisely when 
$$u=(u_{ij})$$
commutes with all the matrices coming from the color-spectral decomposition of $d$.
\end{theorem}

\begin{proof}
This follows by combining Theorem 14.3 and Theorem 14.6, with the ``color-spectral'' decomposition in the statement referring to what comes out by succesively doing the color and spectral decomposition, until the process stabilizes.
\end{proof}

This latter statement is quite interesting, with the color-spectral decomposition there being something quite intriguing. We will be back to this later, when discussing planar algebras, which is the good framework for discussing such things.

\section*{14b. Product operations}

We would like to understand how the operation $X\to G^+(X)$ behaves under taking various products of graphs, in analogy with what we know about $X\to G(X)$, from chapter 10. As a first observation, things can be quite tricky here, as shown by:

\begin{fact}
Although the graph formed by two points $\bullet\,\bullet$ has no quantum symmetry, the graph formed by two copies of it, namely $\bullet\,\bullet\,\bullet\,\,\bullet$\,, does have quantum symmetry.
\end{fact}

Which looks quite scary, but no worries, we will manage to reach to a better understanding of this. Getting to work now, let us recall from chapter 10 that we have:

\index{direct product}
\index{Cartesian product}
\index{lexicographic product}

\begin{definition}
Let $X,Y$ be two finite graphs.
\begin{enumerate}
\item The direct product $X\times Y$ has vertex set $X\times Y$, and edges:
$$(i,\alpha)-(j,\beta)\Longleftrightarrow i-j,\, \alpha-\beta$$

\item The Cartesian product $X\,\square\,Y$ has vertex set $X\times Y$, and edges:
$$(i,\alpha)-(j,\beta)\Longleftrightarrow i=j,\, \alpha-\beta\mbox{ \rm{or} }i-j,\alpha=\beta$$

\item The lexicographic product $X\circ Y$ has vertex set $X\times Y$, and edges:
$$(i,\alpha)-(j,\beta)\Longleftrightarrow \alpha-\beta\mbox{ \rm{or} }\alpha=\beta,\, i-j$$
\end{enumerate}
\end{definition}

The above products are all well-known in graph theory, and we have already studied them in chapter 10, in relation with symmetry groups, with the following conclusion:

\begin{theorem}
We have standard embeddings, as follows,
$$G(X)\times G(Y)\subset G(X \times Y)$$
$$G(X)\times G(Y)\subset G(X\,\square\,Y)$$
$$G(X)\wr G(Y)\subset G(X\circ Y)$$
and under suitable spectral assumptions, these embeddings are isomorphisms.
\end{theorem}

\begin{proof}
This is something that we know from chapter 10, and we refer to the discussion there for both the precise statement of the last assertion, and for the proofs.
\end{proof}

In order to discuss the quantum analogues of these embeddings, we need to introduce first a number of operations on the compact quantum groups, similar to the above operations for the finite groups. Following Wang \cite{wan}, we first have:

\index{product of quantum groups}

\begin{proposition}
Given two compact quantum groups $G,H$, so is their product $G\times H$, constructed according to the following formula:
$$C(G\times H)=C(G)\otimes C(H)$$ 
Equivalently, at the level of the associated discrete duals $\Gamma,\Lambda$, we can set
$$C^*(\Gamma\times\Lambda)=C^*(\Gamma)\otimes C^*(\Lambda)$$
and we obtain the same equality of Woronowicz algebras as above.
\end{proposition}

\begin{proof}
Assume indeed that we have two Woronowicz algebras, $(A,u)$ and $(B,v)$. Our claim is that the following construction produces a Woronowicz algebra:
$$C=A\otimes B\quad,\quad w=diag(u,v)$$

Indeed, the matrix $w$ is unitary, and its coefficients generate $C$. As for the existence of the maps $\Delta,\varepsilon,S$, this follows from the functoriality properties of $\otimes$. But with this claim in hand, the first assertion is clear. As for the second assertion, let us recall that when $G,H$ are classical and abelian, we have the following formula:
$$\widehat{G\times H}=\widehat{G}\times\widehat{H}$$

Thus, our second assertion is simply a reformulation of the first assertion, with the $\times$ symbol used there being justified by this well-known group theory formula. 
\end{proof}

Another standard operation, again from \cite{wan}, is that of taking subgroups:

\index{quantum subgroup}
\index{quotient quantum group}
\index{Hopf ideal}

\begin{proposition}
Let $G$ be compact quantum group, and let $I\subset C(G)$ be a closed $*$-ideal satisfying the following condition: 
$$\Delta(I)\subset C(G)\otimes I+I\otimes C(G)$$
We have then a closed quantum subgroup $H\subset G$, constructed as follows:
$$C(H)=C(G)/I$$
At the dual level we obtain a quotient of discrete quantum groups, $\widehat{\Gamma}\to\widehat{\Lambda}$.
\end{proposition}

\begin{proof}
This follows indeed from the above conditions on $I$, which are designed precisely as for $\Delta,\varepsilon,S$ to factorize through the quotient. As for the last assertion, this is just a reformulation, coming from the functoriality properties of the Pontrjagin duality.
\end{proof}

Regarding now taking quotients, the result here, again from \cite{wan}, is as follows:

\index{quantum subgroup}
\index{quotient quantum group}

\begin{proposition}
Let $G$ be a compact quantum group, and $v=(v_{ij})$ be a corepresentation of $C(G)$. We have then a quotient quantum group $G\to H$, given by:
$$C(H)=<v_{ij}>$$
At the dual level we obtain a discrete quantum subgroup, $\widehat{\Lambda}\subset\widehat{\Gamma}$.
\end{proposition}

\begin{proof}
Here the first assertion follows from definitions, and the second assertion is a reformulation of it, using the basic properties of Pontrjagin duality.
\end{proof}

Finally, we will need the notion of free wreath product, from \cite{bi1}, as follows:

\index{free wreath product}

\begin{proposition}
Given closed subgroups $G\subset S_N^+$, $H\subset S_k^+$, with magic corepresentations $u,v$, the following construction produces a closed subgroup of $S_{Nk}^+$:
$$C(G\wr_*H)=(C(G)^{*k}*C(H))/<[u_{ij}^{(a)},v_{ab}]=0>$$
When $G,H$ are classical, the classical version of $G\wr_*H$ is the wreath product $G\wr H$.
\end{proposition}

\begin{proof}
Consider indeed the matrix $w_{ia,jb}=u_{ij}^{(a)}v_{ab}$, over the quotient algebra in the statement. Our claim is that this matrix is magic. Indeed, the entries are projections, because they appear as products of commuting projections, and the row sums are:
$$\sum_{jb}w_{ia,jb}
=\sum_{jb}u_{ij}^{(a)}v_{ab}
=\sum_bv_{ab}\sum_ju_{ij}^{(a)}
=1$$

As for the column sums, these are as follows:
$$\sum_{ia}w_{ia,jb}
=\sum_{ia}u_{ij}^{(a)}v_{ab}
=\sum_av_{ab}\sum_iu_{ij}^{(a)}
=1$$

With this in hand, it is routine to check that $G\wr_*H$ is indeed a quantum permutation group. Finally, the last assertion is standard as well. See \cite{bi1}.
\end{proof}

With the above discussed, we can now go back to graphs. Following \cite{bbi}, the standard embeddings from Theorem 14.10 have the following quantum analogues:

\begin{theorem}
We have embeddings as follows,
$$G^+(X)\times G^+(Y)\subset G^+(X \times Y)$$
$$G^+(X)\times G^+(Y)\subset G^+(X\,\square\,Y)$$
$$G^+(X)\wr_*G^+(Y)\subset G^+(X\circ Y)$$
with the operation $\wr_*$ being a free wreath product. 
\end{theorem}

\begin{proof}
We use the following identification, given by $\delta_{(i,\alpha)}=\delta_i\otimes\delta_\alpha$:
$$C(X\times Y)=C(X)\otimes C(Y)$$

(1) The adjacency matrix of the direct product is given by:
$$d_{X\times Y}=d_X\otimes d_Y$$

Thus if $u$ commutes with $d_X$ and $v$ commutes with $d_Y$, then $u\otimes v=(u_{ij}v_{\alpha\beta})_{(i\alpha,j\beta)}$ is a magic unitary that commutes with $d_{X\times Y}$. But this gives a morphism as follows:
$$C(G^+(X\times Y))\to C(G^+(X)\times G^+(Y))$$ 

Finally, the surjectivity of this morphism follows by summing over $i$ and $\beta$.

\medskip

(2) The adjacency matrix of the Cartesian product is given by:
$$d_{X\,\square\,Y}=d_X\otimes1+1\otimes d_Y$$

Thus if $u$ commutes with $d_X$ and $v$ commutes with $d_Y$, then $u\otimes v=(u_{ij}v_{\alpha\beta})_{(i\alpha,j\beta)}$ is a magic unitary that commutes with $d_{X\,\square\,Y}$, and this gives the result.

\medskip

(3) The adjacency matrix of the lexicographic product $X\circ Y$ is given by:
$$d_{X\circ Y}=d_X\otimes1+\mathbb I\otimes d_Y$$

Now let $u,v$ be the magic unitary matrices of $G^+(X),G^+(Y)$. The magic unitary matrix of $G^+(X)\wr_*G^+(Y)$ is then given by the following formula:
$$w_{ia,jb}= u_{ij}^{(a)}v_{ab}$$

Since $u,v$ commute with $d_X,d_Y$, we get that $w$ commutes with $d_{X\circ Y}$. But this gives the desired morphism, and the surjectivity follows by summing over $i$ and $b$. 
\end{proof}

The problem now is that of deciding when the embeddings in Theorem 14.15 are isomorphisms. Following \cite{bbi}, we first have the following result: 

\index{connected graph}
\index{regular graph}

\begin{theorem}
Let $X$ and $Y$ be finite connected regular graphs. If their spectra $\{\lambda\}$ and $\{\mu\}$ do not contain $0$ and satisfy
$$\big\{\lambda_i/\lambda_j\big\}\cap\big\{\mu_k/\mu_l\big\}=\{1\}$$
then $G^+(X\times Y)=G^+(X)\times G^+(Y)$. Also, if their spectra satisfy
$$\big\{\lambda_i-\lambda_j\big\}\cap\big\{\mu_k-\mu_l\big\}=\{0\}$$
then $G^+(X\,\square\,Y)=G^+(X)\times G^+(Y)$.
\end{theorem}

\begin{proof}
This is something quite standard, the idea being as follows:

\medskip

(1) First, we know from Theorem 14.15 that we have embeddings as follows, valid for any two graphs $X,Y$, and coming from definitions:
$$G^+(X)\times G^+(Y)\subset G^+(X\times Y)$$
$$G^+(X)\times G^+(Y)\subset G^+(X\,\square\,Y)$$

(2) Now let $\lambda_1$ be the valence of $X$. Since $X$ is regular we have $\lambda_1\in Sp(X)$, with $1$ as eigenvector, and since $X$ is connected $\lambda_1$ has multiplicity 1. Hence if $P_1$ is the orthogonal projection onto $\mathbb C1$, the spectral decomposition of $d_X$ is of the following form:
$$d_X=\lambda_1P_1+\sum_{i\neq1}\lambda_iP_i$$

We have a similar formula for the adjacency matrix $d_Y$, namely:
$$d_Y=\mu_1Q_1+\sum_{j\neq1}\mu_jQ_j$$

(3) But this gives the following formulae for products:
$$d_{X\times Y}=\sum_{ij}(\lambda_i\mu_j)P_{i}\otimes Q_{j}$$
$$d_{X\,\square\,Y}=\sum_{ij}(\lambda_i+\mu_i)P_i\otimes Q_j$$

Here the projections form partitions of unity, and the scalar are distinct, so these are spectral decompositions. The coactions will commute with any of the spectral projections, and hence with both $P_1\otimes1$, $1\otimes Q_1$. In both cases the universal coaction $v$ is the tensor product of its restrictions to the images of $P_1\otimes1$, $1\otimes Q_1$, which gives the result.
\end{proof}

Regarding now the lexicographic products, the result here is as follows:

\begin{theorem}
Let $X,Y$ be regular graphs, with $X$ connected. If their spectra $\{\lambda_i\}$ and $\{\mu_j\}$ satisfy the condition
$$\big\{\lambda_1-\lambda_i\big|i\neq 1\big\}\cap\big\{-n\mu_j\big\}=\emptyset$$
where $n$ and $\lambda_1$ are the order and valence of $X$, then $G^+(X\circ Y)=G^+(X)\wr_*G^+(Y)$.   
\end{theorem}

\begin{proof}
This is something quite tricky, the idea being as follows:

\medskip

(1) First, we know from Theorem 14.15 that we have an embedding as follows, valid for any two graphs $X,Y$, and coming from definitions:
$$G^+(X)\wr_*G^+(Y)\subset G^+(X\circ Y)$$

(2) We denote by $P_i,Q_j$ the spectral projections corresponding to $\lambda_i,\mu_j$. Since $X$ is connected we have $P_1=\mathbb I/n$, and we obtain:
\begin{eqnarray*}
d_{X\circ Y}
&=&d_X\otimes 1+{\mathbb I}\otimes d_Y\\
&=&\left(\sum_i\lambda_iP_i\right)\otimes\left(\sum_jQ_j\right)+\left(nP_1\right)\otimes \left(\sum_i\mu_jQ_j\right)\\
&=&\sum_j(\lambda_1+n\mu_j)(P_1\otimes Q_j) + \sum_{i\not=1}\lambda_i(P_i\otimes 1)
\end{eqnarray*} 

Since in this formula the projections form a partition of unity, and the scalars are distinct, we conclude that this is the spectral decomposition of $d_{X\circ Y}$. 

\medskip

(3) Now let $W$ be the universal magic matrix for $X\circ Y$. Then $W$ must commute with all spectral projections, and in particular, we have:
$$[W,P_1\otimes Q_j]=0$$

Summing over $j$ gives $[W,P_1\otimes 1]=0$, so $1\otimes C(Y)$ is invariant under the coaction. So, consider the restriction of $W$, which gives a coaction of $G^+(X\circ Y)$ on $1\otimes C(Y)$, that we can denote as follows, with $y$ being a certain magic unitary:
$$W(1\otimes e_a)=\sum_b1\otimes e_b\otimes y_{ba}$$

(4) On the other hand, according to our definition of $W$, we can write:
$$W(e_i\otimes 1)=\sum_{jb}e_j\otimes e_b\otimes x_{ji}^b$$  

By multiplying by the previous relation, found in (3), we obtain:
$$W(e_i\otimes e_a)
=\sum_{jb}e_j\otimes e_b\otimes y_{ba}x_{ji}^b
=\sum_{jb}e_j \otimes e_b\otimes x_{ji}^by_{ba}$$

But this shows that the coefficients of $W$ are of the following form:
$$W_{jb,ia}=y_{ba}x_{ji}^b=x_{ji}^b y_{ba}$$

(5) Consider now the matrix $x^b=(x_{ij}^b)$. Since $W$ is a morphism of algebras, each row of $x^b$ is a partition of unity. Also, by using the antipode, we have:
$$S\left(\sum_jx_{ji}^{b}\right)
=S\left(\sum_{ja}W_{jb,ia}\right)
=\sum_{ja}W_{ia,jb}
=1$$

As a conclusion to this, the matrix $x^b$ constructed above is magic. 

\medskip

(6) We check now that both $x^a,y$ commute with $d_X,d_Y$. We have:
$$(d_{X\circ Y})_{ia,jb} = (d_X)_{ij}\delta_{ab} + (d_Y)_{ab}$$

Thus the two products between $W$ and $d_{X\circ Y}$ are given by:
$$(Wd_{X\circ Y})_{ia,kc}=\sum_j W_{ia,jc} (d_X)_{jk} + \sum_{jb}W_{ia,jb}(d_Y)_{bc}$$
$$(d_{X\circ Y}W)_{ia,kc}=\sum_j (d_X)_{ij} W_{ja,kc} + \sum_{jb}(d_Y)_{ab}W_{jb,kc}$$

(7) Now since the magic matrix $W$ commutes by definition with $d_{X\circ Y}$, the terms on the right in the above equations are equal, and by summing over $c$ we get:
$$\sum_j x_{ij}^a(d_X)_{jk} + \sum_{cb} y_{ab}(d_Y)_{bc}
= \sum_{j} (d_X)_{ij}x_{jk}^a + \sum_{cb} (d_Y)_{ab}y_{bc}$$

The second sums in both terms are equal to the valence of $Y$, so we get:
$$[x^a,d_X]=0$$

Now once again from the formula coming from $[W,d_{X\circ Y}]=0$, we get:
$$[y,d_Y] =0$$

(8) Summing up, the coefficients of $W$ are of the following form, where $x^b$ are magic unitaries commuting with $d_X$, and $y$ is a magic unitary commuting with $d_Y$: 
$$W_{jb,ia}=x_{ji}^by_{ba}$$

But this gives a morphism $C(G^+(X)\wr_*G^+(Y))\to G^+(X\circ Y)$ mapping $u_{ji}^{(b)}\to x_{ji}^b$ and $v_{ba}\to y_{ba}$, which is inverse to the morphism in (1), as desired.
\end{proof}

As a main application of the above result, we have:

\index{connected graph}

\begin{theorem}
Given a connected graph $X$, and $k\in\mathbb N$, we have the formulae
$$G(kX)=G(X)\wr S_k\quad,\quad 
G^+(kX)=G^+(X)\wr_*S_k^+$$
where $kX=X\sqcup\ldots\sqcup X$ is the $k$-fold disjoint union of $X$ with itself.
\end{theorem}

\begin{proof}
The first formula is something well-known, which follows as well from the second formula, by taking the classical version. Regarding now the second formula, it is elementary to check that we have an inclusion as follows, for any finite graph $X$:
$$G^+(X)\wr_*S_k^+\subset G^+(kX)$$

Regarding now the reverse inclusion, which requires $X$ to be connected, this follows by doing some matrix analysis, by using the commutation with $u$. To be more precise, let us denote by $w$ the fundamental corepresentation of $G^+(kX)$, and set:
$$u_{ij}^{(a)}=\sum_bw_{ia,jb}\quad,\quad 
v_{ab}=\sum_iv_{ab}$$

It is then routine to check, by using the fact that $X$ is indeed connected, that we have here magic unitaries, as in the definition of the free wreath products. Thus, we obtain:
$$G^+(kX)\subset G^+(X)\wr_*S_k^+$$

But this gives the result, as a consequence of Theorem 14.17. See \cite{bbi}.
\end{proof}

We refer to \cite{bbi} and related papers for further results, along these lines.

\section*{14c. Reflection groups}

We know that we have results involving free wreath products, which replace the usual wreath products from the classical case. In particular, the quantum symmetry group of the graph formed by $N$ segments is the hyperoctahedral quantum group $H_N^+=\mathbb Z_2\wr_*S_N^+$, appearing as a free analogue of the usual hyperoctahedral group $H_N=\mathbb Z_2\wr S_N$.

\bigskip

The free analogues of the reflection groups $H_N^s$ can be constructed as follows:

\begin{definition}
The algebra $C(H_N^{s+})$ is the universal $C^*$-algebra generated by $N^2$ normal elements $u_{ij}$, subject to the following relations,
\begin{enumerate}
\item $u=(u_{ij})$ is unitary,

\item $u^t=(u_{ji})$ is unitary,

\item $p_{ij}=u_{ij}u_{ij}^*$ is a projection,

\item $u_{ij}^s=p_{ij}$,
\end{enumerate}
with Woronowicz algebra maps $\Delta,\varepsilon,S$ constructed by universality.
\end{definition}

Here we allow the value $s=\infty$, with the convention that the last axiom simply disappears in this case. Observe that at $s<\infty$ the normality condition is actually redundant, because a partial isometry $a$ subject to the relation $aa^*=a^s$ is normal.

\bigskip

Observe also that we have an inclusion of quantum groups $H_N^s\subset H_N^{s+}$
which is a liberation, in the sense that the classical version of $H_N^{s+}$, obtained by dividing by the commutator ideal, is the group $H_N^s$. Indeed, this follows exactly as for $S_N\subset S_N^+$.

\bigskip

In analogy with the results from the real case, we have the following result:

\begin{proposition}
The algebras $C(H_N^{s+})$ with $s=1,2,\infty$, and their presentation relations in terms of the entries of the matrix $u=(u_{ij})$, are as follows:
\begin{enumerate}
\item For $C(H_N^{1+})=C(S_N^+)$, the matrix $u$ is magic: all its entries are projections, summing up to $1$ on each row and column.

\item For $C(H_N^{2+})=C(H_N^+)$ the matrix $u$ is cubic: it is orthogonal, and the products of pairs of distinct entries on the same row or the same column vanish.

\item For $C(H_N^{\infty+})=C(K_N^+)$ the matrix $u$ is unitary, its transpose is unitary, and all its entries are normal partial isometries.
\end{enumerate}
\end{proposition}

\begin{proof}
This is something elementary, the idea being as follows:

\medskip

(1) This follows from definitions, and from some standard operator algebra tricks.

\medskip

(2) This follows again from definitions, and standard operator algebra tricks.

\medskip

(3) This is just a translation of the definition of $C(H_N^{s+})$, at $s=\infty$.
\end{proof}

Let us prove now that $H_N^{s+}$ with $s<\infty$ is a quantum permutation group. For this purpose, we must change the fundamental representation. Let us start with:

\begin{definition}
A $(s,N)$-sudoku matrix is a magic unitary of size $sN$, of the form
$$m=\begin{pmatrix}
a^0&a^1&\ldots&a^{s-1}\\
a^{s-1}&a^0&\ldots&a^{s-2}\\
\vdots&\vdots&&\vdots\\
a^1&a^2&\ldots&a^0
\end{pmatrix}$$
where $a^0,\ldots,a^{s-1}$ are $N\times N$ matrices.
\end{definition}

The basic examples of such matrices come from the group $H_n^s$. Indeed, with $w=e^{2\pi i/s}$, each of the $N^2$ matrix coordinates $u_{ij}:H_N^s\to\mathbb C$ decomposes as follows:
$$u_{ij}=\sum_{r=0}^{s-1}w^ra^r_{ij}$$

Here each $a^r_{ij}$ is a function taking values in $\{0,1\}$, and so a projection in the $C^*$-algebra sense, and it follows from definitions that these projections form a sudoku matrix. Now with this sudoku matrix notion in hand, we have the following result:

\begin{theorem}
The following happen:
\begin{enumerate}
\item The algebra $C(H_N^s)$ is isomorphic to the universal commutative $C^*$-algebra generated by the entries of a $(s,N)$-sudoku matrix.

\item The algebra $C(H_N^{s+})$ is isomorphic to the universal $C^*$-algebra generated by the entries of a $(s,N)$-sudoku matrix.
\end{enumerate}
\end{theorem}

\begin{proof}
The first assertion follows from the second one. In order to prove now the second assertion, consider the universal algebra in the statement:
$$A=C^*\left(a_{ij}^p\ \Big\vert \left(a^{q-p}_{ij}\right)_{pi,qj}=(s,N)-\mbox{sudoku }\right)$$

Consider also the algebra $C(H_N^{s+})$. According to Definition 14.19, this algebra is presented by certain relations $R$, that we will call here level $s$ cubic conditions:
$$C(H_N^{s+})=C^*\left(u_{ij}\ \Big\vert\  u=N\times N\mbox{ level $s$ cubic }\right)$$

We will construct now a pair of inverse morphisms between these algebras.

\medskip

(1) Our first claim is that $U_{ij}=\sum_pw^{-p}a^p_{ij}$ is a level $s$ cubic unitary. Indeed, by using the sudoku condition, the verification of (1-4) in Definition 14.19 is routine.

\medskip

(2) Our second claim is that the elements $A^p_{ij}=\frac{1}{s}\sum_rw^{rp}u^r_{ij}$, with the convention $u_{ij}^0=p_{ij}$, form a level $s$ sudoku unitary. Once again, the proof here is routine.

\medskip

(3) According to the above, we can define a morphism $\Phi:C(H_N^{s+})\to A$ by the formula $\Phi(u_{ij})=U_{ij}$, and a morphism $\Psi:A\to C(H_N^{s+})$ by the formula $\Psi(a^p_{ij})=A^p_{ij}$.

\medskip

(4) It is then routine to check that $\Phi,\Psi$ are inverse morphisms, by a direct computation of their compositions. Thus, we have an isomorphism $C(H_N^{s+})=A$, as claimed.
\end{proof}

In order to further advance, we will need the following simple fact: 

\begin{proposition}
A $sN\times sN$ magic unitary commutes with the matrix
$$\Sigma=
\begin{pmatrix}
0&I_N&0&\ldots&0\\
0&0&I_N&\ldots&0\\
\vdots&\vdots&&\ddots&\\
0&0&0&\ldots&I_N\\
I_N&0&0&\ldots&0
\end{pmatrix}$$
precisely when it is a sudoku matrix in the sense of Definition 14.21.
\end{proposition}

\begin{proof}
This follows from the fact that commutation with $\Sigma$ means that the matrix is circulant. Thus, we obtain the sudoku relations from Definition 14.21.
\end{proof}

Now let $Z_s$ be the oriented cycle with $s$ vertices, and consider the graph $NZ_s$ consisting of $N$ disjoint copies of it. Observe that, with a suitable labeling of the vertices, the adjacency matrix of this graph is the above matrix $\Sigma$. We obtain from this:

\begin{theorem}
We have the following results:
\begin{enumerate}
\item $H_N^s$ is the symmetry group of $NZ_s$.

\item $H_N^{s+}$ is the quantum symmetry group of $NZ_s$.
\end{enumerate}
\end{theorem}

\begin{proof}
This is something elementary, the idea being as follows:

\medskip

(1) This follows indeed from definitions.

\medskip

(2) This follows from Theorem 14.18 and Proposition 14.23, because the algebra $C(H_N^{s+})$ is the quotient of the algebra $C(S_{sN}^+)$ by the relations making the fundamental corepresentation commute with the adjacency matrix of $NZ_s$.
\end{proof}

Next in line, we must talk about wreath products. We have here:

\begin{theorem}
We have isomorphisms as follows,
$$H_N^s=\mathbb Z_s\wr S_N\quad,\quad H_N^{s+}=\mathbb Z_s\wr_*S_N^+$$
with $\wr$ being a wreath product, and $\wr_*$ being a free wreath product.
\end{theorem}

\begin{proof}
This follows from the following formulae, valid for any connected graph $X$, and explained before in this chapter, applied to the graph $Z_s$:
$$G(NX)=G(X)\wr S_N\quad,\quad 
G^+(NX)=G^+(X)\wr_*S_N^+$$

Alternatively, (1) follows from definitions, and (2) can be proved directly, by constructing a pair of inverse morphisms. For details here, we refer to the literature.
\end{proof}

Regarding now the easiness property of $H_N^s,H_N^{s+}$, we have here:

\begin{theorem}
The quantum groups $H_N^s,H_N^{s+}$ are easy, the corresponding categories
$$P^s\subset P\quad,\quad 
NC^s\subset NC$$
consisting of partitions satisfying $\#\circ=\#\bullet(s)$, as a weighted sum, in each block.
\end{theorem}

\begin{proof}
This is something quite routine, the idea being as follows:

\medskip

(1) We already know this for the reflection group $H_N^s$, from chapter 12, and the idea is that the computation there works for $H_N^{s+}$ too, with minimal changes. Indeed, at $s=1$, to start with, this is something that we already know, from chapter 13.

\medskip

(2) At $s=2$ now, we know that $H_N^+\subset O_N^+$ appears via the cubic relations, namely:
$$u_{ij}u_{ik}=u_{ji}u_{ki}=0\quad,\quad\forall j\neq k$$

We conclude, exactly as in chapter 12, that $H_N^+$ is indeed easy, coming from:
$$D=<H>=NC_{even}$$

(3) Regarding now that case $s=\infty$, for the quantum group $K_N^+$, the proof here is similar, leading this time to the category $\mathcal{NC}_{even}$ of noncrossing matching partitions. 

\medskip

(4) Summarizing, we have the result at $s=1,2,\infty$. But the passage to the general case $s\in\mathbb N$ is then routine, by using functoriality, and the result at $s=\infty$.
\end{proof}

All the above is very nice, and at the first glance, it looks like a complete theory of quantum reflection groups. However, there is a skeleton in the closet, coming from:

\begin{fact}
The symmetry groups of the hypercube $\square_N\subset\mathbb R^N$ are
$$G(\square_N)=H_N\quad,\quad G^+(\square_N)\neq H_N^+$$
with the problem coming from the fact that $H_N^+$ does not act on $\square_N$.
\end{fact}

Excited about this? We are here at the heart of quantum algebra, with all sorts of new phenomena, having no classical counterpart, waiting to be explored. In answer, we will prove that we have an equality as follows, with $O_N^{-1}$ being a certain twist of $O_N$:
$$G^+(\square_N)=O_N^{-1}$$

In order to introduce this new quantum group $O_N^{-1}$, we will need:

\begin{proposition}
There is a signature map $\varepsilon:P_{even}\to\{-1,1\}$, given by 
$$\varepsilon(\tau)=(-1)^c$$
where $c$ is the number of switches needed to make $\tau$ noncrossing. In addition:
\begin{enumerate}
\item For $\tau\in S_k$, this is the usual signature.

\item For $\tau\in P_2$ we have $(-1)^c$, where $c$ is the number of crossings.

\item For $\tau\leq\pi\in NC_{even}$, the signature is $1$.
\end{enumerate}
\end{proposition}

\begin{proof}
The fact that $\varepsilon$ is indeed well-defined comes from the fact that the number $c$ in the statement is well-defined modulo 2, which is standard combinatorics. In order to prove now the remaining assertion, observe that any partition $\tau\in P(k,l)$ can be put in ``standard form'', by ordering its blocks according to the appearence of the first leg in each block, counting clockwise from top left, and then by performing the switches as for block 1 to be at left, then for block 2 to be at left, and so on. With this convention:

\medskip

(1) For $\tau\in S_k$ the standard form is $\tau'=id$, and the passage $\tau\to id$ comes by composing with a number of transpositions, which gives the signature. 

\medskip

(2) For a general $\tau\in P_2$, the standard form is of type $\tau'=|\ldots|^{\cup\ldots\cup}_{\cap\ldots\cap}$, and the passage $\tau\to\tau'$ requires $c$ mod 2 switches, where $c$ is the number of crossings. 

\medskip

(3) Assuming that $\tau\in P_{even}$ comes from $\pi\in NC_{even}$ by merging a certain number of blocks, we can prove that the signature is 1 by proceeding by recurrence.
\end{proof}

With the above result in hand, we can now formulate:

\begin{definition}
Associated to a partition $\pi\in P_{even}(k,l)$ is the linear map
$$\bar{T}_\pi(e_{i_1}\otimes\ldots\otimes e_{i_k})=\sum_{j_1\ldots j_l}\bar{\delta}_\pi\begin{pmatrix}i_1&\ldots&i_k\\ j_1&\ldots&j_l\end{pmatrix}e_{j_1}\otimes\ldots\otimes e_{j_l}$$
where the signed Kronecker symbols
$$\bar{\delta}_\pi\in\{-1,0,1\}$$
are given by $\bar{\delta}_\pi=\varepsilon(\tau)$ if $\tau\geq\pi$, and $\bar{\delta}_\pi=0$ otherwise, with $\tau=\ker(^i_j)$.
\end{definition}

In other words, what we are doing here is to add signatures to the usual formula of $T_\pi$. Indeed, observe that the usual formula for $T_\pi$ can be written as folllows:
$$T_\pi(e_{i_1}\otimes\ldots\otimes e_{i_k})=\sum_{j:\ker(^i_j)\geq\pi}e_{j_1}\otimes\ldots\otimes e_{j_l}$$

Now by inserting signs, coming from the signature map $\varepsilon:P_{even}\to\{\pm1\}$, we are led to the following formula, which coincides with the one above:
$$\bar{T}_\pi(e_{i_1}\otimes\ldots\otimes e_{i_k})=\sum_{\tau\geq\pi}\varepsilon(\tau)\sum_{j:\ker(^i_j)=\tau}e_{j_1}\otimes\ldots\otimes e_{j_l}$$

Getting now to quantum groups, we have the following construction:

\begin{theorem}
Given a category of partitions $D\subset P_{even}$, the construction
$$Hom(u^{\otimes k},u^{\otimes l})=span\left(\bar{T}_\pi\Big|\pi\in D(k,l)\right)$$
produces via Tannakian duality a quantum group $G_N^{-1}$, for any $N\in\mathbb N$.
\end{theorem}

\begin{proof}
It is routine to check that the assignement $\pi\to\bar{T}_\pi$ is categorical, in the sense that we have the following formulae, where $c(\pi,\sigma)$ are certain positive integers:
$$\bar{T}_\pi\otimes\bar{T}_\sigma=\bar{T}_{[\pi\sigma]}\quad,\quad 
\bar{T}_\pi \bar{T}_\sigma=N^{c(\pi,\sigma)}\bar{T}_{[^\sigma_\pi]}\quad,\quad  
\bar{T}_\pi^*=\bar{T}_{\pi^*}$$

But with this, the result follows from the Tannakian results from chapter 13.
\end{proof}

We can unify the easy quantum groups, or at least the examples coming from categories $D\subset P_{even}$, with the quantum groups constructed above, as follows:

\index{Schur-Weyl twist}

\begin{definition}
A quantum group $G$ is called $q$-easy, or quizzy, with deformation parameter $q=\pm1$, when its tensor category appears as
$$Hom(u^{\otimes k},u^{\otimes l})=span\left(\dot{T}_\pi\Big|\pi\in D(k,l)\right)$$
for a certain category of partitions $D\subset P_{even}$, where, for $q=-1,1$: 
$$\dot{T}=\bar{T},T$$
The Schur-Weyl twist of $G$ is the quizzy quantum group $G^{-1}$ obtained via $q\to-q$.
\end{definition}

As an illustration for all this, which might seem quite abstract, we can now twist the orthogonal group $O_N$, and the unitary group $U_N$ too. The result here is as follows:

\begin{theorem}
The twists of $O_N,U_N$ are obtained by replacing the commutation relations $ab=ba$ between the coordinates $u_{ij}$ and their adjoints $u_{ij}^*$ with the relations
$$ab=\pm ba$$
with anticommutation on rows and columns, and commutation otherwise.
\end{theorem}

\begin{proof}
The basic crossing, $\ker\binom{ij}{ji}$ with $i\neq j$, comes from the transposition $\tau\in S_2$, so its signature is $-1$. As for its degenerated version $\ker\binom{ii}{ii}$, this is noncrossing, so here the signature is $1$. We conclude that the linear map associated to the basic crossing is:
$$\bar{T}_{\slash\!\!\!\backslash}(e_i\otimes e_j)
=\begin{cases}
-e_j\otimes e_i&{\rm for}\ i\neq j\\
e_j\otimes e_i&{\rm otherwise}
\end{cases}$$

We can proceed now as in the untwisted case, and since the intertwining relations coming from $\bar{T}_{\slash\!\!\!\backslash}$ correspond to the relations defining $O_N^{-1},U_N^{-1}$, we obtain the result.
\end{proof}

Getting back now to graphs, we have the following result, from \cite{bbc}:

\begin{theorem}
The quantum symmetry group of the $N$-hypercube is
$$G^+(\square_N)=O_N^{-1}$$
with the corresponding coaction map on the vertex set being given by
$$\Phi:C^*(\mathbb Z_2^N)\to C^*(\mathbb Z_2^N)\otimes C(O_N^{-1})\quad,\quad 
g_i\to\sum_jg_j\otimes u_{ji}$$
via the standard identification $\square_N=\widehat{\mathbb Z_2^N}$. In particular we have $G^+(\square)=O_2^{-1}$.
\end{theorem} 

\begin{proof}
This follows from a standard algebraic study, done in \cite{bbc}, as follows:

\medskip

(1) Our first claim is that $\square_N$ is the Cayley graph of $\mathbb Z_2^N=<\tau_1,\ldots ,\tau_N>$. Indeed, the vertices of this latter Cayley graph are the products of the following form:
$$g=\tau_1^{i_1}\ldots\tau_N^{i_N}$$

The sequence of exponents defining such an element determines a point of $\mathbb R^N$, which is a vertex of the cube. Thus the vertices of the Cayley graph are the vertices of the cube, and in what regards the edges, this is something that we know too, from chapter 4.

\medskip

(2) Our second claim now, which is something routine, coming from an elementary computation, is that when identifying the vector space spanned by the vertices of $\square_N$ with the algebra $C^*(\mathbb Z_2^N)$, the eigenvectors and eigenvalues of $\square_N$ are given by:
$$v_{i_1\ldots i_N}=\sum_{j_1\ldots j_N} (-1)^{i_1j_1
+\ldots+i_Nj_N}\tau_1^{j_1}\ldots\tau_N^{j_N}$$
$$\lambda_{i_1\ldots i_N}=(-1)^{i_1}+\ldots +(-1)^{i_N}$$

(3) We prove now that the quantum group $O_N^{-1}$ acts on the cube $\square_N$. For this purpose, observe first that we have a map as follows:
$$\Phi :C^*(\mathbb Z_2^N)\to C^*(\mathbb Z_2^N)\otimes C(O_N^{-1})
\quad,\quad\tau_i\to\sum_j\tau_j\otimes u_{ji}$$

It is routine to check that for $i_1\neq i_2\neq\ldots\neq i_l$ we have:
$$\Phi(\tau_{i_1}\ldots\tau_{i_l})=\sum_{j_1\neq\ldots\neq j_l}\tau_{j_1}\ldots\tau_{j_l} \otimes u_{j_1i_1}\ldots u_{j_li_l}$$

In terms of eigenspaces $E_s$ of the adjacency matrix, this gives, as desired:
$$\Phi(E_s)\subset E_s\otimes C(O_N^{-1})$$

(4) Conversely now, consider the universal coaction on the cube:
$$\Psi:C^*(\mathbb Z_2^N)\to C^*(\mathbb Z_2^N)\otimes C(G)\quad,\quad\tau_i\to\sum_j\tau_j\otimes u_{ji}$$

By applying $\Psi$ to the relation $\tau_i\tau_j=\tau_j\tau_i$ we get $u^tu=1$, so the matrix $u=(u_{ij})$ is orthogonal. By applying $\Psi$ to the relation $\tau_i^2=1$ we get:
$$1\otimes\sum_ku_{ki}^2+\sum_{k<l}\tau_k\tau_l\otimes(u_{ki}u_{li}+u_{li}u_{ki})=1\otimes 1$$

This gives $u_{ki}u_{li}=-u_{li}u_{ki}$ for $i\neq j$, $k\neq l$, and by using the antipode we get $u_{ik}u_{il}=-u_{il}u_{ik}$ for $k\neq l$. Also, by applying $\Psi$ to $\tau_i\tau_j=\tau_j\tau_i$ with $i\neq j$ we get:
$$\sum_{k<l}\tau_k\tau_l\otimes(u_{ki}u_{lj}+u_{li}u_{kj})=\sum_{k<l}\tau_k\tau_l\otimes
(u_{kj}u_{li}+u_{lj}u_{ki})$$

Identifying coefficients, it follows that for $i\neq j$ and $k\neq l$, we have:
$$u_{ki}u_{lj}+u_{li}u_{kj}=u_{kj}u_{li}+u_{lj}u_{ki}$$

In other words, we have $[u_{ki},u_{lj}]=[u_{kj},u_{li}]$. By using the antipode we get:
$$[u_{jl},u_{ik}]=[u_{il},u_{jk}]$$

Now by combining these relations we get:
$$[u_{il},u_{jk}]=[u_{ik},u_{jl}]=[u_{jk},u_{il}]=-[u_{il},u_{jk}]$$

Thus $[u_{il},u_{jk}]=0$, so the elements $u_{ij}$ satisfy the relations for $C(O_N^{-1})$, as desired.
\end{proof}

Many other things can be said about twists, reflection groups, and their actions on graphs, and for an introduction to this, we recommend \cite{bbc}, followed by \cite{sc2}.

\section*{14d. Small graphs}

Generally speaking, the graphs having small number of vertices can be investigated by using product operations plus complementation. The first graph which is resistent to such a study is the torus graph $K_3\times K_3$, but we have here, following \cite{bbi}:

\begin{theorem}
The torus graph, obtained as a product of a triangle with itself,
$$X=K_3\times K_3$$
has no quantum symmetry, $G^+(X)=G(X)=S_3\wr\mathbb Z_2$. 
\end{theorem}

\begin{proof}
This is something quite tricky, the idea being as follows:

\medskip

(1) To start with, we have $Sp(X)=\{-2,1,4\}$, the corresponding eigenspaces being as follows, with $\xi_{ij}=\xi^i\otimes \xi^j$, where $\xi=(1,w,w^2)$, with $w=e^{2\pi i/3}$:
$$E_{-2}=\mathbb C\xi_{10}\oplus\mathbb C\xi_{01}\oplus\mathbb C\xi_{20}\oplus\mathbb C\xi_{02}$$
$$E_1=\mathbb C\xi_{11}\oplus\mathbb C\xi_{12}\oplus\mathbb C\xi_{21}\oplus\mathbb C\xi_{22}$$
$$E_4=\mathbb C\xi_{00}$$

(2) Since the universal coaction $v:C(X)\to C(X)\otimes A$ preserves the eigenspaces, we can write formulae as follows, for some $a,b,c,d,\alpha,\beta,\gamma,\delta \in A$:
$$v(\xi_{10})=\xi_{10}\otimes a+\xi_{01}\otimes b+\xi_{20}\otimes c+\xi_{02}\otimes d$$
$$v(\xi_{01})=\xi_{10}\otimes\alpha+\xi_{01}\otimes\beta+\xi_{20}\otimes\gamma+\xi_{02}\otimes\delta$$

Taking the square of $v(\xi_{10})$ gives the following formula:
$$v(\xi_{20})=\xi_{20}\otimes a^2+\xi_{02}\otimes b^2+\xi_{10}\otimes c^2+\xi_{01}\otimes d^2$$

Also, from eigenspace preservation, we have the following relations:
$$ab=-ba\ , \ ad=-da\ , \ bc=-cb\ , \ cd=-dc$$
$$ac+ca=-(bd+db)$$

(3) Now since $a,b$ anticommute, their squares have to commute. On the other hand, by applying $v$ to the equality $\xi_{10}^*=\xi_{20}$, we get the following formulae for adjoints:
$$a^*=a^2\ , \ b^*=b^2\ , \ c^*=c^2\ , \ d^*=d^2$$

The commutation relation $a^2b^2=b^2a^2$ reads now $a^*b^*=b^*a^*$, and by taking adjoints we get $ba=ab$. Together with $ab=-ba$ this gives:
$$ab=ba=0$$

The same method applies to $ad,bc,cd$, and we end up with:
$$ab=ba=0\ ,\ ad=da =0\ , \ bc=cb =0\ , \ cd=dc=0$$

(4) We apply now $v$ to the equality $1=\xi_{10}\xi_{20}$. We get that $1$ is the sum of $16$ terms, all of them of the form $\xi_{ij}\otimes P$, where $P$ are products between $a,b,c,d$ and their squares. Due to the above formulae 8 terms vanish, and the $8$ remaining ones give:
$$1=a^3 +b^3 +c^3 +d^3$$

We have as well the relations coming from eigenspace preservation, namely:
$$ac^2=ca^2=bd^2=db^2=0$$

(5) Now from $ac^2=0$ we get $a^2c^2=0$, and by taking adjoints this gives $ca=0$. The same method applies to $ac,bd,db$, and we end up with:
$$ac=ca=0\ ,\ bd=db=0$$

In the same way we can show that $\alpha,\beta,\gamma,\delta$ pairwise commute:
$$\alpha\beta=\beta\alpha=\ldots =\gamma\delta=\delta\gamma=0$$

(6) In order to finish the proof, it remains to show that $a,b,c,d$ commute with $\alpha,\beta,\gamma,\delta$. For this purpose, we apply $v$ to the following equality:
$$\xi_{10}\xi_{01}=\xi_{01}\xi_{10}$$

We obtain in this way an equality between two sums having 16 terms each, and by using 
again the eigenspace preservation condition we get the following formulae relating the corresponding 32 products $a\alpha,\alpha a$, and so on:
$$a\alpha=\alpha a=0\quad,\quad b\beta=\beta b=0$$
$$c\gamma=\gamma c=0\quad,\quad d\delta=\delta d=0$$
$$a\gamma+c\alpha+b\delta+d\beta=0\quad,\quad 
\alpha c+\gamma a+\beta d+\delta b=0$$
$$a\beta+b\alpha=\alpha b+\beta a\quad,\quad
b\gamma+c\beta=\beta c+\gamma b$$ 
$$c\delta+d\gamma=\gamma d+\delta c\quad,\quad
a\delta+d\alpha=\alpha d+\delta a$$  

(7) Now observe that multiplying the first equality in the third row on the left by $a$ and on the right by $\gamma$ gives $a^2\gamma^2 =0$, and by taking adjoints we get $\gamma a=0$. The same method applies to the other 7 products involved in the third row, so all 8 products involved in the third row vanish. That is, we have the following formulae:
$$a\gamma=c\alpha=b\delta=d\beta=\alpha c=\gamma a=\beta d=\delta b=0$$

(8) We use now the first equality in the fourth row. Multiplying it on the left by $a$ gives $a^2\beta=a\beta a$, and multiplying it on the right by $a$ gives $a\beta a=\beta a^2$. Thus we get $a^2\beta=\beta a^2$. On the other hand from $a^3+b^3+c^3+d^3=1$ we get $a^4=a$, so:
$$a\beta=a^4 \beta=a^2a^2 \beta=\beta a^2a^2=\beta a$$

Finally, one can show in a similar manner that the missing commutation formulae $a\delta = \delta a$ and so on, hold as well. Thus the algebra $A$ is commutative, as desired.    
\end{proof}

As a second graph which is resistent to a routine product study, we have the Petersen graph $P_{10}$. In order to explain the computation here, done by Schmidt in \cite{sc1}, we will need a number of preliminaries. Let us start with the following notion, from \cite{bi1}:

\index{Bichon construction}
\index{reduced quantum symmetry group}

\begin{definition}
The reduced quantum automorphism group of $X$ is given by
$$C(G^*(X))=C(G^+(X))\Big/\left<u_{ij}u_{kl}=u_{kl}u_{ij}\Big|\forall i-k,j-l\right>$$
with $i-j$ standing as usual for the fact that $i,j$ are connected by an edge.
\end{definition}

As explained by Bichon in \cite{bi1}, the above construction produces indeed a quantum group $G^*(X)$, which sits as an intermediate subgroup, as follows:
$$G(X)\subset G^*(X)\subset G^+(X)$$

There are many things that can be said about this construction, but in what concerns us, we will rather use it as a technical tool. Following Schmidt \cite{sc1}, we have:

\index{strongly regular}

\begin{proposition}
Assume that a regular graph $X$ is strongly regular, with parameters $\lambda=0$ and $\mu=1$, in the sense that:
\begin{enumerate}
\item $i-j$ implies that $i,j$ have $\lambda$ common neighbors.

\item $i\not\!\!-\,j$ implies that $i,j$ have $\mu$ common neighbors.
\end{enumerate}
The quantum group inclusion $G^*(X)\subset G^+(X)$ is then an isomorphism.
\end{proposition}

\begin{proof}
This is something quite tricky, the idea being as follows:

\medskip

(1) First of all, regarding the statement, a graph is called regular, with valence $k$, when each vertex has exactly $k$ neighbors. Then we have the notion of strong regularity, given by the conditions (1,2) in the statement. And finally we have the notion of strong regularity with parameters $\lambda=0,\mu=1$, that the statement is about, and with as main example here $P_{10}$, which is 3-regular, and strongly regular with $\lambda=0,\mu=1$.

\medskip

(2) Regarding now the proof, we must prove that the following commutation relation holds, with $u$ being the magic unitary of the quantum group $G^+(X)$:
$$u_{ij}u_{kl}=u_{kl}u_{ij}\ ,\ \forall i-k,j-l$$

(3) But for this purpose, we can use the $\lambda=0,\mu=1$ strong regularity of our graph, by inserting some neighbors into our computation. To be more precise, we have:
\begin{eqnarray*}
u_{ij}u_{kl}
&=&u_{ij}u_{kl}\sum_{s-l}u_{is}\\
&=&u_{ij}u_{kl}u_{ij}+\sum_{s-l,s\neq j}u_{ij}u_{kl}u_{is}\\
&=&u_{ij}u_{kl}u_{ij}+\sum_{s-l,s\neq j}u_{ij}\left(\sum_au_{ka}\right)u_{is}\\
&=&u_{ij}u_{kl}u_{ij}+\sum_{s-l,s\neq j}u_{ij}u_{is}\\
&=&u_{ij}u_{kl}u_{ij}
\end{eqnarray*}

(4) But this gives the result. Indeed, we conclude from this that $u_{ij}u_{kl}$ is self-adjoint, and so, by conjugating, that we have $u_{ij}u_{kl}=u_{kl}u_{ij}$, as desired.
\end{proof}

In the particular case of the Petersen graph $P_{10}$, which in addition is 3-regular, we can further build on the above result, and still following Schmidt \cite{sc1}, we have:

\begin{theorem}
The Petersen graph has no quantum symmetry,
$$G^+(P_{10})=G(P_{10})=S_5$$
with $S_5$ acting in the obvious way.
\end{theorem}

\begin{proof}
In view of Proposition 14.36, we must prove that the following commutation relation holds, with $u$ being the magic unitary of the quantum group $G^+(P_{10})$:
$$u_{ij}u_{kl}=u_{kl}u_{ij}\ ,\ \forall i\not\!\!-\,k,j\not\!\!-\,l$$

We can assume $i\neq k$, $j\neq l$. Now if we denote by $s,t$ the unique vertices having the property $i-s,k-s$ and $j-t,l-t$, a routine study shows that we have:
$$u_{ij}u_{kl}=u_{ij}u_{st}u_{kl}$$

With this in hand, if we denote by $q$ the third neighbor of $t$, we obtain:
\begin{eqnarray*}
u_{ij}u_{kl}
&=&u_{ij}u_{st}u_{kl}(u_{ij}+u_{il}+u_{iq})\\
&=&u_{ij}u_{st}u_{kl}u_{ij}+0+0\\
&=&u_{ij}u_{st}u_{kl}u_{ij}\\
&=&u_{ij}u_{kl}u_{ij}
\end{eqnarray*}

Thus the element $u_{ij}u_{kl}$ is self-adjoint, and we obtain, as desired:
$$u_{ij}u_{kl}=u_{kl}u_{ij}$$

As for the fact that the usual symmetry group is $S_5$, this is something that we know well from chapter 10, coming from the Kneser graph picture of $P_{10}$.
\end{proof}

As an application of this, we have the following classification table from \cite{bbi}, improved by using \cite{sc1}, containing all the vertex-transitive graphs of order $\leq 11$ modulo complementation, with their classical and quantum symmetry groups:
\begin{center}\begin{tabular}[t]{|l|l|l|l|}
\hline Order&Graph&Classical group&Quantum group\\ 
\hline\hline 2&$K_2$&$\mathbb Z_2$&$\mathbb  Z_2$\\
\hline\hline 3&$K_3$&$S_3$&$S_3$\\ 
\hline\hline 4&$2K_2$&$H_2$&$H_2^+$\\
\hline 4&$K_4$&$S_4$&$S_4^+$\\
\hline\hline 5&$C_5$&$D_5$&$D_5$\\ 
\hline5&$K_5$&$S_5$&$S_5^+$\\ 
\hline\hline 6&$C_6$&$D_6$&$D_6$\\ 
\hline 6&$2K_3$&$S_3\wr\mathbb Z_2$&$S_3{\,\wr_*\,}\mathbb Z_2$\\
\hline 6&$3K_2$&$H_3$&$H_3^+$\\ 
\hline 6&$K_6$&$S_6$&$S_6^+$\\
\hline\hline 7&$C_7$&$D_7$&$D_7$\\ 
\hline7&$K_7$&$S_7$&$S_7^+$\\ 
\hline\hline 8&$C_8$, $C_8^+$&$D_8$&$D_8$\\
\hline 8&$P(C_4)$& $H_3$&$S_4^+\times \mathbb Z_2$\\ 
\hline 8&$2K_4$&$S_4\wr \mathbb Z_2$&$S_4^+{\,\wr_*\,}\mathbb Z_2$\\
\hline 8&$2C_4$& $H_2\wr\mathbb Z_2$ & $H_2^+{\,\wr_*\,}\mathbb Z_2$\\ 
\hline 8&$4K_2$&$H_4$&$H_4^+$ \\ \hline 8&$K_8$&$S_8$&$S_8^+$\\
\hline\hline 9&$C_9$, $C_9^3$&$D_9$&$D_9$\\ 
\hline 9 & $K_3\times K_3$&$S_3\wr\mathbb Z_2$&$S_3\wr\mathbb Z_2$\\ 
\hline 9&$3K_3$&$S_3\wr S_3$&$S_3{\,\wr_*\,}S_3$\\ 
\hline9&$K_9$&$S_9$&$S_9^+$\\ 
\hline\hline 10&$C_{10}$, $C_{10}^2$, $C_{10}^+$, $P(C_5)$&$D_{10}$&$D_{10}$\\ 
\hline 10 &$P(K_5)$&$S_5\times\mathbb  Z_2$&$S_5^+\times\mathbb  Z_2$\\ 
\hline10&$C_{10}^4$&$\mathbb Z_2\wr D_5$&$\mathbb Z_2{\,\wr_*\,}D_5$\\ 
\hline10&$2C_5$&$D_5\wr\mathbb  Z_2$&$D_5{\,\wr_*\,}\mathbb Z_2$\\ 
\hline10&$2K_{5}$&$S_5\wr\mathbb  Z_2$&$S_5^+{\,\wr_*\,}\mathbb Z_2$\\ 
\hline10&$5K_2$&$H_5$&$H_5^+$\\ 
\hline10&$K_{10}$&$S_{10}$&$S_{10}^+$\\ 
\hline10&$P_{10}$&$S_5$&$S_5$\\
\hline\hline 11&$C_{11}$, $C_{11}^2$, $C_{11}^3$&$D_{11}$&$D_{11}$\\ 
\hline11&$K_{11}$&$S_{11}$&$S_{11}^+$\\ 
\hline
\end{tabular}\end{center}

Here $K$ denote the complete graphs, $C$ the cycles with chords, and $P$ stands for prisms. Moreover, by using more advanced techniques, the above table can be considerably extended. For more on all this, we refer to Schmidt's papers \cite{sc1}, \cite{sc2}, \cite{sc3}.

\section*{14e. Exercises}

We had a quite technical algebraic chapter here, in the hope that you survived all this, and here are a few exercises, in relation with the above:

\begin{exercise}
Clarify all the details in relation with free wreath products.
\end{exercise}

\begin{exercise}
Work out some basic applications of our usual product results.
\end{exercise}

\begin{exercise}
Work out basic applications of our lexicographic product results.
\end{exercise}

\begin{exercise}
Work out colored versions of our lexicographic product results.
\end{exercise}

\begin{exercise}
Clarify the easiness property of the complex reflection groups.
\end{exercise}

\begin{exercise}
Discuss the representations of the complex reflection groups.
\end{exercise}

\begin{exercise}
Discuss the character laws for the complex reflection groups.
\end{exercise}

\begin{exercise}
Come up with new twisting results, inspired by our $O_N^{-1}$ result.
\end{exercise}

As a bonus exercise, learn more about the usual complex reflection groups, and their classification into series and exceptional groups, from Shephard and Todd \cite{sto}.

\chapter{Advanced results}

\section*{15a. Orbits, orbitals}

We have seen in chapter 14 how to develop the basic theory of quantum symmetry groups of graphs $G^+(X)$, with this corresponding more or less to the knowledge of the subject from the mid to the end 00s. Many things have happened since, and in the remainder of this book we will attempt to navigate the subject, with some basics.

\bigskip

In the present chapter we would like to discuss a number of more advanced techniques for the computation of $G^+(X)$. More specifically, we would like to talk about:

\bigskip

(1) The orbits and orbitals of subgroups $G\subset S_N^+$, the point here being that an action $G\curvearrowright X$ is possible when the adjacency matrix of $X$ is constant on the orbitals of $G$. This is a key observation of Lupini, Man\v cinska and Roberson \cite{lmr}, heavily used ever since.

\bigskip

(2) The study of dual quantum groups $\widehat{\Gamma}\subset S_N^+$, and their actions on graphs $\widehat{\Gamma}\curvearrowright X$, when $\Gamma$ is finite, or classical, or arbitrary. Here the main ideas, which are actually related to the orbit and orbital problematics, are due to Bichon \cite{bi2} and McCarthy \cite{mcc}.

\bigskip

(3) The quantum group actions $G\curvearrowright X$ on graphs which are circulant, $\mathbb Z_N\curvearrowright X$, or more generally which admit an action of a finite abelian group, $A\curvearrowright X$. There is a lot of Fourier magic here, first discussed in my paper with Bichon and Chenevier \cite{bbg}.

\bigskip

(4) The construction and main properties of the quantum semigroup of quantum partial isometries $\widetilde{G}^+(X)\subset\widetilde{S}_N^+$, in analogy with what we know about $\widetilde{G}(X)\subset\widetilde{S}_N$. This is something more recent, not yet truly developed, and believed to be of interest.

\bigskip

So, this will be the plan for the present chapter, obviously many interesting things to be discussed, and I can only imagine that you are quite excited about this. Unfortunately, there is no easy way of compacting 20 years of math into one chapter, and we have:

\begin{disclaimer}
Contrary to the previous chapters, where new theory was accompanied by decent theorems using it, here we will be rather theoretical, explaining what is needed to know, in order to read the recent papers on the subject, and their theorems. 
\end{disclaimer}

And as a second disclaimer, this will be just half of the story, because we will still have to talk afterwards about planar algebra methods, which are actually things known since the early and mid 00s. We will do this in chapter 16 below.

\bigskip

But probably too much talking, leaving aside all this arithmetics of knowledge, let us get to work, according to our (1,2,3,4) plan above, for the present chapter. 

\bigskip

To start with, a useful tool for the study of the permutation groups $G\subset S_N$ are the orbits of the action $G\curvearrowright\{1,\ldots,N\}$. In the quantum permutation group case, $G\subset S_N^+$, following Bichon \cite{bi2}, the orbits can be introduced as follows:

\index{orbits}

\begin{theorem}
Given a closed subgroup $G\subset S_N^+$, with standard coordinates denoted $u_{ij}\in C(G)$, the following defines an equivalence relation on $\{1,\ldots,N\}$,
$$i\sim j\iff u_{ij}\neq0$$
that we call orbit decomposition associated to the action $G\curvearrowright\{1,\ldots,N\}$. In the classical case, $G\subset S_N$, this is the usual orbit equivalence.
\end{theorem}

\begin{proof}
We first check the fact that we have indeed an equivalence relation. The reflexivity axiom $i\sim i$ follows by using the counit, as follows:
$$\varepsilon(u_{ii})=1
\implies u_{ii}\neq0$$

The symmetry axiom $i\sim j\implies j\sim i$ follows by using the antipode:
$$S(u_{ji})=u_{ij}\implies[u_{ij}\neq0\implies u_{ji}\neq0]$$

As for the transitivity axiom $i\sim k,k\sim j\implies i\sim j$, this follows by using the comultiplication, and positivity. Consider indeed the following formula:
$$\Delta(u_{ij})=\sum_ku_{ik}\otimes u_{kj}$$

On the right we have a sum of projections, and we obtain from this, as desired:
\begin{eqnarray*}
u_{ik}\neq0,u_{kj}\neq0
&\implies&u_{ik}\otimes u_{kj}>0\\
&\implies&\Delta(u_{ij})>0\\
&\implies&u_{ij}\neq0
\end{eqnarray*}

Finally, in the classical case, where $G\subset S_N$, the standard coordinates are:
$$u_{ij}=\chi\left(\sigma\in G\Big|\sigma(j)=i\right)$$

Thus $u_{ij}\neq0$ means that $i,j$ must be in the same orbit, as claimed.
\end{proof}

Generally speaking, the theory from the classical case extends well to the quantum group setting, and we have in particular the following result, also from \cite{bi2}:

\index{fixed points}

\begin{theorem}
Given a closed subgroup $G\subset S_N^+$, with magic matrix $u=(u_{ij})$, consider the associated coaction map, on the space $X=\{1,\ldots,N\}$:
$$\Phi:C(X)\to C(X)\otimes C(G)\quad,\quad e_i\to\sum_je_j\otimes u_{ji}$$
The following three subalgebras of $C(X)$ are then equal,
$$Fix(u)=\left\{\xi\in C(X)\Big|u\xi=\xi\right\}$$
$$Fix(\Phi)=\left\{\xi\in C(X)\Big|\Phi(\xi)=\xi\otimes1\right\}$$
$$Fix(\sim)=\left\{\xi\in C(X)\Big|i\sim j\implies \xi_i=\xi_j\right\}$$
where $\sim$ is the orbit equivalence relation constructed in Theorem 15.2.
\end{theorem}

\begin{proof}
The fact that we have $Fix(u)=Fix(\Phi)$ is standard, with this being valid for any corepresentation $u=(u_{ij})$. Indeed, we first have the following computation:
\begin{eqnarray*}
\xi\in Fix(u)
&\iff&u\xi=\xi\\
&\iff&(u\xi)_j=\xi_j,\forall j\\
&\iff&\sum_iu_{ji}\xi_i=\xi_j,\forall j
\end{eqnarray*}

On the other hand, we have as well the following computation:
\begin{eqnarray*}
\xi\in Fix(\Phi)
&\iff&\Phi(\xi)=\xi\otimes1\\
&\iff&\sum_i\Phi(e_i)\xi_i=\xi\otimes1\\
&\iff&\sum_{ij}e_j\otimes u_{ji}\xi_i=\sum_je_j\otimes\xi_j\\
&\iff&\sum_iu_{ji}\xi_i=\xi_j,\forall j
\end{eqnarray*}

Thus we have $Fix(u)=Fix(\Phi)$, as claimed. Regarding now the equality of this algebra with $Fix(\sim)$, observe first that given a vector $\xi\in Fix(\sim)$, we have:
\begin{eqnarray*}
\sum_iu_{ji}\xi_i
&=&\sum_{i\sim j}u_{ji}\xi_i\\
&=&\sum_{i\sim j}u_{ji}\xi_j\\
&=&\sum_iu_{ji}\xi_j\\
&=&\xi_j
\end{eqnarray*}

Thus $\xi\in Fix(u)=Fix(\Phi)$. Finally, for the reverse inclusion, we know from Theorem 15.2 that the magic unitary $u=(u_{ij})$ is block-diagonal, with respect to the orbit decomposition there. But this shows that the algebra $Fix(u)=Fix(\Phi)$ decomposes as well with respect to the orbit decomposition, so in order to prove the result, we are left with a study in the transitive case. More specifically we must prove that if the action is transitive, then $u$ is irreducible, and this being clear, we obtain the result. See \cite{bi2}.
\end{proof}

We have as well a useful analytic result, as follows:

\begin{theorem}
Given a closed subgroup $G\subset S_N^+$, the matrix
$$P_{ij}=\int_Gu_{ij}$$
is the orthogonal projection onto $Fix(\sim)$, and determines the orbits of $G\curvearrowright\{1,\ldots,N\}$.
\end{theorem}

\begin{proof}
This follows from Theorem 15.3, and from the standard fact, coming from Peter-Weyl theory, that $P$ is the orthogonal projection onto $Fix(u)$. 
\end{proof}

As an application of the above, let us discuss now the notion of transitivity. We have here the following result, once again coming from \cite{bi2}:

\begin{theorem}
For a closed subgroup $G\subset S_N^+$, the following are equivalent:
\begin{enumerate}
\item $G$ is transitive, in the sense that $i\sim j$, for any $i,j$.

\item $Fix(u)=\mathbb C\xi$, where $\xi$ is the all-one vector.

\item $\int_Gu_{ij}=\frac{1}{N}$, for any $i,j$.
\end{enumerate}
In the classical case, $G\subset S_N$, this is the usual notion of transitivity.
\end{theorem}

\begin{proof}
This is well-known in the classical case. In general, the proof is as follows:

\medskip

$(1)\iff(2)$ We use the standard fact that the fixed point space of a corepresentation coincides with the fixed point space of the associated coaction:
$$Fix(u)=Fix(\Phi)$$

As explained in the beginning of this chapter, the fixed point space of the magic corepresentation $u=(u_{ij})$ has the following interpretation, in terms of orbits:
$$Fix(u)=\left\{\xi\in C(X)\Big|i\sim j\implies \xi(i)=\xi(j)\right\}$$

In particular, the transitivity condition corresponds to $Fix(u)=\mathbb C\xi$, as stated.

\medskip

$(2)\iff(3)$ This is clear from the general properties of the Haar integration, and more precisely from the fact that $(\int_Gu_{ij})_{ij}$ is the projection onto $Fix(u)$.
\end{proof}

Following Lupini, Man\v cinska and Roberson \cite{lmr}, let us discuss now the higher orbitals. Things are quite tricky here, and we have the following result, to start with:

\index{higher orbitals}

\begin{theorem}
For a closed aubgroup $G\subset S_N^+$, with magic unitary $u=(u_{ij})$, and $k\in\mathbb N$, the relation 
$$(i_1,\ldots,i_k)\sim(j_1,\ldots,j_k)\iff u_{i_1j_1}\ldots u_{i_kj_k}\neq0$$
is reflexive and symmetric, and is transitive at $k=1,2$. In the classical case, $G\subset S_N$, this relation is transitive at any $k\in\mathbb N$, and is the usual $k$-orbital equivalence.
\end{theorem}

\begin{proof}
This is known from \cite{lmr}, the proof being as follows:

\medskip

(1) The reflexivity of $\sim$ follows by using the counit, as follows:
\begin{eqnarray*}
\varepsilon(u_{i_ri_r})=1,\forall r
&\implies&\varepsilon(u_{i_1i_1}\ldots u_{i_ki_k})=1\\
&\implies&u_{i_1i_1}\ldots u_{i_ki_k}\neq0\\
&\implies&(i_1,\ldots,i_k)\sim(i_1,\ldots,i_k)
\end{eqnarray*}

(2) The symmetry follows by applying the antipode, and then the involution:
\begin{eqnarray*}
(i_1,\ldots,i_k)\sim(j_1,\ldots,j_k)
&\implies&u_{i_1j_1}\ldots u_{i_kj_k}\neq0\\
&\implies&u_{j_ki_k}\ldots u_{j_1i_1}\neq0\\
&\implies&u_{j_1i_1}\ldots u_{j_ki_k}\neq0\\
&\implies&(j_1,\ldots,j_k)\sim(i_1,\ldots,i_k)
\end{eqnarray*}

(3) The transitivity at $k=1,2$ is more tricky. Here we need to prove that:
$$u_{i_1j_1}\ldots u_{i_kj_k}\neq0\ ,\ u_{j_1l_1}\ldots u_{j_kl_k}\neq0\implies u_{i_1l_1}\ldots u_{i_kl_k}\neq0$$

In order to do so, we use the following formula:
$$\Delta(u_{i_1l_1}\ldots u_{i_kl_k})=\sum_{s_1\ldots s_k}u_{i_1s_1}\ldots u_{i_ks_k}\otimes u_{s_1l_1}\ldots u_{s_kl_k}$$

At $k=1$, we already know this. At $k=2$ now, we can use the following trick:
\begin{eqnarray*}
(u_{i_1j_1}\otimes u_{j_1l_1})\Delta(u_{i_1l_1}u_{i_2l_2})(u_{i_2j_2}\otimes u_{j_2l_2})
&=&\sum_{s_1s_2}u_{i_1j_1}u_{i_1s_1}u_{i_2s_2}u_{i_2j_2}\otimes u_{j_1l_1}u_{s_1l_1}u_{s_2l_2}u_{j_2l_2}\\
&=&u_{i_1j_1}u_{i_2j_2}\otimes u_{j_1l_1}u_{j_2l_2}
\end{eqnarray*}

Indeed, we obtain from this the following implication, as desired:
$$u_{i_1j_1}u_{i_2j_2}\neq0,u_{j_1l_1}u_{j_2l_2}\neq0\implies u_{i_1l_1}u_{i_2l_2}\neq0$$

(4) Finally, assume that we are in the classical case, $G\subset S_N$. We have:
$$u_{ij}=\chi\left(\sigma\in G\Big|\sigma(j)=i\right)$$

But this formula shows that we have the following equivalence:
$$u_{i_1j_1}\ldots u_{i_kj_k}\neq0\iff\exists \sigma\in G,\ \sigma(i_1)=j_1,\ldots,\sigma(i_k)=j_k$$

In other words, $(i_1,\ldots,i_k)\sim(j_1,\ldots,j_k)$ happens precisely when $(i_1,\ldots,i_k)$ and $(j_1,\ldots,j_k)$ are in the same $k$-orbital of $G$, and this gives the last assertion.
\end{proof}

The above result raises the question about what exactly happens at $k=3$, in relation with transitivity, and the answer here is negative in general. To be more precise, as explained by McCarthy in \cite{mcc}, there are closed subgroups $G\subset S_N^+$, as for instance the Kac-Paljutkin quantum group $G\subset S_4^+$, for which $\sim$ is not transitive at $k=3$.

\bigskip

In view of all this, we can only formulate a modest definition, as follows:

\begin{definition}
Given a closed subgroup $G\subset S_N^+$, consider the relation defined by: 
$$(i_1,\ldots,i_k)\sim(j_1,\ldots,j_k)\iff u_{i_1j_1}\ldots u_{i_kj_k}\neq0$$
\begin{enumerate}
\item At $k=1$, the equivalence classes of $\sim$ are called orbits of $G$.

\item At $k=2$, the equivalence classes of $\sim$ are called orbitals of $G$.

\item At $k\geq3$, if $\sim$ is transitive, we call its equivalence classes $k$-orbitals of $G$.
\end{enumerate}
\end{definition}

It is possible to say more things here, but generally speaking, the good theory remains at $k=1,2$. In what follows we will focus on the case $k=2$, where $\sim$ is given by:
$$(i,k)\sim (j,l)\iff u_{ij}u_{kl}\neq0$$

As a key theoretical result on the subject, again from \cite{lmr}, we have the following key analogue of Theorem 15.3, which makes a connection with the graph problematics:

\begin{theorem}
Given a closed subgroup $G\subset S_N^+$, with magic matrix $u=(u_{ij})$, consider the following vector space coaction map, where $X=\{1,\ldots,N\}$:
$$\Phi:C(X\times X)\to C(X\times X)\otimes C(G)\quad,\quad e_{ik}\to\sum_{jl}e_{jl}\otimes u_{ji}u_{lk}$$
The following three algebras are then isomorphic,
$$End(u)=\left\{d\in M_N(\mathbb C)\Big|du=ud\right\}$$
$$Fix(\Phi)=\left\{\xi\in C(X\times X)\Big|\Phi(\xi)=\xi\otimes1\right\}$$
$$Fix(\sim)=\left\{\xi\in C(X\times X)\Big|(i,k)\sim(j,l)\implies \xi_{ik}=\xi_{jl}\right\}$$
where $\sim$ is the orbital equivalence relation from Definition 15.7 (2).
\end{theorem}

\begin{proof}
This follows by doing some computations which are quite similar to those from the proof of Theorem 15.3, and we refer here to \cite{lmr}, for the details.
\end{proof}

As already mentioned, the above result makes a quite obvious connection with the graph problematics, the precise statement here being as follows:

\begin{theorem}
In order for a quantum permutation group $G\subset S_N^+$ to act on a graph $X$, having $N$ vertices, the adjacency matrix $d\in M_N(0,1)$ of the graph must be, when regarded as function on the set $\{1,\ldots,N\}^2$, constant on the orbitals of $G$. 
\end{theorem}

\begin{proof}
This follows indeed from the following isomorphism, from Theorem 15.8:
$$End(u)\simeq Fix(\sim)$$

For more on all this, details, examples, and applications too, we refer to \cite{lmr}.
\end{proof}

Finally, as a theoretical application of the theory of orbitals, as developed above, let us discuss now the notion of double transitivity. Following \cite{lmr}, we have:

\begin{definition}
Let $G\subset S_N^+$ be a closed subgroup, with magic unitary $u=(u_{ij})$, and consider as before the equivalence relation on $\{1,\ldots,N\}^2$ given by:
$$(i,k)\sim (j,l)\iff u_{ij}u_{kl}\neq0$$ 
\begin{enumerate}
\item The equivalence classes under $\sim$ are called orbitals of $G$.

\item $G$ is called doubly transitive when the action has two orbitals. 
\end{enumerate}
In other words, we call $G\subset S_N^+$ doubly transitive when $u_{ij}u_{kl}\neq0$, for any $i\neq k,j\neq l$.
\end{definition}

To be more precise, it is clear from definitions that the diagonal $D\subset\{1,\ldots,N\}^2$ is an orbital, and that its complement $D^c$ must be a union of orbitals. With this remark in hand, the meaning of (2) is that the orbitals must be $D,D^c$. 

\bigskip

In more analytic terms, we have the following result, also from \cite{lmr}:

\begin{theorem}
For a doubly transitive subgroup $G\subset S_N^+$, we have:
$$\int_Gu_{ij}u_{kl}=\begin{cases}
\frac{1}{N}&{\rm if}\ i=k,j=l\\
0&{\rm if}\ i=k,j\neq l\ {\rm or}\ i\neq k,j=l\\
\frac{1}{N(N-1)}&{\rm if}\ i\neq k,j\neq l
\end{cases}$$
Moreover, this formula characterizes the double transitivity.
\end{theorem}

\begin{proof}
We use the standard fact, from \cite{wo1}, that the integrals in the statement form the projection onto $Fix(u^{\otimes 2})$. Now if we assume that $G$ is doubly transitive, $Fix(u^{\otimes 2})$ has dimension 2, and therefore coincides with $Fix(u^{\otimes 2})$ for the usual symmetric group $S_N$. Thus the integrals in the statement coincide with those for the symmetric group $S_N$, which are given by the above formula. Finally, the converse is clear as well.
\end{proof}

So long for orbits, orbitals, and transitivity. As already mentioned, Theorem 15.9, which makes the connection with the actions on graphs $G\curvearrowright X$, is something quite far-reaching, and for the continuation of this, we refer to \cite{lmr} and subsequent papers.

\section*{15b. Dual formalism}

As another application of the orbit theory developed above, following Bichon \cite{bi2}, let us discuss now the group duals $\widehat{\Gamma}\subset S_N^+$. We first have the following result:

\begin{theorem}
Given a quotient group $\mathbb Z_{N_1}*\ldots*\mathbb Z_{N_k}\to\Gamma$, we have an embedding $\widehat{\Gamma}\subset S_N^+$, with $N=N_1+\ldots+N_k$, having the following properties:
\begin{enumerate}
\item This embedding appears by diagonally joining the embeddings $\widehat{\mathbb Z_{N_k}}\subset S_{N_k}^+$, and the corresponding magic matrix has blocks of sizes $N_1,\ldots,N_k$.

\item The equivalence relation on $X=\{1,\ldots,N\}$ coming from the orbits of the action $\widehat{\Gamma}\curvearrowright X$ appears by refining the partition $N=N_1+\ldots+N_k$.
\end{enumerate}
\end{theorem}

\begin{proof}
This is something elementary, the idea being as follows:

\medskip

(1) Given a quotient group $\mathbb Z_{N_1}*\ldots*\mathbb Z_{N_k}\to\Gamma$, we have indeed a standard embedding as follows, with $N=N_1+\ldots+N_k$, that we actually know well since chapter 13:
\begin{eqnarray*}
\widehat{\Gamma}
&\subset&\widehat{\mathbb Z_{N_1}*\ldots*\mathbb Z_{N_k}}
=\widehat{\mathbb Z_{N_1}}\,\hat{*}\,\ldots\,\hat{*}\,\widehat{\mathbb Z_{N_k}}\\
&\simeq&\mathbb Z_{N_1}\,\hat{*}\,\ldots\,\hat{*}\,\mathbb Z_{N_k}
\subset S_{N_1}\,\hat{*}\,\ldots\,\hat{*}\,S_{N_k}\\
&\subset&S_{N_1}^+\,\hat{*}\,\ldots\,\hat{*}\,S_{N_k}^+
\subset S_N^+
\end{eqnarray*}

(2) Regarding the magic matrix, our claim is that this is as follows, $F_N=\frac{1}{\sqrt{N}}(w_N^{ij})$ with $w_N=e^{2\pi i/N}$ being Fourier matrices, and $g_l$ being the standard generator of $\mathbb Z_{N_l}$:
$$u=\begin{pmatrix}
F_{N_1}I_1F_{N_1}^*\\
&\ddots\\
&&F_{N_k}I_kF_{N_k}^*
\end{pmatrix}
\quad,\quad 
I_l=\begin{pmatrix}
1\\
&g_l\\
&&\ddots\\
&&&g_l^{N_l-1}
\end{pmatrix}$$

(3) Indeed, let us recall that the magic matrix for $\mathbb Z_N\subset S_N\subset S_N^+$ is given by:
$$v_{ij}
=\chi\left(\sigma\in\mathbb Z_N\Big|\sigma(j)=i\right)
=\delta_{i-j}$$

Let us apply now the Fourier transform. According to our usual Pontrjagin duality conventions, we have a pair of inverse isomorphisms, as follows:
$$\Phi:C(\mathbb Z_N)\to C^*(\mathbb Z_N)\quad,\quad\delta_i\to\frac{1}{N}\sum_kw^{ik}g^k$$
$$\Psi:C^*(\mathbb Z_N)\to C(\mathbb Z_N)\quad,\quad g^i\to\sum_kw^{-ik}\delta_k$$

Here $w=e^{2\pi i/N}$, and we use the standard Fourier analysis convention that the indices are $0,1,\ldots,N-1$. With $F=\frac{1}{\sqrt{N}}(w^{ij})$ and $I=diag(g^i)$ as above, we have:
\begin{eqnarray*}
u_{ij}
&=&\Phi(v_{ij})\\
&=&\frac{1}{N}\sum_kw^{(i-j)k}g^k\\
&=&\frac{1}{N}\sum_kw^{ik}g^kw^{-jk}\\
&=&(FIF^*)_{ij}
\end{eqnarray*}

Thus, the magic matrix that we are looking for is $u=FIF^*$, as claimed.

\medskip

(4) Finally, the second assertion in the statement is clear from the fact that $u$ is block-diagonal, with blocks corresponding to the partition $N=N_1+\ldots+N_k$.
\end{proof}

As a first comment on the above result, not all group dual subgroups $\widehat{\Gamma}\subset S_N^+$ appear as above, a well-known counterexample here being the Klein group:
$$K=\mathbb Z_2\times\mathbb Z_2\subset S_4\subset S_4^+$$

Indeed, with $K=\{1,a,b,c\}$, where $c=ab$, consider the embedding $K\subset S_4$ given by $a=(12)(34)$, $b=(13)(24)$, $c=(14)(23)$. The corresponding magic matrix is:
$$u=\begin{pmatrix}
\delta_1&\delta_a&\delta_b&\delta_c\\
\delta_a&\delta_1&\delta_c&\delta_b\\
\delta_b&\delta_c&\delta_1&\delta_a\\
\delta_c&\delta_b&\delta_a&\delta_1
\end{pmatrix}\in M_4(C(K))$$

Now since this matrix is not block-diagonal, the only choice for $K=\widehat{K}$ to appear as in Theorem 15.12 would be via a quotient map $\mathbb Z_4\to K$, which is impossible. As a second comment now on Theorem 15.12, in the second assertion there we really have a possible refining operation, as shown by the example provided by the trivial group, namely:
$$\mathbb Z_{N_1}*\ldots*\mathbb Z_{N_k}\to\{1\}$$

In order to further discuss all this, let us first enlarge the attention to the group dual subgroups $\widehat{\Gamma}\subset G$ of an arbitrary closed subgroup $G\subset U_N^+$. We have here: 

\begin{theorem}
Given a closed subgroup $G\subset U_N^+$ and a matrix $Q\in U_N$, we let $T_Q\subset G$ be the diagonal torus of $G$, with fundamental representation spinned by $Q$:
$$C(T_Q)=C(G)\Big/\left<(QuQ^*)_{ij}=0\Big|\forall i\neq j\right>$$
This torus is then a group dual, $T_Q=\widehat{\Lambda}_Q$, where $\Lambda_Q=<g_1,\ldots,g_N>$ is the discrete group generated by the elements $g_i=(QuQ^*)_{ii}$, which are unitaries inside $C(T_Q)$.
\end{theorem}

\begin{proof}
Since $v=QuQ^*$ is a unitary corepresentation, its diagonal entries $g_i=v_{ii}$, when regarded inside $C(T_Q)$, are unitaries, and satisfy:
$$\Delta(g_i)=g_i\otimes g_i$$

Thus $C(T_Q)$ is a group algebra, and more specifically we have $C(T_Q)=C^*(\Lambda_Q)$, where $\Lambda_Q=<g_1,\ldots,g_N>$ is the group in the statement, and this gives the result.
\end{proof}

Generally speaking, the above family $\{T_Q|Q\in U_N\}$ can be thought of as being a kind of ``maximal torus'' for $G\subset U_N^+$. Now back to quantum permutations, we have:

\begin{theorem}
For the quantum permutation group $S_N^+$, the discrete group quotient $F_N\to\Lambda_Q$ with $Q\in U_N$ comes from the following relations:
$$\begin{cases}
g_i=1&{\rm if}\ \sum_lQ_{il}\neq 0\\
g_ig_j=1&{\rm if}\ \sum_lQ_{il}Q_{jl}\neq 0\\ 
g_ig_jg_k=1&{\rm if}\ \sum_lQ_{il}Q_{jl}Q_{kl}\neq 0
\end{cases}$$
Also, given a decomposition $N=N_1+\ldots+N_k$, for the matrix $Q=diag(F_{N_1},\ldots,F_{N_k})$, where $F_N=\frac{1}{\sqrt{N}}(\xi^{ij})_{ij}$ with $\xi=e^{2\pi i/N}$ is the Fourier matrix, we obtain
$$\Lambda_Q=\mathbb Z_{N_1}*\ldots*\mathbb Z_{N_k}$$
with dual embedded into $S_N^+$ in a standard way, as in Theorem 15.12.
\end{theorem}

\begin{proof}
This can be proved by a direct computation, as follows:

\medskip

(1) Fix a unitary matrix $Q\in U_N$, and consider the following quantities:
$$\begin{cases}
c_i=\sum_lQ_{il}\\
c_{ij}=\sum_lQ_{il}Q_{jl}\\
d_{ijk}=\sum_l\bar{Q}_{il}\bar{Q}_{jl}Q_{kl}
\end{cases}$$

We write $w=QvQ^*$, where $v$ is the fundamental corepresentation of $C(S_N^+)$. Assume $X\simeq\{1,\ldots,N\}$, and let $\alpha$ be the coaction of $C(S_N^+)$ on $C(X)$. Let us set:
$$\varphi_i=\sum_l\bar{Q}_{il}\delta_l\in C(X)$$

Also, let $g_i=(QvQ^*)_{ii}\in C^*(\Lambda_Q)$. If $\beta$ is the restriction of $\alpha$ to $C^*(\Lambda_Q)$, then:
$$\beta(\varphi_i)=\varphi_i\otimes g_i$$

(2) Now recall that $C(X)$ is the universal $C^*$-algebra generated by elements $\delta_1,\ldots,\delta_N$ which are pairwise orthogonal projections. Writing these conditions in terms of the linearly independent elements $\varphi_i$ by means of the formulae $\delta_i=\sum_lQ_{il}\varphi_l$, we find that the universal relations for $C(X)$ in terms of the elements $\varphi_i$ are as follows:
$$\begin{cases}
\sum_ic_i\varphi_i=1\\
\varphi_i^*=\sum_jc_{ij}\varphi_j\\
\varphi_i\varphi_j=\sum_kd_{ijk}\varphi_k
\end{cases}$$

(3) Let $\widetilde{\Lambda}_Q$ be the group in the statement. Since $\beta$ preserves these relations, we get:
$$\begin{cases}
c_i(g_i-1)=0\\
c_{ij}(g_ig_j-1)=0\\
d_{ijk}(g_ig_j-g_k)=0
\end{cases}$$

We conclude from this that $\Lambda_Q$ is a quotient of $\widetilde{\Lambda}_Q$. On the other hand, it is immediate that we have a coaction map as follows:
$$C(X)\to C(X)\otimes C^*(\widetilde{\Lambda}_Q)$$

Thus $C(\widetilde{\Lambda}_Q)$ is a quotient of $C(S_N^+)$. Since $w$ is the fundamental corepresentation of $S_N^+$ with respect to the basis $\{\varphi_i\}$, it follows that the generator $w_{ii}$ is sent to $\widetilde{g}_i\in\widetilde{\Lambda}_Q$, while $w_{ij}$ is sent to zero. We conclude that $\widetilde{\Lambda}_Q$ is a quotient of $\Lambda_Q$. Since the above quotient maps send generators on generators, we conclude that $\Lambda_Q=\widetilde{\Lambda}_Q$, as desired.

\medskip

(4) We apply the result found in (3), with the $N$-element set $X$ there being:
$$X=\mathbb Z_{N_1}\sqcup\ldots\sqcup\mathbb Z_{N_k}$$

With this choice, we have $c_i=\delta_{i0}$ for any $i$. Also, we have $c_{ij}=0$, unless $i,j,k$ belong to the same block to $Q$, in which case $c_{ij}=\delta_{i+j,0}$, and also $d_{ijk} =0$, unless $i,j,k$ belong to the same block of $Q$, in which case $d_{ijk}=\delta_{i+j,k}$. We conclude from this that $\Lambda_Q$ is the free product of $k$ groups which have generating relations as follows:
$$g_ig_j=g_{i+j}\quad,\quad g_i^{-1}=g_{-i}$$

But this shows that our group is $\Lambda_Q=\mathbb Z_{N_1}*\ldots*\mathbb Z_{N_k}$, as stated.
\end{proof}

In relation now with actions on graphs, let us start with the following simple fact:

\begin{proposition}
In order for a closed subgroup $G\subset U_K^+$ to appear as $G=G^+(X)$, for a certain graph $X$ having $N$ vertices, the following must happen:
\begin{enumerate}
\item We must have a representation $G\subset U_N^+$.

\item This representation must be magic, $G\subset S_N^+$.

\item We must have a graph $X$ having $N$ vertices, such that $d\in End(u)$.

\item $X$ must be in fact such that the Tannakian category of $G$ is precisely $<d>$.
\end{enumerate}
\end{proposition}

\begin{proof}
This is more of an empty statement, coming from the definition of the quantum automorphism group $G^+(X)$, as formulated in chapter 14.
\end{proof}

The above result remains something quite philosophical. In the group dual case, that we will be interested in now, we can combine it with the following result:

\begin{proposition}
Given a discrete group $\Gamma=<g_1,\ldots,g_N>$, embed diagonally $\widehat{\Gamma}\subset U_N^+$, via the unitary matrix $u=diag(g_1,\ldots,g_N)$. We have then the formula
$$Hom(u^{\otimes k},u^{\otimes l})=\left\{T=(T_{j_1\ldots j_l,i_1\ldots i_k})\Big|g_{i_1}\ldots g_{i_k}\neq g_{j_1}\ldots g_{j_l}\implies T_{j_1\ldots j_l,i_1\ldots i_k}=0\right\}$$
and in particular, with $k=l=1$, we have the formula
$$End(u)=\left\{T=(T_{ji})\Big|g_i\neq g_j\implies T_{ji}=0\right\}$$
with the linear maps being identified with the corresponding scalar matrices.
\end{proposition}

\begin{proof}
This is indeed elementary, with the first assertion coming from:
\begin{eqnarray*}
T\in Hom(u^{\otimes k},u^{\otimes l})
&\iff&Tu^{\otimes k}=u^{\otimes l}T\\
&\iff&(Tu^{\otimes k})_{j_1\ldots j_l,i_1\ldots i_k}=(u^{\otimes l}T)_{j_1\ldots j_l,i_1\ldots i_k}\\
&\iff&T_{j_1\ldots j_l,i_1\ldots i_k}g_{i_1}\ldots g_{i_k}=g_{j_1}\ldots g_{j_l}T_{j_1\ldots j_l,i_1\ldots i_k}\\
&\iff&T_{j_1\ldots j_l,i_1\ldots i_k}(g_{i_1}\ldots g_{i_k}-g_{j_1}\ldots g_{j_l})=0
\end{eqnarray*}

As for the second assertion, this follows from the first one.
\end{proof}

Still in the group dual setting, we have now, refining Proposition 15.15:

\begin{theorem}
In order for a group dual $\widehat{\Gamma}$ as above to appear as $G=G^+(X)$, for a certain graph $X$ having $N$ vertices, the following must happen:
\begin{enumerate}
\item First, we need a quotient map $\mathbb Z_{N_1}*\ldots*\mathbb Z_{N_k}\to\Gamma$.

\item Let $u=diag(I_1,\ldots,I_k)$, with $I_l=diag(\mathbb Z_{N_l})$, for any $l$.

\item Consider also the matrix $F=diag(F_{N_1},\ldots,F_{N_k})$.

\item We must then have a graph $X$ having $N$ vertices.

\item This graph must be such that $F^*dF\neq0\implies I_i=I_j$.

\item In fact, $<F^*dF>$ must be the category in Proposition 15.16.
\end{enumerate}
\end{theorem}

\begin{proof}
This is something rather informal, which follows from the above.
\end{proof}

Going ahead, in connection with the Fourier tori from Theorem 15.14, we have:

\index{Fourier torus}
\index{Fourier liberation}

\begin{proposition}
The Fourier tori of $G^+(X)$ are the biggest quotients
$$\mathbb Z_{N_1}*\ldots*\mathbb Z_{N_k}\to\Gamma$$
whose duals act on the graph, $\widehat{\Gamma}\curvearrowright X$.
\end{proposition}

\begin{proof}
We have indeed the following computation, at $F=1$:
\begin{eqnarray*}
C(T_1(G^+(X)))
&=&C(G^+(X))/<u_{ij}=0,\forall i\neq j>\\
&=&[C(S_N^+)/<[d,u]=0>]/<u_{ij}=0,\forall i\neq j>\\
&=&[C(S_N^+)/<u_{ij}=0,\forall i\neq j>]/<[d,u]=0>\\
&=&C(T_1(S_N^+))/<[d,u]=0>
\end{eqnarray*}

Thus, we obtain the result, with the remark that the quotient that we are interested in appears via relations of type $d_{ij}=1\implies g_i=g_j$. The proof in general is similar.
\end{proof}

In order to formulate our main result, let us call $G\subset G^+$ a Fourier liberation when $G^+$ is generated by $G$, and its Fourier tori. With this convention, we have:

\begin{theorem}
Consider the following conditions:
\begin{enumerate}
\item We have $G(X)=G^+(X)$.

\item $G(X)\subset G^+(X)$ is a Fourier liberation.

\item $\widehat{\Gamma}\curvearrowright X$ implies that $\Gamma$ is abelian.
\end{enumerate}
We have then the equivalence $(1)\iff(2)+(3)$.
\end{theorem}

\begin{proof}
This is something elementary. Indeed, the implications $(1)\implies(2,3)$ are trivial. As for $(2,3)\implies(1)$, assuming $G(X)\neq G^+(X)$, from (2) we know that $G^+(X)$ has at least one non-classical Fourier torus, and this contradicts (3), as desired.
\end{proof}

As an application of this, we have the following elementary result, which is a particular case of more general and advanced results regarding the random graphs, from \cite{lmr}:

\begin{theorem}
For a finite graph $X$, the probability for having an action
$$\widehat{\Gamma}\curvearrowright X$$
with $\Gamma$ being a non-abelian group goes to $0$ with $|X|\to\infty$.
\end{theorem}

\begin{proof}
Observe first that in the cyclic case, where $F=F_N$ is a usual Fourier matrix, associated to a cyclic group $\mathbb Z_N$, we have the following formula, with $w=e^{2\pi i/N}$:
$$(F^*dF)_{ij}
=\sum_{kl}(F^*)_{ik}d_{kl}F_{lj}
=\sum_{kl}w^{lj-ik}d_{kl}
=\sum_{k\sim l}w^{lj-ik}$$

In the general setting now, where we have a quotient map $\mathbb Z_{N_1}*\ldots*\mathbb Z_{N_k}\to\Gamma$, with $N_1+\ldots+N_k=N$, the computation is similar, as follows, with $w_i=e^{2\pi i/N_i}$:
$$(F^*dF)_{ij}
=\sum_{kl}(F^*)_{ik}d_{kl}F_{lj}
=\sum_{k\sim l}(F^*)_{ik}F_{lj}
=\sum_{k:i,l:j,k\sim l}(w_{N_i})^{-ik}(w_{N_j})^{lj}$$

The point now is that with the partition $N_1+\ldots+N_k=N$ fixed, and with $d\in M_N(0,1)$ being random, we have $(F^*dF)_{ij}\neq 0$ almost everywhere in the $N\to\infty$ limit, and so we obtain $I_i=I_j$ almost everywhere, and so abelianity of $\Gamma$, with $N\to\infty$.
\end{proof}

Many other things can be said, as a continuation of the above, and we refer here to \cite{lmr}, \cite{mcc} for an introduction, and to the recent literature on the subject, for more.

\section*{15c. Circulant graphs}

Changing topics, but still obsessed by Fourier analysis, let us discuss now, following \cite{bbg}, some sharp results in the circulant graph case. Let us start with:

\index{type of circulant graph}

\begin{definition}
Associated to any circulant graph $X$ having $N$ vertices are:
\begin{enumerate}
\item The set $S\subset\mathbb Z_N$ given by $i-_Xj\iff j-i\in S$.

\item The group $E\subset\mathbb Z_N^*$ consisting of elements $g$ such that $gS=S$.

\item The number $k=|E|$, called type of $X$.
\end{enumerate}
\end{definition}

The interest in the type $k$ is that this is the good parameter measuring the complexity of the spectral theory of $X$, as we will soon see. To start with, here are a few basic examples, and properties of the type, with $\varphi$ being the Euler function:

\medskip

(1) The type can be $2,4,6,8,\ldots$ This is because $\{\pm 1\}\subset E$.

\medskip

(2) $C_N$ is of type $2$. Indeed, we have $S=\{\pm 1\}$, $E=\{\pm 1\}$.

\medskip

(3) $K_N$ is of type $\varphi(N)$. Indeed, here $S=\emptyset$, $E=\mathbb Z_N^*$.

\medskip

Let us first discuss the spectral theory of the circulant graphs. In what follows $X$ will be a circulant graph having $p$ vertices, with $p$ prime. We denote by $\xi$ the column vector $(1,w,w^2,\ldots ,w^{p-1})$, where $w=e^{2\pi i/p}$. With this convention, we have:

\index{Dedekind theorem}
\index{Galois extension}

\begin{theorem}
The eigenspaces of $d$ are given by $V_0={\mathbb C}1$ and
$$V_x=\bigoplus_{a\in E}{\mathbb C}\,\xi^{xa}$$
with $x\in \mathbb Z_p^*$. Moreover, we have $V_x=V_y$ if and only if $xE=yE$.
\end{theorem}

\begin{proof}
Since $d$ is circulant, we have $d(\xi^x)=f(x)\xi^x$, with $f:{\mathbb Z}_p\to{\mathbb C}$ being:
$$f(x)=\sum_{t\in S}w^{xt}$$

Let $K={\mathbb Q}(w)$ and let $H$ be the Galois group of the Galois extension $\mathbb Q \subset K$. We have then a well-known group isomorphism, as follows:
$$\mathbb Z_p^*\simeq H\quad,\quad x\to s_x=[w\to w^x]$$

Also, we know from a standard theorem of Dedekind that the family $\{s_x\mid x\in{\mathbb Z}_p^*\}$ is free in ${\rm End}_{\mathbb Q}(K)$. Now for $x,y\in \mathbb Z_p^*$ consider the following operator:
$$L = \sum_{t \in S} s_{xt} - \sum_{t \in S} s_{yt} \in
End_{\mathbb Q}(K)$$

We have $L({w}) = f(x)-f(y)$, and since $L$ commutes with the action of $H$, we have:
$$L=0 \iff L({w}) =0 \iff f(x)=f(y)$$

By linear independence of the family $\{s_x\mid x\in \mathbb Z_p^*\}$ we get:
$$f(x) = f(y) \iff xS=yS \iff xE=yE$$

It follows that $d$ has precisely $1+(p-1)/k$ distinct eigenvalues, the corresponding
eigenspaces being those in the statement.
\end{proof}

Consider now a commutative ring $(R,+,\cdot)$. We denote by $R^*$ the group of invertibles, and we assume $2\in R^*$. A subgroup $G\subset R^*$ is called even if $-1\in G$. We have:

\index{even subgroup}
\index{2-maximal}
\index{hexagonal solution}

\begin{definition}
An even subgroup $G\subset R^*$ is called $2$-maximal if, inside $G$:
$$a-b=2(c-d)\implies a=\pm b$$
We call $a=b,c=d$ trivial solutions, and $a=-b=c-d$ hexagonal solutions. 
\end{definition}

To be more precise, in what regards our terminology, consider the group $G\subset{\mathbb C}$ formed by $k$-th roots of unity, with $k$ even. An equation of the form $a-b=2(c-d)$ with $a,b,c,d\in G$ says that the diagonals $a-b$ and $c-d$ must be parallel, and that the first one is twice as much as the second one. But this can happen only when $a,c,d,b$ are consecutive vertices of a regular hexagon, and here we have $a+b=0$.

\bigskip

The relation with our quantum symmetry considerations will come from:

\begin{proposition}
Assume that $R$ has the property $3\neq 0$, and consider a $2$-maximal subgroup $G\subset R^*$. Then, the following happen:
\begin{enumerate}
\item $2,3\not\in G$. 

\item $a+b=2c$ with $a,b,c\in G$ implies $a=b=c$. 

\item $a+2b=3c$ with $a,b,c\in G$ implies $a=b=c$.
\end{enumerate}
\end{proposition}

\begin{proof}
All these assertions are elementary, as follows:

\medskip

(1) This follows from the following formulae, which cannot hold in
$G$:
$$4-2=2(2-1)\quad,\quad 
3-(-1)=2(3-1)$$

Indeed, the first one would imply $4=\pm 2$, and the second one
would imply $3=\pm 1$. But from $2\in R^*$ and $3\neq 0$ we get
$2,4,6\neq 0$, contradiction.

\medskip

(2) We have $a-b=2(c-b)$. For a trivial solution we have
 $a=b=c$, and for a hexagonal
 solution we have $a+b=0$,
hence $c=0$, hence $0\in{G}$, contradiction.

\medskip

(3) We have $a-c=2(c-b)$. For a trivial solution we have
 $a=b=c$, and for a hexagonal
 solution we have $a+c=0$,
hence $b=-2a$, hence $2\in{G}$, contradiction.
\end{proof}

As a first result now, coming from this study, we have:

\index{no quantum symmetry}

\begin{theorem}
A circulant graph $X$, on $p\geq 5$ prime vertices, such that
$$E\subset {\mathbb Z}_p$$
is $2$-maximal, has no quantum symmetry, $G^+(X)=G(X)$.
\end{theorem}

\begin{proof}
This comes from the above results, via a long algebraic study, as follows:

\medskip

(1) We use Proposition 15.24, which ensures that $V_1,V_2,V_3$ are eigenspaces of $d$. By $2$-maximality of $E$, these eigenspaces are different. From the eigenspace preservation in Theorem 15.22 we get formulae of the following type, with $r_a,r_a',r_a''\in{\mathcal A}$:
$$\alpha (\xi)=\sum_{a\in E}\xi^a\otimes r_a\quad,\quad 
\alpha (\xi^2)=\sum_{a\in E}\xi^{2a}\otimes r_a'\quad,\quad
\alpha (\xi^3)=\sum_{a\in E}\xi^{3a}\otimes r_a''$$

(2) We take the square of the first relation, we compare with the formula of $\alpha(\xi^2)$, and we use $2$-maximality. We obtain in this way the following formula:
$$\alpha(\xi^2)=\sum_{c\in E}\xi^{2c}\otimes r_c^2$$

(3) We multiply now this relation by the formula of $\alpha(\xi)$, we compare with the formula of $\alpha(\xi^3)$, and we use $2$-maximality. We obtain the following formula:
$$\alpha(\xi^3)=\sum_{b\in E}\xi^{3b}\otimes r_b^3$$

(4) As a conclusion, the three formulae in the beginning are in fact as follows:
$$\alpha (\xi)=\sum_{a\in E}\xi^a\otimes r_a\quad,\quad 
\alpha(\xi^2)=\sum_{a\in E}\xi^{2a}\otimes r_a^2\quad,\quad 
\alpha(\xi^3)=\sum_{a\in E}\xi^{3a}\otimes r_a^3$$

(5) Our claim now is that for $a\neq b$, we have the following ``key formula'':
$$r_ar_b^3=0$$

Indeed, in order to prove this claim, consider the following equality:
$$\left(\sum_{a\in E}\xi^a\otimes r_a\right)
\left(\sum_{b\in E}\xi^{2b}\otimes r_b^2\right)=\sum_{c\in
E}\xi^{3c}\otimes r_c^3$$

By eliminating all $a=b$ terms, which produce the sum on the right, we get:
$$\sum\left\{r_ar_b^2\Big| a,b\in E,\,a\neq b,\,a+2b=x\right\}=0$$

(6) We fix now elements $a,b\in E$ satisfying $a\neq b$. We know from $2$-maximality that the equation $a+2b=a'+2b'$ with $a',b'\in E$ has at most one non-trivial solution, namely the hexagonal one, given by $a'=-a$ and $b'=a+b$. Now with $x=a+2b$, we get that the
above equality is in fact one of the following two equalities:
$$r_ar_b^2=0\quad,\quad r_ar_b^2+r_{-a}r_{a+b}^2=0$$

(7) In the first case, we are done. In the second case, we know that $a_1=b$ and $b_1=a+b$ are distinct elements of $E$. So, consider the following equation, over $E$:  
$$a_1+2b_1=a_1'+2b_1'$$

The hexagonal solution of this equation, given by $a_1'=-a_1$ and $b_1'=a_1+b_1$, cannot appear, because $b_1'=a_1+b_1$ can be written as $b_1'=a+2b$, and by $2$-maximality we get $b_1'=-a=b$, which contradicts $a+b\in E$. Thus the equation $a_1+2b_1=a_1'+2b_1'$ with $a_1',b_1'\in E$ has only trivial solutions, and with $x=a_1+2b_1$ in the above, we get:
$$r_{a_1}r_{b_1}^2=0$$

Now remember that this follows by identifying coefficients in $\alpha(\xi)\alpha(\xi^2)=\alpha(\xi^3)$. The same method applies to the formula $\alpha(\xi^2)\alpha(\xi)=\alpha(\xi^3)$, and we get as well:
$$r_{b_1}^2r_{a_1}=0$$

We have now all ingredients for finishing the proof of the key formula, as follows:
$$r_ar_b^3
=r_ar_b^2r_b
=-r_{-a}r_{a+b}^2r_b
=-r_{-a}r_{b_1}^2r_{a_1}
=0$$

(8) We come back now to the following formula, proved for $s=1,2,3$:
$$\alpha(\xi^s)=\sum_{a\in E}\xi^{sa}\otimes r_a^s$$

By using the key formula, we get by recurrence on $s$ that
this holds in general.

\medskip

(9) In particular with $s=p-1$ we get the following formula:
$$\alpha(\xi^{-1})=\sum_{a\in E}\xi^{-a}\otimes r_a^{p-1}$$

On the other hand, from $\xi^*=\xi^{-1}$ we get the following formula:
$$\alpha(\xi^{-1})=\sum_{a\in E}\xi^{-a}\otimes r_a^*$$

But this gives $r_a^* = r_a^{p-1}$ for any $a$. Now by using the key
formula we get:
$$(r_ar_b)(r_ar_b)^*
=r_ar_br_b^*r_a^*
=r_ar_b^pr_a^*
=(r_ar_b^3)(r_b^{p-3}r_a^*)
=0$$

(10) But this gives $r_ar_b=0$. Thus, we have the following equalities:
$$r_ar_b=r_br_a=0$$

On the other hand, ${\mathcal A}$ is generated by coefficients of $\alpha$, which are in turn powers of elements $r_a$. It follows that ${\mathcal A}$ is commutative, and we are done.
\end{proof}

Still following \cite{bbg}, we can now formulate a main result, as follows:

\begin{theorem}
A type $k$ circulant graph having $p>>k$ vertices, with $p$ prime,
has no quantum symmetry.
\end{theorem}

\begin{proof}
This follows from Theorem 15.25 and some arithmetics, as follows:

\medskip

(1) Let $k$ be an even number, and consider the group of $k$-th roots of unity $G=\{1,w,\ldots,w^{k-1}\}$, where $w=e^{2\pi i/k}$. By standard arithmetics, $G$ is $2$-maximal in ${\mathbb C}$.

\medskip

(2) As a continuation of this, again by some standard arithmetics, for $p>6^{\varphi(k)}$, with $\varphi$ being the Euler function, any subgroup $E\subset{\mathbb Z}_p^*$ of order $k$ is $2$-maximal.

\medskip

(3) But this proves our result. Indeed, by using (2), we can apply Theorem 15.25 provided that we have $p>6^{{\varphi(k)}}$, and our graph has no quantum symmetry, as desired.
\end{proof}

We should mention that the above result, from \cite{bbg}, is quite old. The challenge is to go beyond this, with results for the graphs having an abelian group action, $A\curvearrowright X$. 

\section*{15d. Partial symmetries}

As a last topic for this chapter, let us discuss the construction and main properties of the quantum semigroup $\widetilde{G}^+(X)\subset\widetilde{S}_N^+$, in analogy with what we know about $\widetilde{G}(X)\subset\widetilde{S}_N$. This is something more recent, and potentially interesting too. Let us start with:

\begin{definition}
$C(\widetilde{S}_N^+)$ is the universal $C^*$-algebra generated by the entries of a $N\times N$ submagic matrix $u$, with comultiplication and counit maps given by
$$\Delta(u_{ij})=\sum_ku_{ik}\otimes u_{kj}\quad,\quad 
\varepsilon(u_{ij})=\delta_{ij}$$
where ``submagic'' means formed of projections, which are pairwise orthogonal on rows and columns. We call $\widetilde{S}_N^+$ semigroup of quantum partial permutations of $\{1,\ldots,N\}$.
\end{definition}

Observe that the morphisms $\Delta,\varepsilon$ constructed above satisfy the usual axioms for a comultiplication and an antipode, in the bialgebra setting, namely:
$$(\Delta\otimes id)\Delta=(id\otimes \Delta)\Delta$$
$$(\varepsilon\otimes id)\Delta=(id\otimes\varepsilon)\Delta=id$$

As a conclusion to this, the basic properties of the quantum semigroup $\widetilde{S}_N^+$ that we constructed in Definition 15.27 can be summarized as follows:

\begin{proposition}
We have maps as follows,
$$\begin{matrix}
C(S_N^+)&\leftarrow&C(\widetilde{S}_N^+)\\
\\
\downarrow&&\downarrow\\
\\
C(S_N)&\leftarrow&C(\widetilde{S}_N)
\end{matrix}
\quad \quad \quad:\quad \quad\quad
\begin{matrix}
S_N^+&\subset&\widetilde{S}_N^+\\
\\
\cup&&\cup\\
\\
S_N&\subset&\widetilde{S}_N
\end{matrix}$$
with the bialgebras at left corresponding to the quantum semigroups at right.
\end{proposition}

\begin{proof}
This is clear from the above discussion, and from the well-known fact that projections which sum up to $1$ are pairwise orthogonal. 
\end{proof}

We recall from chapter 13 that we have $S_N^+\neq S_N$, starting from $N=4$. At the semigroup level things get interesting starting from $N=2$, where we have:

\begin{proposition}
We have an isomorphism as follows,
$$C(\widetilde{S}_2^+)\simeq\left\{(x,y)\in C^*(D_\infty)\oplus C^*(D_\infty)\Big|\varepsilon(x)=\varepsilon(y)\right\}$$
where $\varepsilon:C^*(D_\infty)\to\mathbb C1$ is the counit, given by the formula
$$u=\begin{pmatrix}p\oplus 0&0\oplus r\\0\oplus s&q\oplus 0\end{pmatrix}$$ where $p,q$ and $r,s$ are the standard generators of the two copies of $C^*(D_\infty)$.
\end{proposition}

\begin{proof}
Consider an arbitrary $2\times 2$ matrix formed by projections:
$$u=\begin{pmatrix}P&R\\S&Q\end{pmatrix}$$

This matrix is submagic when the following conditions are satisfied:
$$PR=PS=QR=QS=0$$

Now observe that these conditions tell us that the non-unital algebras $X=<P,Q>$ and $Y=<R,S>$ must commute, and must satisfy $xy=0$, for any $x\in X,y\in Y$. Thus, if we denote by $Z$ the universal algebra generated by two projections, we have:
$$C(\widetilde{S}_2^+)\simeq\mathbb C1\oplus Z\oplus Z$$

Now since we have $C^*(D_\infty)=\mathbb C1\oplus Z$, we obtain an isomorphism as follows:
$$C(\widetilde{S}_2^+)\simeq\left\{(\lambda+a,\lambda+b)\Big|\lambda\in\mathbb C, a,b\in Z\right\}$$

Thus, we are led to the conclusion in the statement.
\end{proof}

We recall now from chapter 9 that given a graph $X$ with $N$ vertices, and adjacency matrix $d\in M_N(0,1)$, its partial automorphism semigroup is given by:
$$\widetilde{G}(X)=\left\{\sigma\in\widetilde{S}_N\Big|d_{ij}=d_{\sigma(i)\sigma(j)},\ \forall i,j\in Dom(\sigma)\right\}$$

We have the following formula, from chapter 11, with $R=diag(R_i)$, $C=diag(C_j)$, with $R_i,C_j$ being the row and column sums of the associated submagic matrix $u$:
$$C(\widetilde{G}(X))=C(\widetilde{S}_N)\Big/\Big<R(du-ud)C=0\Big>$$

With these results in hand, we are led to the following statement:

\begin{theorem}
The following construction, with $R,C$ being the diagonal matrices formed by the row and column sums of $u$, produces a subsemigroup $\widetilde{G}^+(X)\subset\widetilde{S}_N^+$,
$$C(\widetilde{G}^+(X))=C(\widetilde{S}_N^+)\Big/\Big<R(du-ud)C=0\Big>$$
called semigroup of quantum partial automorphisms of $X$, whose classical version is $\widetilde{G}(X)$. When using $du=ud$, we obtain a semigroup $\bar{G}^+(X)\subset\widetilde{G}^+(X)$ which can be smaller.
\end{theorem}

\begin{proof}
All this is elementary, the idea being as follows:

\medskip

(1) In order to construct the comultiplication $\Delta$, consider the following elements:
$$U_{ij}=\sum_ku_{ik}\otimes u_{kj}$$

By using the fact that $u$ is submagic, we deduce that we have:
$$R^U_i(dU)_{ij}C^U_j=\Delta(R_i(du)_{ij}C_j)$$
$$R^U_i(Ud)_{ij}C^U_j=\Delta(R_i(ud)_{ij}C_j)$$

Thus we can define $\Delta$ by mapping $u_{ij}\to U_{ij}$, as desired.

\medskip

(2) Regarding now $\varepsilon$, the algebra in the statement has indeed a morphism $\varepsilon$ defined by $u_{ij}\to\delta_{ij}$, because the following relations are trivially satisfied:
$$R_i(d1_N)_{ij}C_j=R_i(1_Nd)_{ij}C_j$$

(3) Regarding now $S$, we must prove that we have a morphism $S$ given by $u_{ij}\to u_{ji}$. For this purpose, we know that with $R=diag(R_i)$ and $C=diag(C_j)$, we have:
$$R(du-ud)C=0$$

Now when transposing this formula, we obtain:
$$C^t(u^td-du^t)R^t=0$$

Since $C^t,R^t$ are respectively the diagonal matrices formed by the row sums and column sums of $u^t$, we conclude that the relations $R(du-ud)C=0$ are satisfied by the transpose matrix $u^t$, and this gives the existence of the subantipode map $S$.

\medskip

(4) The fact that we have $\widetilde{G}^+(X)_{class}=\widetilde{G}(X)$ follows from $(S_N^+)_{class}=S_N$.

\medskip

(5) Finally, the last assertion follows from our similar results from chapter 7, by taking classical versions, the simplest counterexample being the simplex. 
\end{proof}

As a first result now regarding the correspondence $X\to\widetilde{G}^+(X)$, we have:

\begin{proposition}
For any finite graph $X$ we have
$$\widetilde{G}^+(X)=\widetilde{G}^+(X^c)$$
where $X^c$ is the complementary graph.
\end{proposition}

\begin{proof}
The adjacency matrices of a graph $X$ and of its complement $X^c$ are related by the following formula, where $\mathbb I_N$ is the all-1 matrix:
$$d_X+d_{X^c}=\mathbb I_N-1_N$$

Thus, in order to establish the formula in the statement, we must prove that:
$$R_i(\mathbb I_Nu)_{ij}C_j=R_i(u\mathbb I_N)_{ij}C_j$$

For this purpose, let us recall that, the matrix $u$ being submagic, its row sums and column sums $R_i,C_j$ are projections. By using this fact, we have:
$$R_i(\mathbb I_Nu)_{ij}C_j=R_iC_jC_j=R_iC_j$$
$$R_i(u\mathbb I_N)_{ij}C_j=R_iR_iC_j=R_iC_j$$

Thus we have proved our equality, and the conclusion follows.
\end{proof}

In order to discuss now various aspects of the correspondence $X\to\widetilde{G}^+(X)$, it is technically convenient to slightly enlarge our formalism, as follows:

\begin{definition}
Associated to any complex-colored oriented graph $X$, with adjacency matrix $d\in M_N(\mathbb C)$, is its semigroup of partial automorphisms, given by
$$\widetilde{G}(X)=\left\{\sigma\in\widetilde{S}_N\Big|d_{ij}=d_{\sigma(i)\sigma(j)},\ \forall i,j\in Dom(\sigma)\right\}$$
as well as its quantum semigroup of quantum partial automorphisms, given by 
$$C(\widetilde{G}^+(X))=C(\widetilde{S}_N^+)\Big/\Big<R(du-ud)C=0\Big>$$
where $R=diag(R_i)$, $C=diag(C_j)$, with $R_i,C_j$ being the row and column sums of $u$.
\end{definition}

With this notion in hand, following the material in chapter 14, let us discuss now the color independence. Let $m,\gamma$ be the multiplication and comultiplication of $\mathbb C^N$:
$$m(e_i\otimes e_j)=\delta_{ij}e_i\quad,\quad 
\gamma(e_i)=e_i\otimes e_i$$ 

We denote by $m^{(p)},\gamma^{(p)}$ their iterations, given by the following formulae:
$$m^{(p)}(e_{i_1}\otimes\ldots\otimes e_{i_1})=\delta_{i_1\ldots i_p}e_{i_1}\quad,\quad 
\gamma^{(p)}(e_i)=e_i\otimes\ldots\otimes e_i$$ 

Our goal is to use these iterations in the semigroup case, exactly as we did in chapter 14, in  the quantum group case. We will need some technical results. Let us start with:

\begin{proposition}
We have the following formulae,
$$m^{(p)}u^{\otimes p}=um^{(p)}\quad,\quad 
u^{\otimes p}\gamma^{(p)}=\gamma^{(p)}u$$
valid for any submagic matrix $u$.
\end{proposition}

\begin{proof}
We have the following computations, which prove the first formula:
$$m^{(p)}u^{\otimes p}(e_{i_1}\otimes\ldots\otimes e_{i_p})
=\sum_je_j\otimes u_{ji_1}\ldots u_{ji_p}
=\delta_{i_1\ldots i_p}\sum_je_j\otimes u_{ji_1}$$
$$um^{(p)}(e_{i_1}\otimes\ldots\otimes e_{i_p})
=\delta_{i_1\ldots i_p}u(e_{i_1})
=\delta_{i_1\ldots i_p}\sum_je_j\otimes u_{ji_1}$$

We have as well the following computations, which prove the second formula:
$$u^{\otimes p}\gamma^{(p)}(e_i)
=u^{\otimes p}(e_i\otimes\ldots\otimes e_i)
=\sum_je_j\otimes\ldots\otimes e_j\otimes u_{ji}$$
$$\gamma^{(p)}u(e_i)
=\gamma^{(p)}\left(\sum_je_j\otimes u_{ji}\right)
=\sum_je_j\otimes\ldots\otimes e_j\otimes u_{ji}$$

Summarizing, we have proved both formulae in the statement.
\end{proof}

We will need as well a second technical result, as follows:

\begin{proposition}
We have the following formulae, with $u,m,\gamma$ being as before,
$$m^{(p)}R^{\otimes p}d^{\otimes p}\gamma^{(p)}=Rd^{\times p}\quad,\quad
m^{(p)}d^{\otimes p}C^{\otimes p}\gamma^{(p)}=d^{\times p}C$$
and with $\times$ being the componentwise, or Hadamard, product of matrices.
\end{proposition}

\begin{proof}
We have the following computations, which prove the first formula:
$$m^{(p)}R^{\otimes p}d^{\otimes p}\gamma^{(p)}(e_i)
=m^{(p)}R^{\otimes p}d^{\otimes p}(e_i\otimes\ldots\otimes e_i)
=\sum_je_j\otimes R_jd_{ji}^p$$
$$Rd^{\times p}(e_i)
=R\left(\sum_je_j\otimes d_{ji}^p\right)
=\sum_je_j\otimes R_jd_{ji}^p$$

We have as well the following computations, which prove the second formula:
$$m^{(p)}d^{\otimes p}C^{\otimes p}\gamma^{(p)}(e_i)
=m^{(p)}d^{\otimes p}(e_i\otimes\ldots\otimes e_i\otimes C_i)
=\sum_je_j\otimes d_{ji}^pC_i$$
$$d^{\times p}C(e_i)
=d^{\times p}(e_i\otimes C_i)
=\sum_je_j\otimes d_{ji}^pC_i$$

Thus, we have proved both formulae in the statement.
\end{proof}

We can now prove a key result, as follows:

\begin{proposition}
We have the following formulae, with $u,m,\gamma$ being as before,
$$m^{(p)}(Rdu)^{\otimes p}\gamma^{(p)}=Rd^{\times p}u\quad,\quad 
m^{(p)}(udC)^{\otimes p}\gamma^{(p)}=ud^{\times p}C$$
and with $\times$ being the componentwise product of matrices.
\end{proposition}

\begin{proof}
By using the formulae in Propositions 15.33 and 15.34, we get:
\begin{eqnarray*}
m^{(p)}(Rdu)^{\otimes p}\gamma^{(p)}
&=&m^{(p)}R^{\otimes p}d^{\otimes p}u^{\otimes p}\gamma^{(p)}\\
&=&m^{(p)}R^{\otimes p}d^{\otimes p}\gamma^{(p)}u\\
&=&Rd^{\times p}u
\end{eqnarray*}

Once again by using Proposition 15.33 and Proposition 15.34, we have:
\begin{eqnarray*}
m^{(p)}(udC)^{\otimes p}\gamma^{(p)}
&=&m^{(p)}u^{\otimes p}d^{\otimes p}C^{\otimes p}\gamma^{(p)}\\
&=&um^{(p)}d^{\otimes p}C^{\otimes p}\gamma^{(p)}\\
&=&ud^{\times p}C
\end{eqnarray*}

Thus, we have proved both formulae in the statement.
\end{proof}

We can now prove the color independence result, as follows:

\begin{theorem}
The quantum semigroups of quantum partial isomorphisms of finite graphs are subject to the ``independence on the colors'' formula
$$\Big[d_{ij}=d_{kl}\iff d'_{ij}=d'_{kl}\Big]\implies\widetilde{G}^+(X)=\widetilde{G}^+(X')$$
valid for any graphs $X,X'$, having adjacency matrices $d,d'$.
\end{theorem}

\begin{proof}
Given a matrix $d\in M_N(\mathbb C)$, consider its color decomposition, which is as follows, with the color components $d_c$ being by definition 0-1 matrices:
$$d=\sum_{c\in\mathbb C}c\cdot d_c$$

We want to prove that a given quantum semigroup $G$ acts on $(X,d)$ if and only if it acts on $(X,d_c)$, for any $c\in\mathbb C$. For this purpose, consider the following linear space:
$$E_u=\left\{f\in M_N(\mathbb C)\Big|Rfu=ufC\right\}$$

In terms of this space, we want to prove that we have:
$$d\in E_u\implies d_c\in E_u,\forall c\in\mathbb C$$

For this purpose, observe that we have the following implication, as a consequence of the formulae established in Proposition 15.35:
$$Rdu=udC\implies Rd^{\times p}u=ud^{\times p}C$$

We conclude that we have the following implication:
$$d\in E_u\implies d^{\times p}\in E_u,\forall p\in\mathbb N$$

But this gives the result, exactly as in chapter 14, via the standard fact that the color components $d_c$ can be obtained from the componentwise powers $d^{\times p}$.
\end{proof}

In contrast with what happens for the groups or quantum groups, in the semigroup setting we do not have a spectral decomposition result as well. To be more precise, consider as before the following linear space, associated to a submagic matrix $u$:
$$E_u=\left\{d\in M_N(\mathbb C)\Big|Rdu=udC\right\}$$

It is clear that $E_u$ is a linear space, containing 1, and stable under the adjoint operation $*$ too. We also know from Theorem 15.36 that $E_u$ is stable under color decomposition. However, $E_u$ is not stable under taking products, and so is not an algebra, in general.

\bigskip

In general, the computation of $\widetilde{G}^+(X)$ remains a very interesting question. Interesting as well is the question of generalizing all this to the infinite graph case, $|X|=\infty$, with the key remark that this might be simpler than talking about $G^+(X)$ with $|X|=\infty$.

\section*{15e. Exercises}

We had a quite technical chapter, with a lot of tricky noncommutative algebra, and as exercises we have, and no surprise here, more noncommutative algebra, as follows:

\begin{exercise}
Learn about the Kac-Paljutkin quantum group.
\end{exercise}

\begin{exercise}
Learn about exotic subgroups $S_N\subset G\subset S_N^+$, and their orbitals.
\end{exercise}

\begin{exercise}
Learn more about random graphs, in the classical case.
\end{exercise}

\begin{exercise}
Learn about quantum information aspects of all the above.
\end{exercise}

\begin{exercise}
Clarify the arithmetic input that we used for circulant graphs.
\end{exercise}

\begin{exercise}
Extend our results on circulant graphs, using abelian groups. 
\end{exercise}

\begin{exercise}
Compute some further quantum semigroups of type $\widetilde{G}^+(X)$.
\end{exercise}

\begin{exercise}
Construct a quantum semigroup $\widetilde{G}^+(X)$, when $X$ is infinite.
\end{exercise}

As bonus exercise, read in detail the papers \cite{lmr} and \cite{sc3}. These will bring your knowledge to the state of the art of the subject, as of the late 2010s.

\chapter{Planar algebras}

\section*{16a. Operator algebras}

We discuss here a question that we met several times in this book, in relation with our graph theory investigations, be them classical of quantum, namely the computation of the algebra generated by the adjacency matrix $d\in M_N(0,1)$ of a graph $X$, under the operations consisting in taking the spectral decomposition, and the color decomposition. This question makes in fact sense for any complex matrix, as follows:

\begin{question}
What is the algebra generated by a matrix $d\in M_N(\mathbb C)$,
$$\lhd\, d\,\rhd=?$$
with respect to the spectral decomposition, and the color decomposition?
\end{question}

Here we use the above symbols in the lack of something known and standard, regarding this seemingly alien operation. However, as we will soon discover, there is nothing that alien regarding that operation, which is in fact something very familiar in operator algebras, knot theory, quantum field theory, and many more, namely the operation which consists in computing the associated planar algebra, in the sense of Jones \cite{jo6}.

\bigskip

So, this will be our plan, and the whole discussion will bring us on a long trip into modern mathematics, featuring advanced operator algebras, advanced linear algebra, and other advanced things, quite often with a flavor of modern topology and modern physics too, and as a bonus, we will reach to an answer to the following question too:

\begin{question}
What is the general mathematical theory behind the advanced topics that we saw in Parts I-II, namely walks on ADE graphs, and knot invariants?
\end{question}

With the answer to this latter question being, and you guessed right, again planar algebras in the sense of Jones \cite{jo6}. But probably enough advertisement, let us get to work. We will need some basic von Neumann algebra theory, coming as a complement to the basic $C^*$-algebra theory developed in chapter 13, and we have here:

\index{weak topology}
\index{bicommutant}

\begin{theorem}
For a $*$-algebra of operators $A\subset B(H)$ the following conditions are equivalent, and if satisfied, we say that $A$ is a von Neumann algebra:
\begin{enumerate}
\item $A$ is closed under the weak topology, making each $T\to Tx$ continuous.

\item $A$ equals its bicommutant, $A=A''$, computed inside $B(H)$.
\end{enumerate}
\end{theorem}

\begin{proof}
This is von Neumann's bicommutant theorem, that we actually invoked a few times already, when talking Tannakian duality, in its finite dimensional particular case, which is elementary, with the discussion here, and then the proof, being as follows:

\medskip

(1) As a first comment, the weak topology on $B(H)$, making each $T\to Tx$ with $x\in H$ continuous, is indeed weaker than the norm topology, in the sense that we have:
$$||T_n-T||\to0\implies ||T_nx-Tx||\to0,\forall x\in H$$

In particular, we see that a von Neumann algebra in the sense of (1), that is, closed under the weak topology, must be a $C^*$-algebra, that is, closed under the norm.

\medskip

(2) Before getting further, let us see if the converse of this fact is true. This is certainly true in finite dimensions, $H=\mathbb C^N$, where we have $B(H)=M_N(\mathbb C)$, and where the operator $*$-algebras $A\subset B(H)$ are as follows, automatically closed both for the norm topology, and the weak topology, and with these two topologies actually coinciding:
$$A=M_{n_1}(\mathbb C)\oplus\ldots\oplus M_{n_k}(\mathbb C)$$

(3) In infinite dimensions, however, things change. Indeed, let us first take a look at the most basic examples of commutative $C^*$-algebras that we know, the commutative ones. These naturally appear from compact measured spaces $X$, as follows:
$$C(X)\subset B(L^2(X))\quad,\quad f\to[g\to fg]$$

(4) But, it is pretty much clear that such an algebra will not be weakly closed, unless $X$ is discrete, with the details here being left to you. So, in infinite dimensions, there are far less von Neumann algebras than $C^*$-algebras, with this being good to know.

\medskip

(5) Still talking about this, the following natural question appears, what happens if we take the weak closure of the algebra $C(X)\subset B(L^2(X))$ considered above? And the answer here, obtained via some basic measure theory and functional analysis, that we will leave as an exercise, is that we obtain the following algebra:
$$L^\infty(X)\subset B(L^2(X))\quad,\quad f\to[g\to fg]$$

(6) But this is quite interesting, because forgetting now about $C^*$-algebras, what we have here is a nice method of producing von Neumann algebras, in the weakly closed sense, and with the measured space $X$ being no longer required to be compact.

\medskip

(7) As a conclusion to all this, ``von Neumann algebras have something to do with measure theory, in the same way as $C^*$-algebras have something to do with topology''. Which sounds quite deep, so good, and time to stop here. More on this later.

\medskip

(8) Hang on, we are not done yet with the preliminaries, because all the above was in relation with the condition (1) in the statement, and we still have the condition (2) in the statement to comment on. So, here we go again, with a basic exploration, of that condition. To start with, given a subalgebra $A\subset B(H)$, or even a subset $A\subset B(H)$, we can talk about its commutant inside $B(H)$, constructed as follows:
$$A'=\left\{T\in B(H)\Big|TS=ST,\forall S\in A\right\}$$

Now if we take the commutant $A''$ of this commutant $A'$, it is obvious that the elements of the original algebra or set $A$ will be inside. Thus, we have an inclusion as follows:
$$A\subset A''$$

(9) The question is now, why $A=A''$ should be equivalent to $A$ being weakly closed, and why should we care about this? These are both good questions, so let us start with the first one. As a first observation, in finite dimensions the bicommutant condition is automatic, because with $A\subset M_N(\mathbb C)$ being as in (2) above, its commutant is:
$$A'=\mathbb C\oplus\ldots\oplus \mathbb C$$

But now, by taking again the commutant, we obtain the original algebra $A$:
$$A''=M_{n_1}(\mathbb C)\oplus\ldots\oplus M_{n_k}(\mathbb C)$$

(10) Moving now to infinite dimensions, the first thought goes into taking the commutant of the basic examples of $C^*$-algebras, $C(X)\subset B(L^2(X))$. But here, up to some mesure theory and functional analysis work, that we will leave as an exercise, we are led to the following conclusion, which proves the bicommutant theorem in this case:
$$C(X)''=L^\infty(X)$$

(11) Summarizing, we have some intuition on the condition $A=A''$ from the statement, and we can also say, based on the above, that the method for proving the bicommutant theorem would be that of establishing the following equality, for any $*$-subalgebra $A\subset B(H)$, with on the right being the closure with respect to the weak topology:
$$A''=\overline{A}^{\,w}$$

(12) Before getting to work, however, we still have a question to be answered, namely, why should we care about all this? I mean, the condition (1) in the statement, weak closedness, looks very nice and mathematical, that would be a good axiom for the von Neumann algebras, so why bothering with commutants, and with the condition (2).

\medskip

(13) In answer, at the elementary level, and with my apologies for calling these damn things ``elementary'', we have seen in chapters 12-13, when struggling with Tannakian duality, that the bicommutant operation and theorem can be something very useful.

\medskip

(14) In answer too, at the advanced level now, in abstract quantum mechanics the vectors of the Hilbert space $x\in H$ are the states of the system, and the linear self-adjoint operators $T:H\to H$ are the observables, and taking the commutant of a set or algebra of observables is something extremely natural. And this is how von Neumann came upon such things, back in the 1930s, and looking now retrospectively, we can even say that his bicommutant theorem is not only important in the context of quantum mechanics, but even ``makes abstract quantum mechanics properly work''. So, in short, trust me, with the present bicommutant theorem we are into first-class mathematics and physics.

\medskip

(15) Time perhaps for the proof? We recall from (11) that we would like to prove the following equality, for any $*$-algebra of operators $A\subset B(H)$:
$$A''=\overline{A}^{\,w}$$

(16) Let us first prove $\supset$. Since we have $A\subset A''$, we just have to prove that $A''$ is weakly closed. But, assuming $T_i\to T$ weakly, we have indeed:
\begin{eqnarray*}
T_i\in A''
&\implies&ST_i=T_iS,\ \forall S\in A'\\
&\implies&ST=TS,\ \forall S\in A'\\
&\implies&T\in A
\end{eqnarray*}

(17) Let us prove now $\subset$. Here we must establish the following implication:
$$T\in A''\implies T\in\overline{A}^{\,w}$$

For this purpose, we use an amplification trick. Consider indeed the Hilbert space $K$ obtained by summing $n$ times $H$ with itself:
$$K=H\oplus\ldots\oplus H$$

The operators over $K$ can be regarded as being square matrices with entries in $B(H)$, and in particular, we have a representation $\pi:B(H)\to B(K)$, as follows:
$$\pi(T)=\begin{pmatrix}
T\\
&\ddots\\
&&T
\end{pmatrix}$$

(18) The idea will be that of doing the computations in this representation. First, in this representation, the image of our algebra $A\subset B(H)$ is given by:
$$\pi(A)=\left\{\begin{pmatrix}
T\\
&\ddots\\
&&T
\end{pmatrix}\Big|T\in A\right\}$$

We can compute the commutant of this image, exactly as in the usual scalar matrix case, and we obtain the following formula:
$$\pi(A)'=\left\{\begin{pmatrix}
S_{11}&\ldots&S_{1n}\\
\vdots&&\vdots\\
S_{n1}&\ldots&S_{nn}
\end{pmatrix}\Big|S_{ij}\in A'\right\}$$

(19) We conclude from this that, given an operator $T\in A''$ as above, we have:
$$\begin{pmatrix}
T\\
&\ddots\\
&&T
\end{pmatrix}\in\pi(A)''$$

In other words, the conclusion of all this is that we have:
$$T\in A''\implies \pi(T)\in\pi(A)''$$

(20) Now given a vector $x\in K$, consider the orthogonal projection $P\in B(K)$ on the norm closure of the vector space $\pi(A)x\subset K$. Since the subspace $\pi(A)x\subset K$ is invariant under the action of $\pi(A)$, so is its norm closure inside $K$, and we obtain from this:
$$P\in\pi(A)'$$

By combining this with what we found above, we conclude that we have:
$$T\in A''\implies \pi(T)P=P\pi(T)$$

Now since this holds for any vector $x\in K$, we conclude that any operator $T\in A''$ belongs to the weak closure of $A$. Thus, we have $A''\subset\overline{A}^{\,w}$, as desired.  
\end{proof}

Very nice all this, but as you can see, the von Neumann algebras are far more subtle objects than the $C^*$-algebras, and their proper understanding, even at the very basic level, is a far more complicated business than what we quickly did in chapter 13, for the $C^*$-algebras. Welcome to the real quantum, the quantum mechanics one.

\bigskip

Moving ahead, the continuation of the story involves an accumulation of non-trivial results, due to Murray and von Neumann, from the 1930s and 1940s, and then due to Connes, much later, in the 1970s, the conclusions being as follows:

\index{von Neumann algebra}
\index{factor}
\index{reduction theory}

\begin{theorem}
The von Neumann algebras are as follows:
\begin{enumerate}
\item In the commutative case, these are the algebras $A=L^\infty(X)$, with $X$ measured space, represented on $H=L^2(X)$, up to a multiplicity.

\item If we write the center as $Z(A)=L^\infty(X)$, then we have a decomposition of type $A=\int_XA_x\,dx$, with the fibers $A_x$ having trivial center,  $Z(A_x)=\mathbb C$.

\item The factors, $Z(A)=\mathbb C$, can be fully classified in terms of ${\rm II}_1$ factors, which are those satisfying $\dim A=\infty$, and having a faithful trace $tr:A\to\mathbb C$.
\end{enumerate}
\end{theorem}

\begin{proof}
This is something quite heavy, the idea being as follows:

\medskip

(1) As already discussed above, it is clear that $L^\infty(X)$ is indeed a von Neumann algebra on $H=L^2(X)$. The converse can be proved as well, by using spectral theory, one way of viewing this being by saying that, given a commutative von Neumann algebra $A\subset B(H)$, its elements $T\in A$ are commuting normal operators, so the Spectral Theorem for such operators applies, and gives $A=L^\infty(X)$, for some measured space $X$.

\medskip

(2) This is von Neumann's reduction theory main result, whose statement is already quite hard to understand, and whose proof uses advanced functional analysis. To be more precise, in finite dimensions this is something that we know well, with the formula $A=\int_XA_x\,dx$ corresponding to our usual direct sum decomposition, namely:
$$A=M_{n_1}(\mathbb C)\oplus\ldots\oplus M_{n_k}(\mathbb C)$$

In infinite dimensions, things are more complicated, but the idea remains the same, namely using (1) for the commutative von Neumann algebra $Z(A)$, as to get a measured space $X$, and then making your way towards a decomposition of type $A=\int_XA_x\,dx$.

\medskip

(3) This is something fairly heavy, due to Murray-von Neumann and Connes, the idea being that the other factors can be basically obtained via crossed product constructions. To be more precise, the various type of factors can be classified as follows:

\medskip

-- Type I. These are the matrix algebras $M_N(\mathbb C)$, called of type ${\rm I}_N$, and their infinite generalization, $B(H)$ with $H$ infinite dimensional, called of type ${\rm I}_\infty$. Although these factors are very interesting and difficult mathematical objects, from the perspective of the general von Neumann algebra classification work, they are dismissed as ``trivial''.

\medskip

-- Type II. These are the infinite dimensional factors having a trace, which is a usual trace $tr:A\to\mathbb C$ in the type ${\rm II}_1$ case, and is something more technical, possibly infinite, in the remaining case, the type ${\rm II}_\infty$ one, with these latter factors being of the form $B(H)\otimes A$, with $A$ being a ${\rm II}_1$ factor, and with $H$ being an infinite dimensional Hilbert space.

\medskip

-- Type III. These are the factors which are infinite dimensional, and do not have a trace $tr:A\to\mathbb C$. Murray and von Neumann struggled a lot with such beasts, with even giving an example being a non-trivial task, but later Connes came and classified them, basically showing that they appear from ${\rm II}_1$ factors, via crossed product constructions.
\end{proof}

So long for basic, or rather advanced but foundational, von Neumann algebra theory. In what follows we will focus on the ${\rm II}_1$ factors, according to the following principle:

\index{factor}

\begin{principle}
The building blocks of the von Neumann algebra theory are the ${\rm II}_1$ factors, which are the von Neumann algebras having the following properties:
\begin{enumerate}
\item They are infinite dimensional, $\dim A=\infty$.

\item They are factors, their center being $Z(A)=\mathbb C$.

\item They have a faithful trace $tr:A\to\mathbb C$.
\end{enumerate}
\end{principle}

But you might perhaps ask, is it even clear that such beasts exist? Good point, and in answer, given a discrete group $\Gamma$, you can talk about its von Neumann algebra, obtained by talking the weak closure of the usual group algebra, or group $C^*$-algebra:
$$L(\Gamma)\subset B(l^2(\Gamma))$$

This algebra is then infinite dimensional when $\Gamma$ is infinite, and also has a trace, given on group elements by $tr(g)=\delta_{g1}$. As for the center, this consists of the functions on $\Gamma$ which are constant on the conjugacy classes, so when $\Gamma$ has infinite conjugacy classes, called ICC property, what we have is a factor. Thus, as a conclusion, when $\Gamma$ is infinite and has the ICC property, its von Neumann algebra $L(\Gamma)$ is a ${\rm II}_1$ factor.

\bigskip

Let us summarize this finding, along with a bit more, as follows:

\index{hyperfinite factor}

\begin{theorem}
We have the following examples of ${\rm II}_1$ factors:
\begin{enumerate}
\item The group von Neumann algebras $L(\Gamma)$, with $\Gamma$ being an infinite group, having the infinite conjugacy class (ICC) property.

\item The Murray-von Neumann hyperfinite factor $R=\overline{\cup_kM_{n_k}(\mathbb C)}^{\,w}$, with the limit being independent on the summands, and on the inclusions between them.

\item With the remark that when $\Gamma$ as above is assumed to be amenable, its associated ${\rm II}_1$ factor $L(\Gamma)$ is the Murray-von Neumann hyperfinite factor $R$.
\end{enumerate}
\end{theorem}

\begin{proof}
Here the first assertion comes from the above discussion, and the rest, regarding the factor $R$, is due to Murray and von Neumann, using standard functional analysis. With the remark however that the notion of hyperfiniteness can be plugged into the general considerations from Theorem 16.4, and with the resulting questions, which are of remarkable difficulty, having been solved only relatively recently, basically by Connes in the 1970s, and with a last contribution by Haagerup in the 1980s, the general idea being that, in the end, everything hyperfinite can be reconstructed from $R$.
\end{proof}

Many other things can be said, along these lines, and if truly interested in theoretical physics, be that quantum mechanics, or statistical mechanics, or other, have a look at all this, von Neumann algebras, this is first-class mathematical technology.

\section*{16b. Subfactor theory}

In view of Principle 16.5, and its quantum mechanics ramifications, it looks reasonable to forget about the Hilbert space $H$, about operators $T\in B(H)$, about other von Neumann algebras and factors $A\subset B(H)$ that might appear, about other mathematics and physics too, why not about your friends, spouse and hobbies too, but please keep teaching some calculus, that is first class mathematics, and focus on the ${\rm II}_1$ factors.

\bigskip

With this idea in mind, we have our objects, the ${\rm II}_1$ factors, but what about morphisms. And here, a natural idea is that of looking at the inclusions of such factors:

\begin{definition}
A subfactor is an inclusion of ${\rm II}_1$ factors $A_0\subset A_1$.
\end{definition}

So, these will be the objects that we will be interested in, in what follows. With the comment that, while quantum mechanics and von Neumann algebras have been around for a while, since the 1920s, and Definition 16.7 is something very natural emerging from this, it took mathematics and physics a lot of time to realize this, with Definition 16.7 basically dating back to the late 1970s, with the beginning of the work of Jones, on it. Moral of the story, sometimes it takes a lot of skill, to come up with simple things. 

\bigskip

Now given a subfactor $A_0\subset A_1$, a first question is that of defining its index, measuring how big $A_1$ is, when compared to $A_0$. But this can be done as follows:

\index{factor}
\index{subfactor}
\index{index of subfactor}
\index{coupling constant}

\begin{theorem}
Given a subfactor $A_0\subset A_1$, the number
$$N=\frac{\dim_{A_0}H}{\dim_{A_1}H}$$
is independent of the ambient Hilbert space $H$, and is called index.
\end{theorem}

\begin{proof}
This is something quite standard, the idea being as follows:

\medskip

(1) To start with, given a representation of a ${\rm II}_1$ factor $A\subset B(H)$, we can talk about the corresponding coupling constant, as being a number as follows:
$$\dim_AH\in(0,\infty]$$

To be more precise, we can construct this coupling constant in the following way, with $u:H\to L^2(A)\otimes l^2(\mathbb N)$ being an isometry satisfying $ux=(x\otimes1)u$:
$$\dim_AH=tr(uu^*)$$

(2) Alternatively, we can use the following formula, after proving first that the number on the right is indeed independent of the choice on a nonzero vector $x\in H$:
$$\dim_AH=\frac{tr_A(P_{A'x})}{tr_{A'}(P_{Ax})}$$

This latter formula was in fact the original definition of the coupling constant, by Murray and von Neumann. However, technically speaking, it is better to use (1).

\medskip

(3) Now with this in hand, given a subfactor $A_0\subset A_1$, the fact that the index as defined above is indeed independent of the ambient Hilbert space $H$ comes from the various basic properties of the coupling constant, established by Murray and von Neumann.
\end{proof}

There are many examples of subfactors coming from groups, and every time we obtain the intuitive index. In general now, following Jones \cite{jo1}, let us start with:

\index{conditional expectation}

\begin{proposition}
Given a subfactor $A_0\subset A_1$, there is a unique linear map
$$E:A_1\to A_0$$
which is positive, unital, trace-preserving and which is such that, for any $a_1,a_2\in A_0$:
$$E(a_1ba_2)=a_1E(b)a_2$$
This map is called conditional expectation from $A_1$ onto $A_0$.
\end{proposition}

\begin{proof}
We make use of the standard representation of the ${\rm II}_1$ factor $A_1$, with respect to its unique trace $tr:A_1\to\mathbb C$, namely:
$$A_1\subset L^2(A_1)$$

If we denote by $\Omega$ the standard cyclic and separating vector of $L^2(A_1)$, we have an identification of vector spaces $A_0\Omega=L^2(A_0)$. Consider now the following projection:
$$e:L^2(A_1)\to L^2(A_0)$$

It follows from definitions that we have an inclusion $e(A_1\Omega)\subset A_0\Omega$. Thus the above projection $e$ induces by restriction a certain linear map, as follows:
$$E:A_1\to A_0$$

This linear map $E$ and the orthogonal projection $e$ are related by:
$$exe=E(x)e$$

But this shows that the linear map $E$ satisfies the various conditions in the statement, namely positivity, unitality, trace preservation and bimodule property. As for the uniqueness assertion, this follows by using the same argument, applied backwards, the idea being that a map $E$ as in the statement must come from a projection $e$.
\end{proof}

We will be interested in what follows in the orthogonal projection $e:L^2(A_1)\to L^2(A_0)$ producing the expectation $E:A_1\to A_0$, rather than in $E$ itself:

\index{Jones projection}

\begin{definition}
Associated to any subfactor $A_0\subset A_1$ is the orthogonal projection
$$e:L^2(A_1)\to L^2(A_0)$$
producing the conditional expectation $E:A_1\to A_0$ via the following formula:
$$exe=E(x)e$$
This projection is called Jones projection for the subfactor $A_0\subset A_1$.
\end{definition}

Quite remarkably, the subfactor $A_0\subset A_1$, as well as its commutant, can be recovered from the knowledge of this projection, in the following way:

\begin{proposition}
Given a subfactor $A_0\subset A_1$, with Jones projection $e$, we have
$$A_0=A_1\cap\{e\}'\quad,\quad 
A_0'=(A_1'\cap\{e\})''$$
as equalities of von Neumann algebras, acting on the space $L^2(A_1)$.
\end{proposition}

\begin{proof}
The above two formulae both follow from $exe=E(x)e$, via some elementary computations, and for details here, we refer to Jones' paper \cite{jo1}.
\end{proof}

We are now ready to formulate a key definition, as follows:

\index{basic construction}

\begin{definition}
Associated to any subfactor $A_0\subset A_1$ is the basic construction
$$A_0\subset_eA_1\subset A_2$$
with $A_2=<A_1,e>$ being the algebra generated by $A_1$ and by the Jones projection
$$e:L^2(A_1)\to L^2(A_0)$$
acting on the Hilbert space $L^2(A_1)$.
\end{definition}

The idea now, following as before Jones \cite{jo1}, will be that the inclusion $A_1\subset A_2$ appears as a kind of ``reflection'' of the original inclusion $A_0\subset A_1$, and also that the basic construction can be iterated, with all this leading to non-trivial results. We first have:

\begin{proposition}
Given a subfactor $A_0\subset A_1$ having finite index, 
$$[A_1:A_0]<\infty$$
the basic construction $A_0\subset_eA_1\subset A_2$ has the following properties:
\begin{enumerate}
\item $A_2=JA_0'J$.

\item $A_2=\overline{A_1+A_1eb}$.

\item $A_2$ is a ${\rm II}_1$ factor.

\item $[A_2:A_1]=[A_1:A_0]$.

\item $eA_2e=A_0e$.

\item $tr(e)=[A_1:A_0]^{-1}$.

\item $tr(xe)=tr(x)[A_1:A_0]^{-1}$, for any $x\in A_1$.
\end{enumerate}
\end{proposition}

\begin{proof}
All this is standard, by using the same type of mathematics as in the proof of Proposition 16.9, and for details here, we refer to Jones' paper \cite{jo1}.
\end{proof}

Let us perform now twice the basic construction, and see what we get. The result here, which is something more technical, at least at the first glance, is as follows:

\begin{proposition}
Associated to $A_0\subset A_1$ is the double basic construction
$$A_0\subset_eA_1\subset_fA_2\subset A_3$$
with $e:L^2(A_1)\to L^2(A_0)$ and $f:L^2(A_2)\to L^2(A_1)$ having the following properties:
$$fef=[A_1:A_0]^{-1}f\quad,\quad 
efe=[A_1:A_0]^{-1}e$$
\end{proposition}

\begin{proof}
We have two formulae to be proved, the idea being as follows:

\medskip

(1) The first formula in the statement is clear, because we have:
$$fef
=E(e)f
=tr(e)f
=[A_1:A_0]^{-1}f$$

(2) Regarding now the second formula, it is enough to check this on the dense subset $(A_1+A_1eA_1)\Omega$. Thus, we must show that for any $x,y,z\in A_1$, we have:
$$efe(x+yez)\Omega=[A_1:A_0]^{-1}e(x+yez)\Omega$$

But this is something which is routine as well. See Jones \cite{jo1}.
\end{proof}

We can in fact perform the basic construction by recurrence, and we obtain:

\index{Jones tower}
\index{Jones projections}

\begin{theorem}
Associated to any subfactor $A_0\subset A_1$ is the Jones tower
$$A_0\subset_{e_1}A_1\subset_{e_2}A_2\subset_{e_3}A_3\subset\ldots\ldots$$
with the Jones projections having the following properties:
\begin{enumerate}
\item $e_i^2=e_i=e_i^*$.

\item $e_ie_j=e_je_i$ for $|i-j|\geq2$.

\item $e_ie_{i\pm1}e_i=[A_1:A_0]^{-1}e_i$.

\item $tr(we_{n+1})=[A_1:A_0]^{-1}tr(w)$, for any word $w\in<e_1,\ldots,e_n>$.
\end{enumerate}
\end{theorem}

\begin{proof}
This follows from Proposition 16.13 and Proposition 16.14, because the triple basic construction does not need in fact any further study. See \cite{jo1}.
\end{proof}

The relations found in Theorem 16.15 are in fact well-known, from the standard theory of the Temperley-Lieb algebra. This algebra, discovered by Temperley and Lieb in the context of statistical mechanics \cite{tli}, has a very simple definition, as follows:

\index{Temperley-Lieb algebra}
\index{noncrossing pairings}

\begin{definition}
The Temperley-Lieb algebra of index $N\in[1,\infty)$ is defined as
$$TL_N(k)=span(NC_2(k,k))$$
with product given by vertical concatenation, with the rule
$$\bigcirc=N$$
for the closed circles that might appear when concatenating.
\end{definition}

In other words, the algebra $TL_N(k)$, depending on parameters $k\in\mathbb N$ and $N\in[1,\infty)$, is the linear span of the pairings $\pi\in NC_2(k,k)$. The product operation is obtained by linearity, for the pairings which span $TL_N(k)$ this being the usual vertical concatenation, with the conventions that things go ``from top to bottom'', and that each circle that might appear when concatenating is replaced by a scalar factor, equal to $N$.

\bigskip

In what concerns us, we will just need some elementary results. First, we have:

\index{Jones projection}

\begin{proposition}
The Temperley-Lieb algebra $TL_N(k)$ is generated by the diagrams
$$\varepsilon_1={\ }^\cup_\cap\quad,\quad 
\varepsilon_2=|\!{\ }^\cup_\cap\quad,\quad
\varepsilon_3=||\!{\ }^\cup_\cap\quad,\quad 
\ldots$$
which are all multiples of projections, in the sense that their rescaled versions
$$e_i=N^{-1}\varepsilon_i$$
satisfy the abstract projection relations $e_i^2=e_i=e_i^*$.
\end{proposition}

\begin{proof}
We have two assertions here, the idea being as follows:

\medskip

(1) The fact that the Temperley-Lieb algebra $TL_N(k)$ is indeed generated by the sequence $\varepsilon_1,\varepsilon_2,\ldots$ follows by drawing pictures, and more specifically by decomposing each basis element $\pi\in NC_2(k,k)$ as a product of such elements $\varepsilon_i$.

\medskip

(2) Regarding now the projection assertion, when composing $\varepsilon_i$ with itself we obtain $\varepsilon_i$ itself, times a circle. Thus, according to our multiplication convention, we have:
$$\varepsilon_i^2=N\varepsilon_i$$

Also, when turning upside-down $\varepsilon_i$, we obtain $\varepsilon_i$ itself. Thus, according to our involution convention for the Temperley-Lieb algebra, we have the following formula:
$$\varepsilon_i^*=\varepsilon_i$$

We conclude that the rescalings $e_i=N^{-1}\varepsilon_i$ satisfy $e_i^2=e_i=e_i^*$, as desired.
\end{proof}

As a second result now, making the link with Theorem 16.15, we have:

\begin{proposition}
The standard generators $e_i=N^{-1}\varepsilon_i$ of the Temperley-Lieb algebra $TL_N(k)$ have the following properties, where $tr$ is the trace obtained by closing:
\begin{enumerate}
\item $e_ie_j=e_je_i$ for $|i-j|\geq2$.

\item $e_ie_{i\pm1}e_i=N^{-1}e_i$.

\item $tr(we_{n+1})=N^{-1}tr(w)$, for any word $w\in<e_1,\ldots,e_n>$.
\end{enumerate}
\end{proposition}

\begin{proof}
This follows indeed by doing some elementary computations with diagrams, in the spirit of those performed in the proof of Proposition 16.17.
\end{proof}

With the above results in hand, and still following Jones' paper \cite{jo1}, we can now reformulate Theorem 16.15 into something more conceptual, as follows:

\index{subfactor}
\index{Temperley-Lieb algebra}

\begin{theorem}
Given a subfactor $A_0\subset A_1$, construct its the Jones tower:
$$A_0\subset_{e_1}A_1\subset_{e_2}A_2\subset_{e_3}A_3\subset\ldots\ldots$$
The rescaled sequence of projections $e_1,e_2,e_3,\ldots\in B(H)$ produces then a representation 
$$TL_N\subset B(H)$$
of the Temperley-Lieb algebra of index $N=[A_1:A_0]$.
\end{theorem}

\begin{proof}
We know from Theorem 16.15 that the rescaled sequence of Jones projections $e_1,e_2,e_3,\ldots\in B(H)$ behaves algebrically exactly as the following $TL_N$ diagrams:
$$\varepsilon_1={\ }^\cup_\cap\quad,\quad 
\varepsilon_2=|\!{\ }^\cup_\cap\quad,\quad
\varepsilon_3=||\!{\ }^\cup_\cap\quad,\quad 
\ldots$$

But these diagrams generate $TL_N$, and so we have an embedding $TL_N\subset B(H)$, where $H$ is the Hilbert space where our subfactor $A_0\subset A_1$ lives, as claimed.
\end{proof}

Let us make the following key observation, also from \cite{jo1}:

\index{relative commutant}
\index{higher commutant}
\index{planar algebra}

\begin{theorem} 
Given a finite index subfactor $A_0\subset A_1$, the graded algebra $P=(P_k)$ formed by the sequence of higher relative commutants
$$P_k=A_0'\cap A_k$$
contains the copy of the Temperley-Lieb algebra constructed above, $TL_N\subset P$. This graded algebra $P=(P_k)$ is called ``planar algebra'' of the subfactor.
\end{theorem}

\begin{proof}
As a first observation, since the Jones projection $e_1:A_1\to A_0$ commutes with $A_0$, we have $e_1\in P_2$. By translation we obtain, for any $k\in\mathbb N$:
$$e_1,\ldots,e_{k-1}\in P_k$$

Thus we have indeed an inclusion of graded algebras $TL_N\subset P$, as claimed.
\end{proof}

As an interesting consequence of the above results, also from \cite{jo1}, we have:

\index{index of subfactors}
\index{small index}

\begin{theorem}
The index of subfactors $A\subset B$ is ``quantized'' in the $[1,4]$ range,
$$N\in\left\{4\cos^2\left(\frac{\pi}{n}\right)\Big|n\geq3\right\}\cup[4,\infty]$$
with the obstruction coming from the existence of the representation $TL_N\subset B(H)$.
\end{theorem}

\begin{proof}
This comes from the basic construction, and more specifically from the combinatorics of the Jones projections $e_1,e_2,e_3,\ldots$, the idea being as folows:

\medskip

(1) In order to best comment on what happens, when iterating the basic construction, let us record the first few values of the numbers in the statement:
$$4\cos^2\left(\frac{\pi}{3}\right)=1\quad,\quad 
4\cos^2\left(\frac{\pi}{4}\right)=2$$
$$4\cos^2\left(\frac{\pi}{5}\right)=\frac{3+\sqrt{5}}{2}\quad,\quad 
4\cos^2\left(\frac{\pi}{6}\right)=3$$
$$\ldots$$

(2) When performing a basic construction, we obtain, by trace manipulations on $e_1$:
$$N\notin(1,2)$$

With a double basic construction, we obtain, by trace manipulations on $<e_1,e_2>$:
$$N\notin\left(2,\frac{3+\sqrt{5}}{2}\right)$$

With a triple basic construction, we obtain, by trace manipulations on $<e_1,e_2,e_3>$:
$$N\notin\left(\frac{3+\sqrt{5}}{2},3\right)$$

Thus, we are led to the conclusion in the statement, by a kind of recurrence, involving a certain family of orthogonal polynomials.

\medskip

(3) In practice now, the most elegant way of proving the result is by using the fundamental fact, explained in Theorem 16.19, that that sequence of Jones projections $e_1,e_2,e_3,\ldots\subset B(H)$ generate a copy of the Temperley-Lieb algebra of index $N$:
$$TL_N\subset B(H)$$

With this result in hand, we must prove that such a representation cannot exist in index $N<4$, unless we are in the following special situation:
$$N=4\cos^2\left(\frac{\pi}{n}\right)$$

But this can be proved by using some suitable trace and positivity manipulations on $TL_N$, as in (2) above. For full details here, we refer to \cite{jo1}.
\end{proof}

So long for basic subfactor theory. As a continuation of the story, the subfactors of index $N\leq4$ are classified by the ADE graphs that we met in chapter 3. See \cite{jo5}.

\section*{16c. Planar algebras}

Quite remarkably, the planar algebra structure of $TL_N$, taken in an intuitive sense, of composing diagrams, extends to a planar algebra structure on $P$. In order to discuss this, let us start with axioms for the planar algebras. Following Jones \cite{jo6}, we have:

\index{planar tangle}
\index{planar algebra}

\begin{definition}
The planar algebras are defined as follows:
\begin{enumerate}
\item We consider rectangles in the plane, with the sides parallel to the coordinate axes, and taken up to planar isotopy, and we call such rectangles boxes.

\item A labeled box is a box with $2n$ marked points on its boundary, $n$ on its upper side, and $n$ on its lower side, for some integer $n\in\mathbb N$.

\item A tangle is labeled box, containing a number of labeled boxes, with all marked points, on the big and small boxes, being connected by noncrossing strings.

\item A planar algebra is a sequence of finite dimensional vector spaces $P=(P_n)$, together with linear maps $P_{n_1}\otimes\ldots\otimes P_{n_k}\to P_n$, one for each tangle, such that the gluing of tangles corresponds to the composition of linear maps.
\end{enumerate}
\end{definition}

In this definition we are using rectangles, but everything being up to isotopy, we could have used instead circles with marked points, as in \cite{jo6}. Our choice for using rectangles comes from the main examples that we have in mind, to be discussed below, where the planar algebra structure is best viewed by using rectangles, as above.

\bigskip

Let us also mention that Definition 16.22 is something quite simplified, based on \cite{jo6}. As explained in \cite{jo6}, in order for subfactors to produce planar algebras and vice versa, there are quite a number of supplementary axioms that must be added, and in view of this, it is perhaps better to start with something stronger than Definition 16.22, as basic axioms. However, as before with rectangles vs circles, our axiomatic choices here are mainly motivated by the concrete examples that we have in mind. More on this later.

\bigskip

As a basic example of a planar algebra, we have the Temperley-Lieb algebra:

\begin{theorem}
The Temperley-Lieb algebra $TL_N$, viewed as graded algebra
$$TL_N=(TL_N(n))_{n\in\mathbb N}$$
is a planar algebra, with the corresponding linear maps associated to the planar tangles
$$TL_N(n_1)\otimes\ldots\otimes TL_N(n_k)\to TL_N(n)$$
appearing by putting the various $TL_N(n_i)$ diagrams into the small boxes of the given tangle, which produces a $TL_N(n)$ diagram.
\end{theorem}

\begin{proof}
This is something trivial, which follows from definitions:

\medskip

(1) Assume indeed that we are given a planar tangle $\pi$, as in Definition 16.22, consisting of a box having $2n$ marked points on its boundary, and containing $k$ small boxes, having respectively $2n_1,\ldots,2n_k$ marked points on their boundaries, and then a total of $n+\Sigma n_i$ noncrossing strings, connecting the various $2n+\Sigma 2n_i$ marked points.

\medskip

(2) We want to associate to this tangle $\pi$ a linear map as follows:
$$T_\pi:TL_N(n_1)\otimes\ldots\otimes TL_N(n_k)\to TL_N(n)$$

For this purpose, by linearity, it is enough to construct elements as follows, for any choice of Temperley-Lieb diagrams $\sigma_i\in TL_N(n_i)$, with $i=1,\ldots,k$:
$$T_\pi(\sigma_1\otimes\ldots\otimes\sigma_k)\in TL_N(n)$$

(3) But constructing such an element is obvious, just by putting the various diagrams $\sigma_i\in TL_N(n_i)$ into the small boxes the given tangle $\pi$. Indeed, this procedure produces a certain diagram in $TL_N(n)$, that we can call $T_\pi(\sigma_1\otimes\ldots\otimes\sigma_k)$, as above.

\medskip

(4) Finally, we have to check that everything is well-defined up to planar isotopy, and that the gluing of tangles corresponds to the composition of linear maps. But both these checks are trivial, coming from the definition of $TL_N$, and we are done.
\end{proof}

As a conclusion to all this, $P=TL_N$ is indeed a planar algebra, but of somewhat ``trivial'' type, with the triviality coming from the fact that, in this case, the elements of $P$ are planar diagrams themselves, and so the planar structure appears trivially.

\bigskip

The Temperley-Lieb planar algebra $TL_N$ is however an important planar algebra, because it is the ``smallest'' one, appearing inside the planar algebra of any subfactor. But more on this later, when talking about planar algebras and subfactors.

\bigskip

Moving ahead now, here is our second basic example of a planar algebra, which is also ``trivial'' in the above sense, with the elements of the planar algebra being planar diagrams themselves, but which appears in a bit more complicated way:

\index{Fuss-Catalan algebra}
\index{colored Temperley-Lieb}

\begin{theorem}
The Fuss-Catalan algebra $FC_{N,M}$, which appears by coloring the Temperley-Lieb diagrams with black/white colors, clockwise, as follows 
$$\circ\bullet\bullet\circ\circ\bullet\bullet\circ\ldots\ldots\ldots\circ\bullet\bullet\circ$$
and keeping those diagrams whose strings connect either $\circ-\circ$ or $\bullet-\bullet$, is a planar algebra, with again the corresponding linear maps associated to the planar tangles
$$FC_{N,M}(n_1)\otimes\ldots\otimes FC_{N,M}(n_k)\to FC_{N,M}(n)$$
appearing by putting the various $FC_{N,M}(n_i)$ diagrams into the small boxes of the given tangle, which produces a $FC_{N,M}(n)$ diagram.
\end{theorem}

\begin{proof}
The proof here is nearly identical to the proof of Theorem 16.23, with the only change appearing at the level of the colors. To be more precise:

\medskip

(1) Forgetting about upper and lower sequences of points, which must be joined by strings, a Temperley-Lieb diagram can be thought of as being a collection of strings, say black strings, which compose in the obvious way, with the rule that the value of the circle, which is now a black circle, is $N$. And it is this obvious composition rule that gives the planar algebra structure, as explained in the proof of Theorem 16.23. 

\medskip

(2) Similarly, forgetting about points, a Fuss-Catalan diagram can be thought of as being a collection of strings, which come now in two colors, black and white. These Fuss-Catalan diagrams compose then in the obvious way, with the rule that the value of the black circle is $N$, and the value of the white circle is $M$. And it is this obvious composition rule that gives the planar algebra structure, as before for $TL_N$.
\end{proof}

Getting back now to generalities, and to Definition 16.22, that of a general planar algebra, we have so far two illustrations for it, which, while both important, are both ``trivial'', with the planar structure simply coming from the fact that, in both these cases, the elements of the planar algebra are planar diagrams themselves.

\bigskip

In general, the planar algebras can be more complicated than this, and we will see some further examples in a moment. However, the idea is very simple, namely ``the elements of a planar algebra are not necessarily diagrams, but they behave like diagrams".

\bigskip

In relation now with subfactors, the result, which extends Theorem 16.20, and which was found by Jones in \cite{jo6}, almost 20 years after \cite{jo1}, is as follows:

\index{higher commutant}
\index{planar algebra}

\begin{theorem} 
Given a subfactor $A_0\subset A_1$, the collection $P=(P_n)$ of linear spaces 
$$P_n=A_0'\cap A_n$$
has a planar algebra structure, extending the planar algebra structure of $TL_N$.
\end{theorem}

\begin{proof}
We know from Theorem 16.20 that we have an inclusion as follows, coming from the basic construction, and with $TL_N$ itself being a planar algebra:
$$TL_N\subset P$$

Thus, the whole point is that of proving that the trivial planar algebra structure of $TL_N$ extends into a planar algebra structure of $P$. But this can be done via a long algebraic study, and for the full computation here, we refer to Jones' paper \cite{jo6}.
\end{proof}

As a first illustration for the above result, we have:

\index{Temperley-Lieb}
\index{Fuss-Catalan algebra}

\begin{theorem}
We have the following universality results:
\begin{enumerate}
\item The Temperley-Lieb algebra $TL_N$ appears inside the planar algebra of any subfactor $A_0\subset A_1$ having index $N$.

\item The Fuss-Catalan algebra $FC_{N,M}$ appears inside the planar algebra of any subfactor $A_0\subset A_1$, in the presence of an intermediate subfactor $A_0\subset B\subset A_1$.
\end{enumerate}
\end{theorem}

\begin{proof}
Here the first assertion is something that we already know, from Theorem 16.20, and the second assertion is something quite standard as well, by carefully working out the basic construction for $A_0\subset A_1$, in the presence of an intermediate subfactor $A_0\subset B\subset A_1$. For details here, we refer to Bisch and Jones \cite{bj1}.
\end{proof}

The above results raise the question on whether any planar algebra produces a subfactor. The answer here is yes, but with many subtleties, as follows:

\index{amenable subfactor}
\index{Popa theorem}

\begin{theorem}
We have the following results:
\begin{enumerate}
\item Any planar algebra with positivity produces a subfactor.

\item In particular, we have $TL$ and $FC$ type subfactors.

\item In the amenable case, and with $A_1=R$, the correspondence is bijective.

\item In general, we must take $A_1=L(F_\infty)$, and we do not have bijectivity.

\item The axiomatization of $P$, in the case $A_1=R$, is not known.
\end{enumerate}
\end{theorem}

\begin{proof}
All this is quite heavy, basically coming from the work of Popa in the 90s, using heavy functional analysis, the idea being as follows:

\medskip

(1) As already mentioned in the comments after Definition 16.22, our planar algebra axioms here are something quite simplified, based on \cite{jo6}. However, when getting back to Theorem 16.25, the conclusion is that the subfactor planar algebras there satisfy a number of supplementary ``positivity'' conditions, basically coming from the positivity of the ${\rm II}_1$ factor trace. And the point is that, with these positivity conditions axiomatized, we reach to something which is equivalent to Popa's axiomatization of the lattice of higher relative commutants $A_i'\cap A_j$ of the finite index subfactors, obtained in the 90s via heavy functional analysis. For the full story here, and details, we refer to Jones' paper \cite{jo6}.

\medskip

(2) The existence of the $TL_N$ subfactors, also known as ``$A_\infty$ subfactors'', is something which was known for some time, since some early work of Popa on the subject. As for the existence of the $FC_{N,M}$ subfactors, this can be shown by using the intermediate subfactor picture, $A_0\subset B\subset A_1$, by composing two $A_\infty$ subfactors of suitable indices, $A_0\subset B$ and $B\subset A_1$. For the full story here, we refer as before to Jones \cite{jo6}.

\medskip

(3) This is something fairly heavy, as it is always the case with operator algebra results regarding hyperfiniteness and amenability, due to Popa. For the story here, see \cite{jo6}.

\medskip

(4) This is something a bit more recent, obtained by further building on the above-mentioned constructions of Popa. Again, we refer here to \cite{jo6} and related work.

\medskip

(5) This is the big open question in subfactors. The story here goes back to Jones' original  paper \cite{jo1}, which contains at the end the question, due to Connes, of finding the possible values of the index for the irreducible subfactors of $R$. This question, which certainly looks much easier than (5) in the statement, is in fact still open, now 40 years after its formulation, and with no one having any valuable idea in dealing with it.
\end{proof}

\section*{16d. Graph symmetries} 

Getting back to quantum groups, all this machinery is interesting for us. We will need the construction of the tensor and spin planar algebras $\mathcal T_N,\mathcal S_N$. Let us start with:

\begin{definition}
The tensor planar algebra $\mathcal T_N$ is the sequence of vector spaces 
$$P_k=M_N(\mathbb C)^{\otimes k}$$
with the multilinear maps $T_\pi:P_{k_1}\otimes\ldots\otimes P_{k_r}\to P_k$
being given by the formula
$$T_\pi(e_{i_1}\otimes\ldots\otimes e_{i_r})=\sum_j\delta_\pi(i_1,\ldots,i_r:j)e_j$$
with the Kronecker symbols $\delta_\pi$ being $1$ if the indices fit, and being $0$ otherwise.
\end{definition}

In other words, we put the indices of the basic tensors on the marked points of the small boxes, in the obvious way, and the coefficients of the output tensor are then given by Kronecker symbols, exactly as in the easy quantum group case.

\bigskip

The fact that we have indeed a planar algebra, in the sense that the gluing of tangles corresponds to the composition of linear maps, as required by Definition 16.22, is something elementary, in the same spirit as the verification of the functoriality properties of the correspondence $\pi\to T_\pi$, from easiness, and we refer here to Jones \cite{jo6}. 

\bigskip

Let us discuss now a second planar algebra of the same type, which is important as well for various reasons, namely the spin planar algebra $\mathcal S_N$. This planar algebra appears somehow as the ``square root'' of the tensor planar algebra $\mathcal T_N$. Let us start with:

\begin{definition}
We write the standard basis of $(\mathbb C^N)^{\otimes k}$ in $2\times k$ matrix form,
$$e_{i_1\ldots i_k}=
\begin{pmatrix}i_1 & i_1 &i_2&i_2&i_3&\ldots&\ldots\\
i_k&i_k&i_{k-1}&\ldots&\ldots&\ldots&\ldots 
\end{pmatrix}$$
by duplicating the indices, and then writing them clockwise, starting from top left.
\end{definition}

Now with this convention in hand for the tensors, we can formulate the construction of the spin planar algebra $\mathcal S_N$, also from \cite{jo6}, as follows:

\begin{definition}
The spin planar algebra $\mathcal S_N$ is the sequence of vector spaces 
$$P_k=(\mathbb C^N)^{\otimes k}$$
written as above, with the multiplinear maps $T_\pi:P_{k_1}\otimes\ldots\otimes P_{k_r}\to P_k$ 
being given by
$$T_\pi(e_{i_1}\otimes\ldots\otimes e_{i_r})=\sum_j\delta_\pi(i_1,\ldots,i_r:j)e_j$$
with the Kronecker symbols $\delta_\pi$ being $1$ if the indices fit, and being $0$ otherwise.
\end{definition}

Here are some illustrating examples for the spin planar algebra calculus:

\medskip

(1) The identity $1_k$ is the $(k,k)$-tangle having vertical strings only. The solutions of $\delta_{1_k}(x,y)=1$ being the pairs of the form $(x,x)$, this tangle $1_k$ acts by the identity:
$$1_k\begin{pmatrix}j_1 & \ldots & j_k\\ i_1 & \ldots & i_k\end{pmatrix}=\begin{pmatrix}j_1 & \ldots & j_k\\ i_1 & \ldots & i_k\end{pmatrix}$$

(2) The multiplication $M_k$ is the $(k,k,k)$-tangle having 2 input boxes, one on top of the other, and vertical strings only. It acts in the following way:
$$M_k\left( 
\begin{pmatrix}j_1 & \ldots & j_k\\ i_1 & \ldots & i_k\end{pmatrix}
\otimes\begin{pmatrix}l_1 & \ldots & l_k\\ m_1 & \ldots & m_k\end{pmatrix}
\right)=
\delta_{j_1m_1}\ldots \delta_{j_km_k}
\begin{pmatrix}l_1 & \ldots & l_k\\ i_1 & \ldots & i_k\end{pmatrix}$$

(3) The inclusion $I_k$ is the $(k,k+1)$-tangle which looks like $1_k$, but has one more vertical string, at right of the input box. Given $x$, the solutions of $\delta_{I_k}(x,y)=1$ are the elements $y$ obtained from $x$ by adding to the right a vector of the form $(^l_l)$, and so:
$${I_k}\begin{pmatrix}j_1 & \ldots & j_k\\ i_1 & \ldots & i_k\end{pmatrix}=
\sum_l\begin{pmatrix}j_1 & \ldots & j_k& l\\ i_1 & \ldots & i_k& l\end{pmatrix}$$

(4) The expectation $U_k$ is the $(k+1,k)$-tangle which looks like $1_k$, but has one more string, connecting the extra 2 input points, both at right of the input box:
$$U_k
\begin{pmatrix}j_1 & \ldots &j_k& j_{k+1}\\ i_1 & \ldots &i_k& i_{k+1}\end{pmatrix}=
\delta_{i_{k+1}j_{k+1}}
\begin{pmatrix}j_1 & \ldots & j_k\\ i_1 & \ldots & i_k\end{pmatrix}$$

(5) The Jones projection $E_k$ is a $(0,k+2)$-tangle, having no input box. There are $k$ vertical strings joining the first $k$ upper points to the first $k$ lower points, counting from left to right. The remaining upper 2 points are connected by a semicircle, and the remaining lower 2 points are also connected by a semicircle. We have:
$$E_k(1)=\sum_{ijl}\begin{pmatrix}i_1 & \ldots &i_k&j&j\\ i_1 & \ldots &i_k&l&l\end{pmatrix}$$

The elements $e_k=N^{-1}E_k(1)$ are then projections, and define a representation of the infinite Temperley-Lieb algebra of index $N$ inside the inductive limit algebra $\mathcal S_N$.

\medskip

(6) The rotation $R_k$ is the $(k,k)$-tangle which looks like $1_k$, but the first 2 input points are connected to the last 2 output points, and the same happens at right:
$$R_k=\begin{matrix}
\hskip 0.3mm\Cap \ |\ |\ |\ |\hskip -0.5mm |\cr
|\hskip -0.5mm |\hskip 10.3mm |\hskip -0.5mm |\cr
\hskip -0.3mm|\hskip -0.5mm |\ |\ |\ |\ \hskip -0.1mm\Cup
\end{matrix}$$

The action of $R_k$ on the standard basis is by rotation of the indices, as follows:
$$R_k(e_{i_1i_2\ldots i_k})=e_{i_2\ldots i_ki_1}$$

There are many other interesting examples of $k$-tangles, but in view of our present purposes, we can actually stop here, due to the following fact:

\begin{theorem}
The multiplications, inclusions, expectations, Jones projections and rotations generate the set of all tangles, via the gluing operation.
\end{theorem}

\begin{proof}
This is something well-known and elementary, obtained by ``chopping'' the various planar tangles into small pieces, as in the above list. See \cite{jo6}.
\end{proof}

Finally, in order for our discussion to be complete, we must talk as well about the $*$-structure of the spin planar algebra. This is constructed as follows:
$$\begin{pmatrix}j_1 & \ldots & j_k\\ i_1 & \ldots & i_k\end{pmatrix}^*
=\begin{pmatrix}i_1 & \ldots & i_k\\ j_1 & \ldots & j_k\end{pmatrix}$$

As before, we refer to Jones' paper \cite{jo6} for more on all this. Getting back now to quantum groups, following \cite{ba3}, we have the following result:

\begin{theorem}
Given $G\subset S_N^+$, consider the tensor powers of the associated coaction map on $C(X)$, where $X=\{1,\ldots,N\}$, which are the folowing linear maps:
$$\Phi^k:C(X^k)\to C(X^k)\otimes C(G)$$
$$e_{i_1\ldots i_k}\to\sum_{j_1\ldots j_k}e_{j_1\ldots j_k}\otimes u_{j_1i_1}\ldots u_{j_ki_k}$$
The fixed point spaces of these coactions, which are by definition the spaces
$$P_k=\left\{ x\in C(X^k)\Big|\Phi^k(x)=1\otimes x\right\}$$
are given by $P_k=Fix(u^{\otimes k})$, and form a subalgebra of the spin planar algebra $\mathcal S_N$.
\end{theorem}

\begin{proof}
Since the map $\Phi$ is a coaction, its tensor powers $\Phi^k$ are coactions too, and at the level of fixed point algebras we have the following formula:
$$P_k=Fix(u^{\otimes k})$$

In order to prove now the planar algebra assertion, we will use Theorem 16.31. Consider the rotation $R_k$. Rotating, then applying $\Phi^k$, and rotating backwards by $R_k^{-1}$ is the same as applying $\Phi^k$, then rotating each $k$-fold product of coefficients of $\Phi$. Thus the elements obtained by rotating, then applying $\Phi^k$, or by applying $\Phi^k$, then rotating, differ by a sum of Dirac masses tensored with commutators in $A=C(G)$:
$$\Phi^kR_k(x)-(R_k\otimes id)\Phi^k(x)\in C(X^k)\otimes [A,A]$$

Now let $\int_A$ be the Haar functional of $A$, and consider the conditional expectation onto the fixed point algebra $P_k$, which is given by the following formula:
$$\phi_k=\left(id\otimes\int_A\right)\Phi^k$$

Since $\int_A$ is a trace, it vanishes on commutators. Thus $R_k$ commutes with $\phi_k$:
$$\phi_kR_k=R_k\phi_k$$

The commutation relation $\phi_kT=T\phi_l$ holds in fact for any $(l,k)$-tangle $T$. These tangles are called annular, and the proof is by verification on generators of the annular category. In particular we obtain, for any annular tangle $T$:
$$\phi_kT\phi_l=T\phi_l$$

We conclude from this that the annular category is contained in the suboperad $\mathcal P'\subset\mathcal P$ of the planar operad consisting of tangles $T$ satisfying the following condition, where $\phi =(\phi_k)$, and where $i(.)$ is the number of input boxes:
$$\phi T\phi^{\otimes i(T)}=T\phi^{\otimes i(T)}$$

On the other hand the multiplicativity of $\Phi^k$ gives $M_k\in\mathcal P'$. Now since the planar operad $\mathcal P$ is generated by multiplications and annular tangles, it follows that we have $\mathcal P'=P$. Thus for any tangle $T$ the corresponding multilinear map between spaces $P_k(X)$ restricts to a multilinear map between spaces $P_k$. In other words, the action of the planar operad $\mathcal P$ restricts to $P$, and makes it a subalgebra of $\mathcal S_N$, as claimed.
\end{proof}

As a second result now, also from \cite{ba3}, completing our study, we have:

\begin{theorem}
We have a bijection between quantum permutation groups and subalgebras of the spin planar algebra,
$$(G\subset S_N^+)\quad\longleftrightarrow\quad (Q\subset\mathcal S_N)$$
given in one sense by the construction in Theorem 16.32, and in the other sense by a suitable modification of Tannakian duality.
\end{theorem}

\begin{proof}
The idea is that this will follow by applying Tannakian duality to the annular category over $Q$. Let $n,m$ be positive integers. To any element $T_{n+m}\in Q_{n+m}$ we associate a linear map $L_{nm}(T_{n+m}):P_n(X)\to P_m(X)$ in the following way:
$$L_{nm}\left(\begin{matrix}|\ |\ |\\ T_{n+m}\\ |\ |\ |\end{matrix}\right):
\left(\begin{matrix}|\\ a_n\\ |\end{matrix}\right)
\to \left(\begin{matrix}
\hskip 1.5mm |\hskip 3.0mm |\hskip 3.0mm \cap\\
\ \ T_{n+m}\hskip 0.0mm  |\\
\hskip 1.9mm |\hskip 1.2mm |\hskip 3.2mm |\hskip2.2mm |\\
a_n|\hskip 3.2mm |\hskip 2.2mm |\\
\hskip 2.1mm\cup \hskip3.5mm |\hskip 2.2mm |
\end{matrix}\right)$$

That is, we consider the planar $(n,n+m,m)$-tangle having an small input $n$-box, a big input $n+m$-box and an output $m$-box, with strings as on the picture of the right. This defines a certain multilinear map, as follows:
$$P_n(X)\otimes P_{n+m}(X)\to P_m(X)$$

If we put the element $T_{n+m}$ in the big input box, we obtain in this way a certain linear map $P_n(X)\to P_m(X)$, that we call $L_{nm}$. With this convention, let us set:
$$Q_{nm}=\left\{ L_{nm}(T_{n+m}):P_n(X)\to P_m(X)\Big| T_{n+m}\in Q_{n+m}\right\}$$

These spaces form a Tannakian category, so by \cite{wo2} we obtain a Woronowicz algebra $(A,u)$, such that the following equalities hold, for any $m,n$:
$$Hom(u^{\otimes m},u^{\otimes n})=Q_{mn}$$

We prove that $u$ is a magic unitary. We have $Hom(1,u^{\otimes 2})=Q_{02}=Q_2$, so the unit of $Q_2$ must be a fixed vector of $u^{\otimes 2}$. But $u^{\otimes 2}$ acts on the unit of $Q_2$ as follows:
\begin{eqnarray*}
u^{\otimes 2}(1)
&=&u^{\otimes 2}\left( \sum_i \begin{pmatrix}i&i\\ i&i\end{pmatrix}\right)\\
&=&\sum_{ikl}\begin{pmatrix}k&k\\ l&l\end{pmatrix}\otimes u_{ki}u_{li}\\
&=&\sum_{kl}\begin{pmatrix}k&k\\ l&l\end{pmatrix}\otimes (uu^t)_{kl}
\end{eqnarray*}

From $u^{\otimes 2}(1)=1\otimes 1$ ve get that $uu^t$ is the identity matrix. Together with the unitarity of $u$, this gives the following formulae:
$$u^t=u^*=u^{-1}$$

Consider the Jones projection $E_1\in Q_3$. After isotoping, $L_{21}(E_1)$ looks as follows:
$$L_{21}\left( \Bigl| \begin{matrix}\cup\\\cap\end{matrix}\right) :
\begin{pmatrix} \,|\ |\\ {\ }^i_j{\ }^i_j\\ \,|\ |\end{pmatrix}\,\to\,
\begin{pmatrix}\hskip -5.8mm |\\ {\ }^i_j{\ }^i_j\supset\\ \hskip -5.8mm |\end{pmatrix}
=\,\delta_{ij}\begin{pmatrix}\,|\\ {\ }^i_i\\ \,|\end{pmatrix}$$

In other words, the linear map $M=L_{21}(E_1)$ is the multiplication $\delta_i\otimes\delta_j\to\delta_{ij}\delta_i$:
$$M\begin{pmatrix}i&i\\ j&j\end{pmatrix}
=\delta_{ij}\begin{pmatrix}i\\ i\end{pmatrix}$$

In order to finish, consider the following element of $C(X)\otimes A$:
$$(M\otimes id)u^{\otimes 2}\left(\begin{pmatrix}i&i\\ j&j\end{pmatrix}\otimes 1\right)
=\sum_k\begin{pmatrix}k\\ k\end{pmatrix}\delta_k\otimes u_{ki}u_{kj}$$

Since $M\in Q_{21}=Hom(u^{\otimes 2},u)$, this equals the following element of $C(X)\otimes A$:
$$u(M\otimes id)\left(\begin{pmatrix}i&i\\ j&j\end{pmatrix}\otimes 1\right)
=\sum_k\begin{pmatrix}k\\ k\end{pmatrix}\delta_k\otimes\delta_{ij}u_{ki}$$

Thus we have $u_{ki}u_{kj}=\delta_{ij}u_{ki}$ for any $i,j,k$, which shows that $u$ is a magic unitary. Now if $P$ is the planar algebra associated to $u$, we have $Hom(1,v^{\otimes n})=P_n=Q_n$, as desired. As for the uniqueness, this is clear from the Peter-Weyl theory.
\end{proof}

All the above might seem a bit technical, but is worth learning, and for good reason, because it is extremely powerful. As an example of application, if you agree with the bijection $G\leftrightarrow Q$ in Theorem 16.33, then $G=S_N^+$ itself, which is the biggest object on the left, must correspond to the smallest object on the right, namely $Q=TL_N$. 

\bigskip

Back now to our usual business, graphs, we have the following result:

\begin{theorem}
The planar algebra associated to $G^+(X)$ is equal to the planar algebra generated by $d$, viewed as a $2$-box in the spin planar algebra $\mathcal S_N$, with $N=|X|$.
\end{theorem}

\begin{proof}
We recall from the above that any quantum permutation group $G\subset S_N^+$ produces a subalgebra $P\subset\mathcal S_N$ of the spin planar algebra, given by:
$$P_k=Fix(u^{\otimes k})$$

In our case, the idea is that $G=G^+(X)$ comes via the relation $d\in End(u)$, but we can view this relation, via Frobenius duality, as a relation of the following type:
$$\xi_d\in Fix(u^{\otimes 2})$$

Indeed, let us view the adjacency matrix $d\in M_N(0,1)$ as a 2-box in $\mathcal S_N$, by using the canonical identification between $M_N(\mathbb C)$ and the algebra of 2-boxes $\mathcal S_N(2)$:
$$(d_{ij})\leftrightarrow \sum_{ij} d_{ij}\begin{pmatrix}i&i\\ j&j\end{pmatrix}$$

Let $P$ be the planar algebra associated to $G^+(X)$ and let $Q$ be the planar algebra generated by $d$. The action of $u^{\otimes 2}$ on $d$ viewed as a 2-box is given by:
$$u^{\otimes 2}\left(\sum_{ij} d_{ij}\begin{pmatrix}i&i\\ j&j\end{pmatrix}\right)
=\sum_{ijkl} d_{ij}\begin{pmatrix}k&k\\ l&l\end{pmatrix}\otimes u_{ki}u_{lj}
=\sum_{kl}\begin{pmatrix}k&k\\ l&l\end{pmatrix}\otimes (udu^t)_{kl}$$

Since $v$ is a magic unitary commuting with $d$ we have:
$$udu^t=duu^t=d$$

But this means that $d$, viewed as a 2-box, is in the algebra $P_2$ of fixed points of $u^{\otimes 2}$. Thus $Q\subset P$. As for $P\subset Q$, this follows from the duality found above.
\end{proof}

Generally speaking, the above material, when coupled with what we did in this book about graphs, leads us into the classification of the subalgebras of the spin planar algebra generated by a 2-box. But this can be regarded as a particular case of the Bisch-Jones question of classifying, in general, the planar algebras generated by a 2-box \cite{bj2}. 

\bigskip

There are many more things that can be said here, notably with the introduction of a general notion of ``quantum graph'', which is something quite interesting on its own, and which brings the whole graph problematics closer to the level of generality of \cite{bj2}. For more on all this, we refer to the book \cite{ba3}, which is more specialized.

\bigskip

Finally, let us mention that, as already explained since the beginning of this book, graphs are not everything, in relation with discrete mathematics. Discrete mathematics is something far wider than graph theory, with all sorts of interesting objects involved, which are not necessarily graphs. But the main principles that we learned here, namely that things are usually encoded by a matrix, and that symmetries and quantum symmetries play a key role, generally apply. For more on all this, other aspects of discrete mathematics, and related algebraic techniques, you can have a look at \cite{ba2}, \cite{dfl}, \cite{sti}.

\section*{16e. Exercises}

Congratulations for having read this book, and no exercises for this final chapter. However, for further reading, you have many possible books, on graphs and related topics. We have referenced some below, and in the hope that you will like some of them.

\baselineskip=14pt

\printindex


\begin{thebibliography}{99}

\baselineskip=14pt

\bibitem{ar1}V.I. Arnold, Ordinary differential equations, Springer (1973).

\bibitem{ar2}V.I. Arnold, Mathematical methods of classical mechanics, Springer (1974).

\bibitem{ar3}V.I. Arnold, Lectures on partial differential equations, Springer (1997).

\bibitem{ati}M.F. Atiyah, The geometry and physics of knots, Cambridge Univ. Press (1990).

\bibitem{ba1}T. Banica, Linear algebra and group theory (2024).

\bibitem{ba2}T. Banica, Invitation to Hadamard matrices (2024).

\bibitem{ba3}T. Banica, Quantum permutation groups (2024).

\bibitem{bbi}T. Banica and J. Bichon, Quantum automorphism groups of vertex-transitive graphs of order $\leq11$, {\em J. Algebraic Combin.} {\bf 26} (2007), 83--105.

\bibitem{bbg}T. Banica, J. Bichon and G. Chenevier, Graphs having no quantum symmetry, {\em Ann. Inst. Fourier} {\bf 57} (2007), 955--971.

\bibitem{bbc}T. Banica, J. Bichon and B. Collins, The hyperoctahedral quantum group, {\em J. Ramanujan Math. Soc.} {\bf 22} (2007), 345--384.

\bibitem{bdb}T. Banica and D. Bisch, Spectral measures of small index principal graphs, {\em Comm. Math. Phys.} {\bf 269} (2007), 259--281.

\bibitem{bpa}H. Bercovici and V. Pata, Stable laws and domains of attraction in free probability theory, {\em Ann. of Math.} {\bf 149} (1999), 1023--1060.

\bibitem{bi1}J. Bichon, Free wreath product by the quantum permutation group, {\em Alg. Rep. Theory} {\bf 7} (2004), 343--362.

\bibitem{bi2}J. Bichon, Algebraic quantum permutation groups, {\em Asian-Eur. J. Math.} {\bf 1} (2008), 1--13.

\bibitem{bj1}D. Bisch and V.F.R. Jones, Algebras associated to intermediate subfactors, {\em Invent. Math.} {\bf 128} (1997), 89--157.

\bibitem{bj2}D. Bisch and V.F.R. Jones, Singly generated planar algebras of small dimension, {\em Duke Math. J.} {\bf 104} (2000), 41--75.

\bibitem{bo1}B. Bollob\'as, Modern graph theory, Springer (1998).

\bibitem{bo2}B. Bollob\'as, Random graphs, Cambridge Univ. Press (1985).

\bibitem{bo3}B. Bollob\'as, Extremal graph theory, Dover (1978).

\bibitem{bce}M. Brannan, A. Chirvasitu, K. Eifler, S. Harris, V. Paulsen, X. Su and M. Wasilewski, Bigalois extensions and the graph isomorphism game, {\em Comm. Math.  Phys.} {\bf 375} (2020), 1777--1809.

\bibitem{bra}R. Brauer, On algebras which are connected with the semisimple continuous groups, {\em Ann. of Math.} {\bf 38} (1937), 857--872.

\bibitem{bha}A. Brouwer and W.H. Haemers, Spectra of graphs, Springer (2011).

\bibitem{cha}A. Chassaniol, Quantum automorphism group of the lexicographic product of finite regular graphs, {\em J. Algebra} {\bf 456} (2016), 23--45.

\bibitem{csn}B. Collins and P. \'Sniady, Integration with respect to the Haar measure on unitary, orthogonal and symplectic groups, {\em Comm. Math. Phys.} {\bf 264} (2006), 773--795.

\bibitem{con}A. Connes, Noncommutative geometry, Academic Press (1994).

\bibitem{cox}H.S.M. Coxeter, Regular polytopes, Dover (1948).

\bibitem{crs}D.M. Cvetkovi\'c, P. Rowlinson and S. Simi\'c, An introduction to the theory of graph spectra, Cambridge Univ. Press (2019).

\bibitem{dfl}W. de Launey and D. Flannery, Algebraic design theory, AMS (2011).

\bibitem{dif}P. Di Francesco, Meander determinants, {\em Comm. Math. Phys.} {\bf 191} (1998), 543--583.

\bibitem{do1}M.P. do Carmo, Differential geometry of curves and surfaces, Dover (1976).

\bibitem{do2}M.P. do Carmo, Riemannian geometry, Birkh\"auser (1992).

\bibitem{dou}S.T. Dougherty, Combinatorics and finite geometry, Springer (2020).

\bibitem{du1}R. Durrett, Probability: theory and examples, Cambridge Univ. Press (1990).

\bibitem{du2}R. Durrett, Random graph dynamics, Cambridge Univ. Press (2006).

\bibitem{eyn}B. Eynard, Counting surfaces, Birkh\"auser (2016).

\bibitem{fel}W. Feller, An introduction to probability theory and its applications, Wiley (1950).

\bibitem{fe1}R.P. Feynman, R.B. Leighton and M. Sands, The Feynman lectures on physics I: mainly mechanics, radiation and heat, Caltech (1963).

\bibitem{fe2}R.P. Feynman, R.B. Leighton and M. Sands, The Feynman lectures on physics II: mainly electromagnetism and matter, Caltech (1964).

\bibitem{fe3}R.P. Feynman, R.B. Leighton and M. Sands, The Feynman lectures on physics III: quantum mechanics, Caltech (1966).

\bibitem{fpi}P. Fima and L. Pittau, The free wreath product of a compact quantum group by a quantum automorphism group, {\em J. Funct. Anal.} {\bf 271} (2016), 1996--2043.

\bibitem{ggr}C. Godsil and G. Royle, Algebraic graph theory, Springer (2001).

\bibitem{gos}D. Goswami, Quantum group of isometries in classical and  noncommutative geometry, {\em Comm. Math. Phys.} {\bf 285} (2009), 141--160.

\bibitem{gr1}D.J. Griffiths, Introduction to electrodynamics, Cambridge Univ. Press (2017).

\bibitem{gr2}D.J. Griffiths and D.F. Schroeter, Introduction to quantum mechanics, Cambridge Univ. Press (2018).

\bibitem{gr3}D.J. Griffiths, Introduction to elementary particles, Wiley (2020).

\bibitem{dg1}D. Gromada, Some examples of quantum graphs, {\em Lett. Math. Phys.} {\bf 112} (2022), 112--122.

\bibitem{dg2}D. Gromada, Quantum symmetries of Cayley graphs of abelian groups, {\em Glasg. Math. J.} {\bf 65} (2023), 655--686.

\bibitem{gtu}J.L. Gross and T.W. Tucker, Topics in topological graph theory, Dover (1987).

\bibitem{har}J. Harris, Algebraic geometry, Springer (1992).

\bibitem{jo1}V.F.R. Jones, Index for subfactors, {\em Invent. Math.} {\bf 72} (1983), 1--25.

\bibitem{jo2}V.F.R. Jones, A polynomial invariant for knots via von Neumann algebras, {\em Bull. Amer. Math. Soc.} {\bf 12} (1985), 103--111.

\bibitem{jo3}V.F.R. Jones, Hecke algebra representations of braid groups and link polynomials, {\em Ann. of Math.} {\bf 126} (1987), 335--388.

\bibitem{jo4}V.F.R. Jones, On knot invariants related to some statistical mechanical models, {\em Pacific J. Math.} {\bf 137} (1989), 311--334.

\bibitem{jo5}V.F.R. Jones, Subfactors and knots, AMS (1991).

\bibitem{jo6}V.F.R. Jones, Planar algebras I (1999).

\bibitem{jo7}V.F.R. Jones, The annular structure of subfactors, {\em Monogr. Enseign. Math.} {\bf 38} (2001), 401--463.

\bibitem{joz}P. J\'oziak, Quantum increasing sequences generate quantum permutation groups, {\em Glasg. Math. J.} {\bf 62} (2020), 631--629.

\bibitem{jsm}L. Junk, S. Schmidt and M. Weber, Almost all trees have quantum symmetry, {\em Arch. Math.} {\bf 115} (2020), 267--278.

\bibitem{kes}H. Kesten, Symmetric random walks on groups, {\em Trans. Amer. Math. Soc.} {\bf 92} (1959), 336--354.

\bibitem{lzv}S.K. Lando and A.K. Zvonkin, Graphs on surfaces and their applications, Springer (2004).

\bibitem{lan}S. Lang, Algebra, Addison-Wesley (1993).

\bibitem{lax}P. Lax, Functional analysis, Wiley (2002).

\bibitem{lta}F. Lemeux and P. Tarrago, Free wreath product quantum groups: the monoidal category, approximation properties and free probability, {\em J. Funct. Anal.} {\bf  270} (2016), 3828--3883.

\bibitem{lin}B. Lindst\"om, Determinants on semilattices, {\em Proc. Amer. Math. Soc.} {\bf 20} (1969), 207--208.

\bibitem{lmr}M. Lupini, L. Man\v cinska and D.E. Roberson, Nonlocal games and quantum permutation groups, {\em J. Funct. Anal.} {\bf 279} (2020), 1--39.

\bibitem{mal}S. Malacarne, Woronowicz's Tannaka-Krein duality and free orthogonal quantum groups, {\em Math. Scand.} {\bf 122} (2018), 151--160.

\bibitem{mro}L. Man\v cinska and D.E. Roberson, Quantum isomorphism is equivalent to equality of homomorphism counts from planar graphs (2019).

\bibitem{mpa}V.A. Marchenko and L.A. Pastur, Distribution of eigenvalues in certain sets of random matrices, {\em Mat. Sb.} {\bf 72} (1967), 507--536.

\bibitem{mcc}J.P. McCarthy, A state-space approach to quantum permutations, {\em Exposition. Math.} {\bf 40} (2022), 628--664.

\bibitem{meh}M.L. Mehta, Random matrices, Elsevier (2004).

\bibitem{mth}B. Mohar and C. Thomassen, Graphs on surfaces, Johns Hopkins Univ. Press (2001).

\bibitem{mrv}B. Musto, D.J. Reutter and D. Verdon, A compositional approach to quantum functions, {\em J. Math. Phys.} {\bf 59} (2018), 1--57.

\bibitem{rwe}S. Raum and M. Weber, The full classification of orthogonal easy quantum groups, {\em Comm. Math. Phys.} {\bf 341} (2016), 751--779.

\bibitem{rsh}L.B. Richmond and J. Shallit, Counting abelian squares, {\em Electron. J. Combin.} {\bf 16} (2009), 1--9.

\bibitem{rs1}D.E. Roberson and S. Schmidt, Quantum symmetry vs nonlocal symmetry (2020).

\bibitem{rs2}D.E. Roberson and S. Schmidt, Solution group representations as quantum symmetries of graphs, {\em J. Lond. Math. Soc.} {\bf 106} (2022), 3379--3410.

\bibitem{ru1}W. Rudin, Principles of mathematical analysis, McGraw-Hill (1964).

\bibitem{ru2}W. Rudin, Real and complex analysis, McGraw-Hill (1966).

\bibitem{sc1}S. Schmidt, The Petersen graph has no quantum symmetry, {\em Bull. Lond. Math. Soc.} {\bf 50} (2018), 395--400.

\bibitem{sc2}S. Schmidt, Quantum automorphisms of folded cube graphs, {\em Ann. Inst. Fourier} {\bf 70} (2020), 949--970.

\bibitem{sc3}S. Schmidt, On the quantum symmetry groups of distance-transitive graphs, {\em Adv. Math.} {\bf 368} (2020), 1--43.

\bibitem{ser}J.P. Serre, Linear representations of finite groups, Springer (1977).

\bibitem{sha}I.R. Shafarevich, Basic algebraic geometry, Springer (1974).

\bibitem{sto}G.C. Shephard and J.A. Todd, Finite unitary reflection groups, {\em Canad. J. Math.} {\bf 6} (1954), 274--304.

\bibitem{sti}D.R. Stinson, Combinatorial designs: constructions and analysis, Springer (2006).

\bibitem{twa}P. Tarrago and J. Wahl, Free wreath product quantum groups and standard invariants of subfactors, {\em Adv. Math.} {\bf 331} (2018), 1--57.

\bibitem{twe}P. Tarrago and M. Weber, Unitary easy quantum groups: the free case and the group case, {\em Int. Math. Res. Not.} {\bf 18} (2017), 5710--5750.

\bibitem{tli}N.H. Temperley and E.H. Lieb, Relations between the ``percolation'' and ``colouring'' problem and other graph-theoretical problems associated with regular planar lattices: some exact results for the ``percolation'' problem, {\em Proc. Roy. Soc. London} {\bf 322} (1971), 251--280.

\bibitem{tru}R.J. Trudeau, Introduction to graph theory, Dover (1993).

\bibitem{vdn}D.V. Voiculescu, K.J. Dykema and A. Nica, Free random variables, AMS (1992).

\bibitem{wan}S. Wang, Quantum symmetry groups of finite spaces, {\em Comm. Math. Phys.} {\bf 195} (1998), 195--211.

\bibitem{we1}S. Weinberg, Foundations of modern physics, Cambridge Univ. Press (2011).

\bibitem{we2}S. Weinberg, Lectures on quantum mechanics, Cambridge Univ. Press (2012).

\bibitem{we3}S. Weinberg, Lectures on astrophysics, Cambridge Univ. Press (2019).

\bibitem{wei}D. Weingarten, Asymptotic behavior of group integrals in the limit of infinite rank, {\em J. Math. Phys.} {\bf 19} (1978), 999--1001.

\bibitem{wey}H. Weyl, The classical groups: their invariants and representations, Princeton Univ. Press (1939).

\bibitem{wig}E. Wigner, Characteristic vectors of bordered matrices with infinite dimensions, {\em Ann. of Math.} {\bf 62} (1955), 548--564. 

\bibitem{wit}E. Witten, Quantum field theory and the Jones polynomial, {\em Comm. Math. Phys.} {\bf 121} (1989), 351--399.

\bibitem{wo1}S.L. Woronowicz, Compact matrix pseudogroups, {\em Comm. Math. Phys.} {\bf 111} (1987), 613--665.

\bibitem{wo2}S.L. Woronowicz, Tannaka-Krein duality for compact matrix pseudogroups. Twisted SU(N) groups, {\em Invent. Math.} {\bf 93} (1988), 35--76.

\end{thebibliography}
\end{document}